\documentclass[b5paper,oneandahalfspaced,palatino,11pt,normalmargins]{thesis-phdNSITunito_me}

\usepackage{babel}
\usepackage{amsmath,amssymb,amsthm,enumerate,array}
\usepackage{caption}
\usepackage{subfig}
\usetikzlibrary{arrows,shapes,matrix,decorations.pathmorphing,decorations.pathreplacing}
\usepackage{mathrsfs}
\usepackage{verbatim}

\usepackage[bookmarks,colorlinks,hyperindex,linkcolor=blue,citecolor=red,menucolor=blue,urlcolor=green,hyperfigures]{hyperref}
\usepackage[all]{hypcap}

 \def\fileversion{0.2}
 \def\filedate{2006/01/25}
\edef\DTAtCode{\the\catcode`\@}
\catcode`\@=11
\long\def\LOOPDIRTREE#1\REPEATDIRTREE{%
  \def\ITERATE{#1\relax\expandafter\ITERATE\fi}%
  \ITERATE
  \let\ITERATE\relax
}
\let\REPEATDIRTREE=\fi
\expandafter\ifx\csname DT@fromsty\endcsname\relax
  \def\@namedef#1{\expandafter\def\csname #1\endcsname}
  \def\@nameuse#1{\csname #1\endcsname}
  \long\def\@gobble#1{}
\fi
\def\@nameedef#1{\expandafter\edef\csname #1\endcsname}
\newdimen\DT@offset \DT@offset=0.2em
\newdimen\DT@width \DT@width=1em
\newdimen\DT@sep \DT@sep=0.2em
\newdimen\DT@all
\DT@all=\DT@offset
\advance\DT@all \DT@width
\advance\DT@all \DT@sep
\newdimen\DT@rulewidth \DT@rulewidth=0.4pt
\newdimen\DT@dotwidth \DT@dotwidth=1.6pt
\newdimen\DTbaselineskip \DTbaselineskip=\baselineskip
\newcount\DT@counti
\newcount\DT@countii
\newcount\DT@countiii
\newcount\DT@countiv
\def\DTsetlength#1#2#3#4#5{%
  \DT@offset=#1\relax
  \DT@width=#2\relax
  \DT@sep=#3\relax
  \DT@all=\DT@offset
  \advance\DT@all by\DT@width
  \advance\DT@all by\DT@sep
  \DT@rulewidth=#4\relax
  \DT@dotwidth=#5\relax
}
\expandafter\ifx\csname DT@fromsty\endcsname\relax
  \def\DTstyle{\tt}
  \def\DTstylecomment{\rm}
\else
  \def\DTstyle{\ttfamily}
  \def\DTstylecomment{\rmfamily}
\fi
\def\DTcomment#1{%
  \kern\parindent\dotfill
  {\DTstylecomment{#1}}%
}
\def\dirtree#1{%
  \let\DT@indent=\parindent
  \parindent=\z@
  \let\DT@parskip=\parskip
  \parskip=\z@
  \let\DT@baselineskip=\baselineskip
  \baselineskip=\DTbaselineskip
  \let\DT@strut=\strut
  \def\strut{\vrule width\z@ height0.7\baselineskip depth0.3\baselineskip}%
  \DT@counti=\z@
  \let\next\DT@readarg
  \next#1\@nil
  \dimen\z@=\hsize
  \advance\dimen\z@ -\DT@offset
  \advance\dimen\z@ -\DT@width
  \setbox\z@=\hbox to\dimen\z@{%
    \hsize=\dimen\z@
    \vbox{\@nameuse{DT@body@1}}%
  }%
  \dimen\z@=\ht\z@
  \advance\dimen0 by\dp\z@
  \advance\dimen0 by-0.7\baselineskip
  \ht\z@=0.7\baselineskip
  \dp\z@=\dimen\z@
  \par\leavevmode
  \kern\DT@offset
  \kern\DT@width
  \box\z@
  \endgraf
  \DT@countii=\@ne
  \DT@countiii=\z@
  \dimen3=\dimen\z@
  \@namedef{DT@lastlevel@1}{-0.7\baselineskip}%
  \loop
  \ifnum\DT@countii<\DT@counti
    \advance\DT@countii \@ne
    \advance\DT@countiii \@ne
    \dimen\z@=\@nameuse{DT@level@\the\DT@countii}\DT@all
    \advance\dimen\z@ by\DT@offset
    \advance\dimen\z@ by-\DT@all
    \leavevmode
    \kern\dimen\z@
    \DT@countiv=\DT@countii
    \count@=\z@
    \LOOPDIRTREE
      \advance\DT@countiv \m@ne
      \ifnum\@nameuse{DT@level@\the\DT@countiv} >
        \@nameuse{DT@level@\the\DT@countii}\relax
      \else
        \count@=\@ne
      \fi
    \ifnum\count@=\z@
    \REPEATDIRTREE
    \edef\DT@hsize{\the\hsize}%
    \count@=\@nameuse{DT@level@\the\DT@countii}\relax
    \dimen\z@=\count@\DT@all
    \advance\hsize by-\dimen\z@
    \setbox\z@=\vbox{\@nameuse{DT@body@\the\DT@countii}}%
    \hsize=\DT@hsize
    \dimen\z@=\ht\z@
    \advance\dimen\z@ by\dp\z@
    \advance\dimen\z@ by-0.7\baselineskip
    \ht\z@=0.7\baselineskip
    \dp\z@=\dimen\z@
    \@nameedef{DT@lastlevel@\the\DT@countii}{\the\dimen3}%
    \advance\dimen3 by\dimen\z@
    \advance\dimen3 by0.7\baselineskip
    \dimen\z@=\@nameuse{DT@lastlevel@\the\DT@countii}\relax
    \advance\dimen\z@ by-\@nameuse{DT@lastlevel@\the\DT@countiv}\relax
    \advance\dimen\z@ by0.3\baselineskip
    \ifnum\@nameuse{DT@level@\the\DT@countiv} <
        \@nameuse{DT@level@\the\DT@countii}\relax
      \advance\dimen\z@ by-0.5\baselineskip
    \fi
    \kern-0.5\DT@rulewidth
    \hbox{\vbox to\z@{\vss\hrule width\DT@rulewidth height\dimen\z@}}%
    \kern-0.5\DT@rulewidth
    \kern-0.5\DT@dotwidth
    \vrule width\DT@dotwidth height0.5\DT@dotwidth depth0.5\DT@dotwidth
    \kern-0.5\DT@dotwidth
    \vrule width\DT@width height0.5\DT@rulewidth depth0.5\DT@rulewidth
    \kern\DT@sep
    \box\z@
    \endgraf
  \repeat
  \parindent=\DT@indent
  \parskip=\DT@parskip
  \DT@baselineskip=\baselineskip
  \let\strut\DT@strut
}
\def\DT@readarg.#1 #2. #3\@nil{%
  \advance\DT@counti \@ne
  \@namedef{DT@level@\the\DT@counti}{#1}%
  \@namedef{DT@body@\the\DT@counti}{\strut{\DTstyle{#2}\strut}}%
  \ifx\relax#3\relax
    \let\next\@gobble
  \fi
  \next#3\@nil
}
\catcode`\@=\DTAtCode\relax
 
\renewcommand*\DTstyle{\rm}
\makeatletter \def\revddots{\mathinner{\mkern1mu\raise\p@
\vbox{\kern7\p@\hbox{.}}\mkern2mu
\raise4\p@\hbox{.}\mkern2mu\raise7\p@\hbox{.}\mkern1mu}} \makeatother

\usepackage{makeidx}  

\newcommand{\K}{\mathbb{K}}
\newcommand{\Kalg}{(\K\textnormal{-algebras})}
\newcommand{\Sets}{(\textnormal{sets})}
\newcommand{\Fun}[2]{\textnormal{Fun}\big({#1},{#2}\big)}
\newcommand{\Ob}{\textnormal{Ob}}
\newcommand{\Mor}{\textnormal{Mor}}
\newcommand{\schm}{(\textnormal{schemes})}
\newcommand{\affschm}{(\textnormal{schemes})^{\textnormal{Aff}}}
\newcommand{\Hom}{\textnormal{Hom}}
\newcommand{\id}{\textnormal{id}}
\newcommand{\hilb}[2]{\mathcal{H}\textnormal{ilb}^{#1}_{#2}}
\newcommand{\PP}{\mathbb{P}}
\newcommand{\QQ}{\mathbb{Q}}
\newcommand{\ZZ}{\mathbb{Z}}
\newcommand{\NN}{\mathbb{N}}
\renewcommand{\AA}{\mathbb{A}}
\newcommand{\Hilb}[2]{\textnormal{\bf Hilb}^{#1}_{#2}}
\newcommand{\grass}[2]{\mathcal{G}\textnormal{r}^{#2}_{#1}}
\newcommand{\Grass}[3]{\textnormal{\bf Gr}_{#3}({#1},{#2})}
\newcommand{\Proj}{\textnormal{Proj}\,}
\newcommand{\Spec}{\textnormal{Spec}\,}
\newcommand{\ind}[1]{\mathrm{#1}}
\newcommand{\ii}{\ind{i}}
\newcommand{\jj}{\ind{j}}
\newcommand{\hh}{\ind{h}}
\newcommand{\kk}{\ind{k}}
\newcommand{\II}{\ind{I}}
\newcommand{\JJ}{\ind{J}}
\newcommand{\HH}{\ind{H}}
\newcommand{\KK}{\ind{K}}
\newcommand{\vv}{\underline{v}}
\newcommand{\ww}{\underline{w}}
\newcommand{\xx}{\underline{x}}
\newcommand{\dv}{\overline{v}}

\newcommand{\PE}{\mathscr{P}}
\newcommand{\GL}{\textnormal{GL}}

\newcommand{\BGL}{\textnormal{B}^{\GL}}
\newcommand{\TGL}{\textnormal{T}^{\GL}}
\newcommand{\OO}{\mathcal{O}}
\newcommand{\rank}{\textnormal{rank}\,}
\newcommand{\IN}{\textnormal{\textsf{in}}}
\newcommand{\GIN}{\textnormal{\textsf{gin}}}
\newcommand{\sat}{\textnormal{sat}}
\newcommand{\reg}{\textnormal{reg}}
\newcommand{\dec}[4]{\langle{#1}\vert{#2}\rangle^{#3}_{#4}}
\newcommand{\syz}{\textnormal{Syz}}
\newcommand{\up}[1]{\textnormal{\textsf{e}}^{+}_{#1}}
\newcommand{\down}[1]{\textnormal{\textsf{e}}^{-}_{#1}}
\newcommand{\pos}[2]{\mathcal{P}({#1},{#2})}
\newcommand{\infpos}[1]{\mathcal{P}({#1})}
\newcommand{\restrict}[2]{#1_{(\geqslant{#2})}}

\newcommand{\HP}[1]{\textnormal{\textsc{hp}}({#1})}
\newcommand{\Borel}[2]{{\mathcal{B}^{#1}_{#2}}}
\newcommand{\NOB}[2]{{\mathcal{N}^{#1}_{#2}}}
\newcommand{\DegLex}{\textnormal{\texttt{DegLex}}}
\newcommand{\RevLex}{\textnormal{\texttt{DegRevLex}}}
\newcommand{\Lex}{\textnormal{\texttt{Lex}}}
\newcommand{\resField}[1]{k(\mathfrak{#1})}
\newcommand{\tail}[3]{\mathcal{T}_{#2}^{#3}(#1)}
\newcommand{\St}[1]{\mathcal{S}\textit{t}_{#1}}
\newcommand{\ed}{\textnormal{ed}\,}
\newcommand{\Supp}{\textnormal{Supp}\,}
\newcommand{\poset}{\mathcal{P}(n,r)}

\newcommand{\Ht}{\textnormal{Ht}}

\newcommand{\Nf}{\textnormal{Nf}}
\newcommand{\cN}{\mathcal{N}}
\newcommand{\Mf}{\mathcal{M}\textit{f}}
\newcommand{\lcm}{\textnormal{lcm}}
\newcommand{\Coeff}{\textnormal{Coeff}}
\newcommand{\lcmHilb}[2]{\textnormal{\bf H}_{{#1},{#2}}}
\newcommand{\smHilb}[2]{\textnormal{\bf H}^{{\tiny\textnormal{sm}}}_{{#1},{#2}}}

\newcommand{\lcmdown}[1]{\underline{\textnormal{\textsf{e}}}^{-}_{#1}}
\newcommand{\maxvar}{\textnormal{max\,var}\,}
\newcommand{\dCoef}{\genfrac{\langle}{\rangle}{0pt}0}
\newcommand{\tCoef}{\genfrac{\langle}{\rangle}{0pt}1}

\newtheorem{theorem}{Theorem}[chapter]
\newtheorem{corollary}[theorem]{Corollary}
\newtheorem{proposition}[theorem]{Proposition}
\newtheorem{lemma}[theorem]{Lemma}

\theoremstyle{definition}
\newtheorem{definition}[theorem]{Definition}
\newtheorem{criterion}[theorem]{Criterion}

\newtheorem{example}{Example}[section]
\newtheorem{remark}[example]{Remark}

\newenvironment{code}{\scriptsize\begin{singlespaced}}{\end{singlespaced}\normalsize\noindent}

\makeindex

\usepackage[sort=use,toc,numberline,style=long3colheader]{glossaries}

\renewcommand{\glossaryname}{Symbols and Notation}

\newglossaryentry{K}{name={\ensuremath{\K}},description={algebraically closed field of characteristic 0}}
\newglossaryentry{KX}{name={\ensuremath{\K[x]}},description={compact notation for the polynomial ring $\K[x_0,\ldots,x_n]$}}
\newglossaryentry{catSets}{name={\ensuremath{\Sets}},description={the category of sets}}
\newglossaryentry{catKalg}{name={\ensuremath{\Kalg}},description={the category of $\K$-algebras}}
\newglossaryentry{catSchmK}{name={\ensuremath{\schm_{\K}}},description={the category of schemes over $\K$}}
\newglossaryentry{catAffSchmK}{name={\ensuremath{\affschm_{\K}}},description={the category of affine schemes over $\K$}}
\newglossaryentry{catFun}{name={\ensuremath{\Fun{\mathcal{C}}{\mathcal{C}'}}},description={the category of functors between the categories $\mathcal{C}$ and $\mathcal{C}'$}}
\newglossaryentry{functorPoints}{name={\ensuremath{h_X}},description={the functor of points of the scheme $X$}}
\newglossaryentry{obCat}{name={\ensuremath{\Ob_\mathcal{C}}},description={the set of objects of the category $\mathcal{C}$}}
\newglossaryentry{morCat}{name={\ensuremath{\Mor_\mathcal{C}}},description={the set of morphisms between objects of the category $\mathcal{C}$}}
\newglossaryentry{extAlg}{name={\ensuremath{\wedge^q V}},description={exterior algebra of the vector space $V$ of order $q$}}
\newglossaryentry{KDelta}{name={\ensuremath{\K[\Delta]}},description={compact notation for the polynomial ring of Pl\"ucker coordinates $\K[\ldots,\Delta_\II,\ldots]$}}
\newglossaryentry{PE}{name={\ensuremath{\PE}},description={Pl\"ucker embedding of the Grassmannian}}
\newglossaryentry{grassFun}{name={\ensuremath{\grass{q}{N}}},description={Grassmann functor}}
\newglossaryentry{genSymb}{name={\ensuremath{\langle S\rangle}},description={if $S$ is a set of vectors, $\langle S\rangle$ denotes the vector space generetad by $S$, whereas if $S$ is a set of polynomials, $\langle S\rangle$ is the ideal generated by such polynomials}}
\newglossaryentry{ZariskiClosed}{name={\ensuremath{Z(I)}},description={the Zariski closed set defined by the ideal $I$}}
\newglossaryentry{transposed}{name={\ensuremath{{}^T M}},description={the transposed matrix of $M$}}
\newglossaryentry{hilbFun}{name={\ensuremath{\hilb{n}{}}},description={Hilbert functor}}
\newglossaryentry{hilbFunPoly}{name={\ensuremath{\hilb{n}{p(t)}}},description={the subfunctor of the Hilbert functor associated to the Hilbert polynomial $p(t)$}}
\newglossaryentry{hilbScheme}{name={\ensuremath{\Hilb{n}{p(t)}}},description={the Hilbert scheme, representing the functor $\hilb{n}{p(t)}$}}
\newglossaryentry{satIdeal}{name={\ensuremath{I(X)}},description={the saturated ideal defining the scheme $X$}}
\newglossaryentry{initIdeal}{name={\ensuremath{\IN_{\sigma}(I)}},description={the initial ideal of the ideal $I$ with respect to the term order $\sigma$}}
\newglossaryentry{gin}{name={\ensuremath{\GIN_{\sigma}(I)}},description={the generic initial ideal of the ideal $I$ with respect to the term order $\sigma$}}
\newglossaryentry{GL}{name={\ensuremath{\GL_{\K}(n+1)}},description={the group of invertible matrices}}
\newglossaryentry{BGL}{name={\ensuremath{\BGL_{\K}(n+1)}},description={the Borel subgroup of $\GL_\K(n+1)$ of the upper triangular matrices}}
\newglossaryentry{TGL}{name={\ensuremath{\TGL_{\K}(n+1)}},description={the torus subgroup of $\GL_\K(n+1)$ of the diagonal matrices}}
\newglossaryentry{saturation}{name={\ensuremath{I^\sat}},description={the saturation of the ideal $I$}}
\newglossaryentry{minMon}{name={\ensuremath{\min x^\alpha}},description={the smallest variable (or its index) dividing $x^\alpha$}}
\newglossaryentry{maxMon}{name={\ensuremath{\max x^\alpha}},description={the greatest variable (or its index) dividing $x^\alpha$}}
\newglossaryentry{borelDec}{name={\ensuremath{\dec{x^{\alpha_i}}{x^\gamma}{I}{}}},description={the canonical decomposition of the monomial $x^\beta = x^{\alpha_i} x^\gamma \in I$, namely $x^\alpha_i$ is a minimal generator of the Borel-fixed ideal $I$ and $\min x^{\alpha_1} \geq \max x^\gamma$}}
\newglossaryentry{decFunc}{name={\ensuremath{\partial_{I}}},description={the decomposition function of the Borel-fixed ideal $I$. For each monomial $x^\beta \in I$, $\partial_I$ maps $x^\beta$ to the minimal generator $x^{\alpha_i}$ of $I$ appearing in its canonical decomposition, i.e. $x^\beta = \dec{x^{\alpha_i}}{x^\gamma}{I}{}$}}
\newglossaryentry{DeltaI}{name={\ensuremath{\Delta^i p(t)}},description={the $i$-th difference of the polynomial $p(t)$: $\Delta^i p(t) = \Delta^{i-1} p(t) - \Delta^{i-1} p(t-1)$}}
\newglossaryentry{Sigma}{name={\ensuremath{\Sigma p(t)}},description={the inverse operation of $\Delta$ on polynomials, namely $\Delta (\Sigma p) (t) = p(t)$ and $\Sigma (\Delta p) (t) = p(t) - c,\ c\geqslant 0$}}
\newglossaryentry{upMove}{name={\ensuremath{\up{i}}},description={the element $\frac{x_{i+1}}{x_i}$ in the field of fractions $\K(x_0,\ldots,x_n)$}}
\newglossaryentry{downMove}{name={\ensuremath{\down{j}}},description={the element $\frac{x_{j-1}}{x_j}$ in the field of fractions $\K(x_0,\ldots,x_n)$}}
\newglossaryentry{BorelOrder}{name={\ensuremath{\leq_B}},description={the Borel partial order}}
\newglossaryentry{pos}{name={\ensuremath{\pos{n}{m}}},description={the poset of monomials in $\K[x_0,\ldots,x_n]_m$ w.r.t. the Borel order}}
\newglossaryentry{BorelSet}{name={\ensuremath{\mathscr{B}}},description={Borel set}}
\newglossaryentry{BorelSetIdeal}{name={\ensuremath{\{I_m\}}},description={the Borel set defined by the monomials in $I_m$}}
\newglossaryentry{orderSet}{name={\ensuremath{\mathscr{N}}},description={order set, the complement of a Borel set}}
\newglossaryentry{restrict}{name={\ensuremath{\restrict{S}{i}}},description={the subset of $S$ containing the monomials with minimum variable greater than or equal to $x_i$}}
\newglossaryentry{BorelNpoly}{name={\ensuremath{\Borel{n}{p(t)}}},description={the set of saturated Borel-fixed ideals of $\K[x_0,\ldots,x_n]$ with Hilbert polynomial $p(t)$}}
\newglossaryentry{NOB}{name={\ensuremath{\NOB{n}{p(t)}}},description={the number of saturated Borel-fixed ideals of $\K[x_0,\ldots,x_n]$ with Hilbert polynomial $p(t)$}}
\newglossaryentry{maxvar}{name={\ensuremath{\maxvar \mathscr{N}}},description={the greatest variable (or its index) dividing a monomial belonging to $\mathscr{N}$}}
\newglossaryentry{NOBsequence}{name={\ensuremath{\NOB{\bullet}{p(t)}}},description={the sequence of the number of saturated Borel-fixed ideals with Hilbert polynomial $p(t)$ in polynomials ring with increasing number of variables}}
\newglossaryentry{DeltaNOB}{name={\ensuremath{\Delta\NOB{i}{s}}},description={the difference between the number of saturated Borel-fixed ideals with Hilbert polynomial $s$ in the polynomial rings $\K[x_0,\ldots,x_{s-i}]$ and $\K[x_0,\ldots,x_{s-i-1}]$}}
\newglossaryentry{DeltaNOBsequence}{name={\ensuremath{\Delta\NOB{\bullet}{s}}},description={the sequence of the differences $\Delta\NOB{i}{s}$ for $i$ varying from $1$ to $s-2$}}
\newglossaryentry{DeltaNOBstable}{name={\ensuremath{\Delta \NOB{i}{}}},description={the constant value of $\Delta \NOB{i}{s}$ for $s \gg 0$}}
\newglossaryentry{DeltaNOBstableSequence}{name={\ensuremath{\Delta \NOB{\bullet}{}}},description={the sequence of integers $\Delta \NOB{i}{}$ for increasing value of $i$}}
\newglossaryentry{tail}{name={\ensuremath{\tail{x^\alpha}{\sigma}{J}}},description={the set of monomials of the same degree of $x^\alpha$ not belonging to the ideal $J$, smaller than $x^\alpha$ w.r.t. the term ordering $\sigma$}}
\newglossaryentry{GBstratum}{name={\ensuremath{\St{\sigma}(J)}},description={the family of homogeneous ideals having $J$ as initial ideal w.r.t. the term ordering $\sigma$ (also the affine scheme describing the family)}}
\newglossaryentry{syz}{name={\ensuremath{\syz(I)}},description={the module of syzygies of the ideal $I$ (also a set of generators of the module)}}
\newglossaryentry{Spoly}{name={\ensuremath{S(f_i,f_j)}},description={the $S$-polynomial between $f_i$ and $f_j$: i.e. $\frac{\IN_\sigma(f_j)}{\textnormal{gcd}(\IN_\sigma(f_i),\IN_\sigma(f_j))}f_i - \frac{\IN_\sigma(f_i)}{\textnormal{gcd}(\IN_\sigma(f_i),\IN_\sigma(f_j))}f_j$ if a term order $\sigma$ is consider, $\frac{\Ht(f_j)}{\textnormal{gcd}(\Ht(f_i),\Ht(f_j))}f_i - \frac{\Ht(f_i)}{\textnormal{gcd}(\Ht(f_i),\Ht(f_j))}f_j$ if $f_i$ and $f_j$ are marked polynomials}}
\newglossaryentry{monomialGroup}{name={\ensuremath{\overline{\mathbb{T}}_x}},description={the multiplicative group of monomials in the fiels of fractions $\K(x)$}}
\newglossaryentry{ed}{name={\ensuremath{\ed \St{\sigma}(J,T)}},description={the embedding dimension of the Gr\"obner stratum $\St{\sigma}(J,T)$}}
\newglossaryentry{supp}{name={\ensuremath{\Supp X}},description={the support of the subscheme $X$}}
\newglossaryentry{sousEscalier}{name={\ensuremath{\mathcal{N}(I)}},description={the sous-escalier of the ideal $I$, that is the set of monomials not belonging to $I$}}
\newglossaryentry{saturateUnderline}{name={\ensuremath{\underline{J}}},description={the saturation of a Borel-fixed ideal $J$}}
\newglossaryentry{suppPoly}{name={\ensuremath{\Supp f}},description={the support of a polynomial $f \in \K[x]$, i.e. the set of monomials appearing in $f$ with non-zero coefficient}}
\newglossaryentry{HT}{name={\ensuremath{\Ht(f)}},description={the specified monomial in $\Supp f$ that makes $f$ a marked polynomial}}
\newglossaryentry{tailMP}{name={\ensuremath{\mathcal{T}(f)}},description={the tail of a marked polynomial $f$, i.e. $\tail{f}{}{} = f - \Ht(f)$}}
\newglossaryentry{markedFamily}{name={\ensuremath{\Mf(J)}},description={the family of ideals $I$ such that the sous-escalier $\cN(J)$ of a Borel-fixed ideal $J$ is a basis as $\K$-vector space of $\K[x]/I$}}
\newglossaryentry{lcm}{name={\ensuremath{\lcm}},description={least common multiple}}
\newglossaryentry{BorelN}{name={\ensuremath{\mathcal{B}^n}},description={the set of (saturated) Borel-fixed ideals of $\K[x_0,\ldots,x_n]$}}
\newglossaryentry{MonIdealN}{name={\ensuremath{\mathcal{M}^n}},description={the set of (saturated) monomial ideals of $\K[x_0,\ldots,x_n]$}}
\newglossaryentry{lcmHilb}{name={\ensuremath{\lcmHilb{d}{g}}},description={the Hilbert scheme of locally Cohen-Macaulay curves in $\PP^3$ of degree $d$ and genus $g$ (contained in $\Hilb{3}{dt+1-g}$)}}
\newglossaryentry{Lex}{name={\ensuremath{\Lex}},description={the lexicographic term order (not graded)}}
\newglossaryentry{DegLex}{name={\ensuremath{\DegLex}},description={the degree lexicographic term order}}
\newglossaryentry{DegRevLex}{name={\ensuremath{\RevLex}},description={the degree reverse lexicographic term order}}

\newglossary{hilbSchemes}{hlb}{hbl}{Hilbert Schemes}
\newglossaryentry{Hilb2P2}{type=hilbSchemes,name={\ensuremath{\Hilb{2}{2}}},description={2 points in the projective plane $\PP^2$}}
\newglossaryentry{Hilb3tp1P3}{type=hilbSchemes,name={\ensuremath{\Hilb{3}{3t+1}}},description={curves of degree 3 and genus 0 in the projective space $\PP^3$ (containing rational normal curves of degree 3)}}
\newglossaryentry{Hilb6tm5P3}{type=hilbSchemes,name={\ensuremath{\Hilb{3}{6t-5}}},description={curves of degree 6 and genus 6 in the projective space $\PP^3$}}
\newglossaryentry{Hilb8P3}{type=hilbSchemes,name={\ensuremath{\Hilb{3}{8}}},description={8 points in the projective space $\PP^3$ (first example of Hilbert scheme of points with reducible components)}}
\newglossaryentry{Hilb4tp1P4}{type=hilbSchemes,name={\ensuremath{\Hilb{4}{4t+1}}},description={curves of degree 4 and genus 0 in the projective space $\PP^4$ (containing rational normal curves of degree 4)}}
\newglossaryentry{Hilb6tm3P3}{type=hilbSchemes,name={\ensuremath{\Hilb{3}{6t-3}}},description={curves of degree 6 and genus 4 in the projective space $\PP^3$ (containing $(2,3)$-complete intersections)}}
\newglossaryentry{Hilb4tP3}{type=hilbSchemes,name={\ensuremath{\Hilb{3}{4t}}},description={curves of degree 4 and genus 1 in the projective space $\PP^3$ (containing $(2,2)$-complete intersections)}}
\newglossaryentry{Hilb3tP3}{type=hilbSchemes,name={\ensuremath{\Hilb{3}{3t}}},description={curves of degree 3 and genus 1 in the projective space $\PP^3$ (containing $(1,3)$-complete intersections)}}
\newglossaryentry{lcmHilb4m3}{type=hilbSchemes,name={\ensuremath{\lcmHilb{4}{-3}}},description={locally Cohen-Macaulay curves of degree 4 and genus $-3$ in the projective space $\PP^3$ (containing disjoint unions of 4 lines on a smooth quadric surface)}}
\newglossaryentry{Hilb7P2}{type=hilbSchemes,name={\ensuremath{\Hilb{2}{7}}},description={7 points in the projective plane $\PP^2$}}

\makeglossaries

\candidate{Paolo Lella}

\title{Computable Hilbert schemes}

\cycle{XXIV}

\department{Matematica}

\academicyears{2009--2011}

\submitdate{Febbraio 2012}

\tutor{Prof.ssa Margherita Roggero}

\coordinator{Prof. Luigi Rodino}

\scientificfield{Geometria MAT/03}

\area{Matematica}

\graphicalabstract{graphicalAbstract.jpg}

\bibfiles{bibFiles/allBib}

\begin{document}

\coverpage

\titlepage

\frontmatter

\begin{dedication}
A mia nonna\qquad
\end{dedication}

\begin{acknowledgements}
Voglio innanzitutto ringraziare la mia relatrice Prof.ssa Margherita Roggero, per avermi seguito in questi tre anni condividendo com me le maggior parte delle sue idee. Ho sempre apprezzato il suo approccio alla ricerca, perch\'e vedo intatto la magia del \lq\lq porsi problemi per poi trovarne la soluzione\rq\rq\ che ha sempre rappresentato per me una fortissima attrazione. 

Inoltre ringrazio tutti i professori e ricercatori con i quali ho avuto il piacere di collaborare e confrontarmi: continuo a rimanere piacevolemente stupito dalla disponibilit\`a con la quale tutti loro si sono confrontati con me, non facendomi sentire mai in soggezione e ascoltando con attenzione ed interesse le mie idee inesperte, mettendo anche in discussione le loro ben pi\`u mature. In particolare, grazie a Maria Grazia Marinari, Francesca Cioffi, Roberto Notari, Enrico Schlesinger, Alberto Albano.

Un\hfill ruolo\hfill di\hfill primo\hfill piano\hfill in\hfill questo\hfill percorso\hfill \`e\hfill stato\hfill sicuramente\hfill ricoperto\\ dall'\lq\lq Ufficio Dottorandi\rq\rq. Non si pu\`o certo dire che sia il miglior posto per raggiungere il massimo della concentrazione, ma \`e sicuramente il miglior posto dove passare i tre anni del dottorato. \`E sempre stato un piacere andare in ufficio e vivere le giornate in un clima amichevole, condividendo gli entusiasmi, confrontandosi sulle difficolt\`a del presente e sulle paure per il futuro e discutendo anche di matematica (grossone docet). Grazie a Davide, Nico e Ube con i quali ho condiviso interamente questa esperienza, ed a tutti quelli che si sono aggiunti negli anni successivi. Un grazie anche a Roberta, Enrico e Marco che da quell'ufficio ci sono passati: rappresentano per me un punto di riferimento e mi fanno guardare con ottimismo al futuro.

Un grazie enorme a Giulia. \`E rassicurante sapere che la sera ci sar\`a lei a casa ed \`e bello svegliarsi la mattina con la colazione pronta. Grazie ai miei genitori per il supporto e pi\`u in generale per l'educazione che mi ha portato ad essere quello che sono oggi. Grazie a mio fratello: fin da piccolo ho sempre ammirato e un p\`o invidiato il suo maggiore talento, al quale ho sempre cercato di sopperire con la determinazione e l'impegno, qualit\`a sulle quali ancora oggi faccio affidamento.

Infine grazie ai miei compagni di squadra, per l'abnegazione con la quale insieme competiamo settimanalemente, non solo sul campo da pallavolo. Alzare la coppa \`e stato un grande onore ed un ricordo che mi accompagner\`a sempre.
\end{acknowledgements}

\tableofcontents

\listoffigures
\listofalgorithms

\mainmatter

\begin{introduction}

The Hilbert scheme $\Hilb{n}{p(t)}$ has been introduced by Grothendieck \cite{GrothendieckFGA} at the beginning of the \rq 60 and belongs to several objects that arose with the schematic re-interpretation of algebraic geometry. It represents the Hilbert functor the associates to any scheme $Z$, over a ground field \gls{K} of characteristic 0, the set of flat families in a projective space $\PP^n$ parametrized by $Z$. For this reason, usually we say that the Hilbert scheme parametrizes all the subschemes and all the (flat) families of subschemes of $\PP^n$ with a fixed Hilbert polynomial $p(t)$. This means that the Hilbert scheme is itself a parameter scheme of a flat family $\mathcal{X}$ of subschemes of $\PP^n$ with Hilbert polynomial $p(t)$ such that any other flat family $\mathcal{Y} \rightarrow S$ can be seen as pullback of $\mathcal{X} \rightarrow \Hilb{n}{p(t)}$ by means of a uniquely defined map $S \rightarrow \Hilb{n}{p(t)}$:
\begin{center}
\begin{tikzpicture}[scale=0.8]
\node (a) at (0,0) [] {$\mathcal{Y} = \mathcal{X} \times_{\Hilb{n}{p(t)}} S$};
\node (b) at (4,0) [] {$\mathcal{X}$};
\node (c) at (0,-2) [] {$S$};
\node (d) at (4,-2) [] {$\Hilb{n}{p(t)}$};

\draw [->] (a) -- (c);
\draw [->] (b) -- (d);
\draw [->,dashed] (a) -- (b);
\draw [->,dashed] (c) --node[rectangle,fill=white,inner sep=1pt]{\footnotesize $\exists!$} (d);
\end{tikzpicture}
\end{center}

The Hilbert scheme is a projective scheme and it is usually defined as subscheme of a suitable Grassmannian. Following the notation used by Gotzmann in \cite{Gotzmann}, given a subscheme $X \subset \PP^n = \Proj \K[x_0,\ldots,x_n]$ and the corresponding saturated ideal $I_X \subset \K[x_0,\ldots,x_n]$ (\gls{KX} for short), we will say Hilbert polynomial of $I_X$ referring to the Hilbert polynomial $p(t)$ of $X$, whereas we will say volume polynomial of $I_X$ referring to the polynomial $q(t) = \binom{n+t}{n} - p(t)$ such that $\dim_{\K} I_t = q(t),\ t \gg 0$. By Gotzmann's Regularity Theorem it is well known that for an integer $r$ large enough, for each subscheme $X \subset \PP^n$ parametrized by $\Hilb{n}{p(t)}$, the saturated ideal $I_X$ is generated in degree lower than or equal to $r$ and that $I_r$ is a $q(r)$-dimensional vector subspace of the base vector space of homogeneous polynomials of $\K[x]$ of degree $r$. Hence any point of $\Hilb{n}{p(t)}$ can be naturally identified with a point of the Grassmannian $\Grass{q(r)}{\K[x]_r}{}$ and then embedded by the Pl\"ucker embedding in the projective space $\PP^{E},\ E = \binom{\binom{n+r}{r}}{q(r)}$. The dimension $E$ of the projective space in which we can embed $\Hilb{n}{p(t)}$ clearly becomes very huge just considering non-trivial cases. This fact reveals immediately the great difficulty of studying explicitly and globally the Hilbert scheme, indeed even Hilbert schemes of easy geometric objects give rise to intractable problems of computational algebra.

Dealing with the problem of finding an ideal defining $\Hilb{n}{p(t)}$ as projective scheme, the study can be oriented towards the equations generating such an ideal and particularly to their degree. Iarrobino and Kleiman (1999) \cite[Appendix C]{IarrobinoKanev} proved that there exists an ideal defining the Hilbert scheme $\Hilb{n}{p(t)}$ as subscheme of $\Grass{q(r)}{\K[x]_r}{}$ generated by polynomials of degree $q(r+1)+1$ in the Pl\"ucker coordinates, and afterwards Haiman and Sturmfels (2004) \cite{HaimanSturmfels}, proving a conjecture by Bayer (1982) \cite{BayerThesis}, defined an ideal generated by polynomials of degree $n+1$.

In\hfill Chapter\hfill \ref{ch:HilbertScheme},\hfill after\hfill having\hfill recalled\hfill some\hfill background\hfill material\hfill about\\ representable functors, Grassmannians and Hilbert schemes, we introduce a set of generators for any exterior power $\wedge^l W$ of a subspace $W \in \Grass{q}{N}{}$ depending \emph{linearly} on the Pl\"ucker coordinates of $W$, given by Pl\"ucker embedding. Exploiting this result in the case of Hilbert schemes, we give new and simpler proofs of the theorems about the degree of the equation defining $\Hilb{n}{p(t)} \subset \Grass{q(r)}{\K[x]_r}{}$ by Iarrobino-Kleiman and Bayer-Haiman-Sturmfels. Furthermore in Chapter \ref{ch:LowDegreeEquations}, we introduce a new ideal defining the Hilbert scheme as subscheme of a  Grassmannian generated by equations of degree smaller than or equal to $\deg p(t) + 2 < n+1 \ll q(r+1)+1$.

\bigskip

The first relevant property proved about the Hilbert scheme is surely its connectedness, proved by Hartshorne (1966) in its PhD thesis \cite{HartshorneThesis}. He used a basic idea widely exploited in the field of commutative algebra in the following years, that is to reduce the study on monomial ideals obtained by flat deformation of any ideal defining a point   on a Hilbert scheme. This idea has to be carefully managed in this context. Using the modern language of Gr\"obner degeneration theory, the problem is that any change of coordinates $g \in \GL(n+1)$ defines an isomorphism of $\Hilb{n}{p(t)}$ that identifies the point $\Proj \K[x]/I \in \Hilb{n}{p(t)}$ defined by an ideal $I$ and the point $\Proj \K[x]/(g\centerdot I)$ defined by $g\centerdot I$, whereas the point defined by $\IN(I)$ could not be mapped to the point defined by $\IN(g\centerdot I)$. To overcome this possible ambiguity, for any ideal $I$, given the equivalence relation $g \sim g' \Leftrightarrow \IN(g\centerdot I) = \IN(g'\centerdot I)$, it was proved that one of the equivalence classes corresponds to an open subset $U$ of $\GL(n+1)$. Hence associating to any ideal $I$ the so-called generic initial ideal $\IN(g\centerdot I)$ computed considering a change of coordinate $g$ in the open equivalence class turns out to be consistent with the properties of the Hilbert scheme. Generic initial ideals belong to the class of monomial ideals called Borel-fixed, because fixed by the action of the Borel subgroup of $\GL(n+1)$ composed by the upper triangular matrices. They are fundamental in the study of Hilbert schemes for two reasons:
\begin{enumerate}
\item each component and each intersection of components of $\Hilb{n}{p(t)}$ contains at least one point defined by a Borel-fixed ideal (roughly speaking they are distributed all over the Hilbert scheme);
\item they have a strong combinatorial characterization that makes them very interesting, also from an algorithmic perspective.
\end{enumerate}
In Chapter \ref{ch:BorelFixedIdeals}, after having recalled the main properties of Borel-fixed ideals and showed several ways to represent their combinatorial structure, we expose an algorithm for computing all the (saturated) Borel-fixed ideals in $\K[x]$ with Hilbert polynomial $p(t)$, that is for computing all the points of $\Hilb{n}{p(t)}$ defined by Borel-fixed ideals. Then we discuss how the number of Borel-fixed ideals varies increasing the number of variables of the polynomial ring and changing the Hilbert polynomial. Furthermore we propose new definitions of ideals that generalize the notion of lexicographic ideal, that we call segment ideals.

\bigskip

The basic idea of Hartshorne's proof of the connectedness of the Hilbert scheme is to construct a sequence of deformations and specializations (through distractions) of Borel-fixed ideals (he called it \emph{balanced} ideals) in order to reach the point defined on $\Hilb{n}{p(t)}$ by the unique (saturated) lexicographic ideal associated to $p(t)$ (see \cite{Macaulay}). A second proof was given by Peeva and Stillman (2005) \cite{PeevaStillman}, who basically rewrote Hartshorne's idea in terms of Gr\"obner degenerations. In both cases affine flat deformations are used. In Chapter \ref{ch:deformations}, we introduce a new type of flat deformations involving Borel-fixed ideals simply relying on their combinatorial structure that lead to rational curves on the Hilbert scheme. The idea of giving a \lq\lq direction\rq\rq\ to the deformations works also in this case, so that a new proof of the connectedness of the Hilbert scheme is proposed. Morever we show that all the points defined by segment ideals can play the same role of the lexicographic point in Hartshorne and Peeva-Stillman's proofs. Finally we are able to define families over $(\PP^1)^{\times s}$ that can be used to detect set of points defined by Borel-fixed ideals lying on a same component of the Hilbert scheme. 

\medskip

As explained above, the explicit and global study of the Hilbert scheme is unachievable because of the huge number of parameters needed to describe these kind of geometric families. An alternative approach is the local one, i.e. considering the open covering induced on $\Hilb{n}{p(t)}$ by the standard affine open covering of $\PP^E$ through the Pl\"ucker embedding. This has been the true starting point of the thesis, indeed the first topic I dealt with in my research activity is the construction of families of ideals sharing the same initial ideal. In \cite{NotariSpreafico} Notari and Spreafico (2000) proposed to cover set-theoretically the Hilbert scheme with such families. In the first part of Chapter \ref{ch:openSubsets}, we prove that the family of ideals $I$ such that $\IN(I) = J$, that we call Gr\"obner stratum of $J$ and denote by $\St{}(J)$, has a well-defined structure of affine scheme and we determine the conditions in order for a Gr\"obner stratum to be an open subset of the corresponding Hilbert scheme, namely a local description of its scheme structure. The ideal $J$ is having to be a segment ideal, truncation of a saturated Borel-fixed ideal in degree equal to the degree $r$ used to define the Grassmannian $\Grass{q(r)}{\K[x]_r}{}$.
$\St{}(J)$ can be viewed quite naturally as a homogeneous variety with respect to (w.r.t. for short) a non-standard positive grading. This property has a great relevance in a computational perspective, because it allows to reduce significantly the number of parameters describing $\St{}(J)$ (and also the associated open subset of $\Hilb{n}{p(t)}$), indeed we prove that in many cases the degree of the truncation can be lowered while obtaining the same family (an isomorphic one) of ideals.

Unfortunately this technique does not solve the problem of the local study of the Hilbert scheme, because not every Borel-fixed ideal is a segment ideal, so that this method can not be always applied and above all the number of open subsets that in principle we would need to consider, i.e. the number of Gr\"obner strata we would have to compute, is still enormous. Giving up the segment hypothesis is really costly, because we can no longer use tools provided by Gr\"obner theory, mainly the \emph{noetherian} Buchberger's algorithm, that are the basis of construction of Gr\"obner strata. Thus the central part of Chapter \ref{ch:openSubsets} is devoted to the ideation and development of a new noetherian algorithm of polynomial reduction based solely on the combinatorial properties of Borel-fixed ideals avoiding any term ordering. With this new procedure we define more general families of ideals, that include Gr\"obner strata: a family of this type is constructed from a Borel-fixed ideal $J$, so we call it $J$-marked family and we denote it by $\Mf(J)$. The property common to each ideal $I \in \Mf(J)$ is that the set of monomials not belonging to $J$ represents a basis of $\K[x]/I$ as $\K$-vector space. 

In the final part of Chapter \ref{ch:openSubsets}, we moreover answer also to the second problem about the local study of the Hilbert scheme, that is the huge number of open affine subsets to be considered. We prove that it is sufficient to study the $J$-marked families $\Mf(J)$, where $J$ is the truncation in some degree $\leqslant r$ of a saturated Borel-fixed ideal defining a point of $\Hilb{n}{p(t)} \subset \Grass{q(r)}{\K[x]_r}{}$, and then to exploit the action of the linear group $\GL(n+1)$ on the Hilbert scheme.

\bigskip

In Chapter \ref{ch:lcmHilb}, we deal with Hilbert schemes of locally Cohen-Macaulay curves in the projective space $\PP^3$. A locally Cohen-Macaulay curve is a curve without embedded or isolated points and the set of points of the Hilbert scheme $\Hilb{3}{dt+1-g}$ corresponding to locally Cohen-Macaulay curves of degree $d$ and genus $g$ turns out to be an open subset denoted by $\lcmHilb{d}{g}$. In turn $\lcmHilb{d}{g}$ contains an open subset corresponding to the set of smooth curves, that we denoted by $\smHilb{d}{g}$, i.e.
\[
\smHilb{d}{g} \subset \lcmHilb{d}{g} \subset \Hilb{3}{dt+1-g}.
\]
As seen in Chapter \ref{ch:deformations}, the full Hilbert scheme is connected, whereas there are known examples of Hilbert schemes of smooth curves in $\PP^3$ which are not connected (for instance $\smHilb{9}{10}$ \cite[Chapter IV Example 6.4.3]{Hartshorne}). For the Hilbert scheme of locally Cohen-Macaulay curves nothing is known, in the sense that there are no examples of non-connected Hilbert schemes and neither there is a proof of the connectedness in the general case.  An approach similar to that one used in the proof of the connectedness of $\Hilb{n}{p(t)}$, for instance a sequence of Gr\"obner deformations in a fixed direction, has always seemed unsuitable, because for any term ordering the initial ideal is a monomial ideal, hence except for some rare cases (ACM curves), a Gr\"obner degeneration would lead to a curve with embedded points.

A Hilbert scheme considered a \lq\lq good\rq\rq\ candidate of being non-connected was $\lcmHilb{4}{-3}$. In Chapter \ref{ch:lcmHilb} we show that this Hilbert scheme is indeed connected, by constructing a Gr\"obner deformation with generic fiber corresponding to four disjoint line on a smooth quadric and special fiber an extremal curve. The key point is to choose appropriately a weight order $\omega$ (which is not a total order on the monomials) such that the initial ideal w.r.t. $\omega$ is not necessarily a monomial ideal and such that the degeneration \lq\lq approaches\rq\rq\ the component of $\lcmHilb{d}{g}$ of extremal curves.

\bigskip

In order to strengthen the algorithmic purpose, I wrote lots of lines of code organized in some libraries to explicitly calculate the objects introduced theoretically in the thesis in many non-trivial example. Appendix \ref{ch:Hilb2P2} contains an handbook for the package \texttt{HilbertSchemesEquations} written in the \textit{Macaulay2} language, that covers the topics explained in Chapter \ref{ch:HilbertScheme}. As seen the number of parameters is in any case too large, but the equations for the Hilbert scheme $\Hilb{2}{2} \subset \Grass{4}{6}{}$ can be computed.

Appendix \ref{ch:HSCpackage} presents a library written in \texttt{java} that provides an implementation of the combinatorial structure of Borel-fixed ideals. All the algorithms described in Chapter \ref{ch:BorelFixedIdeals} and Chapter \ref{ch:deformations} are implemented and made available by means of several java applets working on any web browser with java plugins installed.

Finally\hfill the\hfill algorithms\hfill for\hfill working\hfill on\hfill marked\hfill families\hfill and\hfill Gr\"obner\hfill strata\\ introduced\hfill in\hfill Chapter\hfill \ref{ch:openSubsets}\hfill are\hfill made\hfill available\hfill through\hfill the\hfill \textit{Macaulay2}\hfill package\\ \texttt{MarkedSchemes} described in Appendix \ref{ch:markedSchemes}.

All this material is available at my web page
\begin{center}
\href{http://www.personalweb.unito.it/paolo.lella/HSC/index.html}{www.personalweb.unito.it/paolo.lella/HSC/}.
\end{center}
\end{introduction}

\chapter{The Hilbert scheme}\label{ch:HilbertScheme}

In this first chapter we recall the definitions leading to the introduction of the Hilbert scheme and its usual construction as subscheme of a suitable Grassmannian. 
 
\section{The representability of a functor}

\begin{definition}\label{def:functorOfPOints} For any schemes $X$ over $\K$, we define the contravariant \emph{functor of points}\index{functor!of points} from the category \gls{catSchmK} of schemes over an algebraically closed field \glshyperlink{K} of characteristic 0 to the category \gls{catSets} of sets:
\begin{equation*}
\gls{functorPoints}:\, \schm^{\circ}_{\K} \rightarrow \Sets
\end{equation*}
such that for any object $Z \in \glslink{obCat}{\Ob_{\schm_{\K}}}$
\[
h_X(Z) = \Hom(Z,X)
\]
and for any morphism $f: Z \rightarrow W \in \glslink{morCat}{\Mor_{\schm_{\K}}}$, by the diagram
\begin{center}
\begin{tikzpicture}[]
\node (W) at (2.5,0) [] {$W$};
\node (Z) at (0,0) [] {$Z$};
\node (X) at (2.5,-1.5) [] {$X$};
\draw [->](Z) --node[rectangle,fill=white]{\footnotesize $f$} (W); 
\draw [->](W) --node[rectangle,fill=white]{\footnotesize $a$} (X);
\draw [->,dashed](Z) -- (X);
\node at (1.8,-.5) [] {$\circlearrowright$};
\end{tikzpicture}
\end{center}
we define
\[
\begin{split}
h_X(f): \Hom(W,X) &\rightarrow \Hom(Z,X)\\
        \parbox{1cm}{\centering $a$} & \mapsto h_X(f)(a) = a \circ f.
\end{split}
\]
\end{definition}

By the definition $h_X$ turns out to be a contravariant functor, indeed given $Z \stackrel{f}{\rightarrow} W \stackrel{g}{\rightarrow} T$ and $a \in \Hom(T,W)$
\begin{center}
\begin{tikzpicture}[]
\node (T) at (2.5,0) [] {$T$};
\node (W) at (0,0) [] {$W$};
\node (Z) at (-2.5,0) [] {$Z$};
\node (X) at (0,-1.5) [] {$X$};
\draw [->](Z) --node[rectangle,fill=white]{\footnotesize $f$} (W); 
\draw [->](W) --node[rectangle,fill=white]{\footnotesize $g$} (T); 
\draw [->](T) --node[rectangle,fill=white]{\footnotesize $a$} (X);
\draw [->](W) -- (X);
\draw [->](Z) -- (X);
\node at (0.7,-.5) [] {$\circlearrowright$};
\node at (-0.7,-.5) [] {$\circlearrowright$};
\end{tikzpicture}
\end{center}
$h_X(g \circ f)(a) = a \circ g \circ f = h_X(g)(a) \circ f = h_X(f)\big(h_X(g)(a)\big) = \big(h_X(f)\circ h_X(g)\big)(a)$.

Given another scheme $Y$ over $\K$ and a morphism $\psi: X \rightarrow Y$, it can be defined a natural transformation $h_\psi: h_X \rightarrow h_Y$. For any $Z \in \Ob_{\schm_{\K}}$, we define $h_\psi(Z): h_X(Z) \rightarrow h_Y(Z)$ through the diagram
\begin{center}
\begin{tikzpicture}[]
\node (Z) at (0,-1.5) [] {$Z$};
\node (X) at (0,0) [] {$X$};
\node (Y) at (2.5,0) [] {$Y$};
\draw [->] (Z) --node[rectangle,fill=white]{\footnotesize $a$} (X);
\draw [->,dashed] (Z) -- (Y);
\draw [->] (X) --node[rectangle,fill=white]{\footnotesize $\psi$} (Y);
\node at (0.7,-0.5) [] {$\circlearrowright$};
\end{tikzpicture}
\end{center}
and for any $f: Z \rightarrow W \in \Mor_{\schm_{\K}}$, 
\[
\begin{split}
h_\psi(f) : \parbox{1.8cm}{\centering $h_X(f)$} & \rightarrow \parbox{1.8cm}{\centering $h_Y(f)$}\\
            \big(a \mapsto a \circ f \big) & \mapsto \big(\psi\circ a \mapsto \psi\circ a \circ f\big)
\end{split}
\]
as showed in the following diagram
\begin{center}
\begin{tikzpicture}[]
\node (W) at (2.5,0) [] {$W$};
\node (Z) at (0,0) [] {$Z$};
\node (X) at (2.5,-1.5) [] {$X$};
\node (Y) at (5,-1.5) [] {$Y$};
\draw [->](Z) --node[rectangle,fill=white]{\footnotesize $f$} (W); 
\draw [->](W) --node[rectangle,fill=white]{\footnotesize $a$} (X);
\draw [->](X) --node[rectangle,fill=white]{\footnotesize $\psi$} (Y);
\draw [->](Z) --node[rectangle,fill=white]{\footnotesize $h_X(f)(a)$} (X);
\end{tikzpicture}
\end{center}

Therefore we defined a functor $h$ from the category $\schm_{\K}$ of schemes over $\K$ to the category of functors from $\schm_{\K}$ to $\Sets$
\begin{equation}\label{eq:hFullyFaithful}
\begin{split}
h:\ \schm_{\K} &\rightarrow \glslink{catFun}{\Fun{\schm_{\K}}{\Sets}}\\
 \parbox{2cm}{\centering $X$} &\mapsto \parbox{3cm}{\centering $h_X$}\\
 \parbox{2cm}{\centering\small $\big(\psi:X\rightarrow Y\big)$} &\mapsto \parbox{3cm}{\centering \small$\big(h_{\psi}: h_X \rightarrow h_Y\big)$}\\
\end{split}
\end{equation}

\begin{proposition}
The functor $h$ described in \eqref{eq:hFullyFaithful} is fully faithful\index{functor!fully faithful}.
\end{proposition}
\begin{proof}
We have to prove that the function $\psi \mapsto h_{\psi}$ is a bijection, so starting from a map $f: h_X \rightarrow h_Y$ we will define a morphism between $X$ and $Y$ and we will prove that this correspondence is the inverse function.

$h_X(X)$ obviously contains the identity map $\id_X$ of $X$, so by means of $f(X): h_X(X) \rightarrow h_Y(X)$, we can define the map
\[
\varphi= f(X)(\id_X): X \rightarrow Y,
\]
and we want to show that $h_{\varphi} = f$.

For any scheme $Z \in \Ob_{\schm_{\K}}$ and any morphism $g: Z\rightarrow X$, we have the diagram
\begin{center}
\begin{tikzpicture}[]
\node (A) at (0,0) [] {$h_X(Z)$};
\node (B) at (3,0) [] {$h_Y(Z)$};
\node (C) at (3,-2) [] {$h_Y(X)$};
\node (D) at (0,-2) [] {$h_X(X)$};
\draw [->](A) --node[rectangle,fill=white]{\footnotesize $f(Z)$} (B); 
\draw [->](D) --node[rectangle,fill=white]{\footnotesize $f(X)$} (C);
\draw [<-](A) --node[rectangle,fill=white]{\footnotesize $h_X(g)$} (D);
\draw [<-](B) --node[rectangle,fill=white]{\footnotesize $h_Y(g)$} (C);
\node at (1.5,-1) [] {$\circlearrowright$};
\end{tikzpicture}
\end{center}
Note that $g = h_X(g)(\id_X),\ \forall\ g \in h_{X}(Z)$, so
\[
f(Z)(g) = h_Y(g)\big(f(X)(\id_X)\big) = h_Y(g)(\varphi) = \varphi \circ g = h_{\varphi}(Z)(g). \qedhere
\]
\end{proof}

\begin{theorem}[\textbf{Yoneda's Lemma} {\cite[Lemma VI-1]{EisenbudHarris}}]\label{th:Yoneda}\index{Yoneda's Lemma}
The\hfill functor\hfill of\hfill points\hfill $h_X$\hfill in\\ $\Fun{\schm_{\K}}{\Sets}$ uniquely determines the scheme $X \in \schm_{\K}$ up to isomorphism.
\end{theorem}

Now we can introduce the notion of representable functor.

\begin{definition}\label{def:representableFunctor}\index{functor!representable}
Let $\mathcal{F} : \schm^\circ_{\K} \rightarrow \Sets$ be a contravariant functor. $\mathcal{F}$ is \emph{representable} if there exists a scheme $X \in \Ob_{\schm_{\K}}$ such that $\mathcal{F} \simeq h_X$. By Yoneda's Lemma, we know that such a scheme is uniquely determined, so we will say that the scheme $X$ represents $\mathcal{F}$.
\end{definition}

Let us conclude this section recalling some useful results about the representability of a functor.

\begin{proposition}[{\cite[Proposition VI-2]{EisenbudHarris}}]\label{prop:reductionAffine} A scheme $X$ over $\K$ is completely determined by the restriction of its functor of points to affine schemes over $\K$; in fact
\[
h': \schm_{\K} \rightarrow \Fun{\gls{catAffSchmK}}{\Sets}
\]
is\hfill an\hfill equivalence\hfill of\hfill the\hfill category\hfill of\hfill schemes\hfill over\hfill $\K$\hfill with\hfill a\hfill full\hfill subcategory\hfill of\\ $\Fun{\schm_{\K}}{\Sets}$.
\end{proposition}

From now on, in place of the contravariant functor of points (Definition \ref{def:functorOfPOints}), we will consider the  covariant functor
\begin{equation*}
h_A:\, \gls{catKalg} \rightarrow \Sets
\end{equation*}
such that for any object $B \in \Ob_{\Kalg}$
\[
h_A(B) = \Hom(\Spec B,\Spec A)
\]
and for any morphism $f: B \rightarrow C \in \Mor_{\Kalg}$, by the diagram
\begin{center}
\begin{tikzpicture}[]
\node (W) at (2.5,0) [] {$\Spec B$};
\node (Z) at (-.5,0) [] {$\Spec C$};
\node (X) at (2.5,-1.5) [] {$\Spec A$};
\draw [->](Z) --node[rectangle,fill=white]{\footnotesize $f^\ast$} (W); 
\draw [->](W) --node[rectangle,fill=white]{\footnotesize $a$} (X);
\draw [->,dashed](Z) -- (X);
\node at (1.8,-.5) [] {$\circlearrowright$};
\end{tikzpicture}
\end{center}
we define
\[
\begin{split}
h_A(f):\, h_A(B) &\rightarrow h_A(C)\\
        \parbox{1cm}{\centering $a$} & \mapsto h_A(f)(a) = a \circ f^\ast.
\end{split}
\]

\section{Grassmannians}\label{sec:grassmannians}\index{Grassmannian|(}

Let us consider a $\K$-vector space $V$ of dimension $N$ with a basis $\{\vv_1,\ldots,\vv_N\}$. The Grassmannian $\Grass{q}{N}{\K}$ parametrizes the set of all vector subspaces $W \subset V$ of dimension $q$. Every $q$-dimensional subspace $W$ can be described as a row span of a $q \times N$ matrix $\mathfrak{M}(W)$ of maximal rank. Furthermore the list of all maximal minors of such a matrix up to scale determines uniquely the space $W$ (see \cite[Proposition 14.2]{MillerSturmfels}), that is the matrix $\mathfrak{M}(W)$ representing $W$ is unique up to multiplication by invertible $q \times q$ matrices. 

\index{Pl\"ucker!coordinates|(}By this argument the Grassmannian $\Grass{q}{N}{\K}$  can be embedded in the projective space $\PP_{\K}^{\binom{N}{q}-1}$:\index{Pl\"ucker!embedding}
\begin{equation}\label{eq:pluckerEmbedSubspace}
\begin{split} 
\gls{PE}: \parbox{3cm}{\centering $\Grass{q}{N}{\K}$} &\rightarrow \parbox{2cm}{\centering $\PP_{\K}[\gls{extAlg}]$} \\
       W = \langle \ww_1,\ldots,\ww_q\rangle& \mapsto \ww_1 \wedge \cdots \wedge \ww_q
\end{split}
\end{equation}
Let $\ww_i = \sum\limits_{j=1}^N \delta_{ij}\, \vv_j,\ i = 1,\ldots,q$ be the decomposition of the basis of $W$ with respect the basis of $V$ and let
\[
\left\{ \vv_{\II}^{(q)} = \vv_{\ii_1} \wedge \cdots \wedge \vv_{\ii_q}\quad \big\vert\quad \II=(\ii_1,\ldots,\ii_q),\ 1 \leqslant \ii_1 < \cdots < \ii_q \leqslant N \right\}
\]
the standard basis of $\wedge^q V$. The product $\ww_1 \wedge \cdots \wedge \ww_q$ decomposes as
\[
 \ww_1 \wedge \cdots \wedge \ww_q = \sum_{\vert\II\vert = q} \Delta_{\II}(W)\, \vv_{\II}^{(q)}
\]
where the coefficient $\Delta_{\II}(W)$ is just the determinant of the submatrix $\mathfrak{M}_{\II}(W)$ of $\mathfrak{M}(W)$ composed by the columns of the vectors $\vv_{\ii_1},\ldots,\vv_{\ii_N}$ and the set $[\ldots:\Delta_{\II}(W):\ldots]$ is the set of Pl\"ucker coordinates of $W$, hence we will consider $\PP_{\K}[\wedge^q V] = \Proj \K[\ldots,\Delta_I,\ldots] = \Proj \gls{KDelta}$.

Another way to determine the embedding is considering the short exact sequence
\[
0 \ \longrightarrow\ W\ {\longrightarrow}\ V\ \stackrel{\pi_W}{\longrightarrow}\ V/W \ \longrightarrow\ 0 
\]
and the induced epimorphism of the $p$-th exterior power where $p = N-q$, that is
\begin{equation}\label{eq:epimorphismExtPower}
\wedge^p V\ \stackrel{\pi^{(p)}_W}{\longrightarrow}\ \wedge^p V/W\ \longrightarrow\ 0.
\end{equation}
Since $\dim_{\K} V/W = \dim_{\K} V - \dim_{\K} W = p$, the exterior algebra $\wedge^p W$ is isomorphic to the field $\K$ and the map $\pi_W^{(p)}$ is uniquely determined by $W$ up to multiplication by scalar. So we can identify $W$ with the set  of the images by $\pi_W^{(p)}$ of the vectors of the basis of $\wedge^p V$
\[
\left\{ \vv_{\JJ}^{(p)} = \vv_{\jj_1} \wedge \cdots \wedge \vv_{\jj_p}\quad \big\vert\quad \JJ=(\jj_1,\ldots,\jj_p),\ 1 \leqslant \jj_1 < \cdots < \jj_p \leqslant N \right\}
\]
precisely let us define
\begin{equation}\label{eq:defTheta}
\Theta_{\JJ}(W) = \pi_W^{(p)}(\vv_\JJ)
\end{equation}
and the embedding\index{Pl\"ucker!embedding}
\begin{equation}\label{eq:pluckerEmbedQuotient}
\begin{split} 
\psi: \Grass{q}{N}{\K} &\rightarrow \parbox{2.5cm}{\centering $\PP^{\binom{N}{p}-1}_{\K}$} \\
       \parbox{1.5cm}{\centering $W$} & \mapsto [\ldots:\Theta_{\JJ}(W):\ldots]
\end{split}
\end{equation}
considering $\PP^{\binom{N}{p}-1}_{\K} = \Proj \K[\ldots,\Theta_{\JJ},\ldots] = \Proj \K[\Theta]$ where the variables $\Theta_{\JJ}$ are another set of Pl\"ucker coordinates.

\begin{proposition}
The sets of Pl\"ucker coordinates $[\ldots:\Delta_{\II}:\ldots]$ and $[\ldots:\Theta_{\JJ}:\ldots]$ are equivalent. More precisely, for every $W \in \Grass{q}{N}{\K}$, there exists $c \in \K$ such that
\begin{equation}\label{eq:pluckerCoordEquivalence}
\Delta_{\II}(W) = c\, \varepsilon_{\JJ\vert\II}\, \Theta_{\JJ}(W),\qquad \JJ\cup\II = \{1,\ldots,N\}
\end{equation}
where $\varepsilon_{\JJ\vert\II}$ is the signature of the permutation $\sigma$ that orders $\JJ\vert\II = (\jj_1,\ldots,\jj_p,\ii_1,\ldots,\ii_q)$, i.e. $\sigma(\JJ\vert\II) = (1,\ldots,N)$.
\end{proposition}
\begin{proof}
Let us suppose $W$ contained in the vector space $V$ with basis $\{\vv_1,\ldots,\vv_N\}$ and let $p = N-q$. Since $\wedge^p (V/W) \simeq \K$, the morphism $\pi_W^{(p)}$ described in \eqref{eq:epimorphismExtPower} turns out to be a linear functional over $\wedge^p V$. Let $\{\dv^1,\ldots,\dv^N\}$ be the usual basis of the dual space $V^\ast$, that is
\[
\dv^j(\vv_i) = \left\{\begin{array}{l l}
1, & i=j\\
0, & i \neq j
\end{array}\right..
\]
Directly by the definition of $\Theta_{\JJ}(W)$ given in \eqref{eq:defTheta}, we can identify $\pi^{(p)}_W$ with
\[
\pi^{(p)}_W\ \longleftrightarrow\ \sum_{\vert\JJ\vert = p} \Theta_{\JJ}(W)\, \dv^{\,\JJ}_{(p)} \in \wedge^p (V^\ast) = (\wedge^p V)^\ast
\]
where $\dv^{\,\JJ}_{(p)} = \dv^{\,\jj_1} \wedge \cdots \wedge \dv^{\,\jj_p}$.

Through the standard isomorphism
\[
\begin{split}
\wedge^p (V^\ast) & \stackrel{\sim}{\longrightarrow}  \parbox{2.4cm}{\centering$\wedge^{q} V$}\\
\parbox{1.2cm}{\centering$\overline{t}$} & \longmapsto  \overline{t}(\vv_1 \wedge \cdots \wedge \vv_N)
\end{split}
\]
the image of $\pi^{(p)}_W$ is
\[
\sum_{\vert\JJ\vert = p} \Theta_{\JJ}(W)\, \dv^{\,\JJ}_{(p)} \big(\vv_1 \wedge \cdots \wedge \vv_N\big)
\]
and rewriting in each addend $\vv_1 \wedge \cdots \wedge \vv_N$ as $\vv_{\JJ}^{(p)} \wedge \varepsilon_{\JJ\vert\II}\, \vv_{\II}^{(q)}$ where $\JJ \cup \II = \{1,\ldots,N\}$ and $\varepsilon_{\JJ\vert\II}$ is the signature of the permutation that orders $\JJ\vert\II$ (equal obviously also to the signature of the inverse permutation $\sigma$ such that $\sigma(1,\ldots,N) = \JJ\vert\II = (\jj_1,\ldots,\jj_p,$ $\ii_1,\ldots,\ii_q)$), we finally obtain 
\[
\sum_{\vert\JJ\vert = p} \Theta_{\JJ}(W)\, \dv^{\,\JJ}_{(p)} \big(\vv_\JJ^{(p)} \wedge \varepsilon_{\JJ\vert\II}\, \vv_\II^{(q)}\big) = \sum_{\begin{subarray}{c}\vert\II\vert = q\\ \II \cup \JJ = \{1,\ldots,N\}\end{subarray}} \varepsilon_{\JJ\vert\II}\, \Theta_{\JJ}(W)\, \vv_\II^{(q)}.
\]
By duality the vector subspace generated by this elements has to coincide to the vector subspace identified by the injection
\[
\wedge^q W\ \longrightarrow\ \wedge^q V.
\]
Chosen a basis $\{\ww_1,\ldots,\ww_q\}$ of $W$, the generator of $\wedge^q W$ in $\wedge^q V$ is
\[
\ww_1 \wedge \cdots \wedge \ww_q = \sum_{\vert\II\vert = q} \Delta_I(W)\, \vv_{\II}^{(q)}. \qedhere
\]
\end{proof}
\index{Pl\"ucker!coordinates|)}

\begin{remark}
Given any (even not-ordered) set of indices $H=(\hh_1,\ldots,\hh_q)$ (resp. $H = (\hh_1,\ldots,\hh_p)$), we will denote with $\Delta_{H}(W)$ (resp. $\Theta_{H}(W)$) the determinant of the submatrix of $\mathfrak{M}(W)$ obtained considering the columns corresponding to the vectors $\vv_{\hh_1},\ldots,\vv_{h_q}$ (resp. the image of $\vv_{\hh_1}\wedge\cdots\wedge\vv_{\hh_p}$ using $\pi_W^{(p)}$). It is easy to check that $\Delta_{H} = \varepsilon_{H} \Delta_{\KK}$ (resp. $\Theta_{H} = \varepsilon_{H} \Theta_{\KK}$), where $\varepsilon_{H}$ is the signature of the permutation $\sigma$ that orders $H$, $\KK = \sigma(H)$ is the corresponding ordered set of indices and $\Delta_{\KK}$ (resp. $\Theta_{\KK}$) is a Pl\"ucker coordinate. From now on, if not specified the set of indices are considered in increasing order. 

Given two multi-indices $\KK = (\kk_1,\ldots,\kk_a)$ and $\HH = (\hh_1,\ldots,\hh_b)$, we will denote by $\KK\vert\HH$ the set of indices $(\kk_1,\ldots,\kk_a,\hh_1,\ldots,\hh_b)$, that in general will not be ordered, whereas we will denote with the union $\KK\cup\HH$ the ordered multi-index containing the indices belonging to both $\KK$ and $\HH$.
For instance given $\KK = (1,5)$ and $\HH=(2)$, $\KK\vert\HH = (1,5,2)$, $\HH\vert\KK = (2,1,5)$ and $\KK\cup\HH = \HH\cup\KK = (1,2,5)$. Coming back to Pl\"ucker coordinates the following relation holds:
\[
\Delta_{\KK\vert\HH} = \varepsilon_{\KK\vert\HH} \Delta_{\KK\cup\HH} \qquad (\text{resp. } \Theta_{\KK\vert\HH} = \varepsilon_{\KK\vert\HH} \Theta_{\KK\cup\HH}).
\]
\end{remark}

\bigskip

To determine the equations of the subscheme $\PE\big(\Grass{q}{N}{\K}\big) \subset \PP[\wedge^q V]$, that is the conditions such that an element $\underline{u}^{(q)} \in \wedge^q V$ can be decomposed as exterior product of $q$ vectors
\[
\underline{u}^{(q)} = \ww_1 \wedge \cdots \wedge \ww_q,\qquad \ww_1,\ldots,\ww_q \in V
\]
let us start considering the \emph{contraction operator} (sometimes also called convolution). For any element $\dv^{\, j}$ of the basis of $V^\ast$, let us define the operator
\begin{equation*}
i^{\ast}(\dv^{\, j}):\, \wedge^s V \rightarrow \wedge^{s-1} V
\end{equation*}
as the operator sending the generic element $\vv^{(s)}_{\II} = \vv_{\ii_1} \wedge \cdots \wedge \vv_{\ii_s}$ of the basis of $\wedge^s V$ to $0$ if $j$ does not belong to $\II$ whereas if $j = \ii_k$ for any $k$
\[
i^\ast(\dv^j)(\vv^{(s)}_{\II}) = (-1)^{k-1}\, \vv_{\ii_1} \wedge \cdots \wedge \vv_{\ii_{k-1}} \wedge \vv_{\ii_{k+1}} \wedge \cdots \wedge \vv_{\ii_{s}} = (-1)^{k-1}\, \vv^{(s-1)}_{\II\setminus\{j\}}
\]
so that for the generic element of $\wedge^s V$
\[
i^{\ast}(\vv^j)\left(\sum_{\vert \II \vert = s} a_{\II}\, \vv^{(s)}_{\II}\right) = \sum_{\begin{subarray}{c}\vert \II \vert = s\\ \exists\ k\text{ s.t. } \ii_k = j \end{subarray}} (-1)^{k-1} \, a_{\II} \, \vv^{(s-1)}_{\II\setminus\{j\}}.
\]

Extending by linearity, for any element $\overline{t} = \sum_{j=1}^N b_j\, \dv^j \in V^\ast$ we define $i^\ast(\overline{t}):\wedge^s V \rightarrow \wedge^{s-1} V$ as the operator
\begin{equation}
i^\ast(\overline{t}) = \sum_{j=1}^N b_j\, i^\ast(\dv^j)
\end{equation}

\begin{definition}\label{def:contractionOperator}\index{contraction operator}
For any $\overline{t}_{(r)} = \sum_{\vert\JJ\vert = r} b_{\JJ}\, \dv^{\,\JJ}_{(r)} \in \wedge^r (V^\ast)$ we define the \emph{contraction operator}
\begin{equation}\label{eq:contractionOperator}
i^\ast(\overline{t}_{(r)}) : \wedge^s V  \rightarrow \wedge^{s-r} V
\end{equation}
as
\[
\begin{split}
i^\ast(\overline{t}_{(r)})\left(\sum_{\vert\II\vert=s} a_{\II}\, \vv^{(s)}_{\II}\right) & {}= \left( \sum_{\vert\JJ\vert = r} b_{\JJ}\, i^\ast(\dv^{\,\JJ}_{(r)})\right) \left(\sum_{\vert\II\vert=s} a_{\II}\, \vv^{(s)}_{\II}\right) =\\
&{} = \sum_{\vert\JJ\vert=r}\sum_{\vert\II\vert=s} b_{\JJ} a_{\II} \, i^\ast(\dv^{\,\JJ}_{(r)})\big(\vv^{(s)}_{\II}\big)
\end{split} 
\]
where we define $i^\ast(\dv^{\,\JJ}_{(r)})(\vv^{(s)}_{\II}) \in \wedge^{s-r} V$ through the composition of maps
\[
 \wedge^{r} V\ \xrightarrow{i^\ast(\dv^{\, \jj_r})}\ \wedge^{r-1} V \ \xrightarrow{i^\ast(\dv^{\, \jj_{r-1}})}\quad \cdots \quad \xrightarrow{i^\ast(\dv^{\, \jj_1})}\ \wedge^{s-r} V
\]
as
\begin{equation}
 i^\ast(\dv^{\, \jj_1}) \circ \cdots \circ  i^\ast(\dv^{\, \jj_r})\big( \vv^{(s)}_{\II}\big).
\end{equation}
\end{definition}

\begin{proposition}\index{decomposable vector}
A vector $\underline{u}^{(q)} \in \wedge^q V$ is decomposable, i.e. of the form  $\ww_{\ii_1} \wedge \cdots \wedge \ww_{\ii_q}$, if and only if
\begin{equation}\label{eq:pluckerCondition}
i^\ast(\overline{t}_{(q-1)})\big(\underline{u}^{(q)}\big) \wedge \underline{u}^{(q)} = 0, \qquad \forall\ \overline{t}_{(q-1)} \in \wedge^{q-1} (V^\ast). 
\end{equation}
\end{proposition}
\begin{proof}
See \cite[Chapter I Section 4.1]{Shafarevich1} and \cite[Chapter I Section 5]{GriffithsHarris}.
\end{proof}

To compute the Pl\"ucker relations, we consider the generic element $\sum_{\vert\II\vert=q} \Delta_{\II} \, \vv^{(q)}_{\II}$ of $\wedge^{q} V$ and we impose the condition \eqref{eq:pluckerCondition} considering
\[
i^\ast(\dv^{\, \JJ}_{(q-1)})\left(\sum_{\vert\II\vert=q} \Delta_{\II} \, \vv^{(q)}_{\II}\right) \wedge \sum_{\vert\II\vert=q} \Delta_{\II} \, \vv^{(q)}_{\II}
\]
for all $\dv^{\JJ}_{(q-1)}$ in the standard basis of $\wedge^{q-1} (V^\ast)$. Firstly
\[
i^\ast(\dv^{\, \JJ}_{(q-1)})\left(\sum_{\vert\II\vert=q} \Delta_{\II} \, \vv^{(q)}_{\II}\right) = \sum_{\begin{subarray}{c} \vert\II\vert = q\\ \JJ \subset \II \end{subarray}} \Delta_{\II} \, i^\ast(\dv^{\, \JJ}_{(q-1)})\big(\vv^{(q)}_{\II}\big)
\] 
and supposing that $\II \setminus \JJ = \{\ii_k\}$
\[
i^\ast(\dv^{\, \JJ}_{(q-1)})\big(\vv^{(q)}_{\II}\big) = (-1)^{q-k}\, i^\ast(\dv^{\, \JJ}_{(q-1)})\big(\vv^{(q-1)}_{\JJ}\wedge \vv_{i_k}\big) = (-1)^{\binom{q-1}{2}+q-k}\, \vv_{\ii_k}
\]
so we can rewrite
\[
\begin{split}
i^\ast(\dv^{\, \JJ}_{(q-1)})\left(\sum_{\vert\II\vert=q} \Delta_{\II} \, \vv^{(q)}_{\II}\right) &{}= (-1)^{\binom{q-1}{2}} \sum_{\begin{subarray}{c} \vert\II\vert = q\\ \JJ \subset \II \\ \II\setminus\JJ = \{\ii_k\} \end{subarray}} (-1)^{q-k}\,\Delta_{\II}\, \vv_{\ii_k} =\\
&{} = (-1)^{\binom{q-1}{2}} \sum_{\begin{subarray}{c} 1 \leqslant \jj \leqslant N\\ \jj \notin \JJ \end{subarray}} \varepsilon_{\JJ\vert(\jj)}\, \Delta_{\JJ\cup(\jj)}\, \vv_{\jj}
\end{split}
\]
so that finally
\begin{eqnarray}
\notag &&\left( (-1)^{\binom{q-1}{2}} \sum_{\begin{subarray}{c} 1 \leqslant \jj \leqslant N\\ \jj \notin \JJ \end{subarray}} \varepsilon_{\JJ\vert(\jj)}\, \Delta_{\JJ\cup(\jj)}\, \vv_{\jj}\right) \bigwedge \left(\sum_{\vert\II\vert=q} \Delta_{\II} \, \vv^{(q)}_{\II}\right) = \\
\notag &&\qquad {}= (-1)^{\binom{q-1}{2}} \sum_{\begin{subarray}{c} 1 \leqslant \jj \leqslant N\\ \jj \notin \JJ \end{subarray}} \, \sum_{\vert\II\vert=q} \varepsilon_{\JJ\vert(\jj)}\, \Delta_{\JJ\cup(\jj)} \Delta_{\II} \, \vv_{\jj} \wedge \vv^{(q)}_{\II} =\\
\label{eq:finalSumPluckerRel} &&\qquad\qquad{}= (-1)^{\binom{q-1}{2}} \sum_{\begin{subarray}{c} 1 \leqslant \jj \leqslant N\\ \jj \notin \JJ \end{subarray}} \, \sum_{\vert\II\vert=q} \varepsilon_{\JJ\vert(\jj)}\varepsilon_{(\jj)\vert\II}\, \Delta_{\JJ\cup(\jj)} \Delta_{\II} \, \vv^{(q+1)}_{(\jj)\cup\II}.
\end{eqnarray}
By dimension arguments, each vector $\vv^{(q+1)}_{\II}$ of the basis $\wedge^{q+1} V$ comparing in the previous sum contains at least two vectors $\vv_{\ii}$ such that $\ii \notin \JJ$. Let us decompose such a $\II$ as
\[
\II = (\II \cap \JJ) \cup \KK, \qquad \vert\KK\vert \geqslant 2.
\]
For every $\ii \in \KK$,  $\vv^{(q+1)}_{\II}$ compare in the wedge product $\varepsilon_{\JJ\vert(\ii)} \Delta_{\JJ\cup(\ii)} \vv_{\ii} \wedge \Delta_{\II\setminus(\ii)} \vv^{(q)}_{\II\setminus(\ii)}$ therefore the sum \eqref{eq:finalSumPluckerRel} contains the addend
\begin{equation}\label{eq:genericAddendPluckerRel}
(-1)^{\binom{q-1}{2}} \left(\sum_{\begin{subarray}{c}\ii \in \II\\ \ii \notin \JJ\end{subarray}} \varepsilon_{\JJ\vert(\ii)}\varepsilon_{(\ii)\vert(\II\setminus(\ii))}\, \Delta_{\JJ\cup(\ii)} \Delta_{\II\setminus(\ii)} \right) \vv^{(q+1)}_{\II}
\end{equation}

We remark that if $\vert\KK\vert = 2$, i.e. $\II = \JJ \cup (\ii_1,\ii_2)$ the coefficient of the term $\vv^{(q+1)}_{\II}$ in \eqref{eq:genericAddendPluckerRel} vanish, indeed
\[
\varepsilon_{\JJ\vert(\ii_1)}\varepsilon_{(\ii_1)\vert(\JJ\cup(\ii_2))}\, \Delta_{\JJ\cup(\ii_1)} \Delta_{\JJ\cup(\ii_2)} + \varepsilon_{\JJ\vert(\ii_2)}\varepsilon_{(\ii_2)\vert(\JJ\cup(\ii_1))}\, \Delta_{\JJ\cup(\ii_2)} \Delta_{\JJ\cup(\ii_1)} = 0
\]
because, supposing $\ii_1 < \ii_2$,
\[
\begin{split}
& \varepsilon_{\JJ\vert(\ii_1)}\varepsilon_{(\ii_1)\vert(\JJ\cup(\ii_2))} = (-1)^{p_1} \cdot (-1)^{q-p_1-1} = (-1)^{q-1}, \\
& \varepsilon_{\JJ\vert(\ii_2)}\varepsilon_{(\ii_2)\vert(\JJ\cup(\ii_1))} = (-1)^{p_2} \cdot (-1)^{q-p_2} = (-1)^{q}.
\end{split}
\]

Finally the ideal defining the Grassmannian $\Grass{q}{N}{\K}$ as subscheme of $\PP^{\binom{N}{q}-1}_{\K}$ is generated by the following set of quadrics\index{Pl\"ucker!equations}
\begin{equation}\label{eq:PluckerEquations}
\left\{ \sum_{\begin{subarray}{c}\ii \in \II\\ \ii \notin \JJ\end{subarray}} \varepsilon_{\JJ\vert(\ii)}\varepsilon_{(\ii)\vert(\II\setminus(\ii))}\, \Delta_{\JJ\cup(\ii)} \Delta_{\II\setminus(\ii)} \quad \Bigg\vert\quad \begin{array}{l} \forall\ \JJ,\ \vert\JJ\vert = q-1 \\ \forall\ \II,\ \vert\II\vert = q+1,\ \vert \II\setminus(\II\cap\JJ) \vert > 2 \end{array} \right\}.
\end{equation}

\begin{remark}
Increasing the dimension $N$ of the base vector space and/or the dimension $q$ of the subspaces, the number of Pl\"ucker coordinates and above all the number of quadrics become quickly huge, indeed for describing the Pl\"ucker embedding of $\Grass{q}{N}{\K}$, we need $\binom{N}{q}$ Pl\"ucker coordinates and the Pl\"ucker relations among them are $\binom{N}{q-1}\cdot\left[\binom{N}{q+1}-\binom{N-(q-1)}{2}\right]$. If fact these relations are redundant, so the number of quadrics sufficient to describe the ideal is much lower.
\end{remark}

\begin{example}\label{ex:MainExampleGrass}
Let us consider the Grassmannian $\Grass{4}{6}{\K}$ and given a basis $\{\vv_1,$ $\ldots,\vv_6\}$ of $V \simeq \K^6$ let
\[
\PE : \Grass{4}{6}{\K} \rightarrow \PP^{14}_{\K} = \Proj \K[\Delta_{\II}]
\]
be the Pl\"ucker embedding described in \eqref{eq:pluckerEmbedSubspace} where the 15 Pl\"ucker coordinates are indexed by ordered subsets $\II =(\ii_1,\ii_2,\ii_3,\ii_4)$ of $\{1,2,3,4,5,6\}$.
Let $\{\dv^1,\ldots,\dv^6\}$ be the standard basis of the dual space. 

$\wedge^3 (V^\ast)$ has a basis containing $20$ elements. Let us look for instance at the contraction operator defined by $\dv^2 \wedge \dv^3 \wedge \dv^6$ applied on the vector $\sum_{\vert\II\vert=4} \Delta_{\II}\, \vv^{(4)}_{\II}$:
\[
\begin{split}
i^\ast(\dv^2 \wedge \dv^3 \wedge \dv^6)\left(\sum_{\vert\II\vert = 4} \Delta_{I}\, \vv^{(4)}_{\II}\right)&{} = (-1)^{\binom{3}{2}} \sum_{\jj \in \{1,4,5\}} \varepsilon_{236|\jj}\, \Delta_{236 \cup \jj}\, \vv_{\jj} = \\
&{} = \Delta_{1236}\, \vv_{1} + \Delta_{2346}\, \vv_{4} + \Delta_{2356}\, \vv_{5}.
\end{split}
\]
Then by
\[
\left(\Delta_{1236}\, \vv_{1} + \Delta_{2346}\, \vv_{4} + \Delta_{2356}\, \vv_{5}\right) \wedge \left( \sum_{\vert\II\vert=4} \Delta_{\II}\, \vv^{(4)}_{\II}\right)
\]
we obtain the following 3 Pl\"ucker relations
\[
\begin{split}
& \varepsilon_{236\vert1}\varepsilon_{1\vert2345}\, \Delta_{1236}\Delta_{2345} + \varepsilon_{236\vert4}\varepsilon_{4\vert1235}\, \Delta_{2346}\Delta_{1235} + \varepsilon_{236\vert5}\varepsilon_{5\vert1234}\, \Delta_{2356}\Delta_{1234} = \\
&\qquad{} = -\Delta_{1236}\Delta_{2345} + \Delta_{2346}\Delta_{1235} - \Delta_{2356}\Delta_{1234},\\
& \varepsilon_{236\vert1}\varepsilon_{1\vert2456}\, \Delta_{1236}\Delta_{2456} + \varepsilon_{236\vert4}\varepsilon_{4\vert1256}\, \Delta_{2346}\Delta_{1256} + \varepsilon_{236\vert5}\varepsilon_{5\vert1246}\, \Delta_{2356}\Delta_{1246} = \\
&\qquad{} = -\Delta_{1236}\Delta_{2456} - \Delta_{2346}\Delta_{1256} + \Delta_{2356}\Delta_{1246},\\
& \varepsilon_{236\vert1}\varepsilon_{1\vert3456}\, \Delta_{1236}\Delta_{3456} + \varepsilon_{236\vert4}\varepsilon_{4\vert1356}\, \Delta_{2346}\Delta_{1356} + \varepsilon_{236\vert5}\varepsilon_{5\vert1346}\, \Delta_{2356}\Delta_{1346} = \\
&\qquad{} = -\Delta_{1236}\Delta_{3456} - \Delta_{2346}\Delta_{1356} + \Delta_{2356}\Delta_{1346}.
\end{split}
\]
 Repeating this computation for every element of the basis of $\wedge^{3} (V^\ast)$, we obtain the ideal defining the Grassmannian $\Grass{4}{6}{\K}$ as subscheme of $\PP^{14}_{\K}$ as generated by $3 \cdot 20 = 60$ equations of degree 2. In Example \ref{ex:MainExampleGrassCODE}, we will show that 45 quadrics are redundant, that is we only need 15 polynomials to define this Grassmannian.
\end{example}

\bigskip

Our next task is to understand how we can recover a set of generator of a subspace $W \in \Grass{q}{N}{\K}$ knowing its image by the Pl\"ucker embedding $\PE(W) \in \PE(\Grass{q}{N}{\K}) \subset \PP^{\binom{N}{q}-1}_{\K}$.
It is well known (see \cite[Chapter 4 Section 4.1]{Shafarevich1} and \cite[Section III.2.7]{EisenbudHarris}) that any point $[\ldots:\Delta_{\II}:\ldots] \in \PE(\Grass{q}{N}{\K})$ such that $\Delta_{(1\ldots q)} \neq 0$ corresponds to a subspace $W \subset V$ with a basis of the type
\begin{equation}\label{eq:basisW}
\left\{\ww_{\ii} = \vv_{\ii} + \sum_{\jj > q} \delta_{\ii \jj}\, \vv_{\jj}\ \Big\vert\ \ii = 1,\ldots,q\right\} \qquad \text{where}\qquad \delta_{\ii \jj} = \varepsilon_{(1\ldots q)\setminus(\ii)\vert(\ii)}\dfrac{\Delta_{(1\ldots q)\setminus(\ii)\vert(\jj)}}{\Delta_{(1\ldots q)}},
\end{equation}
in fact after multiplying any matrix $\mathfrak{M}(W)$ by the inverse matrix of $\mathfrak{M}_{(1\ldots q)}(W)$ (invertible because $\Delta_{(1\ldots q)} \neq 0$), the coefficient $\delta_{\ii \jj}$ is equal to the determinant of the matrix composed by the columns with indices in $(1,\ldots,q)\setminus(\ii)\vert(\jj)$. 
The same reasoning can be applied for every vector space $\wedge^s W,\ s \leqslant q$. Starting from the basis of $W$ described in \eqref{eq:basisW}, we can consider as basis for $\wedge^s W$ the set $\{\ww^{(s)}_{\II}\ \vert\ \forall\ \II = (\ii_1,\ldots,\ii_s),\ 1 \leqslant \ii_1 < \cdots < \ii_s \leqslant q\}$. More precisely
\begin{equation}\label{eq:basisWs}
\ww^{(s)}_{\II} = \ww_{\ii_1} \wedge \cdots \wedge \ww_{\ii_s} = \vv_{\II}^{(s)} + \sum_{\begin{subarray}{c} \vert\JJ\vert = s\\ \jj_1 > q \end{subarray}} \left\vert \begin{array}{ccc} \delta_{\ii_1 \jj_1} & \cdots & \delta_{\ii_1 \jj_s}\\ \vdots & \ddots & \vdots \\ \delta_{\ii_s \jj_1} & \cdots & \delta_{\ii_s \jj_s} \end{array}\right\vert\, \vv^{(s)}_{\JJ}.
\end{equation}
and the coefficient of $\vv^{(s)}_{\K}$ is equal to the determinant of the submatrix of $\mathfrak{M}(W)$ composed by the columns with indices in $(1 \ldots q) \setminus \II \vert \JJ$, i.e.
\begin{equation*}
\left\vert\begin{array}{ccc} \delta_{\ii_1 \jj_1} & \cdots & \delta_{\ii_1 \jj_s}\\ \vdots & \ddots & \vdots \\ \delta_{\ii_s \jj_1} & \cdots & \delta_{\ii_s \jj_s} \end{array}\right\vert = \varepsilon_{(1\ldots q) \setminus \II \vert \II}\dfrac{\Delta_{(1\ldots q) \setminus \II \vert \JJ}}{\Delta_{(1\ldots q)}}.
\end{equation*}

Now the final step is to determine the vector space $\wedge^s W$ avoiding the hypothesis on the non-vanishing Pl\"ucker coordinate.

\begin{definition}\label{def:genericGens}
Let $\Grass{q}{N}{\K}$ be the Grassmannian of $q$-dimensional subspaces of $V = \langle\vv_1,\ldots,\vv_N\rangle$ with the usual Pl\"ucker embedding $\PE: \Grass{q}{N}{\K} \rightarrow \PP^{\binom{N}{q}-1}_{\K}$. For any point $[\ldots:\Delta_{I}:\ldots] \in \PP^{\binom{N}{q}-1}_{\K}$ and for any $ s < q$, we associate to the ordered multiindex $\JJ = (\jj_1,\ldots,\jj_{q-s})$ the element of $\wedge^s V$
\begin{equation}\label{eq:LambdaSJJ}
\Lambda_{\JJ}^{(s)} = \sum_{\vert\KK\vert = s} \Delta_{\JJ\vert\KK} \, \vv^{(s)}_{\KK}.
\end{equation}
Moreover we denote by $\Gamma^{(s)}$ the set of all possible $\Lambda^{(s)}_{\JJ}$:
\begin{equation}\label{eq:GammaS}
\Gamma^{(s)} = \left\{ \Lambda_{\JJ}^{(s)}\ \big\vert\ \forall\ \JJ \text{ s.t. } \vert\JJ\vert = q-s\right\}.
\end{equation}
\end{definition}

\begin{proposition}\label{prop:genericGens}
Let $\Grass{q}{N}{\K}$ be the Grassmannian of $q$-dimensional subspaces of $V = \langle\vv_1,\ldots,\vv_N\rangle$ with the usual Pl\"ucker embedding $\PE: \Grass{q}{N}{\K} \rightarrow \PP^{\binom{N}{q}-1}_{\K}$ and let us consider $W \in \Grass{q}{N}{\K}$. The vectors of
\[
\Gamma^{(s)}(W) = \left\{ \Lambda_{\JJ}^{(s)}(W) = \sum_{\vert\KK\vert = s} \Delta_{\JJ\vert\KK}(W) \, \vv^{(s)}_{\KK}\quad \bigg\vert\quad \forall\ \JJ \text{ s.t. } \vert\JJ\vert = q-s \right\}
\]
generate $\wedge^s W$.
\end{proposition}

\begin{proof}
First of all, we show that for every $\JJ$, $\Lambda^{(s)}_{\JJ}(W)$ belongs to $\wedge^{s} W$. Let us suppose $\Lambda^{(s)}_{\JJ}(W) \neq 0$, that is there exists a multiindex $\HH,\ \vert\HH\vert = s$ such that the Pl\"ucker coordinate $\Delta_{\II} = \Delta_{\JJ\cup\HH}$ does not vanish. We can write
\[
\begin{split}
\Lambda^{(s)}_{\JJ}(W)&{} = \Delta_{\JJ\vert\HH}\, \vv^{(s)}_{\HH} + \sum_{\begin{subarray}{c} \vert\KK\vert=s\\ \KK\neq\HH\end{subarray}} \Delta_{\JJ\vert\KK}\, \vv^{(s)}_{\KK} =  \Delta_{\JJ\vert\HH} \left(\vv^{(s)}_{\HH} + \sum_{\begin{subarray}{c} \vert\KK\vert=s\\ \KK\neq\HH\end{subarray}} \dfrac{\Delta_{\JJ\vert\KK}}{\Delta_{\JJ\vert\HH}}\, \vv^{(s)}_{\KK}\right) = \\
&{} = \Delta_{\JJ\vert\HH} \left(\vv^{(s)}_{\HH} + \sum_{\begin{subarray}{c} \vert\KK\vert=s\\ \KK\neq\HH\end{subarray}} \dfrac{\Delta_{\JJ\vert\KK}}{\varepsilon_{\JJ\vert\HH}\,\Delta_{\JJ\cup\HH}}\, \vv^{(s)}_{\KK}\right) =\\
&{} = \Delta_{\JJ\vert\HH} \left(\vv^{(s)}_{\HH} + \sum_{\begin{subarray}{c} \vert\KK\vert=s\\ \KK\neq\HH\end{subarray}} \varepsilon_{\II\setminus\HH\vert\HH}\dfrac{\Delta_{\II\setminus\HH\vert\KK}}{\Delta_{\II}}\, \vv^{(s)}_{\KK}\right)
\end{split}
\]
obtaining an element of the type (except for a factor) described in \eqref{eq:basisWs}. This element belongs to $\wedge^s W$ because it can be obtained as exterior product (up to multiplication by a scalar) of $s$ elements of th basis of $W$
\[
\left\{\Lambda^{(1)}_{\JJ'}(W) \ \big\vert\ \JJ' \subset \II,\ \vert\JJ'\vert = q-1\right\}.
\]
and it is easy to check that this basis is of the same type of that one described in \eqref{eq:basisW}.

At this point we know that $\langle \Gamma^{(s)}(W)\rangle \subset \wedge^s W$, so to prove the equality it suffices to show that $\dim_{\K} \langle \Gamma^{(s)}(W)\rangle = \dim_{\K} \wedge^s W = \binom{q}{s}$. In fact the set
\[
\left\{\Lambda^{(s)}_{\KK}(W) \ \big\vert\ \KK \subset \II,\ \vert\KK\vert = s\right\}
\]
contains $\binom{q}{s}$ elements linearly independents, i.e. it represents a basis for $\wedge^s W$.
\end{proof}

\begin{example}\label{ex:MainExampleGens}
We consider again the Grassmannian $\Grass{4}{6}{\K}$ introduced in Example \ref{ex:MainExampleGrass} and we will use the same notation. Obviously $\Gamma^{(4)}$ contains a single element: the vector defining the point in $\PP^{14}_{\K}$: $\Lambda^{(4)}_{\emptyset} = \sum \Delta_{\II} \, \vv^{(4)}_{\II}$. We compute as example an element of $\Gamma^{(1)}$, $\Gamma^{(2)}$ and $\Gamma^{(3)}$.
\begin{itemize}
\item $\Gamma^{(1)},\ \vert\Gamma^{(1)}\vert = \binom{6}{3} = 20$.
\small
\[
\Lambda^{(1)}_{156} = \Delta_{156\vert2}\, \vv_{2} + \Delta_{156\vert3}\, \vv_{3} + \Delta_{156\vert4}\, \vv_{4} = \Delta_{1256}\, \vv_{2} + \Delta_{1356}\, \vv_{3} + \Delta_{1456}\, \vv_{4}.
\]
\normalsize
\item $\Gamma^{(2)},\ \vert\Gamma^{(2)}\vert = \binom{6}{2} = 15$.
\small
\[
\begin{split}
\Lambda^{(2)}_{24}&{} = \Delta_{24\vert13}\, \vv^{(2)}_{13} + \Delta_{24\vert15}\, \vv^{(2)}_{15} + \Delta_{24\vert16}\, \vv^{(2)}_{16} + \Delta_{24\vert35}\, \vv^{(2)}_{35} + \Delta_{24\vert36}\, \vv^{(2)}_{36} + \Delta_{24\vert56}\, \vv^{(2)}_{56}= \\
&{} = -\Delta_{1234}\, \vv^{(2)}_{13} + \Delta_{1245}\, \vv^{(2)}_{15} + \Delta_{1246}\, \vv^{(2)}_{16} - \Delta_{2345}\, \vv^{(2)}_{35} - \Delta_{2346}\, \vv^{(2)}_{36} + \Delta_{2456}\, \vv^{(2)}_{56}.
\end{split}
\]
\normalsize
\item $\Gamma^{(3)},\ \vert\Gamma^{(3)}\vert = \binom{6}{1} = 6$.
\small
\[
\begin{split}
\Lambda^{(3)}_{3}&{} = \Delta_{3\vert124} \, \vv^{(3)}_{124} + \Delta_{3\vert125} \, \vv^{(3)}_{125} + \Delta_{3\vert126} \, \vv^{(3)}_{126} + \Delta_{3\vert145} \, \vv^{(3)}_{145} + \Delta_{3\vert146} \, \vv^{(3)}_{146} + \\
&\qquad{} + \Delta_{3\vert156} \, \vv^{(3)}_{156} + \Delta_{3\vert245} \, \vv^{(3)}_{245} + \Delta_{3\vert246} \, \vv^{(3)}_{246} + \Delta_{3\vert256} \, \vv^{(3)}_{256} + \Delta_{3\vert456} \, \vv^{(3)}_{456} = \\
&{} = \Delta_{1234} \, \vv^{(3)}_{124} + \Delta_{1235} \, \vv^{(3)}_{125} + \Delta_{1236} \, \vv^{(3)}_{126} - \Delta_{1345} \, \vv^{(3)}_{145} - \Delta_{1346} \, \vv^{(3)}_{146} + \\
&\qquad{} - \Delta_{1356} \, \vv^{(3)}_{156} - \Delta_{2345} \, \vv^{(3)}_{245} - \Delta_{2346} \, \vv^{(3)}_{246} - \Delta_{2356} \, \vv^{(3)}_{256} + \Delta_{3456} \, \vv^{(3)}_{456}.
\end{split}
\]
\normalsize
\end{itemize}
For the complete lists of elements in $\Gamma^{(1)}$, $\Gamma^{(2)}$ and $\Gamma^{(3)}$ see Example \ref{ex:MainExampleGrassCODE}.
\end{example}

\begin{remark}
In the paper \textit{\lq\lq Low degree equations defining the Hilbert scheme\rq\rq} \cite{BLMR}, the property stated in Proposition \ref{prop:genericGens} is proved using a different approach, more abstract, based on some results of a paper \cite{BarnabeiBriniRota} by Barnabei, Brini and Rota. For our computational purpose, we prefer this practical, and in some sense algorithmic, description.
\end{remark}

\begin{remark}
As Example \ref{ex:MainExampleGens} suggests, the number of elements in $\Gamma^{(s)}(W)$ is bigger than the dimension of $\wedge^s W$, indeed
\[
\dim_{\K} \wedge^s W = \binom{q}{s} \quad\text{and}\quad \left\vert \Gamma^{(s)}(W) \right\vert = \binom{N}{q-s}
\]
and, set $p = N-q$,
\[
\begin{split}
\binom{N}{q-s}&{} = \binom{N-1}{q-s} + \binom{N-1}{q-s-1} =\\
&{}= \binom{N-2}{q-s} + \binom{N-2}{q-s-1} + \binom{N-1}{q-s-1} =\\
&{}\qquad \vdots\\
&{}= \binom{N-p}{q-s} + \sum_{i=1}^p \binom{N-i}{q-s-1} = \binom{q}{s} + \sum_{i=1}^p \binom{N-i}{q-s-1}.
\end{split}
\]
\end{remark}
\index{Grassmannian|)}

\subsection{The Grassmann functor}

Classically the same construction is extended to the more general setting of $\K$-algebras and direct summands of $A^{n}$ of rank $q$.

It is therefore natural to introduce the following functor.
\begin{definition}[{\cite[Exercise VI-18]{EisenbudHarris}}]\label{def:grassmanFunctor}\index{Grassmann functor}
For any couple of integer $(q,N),\ 0 < q < N$, the \emph{Grassmann functor} is the functor
\begin{eqnarray*}
& \gls{grassFun}: \Kalg \rightarrow \Sets&\\
& \grass{q}{N}(A) = \left\{\text{rank } q \text{ direct summands of } A^N\right\}&
\end{eqnarray*}
and for any $f: A \rightarrow B \in \Mor_{\Kalg}$ the map $\grass{q}{N}(f)$ is defined as
\[
\begin{split}
 \grass{q}{N}(A) & {} \rightarrow\ \grass{q}{N}(B)\\
\parbox{1.5cm}{\centering $P$} &{} \mapsto\ P \otimes_A B
\end{split}
\]
by means of the extension of scalars (see \cite[Proposition 2.17]{AtiyahMacdonald}) and using the right exactness of the tensor product.
\end{definition}

\begin{theorem}\label{th:representabilityGrassman}\index{Grassmannian}
The Grassmannian $\Grass{q}{N}{\K}$ represents the functor $\grass{q}{N}$.
\end{theorem}

To prove this result we need to say some more words about the Grassmannian. Let us consider the ring homomorphism
\[
f: \K[\ldots,\Delta_{\JJ},\ldots] \rightarrow \K[\ldots,\delta_{\ii\jj},\ldots]
\]
where variables $\Delta_{\JJ}$ are indexed over ordered subsets $\JJ \subset \{1,\ldots,N\},\ \vert\II\vert = q$ and variables $\delta_{\ii\jj}$ over $1\leqslant \ii \leqslant q$, $1\leqslant \jj \leqslant N$, such that $f(\Delta_{\JJ}) = \det(\delta_{\ii\jj})\vert_{\jj\in\JJ}$. $f$ induces the affine homomorphism
\[
f^\ast: \AA^{qN}_{\K} = \Spec \K[\ldots,\delta_{\ii\jj},\ldots] \rightarrow \AA^{\binom{N}{q}}_{\K} = \Spec \K[\ldots,\Delta_{\JJ},\ldots],
\]
and considering its restriction to the Zariski open subset 
\[
U = \AA^{qN}_{\K}\setminus \glslink{ZariskiClosed}{Z(\langle \det(\delta_{\ii\jj})\vert_{\jj\in\JJ},\ \forall\ \vert\JJ\vert=q \rangle)}
\]
of the matrices $(\delta_{\ii\jj})$ of rank $q$, also the morphism
\[
\phi: U \rightarrow \PP^{\binom{N}{q}-1}_{\K} = \Proj \K[\ldots,\Delta_{\JJ},\ldots]
\]
turns out to be well defined.

\begin{proposition}\label{prop:phiFactors}
The morphism $\phi: U \rightarrow \PP^{\binom{N}{q}-1}_{\K}$ factors through $\Grass{q}{N}{\K}$:
\begin{center}
\begin{tikzpicture}[]
\node (U) at (0,0) [] {$U$};
\node (PP) at (3,0) [] {$\PP^{\binom{N}{q}-1}_{\K}$};
\node (G) at (1.5,-1.5) [] {$\Grass{q}{N}{\K}$};
\draw [->] (U) --node[rectangle,fill=white]{\footnotesize $\phi$} (PP);
\draw [right hook->] (G) -- (PP);
\draw [->,dashed] (U) -- (G);
\end{tikzpicture}
\end{center}
\end{proposition}
\begin{proof}
Firstly $\phi$ is equivariant with respect the action of $\GL_{\K}(q)$ over $U$ given by left multiplication $\mu: \GL_{\K}(q) \times U \rightarrow U$, that is
\begin{center}
\begin{tikzpicture}[>=latex]
\node (GL) at (0,0) [] {$\GL_{\K}(q) \times U$};
\node (U1) at (3,0) [] {$U$};
\node (U2) at (0,-2) [] {$U$};
\node (PP) at (3,-2) [] {$\PP_{\K}^{\binom{N}{q}-1}$};
\draw [->] (GL) --node[rectangle,fill=white]{\footnotesize $\mu$} (U1);
\draw [->] (GL) --node[rectangle,fill=white]{\footnotesize $\pi_{U}$} (U2);
\draw [->] (U1) --node[rectangle,fill=white]{\footnotesize $\phi$} (PP);
\draw [->] (U2) --node[rectangle,fill=white]{\footnotesize $\phi$} (PP);
\node at (1.5,-1) [] {$\circlearrowright$};
\end{tikzpicture}
\end{center}

Then for any $\JJ\ (\vert\JJ\vert = q)$, let $U_{\JJ}$ the closed subscheme of $U$ such that the submatrix $(\delta_{\ii\jj})\vert_{\jj\in\JJ}$ is equal to the identity matrix. Obviously $U_{\JJ} \simeq \AA_{\K}^{q(N-q)}$. Moreover consider the standard open covering on $\PP^{\binom{N}{q}-1}_{\K}$ 
\[
\bigcup_{\JJ}\, \mathcal{U}_{\JJ}, \qquad \mathcal{U}_{\JJ} = \PP^{\binom{N}{q}-1}_{\K} \setminus Z(\Delta_{\JJ}).
\]
To prove that $\phi$ factors through $\Grass{q}{N}{\K}$ we remark that
\begin{enumerate}[(i)]
\item $\mu\vert_{U_{\JJ}}: \GL_{\K}(q) \times U_{\JJ} \rightarrow \phi^{-1}(\mathcal{U}_{\JJ})$ is an isomorphism;
\item $\phi\vert_{U_{\JJ}}: U_{\JJ} \rightarrow \mathcal{U}_{\JJ} \cap \Grass{q}{N}{\K}$ is an isomorphism;
\item for\hfill each\hfill Pl\"ucker\hfill relation\hfill $Q$\hfill of\hfill the\hfill set\hfill in\hfill \eqref{eq:PluckerEquations},\hfill the\hfill image\hfill of\hfill $Q/\Delta_{\JJ}^2$\hfill in\\ $\K[\ldots,\delta_{\ii\jj},\ldots]$ belongs to the ideal $I(U_{\JJ})$ defining $U_{\JJ}$. \qedhere
\end{enumerate}
\end{proof}

We rewrite the statement of Theorem \ref{th:representabilityGrassman} through the following proposition.
\begin{proposition}There exists an invertible natural transformation
\begin{equation}
\mathscr{F}: \grass{q}{N} \rightarrow h_{\Grass{q}{N}{\K}}.
\end{equation}
\end{proposition}

\begin{proof}
Let\hfill us\hfill consider\hfill an\hfill element\hfill $P \in \grass{q}{N}(A)$\hfill and\hfill let\hfill us\hfill define\hfill a\hfill morphism\\ $\Spec A \rightarrow \Grass{q}{N}{\K}$, i.e. an element in $h_{\Grass{q}{N}{\K}}(A)$. Being $P$ a direct summand of $A^N$ of rank $q$, the injective map
\[
i_{P} : P \hookrightarrow A^N
\]
is described by a $N\times q$ matrix with coefficients in $A$. Therefore the chain of morphism
\begin{equation}\label{eq:chainRepresentabilityGr}
\K[\ldots,\delta_{\ii\jj},\ldots]\ \longrightarrow\ A[\ldots,\delta_{\ii\jj},\ldots] =\K[\ldots,\delta_{\ii\jj},\ldots]\otimes_{\K} A\ \longrightarrow\ A,
\end{equation}
where the second map is the evaluation map over the transposed matrix ${}^T i_P$, 
induces a morphism $f_P:\Spec A \rightarrow \AA_{\K}^{qN}$ that factors over $U$. Finally $\phi \circ f_{P}$ factors over $\Grass{q}{N}{\K}$ giving a morphism $\phi\circ f_p: \Spec A \rightarrow \Grass{q}{N}{\K} \in h_{\Grass{q}{N}{\K}}(A)$. The transformation between morphisms of the two category follows directly from \eqref{eq:chainRepresentabilityGr} and extension of scalars.

To invert $\mathscr{F}$, we look for the universal family over $\Grass{q}{N}{\K}$, that is we want to construct a sub-bundle $\mathcal{K} \hookrightarrow \OO^N_{\Grass{q}{N}{\K}}$. We consider again the action $\mu: \GL_{\K}(q)  \times U \rightarrow U$ and let $\mathcal{V}_{\II} = \mathcal{U}_{\II} \cap \Grass{q}{N}{\K}$ and $\mathcal{V}_{\JJ} = \mathcal{U}_{\JJ} \cap \Grass{q}{N}{\K}$. By the property (i) exposed in the proof of Proposition \ref{prop:phiFactors}, $\mathcal{V}_{\II} \simeq U_{\II} \simeq \AA_{\K}^{q(N-q)}$ and $\mathcal{V}_{\JJ} \simeq U_{\JJ} \simeq \AA_{\K}^{q(N-q)}$ and using the property (iii) we can define the isomorphism
\[
\rho_{\II\JJ}: \GL_{\K}(q) \times (\mathcal{V}_{\II}\cap \mathcal{V}_{\JJ}) \stackrel{\sim}{\longrightarrow} \phi^{-1}(\mathcal{V}_{\II}\cap\mathcal{V}_{\JJ})  \stackrel{\sim}{\longrightarrow} \GL_{\K}(q) \times (\mathcal{V}_{\II}\cap \mathcal{V}_{\JJ})
\]
where the intersection on the left $\mathcal{V}_{\II}\cap \mathcal{V}_{\JJ}$ is considered as open subset of $\mathcal{V}_{\II}$ whereas the intersection on the right as open subset of $\mathcal{V}_{\JJ}$. Thus the morphism to $\GL_{\K}(q)$ induced by $\rho_{\II\JJ}$ determines a transition function $g_{\II\JJ}: \mathcal{V}_{\II}\cap\mathcal{V}_{\JJ} \rightarrow \GL_{K}(q)$.

For every $\II$, the embedding $\mathcal{V}_{\II} \simeq U_{\II} \simeq \AA^{q(N-q)}_{\K} \hookrightarrow \AA^{qN}_{\K}$ induces the map $\OO^q_{\mathcal{V}_{\II}} \rightarrow \OO^N_{\mathcal{V}_{\II}}$. This collection of maps glue together by means of the transition functions $g_{\II\JJ}$ to the injective map $\mathcal{K} \rightarrow \OO^N_{\Grass{q}{N}{\K}}$, whose cokernel turns out to be a locally free sheaf of rank $N-q$.

Finally for every morphism $f: \Spec A \rightarrow \Grass{q}{N}{\K} \in h_{\Grass{q}{N}{\K}}(A)$, we define an element $\mathscr{F}^{-1}(f) \in \grass{q}{N}(A)$ starting from the family
\[
0\ \longrightarrow\ \mathcal{K}\ \longrightarrow\ \OO^N_{\Grass{q}{N}{\K}}\ \longrightarrow\ \mathcal{Q}\ \longrightarrow\ 0
\]
via pull-back
\[
0\ \longrightarrow\ f^{\ast}\mathcal{K}\ \longrightarrow\ f^\ast\OO^N_{\Grass{q}{N}{\K}} \simeq A^N\ \longrightarrow\ f^\ast\mathcal{Q}\ \longrightarrow\ 0. \qedhere
\]
\end{proof}

By Yoneda's Lemma (Theorem \ref{th:Yoneda}) the schemes representing $\grass{q}{N}$ and $h_{\Grass{q}{N}{\K}}$ have to be isomorphic, so $\Grass{q}{N}{\K}$ represents also $\grass{q}{N}$.

\section{The Hilbert functor}

\begin{definition}\label{def:hilbertFunctor}\index{Hilbert functor}
Let us define the \emph{Hilbert functor} as the covariant functor
\begin{equation*}
\gls{hilbFun} : \schm_{\K} \rightarrow \Sets
\end{equation*}
such that for $X\in \Ob_{\schm_{\K}}$
\begin{equation*}
\hilb{n}{}(X) = \left\{\mathcal{Z} \subset \PP^n \times X \ \vert\ \mathcal{Z} \text{ flat over } X \text{ via } \pi: \mathcal{Z} \hookrightarrow \PP^n \times X \rightarrow X \right\}
\end{equation*}
and for $f: X \rightarrow Y \in \Mor_{\schm_{\K}}$
\[
\begin{split}
\hilb{n}{}(f): \hilb{n}{}(Y) &\rightarrow \hilb{n}{}(X)\\
\parbox{1.5cm}{\centering $\mathcal{Z}$} & \mapsto\ f^\ast(\mathcal{Z})
\end{split}
\]
that is well defined because the pullback preserves the flatness:
\begin{center}
\begin{tikzpicture}[]
\node (fZ) at (0,0) [] {$f^\ast(\mathcal{Z})$};
\node (PX) at (2.5,0) [] {$\PP^n_{\K} \times X$};
\node (X) at (5,0) [] {$X$};
\draw [->] (fZ) -- (PX);
\draw [->] (PX) -- (X);
\node (Z) at (0,-1.8) [] {$\mathcal{Z}$};
\node (PY) at (2.5,-1.8) [] {$\PP^n_{\K} \times Y$};
\node (Y) at (5,-1.8) [] {$Y$};
\draw [->] (Z) -- (PY);
\draw [->] (PY) -- (Y);
\draw [->] (fZ) -- (Z);
\draw [->] (X) --node[rectangle,fill=white,inner sep=1pt]{\footnotesize $f$} (Y);
\draw [->] (PX) --node[rectangle,fill=white,inner sep=1pt]{\footnotesize $\id \times f$} (PY);
\end{tikzpicture}
\end{center}
\end{definition}

For the properties of flatness, we know that the Hilbert polynomial is locally constant for $x \in X$, so we can decompose the functor $\hilb{n}{}$ as
\[
\hilb{n}{} = \coprod_{p(t) \in \QQ[t]} \hilb{n}{p(t)}.
\]
where $p(t) \in \QQ[t]$ is a numerical polynomial and $\hilb{n}{p(t)}$ is the subfunctor of $\hilb{n}{}$ such that
\[
\gls{hilbFunPoly}(X) = \left\{\mathcal{Z} \subset \PP^n_{\K} \times X \ \bigg\vert\ \begin{array}{l}\mathcal{Z} \text{ flat over } X \text{ with fibers}\\ \text{having Hilbert polynomial } p(t)\end{array}\right\}.
\]

\begin{theorem}[{Grothendieck \cite{GrothendieckFGA}}]
The functor $\hilb{n}{p(t)}$ is representable.
\end{theorem}

\begin{definition}\label{def:hilbertScheme}\index{Hilbert scheme}
We call \emph{Hilbert scheme} and denote by \gls{hilbScheme} the scheme representing the functor $\hilb{n}{p(t)}$.
\end{definition}

There is another useful meaning of the representability of a functor.
\begin{proposition}\label{prop:universalFamily}
The Hilbert scheme $\Hilb{n}{p(t)}$ representing $\hilb{n}{p(t)}$ is the parameter scheme of a flat family $\mathcal{X}$ with Hilbert polynomial $p(t)$
\begin{center}
\begin{tikzpicture}[>=latex]
\node (X) at (0,0) [] {$\mathcal{X}$};
\node (PX) at (2.5,0) [] {$\PP^n_{\K} \times \Hilb{n}{p(t)}$};
\node (H) at (2.5,-1.5) [] {$\Hilb{n}{p(t)}$};
\draw [->] (X) -- (PX);
\draw [->] (X) -- (H);
\draw [->] (PX) -- (H);
\end{tikzpicture}
\end{center}
such that any subscheme $\mathcal{Y} \subset \PP^n_{\K} \times S$ flat over $S$ with Hilbert polynomial $p(t)$ coincides with the fiber product $\mathcal{X} \times_{\Hilb{n}{p(t)}} S\subset \PP^n_{\K} \times S$ for a unique map $S \rightarrow \Hilb{n}{p(t)}$:
\begin{center}
\begin{tikzpicture}[>=latex]
\node (Y) at (0,0) [] {$\mathcal{Y} = \mathcal{X} \times_{\Hilb{n}{p(t)}} S$};
\node (X) at (2.5,0) [] {$\mathcal{X}$};
\node (S) at (0,-1.5) [] {$S$};
\node (H) at (2.5,-1.5) [] {$\Hilb{n}{p(t)}$};
\draw [->] (Y) -- (X);
\draw [->] (Y) -- (S);
\draw [->] (X) -- (H);
\draw [->] (S) -- (H);
\end{tikzpicture}
\end{center}
The family $\mathcal{X}\rightarrow\Hilb{n}{p(t)}$ is called \emph{universal family}\index{universal family} and we will refer to this property as \emph{universal property}\index{Hilbert scheme!universal property of the} of the Hilbert scheme.
\end{proposition}

To simplify the discussion, we use again Proposition \ref{prop:reductionAffine} for rewriting the Hilbert functor as a functor over the category of affine schemes: 
\begin{eqnarray*}
& \glshyperlink{hilbFunPoly}: \affschm_{\K} \rightarrow \Sets, &\\
&\hilb{n}{p(t)}(\Spec A) = \left\{\mathcal{Z} \subset \PP^n_{A}\ \vert\ \mathcal{Z} \text{ flat over } \Spec A \text{ with Hilbert polynomial } p(t)\right\}.& 
\end{eqnarray*}

It is well known that saying $\mathcal{Z}$ flat over $\Spec A$ means that the structure sheaf $\OO_{\mathcal{Z}}$ is flat over $\Spec A$. Being in the affine case we can consider the graded module
\[
M = \bigoplus_{t \geqslant 0} H^0 (\mathcal{Z},\OO_{\mathcal{Z}}(t)) \qquad (\text{s.t. } \widetilde{M} = \OO_{\mathcal{Z}})
\]
over $A[x_0,\ldots,x_n]$ and $\OO_{\mathcal{Z}}$ flat over $\Spec A$ is equivalent to $M$ flat over $A$ (see \cite[Chapter III Proposition 9.2]{Hartshorne}). The flatness is preserved by localization (\cite[Chapter III Proposition 9.1]{Hartshorne}), i.e. $M$ is flat over $A$ if and only if $M_{\mathfrak{p}}$ is flat over $A_{\mathfrak{p}}$, for all $\mathfrak{p} \in \Spec A$. Moreover for any local algebra $A_{\mathfrak{p}}$, the flatness of a finite $A_{\mathfrak{p}}$-module is equivalent to its freeness (\cite[Proposition 3.G]{MatsumuraCA}), so denoted by $k(\mathfrak{p})$ the residue field of $A_{\mathfrak{p}}$, the Hilbert polynomial $p(t)$ of $M_{\mathfrak{p}}$ is well defined as the Hilbert polynomial of the module $M_{\mathfrak{p}} \otimes k(\mathfrak{p})$, that is
\begin{equation}
p(t) = \dim_{k(\mathfrak{p})} \left(M_{\mathfrak{p}}\right)_t \otimes k(\mathfrak{p}),\qquad t \gg 0. 
\end{equation}
Finally the flatness ensures that the Hilbert polynomial does not depend on the point $\mathfrak{p} \in \Spec A$ for which we localize (\cite[Chapter III Theorem 9.9]{Hartshorne} and \cite[Exercise 6.11]{EisenbudHarris}).

Hence we can redefine the Hilbert functor as follows
\begin{eqnarray*}
&\hilb{n}{p(t)}: \Kalg \rightarrow \Sets&\\
&\hilb{n}{p(t)}(A) = \left\{ \begin{array}{c}M,\text{ graded module over } A[x_0,\ldots,x_n] \text{ flat over } A\\ \text{with Hilbert polynomial } p(t)\end{array}\right\} &\\
\end{eqnarray*}
and for any $f: A \rightarrow B$, 
\[
\begin{split}
\hilb{n}{p(t)}(f): \hilb{n}{p(t)}(A) & \rightarrow \hilb{n}{p(t)}(B)\\
\parbox{2cm}{\centering $M$} & \mapsto \parbox{2cm}{\centering $M \otimes_{A} B$}
\end{split}
\]
using the fact that the extension of scalars preserves the flatness (see \cite[(3.C)]{MatsumuraCA}).

\begin{definition}\label{def:Avaluedpoints}
The elements of the set $\hilb{n}{p(t)}(A)$ are called $A$-valued points of the Hilbert scheme $\Hilb{n}{p(t)}$. With this terminology, the universal property of $\Hilb{n}{p(t)}$ can be described saying that any flat family $\mathcal{Y} \rightarrow \Spec A$ defines a unique $A$-valued point of $\Hilb{n}{p(t)}$.
\end{definition}

\begin{remark}
The set of $\K$-valued points $\hilb{n}{p(t)}(\K)$ contains exactly the subschemes $X \subset \PP^n_{\K}$ with Hilbert polynomial $p(t)$, being $\Spec \K$ a single point. 
\end{remark}

\section{The Hilbert scheme as subscheme of the Grassmannian}

\begin{definition}
An admissible Hilbert polynomial $p(t)$ has a unique \emph{Gotzmann representation}\index{Hilbert polynomial!Gotzmann's representation of a}
\begin{equation}\label{eq:GotzmannDecomposition}
p(t) = \binom{t+a_1}{a_1} + \binom{t+ a_2 -1}{a_2} + \ldots + \binom{t+a_r - (r-1)}{a_r}, \qquad a_1 \geqslant \cdots \geqslant a_r.
\end{equation}
The number $r$ of terms in this sum is said \emph{Gotzmann number}\index{Gotzmann number}\index{Hilbert polynomial!Gotzmann number of a|see{Gotzmann number}} of $p(t)$.
\end{definition}

\begin{example}
The Hilbert polynomial $p(t) = 4t$ has Gotzmann number equal to 6, indeed
\[
4t = \binom{t+1}{1} + \binom{t}{1} + \binom{t-1}{1} + \binom{t-2}{1} + \binom{t-4}{0} + \binom{t-5}{0}.
 \]
\end{example}

\begin{theorem}[\textbf{Gotzmann's Regularity Theorem} {\cite[Theorem 3.11]{GreenGIN}}]\label{th:RegularityTheorem}\index{Gotzmann's Regularity Theorem}\index{regularity of a sheaf}
Let $A$ be any $\K$-algebra and let $Z \subset \Proj A[x_0,\ldots,x_n]$ be any subscheme with Hilbert polynomial $p(t)$, whose Gotzmann number is $r$. Then the sheaf of ideals $\mathcal{I}_{Z}$ is $r$-regular.
\end{theorem}

As usual to any subscheme $Z \subset \PP^{n}_A$, we can associate the saturated ideal $\glslink{satIdeal}{I(Z)} = \bigoplus_t H^0(Z,\mathcal{I}_Z(t))$. Gotzmann's Regularity Theorem ensures that the truncated ideal $I(Z)_{\geqslant r}$ is generated by its homogenous piece of degree $r$, that is
\begin{equation}
I(Z)_{\geqslant r} = \left\langle I(Z)_r \right\rangle.
\end{equation}

This result suggests to associate to any subscheme $Z \subset \Proj A[x_0,\ldots,x_n]$ the truncation $I(Z)_{\geqslant r}$ instead of the saturated ideal $I(Z)$, with the main advantage of a more uniform description, indeed with this approach any ideal associated to such subschemes is generated by the same number of linearly independent polynomials of the same degree.
 
For any $M \in \hilb{n}{p(t)}(A)$, called $Z$ the affine subscheme flat over $\Spec A$ such that $\widetilde{M} = \mathcal{O}_{Z}$, $M = A[x_0,\ldots,x_n]/I(Z)$. We are interested in the homogeneous piece of degree $r$: $M_r = A[x_0,\ldots,n_n]_r/I(Z)_r$. For each $\mathfrak{p} \in \Spec A$, we have already said that $\left(M_{\mathfrak{p}}\right)_r$ is a free $A_{\mathfrak{p}}$-module and this property implies that $M_{r}$ has to be a projective $A$-module (\cite[Theorem A3.2]{EisenbudHarris}). As explained in Appendix C of \cite[pp. 302-303]{IarrobinoKanev} (see also \cite[Section A3.2]{EisenbudHarris}), being $A[x_0,\ldots,n_n]_r/I(Z)_r$ projective ensures that $I(Z)_r$ is a direct summand of $A[x_0,\ldots,x_n]_r$. 

Therefore, set $N(t) = \binom{n+t}{n}$ and $q(t) = N(t)-p(t)$, we can define a natural transformation of functors from $\hilb{n}{p(t)}$ and $\grass{q(r)}{N(r)}$
\begin{equation}\label{eq:hilbToGrass}
\begin{array}{ccc}
\hilb{n}{p(t)}(A) & \rightarrow & \grass{q(r)}{N(r)}(A) \\
M = A[x_0,\ldots,x_n]/I(Z) & \mapsto & I(Z)_r \subset A^{N(r)} \simeq A[x_0,\ldots,x_n]_r
\end{array}
\end{equation}
that suggests to determine the Hilbert scheme $\Hilb{n}{p(t)}$ representing $\hilb{n}{p(t)}$ as subscheme of the Grassmannian $\Grass{q(r)}{N(r)}{\K}$ (representing $\grass{q(r)}{N(r)}$). To accomplish this purpose, we need to understand under which conditions an ideal $I = \langle I_r\rangle \subset A[x_0,\ldots,x_n]$ generated by $q(r)$ linearly independent homogeneous polynomials of degree $r$ determines a module $A[x_0,\ldots,x_n]/I$ with Hilbert polynomial $p(t)$. One of the best characterization is stated in the following theorem.

\begin{theorem}[\textbf{Gotzmann's Persistence Theorem} {\cite[Theorem C.17]{IarrobinoKanev}}]\label{th:PersistenceTheorem}\index{Gotzmann's Persistence Theorem} Let $I \subset A[x_0,\ldots,x_n]$ be a homogeneous ideal generated by its piece of degree $r$, i.e. $I = \langle I_r\rangle$. If $A[x_0,\ldots,x_n]_r/I_r$ and $A[x_0,\ldots,x_n]_{r+1}/I_{r+1}$ are flat $A$-modules of rank $p(r)$ and $p(r+1)$, then $A[x_0,\ldots,x_n]_t/I_t$ is a $A$-flat module of rank $p(t)$ for all $t\geqslant r$. 
\end{theorem}

There is another important results by Macaulay that further simplifies the characterization.

\begin{theorem}[\textbf{Macaulay's Estimate on the Growth of Ideals} {\cite[Theorem 3.3]{GreenGIN}}]\label{th:MacaulayEstimate}\index{Macaulay's Estimate on the Growth of Ideals}
Let $I \subset A[x_0,\ldots,x_n]$ be a homogeneous ideal and let $p(t)$ be an admissible Hilbert polynomial. If $\rank A[x_0,\ldots,x_n]_r/I_r = p(r)$, then $\rank A[x_0,\ldots,x_n]_{r+1}/I_{r+1} \leqslant p(r+1)$.
\end{theorem}

Putting together Theorem \ref{th:PersistenceTheorem} and Theorem \ref{th:MacaulayEstimate}, the condition we have to impose on $I \subset A[x_0,\ldots,x_n]$ is $\rank A[x_0,\ldots,x_n]_{r+1}/I_{r+1} \geqslant p(r+1)$ or $\rank I_{r+1} \leqslant q(r+1)$.

\section{Known sets of equations}

In this section we will consider a fixed Hilbert polynomial $p(t)$ with Gotzmann number $r$, a fixed projective space $\PP^n_{\K} = \Proj \K[x_0,\ldots,x_n]$. Let $N = N(r)$ and $q = q(r) = N(r)-p(r)$. In the previous section we embedded the Hilbert scheme $\Hilb{n}{p(t)}$ in the Grassmannian $\Grass{q}{N}{\K}$ parametrizing the vector subspaces of dimension $q$ in the base vector space $\K[x_0,\ldots,x_n]_r \simeq \K^{N}$.

To determine equations defining $\Hilb{n}{p(t)}$, we will use tools introduced in Section \ref{sec:grassmannians} and for this reason we slightly modify the notation, adapting it to this special case. We consider $\K[x_0,\ldots,x_n]_r$ equipped with its standard monomial basis $\{x^\beta \text{ s.t. } \vert\beta\vert=r\}$ and the bijective function
\begin{equation}\label{eq:monomialCorrespondence}
\begin{array}{c}
\alpha: \left\{1,\ldots,N\right\} \rightarrow \left\{(a_0,\ldots,a_n) \in \ZZ^{n+1}\ \Big\vert\ a_i \geqslant 0,\ \forall\ i \text{ and } \sum_{i=0}^n a_i = r \right\}\\
x^{\alpha(1)} >_{\texttt{DegRevLex}} x^{\alpha(2)} >_{\texttt{DegRevLex}} \cdots >_{\texttt{DegRevLex}} x^{\alpha(N)}
\end{array}
\end{equation}
so that the Pl\"ucker coordinate $\Delta_{\II}$ corresponds to the vector $\xx^{(s)}_{\II} = x^{\alpha(\ii_1)} \wedge \cdots \wedge x^{\alpha(\ii_{q})}$ of the basis of $\wedge^{q}\, \K[x_0,\ldots,x_n]_r$ and for all $\JJ \subset \{1,\ldots,N\},\ \vert\JJ\vert = q-s$,
\[
\Lambda^{(s)}_{\JJ} = \sum_{\vert\KK\vert=s} \Delta_{\JJ\vert\KK}\, \xx^{(s)}_{\KK} = \sum_{\vert\KK\vert=s} \Delta_{\JJ\vert\KK}\, x^{\alpha(\kk_1)}\wedge\cdots\wedge x^{\alpha(\kk_s)}.
\]
Moreover we define for each variable $x_i$
\begin{eqnarray}
&& x_i\Lambda^{(s)}_{\JJ} = \sum_{\vert\KK\vert=s} \Delta_{\JJ\vert\KK}\, x_i \xx^{(s)}_{\KK} = \sum_{\vert\KK\vert=s} \Delta_{\JJ\vert\KK}\, x_i x^{\alpha(\kk_1)}\wedge\cdots\wedge x_i x^{\alpha(\kk_s)},\\
&& x_i\Gamma^{(s)} = \left\{ x_i \Lambda^{(s)}_{\JJ}\ \big\vert\ \forall\ \JJ \subset \{1,\ldots,N\},\ \vert\JJ\vert = q-s\right\}.
\end{eqnarray}

\subsection{Gotzmann equations}

Moving from the Persistence Theorem (Theorem \ref{th:PersistenceTheorem}), the idea of Gotzmann \cite{Gotzmann} was to consider the natural transformation of functors $\hilb{n}{p(t)} \rightarrow \grass{q(r)}{N(r)} \times \grass{q(r+1)}{N(r+1)}$
\begin{equation}
\begin{array}{ccc}
 \hilb{n}{p(t)}(A) &\rightarrow& \grass{q(r)}{N(r)}(A) \times \grass{q(r+1)}{N(r+1)}(A)\\
 M = A[x]/I(Z) &\mapsto&  I(Z)_r \times I(Z)_{r+1}\\
 \end{array}
\end{equation}
in order to translate the Hilbert functor as
\begin{equation}
\hilb{n}{p(t)}(A) = \left\{ \begin{array}{c}(I,J) \in \grass{q(r)}{N(r)}(A) \times \grass{q(r+1)}{N(r+1)}(A) \text{ s.t.}\\ 
 I \subset A[x]_{r},\ J \subset A[x]_{r+1} \text{ and } I \cdot A[x]_1 = J \end{array}\right\}
\end{equation}
and applying Theorem \ref{th:MacaulayEstimate} we can write
\begin{equation}
\hilb{n}{p(t)}(A) = \left\{ \begin{array}{c}(I,J) \in \grass{q(r)}{N(r)}(A) \times \grass{q(r+1)}{N(r+1)}(A)\text{ s.t.} \\ I \subset A[x]_{r},\ J \subset A[x]_{r+1} \text{ and } I \cdot A[x]_1 \subset J \end{array}\right\}.
\end{equation}

\begin{theorem}[{\cite[Proposition C.28, Theorem C.29]{IarrobinoKanev}}]\label{th:GotzmannEquations}
The Hilbert scheme $\Hilb{n}{p(t)}$ is defined by quadric equations in the product of Grassmannians 
\[
\Grass{q(r)}{N(r)}{\K}\times\Grass{q(r+1)}{N(r+1)}{\K}
\]
 and then can be embedded in $\Grass{q(r)}{N(r)}{\K}$ through the projection on the first factor.
\end{theorem}
\begin{proof}
We want to find closed conditions to describe the set
\[
\big\{ (I,J) \in \Grass{q(r)}{N(r)}{\K}\times\Grass{q(r+1)}{N(r+1)}{\K}\ \vert\ I\cdot \K[x_0,\ldots,x_n]_1 \subset J\big\}.
\]
Through the isomorphism
\[
\begin{split}
\Grass{q(r+1)}{N(r+1)}{\K} &\rightarrow \Grass{p(r+1)}{N(r+1)}{\K}\\
\parbox{1.5cm}{\centering $J$} &\mapsto \parbox{1.5cm}{\centering $J^{\perp}$}
\end{split}
\]
we redefine the condition as
\[
\big\{ (I,J) \in \Grass{q(r)}{N(r)}{\K}\times\Grass{p(r+1)}{N(r+1)}{\K}\ \vert\ I\cdot \K[x]_1 \subset J^{\perp}\big\}.
\]
We consider the Grassmannian $\Grass{q(r)}{N(r)}{\K}$ equipped with the Pl\"ucker coordinates $[\ldots,\Delta_{\II},\ldots]$. The ideal $I$ is generated by the set $\Gamma^{(1)}$ and so
\[
I \cdot \K[x]_1 = \left\langle x_i\Gamma^{(1)} \ \big\vert\ i=0,\ldots,n\right\rangle.
\]
Whereas\hfill for\hfill $\Grass{p(r+1)}{N(r+1)}{\K}$,\hfill we\hfill denote\hfill the\hfill Pl\"ucker\hfill coordinates\hfill by\\ $[\ldots,\nabla_{\JJ},\ldots]$, the generators of $J$ by $\widetilde{\Lambda}^{(1)}_{\KK}$ and the whole set of them by $\widetilde{\Gamma}^{(1)}$.

The condition $I\cdot \K[x]_1 \subset J^{\perp}$ is equivalent to the vanishing of the scalar products
\[
\begin{split}
x_i \Lambda^{(1)}_{\HH} \cdot \widetilde{\Lambda}^{(1)}_{\KK}&{} = \left(\sum_{\hh} \Delta_{\HH\vert(\hh)}\, x_i x^{\alpha(\hh)}\right)\cdot\left(\sum_{\kk}\nabla_{\KK\vert(\kk)}\, x^{\alpha(\kk)}\right) = \\
&{}= \sum_{\begin{subarray}{c}\hh,\ \kk\\x_i x^{\alpha(\hh)} = x^{\alpha(\kk)}\end{subarray}} \Delta_{\HH\vert(\hh)} \nabla_{\KK\vert(\kk)},
\end{split}
\]
for each $i = 0,\ldots,n$ and for all $x_i \Lambda^{(1)}_{\HH} \in x_i\Gamma^{(1)}$ and $\widetilde{\Lambda}^{(1)}_{\KK} \in \widetilde{\Gamma}^{(1)}$.
Denoted by $\mathcal{I}_{\mathcal{H}}$ the ideal generated by the set of quadrics
\[
\left\{\sum_{\begin{subarray}{c}\hh,\ \kk\\x_i x^{\alpha(\hh)} = x^{\alpha(\kk)}\end{subarray}} \Delta_{\HH\vert(\hh)} \nabla_{\KK\vert(\kk)}\quad \Bigg\vert\quad \forall\ \HH,\KK,i \right\}
\]
and by the Pl\"ucker relations of both $\Grass{q(r)}{N(r)}{\K}$ and $\Grass{p(r+1)}{N(r+1)}{\K}$, the Hilbert scheme is defined as
\[
\Hilb{n}{p(t)} \simeq \Proj \big(\K[\ldots,\Delta_{\II},\ldots]\times\K[\ldots,\nabla_{\JJ},\ldots]\big)/\mathcal{I}_{\mathcal{H}}. \qedhere
\]
\end{proof}

\begin{example}\label{ex:MainExampleGotzmann}
Let us compute Gotzmann's equations of the Hilbert scheme \gls{Hilb2P2}. Since the Gotzmann number of the Hilbert polynomial $p(t)=2$ is 2 and $\binom{2+2}{2}=6$, $\binom{2+3}{2}=10$, we embed $\Hilb{2}{2}$ in
\[
\Grass{4}{6}{\K} \times \Grass{8}{10}{\K} \simeq \Grass{4}{6}{\K} \times \Grass{2}{10}{\K}
\]
The first Grassmannian is the same already introduced in Examples \ref{ex:MainExampleGrass} and \ref{ex:MainExampleGens}, but\hfill in\hfill this\hfill case\hfill we\hfill consider\hfill the\hfill monomial\hfill basis\hfill $\{x_0^2,x_0x_1,x_1^2,x_0x_2,x_1x_2,x_2^2\}$\hfill of\\ $\K[x_0,x_1,x_2]_2$. By Example \ref{ex:MainExampleGens}, we know that among the generators of $I \in \Grass{4}{6}{\K}$ there is
\[
\Lambda^{(1)}_{156} = \Delta_{1256}\, x_0x_1 + \Delta_{1356}\, x_1^2 + \Delta_{1456}\, x_0x_2. 
\]
Fixed\hfill the\hfill monomial\hfill basis\hfill $\{x_0^3,x_0^2x_1,x_0x_1^2,x_1^3,x_0^2x_2,x_0x_1x_2,x_1^2x_2,x_0x_2^2,x_1x_2^2,x_2^3\}$\hfill of\\ $\K[x_0,x_1,x_2]_3$ among the generators $\widetilde{\Gamma}^{(1)}$ of $J \in \Grass{2}{10}{\K}$ there is
\[
\begin{split}
\widetilde{\Lambda}^{(1)}_{7} = {}& -\nabla_{1,7}\, x_0^3 -\nabla_{2,7}\, x_0^2x_1 -\nabla_{3,7}\, x_0x_1^2 -\nabla_{4,7}\, x_1^3 -\nabla_{5,7}\, x_0^2x_2 +\\
& -\nabla_{6,7}\, x_0x_1x_2 +\nabla_{7,8}\, x_0x_2^2 +\nabla_{7,9}\, x_1x_2^2 +\nabla_{7,10}\, x_2^3.
\end{split}
\]

Thus among the equations described in Theorem \ref{th:GotzmannEquations}, we will find
\[
\begin{split}
x_0\Lambda^{(1)}_{156} \cdot \widetilde{\Lambda}^{(1)}_{7}&{} = \left(\Delta_{1256}\, x_0^2x_1 + \Delta_{1356}\, x_0x_1^2 + \Delta_{1456}\, x_0^2x_2\right) \cdot \widetilde{\Lambda}^{(1)}_{7} = \\
&{} = -\Delta_{1256}\nabla_{2,7} - \Delta_{1356}\nabla_{3,7} - \Delta_{1456}\nabla_{5,7},\\
x_1\Lambda^{(1)}_{156} \cdot \widetilde{\Lambda}^{(1)}_{7}&{} = \left(\Delta_{1256}\, x_0x_1^2 + \Delta_{1356}\, x_1^3 + \Delta_{1456}\, x_0x_1x_2\right) \cdot \widetilde{\Lambda}^{(1)}_{7} = \\
&{} = -\Delta_{1256}\nabla_{3,7} - \Delta_{1356}\nabla_{4,7} - \Delta_{1456}\nabla_{6,7},\\
x_2\Lambda^{(1)}_{156} \cdot \widetilde{\Lambda}^{(1)}_{7}&{} = \left(\Delta_{1256}\, x_0x_1x_2 + \Delta_{1356}\, x_1^2x_2 + \Delta_{1456}\, x_0x_2^2\right) \cdot \widetilde{\Lambda}^{(1)}_{7} = \\
&{} = -\Delta_{1256}\nabla_{6,7} + \Delta_{1456}\nabla_{7,8}.\\
\end{split}
\]
Since $\vert \Gamma^{(1)}\vert = \binom{6}{3} = 20$ and $\vert \widetilde{\Gamma}^{(1)}\vert= \binom{10}{1} = 10$, the condition $I \subset J^{\perp}$ is represented by $3\cdot20\cdot10 = 600$ bilinear equations (the same number of equation was determined by Haiman and Sturmfels \cite[Example 4.3]{HaimanSturmfels}) to whom we must add $60+840 = 900$ Pl\"ucker relations. See Example \ref{ex:MainExampleGotzmannCODE} in Appendix \ref{ch:Hilb2P2} for further details.
\end{example}

\subsection{Iarrobino-Kleiman equations}
As Example \ref{ex:MainExampleGotzmann} shows the transformation $\hilb{n}{p(t)}\rightarrow \grass{q(r)}{N(r)}\times\grass{q(r+1)}{N(r+1)}$ allows to describe the Hilbert scheme by quadric equation, i.e. of low degree, but requiring lots of Pl\"ucker coordinates and consequently lots of Pl\"ucker relations, even in one of the simplest cases. Therefore  we come back to the transformation of functors described in \eqref{eq:hilbToGrass}, in order to interpret the Hilbert functor as
\[
\hilb{n}{p(t)}(A) = \left\{ \begin{array}{c} I  \in \grass{q(r)}{N(r)}(A) \text{ s.t. }I \subset A[x]_r\\ \text{and }\rank A[x]_1 \cdot I = q(r+1) \end{array}\right\}
\]
and again by Macaulay's Estimate (Theorem \ref{th:MacaulayEstimate})
\[
\hilb{n}{p(t)}(A) = \left\{ \begin{array}{c} I  \in \grass{q(r)}{N(r)}(A) \text{ s.t. }I \subset A[x]_r\\ \text{and }\rank A[x]_1 \cdot I \leqslant q(r+1)\end{array}\right\}.
\]

\begin{theorem}[{\cite[Proposition C.30]{IarrobinoKanev}}]\label{th:IarrobinoKleimanEquations} The Hilbert scheme $\Hilb{n}{p(t)}$ can be defined as subscheme of the Grassmannian $\Grass{q}{N}{\K}$ by equations of degree $q(r+1)+1$.
\end{theorem}
\begin{proof}
Associated to any point of $\Grass{q}{N}{\K}$, we consider the ideal $I$ generated by $\Gamma^{(1)}$, so that the vector space $I_{r+1}$ is spanned by $\left\langle x_i\Gamma^{(1)}\ \vert\ i=0,\ldots,n\right\rangle$. Requiring that $\dim_{\K} I_{r+1} \leqslant q(r+1)$ is equivalent to ask that the $\big(q(r+1)+1\big)$-th exterior power of $I_{r+1}$ vanishes:
\begin{equation}\label{eq:conditionsIarrobinoKleiman}
\dim_{\K} I_{r+1} \leqslant q(r+1)\quad \Longleftrightarrow\quad \bigwedge^{q(r+1)+1} I_{r+1} = 0.
\end{equation}
The vector space $\wedge^{q(r+1)+1} I_{r+1}$ is spanned by the set
\[
\left\{ \bigwedge_{j=1}^{q(r+1)+1} x_{i_j} \Lambda^{(1)}_{\JJ_j}\quad \Bigg\vert\quad \forall\ i_j \in \{0,\ldots,n\},\ \forall\ x_{i_j}\Lambda^{(1)}_{\JJ_j} \in x_{i_j}\Gamma^{(1)} \right\}
\]
whose elements have coefficient in the Pl\"ucker coordinates of $\Grass{q}{N}{\K}$ of degree $q(r+1)+1$. Considering the ideal $\mathcal{I}_{\mathcal{H}}$ generated by all the coefficients of these generators of $\wedge^{q(r+1)+1} I_{r+1}$ and the ideal $\mathcal{Q}$ of the Pl\"ucker relations, we define the Hilbert scheme as
\[
\Hilb{n}{p(t)} \simeq \Proj \K[\ldots,\Delta_{\II},\ldots]/(\mathcal{Q},\mathcal{I}_{\mathcal{H}}). \qedhere
\]
\end{proof}

\begin{remark}
Iarrobino and Kleiman in \cite{IarrobinoKanev} proved this result in affine coordinates. In fact considering the standard affine covering of $\Proj \K[\ldots,\Delta_{\II},\ldots]$, on each affine chart we can consider as basis the ideal $I$ that one described in \eqref{eq:basisW}, thus working with an optimal set of generators, and then we can glue together these affine subschemes via transition maps. 
\end{remark}

\begin{remark}\label{rk:reduceIarrobinoKleimanEquations}
From a computational perspective, the set spanning $\wedge^{q(r+1)+1} I_{r+1}$ considered in the proof of Theorem \ref{th:IarrobinoKleimanEquations} contains many vanishing elements, indeed any set of $s > q$ polynomials, coming from $s$ polynomials of degree $r$ multiplied by the same variable, is surely dependent. A better set spanning $\wedge^{q(r+1)+1} I_{r+1}$ is
\begin{equation}\label{eq:optimalSet}
\left\{ \bigwedge_{j=0}^n \left(\bigwedge_{\ii = 1}^{s_j} x_j \Lambda^{(1)}_{\JJ_{\ii}} \right)\quad\Bigg\vert\quad \begin{array}{c}\forall\ 0 \leqslant s_j \leqslant q \text{ s.t. } \sum_{j=0}^n s_j = q(r+1)+1\\ \forall\ x_j\Lambda^{(1)}_{\JJ_{\ii}} \in x_j\Gamma^{(1)},\ j = 0,\ldots,n\end{array} \right\}
\end{equation}
\end{remark}

\begin{example}\label{ex:MainExampleIarrobinoKleiman}
Let us look at the Iarrobino-Kleiman equation of the Hilbert scheme \gls{Hilb2P2} already introduced in Example \ref{ex:MainExampleGotzmann}. In the proof of Theorem \ref{th:IarrobinoKleimanEquations}, we are asking that any set of $q(r+1)+1$ polynomials in $\{x_i\Gamma^{(1)}\ \vert\ i=0,\ldots,n\}$ is linearly dependent. Hence in this special case we have to impose the dependency of any set of 9 polynomials in $\{x_0\Gamma^{(1)},x_1\Gamma^{(1)},x_2\Gamma^{(1)}\}$. For instance, the polynomials represented in the following matrix
\[
\footnotesize
\begin{array}{r|cccccccccc|}
\multicolumn{1}{c}{} & x_0^3 & x_0^2x_1 & x_0x_1^2 & x_1^3 & x_0^2x_2 & x_0x_1x_2 & x_1^2x_2 & x_0x_2^2 & x_1x_2^2 & \multicolumn{1}{c}{x_2^3}\\
\cline{2-11} 
x_0 \Lambda^{(1)}_{126} & 0 & 0 & -\Delta_{1236} & 0 & -\Delta_{1246} & -\Delta_{1256} & 0 & 0 & 0 & 0  \\
x_0 \Lambda^{(1)}_{156} & 0 & \Delta_{1256} & \Delta_{1356} & 0 & \Delta_{1456} & 0 & 0 & 0 & 0 & 0\\
x_0 \Lambda^{(1)}_{234} & -\Delta_{1234} & 0 & 0 & 0 & 0 & \Delta_{2345} & 0 & \Delta_{2346} & 0 & 0\\
x_0 \Lambda^{(1)}_{356} & -\Delta_{1356} & -\Delta_{2356} & 0 & 0 & \Delta_{3456} & 0 & 0 & 0 & 0 & 0\\
x_1 \Lambda^{(1)}_{123} & 0 & 0 & 0 & 0 & 0 & \Delta_{1234} & \Delta_{1235} & 0 & \Delta_{1236} & 0\\
x_1 \Lambda^{(1)}_{245} & 0 & -\Delta_{1245} & 0 & \Delta_{2345} & 0 & 0 & 0 & 0 & \Delta_{2456} & 0\\
x_2 \Lambda^{(1)}_{146} & 0 & 0 & 0 & 0 & 0 & \Delta_{1246} & \Delta_{1346} & 0 & -\Delta_{1456} & 0\\
x_2 \Lambda^{(1)}_{234} & 0 & 0 & 0 & 0 & -\Delta_{1234}  & 0 & 0 & 0 & \Delta_{2345} & \Delta_{2346}\\
x_2 \Lambda^{(1)}_{456} & 0 & 0 & 0 & 0 & -\Delta_{1456} & -\Delta_{2456} & -\Delta_{3456} & 0 & 0 & 0\\
\cline{2-11}
\end{array}
\]
are linearly dependent if the rank of the matrix is not maximal, that is if the following polynomials of degree 9 in the Pl\"ucker coordinates, corresponding to the minors of dimension 9, vanish
\[
\footnotesize
\begin{split}
\bullet\quad &\Delta_{1235}\Delta_{1236}\Delta_{1246}\Delta_{1256}\Delta_{1356}\Delta_{1456}\Delta_{2345}^2\Delta_{2346}-\Delta_{1234}\Delta_{1236}\Delta_{1256}\Delta_{1346}\Delta_{1356}\Delta_{1456}\Delta_{2345}^2\Delta_{2346}+\\
&+\Delta_{1234}\Delta_{1236}^2\Delta_{1256}\Delta_{1346}\Delta_{1356}\Delta_{2345}\Delta_{2346}\Delta_{2456}+\Delta_{1234}\Delta_{1235}\Delta_{1236}\Delta_{1256}\Delta_{1356}\Delta_{1456}\Delta_{2345}\Delta_{2346}\Delta_{2456}+\\
&-\Delta_{1234}\Delta_{1236}^2\Delta_{1246}\Delta_{1256}\Delta_{1356}\Delta_{2345}\Delta_{2346}\Delta_{3456}-\Delta_{1234}^2\Delta_{1236}\Delta_{1256}\Delta_{1356}\Delta_{1456}\Delta_{2345}\Delta_{2346}\Delta_{3456},\\
\bullet\quad &\Delta_{1235}\Delta_{1236}\Delta_{1246}\Delta_{1256}\Delta_{1356}\Delta_{1456}\Delta_{2345}\Delta_{2346}^2-\Delta_{1234}\Delta_{1236}\Delta_{1256}\Delta_{1346}\Delta_{1356}\Delta_{1456}\Delta_{2345}\Delta_{2346}^2,\\
\bullet\quad &\Delta_{1236}^2\Delta_{1256}\Delta_{1346}\Delta_{1356}\Delta_{1456}\Delta_{2345}^2\Delta_{2346}+\Delta_{1235}\Delta_{1236}\Delta_{1256}\Delta_{1356}\Delta_{1456}^2\Delta_{2345}^2\Delta_{2346}+\\
&+\Delta_{1234}\Delta_{1236}\Delta_{1256}\Delta_{1346}\Delta_{1356}\Delta_{1456}\Delta_{2345}\Delta_{2346}\Delta_{2356}+\Delta_{1234}\Delta_{1235}\Delta_{1256}\Delta_{1356}\Delta_{1456}^2\Delta_{2345}\Delta_{2346}\Delta_{2356}+\\
&-\Delta_{1234}\Delta_{1236}\Delta_{1246}\Delta_{1346}\Delta_{1356}\Delta_{2345}\Delta_{2346}\Delta_{2356}\Delta_{2456}+\Delta_{1234}\Delta_{1236}^2\Delta_{1346}\Delta_{1456}\Delta_{2345}\Delta_{2346}\Delta_{2356}\Delta_{2456}+\\
&-\Delta_{1234}\Delta_{1235}\Delta_{1246}\Delta_{1356}\Delta_{1456}\Delta_{2345}\Delta_{2346}\Delta_{2356}\Delta_{2456}+\Delta_{1234}\Delta_{1235}\Delta_{1236}\Delta_{1456}^2\Delta_{2345}\Delta_{2346}\Delta_{2356}\Delta_{2456}+\\
&+\Delta_{1234}\Delta_{1236}\Delta_{1246}^2\Delta_{1356}\Delta_{2345}\Delta_{2346}\Delta_{2356}\Delta_{3456}-\Delta_{1234}\Delta_{1236}^2\Delta_{1246}\Delta_{1456}\Delta_{2345}\Delta_{2346}\Delta_{2356}\Delta_{3456}+\\
&+\Delta_{1234}^2\Delta_{1246}\Delta_{1356}\Delta_{1456}\Delta_{2345}\Delta_{2346}\Delta_{2356}\Delta_{3456}-\Delta_{1234}^2\Delta_{1236}\Delta_{1456}^2\Delta_{2345}\Delta_{2346}\Delta_{2356}\Delta_{3456}+\\
&+\Delta_{1234}\Delta_{1236}^2\Delta_{1256}\Delta_{1346}\Delta_{2345}\Delta_{2346}\Delta_{2456}\Delta_{3456}+\Delta_{1234}\Delta_{1235}\Delta_{1236}\Delta_{1256}\Delta_{1456}\Delta_{2345}\Delta_{2346}\Delta_{2456}\Delta_{3456}+\\
&-\Delta_{1234}\Delta_{1236}^2\Delta_{1246}\Delta_{1256}\Delta_{2345}\Delta_{2346}\Delta_{3456}^2-\Delta_{1234}^2\Delta_{1236}\Delta_{1256}\Delta_{1456}\Delta_{2345}\Delta_{2346}\Delta_{3456}^2,\\
\bullet\quad & -\Delta_{1236}^2\Delta_{1246}\Delta_{1256}\Delta_{1356}\Delta_{1456}\Delta_{2345}\Delta_{2346}^2-\Delta_{1234}\Delta_{1236}\Delta_{1256}\Delta_{1356}\Delta_{1456}^2\Delta_{2345}\Delta_{2346}^2,\\
\end{split}
\]
\[
\footnotesize
\begin{split}
\bullet\quad &-\Delta_{1236}^2\Delta_{1256}\Delta_{1346}\Delta_{1356}\Delta_{1456}\Delta_{2345}\Delta_{2346}^2-\Delta_{1235}\Delta_{1236}\Delta_{1256}\Delta_{1356}\Delta_{1456}^2\Delta_{2345}\Delta_{2346}^2,\\
\bullet\quad &-\Delta_{1236}^2\Delta_{1256}\Delta_{1346}\Delta_{1356}\Delta_{2345}\Delta_{2346}^2\Delta_{2456}-\Delta_{1235}\Delta_{1236}\Delta_{1256}\Delta_{1356}\Delta_{1456}\Delta_{2345}\Delta_{2346}^2\Delta_{2456}+\\
&+\Delta_{1236}^2\Delta_{1246}\Delta_{1256}\Delta_{1356}\Delta_{2345}\Delta_{2346}^2\Delta_{3456}+\Delta_{1234}\Delta_{1236}\Delta_{1256}\Delta_{1356}\Delta_{1456}\Delta_{2345}\Delta_{2346}^2\Delta_{3456},\\
\bullet\quad & -\Delta_{1236}\Delta_{1256}^2\Delta_{1346}\Delta_{1356}\Delta_{1456}\Delta_{2345}\Delta_{2346}^2-\Delta_{1235}\Delta_{1256}^2\Delta_{1356}\Delta_{1456}^2\Delta_{2345}\Delta_{2346}^2+\\
&+\Delta_{1236}\Delta_{1246}\Delta_{1256}\Delta_{1346}\Delta_{1356}\Delta_{2345}\Delta_{2346}^2\Delta_{2456}+\Delta_{1235}\Delta_{1246}\Delta_{1256}\Delta_{1356}\Delta_{1456}\Delta_{2345}\Delta_{2346}^2\Delta_{2456}+\\
&-\Delta_{1236}\Delta_{1246}^2\Delta_{1256}\Delta_{1356}\Delta_{2345}\Delta_{2346}^2\Delta_{3456}-\Delta_{1234}\Delta_{1246}\Delta_{1256}\Delta_{1356}\Delta_{1456}\Delta_{2345}\Delta_{2346}^2\Delta_{3456},\\
\bullet\quad & -\Delta_{1236}\Delta_{1245}\Delta_{1256}\Delta_{1346}\Delta_{1356}^2\Delta_{1456}\Delta_{2346}^2-\Delta_{1235}\Delta_{1245}\Delta_{1256}\Delta_{1356}^2\Delta_{1456}^2\Delta_{2346}^2+\\
&+\Delta_{1236}\Delta_{1245}\Delta_{1246}\Delta_{1346}\Delta_{1356}^2\Delta_{2346}^2\Delta_{2456}+\Delta_{1235}\Delta_{1236}\Delta_{1246}\Delta_{1256}\Delta_{1356}\Delta_{1456}\Delta_{2346}^2\Delta_{2456}+\\
&-\Delta_{1236}^2\Delta_{1245}\Delta_{1346}\Delta_{1356}\Delta_{1456}\Delta_{2346}^2\Delta_{2456}-\Delta_{1234}\Delta_{1236}\Delta_{1256}\Delta_{1346}\Delta_{1356}\Delta_{1456}\Delta_{2346}^2\Delta_{2456}+\\
&+\Delta_{1235}\Delta_{1245}\Delta_{1246}\Delta_{1356}^2\Delta_{1456}\Delta_{2346}^2\Delta_{2456}-\Delta_{1235}\Delta_{1236}\Delta_{1245}\Delta_{1356}\Delta_{1456}^2\Delta_{2346}^2\Delta_{2456}+\\
&-\Delta_{1236}\Delta_{1245}\Delta_{1246}^2\Delta_{1356}^2\Delta_{2346}^2\Delta_{3456}+\Delta_{1236}^2\Delta_{1245}\Delta_{1246}\Delta_{1356}\Delta_{1456}\Delta_{2346}^2\Delta_{3456}+\\
&-\Delta_{1234}\Delta_{1245}\Delta_{1246}\Delta_{1356}^2\Delta_{1456}\Delta_{2346}^2\Delta_{3456}+\Delta_{1234}\Delta_{1236}\Delta_{1245}\Delta_{1356}\Delta_{1456}^2\Delta_{2346}^2\Delta_{3456},\\
\bullet\quad & -\Delta_{1236}\Delta_{1256}\Delta_{1346}\Delta_{1356}^2\Delta_{1456}\Delta_{2345}\Delta_{2346}^2-\Delta_{1235}\Delta_{1256}\Delta_{1356}^2\Delta_{1456}^2\Delta_{2345}\Delta_{2346}^2+\\
&+\Delta_{1236}\Delta_{1246}\Delta_{1346}\Delta_{1356}^2\Delta_{2345}\Delta_{2346}^2\Delta_{2456}-\Delta_{1236}^2\Delta_{1346}\Delta_{1356}\Delta_{1456}\Delta_{2345}\Delta_{2346}^2\Delta_{2456}+\\
&+\Delta_{1235}\Delta_{1246}\Delta_{1356}^2\Delta_{1456}\Delta_{2345}\Delta_{2346}^2\Delta_{2456}-\Delta_{1235}\Delta_{1236}\Delta_{1356}\Delta_{1456}^2\Delta_{2345}\Delta_{2346}^2\Delta_{2456}+\\
&-\Delta_{1236}\Delta_{1246}^2\Delta_{1356}^2\Delta_{2345}\Delta_{2346}^2\Delta_{3456}+\Delta_{1236}^2\Delta_{1246}\Delta_{1356}\Delta_{1456}\Delta_{2345}\Delta_{2346}^2\Delta_{3456}+\\
&-\Delta_{1234}\Delta_{1246}\Delta_{1356}^2\Delta_{1456}\Delta_{2345}\Delta_{2346}^2\Delta_{3456}+\Delta_{1234}\Delta_{1236}\Delta_{1356}\Delta_{1456}^2\Delta_{2345}\Delta_{2346}^2\Delta_{3456},\\
\bullet\quad & -\Delta_{1236}\Delta_{1256}\Delta_{1346}\Delta_{1356}\Delta_{1456}\Delta_{2345}\Delta_{2346}^2\Delta_{2356}-\Delta_{1235}\Delta_{1256}\Delta_{1356}\Delta_{1456}^2\Delta_{2345}\Delta_{2346}^2\Delta_{2356}+\\
&+\Delta_{1236}\Delta_{1246}\Delta_{1346}\Delta_{1356}\Delta_{2345}\Delta_{2346}^2\Delta_{2356}\Delta_{2456}-\Delta_{1236}^2\Delta_{1346}\Delta_{1456}\Delta_{2345}\Delta_{2346}^2\Delta_{2356}\Delta_{2456}+\\
&+\Delta_{1235}\Delta_{1246}\Delta_{1356}\Delta_{1456}\Delta_{2345}\Delta_{2346}^2\Delta_{2356}\Delta_{2456}-\Delta_{1235}\Delta_{1236}\Delta_{1456}^2\Delta_{2345}\Delta_{2346}^2\Delta_{2356}\Delta_{2456}+\\
&-\Delta_{1236}\Delta_{1246}^2\Delta_{1356}\Delta_{2345}\Delta_{2346}^2\Delta_{2356}\Delta_{3456}+\Delta_{1236}^2\Delta_{1246}\Delta_{1456}\Delta_{2345}\Delta_{2346}^2\Delta_{2356}\Delta_{3456}+\\
&-\Delta_{1234}\Delta_{1246}\Delta_{1356}\Delta_{1456}\Delta_{2345}\Delta_{2346}^2\Delta_{2356}\Delta_{3456}+\Delta_{1234}\Delta_{1236}\Delta_{1456}^2\Delta_{2345}\Delta_{2346}^2\Delta_{2356}\Delta_{3456}+\\
&-\Delta_{1236}^2\Delta_{1256}\Delta_{1346}\Delta_{2345}\Delta_{2346}^2\Delta_{2456}\Delta_{3456}-\Delta_{1235}\Delta_{1236}\Delta_{1256}\Delta_{1456}\Delta_{2345}\Delta_{2346}^2\Delta_{2456}\Delta_{3456}+\\
&+\Delta_{1236}^2\Delta_{1246}\Delta_{1256}\Delta_{2345}\Delta_{2346}^2\Delta_{3456}^2+\Delta_{1234}\Delta_{1236}\Delta_{1256}\Delta_{1456}\Delta_{2345}\Delta_{2346}^2\Delta_{3456}^2.
\end{split}
\]

If we consider any possible subset of 9 elements of $\{x_0\Gamma^{(1)},x_1\Gamma^{(1)},$ $x_2\Gamma^{(1)}\}$, we would need to examine $\binom{60}{9}$ exterior products, each of them containing at most 10 terms ($\dim_{\K} \wedge^9 \K[x_0,x_1,x_2]_3 = 10$) so that an overestimate of the number of equations would be $10\cdot\binom{60}{9} = 147831426600$. Following Remark \ref{rk:reduceIarrobinoKleimanEquations}, it suffices to consider the subsets of 9 polynomials subdivided according to the multiplication variable in 3 subsets of cardinality 4,4,1 or 4,3,2 or 3,3,3.
Of the type (4,4,1), there are $3\cdot\binom{20}{4}^2\cdot\binom{20}{1} = 1408441500$ possibilities, $6\cdot\binom{20}{4}\cdot\binom{20}{3}\cdot\binom{20}{2} = 6296562000$ of the type (4,3,2) and $\binom{20}{3}^3 = 1481544000$ corresponding to (3,3,3), for a total of 9186547500: still a huge number just under the half of $\binom{60}{9}$. It is clear that the set of equations defining the Hilbert scheme provided by Theorem \ref{th:IarrobinoKleimanEquations} can not be use to project an effective algorithm (see Section \ref{sec:HilbEquationsCode} of Appendix \ref{ch:Hilb2P2}).
\end{example}

\subsection{Bayer-Haiman-Sturmfels equations}

Even in the very simple case of $\Hilb{2}{2}$, the equations determined by Iarrobino and Kleiman have a large degree and they are too many to think about using them for a methodical study of Hilbert schemes through computational software. Then one of the first goal is to lower the degree of the equations. The idea introduced by Bayer in his thesis \cite{BayerThesis} is to put together polynomials of degree $r+1$ coming from polynomials of degree $r$ multiplied for the same variable, that is to associate to any subproduct in \eqref{eq:optimalSet} a generator of an exterior power of $I_{r+1}$ having as coefficients linear polynomials in the Pl\"ucker coordinates:
\[
\bigwedge_{\ii = 1}^{s_j} x_j \Lambda^{(1)}_{\JJ_{\ii}}\qquad \stackrel{?}{\longleftrightarrow}\qquad x_j\Lambda^{(s_j)}_{\JJ}.
\]

\begin{theorem}[{\cite[Theorem 4.4]{HaimanSturmfels}}]\label{th:BayerHaimanSturmfelsEquations} The Hilbert scheme $\Hilb{n}{p(t)}$ can be defined as subscheme of the Grassmannian $\Grass{q}{N}{\K}$ by equations of degree equal to or less than $n+1$.
\end{theorem}
\begin{proof}
We\hfill consider\hfill again\hfill the\hfill equivalent\hfill condition\hfill showed\hfill in\hfill \eqref{eq:conditionsIarrobinoKleiman}\hfill to\hfill impose\\ $\dim_{\K} I_{r+1} \leqslant q(r+1)$, but we construct a different set of generators for $\wedge^{q(r+1)+1} I_{r+1}$ taking advantage of the elements of $x_i\Gamma^{(s)},\ \forall\ 1 \leqslant s \leqslant q$, indeed $\wedge^{q(r+1)+1} I_{r+1}$ is spanned by the set
\[
\left\{ \bigwedge_{j=0}^n x_{j} \Lambda^{(s_j)}_{\JJ_j}\quad\Bigg\vert\quad \begin{array}{l}\forall\ s_j \leqslant q \text{ s.t. } \sum_{j=0}^n s_j = q(r+1)+1\\ \forall\ x_0\Lambda^{(s_0)}_{\JJ_0} \in x_0\Gamma^{(s_0)},\ \ldots,\ \forall\ x_n\Lambda^{(s_n)}_{\JJ_n} \in x_n\Gamma^{(s_n)}\end{array} \right\}.
\]
Denoted by $\mathcal{I}_{\mathcal{H}}$ the ideal generated by all the coefficients of the elements spanning $\wedge^{q(r+1)+1} I_{r+1}$, that are polynomials of degree $n+1$ in the Pl\"ucker coordinates, and by $\mathcal{Q}$ the ideal of the Pl\"ucker relations, the Hilbert scheme $\Hilb{n}{p(t)}$ can be defined as
\[
\Hilb{n}{p(t)} \simeq \Proj \K[\ldots,\Delta_{\II},\ldots]/(\mathcal{Q},\mathcal{I}_{\mathcal{H}}). \qedhere
\]
\end{proof}

\begin{example}\label{ex:MainExampleBayerHaimanSturmfels}
Let us consider again \gls{Hilb2P2} as in Examples \ref{ex:MainExampleGotzmann} and \ref{ex:MainExampleIarrobinoKleiman}. Applying Theorem \ref{th:BayerHaimanSturmfelsEquations}, we have to impose the vanishing of the exterior products
\[
x_0 \Lambda^{(s_0)}_{\JJ_0} \wedge x_1 \Lambda^{(s_1)}_{\JJ_1} \wedge x_2 \Lambda^{(s_2)}_{\JJ_2}, \qquad\forall\  \Lambda^{(s_0)}_{\JJ_0} \in \Gamma^{(s_0)},\ \Lambda^{(s_1)}_{\JJ_1} \in \Gamma^{(s_1)},\ \Lambda^{(s_2)}_{\JJ_2} \in \Gamma^{(s_2)},\ 
\]
where $(s_0,s_1,s_2)$ is a triple chosen among the set $\{(4,4,1),(4,3,2),(4,2,3),(4,1,4),$ $(3,4,2),(3,3,3),(3,2,4),(2,4,3),(2,3,4),(1,4,4)\}$. For instance let us compute explicitly the product $x_0 \Lambda^{(4)}_{\emptyset}\wedge  x_1 \Lambda^{(2)}_{25} \wedge x_2 \Lambda^{(3)}_{1}$:
\[
\footnotesize
\begin{split}
x_0 \Lambda^{(4)}_{\emptyset} &{} = \Delta_{1234}\, x_0^3 \wedge x_0^2 x_1 \wedge x_0 x_1^2 \wedge x_0^2x_2 + \Delta_{1235}\, x_0^3 \wedge x_0^2 x_1 \wedge x_0 x_1^2 \wedge x_0x_1x_2 +\\
&{} + \Delta_{1236}\, x_0^3 \wedge x_0^2 x_1 \wedge x_0 x_1^2 \wedge x_0x_2^2 + \Delta_{1245}\, x_0^3 \wedge x_0^2 x_1 \wedge x_0^2x_2 \wedge x_0x_1x_2 +\\
&{} + \Delta_{1246}\, x_0^3 \wedge x_0^2 x_1 \wedge x_0^2x_2 \wedge x_0x_2^2 + \Delta_{1256}\, x_0^3 \wedge x_0^2 x_1 \wedge x_0x_1x_2 \wedge x_0x_2^2 + \\
&{} + \Delta_{1345}\, x_0^3 \wedge x_0 x_1^2 \wedge x_0^2x_2 \wedge x_0x_1x_2 + \Delta_{1346}\, x_0^3 \wedge x_0 x_1^2 \wedge x_0^2x_2 \wedge x_0x_2^2 + \\
&{} + \Delta_{1356}\, x_0^3 \wedge x_0 x_1^2 \wedge x_0x_1x_2 \wedge x_0x_2^2 + \Delta_{1456}\, x_0^3 \wedge x_0^2 x_2 \wedge x_0x_1x_2 \wedge x_0x_2^2 +\\
&{} + \Delta_{2345}\, x_0^2 x_1 \wedge x_0 x_1^2 \wedge x_0^2x_2 \wedge x_0x_1x_2 + \Delta_{2346}\, x_0^2x_1 \wedge x_0 x_1^2 \wedge x_0^2x_2 \wedge x_0x_2^2 + \\
&{} + \Delta_{2356}\, x_0^2x_1 \wedge x_0 x_1^2 \wedge x_0 x_1 x_2 \wedge x_0x_2^2 + \Delta_{2456}\, x_0^2x_1 \wedge x_0^2 x_2 \wedge x_0 x_1 x_2 \wedge x_0x_2^2 +\\
&{} + \Delta_{3456}\, x_0x_1^2 \wedge x_0^2 x_2 \wedge x_0 x_1 x_2 \wedge x_0x_2^2,\\
x_1 \Lambda^{(2)}_{25} &{}= - \Delta_{1235}\, x_0^2x_1 \wedge x_1^3 - \Delta_{1245}\, x_0^2x_1 \wedge x_0x_1x_2 + \Delta_{1256}\, x_0^2x_1 \wedge x_1x_2^2 + \Delta_{2345}\, x_1^3 \wedge x_0x_1x_2 +\\
&{} - \Delta_{2356}\, x_1^3 \wedge x_1x_2^2 - \Delta_{2456}\, x_0x_1x_2 \wedge x_1 x_2^2,\\
x_2 \Lambda^{(3)}_{1}&{} = \Delta_{1234}\, x_0x_1x_2 \wedge x_1^2 x_2 \wedge x_0x_2^2 + \Delta_{1235}\, x_0x_1x_2 \wedge x_1^2 x_2 \wedge x_1x_2^2 + \Delta_{1236}\, x_0x_1x_2 \wedge x_1^2 x_2 \wedge x_2^3 + \\
&{}  + \Delta_{1245}\, x_0x_1x_2 \wedge x_0x_2^2 \wedge x_1x_2^2 + \Delta_{1246}\, x_0x_1x_2 \wedge x_0x_2^2 \wedge x_2^3 + \Delta_{1256}\, x_0x_1x_2 \wedge x_1x_2^2 \wedge x_2^3 +\\
&{} + \Delta_{1345}\, x_1^2x_2 \wedge x_0x_2^2 \wedge x_1x_2^2 + \Delta_{1346}\, x_1^2x_2 \wedge x_0x_2^2 \wedge x_2^3 + \Delta_{1356}\, x_1^2x_2 \wedge x_1x_2^2 \wedge x_2^3 +\\
&{} + \Delta_{1456}\, x_0x_2^2 \wedge x_1x_2^2 \wedge x_2^3.
\end{split}
\]
We obtain the following 10 equations:
\[
\footnotesize
\begin{split}
\bullet\quad& \Delta_{1235} \Delta_{1345}^2-\Delta_{1234} \Delta_{1345} \Delta_{2345}+\Delta_{1235}^2 \Delta_{1346}-\Delta_{1234}^2 \Delta_{2356},\\
\bullet\quad& \Delta_{1235} \Delta_{1345} \Delta_{1346}-\Delta_{1234} \Delta_{2345} \Delta_{1346}+\Delta_{1235} \Delta_{1236} \Delta_{1346},\\
\bullet\quad& \Delta_{1235} \Delta_{1345} \Delta_{1356}-\Delta_{1234} \Delta_{2345} \Delta_{1356}+\Delta_{1234} \Delta_{1236} \Delta_{2356},\\
\bullet\quad& -\Delta_{1235} \Delta_{1346} \Delta_{1256}+\Delta_{1234} \Delta_{1246} \Delta_{2356}+\Delta_{1235} \Delta_{1345} \Delta_{1456}-\Delta_{1234} \Delta_{2345} \Delta_{1456},\\
\bullet\quad& -\Delta_{1235} \Delta_{1346} \Delta_{1356}+\Delta_{1234} \Delta_{1346} \Delta_{2356},\\
\bullet\quad& -\Delta_{2345} \Delta_{1236} \Delta_{1356}-\Delta_{1235} \Delta_{1356}^2+\Delta_{1236}^2 \Delta_{2356}+\Delta_{1235} \Delta_{1346} \Delta_{2356},\\
\bullet\quad& -\Delta_{1345} \Delta_{1346} \Delta_{1256}-\Delta_{1236} \Delta_{1346} \Delta_{1256}+\Delta_{1245} \Delta_{1346} \Delta_{1356}-\Delta_{1234} \Delta_{1346} \Delta_{2456},\\
\bullet\quad& \Delta_{2345} \Delta_{1246} \Delta_{1356}-\Delta_{1236} \Delta_{1246} \Delta_{2356}-\Delta_{1245} \Delta_{1346} \Delta_{2356}+\Delta_{1235} \Delta_{1356} \Delta_{1456},\\
\bullet\quad& \Delta_{2345} \Delta_{1346} \Delta_{1356}-\Delta_{1345} \Delta_{1346} \Delta_{2356}-\Delta_{1236} \Delta_{1346} \Delta_{2356},\\
\bullet\quad& \Delta_{2345} \Delta_{2346} \Delta_{1356}-\Delta_{2345} \Delta_{1346} \Delta_{2356}-\Delta_{1236} \Delta_{2346} \Delta_{2356}-\Delta_{1235} \Delta_{1356} \Delta_{3456}.
\end{split}
\]

Globally for the tern $(4,4,1)$ there are $1\cdot1\cdot 20$ possibilities, there are $1\cdot 6 \cdot 15 = 90$ to $(4,3,2)$ and $3\cdot 6 = 18$ to $(3,3,3)$, so we have to examine $3 \cdot 20 + 6 \cdot 90 + 18 = 618$ products. Since $\dim_{\K} \wedge^9\, \K[x_0,x_1,x_2]_3 = 10$, for each product we could have at most 10 coefficients and an upper bound of the number of equations for $\mathcal{I}_{\mathcal{H}}$ is 6180. 
Comparing this set of equations with that discussed in Example \ref{ex:MainExampleIarrobinoKleiman}, we discover that the number of exterior products to be examined is reduced by about a factor of $10^7$.
\end{example}

Let us examine in the general case the difference between the number of exterior products to be considered according to Theorem \ref{th:IarrobinoKleimanEquations} (Remark \ref{rk:reduceIarrobinoKleimanEquations}) and Theorem \ref{th:BayerHaimanSturmfelsEquations}. We saw that in both cases every exterior product can be subdivided in $n+1$ parts depending on the variable used to move by multiplication from a element of degree $r$ to an element of degree $r+1$, so let us consider any sequence of integers $(s_0,\ldots,s_n)$ such that $0 \leqslant s_i \leqslant q,\ \forall\ i =0,\ldots,n$ and $\sum_{i=0}^n s_i = q(r+1)+1$.

To compute Iarrobino-Kleiman equations, we have to choose in every possible way $s_i$ polynomials among $x_i\Gamma^{(1)}$ for each $i$, hence the possibilities are
\begin{equation*}
\prod_{i=0}^n \binom{\vert x_i \Gamma^{(1)} \vert}{s_i} = \prod_{i=0}^n \binom{\binom{N}{q-1}}{s_i}.
\end{equation*}
Instead for Bayer-Haiman-Sturmfels equations, for every $i$ we need to pick a single element among $x_i\Gamma^{(s_i)}$, so that the possibilities are
\begin{equation*}
\prod_{i=0}^n \left\vert x_i\Gamma^{(s_i)}\right\vert = \prod_{i=0}^n \binom{N}{q-s_i}.
\end{equation*}

It would be interesting to determine a good underestimation of the ratio between these two numbers, independent of $s_i$, i.e. a formula
\begin{equation}
\dfrac{\prod_{i=0}^n \binom{\binom{N}{q-1}}{s_i}}{\prod_{i=0}^n \binom{N}{q-s_i}} = \prod_{i=0}^n \dfrac{ \binom{\binom{N}{q-1}}{s_i}}{\binom{N}{q-s_i}} \geqslant F(q,N).
\end{equation}

\chapter{Borel-fixed ideals}\label{ch:BorelFixedIdeals}

In this chapter we introduce the most important objects used in this thesis, i.e. the Borel-fixed ideals.

\section{Definition}
Let us consider the usual action of the linear group \gls{GL} of the square matrix of dimension $n+1$ on the polynomial ring $\K[x_0,\ldots,n_n]$, i.e. for any invertible matrix $g=(g_{ij}) \in \GL_{\K}(n+1)$, the variable $x_i$ is mapped to the linear form $g\centerdot x_i = \sum_{j} g_{ij} x_j$, so that any polynomial $P(x) \in \K[x]$ is mapped to
\[
g\centerdot P(x) = g\centerdot P(x_0,\ldots,x_n) = P(g\centerdot x_0,\ldots,g \centerdot x_n) = P(g\centerdot x).
\]
Hence given an ideal $I\subset \K[x]$, it is well define the ideal 
\[
g\centerdot I = \langle g\centerdot P\ \vert\ \forall\ P \in I \rangle.
\]

For any ideal $I \subset \K[x]$ and for any term ordering $\sigma$, we can define the following equivalence relation on $\GL_{\K}(n+1)$:
\begin{equation}
g \sim g' \qquad \Longleftrightarrow \qquad \glslink{initIdeal}{\IN_{\sigma}(g\centerdot I)} = \IN_{\sigma}(g'\centerdot I)
\end{equation}
Viewing each matrix $g \in \GL_{\K}(n+1)$ as a element of the affine space $\AA^{(n+1)^2} = \Spec \K[\ldots, g_{ij},\ldots]$, the equivalence classes correspond to a stratification of $\AA^{(n+1)^2}$.

\begin{lemma}[{\cite[Lemma 2.6]{MillerSturmfels}}]\label{lem:finiteClasses}
For a fixed ideal $I$ and a term order $\sigma$, the number of equivalence classes in $\GL_{\K}(n+1)$ is finite and one of these classes is a nonempty Zariski open subset $U$ of $\GL_{\K}(n+1)$
\end{lemma}

\begin{definition}\label{def:gin}\index{generic initial ideal}
Let $I$ be an ideal of $\K[x]$ and let $\sigma$ be a fixed term ordering. The initial ideal $\IN_{\sigma}(g\centerdot I)$ obtained considering a matrix $g$ in the open subset $U \subset \GL_{\K}(n+1)$ corresponding to an equivalence class as in Lemma \ref{lem:finiteClasses} is called \emph{generic initial ideal} of $I$ w.r.t. $\sigma$ and it is denote by \gls{gin}.
\end{definition}

\bigskip

To better understand the interaction between the change of coordinates and the computation of initial ideals, let us now recall the \textsf{LU} decomposition of a matrix.
\begin{theorem}
Given any square matrix $g \in \GL_{\K}(n+1)$, if $g$ can be turned into an upper triangular matrix by means of Gaussian elimination without swapping rows, the $g$ can be decomposed as the product $l\cdot u$, where $l$ is a lower triangular matrix and $u$ is an upper triangular matrix with all entries on the diagonal equal to 1:
\begin{equation*}
\small
\left(\begin{array}{c c c c c}
g_{00} && \ldots && g_{0n}\\
&&&&\\
\vdots && g_{ij} && \vdots\\
&&&&\\
g_{n0} && \ldots && g_{nn}\\
\end{array}\right) =
\left(\begin{array}{c c c c c}
l_{00} && 0 && 0\\
& \ddots && \ddots&\\
l_{i0} && l_{ii} && 0\\
& \ddots && \ddots &\\
l_{n0} && l_{ni} && l_{nn}\\
\end{array}\right)
\left(\begin{array}{c c c c c}
1 && u_{0j} && u_{0n}\\
& \ddots && \ddots&\\
0 && 1 && u_{jn}\\
& \ddots && \ddots &\\
0 && 0 && 1\\
\end{array}\right)
\normalsize 
\end{equation*}
\end{theorem}

We remark that always considering matrices as points of $\AA^{(n+1)^2}$, the matrices having a \textsf{LU} decomposition correspond to an open subset, indeed any matrix that need row interchanges to be made upper triangular satisfies linear realations among the variables.

\begin{lemma}
Let $\sigma$ be any term ordering. For each homogeneous polynomial $P \in \K[x]_r$
\begin{equation}
\IN_\sigma(P) = \IN_{\sigma}(l \centerdot P)
\end{equation}
for each lower triangular matrix $l$.
\end{lemma}
\begin{proof}
It suffices to prove the statement for any monomial $x^{\alpha} = x_n^{\alpha_n}\cdots x_0^{\alpha_0},\ \vert\alpha\vert = r$.
\[
\begin{split}
l\centerdot x^{\alpha} &{}= (l \centerdot x_n)^{\alpha_n} \cdots (l\centerdot x_i)^{\alpha_i} \cdots (l\centerdot x_0^{\alpha_0})= \\
&{} = (l_{nn}x_n + \ldots + l_{n0}x_0)^{\alpha_n}\cdots(l_{ii}x_i + \ldots + u_{j0}x_0)^{\alpha_i}\cdots (l_{00}x_0)^{\alpha_0} = \\
&{} = (l_{nn}^{\alpha_n}x_n^{\alpha_n} + \text{ lower terms})\cdots(l_{ii}^{\alpha_i}x_i^{\alpha_i} + \text{ lower terms}) \cdots (l_{00}^{\alpha_0}x_0^{\alpha_0}) =\\
&{} = \left(\prod_{i=0}^n l_{ii}^{\alpha_i}\right) x^{\alpha} + \text{ lower terms}. \qedhere
\end{split}
\]
\end{proof}

From this lemma, it is clear that to understand which ideals have a nice behaviour under change of coordinates and initial ideal computation the key point is analyzing the action of the Borel subgroup \gls{BGL} of \glshyperlink{GL} of the upper triangular matrices.\index{Borel group}

\begin{definition}\label{def:Borel-fixed}\index{Borel-fixed ideal}
An ideal $I \subset \K[x]$ is called \emph{Borel-fixed} if it is fixed by the action of the Borel subgroup, i.e.
\[
g\centerdot I = I,\qquad\forall\ g \in \BGL_{\K}(n+1).
\]
\end{definition}

First of all, a Borel-fixed ideal has to be a monomial ideal because the Borel subgroup contains the algebraic torus group \gls{TGL} of the diagonal matrices that fixes all (and only) the monomials ideals. In the case of our interest, i.e. with a ground field $\K$ of characteristic 0, there is the following characterization.

\begin{proposition}[{\cite[Proposition 1.25]{GreenGIN}, \cite[Proposition 2.3]{MillerSturmfels}}]\label{prop:BorelIdealProperties}\index{Borel-fixed ideal!combinatorial characterization of a}
Let $I \subset \K[x]$ be a monomial ideal. The following statements are equivalent
\begin{enumerate}[(1)]
\item\label{it:BorelIdealProperties_i} $I$ is Borel-fixed;
\item\label{it:BorelIdealProperties_ii} if $x^\alpha \in I$, then $\frac{x_j}{x_i}x^\alpha \in I,\ \forall\ x_i \mid x^\alpha, j > i$;
\end{enumerate}
\end{proposition}
\begin{proof}
\emph{(\ref{it:BorelIdealProperties_i})$\Rightarrow$(\ref{it:BorelIdealProperties_ii}).} Let us consider the action on $x^\alpha$ of the matrix $g \in \BGL_{\K}(n+1)$ sending the variable $x_i$ to $x_j + x_i$ (i.e. $j > i$) and leaving fixed the others:
\[
\begin{split}
g\centerdot x^\alpha&{} = g\centerdot \prod_{k=0}^n x_k^{\alpha_k} = x_n^{\alpha_n}\cdots x_{i+1}^{\alpha_{i+1}} \cdot (x_j+x_i)^{\alpha_i} \cdot x_{i-1}^{\alpha_{i-1}} \cdots x_0^{\alpha_0} = \\
&{}= \left(\prod_{k\neq i} x_k^{\alpha_k} \right) \cdot \sum_{h=0}^{\alpha_i} \binom{\alpha_i}{h} x_j^{h} x_i^{\alpha_i-h}.
\end{split}
\]
Since $g\centerdot I = I$ and $I$ is a monomial ideal, each monomial appearing in $g\centerdot x^\alpha$ belongs to $I$ and $\left(\prod_{k\neq i} x_k^{\alpha_k} \right) x_j x_i^{\alpha_i-1} = \frac{x_j}{x_i}x^\alpha$.

\emph{(\ref{it:BorelIdealProperties_ii})$\Rightarrow$(\ref{it:BorelIdealProperties_i}).} Let $x^\alpha \in I$. For all $g \in \BGL_{\K}(n+1)$, each monomial in $g\centerdot x^\alpha$ can be obtained from $x^\alpha$ through a sequence of multiplication by $\frac{x_j}{x_i},\ j > i$. By the hypothesis \emph{(\ref{it:BorelIdealProperties_ii})} each monomial belongs to $I$, so that $g\centerdot x^\alpha \in I$.
\end{proof}

Now we can state the theorem providing the link between Borel-fixed ideals and Hilbert schemes due to Galligo \cite{Galligo} in characteristic 0 and generalized by Bayer and Stillman \cite{BayerStillmanGIN} in any characteristic.

\begin{theorem}\index{generic initial ideal}
The generic initial ideal $\GIN_{\sigma}(I)$ is Borel-fixed.
\end{theorem}
\begin{proof}
See \cite[Theorem 15.20]{Eisenbud} or \cite[Theorem 1.27]{GreenGIN}.
\end{proof}

Let $I \subset \K[x]$ be an ideal and $\sigma$ any term ordering. It is well known that the ideal $\IN_{\sigma}(I)$ has the same Hilbert function of $I$ and that there is a family over $\AA^1_{\K}$ having as fibers both $I$ and $\IN_{\sigma}(I)$. In the Hilbert schemes context, this means that the ideals $I$ and $\IN_{\sigma}(I)$ defines two $\K$-rational point of the Hilbert scheme $\Hilb{n}{p(t)}$, where $p(t)$ is the Hilbert polynomial of the subscheme $\Proj \K[x]/I$.
Since there is a flat deformation of $I$ which specializes to $\IN_{\sigma}(I)$, there is necessarily a component of $\Hilb{n}{p(t)}$ containing both points. Being embedded in a suitable Grassmannian, the Hilbert scheme turns out to be invariant under the action of the linear group (Proposition \ref{prop:phiFactors}), so it is interesting to face again with the problem of the relation between change of coordinates and initial ideal on the Hilbert scheme.

Any change of coordinates $g \in \GL_{\K}(n+1)$ of $\K[x]$ induces a linear action on the vector space $\K[x]_r$ of the homogeneous polynomials of degree $r,\ \forall\ r$ so also on the Hilbert scheme $\Hilb{n}{p(t)}$. Fixed any term ordering $\sigma$, the correspondence $I \mapsto \IN_{\sigma}(I)$ in general is not well-defined, indeed $I$ and $g.I$ define the same $\K$-rational point on $\Hilb{n}{p(t)}$ (up to isomorphism) but $\IN_{\sigma}(I)$ and $\IN_{\sigma}(g.I)$ could not.
But by Lemma \ref{lem:finiteClasses}, we know that the generic initial ideal $\GIN_{\sigma}(I)$ is stable for change of coordinates in an open subset $U \subset \GL_{\K}(n+1)$, and the same holds for the corresponding points on the Hilbert schemes. Hence the correspondence $I \mapsto \GIN_{\sigma}(I)$ results to be well-defined also in the context of Hilbert schemes.

 $\GIN_{\sigma}(I)$ lies on the same component of $I$ and if $I$ belongs to an intersection of components then $\GIN_{\sigma}(I)$ does too. So each components and each intersection of components of the Hilbert schemes has to contain at least one $\K$-rational point defined by a Borel-fixed ideal. The key role played by Borel ideals in the study of Hilbert schemes spring out from this remark, indeed we can consider the points defined by Borel ideals as distributed throughout the Hilbert scheme. Each Borel-fixed ideal can be used as the starting point of a local study of the Hilbert scheme (see Chapter \ref{ch:openSubsets}), whereas as a whole they can be used to investigate global properties (see Chapter \ref{ch:deformations}).

\section{Basic properties}

We now recall some of the main properties of Borel-fixed ideals that will be useful hereinafter. For this part we mainly refer to the first two sections of \cite{GreenGIN}.

\begin{definition}\index{saturated ideal}
Denoted by $\mathfrak{m}$ the irrelevant ideal $(x_0,\ldots,x_n)$, a homogeneous ideal $I \subset \K[x]$ is \emph{saturated} if $(I : \mathfrak{m}) = I$. Given a non-saturated ideal $J \subset \K[x]$, its saturation is the ideal
\begin{equation}
\glslink{saturation}{J^{\sat}} = \bigcup_{k \geqslant 0}(J : \mathfrak{m}^k).
\end{equation}
\end{definition}

\begin{proposition}\label{prop:saturationBorelIdeal}
For a Borel-fixed ideal $I \subset \K[x]$, $(I : \mathfrak{m}) = (I : x_0)$.
\end{proposition}
\begin{proof}
The inclusion $(I : \mathfrak{m}) \subseteq (I : x_0)$ is obvious. For any monomial $x^\alpha \in (I : x_0)$, $x^\alpha x_0$ belongs to $I$ and so by Proposition \ref{prop:BorelIdealProperties} $\frac{x_i}{x_0} x^\alpha x_0 = x_i x^\alpha,\ i = 1,\ldots,n$ belongs to $I$ too. Finally $x^\alpha \in (I : \mathfrak{m})$.
\end{proof}

\begin{corollary}\label{cor:saturationBorelIdeal}
A Borel-fixed ideal $I \subset \K[x]$ is saturated if the variable $x_0$ does not appear in any generator of $I$.
\end{corollary}

Therefore from a operative point of view, given a Borel ideal $I$, to compute its saturation we can consider any set of generators and impose $x_0 = 1$, indeed if the monomial $x^\alpha x_0^k$ ($x_0 \nmid x^\alpha$) belongs to $I$, the monomial $x^\alpha$ belongs to $(I : x_0^k) = (I : \mathfrak{m}^k)$ and so to $I^{\sat}$. 

\begin{definition}\label{def:regularityIdeal}\index{regularity of an ideal}
Given an ideal $I \subset \K[x]$ and a minimal free resolution
\[
0\ \rightarrow\ M_{n}\ \rightarrow\ \cdots\ \rightarrow\ M_1 \ \rightarrow\ M_0\ \rightarrow\ I\ \rightarrow 0,
\]
where $M_i = \oplus_j\, \K[x](-m_{ij})$, the \emph{regularity} of $I$ is $\max \{m_{ij}-i\}.$
\end{definition}

\begin{proposition}[{\cite[Proposition 2.11]{GreenGIN}, \cite[Proposition 2.9]{BayerStillmanREG}}]\label{prop:regularityDegreeGenerators}
The regularity of a Borel-fixed ideal $I \subset \K[x]$ is equal to the maximal degree of one of its generators.
\end{proposition}

Now we state a proposition giving a meaning of these properties of ideals in the context of schemes.

\begin{proposition}[{\cite[Proposition 2.6]{GreenGIN}}]\label{prop:regularityVsCMregularity}\index{Castelnuovo-Mumford regularity}
Let $I \subset \K[x]$ be a saturated ideal. The regularity of $I$ (as in Definition \ref{def:regularityIdeal}) is equal to the Castelnuovo-Mumford regularity of the sheaf of ideals $\mathscr{I}$ obtained from the sheafification of $I$.
\end{proposition}

Another important property of the last variable $x_0$ is the following.
\begin{proposition}\label{prop:hyperplaneSection}
Let $I \subset \K[x]$ be a Borel-fixed ideal. The linear form $x_0$ is regular for $I$, i.e. the hyperplane $H = \{x_0 = 0\}$ does not contain any irreducible component of $\Proj \K[x]/I$. Thus there is the short exact sequence induced by the multiplication by $x_0$
\begin{equation}\label{eq:hyperplaneSection}
0\ \longrightarrow\ \dfrac{\K[x]}{I}(t-1) \ \stackrel{\cdot x_0}{\longrightarrow}\ \dfrac{\K[x]}{I}(t) \ \longrightarrow\ \dfrac{\K[x]}{(I,x_0)}(t) \ \longrightarrow\ 0 
\end{equation}
\end{proposition}

Called $p(t)$ the Hilbert polynomial of the scheme $X = \Proj \K[x]/I$ of dimension $d$, it is well-known that the scheme obtained through the generic hyperplane section $X_H = X \cap H$ and defined by the ideal $(I,x_0)$ has dimension $d-1$, indeed its Hilbert polynomial is
\[
p_H(t) = \Delta p(t) = p(t) - p(t-1) \quad\text{and}\quad \deg p(t) = d\ \Rightarrow \deg p_H(t) = d-1.
\]

Moreover we remark that the ideal $(I,x_0) \cap \K[x_1,\ldots,x_n]$ is still Borel-fixed, because its monomials satisfy the characterization of Proposition \ref{prop:BorelIdealProperties}.

\begin{proposition}\label{prop:degreeHPpowerVars}
Let $I\subset \K[x]$ be a Borel-fixed ideal and let $p(t)$ be the Hilbert polynomial of $\Proj \K[x]/I$. 
The degree of $p(t)$ is equal to $\max \{ i\ \vert\ x_i^m \notin I,\ m\geqslant \reg(I)\} = \min \{j\ \vert x_j^m \in I,\ m\geqslant \reg(I)\}-1$.
\end{proposition}
\begin{proof}
Let us proceed by induction on the degree $d$ of the Hilbert polynomial. $x_0^m$ does not belong to the ideal $I$ for all $m$, because if it does, applying repeatedly Proposition \ref{prop:BorelIdealProperties} any other monomial should belong to the ideal, i.e. $I_{\geqslant m} = \K[x]_{\geqslant m}$. Considering the sequence in \eqref{eq:hyperplaneSection}, if $\deg p(t) =0$, the intersection between points and a generic hyperplane section should be empty, so that the ideal $(I,x_0)_{\geqslant m}$ coincides with $\K[x]_{\geqslant m}$. Hence $x_1^m,\ldots,x_n^m$ belong to $I$ and $\deg p(t) = 0 = \max \{i\ \vert\ x_i^m \notin I\} = \min \{j\ \vert\ x_j^m \in I\}-1$.

Let us suppose that the statement is true for $\deg p(t) = d-1$ and let us consider an ideal $I$ defining a scheme of dimension $d$. Again through the exact sequence in \eqref{eq:hyperplaneSection}, in order to obtain the ideal $(I,x_0)$ defining a subscheme of dimension $d-1$. Let us consider the map
\[
\begin{split}
\varphi: \K[x_0,\ldots,x_n]& \rightarrow \K[y_0,\ldots,y_{n-1}]\\
x_0& \mapsto 0\\
x_i& \mapsto y_{i-1},\quad i=1,\ldots,n.
\end{split}
\]
The ideal $\widetilde{I}=\phi((I,x_0))$ is still Borel-fixed and defines a subscheme isomorphic to the subscheme defined by $(I,x_0)$ and $\reg(\widetilde{I}) \leqslant \reg(I)$. By the inductive hypothesis we know that $d-1 = \max \{i\ \vert\ y_{i}^m \notin \widetilde{I}\} = \min \{j\ \vert\ y_j^{m} \in \widetilde{I}\}$ so that $x_{d+1}^m = \varphi^{-1}(y_{d}^m) \in I$ and $x_{d}^m = \varphi^{-1}(y_{d-1}^m) \notin I$.
\end{proof}

\begin{definition}\label{def:minMAXmonomials}
For any (non-constant) monomial $x^\alpha \in \K[x]$ we define
\begin{itemize}
\item $\gls{minMon} = \min \{i\ \text{s.t.}\ x_i\mid x^\alpha\}$;
\item $\gls{maxMon} = \max \{j\ \text{s.t.}\ x_j\mid x^\alpha\}$.
\end{itemize}
As we chose $x_n$ as the greatest variable and $x_0$ as the smallest one, the definition makes sense also in the following way
\begin{itemize}
\item $\gls{minMon} = \min \{x_i\ \text{s.t.}\ x_i\mid x^\alpha\}$;
\item $\gls{maxMon} = \max \{x_j\ \text{s.t.}\ x_j\mid x^\alpha\}$.
\end{itemize}
From now on we will use both definitions interchangeably.
\end{definition}

\begin{lemma}[{\cite[Lemma 1.1]{EliahouKervaire},\cite[Lemma 2.11]{MillerSturmfels}}]\label{lem:monomialDecomposition}
Let $I = ( x^{\alpha_1},\ldots,x^{\alpha_s}) \subset \K[x]$ be a Borel-fixed ideal. Each monomial $x^\beta \in I$ can be written uniquely as $x^{\alpha_i} x^{\gamma}$ so that $\min x^{\alpha_i} \geqslant \max x^\gamma$.
\end{lemma}
\begin{proof}
(Existence) Let us consider any decomposition of $x^\beta = x^{\alpha_i} x^\gamma$ and let us determine one of the type described in the statement. Suppose that $\min x^{\alpha_i} < \max x^{\gamma}$. By Proposition \ref{prop:BorelIdealProperties} we know that also $\dfrac{\max x^\gamma}{\min x^{\alpha_i}} x^{\alpha_i}$ belongs to $I$, so we can consider another minimal generator $x^{\alpha_j}$ of $I$ dividing $\frac{\max x^\gamma}{\min x^{\alpha_i}} x^{\alpha_i}$. By construction $\min x^{\alpha_j} \geqslant \min x^{\alpha_i}$, and either $\min x^{\alpha_j} > \min x^{\alpha_i}$ or $\min x^{\alpha_j} = \min x^{\alpha_i}$ and the degree of $\min x^{\alpha_j}$ in $x^{\alpha_j}$ is lower than the degree of $\min x^{\alpha_j}$ in $x^{\alpha_i}$. This procedure can not be repeated infinitely so at the end we will find a good decomposition.

(Uniqueness) Let us suppose that $x^\beta$ has two decomposition $x^{\alpha_i} x^{\gamma} = x^{\alpha_j} x^{\gamma'}$ such that $\min x^{\alpha_i} \geqslant \max x^\gamma$ and $\min x^{\alpha_j} \geqslant \max x^{\gamma'}$. If we suppese that $\min x^{\alpha_i} \geqslant \min x^{\alpha_j}$, then $x^{\alpha_i}$ and $x^{\alpha_j}$ are divided by the same power of each variable $x_k > \min x^{\alpha_i}$. If $\min x^{\alpha_i}$ divides $x^{\alpha_j}$, then $\min x^{\alpha_i} = \min x^{\alpha_j}$, hence either $x^{\alpha_i}$ divides $x^{\alpha_j}$ or viceversa and being both minimal generators they coincide. If $\max x^{\alpha_i}$ does not divide $x^{\alpha_j}$, then the degree of $\max x^{\alpha_i}$ in $x^{\alpha_i}$ is smaller than the degree of $\max x^{\alpha_i}$ in $x^{\alpha_j}$, i.e. the degree of $\max x^{\alpha_i}$ in $x^\beta$. Again $x^{\alpha_i}$ divides $x^{\alpha_j}$ but being both minimal generators of $I$, $x^{\alpha_i}=x^{\alpha_j}$. 
\end{proof}

\begin{definition}\label{def:canonicalDec}\index{canonical decomposition of a monomial of a Borel-fixed ideal} Given a Borel-fixed ideal $I =(x^{\alpha_1},\ldots,x^{\alpha_s}) \subset \K[x]$ and a monomial $x^{\beta} \in I$, the unique decomposition described in Lemma \ref{lem:monomialDecomposition} is called \emph{canonical decomposition} of $x^{\beta}$ over $I$ or \emph{canonical $I$-decomposition} of $x^\beta$ and we will denote it by
\[
x^\beta = \gls{borelDec},\qquad \min x^{\alpha_i} \geqslant \max x^{\gamma}
\]
\end{definition}

In \cite{EliahouKervaire} Eliahou and Kervaire introduce the \emph{decomposition function} $\gls{decFunc}$ from the set of monomials $M(I)$ of the Borel-fixed ideal $I$ to the minimal set of generators $G(I) = \{x^{\alpha_1},\ldots,x^{\alpha_s}\}$ of $I$:
\[
\begin{split}
\partial_{I}: \parbox{2.4cm}{\centering $M(I)$} & \longrightarrow G(I)\\
x^\beta = \dec{x^{\alpha_i}}{x^\gamma}{I}{} & \longmapsto x^{\alpha_i}
\end{split}
\] 
We recall the main properties of this function.
\begin{proposition}
Let $I \subset \K[x]$ be a Borel-fixed ideal and let $\partial_{I}: M(I) \rightarrow G(I)$ be its decomposition function. 
\begin{enumerate}
\item For any pairs of monomials $x^\beta, x^\gamma \in I$,
\begin{equation}
 \partial_{I} (x^{\beta} x^\gamma) = \partial_I (x^\beta)\quad \Longleftrightarrow\quad \min \partial_I(x^{\beta}) \geqslant \max x^{\gamma}.
\end{equation}
\item For any pairs of monomials $x^\beta \in I$ and $x^\gamma \in \K[x]$,
\begin{eqnarray}
& \partial_I\big(x^\gamma \partial_I(x^\beta)\big) = \partial_I (x^{\beta}x^{\gamma}),&\\
& \min \partial_I (x^{\beta}x^{\gamma}) \geqslant \min \partial_I (x^\beta),&\\
& \partial_I (x^{\beta}x^{\gamma}) \leq_{\mathtt{DegRevLex}} \partial_I (x^{\beta}).&
\end{eqnarray} 
\end{enumerate}
\end{proposition}

We finish this section describing the set of generators for the module of syzygies of a Borel-fixed ideal given in the free resolution constructed by Eliahou and Kervaire.

\begin{theorem}[{\cite[Theorem 2.1]{EliahouKervaire}}]\label{th:EKsyzygies}
Let $I \subset \K[x]$ be a Borel-fixed ideal generated by $(x^{\alpha_1},\ldots,x^{\alpha_s})$. The module of syzygies \gls{syz}, i.e. the kernel of the map
\[
\begin{array}{ccc}
\bigoplus\limits_{i=1,\ldots,s} \K[x]\big(-\vert \alpha_i \vert\big) & \rightarrow & \K[x]\\
         \mathbf{e}_i & \mapsto & x^{\alpha_i}
\end{array}
\]
is generated by the elements
\begin{equation}\label{eq:EKsyzygiesGenerator}
x_k \mathbf{e}_i - x^{\eta} \mathbf{e}_j,\quad\forall\ i=1,\ldots,s,\ \forall\ x_k > \min x^{\alpha_i}\text{ s.t. } x_k x^{\alpha_i} = \dec{x^{\alpha_j}}{x^{\eta}}{I}{}.
\end{equation}
\end{theorem}

\begin{example}
Let us consider the Borel ideal $I = (x_3^2,x_3 x_2,x_3 x_1^2,x_3 x_1 x_0,x_2^5,x_2^4 x_1)$ in $\K[x_0,x_1,x_2,x_3]$. The kernel of the map
\[
\K[x](-2)^2\oplus\K[x](-3)^2\oplus\K[x](-5)^2 \rightarrow \K[x],
\]
sending $\mathbf{e}_{1}\mapsto x_3^2,\mathbf{e}_{2}\mapsto x_3 x_2,\mathbf{e}_{3}\mapsto x_3 x_1^2,\mathbf{e}_{4}\mapsto x_3 x_1 x_0,\mathbf{e}_{5}\mapsto x_2^5,\mathbf{e}_{6}\mapsto x_2^4 x_1$, is generated by 9 elements.
\begin{itemize}
\item $\min x_3^2 = x_3$, so no generators of the type \eqref{eq:EKsyzygiesGenerator} can be found.
\item $\min x_3x_2 = x_2$, so there is the generator
\begin{enumerate}
\item $x_3 \mathbf{e}_2 - x_2 \mathbf{e_1}$, since $x_3 (x_3x_2) = \dec{x_3^2}{x_2}{I}{}$.
\end{enumerate}
\item $\min x_3x_1^2 = x_1$, so there are
\begin{enumerate}\setcounter{enumi}{1}
\item $x_2 \mathbf{e}_3 - x_1^2 \mathbf{e}_2$, since $x_2(x_3x_1^2) = \dec{x_3x_2}{x_1^2}{I}{}$;
\item $x_3 \mathbf{e}_3 - x_1^2 \mathbf{e}_1$, since $x_3(x_3x_1^2) = \dec{x_3^2}{x_1^2}{I}{}$.
\end{enumerate}
\item $\min x_3x_1x_0 = x_0$, so there are
\begin{enumerate}\setcounter{enumi}{3}
\item $x_1 \mathbf{e}_4 - x_1^2 \mathbf{e}_3$, since $x_1(x_3x_1x_0) = \dec{x_3x_1^2}{x_0}{I}{}$;
\item $x_2 \mathbf{e}_4 - x_1^2 \mathbf{e}_2$, since $x_2(x_3x_1x_0) = \dec{x_3x_2}{x_1x_0}{I}{}$;
\item $x_3 \mathbf{e}_4 - x_1^2 \mathbf{e}_1$, since $x_3(x_3x_1x_0) = \dec{x_3^2}{x_1x_0}{I}{}$.
\end{enumerate}
\item $\min x_2^5 = x_2$, so there is
\begin{enumerate}\setcounter{enumi}{6}
\item $x_3 \mathbf{e}_5 - x_2^4 \mathbf{e}_2$, since $x_3(x_2^5) = \dec{x_3x_2}{x_2^4}{I}{}$.
\end{enumerate}
\item $\min x_2^4 x_1 = x_1$, so there are	
\begin{enumerate}\setcounter{enumi}{7}
\item $x_2 \mathbf{e}_6 - x_1 \mathbf{e}_5$, since $x_2(x_2^4x_1) = \dec{x_2^4}{x_1}{I}{}$;
\item $x_3 \mathbf{e}_6 - x_2^3x_1 \mathbf{e}_2$, since $x_3(x_2^4x_1) = \dec{x_3x_2}{x_2^3x_1}{I}{}$.
\end{enumerate}
\end{itemize} 
\end{example}

\subsection{Basic manipulations of Hilbert polynomials}

\begin{definition}\label{def:deltaHilbPoly}
Let $p(t)$ be an admissible Hilbert polynomial. We define $\Delta^0 p(t) = p(t)$ and recursively
\begin{equation}\label{eq:deltaHilbPoly}
\gls{DeltaI} = \Delta^{i-1} p(t) -  \Delta^{i-1} p(t-1).
\end{equation}
\end{definition}

Directly from the definition, $\Delta p(t) = \Delta^1 p(t)$ and if $\deg p(t) = d$, $\Delta^{d+1} p(t) = 0$.

\begin{proposition}\label{prop:gotzmannNumberDelta}\index{Gotzmann number}
Let $p(t)$ be an admissible Hilbert polynomial. Let $r$ be its Gotzmann number and $\overline{r}$ the Gotzmann number of $\Delta p(t)$. Then $\overline{r} \leqslant r$.
\end{proposition}
\begin{proof}
Considered the Gotzmann representation \eqref{eq:GotzmannDecomposition} of $p(t)$, the representation of $\Delta p(t)$ is
\[
\begin{split}
\Delta p(t) = p(t) - p(t-1)&{} = \binom{t+a_1}{a_1} + \ldots + \binom{t+a_r - (r-1)}{a_r} \\
&{}\quad - \left(\binom{t-1+a_1}{a_1} + \ldots + \binom{t-1+a_r - (r-1)}{a_r}\right) =\\
&{} = \binom{t+(a_1-1)}{(a_1-1)} + \ldots + \binom{t+(a_r-1) - (r-1)}{a_r-1}
\end{split}
\]
where all the binomial coefficients with $a_i-1 < 0$ vanish, so the Gotzmann number $\overline{r}$ of $\Delta p(t)$ is equal to the number of coefficient $a_i \geqslant 1$.
\end{proof}

Taking inspiration from the proof of the previous proposition we introduce a new manipulation of Hilbert polynomials.
\begin{definition}
Given an admissible Hilbert polynomial $p(t)$ with Gotzmann representation as in \eqref{eq:GotzmannDecomposition}, we define
\begin{equation}
\gls{Sigma} = \binom{t+(a_1+1)}{a_1+1} + \ldots + \binom{t+(a_r+1) - (r-1)}{a_r+1}.
\end{equation}
\end{definition}

Obviously\hfill $\Sigma p(t)$\hfill and\hfill $p(t)$\hfill have\hfill the\hfill same\hfill Gotzmann\hfill number.\hfill Furthermore\\ $\Delta \big(\Sigma p\big)(t) = p(t)$, whereas $\Sigma\big(\Delta p\big)(t) = p(t) - c,\ c\geqslant 0$, indeed in the second case we lose the constant part corresponding in the Gotzmann decomposition to the binomial coefficients with $a_i = 0$.

\begin{definition}\label{def:minimalPolynomial}
Let $\overline{p}(t)$ be an admissible Hilbert polynomial and let us define the set 
\begin{equation}
\HP{\overline{p}(t)} = \left\{ p(t)\ \vert\ \Delta p(t) = \overline{p}(t)\right\}.
\end{equation}
The polynomial $\Sigma \overline{p}(t)$ belongs to $\HP{\overline{p}(t)}$ and
\[
p(t) = \Sigma \overline{p}(t) + c,\ c \geqslant 0,\qquad \forall\ p(t) \in \HP{\overline{p}(t)};
\]
thus we call $\Sigma \overline{p}(t)$ \emph{minimal polynomial} of $\HP{\overline{p}(t)}$.
\end{definition}

\begin{remark}\label{rk:differenceBetweenHPs}
By Proposition \ref{prop:gotzmannNumberDelta}, we deduce that $\Sigma \overline{p}(t)$ is the Hilbert polynomial in $\HP{\overline{p}(t)}$ with lowest Gotzmann number: let us denote it with $\overline{r}.$ For any other polynomial $p(t) \in \HP{\overline{p}(t)}$ with Gotzmann number $r$
\[
p(t) = \Sigma \overline{p}(t) + r - \overline{r}.
\]
\end{remark}

\section{The combinatorial interpretation}

In this section we provide a more combinatorial way to look at Borel-fixed ideals, that will turn out very useful in our algorithmic perspective.

\begin{definition}[{\cite[Definition 1.24]{GreenGIN}}]\label{def:elementaryMoves}\index{Borel elementary moves}
Let $\K[x_0,\ldots,x_n]$ be a polynomial ring and let $\K(x)$ its field of fraction. We define
\begin{itemize}
\item the \emph{$i$-th increasing elementary move} \gls{upMove} as the element $\dfrac{x_{i+1}}{x_i} \in \K(x),\ \forall\ i < n$;  
\item the \emph{$j$-th decreasing elementary move} \gls{downMove} as the element $\dfrac{x_{j-1}}{x_j} \in \K(x),\ \forall\ j >0$.
\end{itemize}
Given a monomial $x^\alpha \in \K[x]$, we will say that the elementary move $\up{i}$ (resp. $\down{j}$) is \emph{admissible} on $x^\alpha$ if $x_i \mid x^{\alpha}$ (resp. $x_j \mid x^\alpha$), i.e. $\up{i}(x^{\alpha}) = \frac{x_{i+1}}{x_i} x^\alpha \in \K[x]$ (resp. $\down{j}(x^{\alpha}) = \frac{x_{j-1}}{x_j} x^\alpha \in \K[x]$).

In the following, we will use an additive notation to denote the composition of an elementary move with itself, that is 
\begin{equation}
2\up{i} = \up{i}\circ\up{i} = \left(\frac{x_{i+1}}{x_i}\right)^2 \quad\text{and}\quad 2\down{j} := \down{j}\circ\down{j} = \left(\frac{x_{j-1}}{x_j}\right)^2
\end{equation}
and so on.
\end{definition}

In general, a composition of elementary moves turns out to be a \lq\lq monomial\rq\rq\ $x^\gamma\ (\gamma \in \ZZ^{n+1})$ of degree 0 in $\K(x)$ and we will say that it is admissible on $x^\alpha$ if the product $x^\gamma x^\alpha$ belongs to $\K[x]$.
Going by the commutativity of the product, given a composition $F = \lambda_a\up{i_a} \circ \cdots \circ \lambda_1 \up{i_1}$, we can suppose $i_1 < \ldots < i_a$ so that whenever $F$ is admissible, each elementary move in the written order is admissible. Similarly for any composition $G = \mu_b \down{j_b}\circ \cdots\circ \mu_1 \down{j_1}$, we will suppose $j_1 > \ldots > j_b$.

\begin{remark}\index{Borel-fixed ideal!combinatorial characterization of a}
Rewriting the characterization of Borel-fixed ideals given in \ref{prop:BorelIdealProperties}, we can say that an ideal $I$ is Borel-fixed if and only if its set of monomials is closed w.r.t. increasing elementary moves, indeed $\forall\ j > i$
\[
\frac{x_j}{x_i} = \frac{x_j}{x_{j-1}}\cdot\frac{x_{j-1}}{x_{j-2}} \cdots \frac{x_{i+2}}{x_{i+1}}\cdot\frac{x_{i+1}}{x_{i}} = \up{j-1} \circ \cdots \circ \up{i}.
\]
\end{remark}

\begin{definition}\label{def:BorelOrder}\index{Borel partial order}
We call \emph{Borel order}, and we denote it by \gls{BorelOrder}, the partial order defined on the set of monomials of a fixed degree by the transitive closure of the relations
\begin{equation}
\up{i}(x^{\alpha}) >_B x^\alpha >_B \down{j}(x^\alpha).
\end{equation}
\end{definition}

We note that the Borel order can be also obtained imposing the compatibility of the assumption $x_n > \ldots > x_0$ with the multiplication, because for any admissible elementary move $\up{i}$ on $x^\alpha$, set $x^{\overline{\alpha}} = \frac{x^\alpha}{x_i}$, we have that $x^{\alpha} = x_i x^{\overline{\alpha}}$, $\up{i}(x^\alpha) = x_{i+1}x^{\overline{\alpha}}$ and
\[
x_{i+1} > x_i \qquad\Longrightarrow\qquad  x_{i+1}x^{\overline{\alpha}} > x_i x^{\overline{\alpha}} \qquad\Longleftrightarrow\qquad \up{i}(x^{\alpha}) >_B x^{\alpha}.
\]
In the definition of a monomial order (see \cite[Definition 1.4.1]{KreuzerRobbiano1}, \cite[Definition 1]{CLOiva}),
the compatibility  between the order relation and the multiplication is always required, therefore any graded term ordering $\sigma$ is a total order on the monomials of fixed degree that refines the Borel partial order, that is
\[
x^{\alpha} >_B x^\beta \qquad\Longrightarrow\qquad x^{\alpha} >_{\sigma} x^\beta.
\]
We now characterize the Borel order by means of an analysis on the sets of exponents of monomials. Firstly, for any pair of multiindices $\alpha,\beta \in \NN^{n+1},\ \vert\alpha\vert = \vert\beta\vert$ and for any $0 \leqslant i \leqslant n$, we define the integer
\begin{equation}\label{eq:rhoMonomials}
\rho(\alpha,\beta,i) = \sum_{j=i}^n (\alpha_j - \beta_j).
\end{equation}

\begin{lemma}
Let $x^\alpha$ and $x^\beta$ be two monomials in $\K[x]_m$.
\begin{equation}
x^\alpha >_B x^\beta \qquad\Longleftrightarrow\qquad \rho(\alpha,\beta,i) \geqslant 0,\ \forall\ i=0,\ldots,n.
\end{equation}
\end{lemma}
\begin{proof}
($\Rightarrow$) $x^\alpha >_B x^\beta$ means $x^\beta = \mu_{b}\down{j_b} \circ \cdots \circ \mu_1\down{j_1}(x^\alpha)$, $j_1 > \ldots > j_b$. Obviously $\rho(\alpha,\beta,0) = \vert\alpha\vert-\vert\beta\vert = 0$. Moreover $\rho(\alpha,\beta,i) = 0,\ \forall\ i > j_1$ and $\rho(\alpha,\beta,j_1)$ has to be positive because $\alpha_{j_1} > \beta_{j_1}$. Let $x^\gamma = \mu_{1}\down{j_1}(x^\alpha)$. By definition $\gamma = (\alpha_0,\ldots,\alpha_{j_1-1} + (\alpha_{j_1}-\beta_{j_1}),\beta_{j_1},\ldots,\alpha_n)$, i.e. $\mu_1 = \alpha_{j_1}-\beta_{j_1}$, so that $\rho(\alpha,\beta,i) = \rho(\gamma,\beta,i),\ \forall\ i < j_1$. Repeating the reasoning on $x^\gamma >_B x^\beta = \mu_b \down{j_b}\circ \cdots \circ \mu_2\down{j_2}(x^\gamma)$, we prove $\rho(\alpha,\beta,i) \geqslant 0,\ \forall\ i$. 

($\Leftarrow$) It suffices to consider the composition of decreasing moves 
\[
G = \rho(\alpha,\beta,1)\down{1} \circ \cdots \circ \rho(\alpha,\beta,n)\down{n}.\qedhere
\]
\end{proof}

\begin{corollary}
Let $x^\alpha$ and $x^\beta$ be two monomials in $\K[x]_m$. They are not comparable w.r.t. the Borel order if there exists two integer $i,j$ such that $\rho(\alpha,\beta,i) \cdot \rho(\alpha,\beta,j) < 0$.
\end{corollary}

\begin{example}
Consider the polynomial ring $\K[x_0,x_1,x_2,x_3,x_4]$ and the monomial $x_4^3 x_3 x_2^3 x_0$. By definition
\[
x_4^3 x_3 x_2^3 x_0 >_B \down{1}\circ3\down{2}\circ2\down{4}(x_4^3 x_3 x_2^3 x_0) = x_4 x_3^3x_1^2 x_0^2
\]
and, set $\alpha = (1,0,3,1,3)$ and $\beta = (2,2,0,3,1)$,
\[
\rho(\alpha,\beta,0) = 0,\ \rho(\alpha,\beta,1) = 1,\ \rho(\alpha,\beta,2) = 3,\ \rho(\alpha,\beta,3) = 0,\ \rho(\alpha,\beta,4) = 2.
\]
Furthermore the monomials $x_4 x_3^2 x_2^3 x_0$ and $x_4^3 x_2 x_1^3$ are not comparable, indeed, set $\gamma=(1,0,3,2,1)$ and $\delta = (0,3,1,3,0)$,
\[
\rho(\gamma,\delta,0) = 0,\ \rho(\gamma,\delta,1) = -1,\ \rho(\gamma,\delta,2) = 2,\ \rho(\gamma,\delta,3) = 0,\ \rho(\gamma,\delta,4) = 1.
\]
\end{example}

\begin{definition}\label{def:poset}\index{poset}
We denote by \gls{pos} the \emph{Partially Ordered SET} (\emph{poset} for short) of the monomials of degree $m$ in the polynomial ring $\K[x_0,\ldots,x_n]$ with the Borel partial order $\leq_B$. 
\end{definition}

\begin{definition}\label{def:borelSet}\index{Borel set}
Following the characterization of Borel-fixed ideals in terms of elementary moves, we call \emph{Borel set} any subset $\gls{BorelSet} \subset \pos{n}{m}$ closed w.r.t. increasing elementary moves, i.e.
\[
x^\alpha \in \mathscr{B}\qquad\Longrightarrow\qquad \up{i}(x^\alpha) \in \mathscr{B},\ \forall\ \up{i} \text{ admissible on } x^\alpha. 
\]
With the terminology of orderings on sets, a Borel set represents a \emph{filter}\index{filter|see{Borel set}} of $\pos{n}{m}$ for the Borel partial order. Given a Borel-fixed ideal $I$, we will write \gls{BorelSetIdeal} referring to the Borel set defined by the piece of degree $m$ of the ideal $I$ in the poset $\pos{n}{m}$. 

Obviously the complement $\gls{orderSet} = \pos{n}{m}\setminus\mathscr{B}$, that we will also denote by $\mathscr{B}^{\mathcal{C}}$, is closed w.r.t. decreasing elementary moves. We will call such a subset an \emph{order set}\index{order set}, taking inspiration from the definition of order ideals, since the dehomogeneization of the complement of a Borel set (imposing $x_0 = 1$) turns out to be exactly an order ideal.
\end{definition}

Given any subset $S \subset \pos{n}{m}$, we will denote with $\restrict{S}{i}$ the subset of $S$
\begin{equation}
\gls{restrict} = \left\{ x^\alpha \in S\ \vert\ \min x^\alpha \geqslant i\right\}. 
\end{equation}
Obviously $\restrict{S}{0} = S$. Now we introduce some further definitions borrowed from the terminology of ordering on sets.

\begin{definition}\label{def:minimalMaximal}\index{minimal monomial}\index{maximal monomial}
Let $\mathscr{B}\subset \pos{n}{m}$ be a Borel set and let $\mathscr{N} = \mathscr{B}^{\mathcal{C}}$ the corresponding order set.
\begin{itemize}
\item $x^\alpha \in \mathscr{B}$ will be called \emph{minimal} if for any admissible decreasing move $\down{j}$, $\down{j}(x^\alpha)$ does not belong to $\mathscr{B}$.
\item $x^\beta \in \mathscr{N}$ will be called \emph{maximal} if for any admissible increasing move $\up{i}$, $\up{i}(x^\beta)$ belongs to $\mathscr{B}$.
\end{itemize}
Moreover we will say that
\begin{itemize}
\item $x^\alpha \in \mathscr{B}$ is \emph{$k$-minimal} if $\down{j}(x^\alpha) \in \mathscr{N}$ for any admissible $\down{j},\ j > k$; 
\item $x^\beta \in \mathscr{N}$ is \emph{$k$-maximal} if $\up{i}(x^\beta) \in \mathscr{B}$ for any admissible $\up{i},\ i \geqslant k$.
\end{itemize}
\end{definition}

\begin{example}
Let us consider the poset $\pos{2}{3}$ and its Borel subset
\[
\mathscr{B}=\{x_2^3,x_2^2 x_1,x_2x_1^2,x_2^2 x_0,x_2x_1x_0\},
\]
so that
\[
\mathscr{N} = \mathscr{B}^{\mathcal{C}} = \{x_1^3,x_1^2x_0,x_2x_0^2,x_1x_0^2,x_0^3\}.
\]
There is a single minimal element $x_2x_1x_0$ and two maximal elements: $x_2x_0^2$ and $x_1^3$. 
Moreover the $1$-minimal monomials are $x_2x_1^2$ and $x_2x_1x_0$ and the $1$-maximal ones are $x_1^3$ and $x_1^2x_0$.
\end{example}	

\begin{remark}\label{rk:minMax}
For any term ordering $\sigma$, refinement of the Borel order, and for any Borel set $\mathscr{B}$
\begin{itemize}
\item $\min_\sigma \restrict{\mathscr{B}}{k}$ is a minimal element of $\restrict{\mathscr{B}}{k}$;
\item $\max_\sigma \restrict{\mathscr{B}^{\mathcal{C}}}{k}$ is a maximal element of $\restrict{\mathscr{B}^{\mathcal{C}}}{k}$.
\end{itemize}
\end{remark}

\section{Graphical representations}

The combinatorial interpretation of Borel-fixed ideals leads up to nice representations of the posets of monomials of the same degree and their Borel subsets. We now briefly describe some different approaches emphasizing positive and negative aspects.

\paragraph{Green's diagrams}\index{Green's diagram|(}
We mainly refer to Section 4 of \cite{GreenGIN}. Green's diagrams can be used to describe few situations, indeed through this approach we can describe only Borel-fixed ideals defining points in $\PP^2$ or curves in $\PP^3$.

Let us begin with $\PP^2 = \Proj	\K[x_0,x_1,x_2]$, i.e. looking at posets of the type $\pos{2}{m}$. Green arranges the monomials of degree $m$ in a triangle shape with the top vertix corresponding to $x_0^m$ and completing the diagram moving down with the rule described in the following picture
\begin{center}
\begin{tikzpicture}[>=latex,scale=1.1]
 \node (top) at (0,0) [] {$x^\alpha$};
 \node (left) at (-1,-1.7) [] {$x^\beta$};
 \node (right) at (1,-1.7) [] {$x^\gamma$};
 \draw [->] (top) --node[fill=white]{\footnotesize $\up{0}$} (right);
 \draw [->] (top) --node[fill=white]{\footnotesize $\up{1}\circ\up{0}$} (left);
 \draw [->] (right) --node[fill=white]{\footnotesize $\up{1}$} (left);
\end{tikzpicture}
\end{center}
so that the base of the triangle contains the monomials of degree $m$ in $\K[x_1,x_2]$ (see Figure \ref{fig:exampleGreenPosetP2}). 

\begin{figure}[!ht]
\begin{center}
\begin{tikzpicture}[scale=0.55]
\node (r0) at (0,0) [] {\footnotesize $x_0^4$};

\node (r10) at (-1,-1.7) [] {\footnotesize $x_2 x_0^3$};
\node (r11) at (1,-1.7) [] {\footnotesize $x_1 x_0^3$};

\node (r20) at (-2,-3.4) [] {\footnotesize $x_2^2 x_0^2$};
\node (r21) at (0,-3.4) [] {\footnotesize $x_2x_1x_0^2$};
\node (r22) at (2,-3.4) [] {\footnotesize $x_1^2x_0^2$};

\node (r30) at (-3,-5.1) [] {\footnotesize $x_2^3 x_0$};
\node (r31) at (-1,-5.1) [] {\footnotesize $x_2^2 x_1 x_0$};
\node (r32) at (1,-5.1) [] {\footnotesize $x_2 x_1^2 x_0$};
\node (r33) at (3,-5.1) [] {\footnotesize $x_1^3 x_0$};

\node (r40) at (-4,-6.8) [] {\footnotesize $x_2^4$};
\node (r41) at (-2,-6.8) [] {\footnotesize $x_2^3 x_1$};
\node (r42) at (0,-6.8) [] {\footnotesize $x_2^2 x_1^2$};
\node (r43) at (2,-6.8) [] {\footnotesize $x_2 x_1^3$};
\node (r44) at (4,-6.8) [] {\footnotesize $x_1^4$};
\end{tikzpicture}
\caption[The Green's diagram describing the poset $\pos{2}{4}$.]{\label{fig:exampleGreenPosetP2} An example of Green's diagram for $\PP^2$: the poset $\pos{2}{4}$.}
\end{center}
\end{figure}
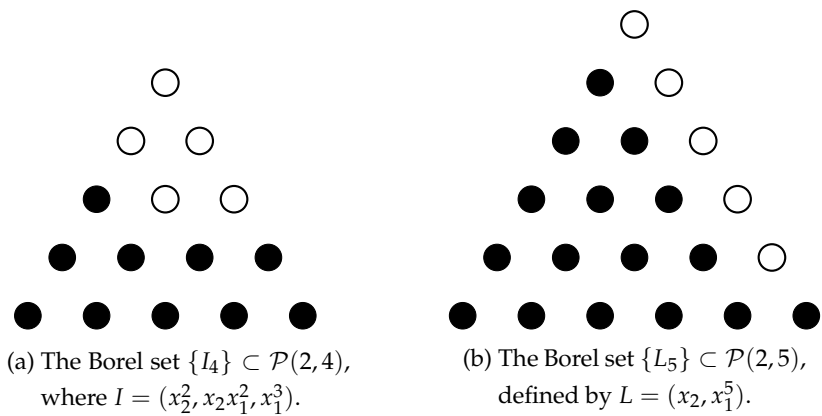
\begin{figure}[!ht]
\begin{center}
\captionsetup[subfloat]{singlelinecheck=false,width=4.5cm,format=hang}
\subfloat[][The Borel set $\{I_{4}\} \subset \pos{2}{4}$, where $I = (x_2^2,x_2 x_1^2,x_1^3).$]{\label{fig:exampleGreenPoints_a}
\begin{tikzpicture}[scale=0.45]
\tikzstyle{ideal}=[circle,draw=black,fill=black,inner sep=3.5pt]
\tikzstyle{quotient}=[circle,draw=black,thick,inner sep=3.5pt]

\node (r0) at (0,0) [quotient] {};

\node (r10) at (-1,-1.7) [quotient] {};
\node (r11) at (1,-1.7) [quotient] {};

\node (r20) at (-2,-3.4) [ideal] {};
\node (r21) at (0,-3.4) [quotient] {};
\node (r22) at (2,-3.4) [quotient] {};

\node (r30) at (-3,-5.1) [ideal] {};
\node (r31) at (-1,-5.1) [ideal] {};
\node (r32) at (1,-5.1) [ideal] {};
\node (r33) at (3,-5.1) [ideal] {};

\node (r40) at (-4,-6.8) [ideal] {};
\node (r41) at (-2,-6.8) [ideal] {};
\node (r42) at (0,-6.8) [ideal] {};
\node (r43) at (2,-6.8) [ideal] {};
\node (r44) at (4,-6.8) [ideal] {};
\end{tikzpicture}
}
\qquad\qquad
\subfloat[][The Borel set $\{L_5\} \subset \pos{2}{5}$, defined by $L = (x_2,x_1^5)$.]{\label{fig:exampleGreenPoints_b}
\begin{tikzpicture}[scale=0.45]
\tikzstyle{ideal}=[circle,draw=black,fill=black,inner sep=3.5pt]
\tikzstyle{quotient}=[circle,draw=black,thick,inner sep=3.5pt]

\node (r0) at (0,0) [quotient] {};

\node (r10) at (-1,-1.7) [ideal] {};
\node (r11) at (1,-1.7) [quotient] {};

\node (r20) at (-2,-3.4) [ideal] {};
\node (r21) at (0,-3.4) [ideal] {};
\node (r22) at (2,-3.4) [quotient] {};

\node (r30) at (-3,-5.1) [ideal] {};
\node (r31) at (-1,-5.1) [ideal] {};
\node (r32) at (1,-5.1) [ideal] {};
\node (r33) at (3,-5.1) [quotient] {};

\node (r40) at (-4,-6.8) [ideal] {};
\node (r41) at (-2,-6.8) [ideal] {};
\node (r42) at (0,-6.8) [ideal] {};
\node (r43) at (2,-6.8) [ideal] {};
\node (r44) at (4,-6.8) [quotient] {};

\node (r50) at (-5,-8.5) [ideal] {};
\node (r51) at (-3,-8.5) [ideal] {};
\node (r52) at (-1,-8.5) [ideal] {};
\node (r53) at (1,-8.5) [ideal] {};
\node (r54) at (3,-8.5) [ideal] {};
\node (r55) at (5,-8.5) [ideal] {};
\end{tikzpicture}
}
\caption[Green's diagrams of Borel sets defined by Borel-fixed ideals of points in $\PP^2$.]{\label{fig:exampleGreenPoints} Example of Green's diagrams of Borel sets defined by Borel-fixed ideals of points in $\K[x_0,x_1,x_2]$.}
\end{center}
\end{figure}

Afterwards he does not write explicitly the monomials, and given an ideal $I\subset \K[x_0,x_1,x_2]$ he uses a black circle to denote a monomial of the ideal and an empty circle to denote a monomial not belonging to $I$. In Figure \ref{fig:exampleGreenPoints}, there are two examples of Borel sets defined by Borel-fixed ideals.

\begin{figure}[!ht]
\begin{center}
\begin{tikzpicture}[scale=0.44]

\node (u0) at (0,0) [] {\footnotesize $x_0^3$};

\node (u10) at (-2,2) [] {\footnotesize $x_2 x_0^2$};
\node (u11) at (2,2) [] {\footnotesize $x_3 x_0^2$};

\node (u20) at (-4,4) [] {\footnotesize $x_2^2 x_0$};
\node (u21) at (0,4) [] {\footnotesize $x_2 x_3 x_0$};
\node (u22) at (4,4) [] {\footnotesize $x_3^2 x_0$};

\node (u30) at (-6,6) [] {\footnotesize $x_2^3$};
\node (u31) at (-2,6) [] {\footnotesize $x_2^2 x_3$};
\node (u32) at (2,6) [] {\footnotesize $x_2 x_3^2$};
\node (u33) at (6,6) [] {\footnotesize $x_3^3$};

\node (f0) at (0,-4) [] {\footnotesize $x_0^3$};
\node (f01) at (0,-6) [] {\footnotesize $x_1 x_0^2$};
\node (f02) at (0,-8) [] {\footnotesize $x_1^2 x_0$};
\node (f03) at (0,-10) [] {\footnotesize $x_1^3$};

\node (f10) at (-2,-4) [black!80] {\footnotesize $x_2 x_0^2$};
\node (f101) at (-2,-6) [black!80] {\footnotesize $x_2 x_1 x_0$};
\node (f102) at (-2,-8) [black!80] {\footnotesize $x_2 x_1^2$};

\node (f11) at (2,-4) [black!80] {\footnotesize $x_3 x_0^2$};
\node (f111) at (2,-6) [black!80] {\footnotesize $x_3 x_1 x_0$};
\node (f112) at (2,-8) [black!80] {\footnotesize $x_3 x_1^2$};

\node (f20) at (-4,-4) [black!60] {\footnotesize $x_2^2 x_0$};
\node (f201) at (-4,-6) [black!60] {\footnotesize $x_2^2 x_1$};

\node (f22) at (4,-4) [black!60] {\footnotesize $x_3^2 x_0$};
\node (f221) at (4,-6) [black!60] {\footnotesize $x_3^2 x_1$};

\node (f30) at (-6,-4) [black!40] {\footnotesize $x_2^3$};

\node (f33) at (6,-4) [black!40] {\footnotesize $x_3^3$};

\node (l0) at (10,-4) [black!40] {\footnotesize $x_0^3$}; 
\node (l01) at (10,-6) [black!40] {\footnotesize $x_1 x_0^3$}; 
\node (l02) at (10,-8) [black!40] {\footnotesize $x_1^2 x_0$}; 
\node (l03) at (10,-10) [black!40] {\footnotesize $x_1^3$}; 

\node (l11) at (12,-4) [black!60] {\footnotesize $x_3 x_0^2$};
\node (l111) at (12,-6) [black!60] {\footnotesize $x_3 x_1 x_0$};
\node (l112) at (12,-8) [black!60] {\footnotesize $x_3 x_1^2$};

\node (l22) at (14,-4) [black!80] {\footnotesize $x_3^2 x_0$};
\node (l221) at (14,-6) [black!80] {\footnotesize $x_3^2 x_1$};

\node (l33) at (16,-4) [] {\footnotesize $x_3^3$};

\draw [thick] (-8,-2) -- (18,-2);
\draw [thick] (8,7) -- (8,-11);

\node at (-8.1,-1.5) [] {\footnotesize \textsc{plane view}};
\node at (-8.1,-2.5) [] {\footnotesize \textsc{front view}};
\node at (18,-2.5) [] {\footnotesize \textsc{side view}};

\draw [very thin,black!20] (u0) -- (f0);
\draw [very thin,black!20] (u10) -- (f10);
\draw [very thin,black!20] (u11) -- (f11);
\draw [very thin,black!20] (u20) -- (f20);
\draw [very thin,black!20] (u22) -- (f22);
\draw [very thin,black!20] (u30) -- (f30);
\draw [very thin,black!20] (u33) -- (f33);

\draw [very thin,black!20] (f33) -- (l0);
\draw [very thin,black!20] (f221) -- (l01);
\draw [very thin,black!20] (f112) -- (l02);
\draw [very thin,black!20] (f03) -- (l03);

\draw [very thin,black!20] (u0) -- (8,0);
\draw [very thin,black!20] (u11) -- (8,2);
\draw [very thin,black!20] (u22) -- (8,4);
\draw [very thin,black!20] (u33) -- (8,6);

\draw [very thin,black!20] (l0) -- (10,-2);
\draw [very thin,black!20] (l11) -- (12,-2);
\draw [very thin,black!20] (l22) -- (14,-2);
\draw [very thin,black!20] (l33) -- (16,-2);

\draw [very thin,black!20] (10,-2) arc (0:90:2);
\draw [very thin,black!20] (12,-2) arc (0:90:4);
\draw [very thin,black!20] (14,-2) arc (0:90:6);
\draw [very thin,black!20] (16,-2) arc (0:90:8);
\end{tikzpicture}
\caption[The Green's trihedron describing the poset $\pos{3}{3}$.]{\label{fig:GreenPyramid} The trihedron describing the poset $\pos{3}{3}$ with orthographic projections.}
\end{center}
\end{figure}
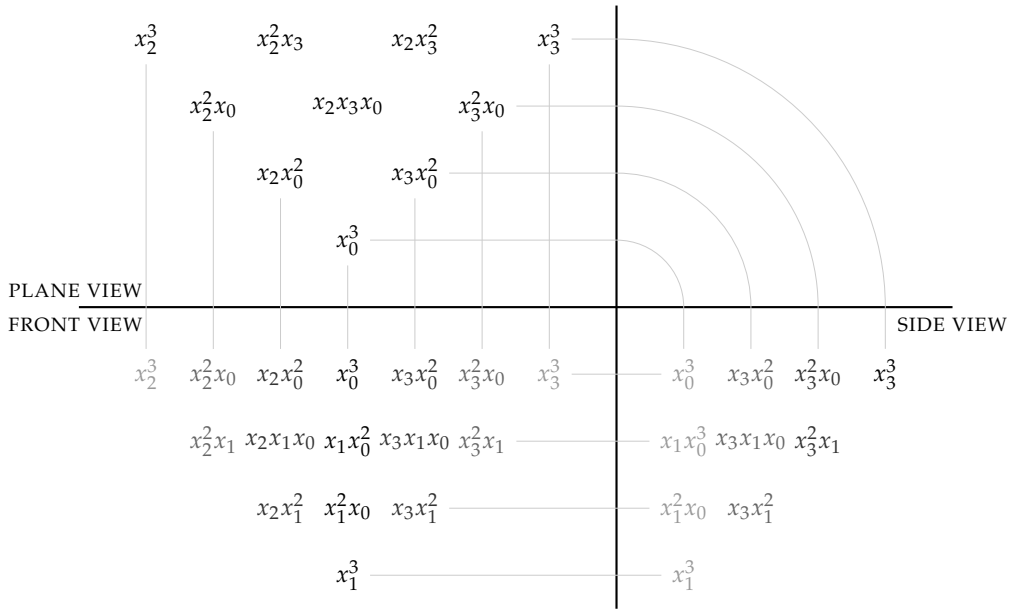

To describe a Borel set $\mathscr{B} \subset \pos{3}{m}$, Green thinks a trihedron (drawn with orthographic projections in Figure \ref{fig:GreenPyramid}) described looking at its plane view with top vertix corresponding to the monomial $x_0^m$ and completed with the following rule
\begin{center}
\begin{tikzpicture}[>=latex,scale=1.1]
 \node (top) at (0,0) [] {$x^\alpha$};
 \node (left) at (-1.5,-1.7) [] {$x^\beta$};
 \node (right) at (1.5,-1.7) [] {$x^\gamma$};
 \draw [->] (top) --node[fill=white]{\footnotesize $\ \up{1}\circ\up{0}$} (right);
 \draw [->] (top) --node[fill=white]{\footnotesize $\up{2}\circ\up{1}\circ\up{0}\ $} (left);
 \draw [->] (right) --node[fill=white]{\footnotesize $\up{2}$} (left);
\end{tikzpicture}
\end{center}

Then the monomials are marked according to the following notation:
\begin{itemize}
\item a black circle denotes a monomial in $\mathscr{B}$ such that also all the other monomials under it belong to $\mathscr{B}$;
\item a empty circle denotes a monomial in $\mathscr{B}^\mathcal{C}$ such that all the monomials under it do not belong to $\mathscr{B}$;
\item a empty circle with inside a positive integer $\lambda$ denotes a monomial $x^\alpha$ in $\mathscr{B}^{\mathcal{C}}$ such that the monomial (under it) $\lambda\up{0}(x^\alpha)$ belongs to $\mathscr{B}$ and $(\lambda-1)\up{0}(x^\alpha)$ does not.
\end{itemize}

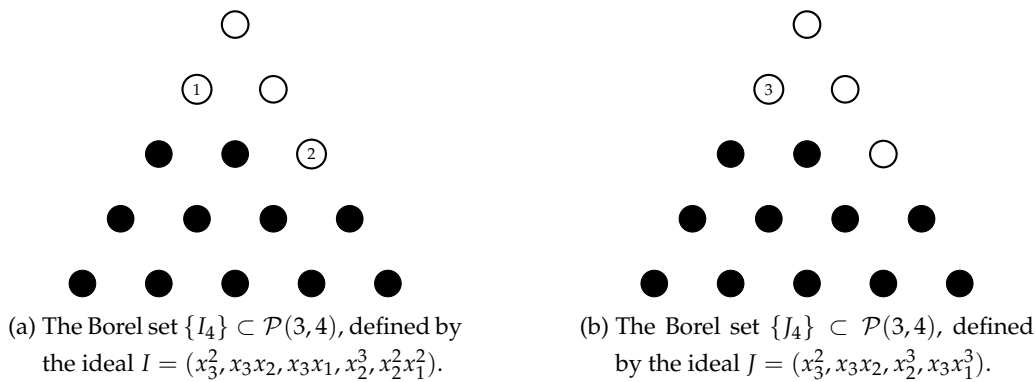
\begin{figure}[!ht]
\begin{center}
\captionsetup[subfloat]{singlelinecheck=false,format=hang}
\subfloat[][The Borel set $\{I_4\} \subset \pos{3}{4}$, defined by the ideal $I = (x_3^2,x_3x_2,x_3x_1,x_2^3,x_2^2x_1^2)$.]{
\begin{tikzpicture}[scale=0.5]
\tikzstyle{ideal}=[circle,draw=black,fill=black,inner sep=3.5pt]
\tikzstyle{quotient}=[circle,draw=black,thick,inner sep=3.5pt]
\tikzstyle{mixed}=[circle,draw=black,thick,inner sep=2pt]
\node at (-5.5,0) [] {};
\node at (5.5,0) [] {};

\node (r0) at (0,0) [quotient] {};

\node (r10) at (-1,-1.7) [mixed] {\tiny $1$};
\node (r11) at (1,-1.7) [quotient] {};

\node (r20) at (-2,-3.4) [ideal] {};
\node (r21) at (0,-3.4) [ideal] {};
\node (r22) at (2,-3.4) [mixed] {\tiny $2$};

\node (r30) at (-3,-5.1) [ideal] {};
\node (r31) at (-1,-5.1) [ideal] {};
\node (r32) at (1,-5.1) [ideal] {};
\node (r33) at (3,-5.1) [ideal] {};

\node (r40) at (-4,-6.8) [ideal] {};
\node (r41) at (-2,-6.8) [ideal] {};
\node (r42) at (0,-6.8) [ideal] {};
\node (r43) at (2,-6.8) [ideal] {};
\node (r44) at (4,-6.8) [ideal] {};
\end{tikzpicture}
}
\qquad\qquad
\subfloat[][The Borel set $\{J_4\} \subset \pos{3}{4}$, defined by the ideal $J = (x_3^2,x_3x_2,x_2^3,x_3x_1^3)$.]{
\begin{tikzpicture}[scale=0.5]
\tikzstyle{ideal}=[circle,draw=black,fill=black,inner sep=3.5pt]
\tikzstyle{quotient}=[circle,draw=black,thick,inner sep=3.5pt]
\tikzstyle{mixed}=[circle,draw=black,thick,inner sep=2pt]

\node at (-5.5,0) [] {};
\node at (5.5,0) [] {};

\node (r0) at (0,0) [quotient] {};

\node (r10) at (-1,-1.7) [mixed] {\tiny $3$};
\node (r11) at (1,-1.7) [quotient] {};

\node (r20) at (-2,-3.4) [ideal] {};
\node (r21) at (0,-3.4) [ideal] {};
\node (r22) at (2,-3.4) [quotient] {};

\node (r30) at (-3,-5.1) [ideal] {};
\node (r31) at (-1,-5.1) [ideal] {};
\node (r32) at (1,-5.1) [ideal] {};
\node (r33) at (3,-5.1) [ideal] {};

\node (r40) at (-4,-6.8) [ideal] {};
\node (r41) at (-2,-6.8) [ideal] {};
\node (r42) at (0,-6.8) [ideal] {};
\node (r43) at (2,-6.8) [ideal] {};
\node (r44) at (4,-6.8) [ideal] {};
\end{tikzpicture}
}
\caption[Green's diagrams of Borel sets defined by Borel-fixed ideals of curves in $\PP^3$.]{\label{fig:exampleGreenCurves} Example of Green's diagrams of Borel sets defined by Borel-fixed ideals of curves in $\K[x_0,x_1,x_2,x_3]$.}
\end{center}
\end{figure}

We remark that in the diagram corresponding to a saturated monomial ideal $I$, a monomial $x^\alpha$ marked with a black circle imposes that also every other monomial under it and contained in a triangle with $x^\alpha$ as top vertix belongs to $I$. Indeed splitting $x^\alpha$ as $x^{\overline{\alpha}} x_0^a$ ($x_0 \nmid x^{\overline{\alpha}}$), the black circle means $x^{\overline{\alpha}} \in I$ and any monomial in the triangle under it is divided by $x^{\overline{\alpha}}$. From now on, we will only draw the black circles defining the saturation of an ideal. Moreover if the ideal $I$ is Borel-fixed also any monomial at the left of $x^\alpha$ has to be marked with a black circle. Thinking about the quotient, any monomial above or at the right of a monomial marked with a empty circle does not belong to the ideal. We can summarize this characterization with the following diagram
\begin{center}
\begin{tikzpicture}[scale=0.5]
\tikzstyle{ideal}=[circle,draw=black,fill=black,inner sep=3.5pt]
\tikzstyle{quotient}=[circle,draw=black,thick,inner sep=3.5pt]

\filldraw[fill=black!20,draw=white] (0,0) -- (-4,0) arc (180:301:4) -- cycle;
\filldraw[fill=black!10,draw=white] (7,-3.4) -- (11,-3.4) arc (0:121:4) -- cycle;

\node at (0,0) [ideal] {};
\draw [-] (0,0) -- (-4,0);
\draw [-] (0,0) -- (2,-3.4);

\node at (7,-3.4) [quotient] {};
\draw [-] (7,-3.4) -- (11,-3.4);
\draw [-] (7,-3.4) -- (5,0);

\node at (-1,-2) [] {ideal}; 
\node at (8,-1.4) [] {quotient};
\end{tikzpicture}
\end{center}

This type of diagram works very well in the context of curves in $\PP^3$, because we can understand many geometrical information about a curve simply looking at its diagram. For instance the number of empty circles corresponds to the degree of the curve, indeed it is easy to check that from the diagram of a curve in $\PP^3$, the diagram of the hyperplane section with $H\vert_{x_0 = 0}$, i.e. points in $\PP^2$, can be obtained substituting the empty circles with an integer inside with black circles.

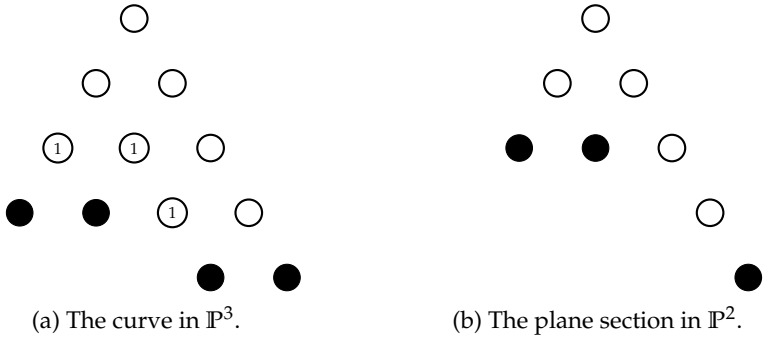
\begin{figure}[!ht]
\begin{center}
\captionsetup[subfloat]{justification=centering}
\subfloat[][The curve in $\PP^3$.]{\label{fig:curveSection_a}
\begin{tikzpicture}[scale=0.5]
\tikzstyle{ideal}=[circle,draw=black,fill=black,inner sep=3.5pt]
\tikzstyle{quotient}=[circle,draw=black,thick,inner sep=3.5pt]
\tikzstyle{mixed}=[circle,draw=black,thick,inner sep=2pt]

\node (r0) at (0,0) [quotient] {};

\node (r10) at (-1,-1.7) [quotient] {};
\node (r11) at (1,-1.7) [quotient] {};

\node (r20) at (-2,-3.4) [mixed] {\tiny $1$};
\node (r21) at (0,-3.4) [mixed] {\tiny $1$};
\node (r22) at (2,-3.4) [quotient] {};

\node (r30) at (-3,-5.1) [ideal] {};
\node (r31) at (-1,-5.1) [ideal] {};
\node (r32) at (1,-5.1) [mixed] {\tiny $1$};
\node (r33) at (3,-5.1) [quotient] {};

\node at (-4,-6.8) [] {};
\node (r43) at (2,-6.8) [ideal] {};
\node (r44) at (4,-6.8) [ideal] {};
\end{tikzpicture}
}
\qquad\qquad
\subfloat[][The plane section in $\PP^2$.]{\label{fig:curveSection_b}
\begin{tikzpicture}[scale=0.5]
\tikzstyle{ideal}=[circle,draw=black,fill=black,inner sep=3.5pt]
\tikzstyle{quotient}=[circle,draw=black,thick,inner sep=3.5pt]
\tikzstyle{mixed}=[circle,draw=black,thick,inner sep=2pt]

\node (r0) at (0,0) [quotient] {};

\node (r10) at (-1,-1.7) [quotient] {};
\node (r11) at (1,-1.7) [quotient] {};

\node (r20) at (-2,-3.4) [ideal] {};
\node (r21) at (0,-3.4) [ideal] {};
\node (r22) at (2,-3.4) [quotient] {};

\node (r33) at (3,-5.1) [quotient] {};

\node at (-4,-6.8) [] {};
\node (r44) at (4,-6.8) [ideal] {};
\end{tikzpicture}
}
\caption[Green's diagram of a curve and of its plane section.]{\label{fig:curveSection} Green's diagrams of the curve in $\PP^3$ and its plane section in $\PP^2$ described in Example \ref{ex:curveSection}.}
\end{center}
\end{figure}

\begin{example}\label{ex:curveSection}
Let us consider the curve defined by the ideal 
\[
I = (x_3^3,x_3^2x_2,x_3^2x_1,x_3x_2^3,x_2^4,x_3x_2^2x_1,x_3x_2x_1^2) \subset \K[x_0,x_1,x_2,x_3].
\]
Its Hilbert polynomial is $p(t) = 5t+2$, i.e. the curve has degree $5$ and genus $-1$.
The ideal defining the plane section of the curve with the plane $H\vert_{x_0=0}$ turns out to be
\[
J = (x_3^2,x_3x_2,x_2^4) \subset \K[x_1,x_2,x_3]
\]
and defines 5 points in $\PP^2$ as expected.
The two diagram are drawn in Figure \ref{fig:curveSection}.
\end{example}

\index{Green's diagram|)}

\paragraph{Marinari's lattices}\index{Marinari's lattice|(}

Another way to represent posets in 3 or 4 variables was taught to me by Maria Grazia Marinari, that with some collegues has worked extensively on Borel-fixed ideals (see \cite{MarinariBUMI,MarinariRamella1999,MarinariRamella2005,MarinariRamella2006}). 
As for Green's diagrams, this approach works only for ideals defining subschemes in $\PP^2$ and $\PP^3$, but without any restriction on the degree of the Hilbert polynomial.

The poset $\pos{2}{m}$ is described again through a triangle shape with the bottom right vertix corresponding to the monomial $x_0^{m}$ and then moving up and left with the following rule
\begin{center}
\begin{tikzpicture}[>=latex,scale=1]
\node (1) at (0,0) [] {$x^\alpha$};
\node (2) at (0,2) [] {$x^\beta$};
\node (3) at (-2,2) [] {$x^\gamma$};
\draw [->] (1) -- node[fill=white]{$\up{0}$} (2);
\draw [->] (2) -- node[fill=white]{$\up{1}$} (3);
\end{tikzpicture}
\end{center}

\begin{figure}
\begin{center}
\begin{tikzpicture}[>=latex,scale=0.75]

\node (1) at (0,0) [] {\footnotesize $x_2^4$};
\node (2) at (2,0) [] {\footnotesize $x_2^3 x_1$};
\node (3) at (4,0) [] {\footnotesize $x_2^2 x_1^2$};
\node (4) at (6,0) [] {\footnotesize $x_2 x_1^3$};
\node (5) at (8,0) [] {\footnotesize $x_1^4$};

\node (6) at (2,-1.5) [] {\footnotesize $x_2^3 x_0$};
\node (7) at (4,-1.5) [] {\footnotesize $x_2^2 x_1 x_0$};
\node (8) at (6,-1.5) [] {\footnotesize $x_2 x_1^2 x_0$};
\node (9) at (8,-1.5) [] {\footnotesize $x_1^3 x_0$};

\node (10) at (4,-3) [] {\footnotesize $x_2^2 x_0^2$};
\node (11) at (6,-3) [] {\footnotesize $x_2 x_1 x_0^2$};
\node (12) at (8,-3) [] {\footnotesize $x_1^2 x_0^2$};

\node (13) at (6,-4.5) [] {\footnotesize $x_2 x_0^3$};
\node (14) at (8,-4.5) [] {\footnotesize $x_1 x_0^3$};

\node (15) at (8,-6) [] {\footnotesize $x_0^4$};

\draw [->] (15) -- (14);
\draw [->] (14) -- (13);
\draw [->] (14) -- (12);
\draw [->] (13) -- (11);
\draw [->] (12) -- (11);
\draw [->] (12) -- (9);
\draw [->] (11) -- (10);
\draw [->] (11) -- (8);
\draw [->] (10) -- (7);
\draw [->] (9) -- (8);
\draw [->] (9) -- (5);
\draw [->] (8) -- (7);
\draw [->] (8) -- (4);
\draw [->] (7) -- (6);
\draw [->] (7) -- (3);
\draw [->] (6) -- (2);
\draw [->] (5) -- (4);
\draw [->] (4) -- (3);
\draw [->] (3) -- (2);
\draw [->] (2) -- (1);
\end{tikzpicture}
\caption[The Marinari's lattice describing the poset $\pos{2}{4}$.]{\label{fig:exampleMarinariPosetP2} An example of Marinari's lattice for $\PP^2$: the poset $\pos{2}{4}$ (cf. Figure \ref{fig:exampleGreenPosetP2}).}
\end{center}
\end{figure}

In the following we will not write explicitly the monomials and we will consider the lattice without the verse of the arrows. Given a Borel set $\mathscr{B} \subset \pos{2}{m}$, we will denote again with a black circle a monomial belonging to $\mathscr{B}$ and with a empty circle a monomial not belonging.

\begin{figure}[!ht]
\begin{center}
\captionsetup[subfloat]{singlelinecheck=false,format=hang}
\subfloat[][The Borel set $\{I_{4}\} \subset \pos{2}{4}$, where $I = (x_2^2,x_2 x_1^2,x_1^3).$]{\label{fig:exampleMarinariPoints_a}
\begin{tikzpicture}[scale=0.55]
\tikzstyle{ideal}=[circle,draw=black,fill=black,inner sep=1.5pt]
\tikzstyle{quotient}=[circle,draw=black,thick,inner sep=1.5pt]

\node at (-2,0) [] {};
\node at (6,0) [] {};

\node (00) at (0,0) [ideal] {};
\node (01) at (1,0) [ideal] {};
\node (02) at (2,0) [ideal] {};
\node (03) at (3,0) [ideal] {};
\node (04) at (4,0) [ideal] {};

\node (11) at (1,-1) [ideal] {};
\node (12) at (2,-1) [ideal] {};
\node (13) at (3,-1) [ideal] {};
\node (14) at (4,-1) [ideal] {};

\node (22) at (2,-2) [ideal] {};
\node (23) at (3,-2) [quotient] {};
\node (24) at (4,-2) [quotient] {};

\node (33) at (3,-3) [quotient] {};
\node (34) at (4,-3) [quotient] {};

\node (44) at (4,-4) [quotient] {};

\node at (4,-5) [] {};

\draw [-,black!50] (00) -- (01);
\draw [-,black!50] (01) -- (02);
\draw [-,black!50] (02) -- (03);
\draw [-,black!50] (03) -- (04);

\draw [-,black!50] (01) -- (11);
\draw [-,black!50] (02) -- (12);
\draw [-,black!50] (03) -- (13);
\draw [-,black!50] (04) -- (14);

\draw [-,black!50] (11) -- (12);
\draw [-,black!50] (12) -- (13);
\draw [-,black!50] (13) -- (14);

\draw [-,black!50] (12) -- (22);
\draw [-,black!50] (13) -- (23);
\draw [-,black!50] (14) -- (24);

\draw [-,black!50] (22) -- (23);
\draw [-,black!50] (23) -- (24);

\draw [-,black!50] (23) -- (33);
\draw [-,black!50] (24) -- (34);

\draw [-,black!50] (33) -- (34);
\draw [-,black!50] (34) -- (44);
\end{tikzpicture}
}
\qquad\qquad
\subfloat[][The Borel set $\{L_5\} \subset \pos{2}{5}$, defined by $L = (x_2,x_1^5)$.]{\label{fig:exampleMarinariPoints_b}
\begin{tikzpicture}[scale=0.55]
\tikzstyle{ideal}=[circle,draw=black,fill=black,inner sep=1.5pt]
\tikzstyle{quotient}=[circle,draw=black,thick,inner sep=1.5pt]

\node at (-3,0) [] {};
\node at (6,0) [] {};

\node (n)  at (-1,1) [ideal] {};
\node (n0) at (0,1) [ideal] {};
\node (n1) at (1,1) [ideal] {};
\node (n2) at (2,1) [ideal] {};
\node (n3) at (3,1) [ideal] {};
\node (n4) at (4,1) [ideal] {};

\node (00) at (0,0) [ideal] {};
\node (01) at (1,0) [ideal] {};
\node (02) at (2,0) [ideal] {};
\node (03) at (3,0) [ideal] {};
\node (04) at (4,0) [quotient] {};

\node (11) at (1,-1) [ideal] {};
\node (12) at (2,-1) [ideal] {};
\node (13) at (3,-1) [ideal] {};
\node (14) at (4,-1) [quotient] {};

\node (22) at (2,-2) [ideal] {};
\node (23) at (3,-2) [ideal] {};
\node (24) at (4,-2) [quotient] {};

\node (33) at (3,-3) [ideal] {};
\node (34) at (4,-3) [quotient] {};

\node (44) at (4,-4) [quotient] {};

\draw [-,black!50] (n) -- (n0);
\draw [-,black!50] (n0) -- (n1);
\draw [-,black!50] (n1) -- (n2);
\draw [-,black!50] (n2) -- (n3);
\draw [-,black!50] (n3) -- (n4);

\draw [-,black!50] (n0) -- (00);
\draw [-,black!50] (n1) -- (01);
\draw [-,black!50] (n2) -- (02);
\draw [-,black!50] (n3) -- (03);
\draw [-,black!50] (n4) -- (04);

\draw [-,black!50] (00) -- (01);
\draw [-,black!50] (01) -- (02);
\draw [-,black!50] (02) -- (03);
\draw [-,black!50] (03) -- (04);

\draw [-,black!50] (01) -- (11);
\draw [-,black!50] (02) -- (12);
\draw [-,black!50] (03) -- (13);
\draw [-,black!50] (04) -- (14);

\draw [-,black!50] (11) -- (12);
\draw [-,black!50] (12) -- (13);
\draw [-,black!50] (13) -- (14);

\draw [-,black!50] (12) -- (22);
\draw [-,black!50] (13) -- (23);
\draw [-,black!50] (14) -- (24);

\draw [-,black!50] (22) -- (23);
\draw [-,black!50] (23) -- (24);

\draw [-,black!50] (23) -- (33);
\draw [-,black!50] (24) -- (34);

\draw [-,black!50] (33) -- (34);
\draw [-,black!50] (34) -- (44);
\end{tikzpicture}
}
\caption[Marinari's lattices representing Borel sets defined by Borel-fixed ideals of points in $\PP^2$.]{\label{fig:exampleMarinariPoints} Example of Marinari's lattices representing Borel sets defined by Borel-fixed ideals of points in $\K[x_0,x_1,x_2]$ (cf. Figure \ref{fig:exampleGreenPoints}).}
\end{center}
\end{figure}

To describe posets in 4 variables, the idea is to decompose $\pos{3}{m}$ as
\[
\pos{3}{m} = \pos{2}{m} \cup x_0 \cdot \pos{2}{m-1} \cup \cdots \cup x_{0}^{m-1} \cdot \pos{2}{1} \cup \{x_0^m\}
\]
and to draw each $\pos{2}{i}$ one above the other, so that a monomial $x^\alpha \in x_0^{m-i}\pos{2}{i}$ has above it the monomial $\up{0}(x^\alpha) \in x_0^{m-i-1}\pos{2}{i+1}$ and beaneath it $\down{0}(x^\alpha) \in x_0^{m-i+1}\pos{2}{i-1}$. 

\begin{figure}[!ht]
\begin{center}
\begin{tikzpicture}[scale=0.55]
\node (0r00) at (0,0) [] {\footnotesize $x_3^3$};
\node (0r01) at (2,0) [] {\footnotesize $x_3^2 x_2$};
\node (0r02) at (4,0) [] {\footnotesize $x_3 x_2^2$};
\node (0r03) at (6,0) [] {\footnotesize $x_2^3$};

\node (0r11) at (3,-1.5) [] {\footnotesize $x_3^2 x_1$};
\node (0r12) at (5,-1.5) [] {\footnotesize $x_3 x_2 x_1$};
\node (0r13) at (7,-1.5) [] {\footnotesize $x_2^2 x_1$};

\node (0r22) at (6,-3) [] {\footnotesize $x_3 x_1^2$};
\node (0r23) at (8,-3) [] {\footnotesize $x_2 x_1^2$};

\node (0r33) at (9,-4.5) [] {\footnotesize $x_1^3$};

\draw [-,black!50] (0r00) -- (0r01);
\draw [-,black!50] (0r01) -- (0r02);
\draw [-,black!50] (0r02) -- (0r03);

\draw [-,black!50] (0r01) -- (0r11);
\draw [-,black!50] (0r02) -- (0r12);
\draw [-,black!50] (0r03) -- (0r13);

\draw [-,black!50] (0r11) -- (0r12);
\draw [-,black!50] (0r12) -- (0r13);

\draw [-,black!50] (0r12) -- (0r22);
\draw [-,black!50] (0r13) -- (0r23);

\draw [-,black!50] (0r22) -- (0r23);

\draw [-,black!50] (0r23) -- (0r33);

\node (1r11) at (3,-4.5) [] {\footnotesize $x_3^2 x_0$};
\node (1r12) at (5,-4.5) [] {\footnotesize $x_3 x_2 x_0$};
\node (1r13) at (7,-4.5) [] {\footnotesize $x_2^2 x_0$};

\node (1r22) at (6,-6) [] {\footnotesize $x_3 x_1 x_0$};
\node (1r23) at (8,-6) [] {\footnotesize $x_2 x_1 x_0$};

\node (1r33) at (9,-7.5) [] {\footnotesize $x_1^2 x_0$};

\draw [-,black!50] (1r11) -- (1r12);
\draw [-,black!50] (1r12) -- (1r13);

\draw [-,black!50] (1r12) -- (1r22);
\draw [-,black!50] (1r13) -- (1r23);

\draw [-,black!50] (1r22) -- (1r23);

\draw [-,black!50] (1r23) -- (1r33);

\node (2r22) at (6,-8.5) [] {\footnotesize $x_3 x_0^2$};
\node (2r23) at (8,-8.5) [] {\footnotesize $x_2 x_0^2$};

\node (2r33) at (9,-10) [] {\footnotesize $x_1 x_0^2$};

\draw [-,black!50] (2r22) -- (2r23);

\draw [-,black!50] (2r23) -- (2r33);

\node (3r33) at (9,-12.5) [] {\footnotesize $x_0^3$};

\draw [-,black,thick,loosely dotted] (0r11) -- (1r11);
\draw [-,black,thick,loosely dotted] (0r12) -- (1r12);
\draw [-,black,thick,loosely dotted] (0r13) -- (1r13);

\draw [-,black,thick,loosely dotted] (0r22) -- (1r22);
\draw [-,black,thick,loosely dotted] (0r23) -- (1r23);

\draw [-,black,thick,loosely dotted] (0r33) -- (1r33);

\draw [-,black,thick,loosely dotted] (1r22) -- (2r22);
\draw [-,black,thick,loosely dotted] (1r23) -- (2r23);

\draw [-,black,thick,loosely dotted] (1r33) -- (2r33);

\draw [-,black,thick,loosely dotted] (2r33) -- (3r33);

\end{tikzpicture}
\caption[The Marinari's lattice describing the poset {$\pos{3}{3}$}.]{\label{fig:MarinariPyramid} The Marinari's lattice describing the poset $\pos{3}{3}$ (cf. Figure \ref{fig:GreenPyramid}).}
\end{center}
\end{figure}

With this representation, given an ideal defining a curve, the ideal of its plane section (with $H\vert_{x_0 = 0}$) in $\K[x_1,x_2,x_3]$ is directly described in the highest triangle.

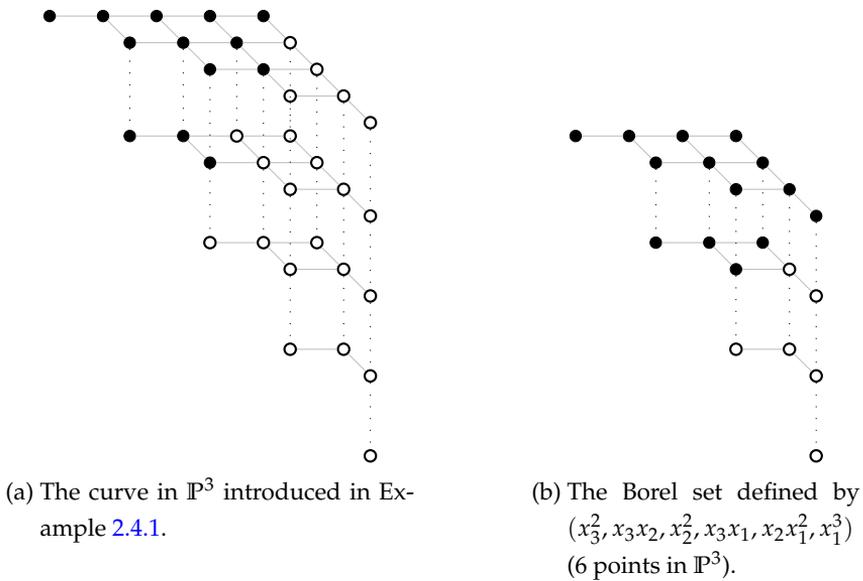
\begin{figure}[!ht]
\begin{center}
\captionsetup[subfloat]{singlelinecheck=false,format=hang}
\subfloat[][The curve in $\PP^3$ introduced in Example \ref{ex:curveSection}.]{\label{fig:curveMarinari_a}
\begin{tikzpicture}[scale=0.35]
\tikzstyle{ideal}=[circle,draw=black,fill=black,inner sep=1.5pt]
\tikzstyle{quotient}=[circle,draw=black,thick,inner sep=1.5pt]

\node at (-4,0) [] {};
\node at (10,0) [] {};

\node (0rnn) at (-3,1) [ideal] {};
\node (0rn0) at (-1,1) [ideal] {};
\node (0rn1) at (1,1) [ideal] {};
\node (0rn2) at (3,1) [ideal] {};
\node (0rn3) at (5,1) [ideal] {};

\node (0r00) at (0,0) [ideal] {};
\node (0r01) at (2,0) [ideal] {};
\node (0r02) at (4,0) [ideal] {};
\node (0r03) at (6,0) [quotient] {};

\node (0r11) at (3,-1) [ideal] {};
\node (0r12) at (5,-1) [ideal] {};
\node (0r13) at (7,-1) [quotient] {};

\node (0r22) at (6,-2) [quotient] {};
\node (0r23) at (8,-2) [quotient] {};

\node (0r33) at (9,-3) [quotient] {};

\draw [-,black!25] (0rnn) -- (0rn0);
\draw [-,black!25] (0rn0) -- (0rn1);
\draw [-,black!25] (0rn1) -- (0rn2);
\draw [-,black!25] (0rn2) -- (0rn3);

\draw [-,black!25] (0rn0) -- (0r00);
\draw [-,black!25] (0rn1) -- (0r01);
\draw [-,black!25] (0rn2) -- (0r02);
\draw [-,black!25] (0rn3) -- (0r03);

\draw [-,black!25] (0r00) -- (0r01);
\draw [-,black!25] (0r01) -- (0r02);
\draw [-,black!25] (0r02) -- (0r03);

\draw [-,black!25] (0r01) -- (0r11);
\draw [-,black!25] (0r02) -- (0r12);
\draw [-,black!25] (0r03) -- (0r13);

\draw [-,black!25] (0r11) -- (0r12);
\draw [-,black!25] (0r12) -- (0r13);

\draw [-,black!25] (0r12) -- (0r22);
\draw [-,black!25] (0r13) -- (0r23);

\draw [-,black!25] (0r22) -- (0r23);

\draw [-,black!25] (0r23) -- (0r33);

\node (1r00) at (0,-3.5) [ideal] {};
\node (1r01) at (2,-3.5) [ideal] {};
\node (1r02) at (4,-3.5) [quotient] {};
\node (1r03) at (6,-3.5) [quotient] {};

\node (1r11) at (3,-4.5) [ideal] {};
\node (1r12) at (5,-4.5) [quotient] {};
\node (1r13) at (7,-4.5) [quotient] {};

\node (1r22) at (6,-5.5) [quotient] {};
\node (1r23) at (8,-5.5) [quotient] {};

\node (1r33) at (9,-6.5) [quotient] {};

\draw [-,black!25] (1r00) -- (1r01);
\draw [-,black!25] (1r01) -- (1r02);
\draw [-,black!25] (1r02) -- (1r03);

\draw [-,black!25] (1r01) -- (1r11);
\draw [-,black!25] (1r02) -- (1r12);
\draw [-,black!25] (1r03) -- (1r13);

\draw [-,black!25] (1r11) -- (1r12);
\draw [-,black!25] (1r12) -- (1r13);

\draw [-,black!25] (1r12) -- (1r22);
\draw [-,black!25] (1r13) -- (1r23);

\draw [-,black!25] (1r22) -- (1r23);

\draw [-,black!25] (1r23) -- (1r33);

\node (2r11) at (3,-7.5) [quotient] {};
\node (2r12) at (5,-7.5) [quotient] {};
\node (2r13) at (7,-7.5) [quotient] {};

\node (2r22) at (6,-8.5) [quotient] {};
\node (2r23) at (8,-8.5) [quotient] {};

\node (2r33) at (9,-9.5) [quotient] {};

\draw [-,black!25] (2r11) -- (2r12);
\draw [-,black!25] (2r12) -- (2r13);

\draw [-,black!25] (2r12) -- (2r22);
\draw [-,black!25] (2r13) -- (2r23);

\draw [-,black!25] (2r22) -- (2r23);

\draw [-,black!25] (2r23) -- (2r33);

\node (3r22) at (6,-11.5) [quotient] {};
\node (3r23) at (8,-11.5) [quotient] {};

\node (3r33) at (9,-12.5) [quotient] {};

\draw [-,black!25] (3r22) -- (3r23);

\draw [-,black!25] (3r23) -- (3r33);

\node (4r44) at (9,-15.5) [quotient] {};

\draw [-,black,loosely dotted] (0r00) -- (1r00);
\draw [-,black,loosely dotted] (0r01) -- (1r01);
\draw [-,black,loosely dotted] (0r02) -- (1r02);
\draw [-,black,loosely dotted] (0r03) -- (1r03);

\draw [-,black,loosely dotted] (0r11) -- (1r11);
\draw [-,black,loosely dotted] (0r12) -- (1r12);
\draw [-,black,loosely dotted] (0r13) -- (1r13);

\draw [-,black,loosely dotted] (0r22) -- (1r22);
\draw [-,black,loosely dotted] (0r23) -- (1r23);

\draw [-,black,loosely dotted] (0r33) -- (1r33);

\draw [-,black,loosely dotted] (1r11) -- (2r11);
\draw [-,black,loosely dotted] (1r12) -- (2r12);
\draw [-,black,loosely dotted] (1r13) -- (2r13);

\draw [-,black,loosely dotted] (1r22) -- (2r22);
\draw [-,black,loosely dotted] (1r23) -- (2r23);

\draw [-,black,loosely dotted] (1r33) -- (2r33);

\draw [-,black,loosely dotted] (2r22) -- (3r22);
\draw [-,black,loosely dotted] (2r23) -- (3r23);

\draw [-,black,loosely dotted] (2r33) -- (3r33);

\draw [-,black,loosely dotted] (3r33) -- (4r44);

\end{tikzpicture}
}
\qquad\qquad
\subfloat[][The Borel set defined by $(x_3^2,x_3x_2,x_2^2,x_3x_1,x_2x_1^2,x_1^3)$ (6 points in $\PP^3$).]{\label{fig:curveMarinari_b}
\begin{tikzpicture}[scale=0.35]
\tikzstyle{ideal}=[circle,draw=black,fill=black,inner sep=1.5pt]
\tikzstyle{quotient}=[circle,draw=black,thick,inner sep=1.5pt]

\node at (-1,-5) [] {};
\node at (10,-5) [] {};

\node (1r00) at (0,-3.5) [ideal] {};
\node (1r01) at (2,-3.5) [ideal] {};
\node (1r02) at (4,-3.5) [ideal] {};
\node (1r03) at (6,-3.5) [ideal] {};

\node (1r11) at (3,-4.5) [ideal] {};
\node (1r12) at (5,-4.5) [ideal] {};
\node (1r13) at (7,-4.5) [ideal] {};

\node (1r22) at (6,-5.5) [ideal] {};
\node (1r23) at (8,-5.5) [ideal] {};

\node (1r33) at (9,-6.5) [ideal] {};

\draw [-,black!25] (1r00) -- (1r01);
\draw [-,black!25] (1r01) -- (1r02);
\draw [-,black!25] (1r02) -- (1r03);

\draw [-,black!25] (1r01) -- (1r11);
\draw [-,black!25] (1r02) -- (1r12);
\draw [-,black!25] (1r03) -- (1r13);

\draw [-,black!25] (1r11) -- (1r12);
\draw [-,black!25] (1r12) -- (1r13);

\draw [-,black!25] (1r12) -- (1r22);
\draw [-,black!25] (1r13) -- (1r23);

\draw [-,black!25] (1r22) -- (1r23);

\draw [-,black!25] (1r23) -- (1r33);

\node (2r11) at (3,-7.5) [ideal] {};
\node (2r12) at (5,-7.5) [ideal] {};
\node (2r13) at (7,-7.5) [ideal] {};

\node (2r22) at (6,-8.5) [ideal] {};
\node (2r23) at (8,-8.5) [quotient] {};

\node (2r33) at (9,-9.5) [quotient] {};

\draw [-,black!25] (2r11) -- (2r12);
\draw [-,black!25] (2r12) -- (2r13);

\draw [-,black!25] (2r12) -- (2r22);
\draw [-,black!25] (2r13) -- (2r23);

\draw [-,black!25] (2r22) -- (2r23);

\draw [-,black!25] (2r23) -- (2r33);

\node (3r22) at (6,-11.5) [quotient] {};
\node (3r23) at (8,-11.5) [quotient] {};

\node (3r33) at (9,-12.5) [quotient] {};

\draw [-,black!25] (3r22) -- (3r23);

\draw [-,black!25] (3r23) -- (3r33);

\node (4r44) at (9,-15.5) [quotient] {};

\draw [-,black,loosely dotted] (1r11) -- (2r11);
\draw [-,black,loosely dotted] (1r12) -- (2r12);
\draw [-,black,loosely dotted] (1r13) -- (2r13);

\draw [-,black,loosely dotted] (1r22) -- (2r22);
\draw [-,black,loosely dotted] (1r23) -- (2r23);

\draw [-,black,loosely dotted] (1r33) -- (2r33);

\draw [-,black,loosely dotted] (2r22) -- (3r22);
\draw [-,black,loosely dotted] (2r23) -- (3r23);

\draw [-,black,loosely dotted] (2r33) -- (3r33);

\draw [-,black,loosely dotted] (3r33) -- (4r44);

\end{tikzpicture}
}
\caption[Marinari's lattices of a curve and of points in $\PP^3$.]{\label{fig:curveMarinari} Borel sets of a curve and points in $\PP^3$ in the Marinari's lattice.}
\end{center}
\end{figure}
\index{Marinari's lattice|)}

\pagebreak

\begin{figure}[!ht]
\begin{center}
\includegraphics[height=5cm]{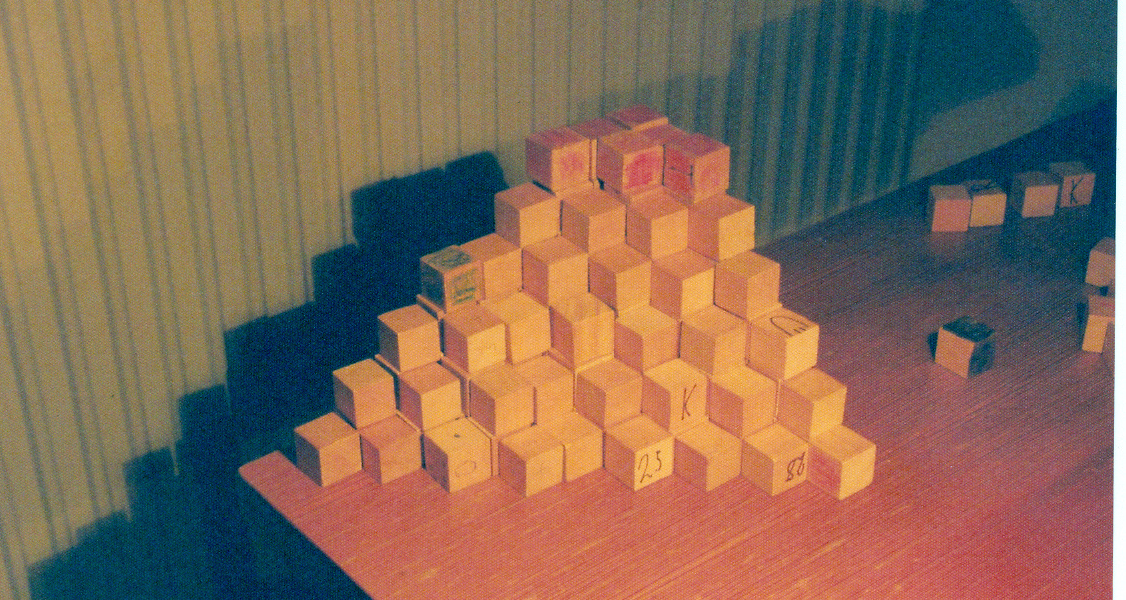}
\caption{\label{fig:gotzmannPhoto} A photo of a Borel set sent to me by Gotzmann.}
\end{center}
\end{figure}

\paragraph{Gotzmann's pyramids}\index{Gotzmann's pyramid|(}
Gotzmann showed to me another slightly different way to draw posets of the type $\pos{2}{m}$ and $\pos{3}{m}$ using two and three dimensonal spaces. He (literally) builds two or three dimensional pyramids (see Figure \ref{fig:gotzmannPhoto})  representing monomials as bricks (squares for $\pos{2}{m}$ and cubes for $\pos{3}{m}$) starting from the monomial $x_2^m$ in the case of $\pos{2}{m}$ and from $x_3^m$ in the case of $\pos{3}{m}$ with the following expansion rules:
\begin{center}
\begin{tikzpicture}[>=latex,scale=0.5]
\draw [thick,->] (0,0) -- (4,0) -- node[above]{\footnotesize $\down{2}$} (6,0);
\draw (0,2) -- (4,2);
\draw [thick,->] (0,0) -- (0,4) -- node[right]{\footnotesize $\down{1}\circ\down{2}$} (0,6);
\draw (2,0) -- (2,4);
\draw (0,4) -- (2,4);
\draw (4,0) -- (4,2);

\end{tikzpicture}
\qquad\qquad
\begin{tikzpicture}[>=latex,scale=0.4]
\draw [thick,->] (0,4) --node[right]{\footnotesize $\down{1}\circ\down{2}\circ\down{3}$} (0,6);

\draw (0,-2) -- (2,-3);
\draw (0,-2) -- (-2,-3);

\draw (-4,-2) -- (-2,-3);
\draw (4,-2) -- (2,-3);

\draw (-2,1) -- (-4,0);
\draw (2,1) -- (4,0);

\draw (-2,1) -- (2,-1);
\draw (2,1) -- (-2,-1);

\draw (-4,0) -- (-2,-1);
\draw (4,0) -- (2,-1);

\draw (0,2) -- (0,-2);
\draw (-2,-1) -- (-2,-3);
\draw (2,-1) -- (2,-3);
\draw (-4,0) -- (-4,-2);
\draw (4,0) -- (4,-2);

\draw (-2,1) -- (-2,3);
\draw (2,1) -- (2,3);

\draw (0,2) -- (2,3) -- (0,4) -- (-2,3) -- cycle;

\draw [thick,->] (-4,-2) --node[above left]{\footnotesize $\down{3}$} (-6,-3);
\draw [thick,->] (4,-2) --node[above right]{\footnotesize $\down{2}\circ\down{3}$} (6,-3);
 
\end{tikzpicture}
\end{center}

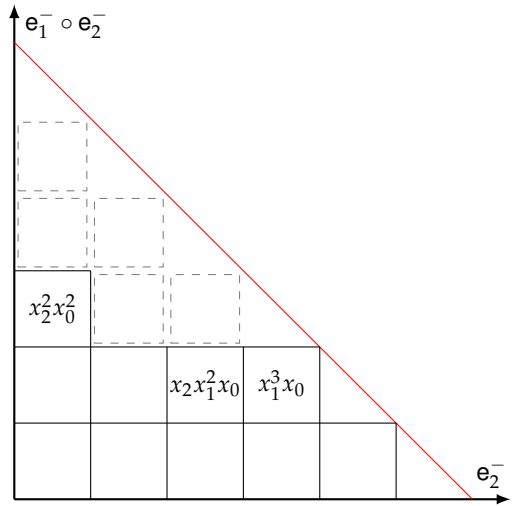
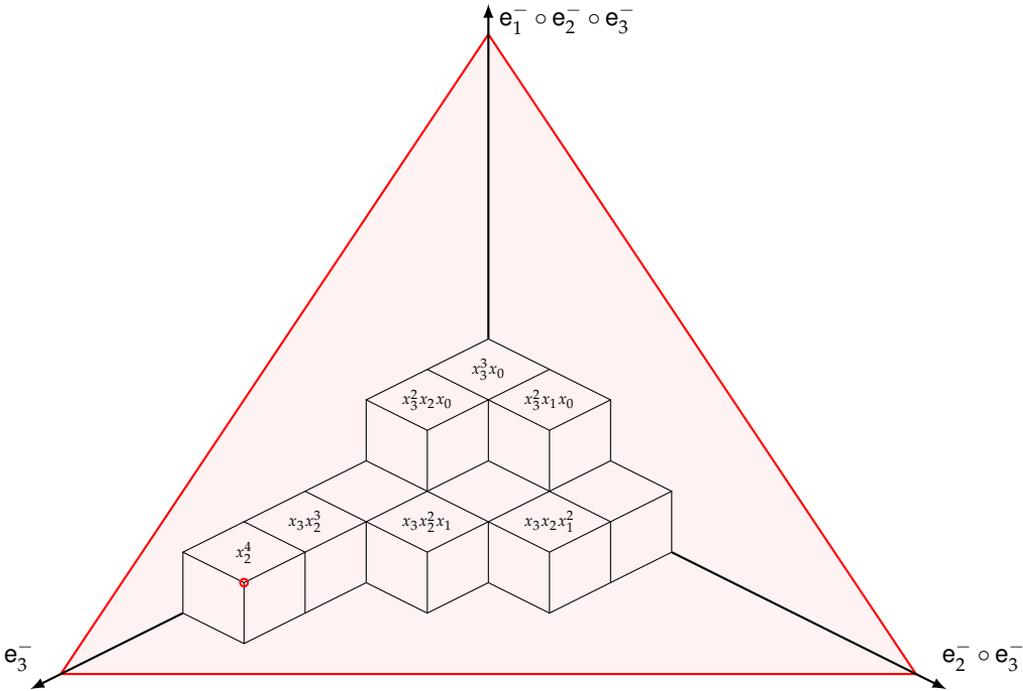
\begin{figure}[H]
\begin{center}
\captionsetup[subfloat]{singlelinecheck=false,format=hang,width=\linewidth,justification=centering}
\subfloat[][The Borel set defined by $(x_2^2,x_2x_1^2,x_1^3)$ in $\pos{2}{4}$ (cf. Figure \ref{fig:exampleMarinariPoints_a} and Figure \ref{fig:exampleGreenPoints_a}).]{\label{fig:pointsCurvesGotzmann_a}
\begin{tikzpicture}[>=latex,scale=0.5]
\draw [->,thick] (0,0) -- (12,0) --node[above]{\footnotesize $\down{2}$} (13,0);
\draw [->,thick] (0,0) -- (0,12) --node[right]{\footnotesize $\down{1}\circ\down{2}$} (0,13);

\draw [red] (0,12) -- (12,0);
\draw (0,2) -- (10,2);
\draw (0,4) -- (8,4);
\draw (0,6) -- (2,6);
\draw (2,0) -- (2,6);
\draw (4,0) -- (4,4);
\draw (6,0) -- (6,4);
\draw (8,0) -- (8,4);
\draw (10,0) -- (10,2);

\draw [dashed,black!50] (0.1,6.1) -- (1.9,6.1) -- (1.9,7.9) -- (0.1,7.9) -- cycle;
\draw [dashed,black!50] (0.1,8.1) -- (1.9,8.1) -- (1.9,9.9) -- (0.1,9.9) -- cycle;
\draw [dashed,black!50] (2.1,6.1) -- (3.9,6.1) -- (3.9,7.9) -- (2.1,7.9) -- cycle;
\draw [dashed,black!50] (2.1,4.1) -- (3.9,4.1) -- (3.9,5.9) -- (2.1,5.9) -- cycle;
\draw [dashed,black!50] (4.1,4.1) -- (5.9,4.1) -- (5.9,5.9) -- (4.1,5.9) -- cycle;

\node at (1,5) [] {\footnotesize $x_2^2 x_0^2$};
\node at (5,3) [] {\footnotesize $x_2 x_1^2 x_0$};
\node at (7,3) [] {\footnotesize $x_1^3 x_0$};
\end{tikzpicture}
}
\\
\subfloat[][The Borel set corresponding to the curve defined in Example  \ref{ex:curveSection} (cf. Figure \ref{fig:curveSection_a} and Figure \ref{fig:curveMarinari_a}).]{\label{fig:pointsCurvesGotzmann_b}
\begin{tikzpicture}[>=latex,scale=0.4]
\filldraw [fill=red!5,draw=red,thick] (0,14) -- (-14,-7) -- (14,-7) -- cycle;
\draw (2,3) -- (0,4) -- (-2,3);
\draw (-2,3) -- (-4,2) -- (-2,1) -- (0,2) -- cycle;
\draw (2,3) -- (4,2) -- (2,1) -- (0,2) -- cycle;

\draw (0,2) -- (0,0);
\draw (-4,2) -- (-4,0);
\draw (4,2) -- (4,0);
\draw (-2,1) -- (-2,-1);
\draw (2,1) -- (2,-1);


\draw (-4,0) -- (-6,-1) -- (-4,-2) -- (-2,-1) -- cycle;
\draw (-6,-1) -- (-8,-2) -- (-6,-3) -- (-4,-2);
\draw (-8,-2) -- (-10,-3) -- (-8,-4) -- (-6,-3);

\draw (-4,-2) -- (-2,-3) -- (0,-2);
\draw (0,0) -- (-2,-1) -- (0,-2) -- (2,-1) -- cycle;
\draw (0,-2) -- (2,-3) -- (4,-2);
\draw (4,0) -- (2,-1) -- (4,-2) -- (6,-1) -- cycle;

\draw (-10,-3) -- (-10,-5);
\draw (-8,-4) -- (-8,-6);
\draw (-6,-3) -- (-6,-5);
\draw (-4,-2) -- (-4,-4);
\draw (-2,-3) -- (-2,-5);
\draw (0,-2) -- (0,-4);
\draw (2,-3) -- (2,-5);
\draw (4,-2) -- (4,-4);
\draw (6,-1) -- (6,-3);

\draw (-10,-5) -- (-8,-6) -- (-4,-4) -- (-2,-5) -- (0,-4) -- (2,-5) -- (6,-3);
\draw [->,thick] (-10,-5) -- (-14,-7) --node[above left]{\footnotesize $\down{3}$} (-15,-7.5);
\draw [->,thick] (6,-3) -- (14,-7) --node[above right]{\footnotesize $\down{2}\circ\down{3}$} (15,-7.5);
\draw [->,thick] (0,4) -- (0,14) --node[right]{\footnotesize $\down{1}\circ\down{2}\circ\down{3}$} (0,15);

\node at (0,3) [] {\tiny $x_3^3 x_0$};
\node at (-2,2) [] {\tiny $x_3^2 x_2 x_0$};
\node at (2,2) [] {\tiny $x_3^2 x_1 x_0$};
\node at (-6,-2) [] {\tiny $x_3 x_2^3$};
\node at (-8,-3) [] {\tiny $x_2^4$};
\node at (-2,-2) [] {\tiny $x_3 x_2^2 x_1$};
\node at (2,-2) [] {\tiny $x_3 x_2 x_1^2$};

\node at (-8,-4) [circle,thick,draw=red,inner sep=1pt] {};
\end{tikzpicture}
}
\caption{\label{fig:pointsCurvesGotzmann} An example of Gotzmann's pyramids.}
\end{center}
\end{figure}

With this technique, the monomials not belonging to the Borel set are not drawn. To understand which monomial are missing, i.e. to have information about the subscheme defined, we can consider for $\pos{2}{m}$ the line touching the monomials $x_1^m$ and $x_0^m$, for $\pos{3}{m}$ the plane touching the monomials $x_2^m$, $x_1^m$ and $x_0^m$ and look at  the \lq\lq bricks\rq\rq\ that we would need to fill the empty space between the pyramid and the line or plane (see Figure \ref{fig:pointsCurvesGotzmann}).

\index{Gotzmann's pyramid|)}

\begin{figure}[H]
\begin{center}
\captionsetup[subfloat]{singlelinecheck=false}
\subfloat[][The Borel set defined by $(x_2^2,x_2x_1^2,x_1^3)$ in $\pos{2}{4}$ (cf. Figure \ref{fig:exampleMarinariPoints_a}, Figure \ref{fig:exampleGreenPoints_a} and Figure \ref{fig:pointsCurvesGotzmann_a}).]{
\begin{tikzpicture}[>=latex,line join=bevel,scale=0.7]
\tikzstyle{ideal}=[ellipse,thick,draw=black,inner sep=1.5pt]
\tikzstyle{quotient}=[inner sep=1.5pt]
\node (N_8) at (40bp,166bp) [ideal] {\footnotesize $x_{1}^{3}x_0$};
  \node (N_9) at (131bp,217bp) [ideal] {\footnotesize $x_2^{2}x_0^{2}$};
  \node (N_4) at (12bp,217bp) [ideal] {\footnotesize $x_{1}^{4}$};
  \node (N_5) at (103bp,319bp) [ideal] {\footnotesize $x_2^{3}x_0$};
  \node (N_6) at (103bp,268bp) [ideal] {\footnotesize $x_2^{2}x_{1}x_0$};
  \node (N_7) at (68bp,217bp) [ideal] {\footnotesize $x_2x_{1}^{2}x_0$};
  \node (N_0) at (75bp,421bp) [ideal] {\footnotesize $x_2^{4}$};
  \node (N_1) at (75bp,370bp) [ideal] {\footnotesize $x_2^{3}x_{1}$};
  \node (N_2) at (47bp,319bp) [ideal] {\footnotesize $x_2^{2}x_{1}^{2}$};
  \node (N_3) at (40bp,268bp) [ideal] {\footnotesize $x_2x_{1}^{3}$};
  \node (N_14) at (75bp,13bp) [quotient] {\footnotesize $x_0^{4}$};
  \node (N_12) at (103bp,115bp) [quotient] {\footnotesize $x_2x_0^{3}$};
  \node (N_13) at (75bp,64bp) [quotient] {\footnotesize $x_{1}x_0^{3}$};
  \node (N_10) at (103bp,166bp) [quotient] {\footnotesize $x_2x_{1}x_0^{2}$};
  \node (N_11) at (47bp,115bp) [quotient] {\footnotesize $x_{1}^{2}x_0^{2}$};
  \draw [->,thin,black!75] (N_9) --  (N_10);
  \draw [->,thin,black!75] (N_6) --  (N_9);
  \draw [->,thin,black!75] (N_2) --  (N_6);
  \draw [->,thin,black!75] (N_6) --  (N_7);
  \draw [->,thin,black!75] (N_7) --  (N_10);
  \draw [->,thin,black!75] (N_3) --  (N_7);
  \draw [->,thin,black!75] (N_0) --  (N_1);
  \draw [->,thin,black!75] (N_11) --  (N_13);
  \draw [->,thin,black!75] (N_5) --  (N_6);
  \draw [->,thin,black!75] (N_1) --  (N_2);
  \draw [->,thin,black!75] (N_10) --  (N_12);
  \draw [->,thin,black!75] (N_2) --  (N_3);
  \draw [->,thin,black!75] (N_10) --  (N_11);
  \draw [->,thin,black!75] (N_1) --  (N_5);
  \draw [->,thin,black!75] (N_7) --  (N_8);
  \draw [->,thin,black!75] (N_3) --  (N_4);
  \draw [->,thin,black!75] (N_8) --  (N_11);
  \draw [->,thin,black!75] (N_12) --  (N_13);
  \draw [->,thin,black!75] (N_13) --  (N_14);
  \draw [->,thin,black!75] (N_4) --  (N_8);
\end{tikzpicture}
}
\qquad\qquad
\subfloat[][The Borel set corresponding to the curve defined in Example  \ref{ex:curveSection} (cf. Figure \ref{fig:curveSection_a}, Figure \ref{fig:curveMarinari_a} and Figure \ref{fig:pointsCurvesGotzmann_b}).]{
\begin{tikzpicture}[>=latex,line join=bevel,scale=0.6]
\tikzstyle{ideal}=[ellipse,thick,draw=black,inner sep=1.5pt]
\tikzstyle{quotient}=[inner sep=1.5pt]

  \node (N_23) at (129bp,217bp) [quotient] {\footnotesize $x_2x_1^{2}x_0$};
  \node (N_22) at (201bp,268bp) [quotient] {\footnotesize $x_3x_1^{2}x_0$};
  \node (N_21) at (131bp,268bp) [quotient] {\footnotesize $x_2^{2}x_1x_0$};
  \node (N_20) at (201bp,319bp) [quotient] {\footnotesize $x_3x_2x_1x_0$};
  \node (N_27) at (192bp,217bp) [quotient] {\footnotesize $x_2^{2}x_0^{2}$};
  \node (N_26) at (271bp,268bp) [quotient] {\footnotesize $x_3x_2x_0^{2}$};
  \node (N_25) at (271bp,319bp) [quotient] {\footnotesize $x_3^{2}x_0^{2}$};
  \node (N_24) at (129bp,166bp) [quotient] {\footnotesize $x_1^{3}x_0$};
  \node (N_29) at (192bp,166bp) [quotient] {\footnotesize $x_2x_1x_0^{2}$};
  \node (N_28) at (255bp,217bp) [quotient] {\footnotesize $x_3x_1x_0^{2}$};
  \node (N_8) at (19bp,370bp) [quotient] {\footnotesize $x_2^{3}x_1$};
  \node (N_9) at (143bp,421bp) [ideal] {\footnotesize $x_3^{2}x_1^{2}$};
  \node (N_4) at (21bp,421bp) [ideal] {\footnotesize $x_2^{4}$};
  \node (N_5) at (171bp,523bp) [ideal] {\footnotesize $x_3^{3}x_1$};
  \node (N_6) at (143bp,472bp) [ideal] {\footnotesize $x_3^{2}x_2x_1$};
  \node (N_7) at (80bp,421bp) [ideal] {\footnotesize $x_3x_2^{2}x_1$};
  \node (N_0) at (143bp,625bp) [ideal] {\footnotesize $x_3^{4}$};
  \node (N_1) at (143bp,574bp) [ideal] {\footnotesize $x_3^{3}x_2$};
  \node (N_2) at (115bp,523bp) [ideal] {\footnotesize $x_3^{2}x_2^{2}$};
  \node (N_3) at (80bp,472bp) [ideal] {\footnotesize $x_3x_2^{3}$};
  \node (N_30) at (164bp,115bp) [quotient] {\footnotesize $x_1^{2}x_0^{2}$};
  \node (N_31) at (255bp,166bp) [quotient] {\footnotesize $x_3x_0^{3}$};
  \node (N_32) at (220bp,115bp) [quotient] {\footnotesize $x_2x_0^{3}$};
  \node (N_33) at (192bp,64bp) [quotient] {\footnotesize $x_1x_0^{3}$};
  \node (N_34) at (192bp,13bp) [quotient] {\footnotesize $x_0^{4}$};
  \node (N_18) at (131bp,319bp) [quotient] {\footnotesize $x_2^{3}x_0$};
  \node (N_19) at (222bp,370bp) [ideal] {\footnotesize $x_3^{2}x_1x_0$};
  \node (N_16) at (206bp,421bp) [ideal] {\footnotesize $x_3^{2}x_2x_0$};
  \node (N_17) at (152bp,370bp) [quotient] {\footnotesize $x_3x_2^{2}x_0$};
  \node (N_14) at (71bp,217bp) [quotient] {\footnotesize $x_1^{4}$};
  \node (N_15) at (206bp,472bp) [ideal] {\footnotesize $x_3^{3}x_0$};
  \node (N_12) at (75bp,319bp) [quotient] {\footnotesize $x_3x_1^{3}$};
  \node (N_13) at (68bp,268bp) [quotient] {\footnotesize $x_2x_1^{3}$};
  \node (N_10) at (82bp,370bp) [ideal] {\footnotesize $x_3x_2x_1^{2}$};
  \node (N_11) at (19bp,319bp) [quotient] {\footnotesize $x_2^{2}x_1^{2}$};
  \draw [->,thin,black!75] (N_2) --  (N_6);
  \draw [->,thin,black!75] (N_6) --  (N_7);
  \draw [->,thin,black!75] (N_14) --  (N_24);
  \draw [->,thin,black!75] (N_23) --  (N_24);
  \draw [->,thin,black!75] (N_11) --  (N_13);
  \draw [->,thin,black!75] (N_2) --  (N_3);
  \draw [->,thin,black!75] (N_23) --  (N_29);
  \draw [->,thin,black!75] (N_29) --  (N_30);
  \draw [->,thin,black!75] (N_24) --  (N_30);
  \draw [->,thin,black!75] (N_30) --  (N_33);
  \draw [->,thin,black!75] (N_22) --  (N_28);
  \draw [->,thin,black!75] (N_22) --  (N_23);
  \draw [->,thin,black!75] (N_31) --  (N_32);
  \draw [->,thin,black!75] (N_19) --  (N_20);
  \draw [->,thin,black!75] (N_20) --  (N_21);
  \draw [->,thin,black!75] (N_32) --  (N_33);
  \draw [->,thin,black!75] (N_11) --  (N_21);
  \draw [->,thin,black!75] (N_10) --  (N_12);
  \draw [->,thin,black!75] (N_0) --  (N_1);
  \draw [->,thin,black!75] (N_17) --  (N_20);
  \draw [->,thin,black!75] (N_21) --  (N_23);
  \draw [->,thin,black!75] (N_16) --  (N_17);
  \draw [->,thin,black!75] (N_15) --  (N_16);
  \draw [->,thin,black!75] (N_18) --  (N_21);
  \draw [->,thin,black!75] (N_13) --  (N_23);
  \draw [->,thin,black!75] (N_10) --  (N_20);
  \draw [->,thin,black!75] (N_25) --  (N_26);
  \draw [->,thin,black!75] (N_20) --  (N_26);
  \draw [->,thin,black!75] (N_28) --  (N_29);
  \draw [->,thin,black!75] (N_7) --  (N_8);
  \draw [->,thin,black!75] (N_7) --  (N_10);
  \draw [->,thin,black!75] (N_27) --  (N_29);
  \draw [->,thin,black!75] (N_12) --  (N_13);
  \draw [->,thin,black!75] (N_19) --  (N_25);
  \draw [->,thin,black!75] (N_20) --  (N_22);
  \draw [->,thin,black!75] (N_26) --  (N_27);
  \draw [->,thin,black!75] (N_5) --  (N_6);
  \draw [->,thin,black!75] (N_6) --  (N_16);
  \draw [->,thin,black!75] (N_10) --  (N_11);
  \draw [->,thin,black!75] (N_8) --  (N_18);
  \draw [->,thin,black!75] (N_8) --  (N_11);
  \draw [->,thin,black!75] (N_9) --  (N_19);
  \draw [->,thin,black!75] (N_33) --  (N_34);
  \draw [->,thin,black!75] (N_3) --  (N_4);
  \draw [->,thin,black!75] (N_29) --  (N_32);
  \draw [->,thin,black!75] (N_9) --  (N_10);
  \draw [->,thin,black!75] (N_6) --  (N_9);
  \draw [->,thin,black!75] (N_5) --  (N_15);
  \draw [->,thin,black!75] (N_12) --  (N_22);
  \draw [->,thin,black!75] (N_1) --  (N_2);
  \draw [->,thin,black!75] (N_17) --  (N_18);
  \draw [->,thin,black!75] (N_26) --  (N_28);
  \draw [->,thin,black!75] (N_3) --  (N_7);
  \draw [->,thin,black!75] (N_28) --  (N_31);
  \draw [->,thin,black!75] (N_16) --  (N_19);
  \draw [->,thin,black!75] (N_1) --  (N_5);
  \draw [->,thin,black!75] (N_7) --  (N_17);
  \draw [->,thin,black!75] (N_13) --  (N_14);
  \draw [->,thin,black!75] (N_21) --  (N_27);
  \draw [->,thin,black!75] (N_4) --  (N_8);
\end{tikzpicture}
}
\caption[Borel sets drawn as planar graphs.]{\label{fig:graphPointsCurves} An example of Borel sets drawn as planar graphs.}
\end{center}
\end{figure}
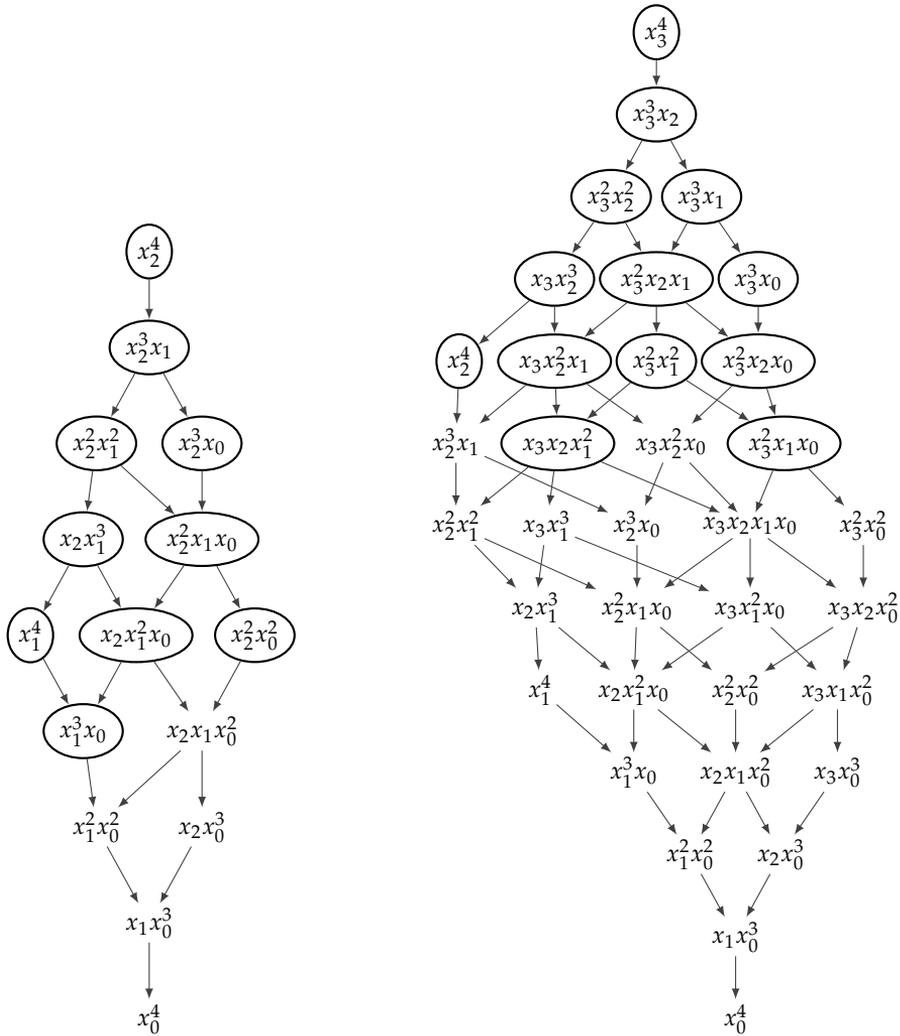

\paragraph{Planar graphs}\index{Borel set!planar graph representing a|(}\index{planar graph|see{Borel set}}
All previous approches result to be intrinsically limitated to the projective $3$-space. We would like to overcome this limit and to find a nice representations available for any kind of poset. A very natural way to describe $\pos{n}{m}$ comes directly from the definition of partially ordered set: we associate to it the graph whose vertices are the monomials of $\K[x_0,\ldots,x_n]_m$ and whose edges correspond to elementary decreasing moves (see Figure \ref{fig:graphPointsCurves}). Given a Borel set $\mathscr{B} \subset \pos{n}{m}$, we will represent its monomials with vertices with elliptic black boundary and without boundaries the monomials outside $\mathscr{B}$.

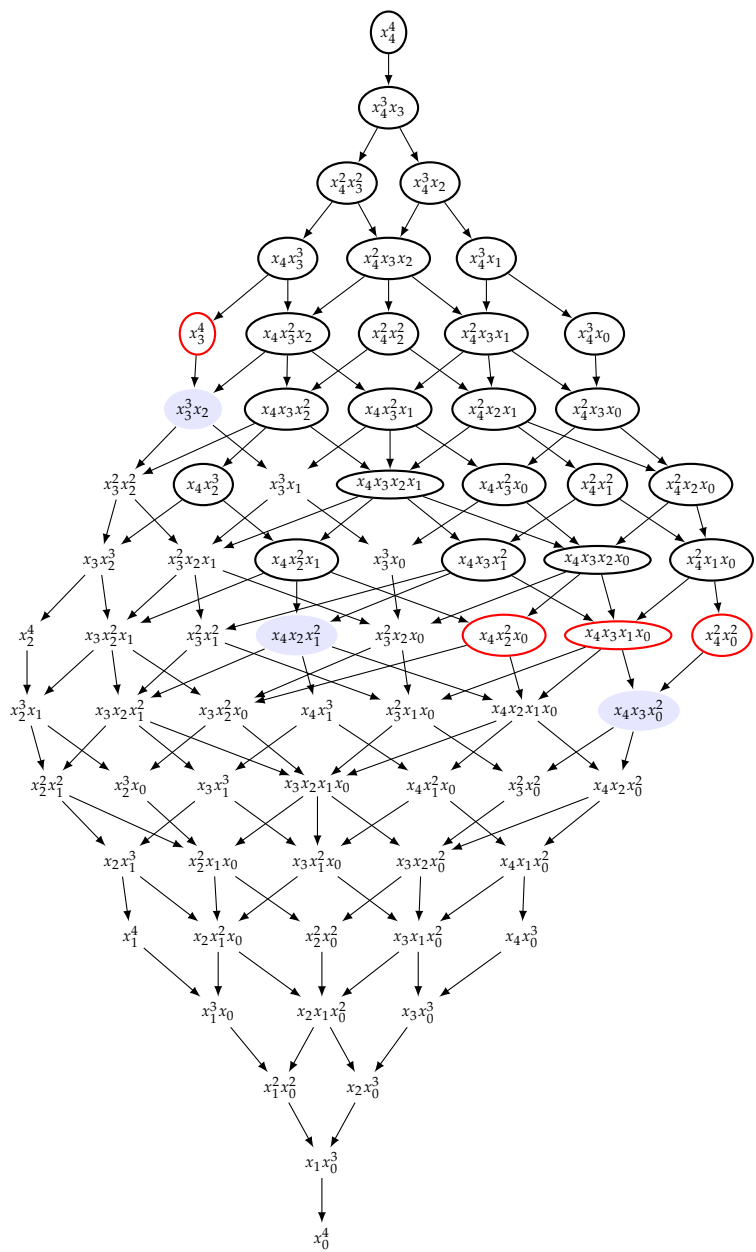
\begin{figure}[H]
\begin{center}
\begin{tikzpicture}[>=latex,line join=bevel,scale=0.55]
\tikzstyle{ideal}=[ellipse,thick,draw=black,inner sep=1.5pt]
\tikzstyle{quotient}=[inner sep=1.5pt]

\tikzstyle{min}=[ellipse,thick,draw=red,inner sep=1.5pt]
\tikzstyle{max}=[ellipse,draw=white,fill=blue!10,inner sep=1.5pt]

\node (N_58) at (418bp,319bp) [quotient] {\tiny $x_4x_{2}x_0^{2}$};
  \node (N_59) at (285bp,268bp) [quotient] {\tiny $x_3x_{2}x_0^{2}$};
  \node (N_52) at (215bp,268bp) [quotient] {\tiny $x_3x_1^{2}x_0$};
  \node (N_53) at (148bp,217bp) [quotient] {\tiny $x_{2}x_1^{2}x_0$};
  \node (N_50) at (145bp,268bp) [quotient] {\tiny $x_{2}^{2}x_1x_0$};
  \node (N_51) at (292bp,319bp) [quotient] {\tiny $x_4x_1^{2}x_0$};
  \node (N_56) at (432bp,370bp) [max] {\tiny $x_4x_3x_0^{2}$};
  \node (N_57) at (355bp,319bp) [quotient] {\tiny $x_3^{2}x_0^{2}$};
  \node (N_54) at (148bp,166bp) [quotient] {\tiny $x_1^{3}x_0$};
  \node (N_55) at (488bp,421bp) [min] {\tiny $x_4^{2}x_0^{2}$};
  \node (N_46) at (418bp,421bp) [min] {\tiny $x_4x_3x_1x_0$};
  \node (N_45) at (481bp,472bp) [ideal] {\tiny $x_4^{2}x_1x_0$};
  \node (N_44) at (89bp,319bp) [quotient] {\tiny $x_{2}^{3}x_0$};
  \node (N_47) at (278bp,370bp) [quotient] {\tiny $x_3^{2}x_1x_0$};
  \node (N_42) at (341bp,421bp) [min] {\tiny $x_4x_{2}^{2}x_0$};
  \node (N_23) at (75bp,421bp) [quotient] {\tiny $x_3x_{2}^{2}x_1$};
  \node (N_22) at (201bp,472bp) [ideal] {\tiny $x_4x_{2}^{2}x_1$};
  \node (N_21) at (131bp,472bp) [quotient] {\tiny $x_3^{2}x_{2}x_1$};
  \node (N_20) at (264bp,523bp) [ideal] {\tiny $x_4x_3x_{2}x_1$};
  \node (N_27) at (138bp,421bp) [quotient] {\tiny $x_3^{2}x_1^{2}$};
  \node (N_26) at (327bp,472bp) [ideal] {\tiny $x_4x_3x_1^{2}$};
  \node (N_25) at (404bp,523bp) [ideal] {\tiny $x_4^{2}x_1^{2}$};
  \node (N_24) at (19bp,370bp) [quotient] {\tiny $x_{2}^{3}x_1$};
  \node (N_41) at (271bp,421bp) [quotient] {\tiny $x_3^{2}x_{2}x_0$};
  \node (N_29) at (82bp,370bp) [quotient] {\tiny $x_3x_{2}x_1^{2}$};
  \node (N_28) at (201bp,421bp) [max] {\tiny $x_4x_{2}x_1^{2}$};
  \node (N_49) at (215bp,319bp) [quotient] {\tiny $x_3x_{2}x_1x_0$};
  \node (N_48) at (355bp,370bp) [quotient] {\tiny $x_4x_{2}x_1x_0$};
  \node (N_40) at (404bp,472bp) [ideal] {\tiny $x_4x_3x_{2}x_0$};
  \node (N_64) at (190bp,115bp) [quotient] {\tiny $x_1^{2}x_0^{2}$};
  \node (N_36) at (404bp,574bp) [ideal] {\tiny $x_4^{2}x_3x_0$};
  \node (N_69) at (218bp,13bp) [quotient] {\tiny $x_0^{4}$};
  \node (N_68) at (218bp,64bp) [quotient] {\tiny $x_1x_0^{3}$};
  \node (N_67) at (246bp,115bp) [quotient] {\tiny $x_{2}x_0^{3}$};
  \node (N_66) at (283bp,166bp) [quotient] {\tiny $x_3x_0^{3}$};
  \node (N_65) at (353bp,217bp) [quotient] {\tiny $x_4x_0^{3}$};
  \node (N_37) at (341bp,523bp) [ideal] {\tiny $x_4x_3^{2}x_0$};
  \node (N_63) at (218bp,166bp) [quotient] {\tiny $x_{2}x_1x_0^{2}$};
  \node (N_62) at (283bp,217bp) [quotient] {\tiny $x_3x_1x_0^{2}$};
  \node (N_61) at (355bp,268bp) [quotient] {\tiny $x_4x_1x_0^{2}$};
  \node (N_60) at (218bp,217bp) [quotient] {\tiny $x_{2}^{2}x_0^{2}$};
  \node (N_38) at (264bp,472bp) [quotient] {\tiny $x_3^{3}x_0$};
  \node (N_39) at (467bp,523bp) [ideal] {\tiny $x_4^{2}x_{2}x_0$};
  \node (N_43) at (152bp,370bp) [quotient] {\tiny $x_3x_{2}^{2}x_0$};
  \node (N_8) at (131bp,574bp) [max] {\tiny $x_3^{3}x_{2}$};
  \node (N_9) at (263bp,625bp) [ideal] {\tiny $x_4^{2}x_{2}^{2}$};
  \node (N_4) at (134bp,625bp) [min] {\tiny $x_3^{4}$};
  \node (N_5) at (291bp,727bp) [ideal] {\tiny $x_4^{3}x_{2}$};
  \node (N_6) at (263bp,676bp) [ideal] {\tiny $x_4^{2}x_3x_{2}$};
  \node (N_7) at (195bp,625bp) [ideal] {\tiny $x_4x_3^{2}x_{2}$};
  \node (N_0) at (263bp,829bp) [ideal] {\tiny $x_4^{4}$};
  \node (N_1) at (263bp,778bp) [ideal] {\tiny $x_4^{3}x_3$};
  \node (N_2) at (235bp,727bp) [ideal] {\tiny $x_4^{2}x_3^{2}$};
  \node (N_3) at (195bp,676bp) [ideal] {\tiny $x_4x_3^{3}$};
  \node (N_30) at (33bp,319bp) [quotient] {\tiny $x_{2}^{2}x_1^{2}$};
  \node (N_31) at (215bp,370bp) [quotient] {\tiny $x_4x_1^{3}$};
  \node (N_32) at (145bp,319bp) [quotient] {\tiny $x_3x_1^{3}$};
  \node (N_33) at (82bp,268bp) [quotient] {\tiny $x_{2}x_1^{3}$};
  \node (N_34) at (89bp,217bp) [quotient] {\tiny $x_1^{4}$};
  \node (N_35) at (402bp,625bp) [ideal] {\tiny $x_4^{3}x_0$};
  \node (N_18) at (194bp,523bp) [quotient] {\tiny $x_3^{3}x_1$};
  \node (N_19) at (334bp,574bp) [ideal] {\tiny $x_4^{2}x_{2}x_1$};
  \node (N_16) at (329bp,625bp) [ideal] {\tiny $x_4^{2}x_3x_1$};
  \node (N_17) at (264bp,574bp) [ideal] {\tiny $x_4x_3^{2}x_1$};
  \node (N_14) at (19bp,421bp) [quotient] {\tiny $x_{2}^{4}$};
  \node (N_15) at (329bp,676bp) [ideal] {\tiny $x_4^{3}x_1$};
  \node (N_12) at (138bp,523bp) [ideal] {\tiny $x_4x_{2}^{3}$};
  \node (N_13) at (68bp,472bp) [quotient] {\tiny $x_3x_{2}^{3}$};
  \node (N_10) at (194bp,574bp) [ideal] {\tiny $x_4x_3x_{2}^{2}$};
  \node (N_11) at (82bp,523bp) [quotient] {\tiny $x_3^{2}x_{2}^{2}$};
  \draw [->] (N_47) -- (N_57);
  \draw [->] (N_42) -- (N_43);
  \draw [->] (N_61) -- (N_62);
  \draw [->] (N_60) -- (N_63);
  \draw [->] (N_11) -- (N_13);
  \draw [->] (N_2) -- (N_3);
  \draw [->] (N_19) -- (N_20);
  \draw [->] (N_20) -- (N_21);
  \draw [->] (N_26) -- (N_28);
  \draw [->] (N_68) -- (N_69);
  \draw [->] (N_66) -- (N_67);
  \draw [->] (N_40) -- (N_46);
  \draw [->] (N_53) -- (N_54);
  \draw [->] (N_19) -- (N_25);
  \draw [->] (N_20) -- (N_22);
  \draw [->] (N_29) -- (N_49);
  \draw [->] (N_50) -- (N_60);
  \draw [->] (N_32) -- (N_33);
  \draw [->] (N_8) -- (N_18);
  \draw [->] (N_45) -- (N_55);
  \draw [->] (N_39) -- (N_40);
  \draw [->] (N_22) -- (N_23);
  \draw [->] (N_33) -- (N_34);
  \draw [->] (N_17) -- (N_37);
  \draw [->] (N_9) -- (N_10);
  \draw [->] (N_6) -- (N_9);
  \draw [->] (N_49) -- (N_59);
  \draw [->] (N_1) -- (N_2);
  \draw [->] (N_37) -- (N_40);
  \draw [->] (N_48) -- (N_58);
  \draw [->] (N_22) -- (N_42);
  \draw [->] (N_32) -- (N_52);
  \draw [->] (N_26) -- (N_46);
  \draw [->] (N_13) -- (N_14);
  \draw [->] (N_21) -- (N_41);
  \draw [->] (N_39) -- (N_45);
  \draw [->] (N_6) -- (N_7);
  \draw [->] (N_59) -- (N_62);
  \draw [->] (N_52) -- (N_53);
  \draw [->] (N_28) -- (N_31);
  \draw [->] (N_29) -- (N_30);
  \draw [->] (N_58) -- (N_61);
  \draw [->] (N_3) -- (N_7);
  \draw [->] (N_56) -- (N_58);
  \draw [->] (N_26) -- (N_27);
  \draw [->] (N_21) -- (N_23);
  \draw [->] (N_0) -- (N_1);
  \draw [->] (N_10) -- (N_20);
  \draw [->] (N_43) -- (N_49);
  \draw [->] (N_50) -- (N_53);
  \draw [->] (N_20) -- (N_26);
  \draw [->] (N_54) -- (N_64);
  \draw [->] (N_33) -- (N_53);
  \draw [->] (N_41) -- (N_43);
  \draw [->] (N_15) -- (N_35);
  \draw [->] (N_46) -- (N_47);
  \draw [->] (N_6) -- (N_16);
  \draw [->] (N_3) -- (N_4);
  \draw [->] (N_67) -- (N_68);
  \draw [->] (N_1) -- (N_5);
  \draw [->] (N_16) -- (N_36);
  \draw [->] (N_9) -- (N_19);
  \draw [->] (N_28) -- (N_48);
  \draw [->] (N_49) -- (N_50);
  \draw [->] (N_48) -- (N_51);
  \draw [->] (N_7) -- (N_17);
  \draw [->] (N_21) -- (N_27);
  \draw [->] (N_2) -- (N_6);
  \draw [->] (N_24) -- (N_30);
  \draw [->] (N_63) -- (N_64);
  \draw [->] (N_65) -- (N_66);
  \draw [->] (N_47) -- (N_49);
  \draw [->] (N_45) -- (N_46);
  \draw [->] (N_31) -- (N_32);
  \draw [->] (N_62) -- (N_63);
  \draw [->] (N_24) -- (N_44);
  \draw [->] (N_16) -- (N_17);
  \draw [->] (N_15) -- (N_16);
  \draw [->] (N_18) -- (N_21);
  \draw [->] (N_41) -- (N_47);
  \draw [->] (N_34) -- (N_54);
  \draw [->] (N_25) -- (N_26);
  \draw [->] (N_19) -- (N_39);
  \draw [->] (N_48) -- (N_49);
  \draw [->] (N_7) -- (N_8);
  \draw [->] (N_18) -- (N_38);
  \draw [->] (N_7) -- (N_10);
  \draw [->] (N_12) -- (N_13);
  \draw [->] (N_42) -- (N_48);
  \draw [->] (N_5) -- (N_6);
  \draw [->] (N_51) -- (N_52);
  \draw [->] (N_10) -- (N_11);
  \draw [->] (N_27) -- (N_29);
  \draw [->] (N_37) -- (N_38);
  \draw [->] (N_29) -- (N_32);
  \draw [->] (N_64) -- (N_68);
  \draw [->] (N_36) -- (N_39);
  \draw [->] (N_5) -- (N_15);
  \draw [->] (N_12) -- (N_22);
  \draw [->] (N_35) -- (N_36);
  \draw [->] (N_23) -- (N_43);
  \draw [->] (N_53) -- (N_63);
  \draw [->] (N_20) -- (N_40);
  \draw [->] (N_17) -- (N_18);
  \draw [->] (N_4) -- (N_8);
  \draw [->] (N_36) -- (N_37);
  \draw [->] (N_30) -- (N_50);
  \draw [->] (N_14) -- (N_24);
  \draw [->] (N_31) -- (N_51);
  \draw [->] (N_23) -- (N_24);
  \draw [->] (N_46) -- (N_56);
  \draw [->] (N_23) -- (N_29);
  \draw [->] (N_52) -- (N_62);
  \draw [->] (N_30) -- (N_33);
  \draw [->] (N_22) -- (N_28);
  \draw [->] (N_40) -- (N_41);
  \draw [->] (N_11) -- (N_21);
  \draw [->] (N_46) -- (N_48);
  \draw [->] (N_17) -- (N_20);
  \draw [->] (N_25) -- (N_45);
  \draw [->] (N_49) -- (N_52);
  \draw [->] (N_27) -- (N_47);
  \draw [->] (N_59) -- (N_60);
  \draw [->] (N_28) -- (N_29);
  \draw [->] (N_61) -- (N_65);
  \draw [->] (N_43) -- (N_44);
  \draw [->] (N_40) -- (N_42);
  \draw [->] (N_55) -- (N_56);
  \draw [->] (N_51) -- (N_61);
  \draw [->] (N_62) -- (N_66);
  \draw [->] (N_8) -- (N_11);
  \draw [->] (N_56) -- (N_57);
  \draw [->] (N_13) -- (N_23);
  \draw [->] (N_38) -- (N_41);
  \draw [->] (N_44) -- (N_50);
  \draw [->] (N_57) -- (N_59);
  \draw [->] (N_10) -- (N_12);
  \draw [->] (N_16) -- (N_19);
  \draw [->] (N_58) -- (N_59);
  \draw [->] (N_63) -- (N_67);
\end{tikzpicture}
\caption[An example of a Borel set drawn as a planar graph, in which the minimal and maximal elements are highlighted.]{\label{fig:minMaxMonomials} The Borel set defined by the ideal $(x_4^2,x_4x_3^2,x_4x_3x_2,x_4x_2^2,x_4x_3x_1,x_3^4)$ in $\pos{4}{4}$. The monomials with red boundary are minimal elements of the Borel set, whereas the monomials with a light blue background are maximal elements in the complement.}
\end{center}
\end{figure}

This representation, as well as allowing to manage any poset, is very advantageous in an algorithmic perspective, indeed there are many tools to work on graphs. We underline that with this description is very easy to detect minimal and maximal monomials (see Figure \ref{fig:minMaxMonomials}).

In the following we will use Marinari's lattices and planar graphs: the pictures of lattices turn out to be very helpful to understand the main ideas and planar graphs allow to generalize such ideas to posets in any number of variables and to concretely project algorithms, indeed the java class \texttt{PosetGraph} of the package \texttt{HSC} implements the poset by means of the associated direct graph (see Appendix \ref{ch:HSCpackage}).

\newpage

\index{Borel set!planar graph representing a|)}

\section{An algorithm computing Borel-fixed ideals}

In this section and in Section \ref{sec:segments}, we will expose some of the results contained in the paper \cite{CLMR} \lq\lq Segments and Hilbert schemes of points\rq\rq\ written in collaboration with Francesca Cioffi, Maria Grazia Marinari and Margherita Roggero.

\medskip 

Let $J\subset \K[x]$ be a Borel-fixed ideal. In this section we denote by $J_{x_0}$ the ideal obtained from $J$ setting $x_0 = 1$. Keeping in mind Corollary \ref{cor:saturationBorelIdeal}, we know that $J_{x_0}=J^{\sat}$. We extend this notation denoting by $J_{x_0x_1}$ the ideal obtained from $J$ setting both $x_0$ and $x_1$ equal to 1. We call $J_{x_0x_1}$ the $x_1$\emph{-saturation} of $J$ and say that $J$ is $x_1$\emph{-saturated} if $J=J_{x_0x_1}$. Hence an ideal $J$ that is $x_1$-saturated is also saturated.

\begin{remark}
A $x_1$-saturated Borel-fixed ideal $J\subset \K[x]$ defining a subscheme with Hilbert polynomial $p(t)$ has the same minimal generators as the saturated Borel ideal $J\cap \K[x_1,\ldots,x_n]\subset \K[x_1,\ldots,x_n]$, that defines a subscheme of $\PP^{n-1}$ with Hilbert polynomial $\Delta p(t)$. 
\end{remark}

The following result is analogous to Theorem 3 of \cite{Reeves}, where the notion of \lq\lq fan\rq\rq\ is used. 
Here we apply the combinatorial properties of Borel ideals only.

\begin{proposition}\label{prop:difference}
Let $J\subset \K[x]$ be a saturated Borel-fixed ideal defining a subscheme with Hilbert polynomial $p(t)$ whose Gotzmann number is $r$. Let $I=J_{x_0x_1}$ be its $x_1$-saturation and let $\overline{p}(t)$ be the Hilbert polynomial of the subscheme defined by $I$ in $\PP^n$. Set $q =\dim_\K I_r- \dim_\K J_r$,
\begin{enumerate}[(i)]
\item\label{it:difference_i} $\overline{p}(t)=p(t)-q$;
\item\label{it:difference_ii} $q$ is equal to the sum of the exponents of $x_1$ in the minimal generators of $J$.
\end{enumerate}
\end{proposition}

\begin{proof}
\emph{(\ref{it:difference_i})} We show that if $q=\dim_\K I_s-\dim_\K J_s$ then $q=\dim_\K I_{s+1}-\dim_\K J_{s+1}$, for every $s\geq r$. 
Let $x^{\beta_1},\ldots,x^{\beta_q}$ be the terms of $I_s\setminus J_s$. Thus, $x_0x^{\beta_1},\ldots,x_0x^{\beta_q}$ are terms of $I_{s+1}\setminus J_{s+1}$ and so $\dim_\K I_{s+1}-\dim_\K J_{s+1}\geqslant q$, since $x_0x^{\beta_i}$ would belong to $J_{s+1}$ if and only if $x^{\beta_i}$ belong to $J_{s}$, being $J$ saturated. Now, to obtain the opposite inequality, it is enough to show that every term of $I_{s+1}\setminus J_{s+1}$ is divisible by $x_0$. Let $x^{\gamma}\in I_{s+1}$ be such that $\min x^{\gamma}\geqslant 1$ and let $x^{\alpha}$ be a minimal generator of $I$ such 
that $x^{\gamma}=x^{\alpha}x^{\delta}$. Since $J$ is saturated and $I$ is the $x_1$-saturation of $J$, 
$x^{\alpha}x_1^a$ is a minimal generator of $J$ for some non negative integer $a$. Hence, for every $x^{\delta'}$ of degree $s+1-\vert\alpha\vert$ and with $\min x^{\delta'} \geqslant 1$, $x^{\delta'} >_B x_1^{a}$ implies $x^{\alpha}x^{\delta'} \in J_{s+1}$.  In particular, $x^{\gamma}\in J_{s+1}$.

\emph{(\ref{it:difference_ii})} Let $x^{\alpha_1}x_1^{s_1}, \ldots, x^{\alpha_h}x_1^{s_h}$ be the minimal generators of $J$, with $\min x^{\alpha_i} > 1$, for every $1\leq i\leq h$. Since the $\sum s_i$ terms $x^{\alpha_i}x_1^{s_i-t}x_0^{r-\vert\alpha_i\vert-s_i+t}, 1\leq t\leq s_i$, are in $I_r\setminus J_r$, one has $q\geq \sum s_i$. Vice versa, we show that each term $x^{\delta}$ in $I_r\setminus J_r$ is of the previous type. We can write $x^{\delta}=x^{\beta}x_1^u x_0^{r-\vert\beta\vert-u}$, with $\min x^{\beta}> 1$ and $u<s_i$. Let $s$ be the minimum non negative integer such that $x^{\beta} x_1^s$ is in $J$. Then there exists $i$ such that $x^{\alpha_i}x_1^{s_i} \vert x^{\beta}x_1^s$, i.e. $x^{\alpha_i} \vert x^{\beta}$ and $s_i\leq s$. By the definition of $s$, we get $s_i=s$ and there exists $x^{\gamma}$ with  $\min x^{\gamma} > 1$ such that $x^{\beta}=x^{\alpha_i}x^{\gamma}$. Since $x^{\beta}$ does not belong to $J$ we have $\vert\gamma\vert<s_i=s$, or otherwise $x^{\alpha_i}x_1^{\vert\gamma\vert}$ and hence, by the Borel property, $x^{\beta}=x^{\alpha_i}x^{\gamma}$ should belong to $J$. Now we can take $x^{\beta}x_1^{s-\vert\gamma\vert}$ and observe that this term belongs to $J$ because it follows $x^{\alpha_i}x_1^s$ in the Borel relation. Thus $s\leq s-\vert \gamma \vert$, so that 
$\gamma=0$, i.e. $x^{\beta}=x^{\alpha_i}$ as claimed.
\end{proof}

\begin{lemma}\label{lem:removingMinimal}
Let $J\subset \K[x]$ be a saturated Borel-fixed ideal such that $\K[x]/J$ has Hilbert polynomial $p(t)$ whose Gotzmann number is $r$. 
Let $x^{\beta}$ be a minimal monomial of $\{J_s\} \subset \pos{n}{s}$ of degree $s\geqslant r$ such that $\min x^{\beta} = x_0$. Then the ideal  $I =\langle \{J_s\}\setminus \{x^{\beta}\}\rangle$ is Borel-fixed and $\K[x]/I$ has Hilbert polynomial $\overline{p}(t) = p(t)+1$.
\end{lemma}
\begin{proof}
First, note that by definition of minimal monomial, $\{I_s\}$ is still a Borel set. Called $q(t)$ the volume polynomial of $J$, we show that $I$ has volume polynomial $\overline{q}(t) = q(t)-1$ applying Gotzmann's Persistence Theorem (Theorem \ref{th:PersistenceTheorem}), i.e. proving that $\dim_{\K} J_s - \dim_{\K} I_s = \dim_{\K} J_{s+1} - \dim_{\K} I_{s+1} = 1$.  
By construction $\dim_{\K} J_s - \dim_{\K} I_s = 1$. The Borel considition ensures that $x^\beta x_0 \in J_{s+1} \setminus I_{s+1}$ and there are no other elements, because $x^\beta x_0 $ is the only monomial that cannot be generated from the monomials in $I_{s}$ by multiplication of a single variable.  In fact let us consider the monomial $x_i x^\beta$, $i > 0$. Since $x_0 \mid x^\beta$ the following identity holds:
\[
x_i x^\beta = \dfrac{x_i}{x_{i-1}} \cdot \ldots \cdot \dfrac{x_1}{x_0} x_0 x^\beta = \up{i-1}\circ \cdots \circ \up{0}(x^\beta) x_0
\]
and for each $i$, $\up{i-1}\circ \cdots \circ \up{0}(x^\beta)$ belongs to $I_s$, by the minimality of $x^\beta$.
\end{proof}

\begin{proposition}\label{$x_1$-saturation}
Let $I$ and $J$ be Borel-fixed ideals of $\K[x]$. If for every $s\gg 0$ we have $I_s \subset J_s$ and 
$p_{\K[x]/I}(t)=p_{\K[x]/J}(t)+a$, with $a\in \mathbb N$, then $I$ and $J$ have the same $x_1$-saturation. 
\end{proposition}
\begin{proof}
Let $s\geq \max\{\reg(I),\reg(J)\}$. In case $a=1$, there exists a unique term in $J_{s+t}\setminus I_{s+t}$, 
for every $t\geq 0$. Let $x^{\alpha}$ be the unique term in $J_s\setminus I_s$. Then, both $x^{\alpha}x_0$ and 
$x^{\alpha}x_1$ belong to $J_{s+1}$. By the Borel property, $x^{\alpha}x_1$ must be in 
$I_{s+1}$ and so the unique term in $J_{s+t}\setminus I_{s+t}$ is $x^{\alpha}{x_0}^t$. This is enough to say that $I$ and $J$ have the same $x_1$-saturation. If $a>1$, the thesis follows 
by induction applying Lemma \ref{lem:removingMinimal}.
\end{proof}

\begin{theorem}\label{th:bijectionIdealsSets}
Let $p(t)$ be an admissible Hilbert polynomial in $\PP^n$. For any $s$, there is a bijective function
\begin{equation}\label{eq:bijectionIdealsSets}
\begin{array}{rcl}
\left\{ \begin{array}{c}J \subset \K[x] \text{ saturated Borel-fixed}\\ \text{ideal s.t. } \reg(J) \leqslant s \text{ and } \K[x]/J \text{ has}\\ \text{Hilbert polynomial } p(t)\end{array}\right\} & \stackrel{1:1}{\longleftrightarrow} & \left\{\begin{array}{c} \mathscr{B} \subset \pos{n}{s} \text{ Borel set s.t.}\\ \text{set } \mathscr{N} = \pos{n}{s}\setminus \mathscr{B}\\ \Big\vert \restrict{\mathscr{N}}{i} \Big\vert = \Delta^i p(s),\ \forall\ i\end{array}\right\}\\
J & \longrightarrow & \{J_s\}\\
\langle\mathscr{B}\rangle^{\sat} & \longleftarrow & \mathscr{B}
\end{array}
\end{equation}
\end{theorem}
\begin{proof}
First of all, note that if the two maps are well-defined, keeping in mind \ref{prop:regularityDegreeGenerators}, i.e. for each $J$, $J_{\geqslant s} = \langle J_s \rangle$, 
\begin{eqnarray*}
& J\ \longrightarrow\ \{J_s\}\ \longrightarrow\ \big\langle\{J_s\}\big\rangle^{\sat} = J, &\\
& \mathscr{B}\ \longrightarrow\ \langle \mathscr{B} \rangle^{\sat}\ \longrightarrow\  \big\{\langle \mathscr{B} \rangle^{\sat}_s\big\} = \mathscr{B}.
\end{eqnarray*}

Let $J\subset \K[x]$ be a Borel-fixed ideal such that the Hilbert polynomial of $\K[x]/J$ is equal to $p(t)$ and let $\mathscr{N}= \pos{n}{s}\setminus \{J_s\}$. Obviously $\big\vert\restrict{\mathscr{N}}{0}\big\vert = \vert\mathscr{N}\vert = p(s)= \Delta^0 p(s)$. Using the short exact sequence \eqref{eq:hyperplaneSection}, we determine the Borel ideal $I = (J,x_0) \cap \K[x_1,\ldots,x_n]$ with module $\K[x_1,\ldots,x_n]/I$ having Hilbert polynomial $\Delta p(t)$.
Thus being $\{I_s\} = \restrict{\{J_s\}}{1} \subset \pos{n-1}{s}$, $\left\vert \restrict{\mathscr{N}}{1} \right\vert = \left\vert \{I_s\}^{\mathcal{C}} \right\vert = \Delta p(s)$. Since $I$ is Borel-fixed in the polynomial ring $\K[x_1,\ldots,x_n]$ we can repeat the reasoning with the hyperplane section defined by $x_1 = 0$ and so on.

Let us now consider a Borel set $\mathscr{B} \subset \pos{n}{s}$, such that the complement $\mathscr{N} = \mathscr{B}^{\mathcal{C}}$ satisfies the condition $\left\vert\restrict{\mathscr{N}}{i}\right\vert = \Delta^{i} p(s)$ for every $i$. Firstly $\reg \langle\mathscr{B}\rangle^{\sat} \leqslant s$ by Proposition \ref{prop:regularityDegreeGenerators}, so let us prove that $\K[x]/\langle\mathscr{B}\rangle$ has Hilbert polynomial $p(t)$. We proceed by induction on the degree $d$ of the Hilbert polynomial. For any $n$, if $\deg p(t) = 0$, then $\restrict{\mathscr{N}}{i} = \emptyset$, for every $i \geqslant 1$, since $\Delta p(t) = 0$, that is for any $x^\beta \in \mathscr{N}$, $\min x^\beta = 0$. Applying repeatedly Lemma \ref{lem:removingMinimal} starting from the Hilbert polynomial $\overline{p}(t) = 0$ (corresponding to the ideal $(1)$), we obtain that $\langle \mathscr{B}\rangle^{\sat}$ defines a module $\K[x]/\langle \mathscr{B}\rangle^{\sat}$ having constant Hilbert polynomial $p(t) (= p(s))$. Let us know suppose that the map $\mathscr{B} \rightarrow \langle\mathscr{B}\rangle^{\sat}$ is well-defined for any Hilbert polynomial of degree $d-1$ and let $p(t)$ be a Hilbert polynomial of degree $d$. $\overline{\mathscr{B}} = \restrict{\mathscr{B}}{1} \subset \pos{n-1}{s}$ realizes the condition of the theorem w.r.t. the Hilbert polynomial $\Delta p(t)$ and $\deg \Delta p(t) = d-1$. Hence by the inductive hypothesis the ideal $\langle \overline{\mathscr{B}}\rangle^\sat \subset \K[x_1,\ldots,x_n]$ defines the module $\K[x_1,\ldots,x_n]/\langle \overline{\mathscr{B}}\rangle^\sat$ with Hilbert polynomial $\Delta p(t)$. Let $\overline{p}(t)$ be the Hilbert polynomial of $\K[x_0,\ldots,x_n]/\langle \overline{\mathscr{B}}\rangle^\sat$: $\overline{p}(t) = p(t) + a$, because $\Delta \overline{p}(t) = \Delta p(t)$. $\langle \overline{\mathscr{B}}\rangle^\sat$ turns out to be the $x_1$-saturation of $\langle \mathscr{B}\rangle^{\sat}$, so by Proposition \ref{prop:difference} the Hilbert polynomial of $\K[x_0,\ldots,x_n]/\langle \mathscr{B}\rangle^{\sat}$ differs by a constant from $\overline{p}(t)$ and since $\vert \mathscr{N} \vert = \vert \restrict{\mathscr{N}}{0}\vert = p(r)$ it coincides with $p(t)$.
\end{proof}

\begin{corollary}\label{cor:bijectionAllIdeals}
Let $p(t)$ be an admissible Hilbert polynomial in $\PP^n$ whose Gotzmann number is $r$. There is a bijective function
\begin{equation}
\begin{array}{rcl}
\left\{ \begin{array}{c}J \subset \K[x] \text{ saturated Borel-fixed}\\ \text{ideal s.t. } \K[x]/J \text{ has}\\ \text{Hilbert polynomial } p(t)\end{array}\right\} & \stackrel{1:1}{\longleftrightarrow} & \left\{\begin{array}{c} \mathscr{B} \subset \pos{n}{r} \text{ Borel set s.t.}\\ \text{set } \mathscr{N} = \pos{n}{r}\setminus \mathscr{B}\\ \Big\vert \restrict{\mathscr{N}}{i} \Big\vert = \Delta^i p(r),\ \forall\ i\end{array}\right\}\\
J & \longrightarrow & \{J_r\}\\
\langle\mathscr{B}\rangle^{\sat} & \longleftarrow & \mathscr{B}
\end{array}
\end{equation}
\end{corollary}
\begin{proof}
By Proposition \ref{prop:regularityDegreeGenerators}, Proposition \ref{prop:regularityVsCMregularity} ant Theorem \ref{th:RegularityTheorem}, any saturated Borel-fixed ideal $J$ defining a module $\K[x]/J$ with Hilbert polynomial $p(t)$ has regularity lower than or equal to the Gotzmann number of $p(t)$.
\end{proof}

Therefore to compute the saturated Borel-fixed ideals we can construct Borel sets with the prescribed property. The proof of Theorem \ref{th:bijectionIdealsSets} suggests to use a recursive algorithm: i.e. to determine the Borel sets in $\pos{n}{r}$ corresponding to the Hilbert polynomial $p(t)$, we begin computing the Borel sets in $\pos{n-1}{r}$ corresponding to the Hilbert polynomial $\Delta p(t)$. 

Let us examine more precisely this idea. Let $\overline{\mathscr{B}} \subset \pos{n-1}{r}$ a Borel set corresponding to the Hilbert polynomial $\Delta p(t)$ and let $\overline{\mathscr{N}} = \overline{\mathscr{B}}^{\mathcal{C}}$. In order for $\overline{\mathscr{B}}$ to be the restriction $\restrict{\mathscr{B}}{1}$ of a Borel set $\mathscr{B} \subset \pos{n}{r}$ (where $\pos{n}{r}$ contains one more variable smaller than variables in $\pos{n-1}{r}$), each monomial that can be obtained by decreasing moves from a monomial in $\overline{\mathscr{N}}$ has to belong to $\mathscr{N} = \mathscr{B}^{\mathcal{C}}$. This extension of an order set $\overline{\mathscr{N}} \subset \pos{n-1}{r}$ to an order set $\mathscr{N} \subset \pos{n}{r}$ has an ideal interpretation.

\begin{lemma}\label{lem:partition} Let $\overline{\mathscr{B}} \subset \pos{n-1}{r}$ be a Borel set and let $\overline{\mathscr{N}} = \overline{\mathscr{B}}^{\mathcal{C}}$. Moreover let $\mathscr{N} \subset \pos{n}{r}$ be the order set containing the monomials in $\overline{\mathscr{N}}$ and all those obtained by descreasing moves from them. Then,
\begin{equation}
\mathscr{N} = \pos{n}{r} \setminus \Big\{ \big(\langle\overline{\mathscr{B}}\rangle^{\sat} \cdot \K[x_0,\ldots,x_n]\big)_r \Big\}.
\end{equation}
\end{lemma} 
\begin{proof}
Let us call $\mathscr{B}$ the Borel set $\{ (\langle\overline{\mathscr{B}}\rangle^{\sat} \cdot \K[x_0,\ldots,x_n])_r\}$.
Let $x^\alpha = x_n^{\alpha_n}\cdots x_0^{\alpha_0}$ be a monomial of $\pos{n}{r}$ and suppose $\min x^\alpha = 0$, i.e. $\alpha_0 > 0$. The monomial $\alpha_0\up{0}(x^\alpha) = x^{\alpha_n}\cdots x_1^{\alpha_1 + \alpha_0}$\hfill belongs\hfill to\hfill $\pos{n-1}{r}$,\hfill so\hfill either\hfill $\alpha_0\up{0}(x^\alpha) \in \overline{\mathscr{B}}$\hfill or\\ $\alpha_0\up{0}(x^\alpha) \in \overline{\mathscr{N}}$.
If $x^{\alpha_n}\cdots x_1^{\alpha_1 + \alpha_0} \in \overline{\mathscr{B}}$, then $x_n^{\alpha_n}\cdots x_2^{\alpha_2}$ is in $\langle \overline{\mathscr{B}} \rangle^{\sat}$ and so $x^\alpha \in \langle\overline{\mathscr{B}}\rangle^{\sat} \cdot \K[x_0,\ldots,x_n]$, otherwise $x^{\alpha_n}\cdots x_1^{\alpha_1 + \alpha_0} \in \overline{\mathscr{N}}$ implies $x^{\alpha} \in \mathscr{N}$.
\end{proof}

At this point, by Proposition \ref{prop:difference}, we know that the Hilbert polynomial corresponding to a Borel set $\mathscr{B}$ of the type $\{(\langle\overline{\mathscr{B}}\rangle^{\sat}\cdot \K[x])_r\}$ differs from the target Hilbert polynomial by a constant: to determine this constant we compare the value $p(r)$ of the Hilbert polynomial $p(t)$ in degree $r$ with the cardinality of the order set $\mathscr{N}$ obtained by decreasing moves from $\overline{\mathscr{N}}$.

\begin{lemma}\label{lem:countingConditions}
Let $\overline{\mathscr{N}} \subset \pos{n-k}{r}$ be an order set and let $\mathscr{N} \subset \pos{n}{r}$ be the order set defined from $\overline{\mathscr{N}}$ by decreasing moves. Then,
\begin{equation}
\left\vert \mathscr{N} \right\vert = \sum_{\begin{subarray}{c} x^{\alpha} \in \overline{\mathscr{N}} \\ x^{\alpha}=x_n^{\alpha_n}\cdots x_{k}^{\alpha_k}\end{subarray}} \binom{\alpha_k+k}{k}.
\end{equation}
\end{lemma}
\begin{proof}
Each monomial $x^{\alpha} \in \overline{\mathscr{N}}$ imposes the belonging to $\mathscr{N}$ of any monomials obtained from it applying a composition of decreasing moves $\lambda_1 \down{1} \circ \cdots \circ \lambda_k \down{k}$. These type of moves act on the maximal power of $x_k$ in $x^\alpha$ and so they describe a poset isomoprhic to $\pos{k}{\alpha_k}$ and
\[
\left\vert \pos{k}{\alpha_k} \right\vert = \dim_{\K} \K[x_0,\ldots,x_k]_{\alpha_k} = \binom{k+\alpha_k}{k}. \qedhere
\]
\end{proof}

There are three possibilities:
\begin{itemize}
\item $p(r) - \vert \mathscr{N} \vert < 0$, $\overline{\mathscr{N}}$ imposes too many monomials ouside the ideal, so the hyperplane section defined by $\langle \overline{\mathscr{B}} \rangle^{\sat}$ has to be discarded (there exist no Borel-fixed ideals corresponding to $p(t)$ with such a hyperplane section);
\item $p(r) - \vert \mathscr{N} \vert = 0$, $\langle\overline{\mathscr{B}} \rangle^{\sat}  \subset \K[x_0,\ldots,x_n]$ is one of the ideals sought;
\item $p(r) - \vert \mathscr{N} \vert > 0$, applying repeatedly Lemma \ref{lem:removingMinimal} we determine the ideals we are looking for.
\end{itemize}

Putting together Remark \ref{rk:differenceBetweenHPs} with Lemma \ref{lem:partition}, we can establish a sharp upper bound of the difference $p(r) - \vert \mathscr{N} \vert$.

\begin{proposition}\label{prop:estimateRemovals}
Let $p(t)$ be an admissible Hilbert polynomial with Gotzmann number $r$ and let $r_1$ be the Gotzmann number of $\Delta p(t)$.
\begin{enumerate}[(i)]
\item\label{it:estimateRemovals_i} Given saturated Borel-fixed ideal $J \subset \K[x_1,\ldots,x_n]$ such that $\K[x_1,\ldots,x_n]/J$ has Hilbert polynomial $\Delta p(t)$, to pass from $\{(J\cdot\K[x_0,\ldots,x_n])_r\}\subset \pos{n}{r}$ to a Borel set corresponding to $p(t)$, we need to remove at most $r-r_1$ monomials.
\item\label{it:estimateRemovals_ii} We need to remove exactly $r-r_1$ monomials if we consider the lexicographic ideal $L \subset \K[x_1,\ldots,x_n]$ corresponding to the polynomial $\Delta p(t)$.
\item\label{it:estimateRemovals_iii} Let $I\subset \K[x]$ be a Borel-fixed ideal, such that the Hilbert polynomial of $\K[x]/I$ is $p(t)$.
To construct the Borel set $\mathscr{B} \subset \pos{r}{n}$ corresponding to $I$,  we need to remove at most $r$ monomials.
\end{enumerate}
\end{proposition}
\begin{proof}\emph{(\ref{it:estimateRemovals_i})}\hfill
Straightforward\hfill from\hfill Remark\hfill \ref{rk:differenceBetweenHPs},\hfill because\hfill the\hfill Hilbert\hfill polynomial\hfill of\\ $\K[x_0,\ldots,x_n]/(J\cdot\K[x_0,\ldots,x_n])$ belongs to $\HP{\Delta p(t)}$.

\emph{(\ref{it:estimateRemovals_ii})}
There are more than one way to prove this point. We exploit Lemma \ref{lem:countingConditions} in order to show that 
the ideal $L \cdot \K[x_0,\ldots,x_n]$ corresponds to the minimal polynomial in $\HP{\Delta p(t)}$. By definition, the order set $\overline{\mathscr{N}} = \{L_r\}^{\mathcal{C}} \subset \pos{n-1}{r}$ contains the smallest $\Delta p(r)$ monomials in $\pos{n-1}{r}$ w.r.t. the degree lexicographic order. To construct $\overline{\mathscr{N}}$ we can think to start from the order set $\{x_1^r\}\subset \pos{n-1}{r}$ and to remove successively the minimum w.r.t. the lexicographic order (see Remark \ref{rk:minMax}) among the minimal monomials of the complement. This minimum can be detected also looking at the elementary increasing move by which we can reach it: it will correspond to the lowest index possible. So whenever it is possible we add a monomials reached with $\up{1}$ so that the number of monomials smaller than it in $\pos{n}{r}$ decreases. Any other choice will generate an order set $\mathscr{N}$ with more elements.

\emph{(\ref{it:estimateRemovals_iii})} It comes directly applying \emph{(\ref{it:estimateRemovals_i})} recursively on the construction of the Borel sets corresponding to any $\Delta^i p(t)$.
\end{proof}

\subsection{The pseudocode description of the algorithm}

In Algorithm \ref{alg:BorelGeneratorBFS}, we give a pseudocode description of the algorithm just designed. Of some auxiliary methods implementing basic operations, we only describe the requirements on the input and the result returned in the output. In Table \ref{tab:exectutionBorelGenerator}, we simulate an execution of \textsc{BorelGenerator} on $\K[x_0,x_1,x_2,x_3]$ and $p(t) = 3t+2$.

\begin{algorithm}[H]	
\begin{algorithmic}
\STATE $\textsc{GotzmannNumber}\big(p(t)\big)$
\REQUIRE $p(t)$, admissible Hilbert polynomial.
\ENSURE the Gotzmann number of $p(t)$.
\end{algorithmic}
\begin{algorithmic}
\STATE \textsc{MinimalElements}($\mathscr{B}$)
\REQUIRE $\mathscr{B}$, Borel set.
\ENSURE the set of minimal elements of $\mathscr{B}$ w.r.t. $\leq_B$.
\end{algorithmic}

\smallskip

\begin{algorithmic}[1]
\STATE\nonumber $\textsc{BorelGenerator}\big(\K[x_k,\ldots,x_h],p(t)\big)$
\REQUIRE $\K[x_k,\ldots,x_h]$, polynomial ring.
\REQUIRE $p(t)$, admissible Hilbert polynomial in $\PP^{h-k} = \Proj \K[x_k,\ldots,x_h]$.
\ENSURE the set of all Borel-fixed ideals in $\K[x_k,\ldots,x_h]$ defining subschemes of $\PP^{h-k}$ with Hilbert polynomial $p(t)$.
\IF{$p(t) = 0$}
\RETURN $\big\{(1)\big\}$;
\ENDIF
\end{algorithmic}
\end{algorithm}

\begin{algorithm}[H]
\caption[Algorithm computing the set of all saturated Borel-fixed ideals in a fixed polynomial ring with a fixed Hilbert polynomial.]{The algorithm computing the set of all saturated Borel-fixed ideals in a fixed polynomial ring with a fixed Hilbert polynomial.}\label{alg:BorelGeneratorBFS}
\begin{algorithmic}[1]\setcounter{ALC@line}{4}
\STATE $\textsf{hyperplaneSections} \leftarrow \textsc{BorelGenerator}\big(\K[x_{k+1},\ldots,x_h],\Delta p(t)\big)$; 
\STATE $\textsf{borelFixedIdeals} \leftarrow \emptyset$;
\STATE $r \leftarrow \textsc{GotzmannNumber}\big(p(t)\big)$;
\FORALL {$J \in \textsf{hyperplaneSections}$}
\STATE $\mathscr{B} \leftarrow \big\{(J \cdot \K[x_k,\ldots,x_h])_r\big\} \subset \pos{h-k}{r}$;
\STATE $q \leftarrow p(r) - \vert \mathscr{B}^{\mathcal{C}}\vert$;
\IF{$q \geqslant 0$}
\STATE $\textsf{borelFixedIdeals} \leftarrow \textsf{borelFixedIdeals} \cup \textsc{Remove}(\mathscr{B},q)$;
\ENDIF
\ENDFOR
\RETURN \textsf{borelFixedIdeals};
\end{algorithmic}

\bigskip

\label{alg:remove}
\begin{algorithmic}[1]
\STATE \textsc{Remove}($\mathscr{B}$,$q$)
\REQUIRE $\mathscr{B}$, a Borel set.
\REQUIRE $q$, the number of monomials to remove from $\mathscr{B}$.
\ENSURE the set of saturated Borel-fixed ideals obtained from Borel sets constructed from $\mathscr{B}$ removing in all the possible ways $q$ monomials.
\IF{$q = 0$}
\RETURN $\{ \langle \mathscr{B} \rangle^{\sat} \}$;
\ELSE
\STATE $\textsf{borelIdeals} \leftarrow \emptyset$;
\STATE $\textsf{minimalMonomials} \leftarrow \textsc{MinimalElements}(\mathscr{B})$;
\FORALL{$x^\alpha \in \textsf{minimalMonomials}$}
 \STATE $\textsf{borelIdeals} \leftarrow \textsf{borelIdeals} \cup \textsc{Remove}\big(\mathscr{B}\setminus\{x^{\alpha}\},q-1\big)$;
\ENDFOR
\RETURN \textsf{borelIdeals};
\ENDIF
\end{algorithmic}
\end{algorithm}

\pagebreak

\begin{table}[H]
\dirtree{%
.1 $\textsc{BorelGenerator}\big(\K[x_0,x_1,x_2,x_3],3t+2\big)$ . 
.2 $\textsc{BorelGenerator}\big(\K[x_1,x_2,x_3],3\big)$ (because $3t+2 \neq 0$) .
.3 $\textsc{BorelGenerator}\big(\K[x_2,x_3],0\big)$ (because $3 \neq 0$) .
.4 $\textbf{return } \{(1)\}$ in $\K[x_2,x_3]$ .
.3 The Gotzmann number of $\Delta p(t) = 3$ is 3 . 
.3 $\mathscr{B}_0 = \big\{\big((1)\cdot\K[x_1,x_2,x_3] \big)_3\big\} = \pos{2}{3}$ and $q_0 = 3 - \vert \mathscr{B}_0^{\mathcal{C}} \vert = 3$ .
.4 $\textsc{Remove}(\mathscr{B}_0,3)$ .
.5 $\textsc{Remove}(\mathscr{B}_0\setminus\{x_1^3\},2)$ . 
.6 $\textsc{Remove}(\mathscr{B}_0\setminus\{x_1^3,x_2x_1^2\},1)$ .
.7 $\textsc{Remove}(\mathscr{B}_0\setminus\{x_1^3,x_2x_1^2,x_2^2 x_1\},0)$ .
.7 $\textsc{Remove}(\mathscr{B}_0\setminus\{x_1^3,x_2x_1^2,x_3x_1^2\},0)$ .
.3 $\textbf{return } \{(x_3,x_2^3),(x_3^2,x_3x_2,x_2^2)\}$ in $\K[x_1,x_2,x_3]$ .
.2 The Gotzmann number of $p(t) = 3t+2$ is 5 . 
.2 $\mathscr{B}_{1,1} = \big\{\big((x_3,x_2^3)\cdot\K[x_0,x_1,x_2,x_3] \big)_5\big\} \subset \pos{3}{5}$ and $q_{1,1} = 17 - \vert \mathscr{B}_{1,1}^{\mathcal{C}} \vert = 2$ .
.3  $\textsc{Remove}(\mathscr{B}_{1,1},2)$ .
.4  $\textsc{Remove}(\mathscr{B}_{1,1}\setminus\{x_2^3x_0^2\},1)$ .
.5  $\textsc{Remove}(\mathscr{B}_{1,1}\setminus\{x_2^3x_0^2,x_2^3x_1x_0\},0)$ .
.5  $\textsc{Remove}(\mathscr{B}_{1,1}\setminus\{x_2^3x_0^2,x_3x_0^4\},0)$ .
.4  $\textsc{Remove}(\mathscr{B}_{1,1}\setminus\{x_3x_0^4\},1)$ .
.5  $\textsc{Remove}(\mathscr{B}_{1,1}\setminus\{x_3x_0^4,x_2^3x_0^2\},0)$ (already found).
.5  $\textsc{Remove}(\mathscr{B}_{1,1}\setminus\{x_3x_0^4,x_3x_1x_0^3\},0)$.
.2 $\mathscr{B}_{1,2} = \big\{\big((x_3^2,x_3x_2,x_2^2)\cdot\K[x_0,x_1,x_2,x_3] \big)_5\big\} \subset \pos{3}{5}$ and\\ \phantom{spa} $q_{1,2} = 17 - \vert \mathscr{B}_{1,2}^{\mathcal{C}} \vert = 1$ .
.3 $\textsc{Remove}(\mathscr{B}_{1,2},1)$ .
.4 $\textsc{Remove}(\mathscr{B}_{1,2}\setminus\{x_2^2 x_0^3\},0)$ .
.2 $\textbf{return } \{(x_3,x_2^4,x_2^3x_1^2),(x_3^2,x_3x_2,x_3x_1,x_2^4,x_2^3x_1),(x_3^2,x_3x_2,x_2^3,x_3x_1^2),$ $(x_3^2,x_3x_2,x_2^3,x_2^2x_1)\}$ in $\K[x_0,x_1,x_2,x_3]$ .
}
\normalsize
\caption{\label{tab:exectutionBorelGenerator} The diagram of the execution of \textsc{BorelGenerator} with as inputs the polynomial ring $\K[x_0,x_1,x_2,x_3]$ and the Hilbert polynomial $p(t)=3t+2$.}
\end{table}

\pagebreak

The example described in Table \ref{tab:exectutionBorelGenerator} shows a first inaccuracy of the strategy, in fact \textsc{BorelGenerator} could compute many times the same ideal (the Borel set $\mathscr{B}_{1,1}\setminus\{x_3x_0^4,x_2^3x_0^2\}$ corresponding to the ideal $(x_3^2,x_3x_2,x_3x_1,x_2^4,x_2^3 x_1)$ is obtained 2 times). To solve this problem, we can use a total order on the monomials, so we fix any term ordering $\sigma$, and then keep trace of the computation: we add as argument of the function $\textsc{Remove}$ a monomial (that usually will be the last monomial removed) and we consider as monomials to remove only those greater than it.

\begin{algorithm}[H]
\caption[Modified version of Algorithm \ref{alg:remove} to avoid repetitions of ideals.]{The modified version of Algorithm \ref{alg:remove} to avoid repetitions of ideals.}
\label{alg:removeWithoutRepetions}
\begin{algorithmic}[1]
\STATE $\textsc{RemoveUniqueness}(\mathscr{B},q,x^\beta)$
\REQUIRE $\mathscr{B}$, a Borel set.
\REQUIRE $q$, the number of monomials to remove from $\mathscr{B}$.
\REQUIRE $x^\beta$, monomial (usually it will be a maximal element of $\mathscr{B}$).
\ENSURE the set of saturated Borel-fixed ideals obtained from Borel sets constructed from $\mathscr{B}$ removing in all the possible ways $q$ monomials without repetitions.
\IF{$q = 0$}
\RETURN $\{ \langle \mathscr{B} \rangle^{\sat} \}$;
\ELSE
\STATE $\textsf{borelIdeals} \leftarrow \emptyset$;
\STATE $\textsf{minimalMonomials} \leftarrow \textsc{MinimalElements}(\mathscr{B})$;
\FORALL{$x^\alpha \in \textsf{minimalMonomials}$}
\IF{$x^\alpha >_{\texttt{DegLex}} x^\beta$}
 \STATE $\textsf{borelIdeals} \leftarrow \textsf{borelIdeals} \cup \textsc{RemoveUniqueness}\big(\mathscr{B}\setminus\{x^{\alpha}\},q-1,x^{\alpha}\big)$;
 \ENDIF
\ENDFOR
\RETURN \textsf{borelIdeals};
\ENDIF
\end{algorithmic}
\end{algorithm}

We remark that Algorithm \ref{alg:BorelGeneratorBFS} could be naturally interpreted as an algorithm visiting a tree. Let us consider any Hilbert polynomial $p(t)$ of degree $d$, admissible for the projective space $\PP^n$ (i.e. $d<n$). We can associate to the pair $(p(t),\PP^n)$ the rooted tree defined as follows:
\begin{itemize}
 \item the nodes are all Borel-fixed ideals of $\K[x_i,\ldots,x_n]$ with Hilbert polynomial $\Delta^i p(t),\ \forall\ i = 0,\ldots,d+1$;
 \item the father of $I \subset \K[x_i,\ldots,x_n]$ is the ideal $J \subset \K[x_{i+1},\ldots,x_n]$ such that $J = \big(I\vert_{x_i=1} \cap \K[x_{i+1},\ldots,x_n]\big)^\sat$, that is $J$ represents the hyperplane section of $I$ w.r.t. $x_i$.
\end{itemize}

\begin{figure}[H]
\begin{center}
\begin{tikzpicture}[>=latex,line join=bevel,scale=0.55]
 \node (a) at (14bp,86.007bp) [] {\tiny $(1)$};
  \node (d5) at (537bp,385.01bp) [] {\tiny $(x_4^2,x_4x_3,x_4x_2,x_4x_1,x_3^7,x_3^6x_2,x_3^6x_1,x_3^5x_2^3)$};
  \node (d8) at (537bp,185.01bp) [] {\tiny $(x_4^2,x_4x_3,x_4x_2,x_4x_1,x_3^7,x_3^6x_2,x_3^5x_2^2)$};
  \node (d2) at (537bp,584.01bp) [] {\tiny $(x_4,x_3^7,x_3^6x_2,x_3^6x_1,x_3^5x_2^4,x_3^5x_2^3x_1)$};
  \node (d3) at (537bp,518.01bp) [] {\tiny $(x_4^2,x_4x_3,x_4x_2,x_4x_1,x_3^6,x_3^5x_2^4,x_3^5x_2^3x_1)$};
  \node (d6) at (537bp,318.01bp) [] {\tiny $(x_4^2,x_4x_3,x_4x_2,x_4x_1^2,x_3^6,x_3^5x_2^3)$};
  \node (d4) at (537bp,451.01bp) [] {\tiny $(x_4,x_3^7,x_3^6x_2,x_3^5x_2^3,x_3^6x_1^2)$};
  \node (d7) at (537bp,252.01bp) [] {\tiny $(x_4,x_3^7,x_3^6x_2,x_3^5x_2^3,x_3^5x_2^2x_1)$};
  \node (b1) at (110bp,219.01bp) [] {\tiny $(x_4,x_3^5)$};
  \node (b2) at (110bp,86.007bp) [] {\tiny $(x_4^2,x_4x_3,x_3^4)$};
  \node (b3) at (110bp,19.007bp) [] {\tiny $(x_4^2,x_4x_3^2,x_3^3)$};
  \node (c3) at (280bp,285.01bp) [] {\tiny $(x_4^2,x_4x_3,x_4x_2,x_3^6,x_3^5x_2^2)$};
  \node (c2) at (280bp,352.01bp) [] {\tiny $(x_4,x_3^7,x_3^6x_2,x_3^5x_2^2)$};
  \node (c1) at (280bp,418.01bp) [] {\tiny $(x_4,x_3^6,x_3^5x_2^3)$};
  \node (d1) at (537bp,651.01bp) [] {\tiny $(x_4,x_3^6,x_3^5x_2^4,x_3^5x_2^3x_1^2)$};
  \node (c7) at (280bp,19.007bp) [] {\tiny $(x_4^2,x_4x_3,x_3^4)$};
  \node (c6) at (280bp,86.007bp) [] {\tiny $(x_4^2,x_4x_3^2,x_4x_3x_2,x_4x_2^2,x_3^5)$};
  \node (c5) at (280bp,152.01bp) [] {\tiny $(x_4^2,x_4x_3,x_4x_2^3,x_3^5)$};
  \node (c4) at (280bp,219.01bp) [] {\tiny $(x_4^2,x_4x_3,x_4x_2^2,x_3^6,x_3^5x_2)$};
  \draw [-] (a) -- (b2);
  \draw [-] (c1) -- (d8);
  \draw [-] (a) -- (b1);
  \draw [-] (b1) -- (c2);
  \draw [-] (c1) -- (d3);
  \draw [-] (c1) -- (d4);
  \draw [-] (b1) -- (c6);
  \draw [-] (c1) -- (d2);
  \draw [-] (b1) -- (c5);
  \draw [-] (c1) -- (d6);
  \draw [-] (b1) -- (c4);
  \draw [-] (a) -- (b3);
  \draw [-] (b2) -- (c7);
  \draw [-] (c1) -- (d5);
  \draw [-] (b1) -- (c1);
  \draw [-] (c1) -- (d1);
  \draw [-] (c1) -- (d7);
  \draw [-] (b1) -- (c3);
  \node at (14bp,-25bp) [] {\footnotesize $\begin{array}{c}\Delta^3 p(t) = 0\\ \K[x_3,x_4]\end{array}$};
  \node at (110bp,-25bp) [] {\footnotesize $\begin{array}{c}\Delta^2 p(t) = 5\\ \K[x_2,x_3,x_4]\end{array}$};
  \node at (280bp,-25bp) [] {\footnotesize $\begin{array}{c}\Delta^1 p(t) = 5t-2\\ \K[x_1,x_2,x_3,x_4]\end{array}$};
  \node at (537bp,-25bp) [] {\footnotesize $\begin{array}{c}p(t) = \frac{5}{2}t^2 + \frac{1}{2}t - 8\\ \K[x_0,x_1,x_2,x_3,x_4]\end{array}$};
\end{tikzpicture}
\caption[The tree of Borel-fixed ideals defining surfaces in $\PP^4$ with Hilbert polynomial $p(t)=\frac{5}{2}t^2 + \frac{1}{2}t - 8$.]{\label{fig:treeBorelIdeals} The tree of Borel-fixed ideals associated to $\PP^4$ and $p(t)=\frac{5}{2}t^2 + \frac{1}{2}t - 8$.}
\end{center}
\end{figure}
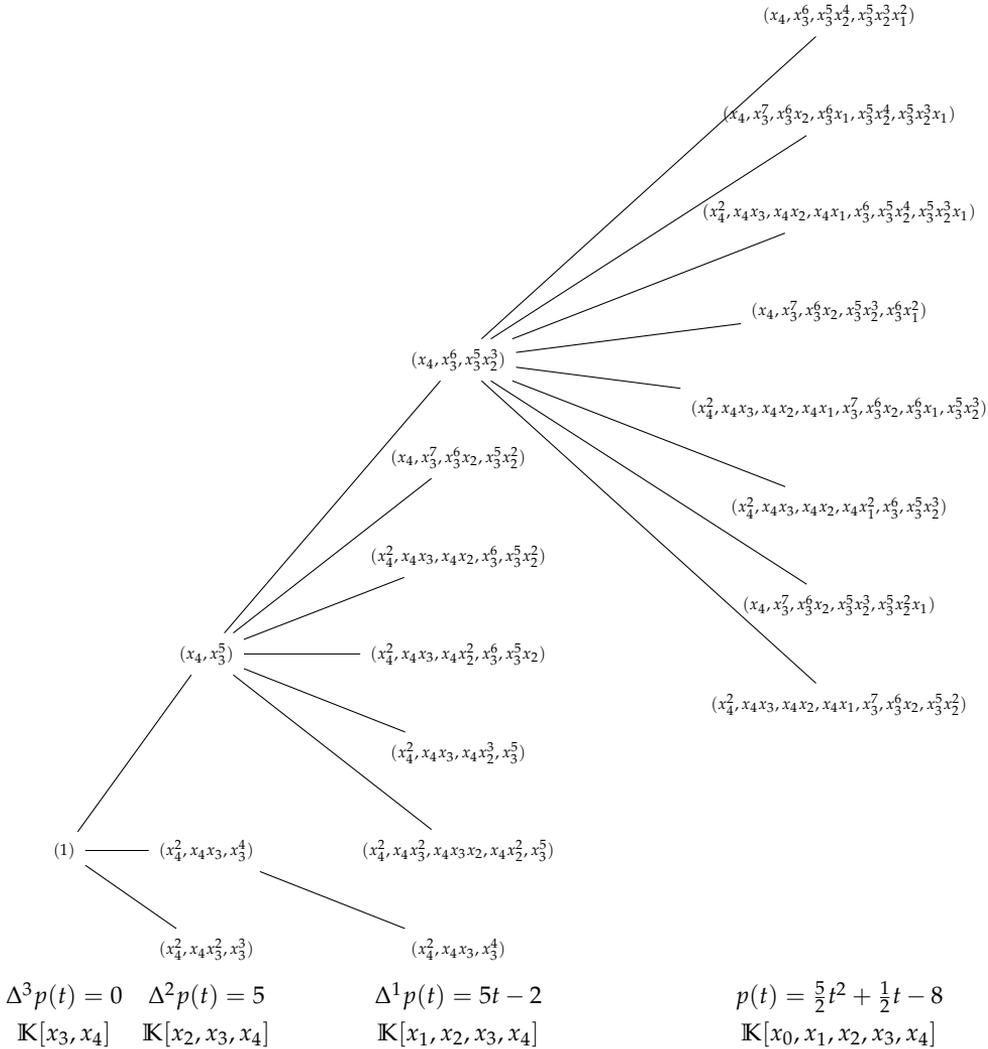

With such a definition, we have that the root of the tree is the ideal $(1) \subset \K[x_{d+1},\ldots,x_n]$ if $d<n-1$, or the ideal $\big(x_{n}^{\Delta^{d} p(t)}\big) \subset \K[x_{n-1},x_n]$, if $d = n-1$, and the Borel-fixed ideals defining subschemes of $\PP^n$ with Hilbert polynomial $p(t)$ are represented by the leaves at maximal distance from the root (see for an example Figure \ref{fig:treeBorelIdeals}). Algorithm \ref{alg:BorelGeneratorBFS} turns out to be a BFS (Breadth First Search) on the tree, indeed to determine the leaves at maximal distance we have to examine before all the nodes closer to the root, that from a computational point of view means that we need to store in the memory of a computer all the intermediate steps. Figure \ref{fig:treeBorelIdeals} clearly shows that this approach could not be optimal also because generally there are many ideals that will be finally discarded by the algorithm (because imposing too many monomials outside the ideal) but that we keep in mind for a long time before examing them.

Therefore a better approach is to visit the nodes of the tree of Borel-fixed ideals by means of a DFS (Depth First Search) visiting algorithm, so that the algorithm discards an ideal (if necessary) immediatly after having determined it. In Algorithm \ref{alg:BorelGeneratorDFS} and Algorithm \ref{alg:coreDFS}, there is the description of this strategy, which is used in the effective implementation \texttt{BorelIdeals} of the algorithm in the package \texttt{HSC} (see Appendix \ref{ch:HSCpackage}).

\begin{algorithm}[H]
\begin{algorithmic}[1]
\STATE $\textsc{BorelGeneratorHS}\big(\K[x_0,\ldots,x_n],p(t),k,I\big)$
\REQUIRE $\K[x_0,\ldots,x_n]$, polynomial ring.
\REQUIRE $p(t)$, admissible Hilbert polynomial in $\PP^{n} = \Proj \K[x_0,\ldots,x_n]$.
\REQUIRE $k$, integer s.t. $0\leqslant k \leqslant \deg p(t)$.
\REQUIRE $I$, Borel-fixed ideal in $\K[x_{k},\ldots,x_n]$ s.t. $\K[x_{k},\ldots,x_n]/I$ has Hilbert polynomial $\Delta^{k} p(t)$.
\ENSURE the set of all Borel-fixed ideals $J$ in $\K[x_0,\ldots,x_n]$ defining subschemes of $\PP^{n}$ with Hilbert polynomial $p(t)$ s.t. $\big(J\vert_{x_0 = \ldots = x_{k} = 1}\big)^{\sat} = I$.
\IF{$k = 0$}
\RETURN $\{I\}$;
\ENDIF
\STATE $r \leftarrow \textsc{GotzmannNumber}\big(\Delta^{k-1}p(t)\big)$;
\end{algorithmic}
\end{algorithm}

\begin{algorithm}[H]
\caption[Core of the DFS strategy to compute Borel-fixed ideals defining subschemes of $\PP^n$ with Hilbert polynomial $p(t)$.]{The core of the DFS strategy to compute Borel-fixed ideals defining subschemes of $\PP^n$ with Hilbert polynomial $p(t)$. This function visits a node and then calls itself on the children of the node.}
\label{alg:coreDFS}
\begin{algorithmic}[1]\setcounter{ALC@line}{5}
\STATE $\mathscr{B} \leftarrow \big\{(I\cdot\K[x_{k-1},\ldots,x_n])_r\big\} \subset \pos{n-k+1}{r}$;
\STATE $q \leftarrow \Delta^{k-1} p(r) - \left\vert\mathscr{B}^{\mathcal{C}}\right\vert$;
\IF{$q\geqslant 0$}
\STATE $\textsc{HS} \leftarrow \textsc{RemoveUniqueness}(\mathscr{B},q,0)$;
\STATE $\textsc{BorelIdeals} \leftarrow \emptyset$;
\FORALL{$\widetilde{I} \in \textsc{HS}$}
\STATE $\textsc{BorelIdeals} \leftarrow {}$ \parbox[t]{10cm}{$\textsc{BorelIdeals} \cup {}$\\ $\textsc{BorelGeneratorHS}\big(\K[x_0,\ldots,x_n],p(t),k-1,\widetilde{I}\big)$;}
\ENDFOR
\RETURN $\textsc{BorelIdeals}$;
\ELSE
\RETURN $\emptyset$;
\ENDIF
\end{algorithmic}
\end{algorithm}

\begin{algorithm}[H]
\caption[Function detecting the root of the tree associated to ${\big(\K[x],p(t)\big)}$ and then starting the DFS visit.]{This function detects the root of the tree associated to $\big(\K[x],p(t)\big)$ and then starts the DFS visit.}
\label{alg:BorelGeneratorDFS}
\begin{algorithmic}[1]
\STATE $\textsc{BorelGeneratorDFS}\big(\K[x_0,\ldots,x_n],p(t)\big)$
\REQUIRE $\K[x_0,\ldots,x_n]$, polynomial ring.
\REQUIRE $p(t)$, admissible Hilbert polynomial in $\PP^{n} = \Proj \K[x_0,\ldots,x_n]$.
\ENSURE the set of all Borel-fixed ideals in $\K[x_0,\ldots,x_n]$ defining subschemes of $\PP^{n}$ with Hilbert polynomial $p(t)$.
\STATE $d \leftarrow \deg p(t)$;
\IF{$d = n-1$}
\STATE $c \leftarrow \Delta^d p(t)$;
\RETURN $\textsc{BorelGeneratorHS}\big(\K[x_0,\ldots,x_n],p(t),d,(x_{n}^{c})\big)$;
\ELSE
\RETURN $\textsc{BorelGeneratorHS}\big(\K[x_0,\ldots,x_n],p(t),d+1,(1)\big)$
\ENDIF
\end{algorithmic}
\end{algorithm}

\begin{table}[H]
\dirtree{%
.1 $\textsc{BorelGeneratorDFS}\big(\K[x_0,x_1,x_2,x_3],3t+2\big)$ . 
.2 $d = \deg 3t+2 = 1 < 2 = 3-1$ .
.2 $\textsc{BorelGeneratorHS}\big(\K[x_0,x_1,x_2,x_3],3t+2,2,(1)\big)$ .
.3 $\Delta^1 (3t+2) = 3 \Rightarrow r = 3$, $\mathscr{B}_0 = \pos{2}{3}$ and $q_0 = 3$ .
.4 $\textsc{RemoveUniqueness}(\mathscr{B}_0,3,0)$ .
.5 $\textsc{RemoveUniqueness}(\mathscr{B}_0\setminus\{x_1^3\},2,x_1^3)$ .
.6 $\textsc{RemoveUniqueness}(\mathscr{B}_0\setminus\{x_1^3,x_2x_1^2\},1,x_2x_1^2)$ .
.7 $\textsc{RemoveUniqueness}(\mathscr{B}_0\setminus\{x_1^3,x_2x_1^2,x_2^2 x_1\},0,x_2^2x_1)$ .
.7 $\textsc{RemoveUniqueness}(\mathscr{B}_0\setminus\{x_1^3,x_2x_1^2,x_3x_1^2\},0,x_3x_1^2)$ .
.4 $\textbf{return } \{(x_3,x_2^3),(x_3^2,x_3x_2,x_2^2)\}$ . 
.3 $\textsc{BorelGeneratorHS}\big(\K[x_0,x_1,x_2,x_3],3t+2,1,(x_3,x_2^3)\big)$ .
.4 $\Delta^0 (3t+2) = 3t+2\Rightarrow r = 5$,\\ \phantom{---}$\mathscr{B}_{1,1} = \big\{\big((x_3,x_2^3)\cdot\K[x_0,x_1,x_2,x_3]\big)_5\big\}$ and $q_{1,1} = 17-15=2$ .
.5 $\textsc{RemoveUniqueness}(\mathscr{B}_{1,1},2,0)$ .
.6 $\textsc{RemoveUniqueness}(\mathscr{B}_{1,1}\setminus\{x_2^3x_0^2\},1,x_2^3x_0^2)$ .
.7 $\textsc{RemoveUniqueness}(\mathscr{B}_{1,1}\setminus\{x_2^3x_0^2,x_2^3 x_1x_0\},0,x_2^3 x_1x_0)$ .
.7 $\textsc{RemoveUniqueness}(\mathscr{B}_{1,1}\setminus\{x_2^3x_0^2,x_3x_0^4\},0,x_3x_0^4)$ .
.6 $\textsc{RemoveUniqueness}(\mathscr{B}_{1,1}\setminus\{x_3x_0^4\},1,x_3x_0^4)$ .
.7 $\textsc{RemoveUniqueness}(\mathscr{B}_{1,1}\setminus\{x_3x_0^4,x_3x_1x_0^3\},0,x_3x_1x_0^3)$ .
.4 $\textbf{return } \{(x_3,x_2^4,x_2^3x_1^2),(x_3^2,x_3x_2,x_3x_1,x_2^4,x_2^3x_1),(x_3^2,x_3x_2,x_2^3,x_3x_1^2)\}$ .
.3 $\textsc{BorelGeneratorHS}\big(\K[x_0,x_1,x_2,x_3],3t+2,1,(x_3^2,x_3x_2,x_2^2)\big)$ .
.4 $\Delta^0 (3t+2) = 3t+2\Rightarrow r = 5$,\\ \phantom{---}$\mathscr{B}_{1,2} = \big\{\big((x_3^2,x_3x_2,x_2^2)\cdot\K[x_0,x_1,x_2,x_3]\big)_5\big\}$ and $q_{1,2} = 17-16=1$ .
.5 $\textsc{RemoveUniqueness}(\mathscr{B}_{1,2},1,0)$ .
.6 $\textsc{RemoveUniqueness}(\mathscr{B}_{1,2}\setminus\{x_3x_1x_0^3\},0,x_3x_1x_0^3)$ .
.4 $\textbf{return } \{(x_3^2,x_3x_2,x_2^3,x_2^2x_1)\}$ .
.2 $\textbf{return } \{(x_3,x_2^4,x_2^3x_1^2),(x_3^2,x_3x_2,x_3x_1,x_2^4,x_2^3x_1),(x_3^2,x_3x_2,x_2^3,x_3x_1^2),$ $(x_3^2,x_3x_2,x_2^3,x_2^2x_1)\}$ in $\K[x_0,x_1,x_2,x_3]$ .
}
\normalsize
\caption{\label{tab:exectutionBorelGeneratorDFS} The diagram of the execution of \textsc{BorelGeneratorDFS} with as arguments $\K[x_0,x_1,x_2,x_3]$ and $3t+2$ (cf. Table \ref{tab:exectutionBorelGenerator}).}
\end{table}

\section{How many Borel-fixed ideals are there?}

An interesting question that naturally arises from the algorithm projected in the previous section is if we are able to predict the number of Borel-fixed ideals associated to any pairs $\big(n,p(t)\big)$ without computing all them. In Table
\ref{tab:exampleN}, there is a summary of the number of Borel-fixed ideals in the case of some constant Hilbert polynomial.

\begin{table}[!ht]
\begin{center}
\footnotesize
\begin{tabular}{c | c | c | c | c | c | c | c | c | c}
 & $n=2$ & $n=3$ & $n=4$ & $n=5$ & $n=6$ & $n=7$ & $n=8$ & $n=9$ & $n=10$ \\
\cline{1-10}\cline{1-10}\cline{1-10}
$p(t)=2$ & 1 & 1 & 1 & 1 & 1 & 1 & 1 & 1 & 1\\
\cline{1-10}
$p(t)=3$ & 2 & 2 & 2 & 2 & 2 & 2 & 2 & 2 & 2\\
\cline{1-10}
$p(t)=4$ & 2 & 3 & 3 & 3 & 3 & 3 & 3 & 3 & 3 \\
\cline{1-10}
$p(t)=5$ & 3 & 4 & 5 & 5 & 5 & 5 & 5 & 5 & 5 \\
\cline{1-10}
$p(t)=6$ & 4 & 6 & 7 & 8 & 8 & 8 & 8 & 8 & 8 \\
\cline{1-10}
$p(t)=7$ & 5 & 9 & 11 & 12 & 13 & 13 & 13 & 13 & 13 \\
\cline{1-10}
$p(t)=8$ & 6 & 12 & 16 & 18 & 19 & 20 & 20 & 20 & 20 \\
\cline{1-10}
$p(t)=9$ & 8 & 17 & 24 & 28 & 30 & 31 & 32 & 32 & 32 \\
\cline{1-10}
$p(t)=10$ & 10 & 24 & 35 & 42 & 46 & 48 & 49 & 50 & 50 \\
\end{tabular}
\normalsize
\end{center}
\caption{\label{tab:exampleN} The number of Borel-fixed ideals in the case of constant Hilbert polynomials $p(t) = s$ in $\PP^n$, for $2 \leqslant n \leqslant 10$ and $2 \leqslant s \leqslant 10$.}
\end{table}

\begin{definition}
Given a projective space $\PP^n$ and an admissible Hilbert polynomial $p(t)$, we denote by \gls{BorelNpoly} the set of all Borel-fixed ideals defining subschemes of $\PP^n$ with Hilbert polynomial $p(t)$ and by \gls{NOB} its cardinality
\begin{equation}
\NOB{n}{p(t)} = \left\vert \Borel{n}{p(t)} \right\vert.
\end{equation}
\end{definition}

The dependence of $\NOB{n}{p(t)}$ on the Hilbert polynomial seems a very hard task to achieve, whereas already the numbers showed in Table \ref{tab:exampleN} suggests that the behavior of $\NOB{n}{p(t)}$ varying $n$ could be more treatable. Thus let us start discussing how the dimension $n$ of the projective space affects $\NOB{n}{p(t)}$.

\begin{remark}\label{rk:linearSpaces}
Let $p(t) = \binom{t+d}{d}$ be the Hilbert polynomial of a $d$-projective space contained in $\PP^n$ for any $0\leqslant d < n$. The Gotzmann number of such a Hilbert polynomial is 1, so that we consider posets of the type $\pos{n}{1}$, that turn out to be totally ordered sets.
\begin{center}
\begin{tikzpicture}[>=latex,scale=1]
\node (1) at (0,0) [draw,ellipse] {\footnotesize $x_n$};
\node (2) at (2,0) [draw,ellipse] {\footnotesize $x_{n-1}$};
\node (3) at (4,0) [draw,ellipse] {\footnotesize$x_{n-2}$};
\node (4) at (9,0) [draw,ellipse] {\footnotesize$x_{1}$};
\node (5) at (11,0) [draw,ellipse] {\footnotesize $x_{0}$};
\draw [->] (1) -- (2);
\draw [->] (2) -- (3);
\draw [->] (3) -- (5.5,0);
\draw [->] (8,0) -- (4);
\draw [->] (4) -- (5);
\draw [loosely dotted,very thick] (6,0) -- (7.5,0);
\end{tikzpicture}
\end{center}
Therefore
\begin{equation}
\Borel{n}{\binom{t+d}{d}} = \big\{(x_n,\ldots,x_{d+1})\big\} \quad\Longrightarrow\quad \NOB{n}{\binom{t+d}{d}} = 1.
\end{equation}
\end{remark}

The inclusion $\K[x_0,\ldots,x_n] \hookrightarrow \K[x_0,\ldots,x_n,x_{x+1}]$ with $x_{n+1} >_B x_n$ naturally extends to the corresponding posets in any degree. Thinking about the characterization of Borel-fixed ideals given in Corollary \ref{cor:bijectionAllIdeals} and looking for a setting that allows to consider simultaneously Borel-fixed ideals in any number of variables we introduce the following setting.
\begin{definition}\index{poset!infinite}
Let $\{x_0,\ldots,x_n,\ldots\}$ an infinite set of variables. We denote by $\infpos{m}$ the poset composed by monomials of degree $m$ in the variables $\{x_0,\ldots,x_n,\ldots\}$ and given by the natural extension to an infinite number of variables of the elementary moves. Moreover for any finite poset $\pos{n}{m}$, we denote by $i_n:\pos{n}{m} \hookrightarrow \infpos{m}$ the inclusion map.
\end{definition}

Now we need to extend the notion of Borel set.
\begin{definition}\index{Borel set}\index{order set}
A subset $\mathscr{B} \subset \infpos{m}$ is called \emph{Borel set} if the complement $\mathscr{B}^{\mathcal{C}} = \infpos{m}\setminus\mathscr{B}$ is a finite order set, i.e. $\mathscr{B}^{\mathcal{C}}$ is closed w.r.t. decreasing elementary moves. 

For any $n$, we define (and denote again with $i_n$) the map
\begin{equation}\label{eq:embeddingBorelSets}
\begin{split}
i_n : \left\{\text{Borel sets of } \pos{n}{m}\right\} &\ \longrightarrow\ \left\{\text{Borel sets of } \infpos{m}\right\}\\
\parbox{4cm}{\centering $\mathscr{B}$}&\ \longmapsto\  \parbox{3.5cm}{\centering $\left(i_n(\mathscr{B}^{\mathcal{C}})\right)^{\mathcal{C}}$}
\end{split}
\end{equation}
so that the order sets defined by $\mathscr{B}$ and $i_n(\mathscr{B})$ are equal, and the map
\begin{equation}\label{eq:restrictionBorelSets}
\begin{split}
s_{n} : \left\{\text{Borel sets of } \infpos{m}\right\} &\ \longrightarrow\ \left\{\text{Borel sets of } \pos{n}{m}\right\}\\
\parbox{3.5cm}{\centering $\mathscr{B}$}&\ \longmapsto\ \parbox{4cm}{\centering $\mathscr{B} \cap \pos{n}{m}$}.
\end{split}
\end{equation}
Obviously $s_n(i_n(\mathscr{B})) = \mathscr{B}$.
\end{definition}

Given a Hilbert polynomial $p(t)$ with Gotzmann number $r$, let us now consider the inclusion
\begin{equation}\label{eq:IdealsToInfinitePos}
\begin{split}
i_n : \left\{\begin{array}{c} \mathscr{B} \subset \pos{n}{r} \text{ Borel set s.t.}\\ \text{set } \mathscr{N} = \pos{n}{r}\setminus \mathscr{B}\\ \Big\vert \restrict{\mathscr{N}}{i} \Big\vert = \Delta^i p(r),\ \forall\ i\end{array}\right\} &\ \longrightarrow\ \left\{\begin{array}{c} \mathscr{B} \subset \infpos{r} \text{ Borel set s.t.}\\ \text{set } \mathscr{N} = \infpos{r}\setminus \mathscr{B}\\ \Big\vert \restrict{\mathscr{N}}{i} \Big\vert = \Delta^i p(r),\ \forall\ i\end{array}\right\}.
\end{split}
\end{equation}
that results still well-defined. To understand how the number of Borel-fixed ideals of $\K[x]$ with fixed Hilbert polynomial is affected by $n$, we will discuss the property of the map \eqref{eq:IdealsToInfinitePos}.

\begin{definition}
Let $\mathscr{N}$ be a finite order set (either in $\pos{n}{m}$ or $\infpos{m}$). We define the maximal variable of $\mathscr{N}$ as
\begin{equation}\label{eq:maxN}
\gls{maxvar} = \max \left\{\max x^\beta\ \big\vert\ x^\beta \in \mathscr{N}\right\}.
\end{equation}
\end{definition}

\begin{proposition}\label{prop:maxOrderSetLinearGen}
Let $I \subset \K[x]$ be a saturated Borel-fixed ideal and let $\mathscr{N} = \{I_r\}^\mathcal{C} \subset \pos{n}{r}$ be the corresponding order set, where $r$ is the Gotzmann number of the Hilbert polynomial of $\K[x]/I$. Then
\[
x_j \text{ is a minimal generator of } I \quad\Longleftrightarrow\quad \maxvar \mathscr{N} < j.
\]
\end{proposition}
\begin{proof}
($\Leftarrow$) If $j > \maxvar \mathscr{N}$, then $x_j x_0^{r-1}$ belongs to $\{I_r\}$, so $x_j$ is a generator of $I$.

($\Rightarrow$) Since $x_j \in I$, $x_j x_0^{r-1}$ does not belong to $\mathscr{N}$. Moreover any other monomial $x^\alpha$ of degree $r$ such that $\max x^\alpha \geqslant j$ cannot belong to $\mathscr{N}$, because $x^{\alpha} \geq_B x_jx_0^{r-1}$. In fact $x^\alpha \geq_B \frac{x^\alpha}{\max x^\alpha} x_j \geq_B x_j x_0^{r-1}$ since $x_0^{r-1}$ is the minimum w.r.t. the Borel order among monomials of degree $r-1$.
\end{proof}

\begin{proposition}
Let $\mathscr{B} \subset \infpos{m}$ be a Borel set. If $\maxvar \mathscr{B}^{\mathcal{C}} = j$, then
\begin{equation}
i_n \left(s_n\left(\mathscr{B}\right)\right) = \mathscr{B},\qquad\forall\ n \geqslant j.
\end{equation}
\end{proposition}
\begin{proof}
$\maxvar \mathscr{B}^{\mathcal{C}} = j$ implies $\mathscr{B}^{\mathcal{C}} \subset \pos{n}{m}$, for all $n \geqslant j$.
\end{proof}

\begin{lemma}\label{lem:firstBound}
Let $\mathscr{N} \subset \infpos{m}$ be any order set such that $\left\vert\mathscr{N}\right\vert = a$. Then
\begin{equation}
\maxvar \mathscr{N} \leqslant a-1.
\end{equation}
\end{lemma}
\begin{proof}
It comes directly from the remark that to reach $x_0^m$ from a monomial $x^\alpha$ s.t. $\max x^\alpha > a-1$ we need at least $a$ decreasing elementary moves, i.e. any order set containing $x^\alpha$ would contain at least $a+1$ elements.

The bound is sharp, indeed for the order set $\overline{\mathscr{N}} = \{x_0^m,x_0^{m-1}x_1,\ldots,$ $x_0^{m-1}x_{a-1}\}$ $\left\vert\overline{\mathscr{N}}\right\vert = a$ and $\maxvar \overline{\mathscr{N}} = a-1$.
\end{proof} 

\begin{proposition}
Let $p(t)$ be an admissible Hilbert polynomial with Gotzmann number $r$. Then for any $n \geqslant p(r)-1$
\begin{equation}
\left\{\begin{array}{c} \mathscr{B} \subset \pos{n}{r} \text{ Borel set s.t.}\\ \text{set } \mathscr{N} = \pos{n}{r}\setminus \mathscr{B}\\ \Big\vert \restrict{\mathscr{N}}{i} \Big\vert = \Delta^i p(r),\ \forall\ i\end{array}\right\} \stackrel{1:1}{\longleftrightarrow} \left\{\begin{array}{c} \mathscr{B} \subset \pos{n+1}{r} \text{ Borel set s.t.}\\ \text{set } \mathscr{N} = \pos{n+1}{r}\setminus \mathscr{B}\\ \Big\vert \restrict{\mathscr{N}}{i} \Big\vert = \Delta^i p(r),\ \forall\ i\end{array}\right\}
\end{equation}
\end{proposition}
\begin{proof}
Using the Borel sets in the infinite poset $\infpos{r}$ as intermediate step, the maps giving the bijection are $s_{n+1}\circ i_n$ and $s_n \circ i_{n+1}$
\end{proof}

From the point of view of saturated Borel-fixed ideals, since $n+1$ is surely greater than $\maxvar \mathscr{N}$, the correspondence turns out to be
\[
\begin{array}{ccc}
\left\{ \begin{array}{c}J \text{ saturated Borel-fixed}\\ \text{ideal s.t. } \K[x_0,\ldots,x_n]/J \text{ has}\\ \text{Hilbert polynomial } p(t)\end{array}\right\}& \longrightarrow &
\left\{ \begin{array}{c}J'\text{ saturated Borel-fixed}\\ \text{ideal s.t. } \K[x_0,\ldots,x_{n+1}]/J' \text{ has}\\ \text{Hilbert polynomial } p(t)\end{array}\right\}\\
J & \longmapsto & (x_{n+1},J).
\end{array}
\]

\begin{corollary}\label{cor:NOBstabilizes}
Let $p(t)$ be an admissible Hilbert polynomial with Gotzmann number $r$. The number $\NOB{n}{p(t)}$ is constant for $n \geqslant p(r)-1$.
\end{corollary}

\begin{definition}
We denote by $\NOB{\bullet}{p(t)}$ the sequence of the the number of Borel-fixed ideals associated to a fixed Hilbert polynomial $p(t)$ and varying number of variables $n$:
\begin{equation}
\NOB{\bullet}{p(t)} = \left( \ldots,\NOB{n}{p(t)},\ldots\right)
\end{equation}
If $d = \deg p(t)$, then $p(t)$ is admissible in $\PP^n$ for $n > d$, so the first integer of the sequence will be always $\NOB{d+1}{p(t)}$. Moreover since by Corollary \ref{cor:NOBstabilizes} the sequence at some point becomes constant, we could write the sequence as a finite list of integer
\[
\gls{NOBsequence} = \left( \NOB{d+1}{p(t)},\ldots,\NOB{A}{p(t)}\right)
\]
meaning that $\NOB{A+k}{p(t)} = \NOB{A}{p(t)},\ \forall\ k \geqslant 0$.
\end{definition}

In general, we expect that the bound $p(r)-1$ is an overestimation of the point of stabilization of  $\NOB{n}{p(t)}$, because the order set $\overline{\mathscr{N}} = \{x_0^r,\ldots,x_{p(r)-1}x_0^{r-1}\}$ constructed in Lemma \ref{lem:firstBound} corresponds to a Borel-fixed ideal in any $\K[x_0,\ldots,x_n],\ n \geqslant p(r)-1$ with \emph{constant} Hilbert polynomial $\overline{p}(t) = p(r)$, indeed the hyperplane section is the ideal $(1) \subset \K[x_1,\ldots,x_n]$. For this reason, we carry on with a more detailed analisys.

Let $p(t)$ be an admissible Hilbert polynomial of degree $d$ with Gotzmann number $r$. Thinking about the recursive strategy of Algorithm \ref{alg:BorelGeneratorDFS}, we want determine the Borel-fixed ideal defining the order set $\mathscr{N}$ with maximum $\maxvar \mathscr{N}$ among the ideals with a given hyperplane section $\overline{I} \subset \K[x_1,\ldots,x_n]$, for some $n$.

Let $\mathscr{\overline{\mathscr{B}}} = \big\{(\overline{I}\cdot\K[x_0,\ldots,x_n])_r\big\} \subset \pos{n}{r}$ and let $\overline{\mathscr{N}}$ be the associated order set, viewed in the infinite poset $\infpos{r}$. We saw that the Hilbert polynomial $\overline{p}(t)$ associated to $\overline{\mathscr{N}}$ differs from the Hilbert polynomial $p(t)$ by a constant (Proposition \ref{prop:difference}). Set $c = p(t) - \overline{p}(t)$, to determine an order set corresponding to $p(t)$ we have to add $c$ monomials to $\overline{\mathscr{N}}$ and we want to achieve it using as many variables as possible.
Let $A = \maxvar \overline{\mathscr{N}}$ ($=\maxvar \{\overline{I}_r\}^{\mathcal{C}}$ by construction); by Proposition \ref{prop:maxOrderSetLinearGen}, we know that $x_{A+1} x_0^{r-1} \notin \overline{\mathscr{N}}$ and moreover this monomial is minimal in $\overline{\mathscr{N}}^{\mathcal{C}}$, so $\overline{\mathscr{N}} \cup \{x_{A+1} x_0^{r-1}\}$ is still an order set and
\[
\maxvar \left(\overline{\mathscr{N}} \cup \{x_{A+1} x_0^{r-1}\}\right) = A+1= \maxvar \overline{\mathscr{N}} + 1.
\]
Repeating the reasoning $c$ times we construct the order set $\mathscr{N}=\overline{\mathscr{N}}\cup\{x_{A+1} x_0^{r-1},$ $\ldots,x_{A+c}x_0^{r-1}\}$ and
\[
\begin{split}
\maxvar \mathscr{N}&{} = \maxvar \left(\overline{\mathscr{N}}\cup\{x_{A+1} x_0^{r-1},\ldots,x_{A+c}x_0^{r-1}\}\right) = \\
&{} = A+c = \maxvar \overline{\mathscr{N}} + c.
\end{split}
\]

We summarize this construction in the following lemma.
\begin{lemma}\label{lem:estimateMaxSingleHS}
Let $p(t)$ be an admissible Hilbert polynomial with Gotzmann number $r$ and let $\overline{I} \subset \K[x_1,\ldots,x_n]$ the saturated Borel-fixed ideal of an admissible hyperplane section of $p(t)$. Moreover let $\overline{p}(t)$ be the Hilbert polynomial associated to $\overline{I}\cdot\K[x_0,\ldots,x_n]$ and $c = p(t) - \overline{p}(t)$. 
\begin{enumerate}[(i)]
\item $\max \left\{ \max \mathscr{N}\text{ s.t. } \vert\mathscr{N}\vert = p(r) \text{ and } \restrict{\mathscr{N}}{1} =  \{\overline{I}_r\}^{\mathcal{C}}  \right\} = \maxvar \{\overline{I}_r\}^{\mathcal{C}}+c$.
\item Any Borel-fixed ideal $I\subset \K[x_0,\ldots,x_n]$ associated to $p(t)$ with hyperplane section $\overline{I}$, contains as minimal generator $x_j,\ \maxvar \{\overline{I}_r\}^{\mathcal{C}}+c < j \leqslant n$.
\end{enumerate}
\end{lemma}

So to estimate a bound, we can consider among the hyperplane sections defining order sets with same maximum the one with smallest Hilbert polynomial.

\begin{proposition}\label{prop:maxLexSection}
Let $p(t)$ be an admissible Hilbert polynomial of degree $d$ with Gotzmann number $r$ and let $L \subset \K[x_1,\ldots,x_n]$ the saturated lexicographic ideal\index{lexicographic ideal} associated to $\Delta p(t)$. For any ideal $I \subset \K[x_0,\ldots,x_n]$ having $L$ as hyperplane section
\begin{equation}
\maxvar \{I_r\}^{\mathcal{C}} \leqslant \begin{cases} d+ p(t) - \Sigma\big(\Delta p\big)(t), & \text{if } p(t) = \binom{t+d}{d} + c,\\ d+ p(t) - \Sigma\big(\Delta p\big)(t)+1, & \text{otherwise.} \end{cases}
\end{equation}
\end{proposition}
\begin{proof}
To determine the smallest number of variables that we need to determine an order set in $\pos{n-1}{r}$ associated to $\Delta p(t)$, we think to Macaulauy's Theorem \cite{Macaulay}. It says that whenever the Hilbert polynomial is admissible a lexicographic ideal $L$ realizing it exists, i.e. the lexicographic ideal has to be the ideal with the order set involving the smallest number of variables. 

If $p(t) = \binom{t+d}{d} + c$, then $\Delta p(t) = \binom{t+d-1}{d-1}$ and by Remark \ref{rk:linearSpaces} the lexicographic ideal is the only Borel-fixed ideal $(x_n,\ldots,x_{d+1}) \subset \K[x_1,\ldots,x_n]$, so that by Proposition \ref{prop:maxOrderSetLinearGen}
\[
\maxvar \{(x_n,\ldots,x_{d+1})_r\}^{\mathcal{C}} = d.
\]
For any other Hilbert polynomial the condition of admissibility is $n>d$, so the order set $\{L_r\}^{\mathcal{C}} \subset \pos{n-1}{r}$ will involve the variables $x_1,\ldots,x_{d+1}$ (regardless of $n$) and
\[
\maxvar \{L_r\}^{\mathcal{C}} = d+1.
\]
Let $D$ be the maximum of the order set defined by the lexicographic ideal.

In the proof of Proposition \ref{prop:estimateRemovals} we showed that the ideal $L\cdot\K[x_0,\ldots,x_n]$ is associated to the minimal polynomial among those with first difference equal to $\Delta p(t)$ and by Definition \ref{def:minimalPolynomial} such minimal polynomial is $\Sigma\big(\Delta p\big)(t)$. Hence called $\mathscr{N}$ the order set defined by $\{(L\cdot\K[x_0,\ldots,x_n])_r\}^{\mathcal{C}}$, to obtain an order set associated to $p(t)$ we need add $p(t) - \Sigma\big(\Delta p\big)(t)$ monomials. With the goal of constructing the order ideal involving as many variables as possible, we begin adding $x_{D+1} x_0^{r-1}$, since $\maxvar \mathscr{N} = D \Rightarrow x_{d+1} x_0^{r-1} \notin \mathscr{N}$. Repeating the reasoning $p(t) - \Sigma\big(\Delta p\big)(t)$ times we can obtain as limit case an order set with maximum variable equal to $D+p(t) - \Sigma\big(\Delta p\big)(t)$.
\end{proof}

Now we compare the lexicographic hyperplane section with any other section defining an order set with maximum greater than the maximum of the order set defined by the lexicographic ideal associated to $\Delta p(t)$ ($d$ or $d+1$).

\begin{lemma}\label{lem:estimateHSs}
Let $\overline{I} \subset \K[x_1,\ldots,x_n]$ be an admissible hyperplane sections for the Hilbert polynomial $p(t)$ of degree $d$ with Gotzmann number $r$ and let $\overline{p}(t)$ the Hilbert polynomial associated to $\overline{I}\cdot\K[x_0,\ldots,x_n]$. Moreover let $D$ be the maximum of the order set defined by the lexicographic ideal that realizes $\Delta p(t)$. If $\maxvar \{\overline{I}_r\}^\mathcal{C} = D + B$, then
\begin{enumerate}[(i)]
\item\label{it:estimateHSs_i} $\overline{p}(t) \geqslant \Sigma\big(\Delta p\big)(t) + B$;
\item\label{it:estimateHSs_ii} $\maxvar \{\overline{I}_r\}^\mathcal{C} + p(t) - \overline{p}(t) \leqslant D + p(t) - \Sigma\big(\Delta p\big)(t)$.
\end{enumerate}
\end{lemma}
\begin{proof} 
\emph{(\ref{it:estimateHSs_i})} Preliminarly we can simplify the problem considering the case $B=1$ and then proceeding iteratively.

Furthermore we suppose that $\{\overline{I}_r\}^{\mathcal{C}}$ contains a single monomial with maximum variable equal to $D+1$, i.e. the monomial $x_{D+1} x_1^{r-1}$, because we want to not move too away from the lexicographic ideal in order to preserve a small polynomial $\overline{p}(t)$ (this idea, at this point almost intuitive, will be clarified in Chapter \ref{ch:deformations}).

Applying Lemma \ref{lem:partition} and Lemma \ref{lem:firstBound} we know that
\[
\overline{p}(t) = p(t) - \left(p(r) - \sum_{x^{\alpha} \in \{\overline{I}_r\}^{\mathcal{C}}} (\alpha_1 + 1)\right).
\]
The order sets defined by $\overline{I}$ and by the lexicographic ideal surely differs for at least 1 monomial: $x_{D+1} x_1^{r-1}$. It is replaced necessarly by a monomial divided by a power of $x_1$ lower than $r-1$ because all the monomials divided by a power of $x_1$ greater than or equal to $r-1$ already belong to $\{\overline{I}\}^{\mathcal{C}}$, so
\[
\sum_{x^{\alpha} \in \{\overline{I}_r\}^{\mathcal{C}}} (\alpha_1 + 1) > \Sigma\big(\Delta p\big)(r) \qquad\Longrightarrow\qquad \text{\emph{(\ref{it:estimateHSs_i})}}\quad \overline{p}(t) > \Sigma\big(\Delta p\big)(t)
\]

\emph{(\ref{it:estimateHSs_ii})} By the inequality \emph{(\ref{it:estimateHSs_i})}
\[
p(t) - \overline{p}(t) < p(t) - \Sigma\big(\Delta p\big)(t)\ \Leftrightarrow\ p(t) - \overline{p}(t) \leqslant p(t) - \Sigma\big(\Delta p\big)(t) - 1 
\]
and finally
\[
\maxvar \{\overline{I}_r\}^{\mathcal{C}} + p(t) - \overline{p}(t) = D + 1 + p(t) - \overline{p}(t) \leqslant D + p(t) - \Sigma\big(\Delta p\big)(t). \qedhere
\]
\end{proof}

\begin{theorem}
Let $p(t)$ be an admissible Hilbert polynomial of degree $d$ with Gotzmann number $r$.
For any $\mathscr{N} \subset \infpos{r}$ order set defined by a Borel-fixed ideal associated to $p(t)$,
\begin{equation}
\maxvar \mathscr{N} \leqslant \begin{cases} d+ p(t) - \Sigma\big(\Delta p\big)(t),& \text{if } p(t) = \binom{t+d}{d} + c \\ d+ p(t) - \Sigma\big(\Delta p\big)(t) +1 ,& \text{otherwise.} \end{cases}
\end{equation}
\end{theorem}
\begin{proof}
It comes directly applying Lemma \ref{lem:estimateHSs} and Proposition \ref{prop:maxLexSection}. Moreover Proposition \ref{prop:maxLexSection} ensures that the bound is sharp.
\end{proof}

\begin{corollary}\label{cor:stabilizationGeneralPoly}
Let $p(t)$ be an admissible Hilbert polynomial of degree $d$ with Gotzmann number $r$. Then,
\begin{equation}
\NOB{n}{p(t)} = \NOB{n+1}{p(t)},\qquad \forall\ n \geqslant \begin{cases} d+ p(t) - \Sigma\big(\Delta p\big)(t),& \text{if } p(t) = \binom{t+d}{d} + c \\ d+ p(t) - \Sigma\big(\Delta p\big)(t) +1 ,& \text{otherwise.} \end{cases}
\end{equation}
\end{corollary}

\begin{example}
Let us check empirically the statement of Corollary \ref{cor:stabilizationGeneralPoly} on some examples, with the help of the java function \texttt{BorelIdeals} of the package \texttt{HSC}.

\noindent$\boldsymbol{p(t) = 5t+1.}$ From the Gotzmann representation
\[
\begin{split}
5t+1&{} = \binom{t+1}{1}+\binom{t}{1}+\binom{t-1}{1}+\binom{t-2}{1}+\binom{t-3}{1}+\\
&{}+\binom{t-5}{0}+\binom{t-6}{0}+\binom{t-7}{0}+\binom{t-8}{0}+\binom{t-9}{0}+\binom{t-10}{0}
\end{split}
\]
we can easily compute
\[
5t+1 - \Sigma\big(\Delta (5t+1)\big) = 5t+1 - (5t-5) = 6
\]
so that the sequence $\NOB{\bullet}{5t+1}$ is supposed to stabilize for $n=2+6 = 8$. In fact
\[
\NOB{\bullet}{5t+1} = (4,38,71,89,95,97,98,98,\ldots).
\]

\noindent$\boldsymbol{p(t) = \frac{3}{2}t^2 + \frac{7}{2}t.}$ The Gotzmann representation is
\[
\begin{split}
\frac{3}{2}t^2 + \frac{7}{2}t&{} = \binom{t+2}{2}+\binom{t+1}{2}+\binom{t}{2}+\binom{t-2}{1}+\binom{t-3}{1}+\\
&{} + \binom{t-5}{0}+\binom{t-6}{0}+\binom{t-7}{0}+\binom{t-8}{0}
\end{split}
\]
so
\[
\frac{3}{2}t^2 + \frac{7}{2}t + 1 - \Sigma\left(\Delta \left(\frac{3}{2}t^2 + \frac{7}{2}t\right)\right) = \frac{3}{2}t^2 + \frac{7}{2}t - \left(\frac{3}{2}t^2 + \frac{7}{2}t -4\right) = 4
\]
and the sequence $\NOB{\bullet}{\frac{3}{2}t^2 + \frac{7}{2}t}$ levels off for $n=3+4= 7$. In fact
\[
\NOB{\bullet}{\frac{3}{2}t^2 + \frac{7}{2}t} = (8,27,36,39,40,40,\ldots).
\]

\noindent$\boldsymbol{p(t) = \frac{1}{3}t^3 + \frac{3}{2}t^2 + \frac{25}{6}t + 3.}$ The Gotzmann representation is
\[
\begin{split}
\frac{1}{3}t^3 + \frac{3}{2}t^2 + \frac{25}{6}t + 3&{} = \binom{t+3}{3}+\binom{t+2}{3}+\binom{t-1}{1}+\binom{t-2}{1}+\\
&{} + \binom{t-4}{0}+ \binom{t-5}{0}+ \binom{t-6}{0}+ \binom{t-7}{0}+ \binom{t-8}{0}
\end{split}
\]
so
\[
\frac{1}{3}t^3 + \frac{3}{2}t^2 + \frac{25}{6}t + 3 - \Sigma\left(\Delta \left(\frac{1}{3}t^3 + \frac{3}{2}t^2 + \frac{25}{6}t + 3\right)\right) = 5
\]
and the sequence $\NOB{\bullet}{\frac{1}{3}t^3 + \frac{3}{2}t^2 + \frac{19}{6}t + 5}$ levels off for $n=4+5= 9$. In fact
\[
\NOB{\bullet}{\frac{1}{3}t^3 + \frac{3}{2}t^2 + \frac{25}{6}t + 3} = (21,84,130,149,155,156,156,\ldots).
\]
\end{example}

\subsection{The special case of constant Hilbert polynomials}

We will now discuss in depth the case of constant Hilbert polynomial. The bound given in Corollary \ref{cor:stabilizationGeneralPoly}, if $p(t) = s = \binom{t+0}{0}+(s-1)$, turns out to be $s-1$.

\begin{definition}
For any constant Hilbert polynomial $p(t) = s$, we define for $1 \leqslant i \leqslant s-2$
\begin{equation}
\gls{DeltaNOB} = \NOB{s-i}{s} - \NOB{s-i-1}{s}
\end{equation}
and
\begin{equation}
\gls{DeltaNOBsequence} = \left(\Delta\NOB{1}{s},\Delta\NOB{2}{s},\ldots,\Delta\NOB{s-2}{s}\right)
\end{equation}
\end{definition}

\begin{example}
Let us consider the Hilbert polynomial $p(t) = 15$. Since
\[
\NOB{\bullet}{15} = (1,27,107,206,287,342,377,398,410,417,421,423,424),
\]
we obtain
\[
\Delta\NOB{\bullet}{15} = (1,1,2,4,7,12,21,35,55,81,99,80,26).
\]
\end{example}

Experimentally we noticed that $\Delta \NOB{i}{s}$ for $s \gg 0$ becomes constant. In the following proposition we explain this behavior. 

\begin{proposition}\label{prop:stabilizationDelta}
$\Delta \NOB{i}{s}$ is constant for $s \geqslant 2i-1$.
\end{proposition}
\begin{proof}
Since $\Delta \NOB{i}{s} = \NOB{s-i}{s} - \NOB{s-i-1}{s}$, we are interested in studying order sets $\mathscr{N} \subset \infpos{s}$, such that $\maxvar \mathscr{N} = s-i$, indeed if $\maxvar \mathscr{N} < s-i$, $\mathscr{N}$ can be defined by an ideal in $\K[x_0,\ldots,x_n]$ with $n \leqslant s-i-1$.

By Proposition \ref{prop:maxOrderSetLinearGen}, if $\maxvar \mathscr{N} = s-i$ then $x_{s-i}x_0^{s-1}$ belongs to $\mathscr{N}$ and being $\mathscr{N}$ closed by decreasing elementary moves
\[
\{x_{s-i}x_0^{s-1},x_{s-i-1}x_0^{s-1},\ldots,x_0^s\} \subset \mathscr{N}.
\]
Thus we have to study the remaining
\[
\left\vert\mathscr{N}\setminus\{x_{s-i}x_0^{s-1},x_{s-i-1}x_0^{s-1},\ldots,x_0^s\}\right\vert = s - (s-i+1) = i-1
\]
monomials (not depending on $s$!). The point is again to determine which is the greatest number of variables involved in a subset $\mathscr{M}$ such that $\mathscr{N} = \mathscr{M} \cup \{x_{s-i}x_0^{s-1},$ $\ldots,x_0^s\}$. The Borel set $\pos{s-i}{s}\setminus \{x_{s-i}x_0^{s-1},\ldots,x_0^s\}$ has only one minimal monomial: $x_1^2 x_0^{s-2}$, so surely $x_1^2 x_0^{s-2} \in \mathscr{M}$. Then the strategy is always the same: we choose the next monomial to add applying succesively the increasing elementary moves with the index as big as possible, i.e.
\[
x_1^2 x_0^{s-2}\ \stackrel{\up{1}}{\longrightarrow}\ x_2x_1x_0^{s-2} \ \stackrel{\up{2}}{\longrightarrow}\ x_3x_1x_0^{s-2} \ \stackrel{\up{3}}{\longrightarrow}\ \cdots
\]
In order for being able to apply this strategy $i-2$ times, that is being able to construct the order set
\[
\{x_{s-i}x_0^{s-1},\ldots,x_0^s\} \cup \{ x_1^2 x_0^{s-2},\ldots,x_{i-1} x_1 x_0^{s-2}\}
\]
we need $i-1 \leqslant s-i\ \Rightarrow\ s \geqslant 2i-1$.
\end{proof}

\begin{definition}
We define
\begin{equation}
\gls{DeltaNOBstable} = \Delta \NOB{i}{s},\qquad s \gg 0
\end{equation}
and we denote by \gls{DeltaNOBstableSequence} the sequence
\begin{equation}
\left(\Delta \NOB{1}{},\Delta \NOB{2}{},\ldots,\Delta \NOB{i}{},\ldots\right).
\end{equation}
\end{definition}

The optimal way to compute $\Delta \NOB{i}{}$ is to consider $\Delta \NOB{i}{2i-1} = \NOB{i-1}{2i-1}-\NOB{i-2}{2i-1}$, that is to count the saturated ideals of $\NOB{i-1}{2i-1}$ not having variables as generator. The first values of the sequence are
\[
1,1,2,4,7,12,21,35,58,96,156,251,403,639,1008,1582,2465,3821,5898,9055,\ldots.
\] 

Proposition \ref{prop:stabilizationDelta} (and its proof) suggests a new strategy for designing an algorithm computing Borel-fixed ideals with constant Hilbert polynomial. Let us consider the Hilbert polynomial $p(t) = s$ and the polynomial ring $\K[x]$. The sequences $\Delta\NOB{\bullet}{s}$ coincides with $\Delta\NOB{\bullet}{}$ for the first $i$ values with $i \leqslant C=\lfloor\frac{s+1}{2}\rfloor$. There are three cases.
\begin{description}
\item[$n \geqslant s$.] By Corollary \ref{cor:stabilizationGeneralPoly}, we know that all the possible order sets involve at most $s-1$ variables, hence we compute the ideals $\Borel{s-1}{s}$ and
\[
\Borel{n}{s} = \left\{(x_n,\ldots,x_s,I)\ \vert\ I \in \Borel{s-1}{s}\right\}.
\]
\item[$n < s - C$.] We compute $\Borel{n}{s}$ using Algorithm \ref{alg:BorelGeneratorDFS}.
\item[$s-C \leqslant n \leqslant s-1$.] We compute $\Borel{s-C-1}{s}$ using Algorithm \ref{alg:BorelGeneratorDFS} and then we apply repeatedly Proposition \ref{prop:stabilizationDelta}. Let us suppose to have computed $\Borel{k}{s}$. The set $\Borel{k+1}{s}$ surely contains the set 
\[
\left\{(x_{k+1},I)\quad \big\vert\quad I \in \Borel{k}{s}\right\}
\]
to which we would add the Borel-ideals defining order set with maximum equal to $k+1$. We are looking at 
\[
\NOB{k+1}{s} - \NOB{k}{s} = \Delta\NOB{i}{s},\qquad i=s-k-1,
\]
new order sets and the best way to study them is considering the Borel ideals in $\Borel{i-1}{2i-1}$. Let $\mathcal{O}$ the set of order sets with maximum equal to $i-1$. The order sets we are looking for are
\[
\overline{\mathcal{O}} = \left\{ \{x_{k+1}x_0^{s-1},\ldots,x_{i}x_0^{s-1}\} \cup x_0^{s-2i+1}\cdot\mathscr{N} \ \vert\ \mathscr{N}\in \mathcal{O}\right\}.
\]
\end{description}

\begin{algorithm}
\caption[Alternative strategy for computing Borel-fixed ideals with constant Hilbert polynomial.]{A new strategy for computing Borel-fixed ideals with constant Hilbert polynomial.}
\label{alg:algorithmConstantHP}
\begin{algorithmic}[1]
\STATE $\textsc{BorelGeneratorConstantHP}(\K[x_0,\ldots,x_n],s)$
\REQUIRE $\K[x_0,\ldots,x_n]$, polynomial ring.
\REQUIRE $s$, positive integer (the Hilbert polynomial).
\ENSURE the set of all Borel-fixed ideals in $\K[x_0,\ldots,x_n]$ with Hilbert polynomial $p(t)=s$.
\STATE $i \leftarrow s-n$;
\IF{$i \leqslant 0$}
\STATE $\textsf{idealsLessVars} \leftarrow \textsc{BorelGeneratorConstantHP}(\K[x_0,\ldots,x_{s-1}],s)$;
\STATE $\textsf{borelIdeals} \leftarrow \emptyset$;
\FORALL{$I \in \textsf{idealsLessVars}$}
\STATE $\textsf{borelIdeals} \leftarrow \textsf{borelIdeals} \cup \left\{(x_n,\ldots,x_s,I)\right\}$;
\ENDFOR
\RETURN $\textsf{borelIdeals}$;
\ELSIF{$i > \lfloor\frac{s+1}{2}\rfloor$}
\RETURN $\textsc{BorelGeneratorDFS}(\K[x_0,\ldots,n_n],s)$;
\ELSE
\STATE $\textsf{idealsLessVars} \leftarrow \textsc{BorelGeneratorConstantHP}(\K[x_0,\ldots,x_{n-1}],s)$; 
\STATE $\textsf{idealsLinGen} \leftarrow \emptyset$;
\FORALL{$I \in \textsf{idealsLessVars}$}
\STATE $\textsf{idealsLinGen} \leftarrow \textsf{idealsLinGen} \cup \left\{(x_n,I)\right\}$;
\ENDFOR
\STATE $\textsf{idealsDelta} \leftarrow \textsc{BorelGeneratorDFS}(\K[x_0,\ldots,x_{i-1}],2i-1)$;
\STATE $\textsf{idealsWithoutLinGen} \leftarrow \emptyset$;
\FORALL{$J \in \textsf{idealsDelta}$}
\IF{$\max \{J_{s}\}^{\mathcal{C}} = i-1$}
\STATE $\mathscr{N} \leftarrow \{J_{s}\}^{\mathcal{C}} \cup \{x_{i}x_0^{s-1},\ldots,x_n x_0^{s-1}\}$;
\STATE $\textsf{idealsWithoutLinGen} \leftarrow \textsf{idealsWithoutLinGen} \cup \left\{\langle\mathscr{N}^{\mathcal{C}}\rangle^{\sat} \right\}$;
\ENDIF
\ENDFOR
\RETURN $\textsf{idealsLinGen} \cup \textsf{idealsWithoutLinGen}$;
\ENDIF
\end{algorithmic}
\end{algorithm}

\begin{example}
Let us see how the strategy works in the concrete case of the polynomial ring $\K[x_0,\ldots,x_6]$ and the Hilbert polynomial $p(t) = 10$. Following Algorithm \ref{alg:algorithmConstantHP}, we have that $i = 10 - 6 = 4$ and since $s=10 > 7 = 2i-1$, we are sure that $\Delta\NOB{4}{10} = \Delta\NOB{4}{}$.

Among the Borel-ideals with Hilbert polynomial $2i-1 = 7$ in the polynomial ring $\K[x_0,\ldots,x_{i-1}] = \K[x_0,x_1,x_2,x_3]$, there are 4 saturated ideals without a linear generator, precisely
\[
\begin{split}
&J_1 = (x_3^2,x_3x_2,x_2^2,x_3x_1,x_2x_1,x_1^5) \ \Rightarrow\ \{(J_1)_{10}\}^{\mathcal{C}} = \mathscr{M}\cup\{x_1^2x_0^{8},x_{1}^3x_0^7,x_1^4 x_0^6\}\\
&J_2 = (x_3^2,x_3x_2,x_2^2,x_3x_1,x_2x_1^2,x_1^4) \ \Rightarrow\ \{(J_2)_{10}\}^{\mathcal{C}} = \mathscr{M}\cup\{x_1^2x_0^{8},x_2 x_1 x_0^8,x_{1}^3x_0^7\}\\
&J_3 = (x_3^2,x_3x_2,x_3x_1,x_2^3,x_2^2x_1,x_2x_1^2,x_1^3) \ \Rightarrow\ \{(J_3)_{10}\}^{\mathcal{C}} = \mathscr{M}\cup\{x_1^2x_0^{8},x_2 x_1 x_0^8,x_{1}^2x_0^8\}\\
&J_4 = (x_3^2,x_3x_2,x_2^2,x_3x_1^2,x_2x_1^2,x_1^3)  \ \Rightarrow\ \{(J_4)_{10}\}^{\mathcal{C}} = \mathscr{M}\cup\{x_1^2x_0^{8},x_2 x_1 x_0^8,x_3x_{1}x_0^8\}\\
\end{split}
\]
where $\mathscr{M}=\{x_3x_0^9,x_2x_0^9,x_1x_0^9,x_0^{10}\}$. Therefore, set $\overline{\mathscr{M}} = \{x_6x_0^9,x_5x_0^9,x_4x_0^9\} \cup \mathscr{M}$, the ideals in $\K[x_0,\ldots,x_6]$ with Hilbert polynomial $p(t) = 10$ without linear generators are
\[
\begin{split}
&\left\langle \left(\overline{\mathscr{M}}\cup\{x_1^2x_0^{8},x_{1}^3x_0^7,x_1^4 x_0^6\} \right)^{\mathcal{C}} \right\rangle^{\sat}\\
&\left\langle \left(\overline{\mathscr{M}}\cup\{x_1^2x_0^{8},x_2 x_1 x_0^8,x_{1}^3x_0^7\} \right)^{\mathcal{C}} \right\rangle^{\sat}\\
&\left\langle \left(\overline{\mathscr{M}}\cup\{x_1^2x_0^{8},x_2 x_1 x_0^8,x_{1}^2x_0^8\}\right)^{\mathcal{C}} \right\rangle^{\sat}\\
&\left\langle \left(\overline{\mathscr{M}}\cup\{x_1^2x_0^{8},x_2 x_1 x_0^8,x_3x_{1}x_0^8\} \right)^{\mathcal{C}} \right\rangle^{\sat}\\
\end{split}
\]

In addition to these 4 ideals, we have to consider the ideals with the same Hilbert polynomial in $\K[x_0,\ldots,x_5]$ to which we add $x_6$ as generator.
\end{example}

\section{Segment ideals}\label{sec:segments}

We conclude this chapter, trying to generalize the notion of lexicographic ideal. We recall that by a theorem of Macaulay \cite{Macaulay}, if a numerical function $f:\NN \rightarrow \NN$ is admissible, i.e. there exists a $\K$-algebra $A$ such that $\text{HF}_A(t) = f(t)$, then the direct sum\index{lexicographic ideal}
\[
L = \bigoplus_{t \in \NN} L_t, \quad L_t = \left\langle \text{biggest } \binom{n+t}{n}-f(t) \text{ monomials of } \K[x]_t \text{ w.r.t. }\texttt{DegLex}\right\rangle
\]
is an ideal of $\K[x]$ and $\K[x]/L$ has $f(t)$ has Hilbert function. $L$ is uniquely determined by $f(t)$, so it is called \emph{lexicographic ideal} associated to $f(t)$. In our context, we are mostly interested in Hilbert polynomials, so we would like to associate uniquely a saturated lexicographic ideal to any Hilbert polynomial $p(t)$. Let us start writing in a slightly different way the Gotzmann representation of $p(t)$ supposing $\deg p(t) = d$:
\begin{equation}\label{eq:gotzmannRepresentationLex}\index{Hilbert polynomial!Gotzmann's representation of a}
\begin{split}
p(t) &{}= \binom{t+d}{d} + \ldots + \binom{t+d-(b_d-1)}{d} +\\
&{}+ \binom{t+d-1-b_d}{d-1} + \ldots + \binom{t+d-1-(b_d+b_{d-1}-1)}{d-1} +\\
&{}+ \ldots + \\
&{}+ \binom{t-(b_d+\ldots+b_1)}{0} + \ldots + \binom{t-(b_d+\ldots+b_1+b_0-1)}{0}.
\end{split}
\end{equation}
Comparing \eqref{eq:gotzmannRepresentationLex} with \eqref{eq:GotzmannDecomposition}, we see that $b_j$ counts the number of $a_i$ equal to $j$.

We have the following characterization of the saturated lexicographic ideal.
\begin{proposition}\label{prop:satLexIdeal}
Let $p(t)$ be an admissible Hilbert polynomial in $\K[x_0,\ldots,x_n]$. The saturated lexicographic ideal $L$, such that $\K[x]/L$ has Hilbert polynomial $p(t)$, is
\begin{equation}\label{eq:satLexIdeal}
L = \left(x_n,\ldots x_{d+2}, x_{d+1}^{b_d+1},x_{d+1}^{b_d} x_{d}^{b_{d-1}+1},\ldots ,x_{d+1}^{b_d}\cdots x_2^{b_1+1},x_{d+1}^{b_d}\cdots x_1^{b_0}\right),
\end{equation}
where the exponents $b_j$ are the integers defined in \eqref{eq:gotzmannRepresentationLex}.
\end{proposition}

\begin{example}
Let us consider the Hilbert polynomial $p(t) = \frac{1}{3}t^3 + \frac{3}{2}t^2 + \frac{25}{6}t + 3$ in $\PP^6$. The Gotzmann representation of $p(t)$ is
\[
\begin{array}{l c l}	
\frac{1}{3}t^3 + \frac{3}{2}t^2 + \frac{25}{6}t + 3 = \binom{t+3}{3}+\binom{t+2}{3}+ & \Rightarrow & b_3 = 2\\
 & \Rightarrow & b_2 = 0\\
\phantom{\frac{1}{3}t^3 + \frac{3}{2}t^2 + \frac{25}{6}t + 3} + \binom{t-1}{1}+\binom{t-2}{1}+ & \Rightarrow & b_1 = 2\\
\phantom{\frac{1}{3}t^3 + \frac{3}{2}t^2 + \frac{25}{6}t + 3} + \binom{t-4}{0}+ \binom{t-5}{0}+ \binom{t-6}{0}+ \binom{t-7}{0}+ \binom{t-8}{0} & \Rightarrow & b_0 = 5\\
\end{array}
\]
so the corresponding saturated lexicographic ideal is
\[
L = (x_6,x_5,x_4^3,x_4^2x_3,x_4^2x_2^3,x_4^2x_2^2x_1^5).
\]
\end{example}

We want to study what happens considering any term ordering $\sigma$ instead of the lexicographic order.

\begin{definition}
A set $S$ of monomials of degree $t$ is called \emph{segment} w.r.t. the term ordering $\sigma$ if for any monomial $x^{\alpha} \in S$, then $x^{\beta} >_{\sigma} x^{\alpha} \Rightarrow x^{\beta} \in S$.
\end{definition}

\begin{definition}\label{def:segments}
Let $I \subset \K[x]$ be an non-null monomial ideal and let $\sigma$ be any term ordering.
\begin{itemize}
\item\index{segment ideal} $I$ is called \emph{segment ideal} w.r.t. $\sigma$, if for every $t \in \NN$, $I_t$ is a segment w.r.t. $\sigma$.
\item\index{segment ideal!hilb-segment ideal}\index{hilb-segment ideal|see{segment ideal}} $I$ is called \emph{hilb-segment ideal} w.r.t. $\sigma$, if $I_r$ is a segment w.r.t. $\sigma$, where $r$ is the Gotzmann number of the Hilbert polynomial of $\K[x]/I$.
\item\index{segment ideal!reg-segment ideal}\index{reg-segment ideal|see{segment ideal}} $I$ is called \emph{reg-segment ideal} w.r.t. $\sigma$, if $I_{\reg(I)}$ is a segment w.r.t. $\sigma$.
\item\index{segment ideal!gen-segment ideal}\index{gen-segment ideal|see{segment ideal}} $I$ is called \emph{gen-segment ideal}  w.r.t. $\sigma$, if for every $t \in \NN$, the generators of $I$ of degree $t$, i.e. the generators of $I_t\setminus\langle I_{t-1}\rangle_t$, are the biggest monomials w.r.t. $\sigma$ among the monomials of degree $t$ non contained in $\langle I_{t-1}\rangle$.
\end{itemize}
\end{definition}

By definition, the lexicographic ideal is a segment ideal w.r.t. the degree lexicographic order. In the following we will call it also lexsegment ideal.

As seen in Definition \ref{def:BorelOrder}, each term ordering $\sigma$ defines a total order on the monomials of a fixed degree that refines the Borel order. Hence it is clear that in order for a set $S$ to be a segment it is necessary to be a Borel set.

\begin{proposition}
\begin{enumerate}[(i)]
\item If $S$ is a segment, then $S$ is a Borel set.
\item If $I$ is an ideal that respects one properties described in Definition \ref{def:segments}, then $I$ is a Borel-fixed ideal.
\end{enumerate}
\end{proposition}

\begin{lemma}\label{lem:segmentLowerDegree}
Let $I\subset \K[x]$ be a saturated Borel-fixed ideal and let $\sigma$ be any term ordering. 
If $I_p$ is a segment then $I_q$, $q < p$, is a segment too. 
\end{lemma}
\begin{proof}
Let $x^{\alpha}$ and $x^\beta$ be two monomials of degree $q$, such that $x^\alpha \in I_q$ and $x^{\beta} >_{\sigma} x^{\alpha}$. $x^{\alpha} x_0^{p-q}$ belongs to $I_p$ and $x^{\beta} x_0^{p-q} >_{\sigma} x^{\alpha}x_0^{p-q}$ implies
 $x^{\beta}x_0^{p-q} \in I_p$, because $I_p$ is a segment. Recalling that $I$ is saturated, $x^{\beta}$ belongs to $I_q$.
\end{proof}

\begin{lemma}\label{lem:notSegment} 
Let $\mathscr{B} \subset \pos{n}{m}$ be a Borel set. If there exists four terms $x^{\alpha}, x^{\beta}\in \mathscr{B}$, $x^\gamma, x^\delta \in \mathscr{B}^{\mathcal{C}}$  such that $x^{\alpha}x^{\beta}=x^{\gamma}x^{\delta}$, then $\mathscr{B}$ is not a segment w.r.t. any term order $\sigma$.
\end{lemma}
\begin{proof}
If $\mathscr{B}$ were a segment w.r.t some $\sigma$, by the given assumptions we would have in particular 
$x^\alpha >_\sigma x^{\gamma}$ and $x^\beta >_{\sigma} x^{\delta}$. From these it would follow 
\[
x^\alpha x^\beta >_\sigma x^{\gamma} x^\beta,\ \quad x^\beta x^\gamma >_{\sigma} x^{\delta} x^\gamma \quad\Rightarrow\quad x^\alpha x^\beta >_\sigma x^\gamma x^{\delta} 
\]
contradicting $x^{\alpha}x^{\beta}=x^{\gamma}x^{\delta}$.
\end{proof}

Clearly this lemma can be used to deduce properties also about ideals. Let $I$ be a Borel-fixed ideal.
\begin{itemize}
\item If $\{I_t\}$ for some $t$ realizes the hypothesis of Lemma \ref{lem:notSegment}, then $I$ could not be a segment ideal.
\item If $\{I_r\}$ realizes the hypothesis of Lemma \ref{lem:notSegment}, then $I$ could not be a hilb-segment ideal.
\item If $\{I_{\reg(I)}\}$ realizes the hypothesis of Lemma \ref{lem:notSegment}, then $I$ could not be a reg-segment ideal.
\item If $\{I_t\}$, for some $t$, and two generators of degree $t$ realize the hypothesis of Lemma \ref{lem:notSegment}, then $I$ could not be a gen-segment ideal.
\end{itemize}

\begin{proposition}\label{prop:segmentsImplications}
Let $I\subset \K[x]$ be a saturated Borel-fixed ideal and let $\sigma$ be a term ordering. Then
\begin{enumerate}[(i)]
\item\label{it:segmentsImplications_i} $I$ segment ideal $\Rightarrow$ $I$ hilb-segment ideal $\Rightarrow$ $I$ reg-segment ideal $\Rightarrow$ $I$ gen-segment ideal. 
\item\label{it:segmentsImplications_ii} $\sigma$ is the lexicographic order $\Leftrightarrow$ the implications in (\ref{it:segmentsImplications_i}) are all equivalences, for every ideal $I$.
\item\label{it:segmentsImplications_iii} If the projective scheme defined by $I$ has constant Hilbert polynomial, then: $I$ segment ideal $\Leftrightarrow$ $I$ hilb-segment ideal $\Leftrightarrow$ $I$ reg-segment ideal.
\end{enumerate}
\end{proposition}
\begin{proof}
\emph{(\ref{it:segmentsImplications_i})} The first implication is obvious. For the second one, it is enough to apply Lemma \ref{lem:segmentLowerDegree}, because the Gotzmann number is greater than or equal to $reg(I)$. For the third implication, recall that $I$ is generated in degrees $\leqslant \reg(I)$, by definition. 
Moreover, if $I$ is a reg-segment ideal, by Lemma \ref{lem:segmentLowerDegree} $I_t$ contains the greatest terms of degree $t$, for every $t\leqslant \reg(I)$. Thus, in particular, minimal generators of $I$ must to be the greatest possible. 

\emph{(\ref{it:segmentsImplications_ii})} First, suppose that $\sigma$ is the lexicographic order. Then, by \emph{(\ref{it:segmentsImplications_i})}, it is enough to show that a gen-segment ideal is also a segment ideal. Indeed, by induction on the degree $s$ of monomials and with $s=0$ as base of induction, 
for $s>0$ suppose that $I_{s-1}$ is a segment. $\langle I_{s-1}\rangle_s$ is still a segment and, since possible minimal generators are always the greatest possible, we are done.

Vice versa, if $\sigma$ is not the lexicographic order, let $s$ be the minimum degree at which the monomials are ordered in  a different way from the lexicographic one. Thus, there exist two terms $x^{\alpha}$ and $x^{\beta}$ with maximum 
variables $x_l$ and $x_h$, respectively, such that $x^{\beta} <_{\sigma} x^{\alpha}$ but $x_h >_{\sigma} x_l$. The ideal 
$I=(x_h,\ldots,x_n)$ is a gen-segment ideal but not a segment ideal, since $x^{\beta}$ belongs to $I$ and 
$x^{\alpha}$ does not.

\emph{(\ref{it:segmentsImplications_iii})} It is enough to show that, in the case of constant Hilbert polynomial, a reg-segment ideal $I$ is also a segment ideal. 
By induction on the degree $s$, if $s\leqslant \reg(I)$, then the thesis follows by the hypothesis and by Lemma 
\ref{lem:segmentLowerDegree}. Suppose that $s>\reg(I)$ and that $I_{s-1}$ is a segment. At degree $s$ there are not minimal generators for $I$ so that a monomial of $I_s$ is always of type $x^{\alpha}x_h$ with $x^{\alpha}$ in $I_{s-1}$. 
Let $x^{\beta}$ be a monomial of degree $s$ such that $x^{\beta}>_\sigma x^{\alpha}x_h$, thus $x^{\beta}>_\sigma x^{\alpha}x_0$. 
By Proposition \ref{prop:degreeHPpowerVars}, we have that $(x_1,\ldots,x_n)^s \subseteq I$. So, if $x^\beta$ is not divided by $x_0$, then $x^\beta$ belongs to $I_s$, otherwise there exists a monomial $x^\gamma$ such that $x^\beta=x^\gamma x_0$. Thus $x^\gamma >_\sigma x^\alpha$ and by induction $x^\gamma$ belongs to $I_{s-1}$ so that $x^\beta=x^\gamma x_0$ belongs to $I_s$.
\end{proof}

In the following chapters, the relevance of the definition of various type of segment ideals will be clearified. At this point, we are interested in explaining how to determine the term ordering $\sigma$ that makes an ideal a segment ideal.

First of all, we recall the characterization of term orderings by means of matrices with rational coefficients.
Let $T \in \GL_{\QQ}(n+1)$ be any invertible matrix. $T$ induces an order relations on the monomials of $\K[x]$ defined as follows
\[
x^{\alpha} >_T x^{\beta} \quad \Longleftrightarrow\quad \text{the first non-zero entry of the vector $T\cdot(\alpha-\beta)^{\text{tr}}$ is positive}. 
\]
\begin{proposition}[{\cite[Proposition 1.4.12]{KreuzerRobbiano1}}]\index{term ordering}
Let $T$ be an invertible matrix in $\GL_{\QQ}(n+1)$. The order induced by $T$ on the monomials of $\K[x]$ is a term ordering if and only if the first non-zero element in each column of $T$ is positive.
\end{proposition}

As well known the matrices representing the \gls{Lex}, \gls{DegLex} and \gls{DegRevLex} are
\[
\left(\begin{array}{cccc}
1 & 0 & \ldots &0 \\
0 & 1 & \ldots  &0 \\
0 & \ddots &\ddots  &0\\
0 & \ldots & 0 & 1
\end{array}\right),\quad
\left(\begin{array}{cccc}
1 & 1 & \ldots &1 \\
1 & 0 & \ldots  &0 \\
0 & \ddots &\ddots  &0\\
0 & \ldots & 1 & 0
\end{array}\right)\quad \text{and}\quad
\left(\begin{array}{cccc}
1 & 1 & \ldots &1 \\
0 & \ldots & 0  &-1 \\
0 & \revddots & \revddots &0\\
0 & -1 & \ldots & 0
\end{array}\right).
\]

We are working in the projective case, i.e. with homogeneous polynomials, so we are interested in matrices with the first rows represented by the vector $(1,\ldots,1)$ to evaluate the total degree of a monomial. We consider as second row a vector $\omega = (\omega_n,\ldots,\omega_0) \in \QQ^{n+1}$ and then we complete the matrix with the same rows of the matrix associated to the \texttt{DegLex} term ordering:
\begin{equation}\label{eq:termOrdering}
T = \left(
\begin{array}{ccccc}
1 & 1 & \ldots & 1 & 1\\
\omega_n & \omega_{n-1} & \ldots & \omega_1 & \omega_0\\
0 & 1 & 0 & \ldots & 0\\
\vdots & \ddots & \ddots & \ddots & \vdots\\
0 & \ldots & 0 & 1 & 0\\
\end{array}
\right).
\end{equation}
We want that the induced term ordering agrees with the hypothesis $x_n >_B \ldots >_B x_0$, so we will always suppose $\omega_n > \ldots > \omega_0$. This fact ensures that the matrix is invertible because the rank of $T$ could be lower than $n+1$ if and only if
\[
\left\vert\begin{array}{cc}
1 & 1\\
\omega_n & \omega_0
\end{array}
\right\vert = \omega_0 - \omega_n = 0.
\]
Since the term ordering induced by $T$ depends on the vector $\omega$, we will say \lq\lq the term ordering $\omega$\rq\rq\ meaning \lq\lq the term ordering induced by the matrix \eqref{eq:termOrdering} with second row equal to the vector $\omega$\rq\rq.

Now let us consider a Borel set $\mathscr{B} \subset \pos{n}{m}$ and let us try to determine the vector $\omega$ in order for $\mathscr{B}$ to be a segment. If $x^{\alpha_1},\ldots,x^{\alpha_a}$ are the minimal monomials of $\mathscr{B}$ and $x^{\beta_1},\ldots,x^{\beta_b}$ are the maximal elements of the complement $\mathscr{B}^{\mathcal{C}}$, we need to impose that $x^{\alpha_i} >_\omega x^{\beta_j},\ \forall\ i = 1,\ldots,a,\ j=1,\ldots,b$. 
The first row of the matrix \eqref{eq:termOrdering} compares the degree of the monomials, therefore it does not affect the ordering on the monomials of $\pos{n}{m}$.
We can look for a vector $\omega$ that orders the monomials as we want, that is we have to solve the following system of inequalitites
\begin{equation}\label{eq:inequalitites}
\begin{cases}
\omega_i > \omega_{i-1}, & i=1,\ldots,n, \\
\omega \cdot (\alpha_i - \beta_j) > 0, & \forall\ x^{\alpha_i},\ \forall\ x^{\beta_j}.
\end{cases}
\end{equation}

\begin{example}\label{ex:segmentInequalities}
Let us consider the Borel set $\{I_4\} \subset \pos{2}{4}$ defined by the ideal $I = (x_2^2,x_2 x_1^2,x_1^4) \subset \K[x_0,x_1,x_2]$. The minimal monomials of $\{I_4\}$ are $x_2^2 x_0^2,x_2x_1^2 x_0,x_1^4$ and the maximal monomials of $\{I_4\}^{\mathcal{C}}$ are $x_2x_1x_0^2, x_1^3x_0$. So the system of inequalities to solve is
\[
\left\{
\begin{array}{ll}
\omega_2 > \omega_1 & (x_2 >_\omega x_1)\\
\omega_1 > \omega_0 & (x_1 >_\omega x_0)\\
2\omega_2 + 2\omega_0 > \omega_2 + \omega_1 + 2\omega_0 & (x_2^2 x_0^2 >_{\omega} x_2x_1x_0^2)\\
2\omega_2 + 2\omega_0 > 3\omega_1 + \omega_0 & (x_2^2 x_0^2 >_{\omega} x_1^3x_0)\\
\omega_2 + 2\omega_1 + \omega_0 > \omega_2 + \omega_1 + 2\omega_0 & (x_2 x_1^2 x_0 >_{\omega} x_2x_1x_0^2)\\
\omega_2 + 2\omega_1 + \omega_0 > 3\omega_1 + \omega_0 & (x_2 x_1^2 x_0 >_{\omega} x_1^3x_0)\\
4\omega_1 > \omega_2 + \omega_1 + 2\omega_0 & (x_1^4 >_{\omega} x_2x_1x_0^2)\\
4\omega_1 > 3\omega_1 + \omega_0 & (x_1^4 >_{\omega} x_1^3x_0)\\
\end{array}
\right.
\
\Rightarrow\
\left\{
\begin{array}{l}
\omega_2 > \omega_1 \\
\omega_1 > \omega_0 \\
\omega_2 > \omega_1 \\
2\omega_2 + \omega_0 > 3\omega_1 \\
\omega_1 > \omega_0 \\
\omega_2 > \omega_1 \\
3\omega_1 > \omega_2 + 2\omega_0 \\
\omega_1 >  \omega_0 \\
\end{array}
\right.
\]
Supposing $\omega_0 = 1$, we obtain
\[
\left\{
\begin{array}{l}
\omega_2 > \omega_1 \\
\omega_1 > 1 \\
2\omega_2 > 3\omega_1 - 1\\
\omega_2 < 3\omega_1 + 2 \\
\end{array}
\right.
\quad\Rightarrow\quad \omega_2 = 3,\ \omega_1 = 2,\ \omega_0 = 1.
\]
\end{example}

The example clearly shows that whenever for the pair $x^{\alpha_i},x^{\beta_j}$ there exists an elementary decreasing move $\down{k}$ such that $\down{k}(x^{\alpha_i}) = x^{\beta_j}$, the inequality given is already determined by the Borel order. Thus we can restrict the conditions to be imposed as follows
\begin{equation}\label{eq:inequalititesRestricted}
\begin{cases}
\omega_i > \omega_{i-1}, & i=1,\ldots,n, \\
\omega \cdot (\alpha_i - \beta_j) > 0, & \forall\ x^{\alpha_i},\ \forall\ x^{\beta_j},\ \nexists\ \down{k} \text{ s.t. } \down{k}(x^{\alpha_i}) = x^{\beta_j}.
\end{cases}
\end{equation}

For the simplest cases, a solution of the system of inequalitites can be found by hand with a little bit of work, but just increasing the number of variables and the degree of monomials of a poset, the system can become much more complicated.  A natural way to face this problem is to use the theory of linear programming and the simplex algorithm (see \cite[Chapter 2]{Vanderbei} for an introduction to the topic).\index{simplex algorithm}

The problem of linear programming that best fits our case is the minimization of a linear expression under some linear constraints. In fact the standard problem of minimization requires:
\begin{itemize}
\item a set of variables $z_0,\ldots,z_n$ with the hypothesis that they do not assume negative values, i.e.
\[
z_0 \geqslant0,\quad \ldots \quad, z_n \geqslant 0,
\]
\item a target function (that we want to minimize)
\[
c_0 z_0 + \ldots + c_n z_n,
\]
\item a set of constraints on the variables expressed by linear inequalities of the type
\[
a_{i0} z_0 + \ldots + a_{in}z_n \geqslant b_i.
\]
\end{itemize}
The simplex algorithm allows to compute the solution $(\overline{z}_0,\ldots,\overline{z}_n)$ of the system of constraints such that the value $c_0 \overline{z}_0 + \ldots + c_n \overline{z}_n$ is minimum.

In our context we are looking for any solution of the system \eqref{eq:inequalititesRestricted}, not the minimal one, so we choose arbitrarily as target function the sum of the vector components $\omega_0 + \ldots + \omega_n$. The next step is to transform our strict inequalities in non-strict ones, so we rewrite the system \eqref{eq:inequalititesRestricted}  as
\begin{equation}\label{eq:inequalititesRestrictedSimplex}
\begin{cases}
\omega_0 \geqslant 1\\ 
\omega_i - \omega_{i-1} \geqslant 1, & i=1,\ldots,n, \\
\omega \cdot (\alpha_i - \beta_j) \geqslant 1, & \forall\ x^{\alpha_i},\ \forall\ x^{\beta_j},\ \nexists\ \down{k} \text{ s.t. } \down{k}(x^{\alpha_i}) = x^{\beta_j}.
\end{cases}
\end{equation}
ensuring also $\omega_0 \geqslant 0,\ldots,\omega_n \geqslant 0$.

\begin{example}
Considering again the Borel set $\{I_4\}$ introduced in Example \ref{eq:inequalititesRestrictedSimplex}, the system of inequalities written so that it can be solved with the simplex algorithm is
\[
\left\{
\begin{array}{l}
\omega_0 \geqslant 1 \\
\omega_1 - \omega_0 \geqslant 1 \\
\omega_2 - \omega_1 \geqslant 1 \\
2\omega_2 - 3\omega_1 + \omega_0 \geqslant 1 \\
- \omega_2 + 3\omega_1 - 2 \omega_0 \geqslant 1\\
\end{array}
\right.
\]
and the solution we obtain with target fucntion $\omega_2+\omega_1+\omega_0$ is again $(3,2,1)$.
\end{example}

In Algorithm \ref{alg:BorelSetSegment}, there is the pseudocode description of the strategy just exposed\hfill and\hfill also\hfill used\hfill for\hfill designing\hfill the\hfill corresponding\hfill functions\hfill of\hfill the\hfill class\\ \texttt{BorelInequalitiesSystem} of the package \texttt{HSC} (see Appendix \ref{ch:HSCpackage}). In Algorithm \ref{alg:HilbRegSegment}, there are the pseudocode descriptions of the methods for determining whenever a Borel-fixed ideal is a hilb-segment or a reg-segment ideal, that can be easily deduced by Algorithm \ref{alg:BorelSetSegment}. 

\begin{algorithm}[H]
\caption*{Some auxiliary methods for Algorithm \ref{alg:BorelSetSegment}.}
\begin{algorithmic}
\STATE \textsc{MinimalElements}($\mathscr{B}$)
\REQUIRE $\mathscr{B}$, Borel set.
\ENSURE the set of minimal elements of $\mathscr{B}$ w.r.t. $\leq_B$.
\end{algorithmic}

\medskip

\begin{algorithmic}
\STATE $\textsc{MaximalElements}(\mathscr{B})$
\REQUIRE $\mathscr{B}$, Borel set.
\ENSURE the set of maximal elements of $\mathscr{B}$ w.r.t. $\leq_B$.
\end{algorithmic}

\medskip

\begin{algorithmic}
\STATE $\textsc{SimplexAlgorithm}(f,S,\omega)$
\REQUIRE $f$, a target function (to be minimized).
\REQUIRE $S$, a set of constraints.
\REQUIRE $\omega$, a vector.
\ENSURE \TRUE\ if the system of constraints is solvable, \FALSE\ otherwise. If \TRUE\, the solution that minimizes $f$ is stored in $\omega$.
\end{algorithmic}
\end{algorithm}

\begin{algorithm}[!ht]
\caption[Algorithm computing a vector defining a term ordering that realizes a Borel set as segment.]{The algorithm computing a vector defining a term ordering that realizes a Borel set as segment.}
\label{alg:BorelSetSegment}
\begin{algorithmic}[1]
\STATE $\textsc{isSegment}(\mathscr{B},\omega)$
\REQUIRE $\mathscr{B} \subset \pos{n}{m}$, Borel set.
\REQUIRE $\omega = (\omega_0,\ldots,\omega_n)$, a vector.
\ENSURE \FALSE, if $\mathscr{B}$ can not be a segment, \TRUE\ otherwise. In this case the term ordering realizing $\mathscr{B}$ as a segment is stored in $\omega$.
\STATE $\textsf{constraints} \leftarrow \{z_0 \geqslant 1\}$;
\FOR{$i = 1,\ldots,n$}
\STATE $\textsf{constraints} \leftarrow \textsf{constraints} \cup \{z_i - z_{i-1} \geqslant 1\}$;
\ENDFOR
\STATE $\textsf{minimalMonomials} \leftarrow \textsc{MinimalElements}(\mathscr{B})$;
\STATE $\textsf{maximalMonomials} \leftarrow \textsc{MaximalElements}(\mathscr{B})$;
\FORALL{$x^{\alpha}\in \textsf{minimalMonomials}$}
\FORALL{$x^\beta \in \textsf{maximalMonomials}$}
\IF{$\nexists\ \down{k} \text{ s.t. } \down{k}(x^{\alpha}) = x^\beta$}
\STATE $\textsf{constraints} \leftarrow \textsf{constraints} \cup \{(z_n,\ldots,z_0)\cdot(\alpha-\beta) \geqslant 1\}$;
\ENDIF
\ENDFOR
\ENDFOR
\RETURN $\textsc{SimplexAlgorithm}(z_0 + \ldots + z_n,\textsf{constraints},\omega)$;
\end{algorithmic}
\end{algorithm}

In the following proposition we use definitions and algorithms just introduced to characterize segment ideals in $\K[x_0,x_1,x_2]$ with constant Hilbert polynomial.

\begin{algorithm}[!ht]
\caption[Methods for computing the term ordering that realizes an Borel-fixed ideal as hilb-segment or reg-segment.]{The methods for computing the term ordering that realizes an Borel-fixed ideal as hilb-segment or reg-segment.}
\label{alg:HilbRegSegment}
\begin{algorithmic}[1]
\STATE $\textsc{isHilbSegment}(I,\omega)$
\REQUIRE $I \subset \K[x_0,\ldots,x_n]$, Borel-fixed ideal.
\REQUIRE $\omega = (\omega_0,\ldots,\omega_n)$, a vector.
\ENSURE \FALSE, if $I$ can not be a hilb-segment, \TRUE\ otherwise. In this case the term ordering realizing $I$ as a hilb-segment is stored in $\omega$.
\STATE $p(t) \leftarrow \textsc{HilbertPolynomial}(\K[x_0,\ldots,x_n]/I)$;
\STATE $r \leftarrow \textsc{GotzmannNumber}\big(p(t)\big)$;
\RETURN $\textsc{isSegment}(\{I_r\},\omega)$;
\end{algorithmic}

\medskip

\begin{algorithmic}[1]
\STATE $\textsc{isRegSegment}(I,\omega)$
\REQUIRE $I \subset \K[x_0,\ldots,x_n]$, Borel-fixed ideal.
\REQUIRE $\omega = (\omega_0,\ldots,\omega_n)$, a vector.
\ENSURE \FALSE, if $I$ can not be a reg-segment, \TRUE\ otherwise. In this case the term ordering realizing $I$ as a reg-segment is stored in $\omega$.
\RETURN $\textsc{isSegment}(\{I_{\reg(I)}\},\omega)$;
\end{algorithmic}
\end{algorithm}

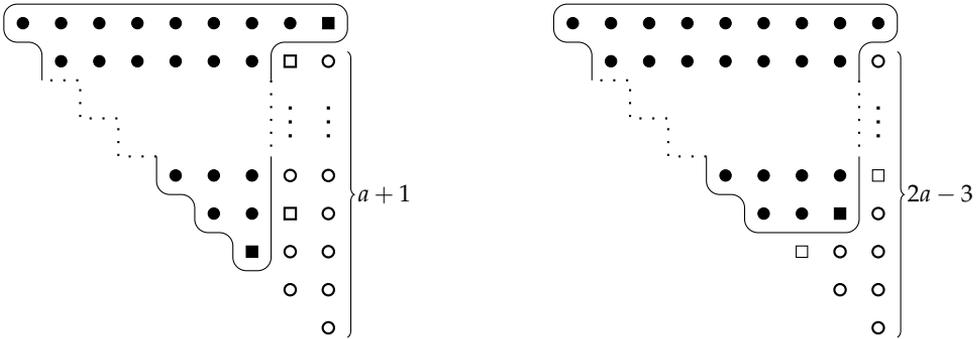
\begin{figure}[!ht]
\begin{center}
\captionsetup[subfloat]{singlelinecheck=false,format=hang}
\subfloat[][The Borel set of $\pos{2}{a+1}$ defined by $(x_2^2,x_2x_1^{a},x_1^{a+1})$.]{\label{fig:notSegmentP2_a}
\begin{tikzpicture}[scale=0.5,decoration=brace]
\tikzstyle{ideal}=[circle,draw=black,fill=black,inner sep=1.5pt]
\tikzstyle{quotient}=[circle,draw=black,thick,inner sep=1.5pt]

\node at (-4,4) [ideal] {};
\node at (-3,4) [ideal] {};
\node at (-2,4) [ideal] {};
\node (n)  at (-1,4) [ideal] {};
\node (n0) at (0,4) [ideal] {};
\node (n1) at (1,4) [ideal] {};
\node (n2) at (2,4) [ideal] {};
\node (n3) at (3,4) [ideal] {};
\node (n4) at (4,4) [regular polygon,regular polygon sides=4,draw=black,fill=black,inner sep=1.5pt] {};

\node at (-3,3) [ideal] {};
\node at (-2,3) [ideal] {};
\node (n)  at (-1,3) [ideal] {};
\node (n0) at (0,3) [ideal] {};
\node (n1) at (1,3) [ideal] {};
\node (n2) at (2,3) [ideal] {};
\node (n3) at (3,3) [regular polygon,regular polygon sides=4,draw=black,thick,inner sep=1.5pt] {};
\node (n4) at (4,3) [quotient] {};

\node (00) at (0,0) [ideal] {};
\node (01) at (1,0) [ideal] {};
\node (02) at (2,0) [ideal] {};
\node (03) at (3,0) [quotient] {};
\node (04) at (4,0) [quotient] {};

\node (11) at (1,-1) [ideal] {};
\node (12) at (2,-1) [ideal] {};
\node (13) at (3,-1) [regular polygon,regular polygon sides=4,draw=black,thick,inner sep=1.5pt] {};
\node (14) at (4,-1) [quotient] {};

\node (22) at (2,-2) [regular polygon,regular polygon sides=4,draw=black,fill=black,inner sep=1.5pt] {};
\node (23) at (3,-2) [quotient] {};
\node (24) at (4,-2) [quotient] {};

\node (33) at (3,-3) [quotient] {};
\node (34) at (4,-3) [quotient] {};

\node (44) at (4,-4) [quotient] {};

\draw [loosely dotted,very thick] (4,1) -- (4,2); 
\draw [loosely dotted,very thick] (3,1) -- (3,2); 

\draw [rounded corners=5pt] (-3.5,2.5) -- (-3.5,3.5) -- (-4.5,3.5) -- (-4.5,4.5) -- (4.5,4.5) -- (4.5,3.5) -- (2.5,3.5) -- (2.5,2.5);
\draw [loosely dotted,thick] (2.5,2.5) -- (2.5,0.5);
\draw [loosely dotted,thick] (-0.5,0.5) -- (-1.5,0.5) -- (-1.5,1.5) -- (-2.5,1.5) -- (-2.5,2.5) -- (-3.5,2.5);
\draw [rounded corners=5pt] (2.5,0.5) -- (2.5,-2.5) -- (1.5,-2.5) -- (1.5,-1.5) -- (0.5,-1.5) -- (0.5,-0.5) -- (-0.5,-0.5) -- (-0.5,0.5);

\draw [decorate] (4.5,3.25) --node[right]{\footnotesize $a+1$} (4.5,-4.25);
\end{tikzpicture}
}
\qquad\qquad
\subfloat[][The Borel set of $\pos{2}{2a-3}$ defined by $(x_2^3,x_2^2x_1,x_2x_1^2, x_1^{2a-3})$.]{\label{fig:notSegmentP2_b}
\begin{tikzpicture}[scale=0.5,decoration=brace]
\tikzstyle{ideal}=[circle,draw=black,fill=black,inner sep=1.5pt]
\tikzstyle{quotient}=[circle,draw=black,thick,inner sep=1.5pt]

\node at (-4,4) [ideal] {};
\node at (-3,4) [ideal] {};
\node at (-2,4) [ideal] {};
\node (n)  at (-1,4) [ideal] {};
\node (n0) at (0,4) [ideal] {};
\node (n1) at (1,4) [ideal] {};
\node (n2) at (2,4) [ideal] {};
\node (n3) at (3,4) [ideal] {};
\node (n4) at (4,4) [ideal] {};

\node at (-3,3) [ideal] {};
\node at (-2,3) [ideal] {};
\node (n)  at (-1,3) [ideal] {};
\node (n0) at (0,3) [ideal] {};
\node (n1) at (1,3) [ideal] {};
\node (n2) at (2,3) [ideal] {};
\node (n3) at (3,3) [ideal] {};
\node (n4) at (4,3) [quotient] {};

\node (00) at (0,0) [ideal] {};
\node (01) at (1,0) [ideal] {};
\node (02) at (2,0) [ideal] {};
\node (03) at (3,0) [ideal]{};
\node (04) at (4,0) [regular polygon,regular polygon sides=4,draw=black,inner sep=1.5pt] {};

\node (11) at (1,-1) [ideal] {};
\node (12) at (2,-1) [ideal] {};
\node (13) at (3,-1) [regular polygon,regular polygon sides=4,draw=black,fill=black,inner sep=1.5pt] {};
\node (14) at (4,-1) [quotient] {};

\node (22) at (2,-2) [regular polygon,regular polygon sides=4,draw=black,inner sep=1.5pt] {};
\node (23) at (3,-2) [quotient] {};
\node (24) at (4,-2) [quotient] {};

\node (33) at (3,-3) [quotient] {};
\node (34) at (4,-3) [quotient] {};

\node (44) at (4,-4) [quotient] {};

\draw [loosely dotted,very thick] (4,1) -- (4,2); 

\draw [rounded corners=5pt] (-3.5,2.5) -- (-3.5,3.5) -- (-4.5,3.5) -- (-4.5,4.5) -- (4.5,4.5) -- (4.5,3.5) -- (3.5,3.5) -- (3.5,2.5);
\draw [loosely dotted,thick] (3.5,2.5) -- (3.5,0.5);
\draw [loosely dotted,thick] (-0.5,0.5) -- (-1.5,0.5) -- (-1.5,1.5) -- (-2.5,1.5) -- (-2.5,2.5) -- (-3.5,2.5);
\draw [rounded corners=5pt] (3.5,0.5) -- (3.5,-1.5) -- (1.5,-1.5) -- (0.5,-1.5) -- (0.5,-0.5) -- (-0.5,-0.5) -- (-0.5,0.5);

\draw [decorate] (4.5,3.25) --node[right]{\footnotesize $2a-3$} (4.5,-4.25);
\end{tikzpicture}
}
\caption[{The Borel sets defined by Borel-fixed ideals of $s \geqslant 7$ points in $\PP^2$, which can not be segment ideals.}]{\label{fig:notSegmentP2} The Borel sets defined by Borel-fixed ideals of $\K[x_0,x_1,x_2]$ with Hilbert polynomial $p(t) = s \geqslant 7$, which can not be segment ideals (see Proposition \ref{prop:pointsP2Segments}). With $\square$ and $\blacksquare$ we denoted the monomials realizing the hypothesis of Lemma \ref{lem:notSegment}.}
\end{center}
\end{figure}

\begin{proposition}\label{prop:pointsP2Segments}
In $\K[x_0,x_1,x_2]$ every saturated Borel-fixed ideal with Hilbert polynomial $p(t)=s \leqslant 6$ is a segment 
ideal. Whereas for every $p(t) = s \geqslant 7$, a saturated Borel-fixed ideal, which is not a segment
for any term order, always exists. 
\end{proposition}  
\begin{proof}
We give a constructive proof of this result, examining the Borel sets defined by the ideals in degree equal to their regularity (Proposition \ref{prop:segmentsImplications}\emph{\ref{it:segmentsImplications_iii}}). 
\begin{description}
\item[$s=1,2.$] There exists a unique saturated Borel-fixed ideal $(x_2,x_1^d)$, which is the lexicographic ideal.
\item[$s=3,4.$] There are two saturated Borel-fixed ideals: the lexsegment ideal $(x_2,x_1^d)$ and $(x_2^2,x_2x_1,x_1^{d-1})$, segment w.r.t. \texttt{DegRevLex}.
\item[$s=5.$] There are three saturated Borel-fixed ideals: the lexsegment ideal $(x_2,x_1^5)$, $(x_2^2,x_2x_1,x_1^4)$\hfill segment\hfill w.r.t.\hfill $(4,2,1)$\hfill and\hfill $(x_2^2,x_2 x_1^2,x_1^3)$\hfill segment\hfill w.r.t.\\ \texttt{DegRevLex}.
\item[$s=6.$] There are four saturated Borel ideals: the lexsegment $(x_2,x_1^6)$,  $(x_2^2,x_2 x_1,x_1^5)$ segment w.r.t. $(5,2,1)$, $(x_2^2,x_2x_1^2,x_1^4)$ segment w.r.t. $(3,2,1)$ and $(x_2^3,x_2^2x_1,x_2x_1^2,$ $x_1^3)$ segment w.r.t. \texttt{DegRevLex}.
\item[$s = 2a+1,\ a \geqslant 3.$] Let us consider the ideal $J = (x_2^2,x_2 x_1^a,x_1^{a+1})$. It has constant Hilbert\hfill polynomial\hfill $p(t)=2a+1$,\hfill because\hfill $x_1^{a+1}$\hfill belongs\hfill to\hfill $J$\hfill and\\ $\dim_{\K} \K[x_0,x_1,x_2]_{a+1}/J_{a+1} = 2a+1$, indeed the monomials of degree $a+1$ not belonging to $J$ are
\[
\left\{x_0^{a+1},x_1 x_0^{a},\ldots,x_1^{a} x_{0}\right\} \cup \left\{ x_2x_0^{a}, x_2x_1x_0^{a-1},\ldots,x_2x_1^{a-1}x_0\right\}.
\]
By Lemma \ref{lem:notSegment}, $J$ cannot be a segment because $x_2^2 x_0^{a-1},x_1^{a+1} \in \{J_{a+1}\}$, $x_2x_1^2 x_0^{a-2},$ $x_2x_1^{a-1}x_0 \in \{J_{a+1}\}^{\mathcal{C}}$ and $x_2^2 x_0^{a-1} \cdot x_1^{a+1} = x_2x_1^2 x_0^{a-2} \cdot x_2x_1^{a-1}x_0$ (see Figure \ref{fig:notSegmentP2_a}).
\item[$s = 2a,\ a \geqslant 4.$] The ideal $J = (x_2^3,x_2^2 x_1,x_2 x_1^2, x_1^{2a-3})$ has Hilbert polynomial $p(t)=2a$, because $x_1^{2a-3} \in J$ and $\K[x_0,x_1,x_2]_{2a-3}/J_{2a-3}$ is spanned by the $2a$ monomials
\[
\left\{x_2^2 x_0^{2a-5},x_2x_1x_0^{2a-5},x_2 x_0^{2a-4}\right\}\cup\left\{x_0^{2a-3},x_1 x_0^{2a-4},\ldots,x_1^{2a-4}x_0\right\}.
\] 
Finally by Lemma \ref{lem:notSegment}, $J$ can not be a segment because $x_2 x_1^2 x_0^{2a-6} \in \{J_{2a-3}\}$, $x_2^2 x_0^{2a-5},$ $x_1^4 x_0^{2a-7} \in \{J_{2a-3}\}^{\mathcal{C}}$ and $(x_2 x_1^2 x_0^{2a-6})^2 = x_2^2 x_0^{2a-5} \cdot x_1^4 x_0^{2a-7}$ (see Figure \ref{fig:notSegmentP2_b}).\qedhere
\end{description}
\end{proof}

A slightly different approach is needed to deal with the gen-segment ideals, because the property required in the definition involes more than a single degree of the monomials. Let us consider the truncation $I_{\geqslant m}$ of a saturated Borel-fixed ideal $I$. In order for $I_{\geqslant m}$ to be a gen-segment ideal, we have to check that $\{I_m\}$ is a segment and that any generator of $I$ of degree $\overline{m} > m$ is bigger than the monomials of $\{I_{\overline{m}}\}^{\mathcal{C}}$:
\begin{equation}\label{eq:inequalititesGenSegment}
\begin{cases}
\omega_0 \geqslant 1, & \\
\omega_i - \omega_{i-1} \geqslant 1, & i=1,\ldots,n, \\
\omega \cdot (\alpha - \beta) \geqslant 1, &
\begin{cases}
\forall\ x^{\alpha} \text{ minimal monomial of } \{J_m\} \\ 
\forall\ x^{\beta} \text{ maximal monomial of } \{J_m\}^{\mathcal{C}}\\
\nexists\ \down{k} \text{ s.t. } \down{k}(x^{\alpha}) = x^{\beta}
\end{cases} \\
\omega \cdot (\gamma - \delta) \geqslant 1, &
\begin{cases}
\forall\ x^{\gamma} \text{ minimal generator of } I,\ \vert\gamma\vert > m \\ 
\forall\ x^{\delta} \text{ maximal monomial of } \{I_{\vert\alpha\vert}\}^{\mathcal{C}}\\
\nexists\ \down{k} \text{ s.t. } \down{k}(x^{\gamma}) = x^{\delta}
\end{cases} \\
\end{cases}
\end{equation}

\begin{example}\label{ex:notRegbutGen}
Let us consider the ideal $I = (x_3^3,x_3^2x_2,x_3x_2^2,x_3^2x_1,x_3x_2x_1,x_3x_1^2,x_2^5) \subset \K[x_0,x_1,x_2,x_3]$. Note that $I$ can not be a reg-segment ideal because
\[
x_3x_2^2 x_0^2 \in I,\ x_2^4 x_0,x_3^2x_0^3 \notin I\quad \text{and}\quad \left(x_3x_2^2 x_0^2\right)^2 = x_2^4 x_0 \cdot x_3^2x_0^3 \qquad \text{(Lemma \ref{lem:notSegment})}
\]
neither a hilb-segment ideal because in degree $9$ (the Gotzmann number of the Hilbert polynomial $5t-1$ of $I$) the same relation multiplied by $x_0^4$ holds.

Applying Algorithm \ref{alg:genSegment} to determine if $I$ could be a gen-segment ideal, we need to impose that $\{I_3\}$ is a segment of $\pos{3}{3}$ and that $x_2^5$ is greater than all the monomials of $\{I_5\}^{\mathcal{C}}$. These conditions lead to the following system of inequalities:
\[
\begin{cases}
\omega_0 \geqslant 1 & \\
\omega_1 - \omega_0 \geqslant 1 & \\
\omega_2 - \omega_1 \geqslant 1 & \\
\omega_3 - \omega_2 \geqslant 1 & \\
\omega_3 - 3\omega_2 + 2\omega_1 \geqslant 1 & (x_3x_1^2 >_{\omega} x_2^3)\\
-\omega_3 2\omega_1 -\omega_0 \geqslant 1 & (x_3x_1^2 >_{\omega} x_3^2 x_0)\\
-2\omega_3 + 5\omega_2 -3\omega_0 \geqslant 1 & (x_2^5 >_{\omega} x_3^2 x_0^3)
\end{cases} \qquad\Rightarrow\qquad \omega = (10,7,6,1).
\]
\end{example}

\begin{algorithm}[!ht]
\caption[Pseudocode description of the algorithm determining if a Borel-fixed ideal is a gen-segment ideal or not.]{The pseudocode description of the algorithm determining if a Borel-fixed ideal is a gen-segment ideal or not.}
\label{alg:genSegment}
\begin{algorithmic}[1]
\STATE $\textsc{isGenSegment}(I,m,\omega)$
\REQUIRE $I \subset \K[x_0,\ldots,x_n]$, saturated Borel-fixed ideal.
\REQUIRE $m$, positive integer.
\REQUIRE $\omega = (\omega_0,\ldots,\omega_n)$, a vector.
\ENSURE \FALSE, if $I_{\geqslant m}$ can not be a gen-segment, \TRUE\ otherwise. In this case the term ordering realizing $I_{\geqslant m}$ as a gen-segment ideal is stored in $\omega$.
\STATE $\textsf{constraints} \leftarrow \{z_0 \geqslant 1\}$;
\FOR{$i = 1,\ldots,n$}
\STATE $\textsf{constraints} \leftarrow \textsf{constraints} \cup \{z_i - z_{i-1} \geqslant 1\}$;
\ENDFOR
\STATE $\textsf{minimalMonomials} \leftarrow \textsc{minimalElements}\big(\{I_m\}\big)$;
\STATE $\textsf{maximalMonomials} \leftarrow \textsc{maximalElements}\big(\{I_m\}^{\mathcal{C}}\big)$;
\FORALL{$x^{\alpha}\in \textsf{minimalMonomials}$}
\FORALL{$x^\beta \in \textsf{maximalMonomials}$}
\IF{$\nexists\ \down{k} \text{ s.t. } \down{k}(x^{\alpha}) = x^\beta$}
\STATE $\textsf{constraints} \leftarrow \textsf{constraints} \cup \{(z_n,\ldots,z_0)\cdot(\alpha-\beta) \geqslant 1\}$;
\ENDIF
\ENDFOR
\ENDFOR
\FOR{$j=m+1,\ldots,\reg(I)$}
\STATE $\textsf{maximalMonomials} \leftarrow \textsf{maximalElements}\big(\{I_{j}\}\big)$;
\FORALL{$x^{\gamma}$ minimal generator of $I_{\geqslant m}$, $\vert\gamma\vert = j$}
\FORALL{$x^\delta \in \textsf{maximalMonomials}$}
\IF{$\nexists\ \down{k} \text{ s.t. } \down{k}(x^{\gamma}) = x^\delta$}
\STATE $\textsf{constraints} \leftarrow \textsf{constraints} \cup \{(z_n,\ldots,z_0)\cdot(\gamma-\delta) \geqslant 1\}$;
\ENDIF
\ENDFOR
\ENDFOR
\ENDFOR
\RETURN $\textsc{SimplexAlgorithm}(z_0 + \ldots + z_n,\textsf{constraints},\omega)$;
\end{algorithmic}
\end{algorithm}

\chapter{Rational curves on the Hilbert scheme}\label{ch:deformations}

In this chapter we expose the results contained in the preprint \cite{LellaDeformations} \lq\lq A network of rational curves on the Hilbert scheme\rq\rq.

\section{Rational deformations of Borel-fixed ideals}

The starting idea is expressed in the following remark.

\begin{remark}\label{rk:AddRemove}\index{Borel set}\index{minimal monomial}\index{maximal monomial}
Let $\mathscr{B} \subset \pos{n}{m}$ be a Borel set and let $x^\alpha$ and $x^\beta$ be a minimal monomial of $\mathscr{B}$ and a maximal monomial of $\mathscr{B}$. By definition both $\mathscr{B}\setminus\{x^\alpha\}$ and $\mathscr{B}\cup\{x^\beta\}$ are still Borel sets. Moreover $x^\alpha$ will be a maximal element of $(\mathscr{B}\setminus\{x^\alpha\})^{\mathcal{C}}$ and $x^\beta$ will be a minimal element of $\mathscr{B}\cup\{x^\beta\}$.
\end{remark}

Let us consider a Borel set $\{I_r\} \subset \pos{n}{r}$, defined by a Borel-fixed ideal $I$ with constant Hilbert polynomial $p(t) = r$, i.e. $\vert \{I_r\}^\mathcal{C}\vert = r$ and $\restrict{\{I_r\}^{\mathcal{C}}}{j} = \emptyset, \forall\ j > 0$. Moreover let us suppose that there exist two monomials $x^\alpha,x^\beta \in \pos{n}{r}$, such that $x^\alpha$ is a minimal element in $\{I_r\}$, $x^\beta$ is a maximal element in $\{I_r\}^{\mathcal{C}}$ and assume $\min x^{\alpha} = \min x^\beta = 0$. If $\mathscr{B} = \{I_r\}\setminus\{x^\alpha\} \cup \{x^\beta\}$ is still a Borel set, the ideal $\langle \mathscr{B}\rangle^{\sat}$ has the same Hilbert polynomial of $I$ (by Corollary \ref{cor:bijectionAllIdeals}).

\begin{example}\label{ex:ExampleSwap}
Let us consider the ideal $I = (x_2^2,x_2 x_1,x_1^3) \subset \K[x_0,x_1,x_2]$ having Hilbert polynomial $p(t) = 4$. The Borel set $\{I_4\} \subseteq \mathcal{P}(2,4)$ has two minimal monomials: $x_2x_1x_0^2$ and $x_1^3x_0$; while the monomials $x_2x_0^3$ and $x_1^2x_0^2$ are maximal elements in the complement $\{I_4\}^{\mathcal{C}}$ (see Figure \ref{fig:ExampleSwap_a}). There are four possibilities of removing a minimal element and adding a maximal one:
\begin{description}
\item[(Figure \ref{fig:ExampleSwap_b})] the set $\mathscr{B}_1 = \{I_4\} \setminus \{x_2x_1x_0^2\} \cup \{x_1^2x_0^2\}$ is not Borel, because $x_1^2 x_0^2 \in \mathscr{B}_1$ and $x_2 x_1 x_0^2 = \up{1}(x_1^2 x_0^2) \notin \mathscr{B}_1$;
\item[(Figure \ref{fig:ExampleSwap_c})] the set $\mathscr{B}_2 = \{I_4\} \setminus \{x_1^3x_0\} \cup \{x_1^2x_0^2\}$ is not Borel, because $x_1^2 x_0^2 \in \mathscr{B}_2$ and $x_1^3 x_0 = \up{0}(x_1^2 x_0^2) \notin \mathscr{B}_2$;
\item[(Figure \ref{fig:ExampleSwap_d})] the set $\mathscr{B}_3 = \{I_4\} \setminus \{x_2x_1x_0^2\} \cup \{x_2x_0^3\}$ is not Borel, because $x_2 x_0^3 \in \mathscr{B}_3$ and $x_2 x_1 x_0^2 = \up{0}(x_2 x_0^3) \notin \mathscr{B}_3$;
\item[(Figure \ref{fig:ExampleSwap_e})]  the set $\mathscr{B}_4 = \{I_4\} \setminus \{x_1^3x_0\} \cup \{x_2x_0^3\}$ is Borel and $\langle\mathscr{B}_4\rangle^{\sat} = (x_2,x_1^4)$.
\end{description}
\end{example}

By Remark \ref{rk:AddRemove} and Example \ref{ex:ExampleSwap}, we see that swapping two monomials $x^\alpha,x^\beta$ fails to preserve the Borel property whenever the minimal monomial is mapped by an elementary move to the maximal monomial, in fact in this case if we first remove the minimal monomial $x^\alpha$ from $\mathscr{B}$, $x^\beta$ is no longer a maximal element of $(\mathscr{B}\setminus\{x^\alpha\})^\mathcal{C}$.

What happens when considering Borel-fixed ideals with Hilbert polynomial $p(t)$ of any degree? The point is to understand how to move in and out monomials for Borel set $\mathscr{B}$ for obtaining another Borel set $\overline{\mathscr{B}}$, with the same Hilbert polynomial, that is by Corollary \ref{cor:bijectionAllIdeals} $\vert \restrict{\mathscr{B}^{\mathcal{C}}}{j}\vert = \vert \restrict{\overline{\mathscr{B}}^{\mathcal{C}}}{j}\vert,\ \forall\ j$. 

The first idea is to exchange a minimal monomial $x^\alpha \in \mathscr{B}$ with a maximal $x^\beta \in \mathscr{B}^{\mathcal{C}}$, but whenever $\min x^\alpha \neq \min x^\beta$ this exchange will not preserve the Hilbert polynomial (see Figure \ref{fig:wrongHP}).

\begin{figure}[H]
\capstart
\begin{center}
\subfloat[][The Borel set $\{I_4\}$ defined by $I=(x_3^2,x_3x_2,x_3x_1,x_2^3) \subset {\K[x_0,x_1,x_2,x_3]}$.]{\label{fig:breakingBorel_a}
\begin{tikzpicture}[scale=0.35]
\tikzstyle{ideal}=[circle,draw=black,fill=black,inner sep=1.5pt]
\tikzstyle{quotient}=[circle,draw=black,thick,inner sep=1.5pt]
\node at (-4,1) [] {};
\node at (10,1) [] {};

\node (0rnn) at (-3,1) [ideal] {};
\node (0rn0) at (-1,1) [ideal] {};
\node (0rn1) at (1,1) [ideal] {};
\node (0rn2) at (3,1) [ideal] {};
\node (0rn3) at (5,1) [ideal] {};

\node (0r00) at (0,0) [ideal] {};
\node (0r01) at (2,0) [ideal] {};
\node (0r02) at (4,0) [ideal] {};
\node (0r03) at (6,0) [ideal] {};

\node (0r11) at (3,-1) [ideal] {};
\node (0r12) at (5,-1) [ideal] {};
\node (0r13) at (7,-1) [regular polygon,regular polygon sides=4,draw=black,inner sep=2pt] {};

\node (0r22) at (6,-2) [regular polygon,regular polygon sides=4,draw=black,fill=black,inner sep=2pt] {};
\node (0r23) at (8,-2) [quotient] {};

\node (0r33) at (9,-3) [quotient] {};

\draw [-,black!25] (0rnn) -- (0rn0);
\draw [-,black!25] (0rn0) -- (0rn1);
\draw [-,black!25] (0rn1) -- (0rn2);
\draw [-,black!25] (0rn2) -- (0rn3);

\draw [-,black!25] (0rn0) -- (0r00);
\draw [-,black!25] (0rn1) -- (0r01);
\draw [-,black!25] (0rn2) -- (0r02);
\draw [-,black!25] (0rn3) -- (0r03);

\draw [-,black!25] (0r00) -- (0r01);
\draw [-,black!25] (0r01) -- (0r02);
\draw [-,black!25] (0r02) -- (0r03);

\draw [-,black!25] (0r01) -- (0r11);
\draw [-,black!25] (0r02) -- (0r12);
\draw [-,black!25] (0r03) -- (0r13);

\draw [-,black!25] (0r11) -- (0r12);
\draw [-,black!25] (0r12) -- (0r13);

\draw [-,black!25] (0r12) -- (0r22);
\draw [-,black!25] (0r13) -- (0r23);

\draw [-,black!25] (0r22) -- (0r23);

\draw [-,black!25] (0r23) -- (0r33);

\node (1r00) at (0,-3.5) [ideal] {};
\node (1r01) at (2,-3.5) [ideal] {};
\node (1r02) at (4,-3.5) [ideal] {};
\node (1r03) at (6,-3.5) [ideal] {};

\node (1r11) at (3,-4.5) [ideal] {};
\node (1r12) at (5,-4.5) [ideal] {};
\node (1r13) at (7,-4.5) [quotient] {};

\node (1r22) at (6,-5.5) [ideal] {};
\node (1r23) at (8,-5.5) [quotient] {};

\node (1r33) at (9,-6.5) [quotient] {};

\draw [-,black!25] (1r00) -- (1r01);
\draw [-,black!25] (1r01) -- (1r02);
\draw [-,black!25] (1r02) -- (1r03);

\draw [-,black!25] (1r01) -- (1r11);
\draw [-,black!25] (1r02) -- (1r12);
\draw [-,black!25] (1r03) -- (1r13);

\draw [-,black!25] (1r11) -- (1r12);
\draw [-,black!25] (1r12) -- (1r13);

\draw [-,black!25] (1r12) -- (1r22);
\draw [-,black!25] (1r13) -- (1r23);

\draw [-,black!25] (1r22) -- (1r23);

\draw [-,black!25] (1r23) -- (1r33);

\node (2r11) at (3,-7.5) [ideal] {};
\node (2r12) at (5,-7.5) [ideal] {};
\node (2r13) at (7,-7.5) [quotient] {};

\node (2r22) at (6,-8.5) [ideal] {};
\node (2r23) at (8,-8.5) [quotient] {};

\node (2r33) at (9,-9.5) [quotient] {};

\draw [-,black!25] (2r11) -- (2r12);
\draw [-,black!25] (2r12) -- (2r13);

\draw [-,black!25] (2r12) -- (2r22);
\draw [-,black!25] (2r13) -- (2r23);

\draw [-,black!25] (2r22) -- (2r23);

\draw [-,black!25] (2r23) -- (2r33);

\node (3r22) at (6,-11.5) [quotient] {};
\node (3r23) at (8,-11.5) [quotient] {};

\node (3r33) at (9,-12.5) [quotient] {};

\draw [-,black!25] (3r22) -- (3r23);

\draw [-,black!25] (3r23) -- (3r33);

\node (4r44) at (9,-15.5) [quotient] {};

\draw [-,black,loosely dotted] (0r00) -- (1r00);
\draw [-,black,loosely dotted] (0r01) -- (1r01);
\draw [-,black,loosely dotted] (0r02) -- (1r02);
\draw [-,black,loosely dotted] (0r03) -- (1r03);

\draw [-,black,loosely dotted] (0r11) -- (1r11);
\draw [-,black,loosely dotted] (0r12) -- (1r12);
\draw [-,black,loosely dotted] (0r13) -- (1r13);

\draw [-,black,loosely dotted] (0r22) -- (1r22);
\draw [-,black,loosely dotted] (0r23) -- (1r23);

\draw [-,black,loosely dotted] (0r33) -- (1r33);

\draw [-,black,loosely dotted] (1r11) -- (2r11);
\draw [-,black,loosely dotted] (1r12) -- (2r12);
\draw [-,black,loosely dotted] (1r13) -- (2r13);

\draw [-,black,loosely dotted] (1r22) -- (2r22);
\draw [-,black,loosely dotted] (1r23) -- (2r23);

\draw [-,black,loosely dotted] (1r33) -- (2r33);

\draw [-,black,loosely dotted] (2r22) -- (3r22);
\draw [-,black,loosely dotted] (2r23) -- (3r23);

\draw [-,black,loosely dotted] (2r33) -- (3r33);

\draw [-,black,loosely dotted] (3r33) -- (4r44);

\end{tikzpicture}
}
\qquad\qquad
\subfloat[][Swapping $x_3x_1^3$ and $x_2^2x_1^2$, we do not obtain a new Borel set, because $x_3x_1^2x_0$ belongs to $\{I_4\}\setminus\{x_3x_1^3\} \cup\{x_2^2x_1^2\}$ and $\up{0}(x_3x_1^2x_0)$ does not.]{\label{fig:breakingBorel_b}
\begin{tikzpicture}[scale=0.35]
\tikzstyle{ideal}=[circle,draw=black,fill=black,inner sep=1.5pt]
\tikzstyle{quotient}=[circle,draw=black,thick,inner sep=1.5pt]
\node at (-4,1) [] {};
\node at (10,1) [] {};

\node (0rnn) at (-3,1) [ideal] {};
\node (0rn0) at (-1,1) [ideal] {};
\node (0rn1) at (1,1) [ideal] {};
\node (0rn2) at (3,1) [ideal] {};
\node (0rn3) at (5,1) [ideal] {};

\node (0r00) at (0,0) [ideal] {};
\node (0r01) at (2,0) [ideal] {};
\node (0r02) at (4,0) [ideal] {};
\node (0r03) at (6,0) [ideal] {};

\node (0r11) at (3,-1) [ideal] {};
\node (0r12) at (5,-1) [ideal] {};
\node (0r13) at (7,-1) [regular polygon,regular polygon sides=4,draw=black,fill=black,inner sep=2pt] {};

\node (0r22) at (6,-2) [regular polygon,regular polygon sides=4,draw=black,inner sep=2pt] {};
\node (0r23) at (8,-2) [quotient] {};

\node (0r33) at (9,-3) [quotient] {};

\draw [-,black!25] (0rnn) -- (0rn0);
\draw [-,black!25] (0rn0) -- (0rn1);
\draw [-,black!25] (0rn1) -- (0rn2);
\draw [-,black!25] (0rn2) -- (0rn3);

\draw [-,black!25] (0rn0) -- (0r00);
\draw [-,black!25] (0rn1) -- (0r01);
\draw [-,black!25] (0rn2) -- (0r02);
\draw [-,black!25] (0rn3) -- (0r03);

\draw [-,black!25] (0r00) -- (0r01);
\draw [-,black!25] (0r01) -- (0r02);
\draw [-,black!25] (0r02) -- (0r03);

\draw [-,black!25] (0r01) -- (0r11);
\draw [-,black!25] (0r02) -- (0r12);
\draw [-,black!25] (0r03) -- (0r13);

\draw [-,black!25] (0r11) -- (0r12);
\draw [-,black!25] (0r12) -- (0r13);

\draw [-,black!25] (0r12) -- (0r22);
\draw [-,black!25] (0r13) -- (0r23);

\draw [-,black!25] (0r22) -- (0r23);

\draw [-,black!25] (0r23) -- (0r33);

\node (1r00) at (0,-3.5) [ideal] {};
\node (1r01) at (2,-3.5) [ideal] {};
\node (1r02) at (4,-3.5) [ideal] {};
\node (1r03) at (6,-3.5) [ideal] {};

\node (1r11) at (3,-4.5) [ideal] {};
\node (1r12) at (5,-4.5) [ideal] {};
\node (1r13) at (7,-4.5) [quotient] {};

\node (1r22) at (6,-5.5) [ideal] {};
\node (1r23) at (8,-5.5) [quotient] {};

\node (1r33) at (9,-6.5) [quotient] {};

\draw [-,black!25] (1r00) -- (1r01);
\draw [-,black!25] (1r01) -- (1r02);
\draw [-,black!25] (1r02) -- (1r03);

\draw [-,black!25] (1r01) -- (1r11);
\draw [-,black!25] (1r02) -- (1r12);
\draw [-,black!25] (1r03) -- (1r13);

\draw [-,black!25] (1r11) -- (1r12);
\draw [-,black!25] (1r12) -- (1r13);

\draw [-,black!25] (1r12) -- (1r22);
\draw [-,black!25] (1r13) -- (1r23);

\draw [-,black!25] (1r22) -- (1r23);

\draw [-,black!25] (1r23) -- (1r33);

\node (2r11) at (3,-7.5) [ideal] {};
\node (2r12) at (5,-7.5) [ideal] {};
\node (2r13) at (7,-7.5) [quotient] {};

\node (2r22) at (6,-8.5) [ideal] {};
\node (2r23) at (8,-8.5) [quotient] {};

\node (2r33) at (9,-9.5) [quotient] {};

\draw [-,black!25] (2r11) -- (2r12);
\draw [-,black!25] (2r12) -- (2r13);

\draw [-,black!25] (2r12) -- (2r22);
\draw [-,black!25] (2r13) -- (2r23);

\draw [-,black!25] (2r22) -- (2r23);

\draw [-,black!25] (2r23) -- (2r33);

\node (3r22) at (6,-11.5) [quotient] {};
\node (3r23) at (8,-11.5) [quotient] {};

\node (3r33) at (9,-12.5) [quotient] {};

\draw [-,black!25] (3r22) -- (3r23);

\draw [-,black!25] (3r23) -- (3r33);

\node (4r44) at (9,-15.5) [quotient] {};

\draw [-,black,loosely dotted] (0r00) -- (1r00);
\draw [-,black,loosely dotted] (0r01) -- (1r01);
\draw [-,black,loosely dotted] (0r02) -- (1r02);
\draw [-,black,loosely dotted] (0r03) -- (1r03);

\draw [-,black,loosely dotted] (0r11) -- (1r11);
\draw [-,black,loosely dotted] (0r12) -- (1r12);
\draw [-,black,loosely dotted] (0r13) -- (1r13);

\draw [-,black,loosely dotted] (0r22) -- (1r22);
\draw [-,black,loosely dotted] (0r23) -- (1r23);

\draw [-,black,loosely dotted] (0r33) -- (1r33);

\draw [-,black,loosely dotted] (1r11) -- (2r11);
\draw [-,black,loosely dotted] (1r12) -- (2r12);
\draw [-,black,loosely dotted] (1r13) -- (2r13);

\draw [-,black,loosely dotted] (1r22) -- (2r22);
\draw [-,black,loosely dotted] (1r23) -- (2r23);

\draw [-,black,loosely dotted] (1r33) -- (2r33);

\draw [-,black,loosely dotted] (2r22) -- (3r22);
\draw [-,black,loosely dotted] (2r23) -- (3r23);

\draw [-,black,loosely dotted] (2r33) -- (3r33);

\draw [-,black,loosely dotted] (3r33) -- (4r44);

\end{tikzpicture}
}
\end{center}
\caption[Example of an exchange of monomials that does not preserves the Borel property.]{\label{fig:breakingBorel} An example of an exchange of monomials that does not preserves the Borel property.}
\end{figure}
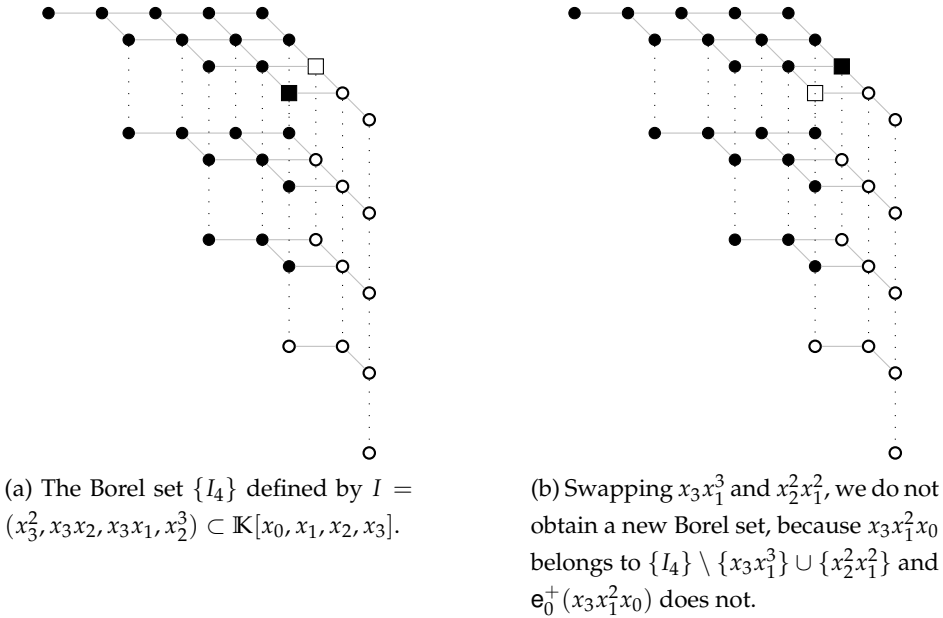

A second idea could be to exchange a minimal monomial $x^\alpha$ in $\restrict{\mathscr{B}}{k}$ and a maximal monomial $x^\beta$ in $\restrict{\mathscr{B}^{\mathcal{C}}}{k}$. In this way, the cardinality of the sets $\restrict{\mathscr{B}}{k}$ and $\restrict{\overline{\mathscr{B}}}{k}$ is preserved, but whenever $\down{k}(x^\alpha)$ belongs to $\mathscr{B}$, swapping $x^\alpha$ and $x^\beta$ we do not obtain a Borel set (see Figure \ref{fig:breakingBorel}). This fact suggests that the general case is more complicated and that we have to swap more monomials than a minimal and a maximal element.

\begin{figure}[H] 
\capstart
\begin{center}
\subfloat[][The Borel set $\{I_4\}$ defined by $I = (x_2^2,x_2x_1,x_1^3)$.]{\label{fig:ExampleSwap_a}
\begin{tikzpicture}[scale=0.9,>=latex]
\tikzstyle{quotient}=[]
\node at (-0.5,5) [] {};
\node at (7.5,5) [] {};
\node (M_1) at (0.5,5) [quotient] {\footnotesize $x_2^4$};
\node (M_2) at (2,5) [quotient] {\footnotesize $x_2^3x_1$};
\node (M_3) at (3.5,5) [quotient] {\footnotesize $x_2^2x_1^2$};
\node (M_4) at (5,5) [quotient] {\footnotesize $x_2x_1^3$};
\node (M_5) at (6.5,5) [quotient] {\footnotesize $x_1^4$};

\node (M_6) at (2,3.75) [quotient] {\footnotesize $x_2^3x_0$};
\node (M_7) at (3.5,3.75) [quotient] {\footnotesize $x_2^2x_1x_0$};
\node (M_8) at (5,3.75) [quotient] {\footnotesize $x_2x_1^2x_0$};
\node (M_9) at (6.5,3.75) [quotient] {\footnotesize $x_1^3x_0$};

\node (M_10) at (3.5,2.5) [quotient] {\footnotesize $x_2^2x_0^2$};
\node (M_11) at (5,2.5) [quotient] {\footnotesize $x_2x_1x_0^2$};
\node (M_12) at (6.5,2.5) [quotient] {\footnotesize $x_1^2x_0^2$};

\node (M_13) at (5,1.25) [quotient] {\footnotesize $x_2x_0^3$};
\node (M_14) at (6.5,1.25) [quotient] {\footnotesize $x_1x_0^3$};

\node (M_15) at (6.5,0) [quotient] {\footnotesize $x_0^4$};

\draw [->] (M_1) -- (M_2);
\draw [->] (M_2) -- (M_3);
\draw [->] (M_3) -- (M_4);
\draw [->] (M_4) -- (M_5);

\draw [->] (M_2) -- (M_6);
\draw [->] (M_3) -- (M_7);
\draw [->] (M_4) -- (M_8);
\draw [->] (M_5) -- (M_9);

\draw [->] (M_6) -- (M_7);
\draw [->] (M_7) -- (M_8);
\draw [->] (M_8) -- (M_9);

\draw [->] (M_7) -- (M_10);
\draw [->] (M_8) -- (M_11);
\draw [->] (M_9) -- (M_12);

\draw [->] (M_10) -- (M_11);
\draw [->] (M_11) -- (M_12);

\draw [->] (M_11) -- (M_13);
\draw [->] (M_12) -- (M_14);

\draw [->] (M_13) -- (M_14);
\draw [->] (M_14) -- (M_15);

\draw [rounded corners=6pt,thick] (4.5,1.85) -- (5.75,1.85) -- (5.75,3.1) -- (7.25,3.1) -- (7.25,5.75) --  (-0.25,5.75) -- (-0.25,4.35) -- (1.25,4.35) -- (1.25,3.1) -- (2.75,3.1) -- (2.75,1.85) -- cycle;
\end{tikzpicture}
}
\\
\subfloat[][$\{I_4\} \setminus \{x_2x_1x_0^2\} \cup \{x_1^2x_0^2\}$.]{\label{fig:ExampleSwap_b}
\begin{tikzpicture}[>=latex,scale=0.8]
\tikzstyle{ideal}=[circle,draw=black,fill=black,inner sep=1.5pt]
\tikzstyle{quotient}=[circle,draw=black,thick,inner sep=1.5pt]

\node at (0,0) [ideal] {};
\node at (1,0) [ideal] {};
\node at (2,0) [ideal] {};
\node at (3,0) [ideal] {};
\node at (4,0) [ideal] {};

\node at (1,-1) [ideal] {};
\node at (2,-1) [ideal] {};
\node at (3,-1) [ideal] {};
\node at (4,-1) [ideal] {};

\node at (2,-2) [ideal] {};
\node (W_5) at (3,-2) [regular polygon,regular polygon sides=4,draw=black,inner sep=1.5pt] {};
\node (W_6) at (4,-2) [regular polygon,regular polygon sides=4,draw=black,fill=black,inner sep=1.5pt] {};

\draw [->] (W_6) -- (W_5);

\node at (3,-3) [quotient] {};
\node at (4,-3) [quotient] {};

\node at (4,-4) [quotient] {};

\draw[rounded corners=5pt,dashed] (-0.5,0.5) -- (4.5,0.5)  -- (4.5,-2.5) -- (3.5,-2.5) -- (3.5,-1.5) -- (2.5,-1.5) -- (2.5,-2.5) -- (1.5,-2.5) -- (1.5,-1.5) -- (0.5,-1.5) -- (0.5,-0.5)  -- (-0.5,-0.5) -- cycle;
\end{tikzpicture}
}
\qquad\qquad\qquad
\subfloat[][$\{I_4\} \setminus \{x_1^3x_0\} \cup \{x_1^2 x_0^2\}$.]{\label{fig:ExampleSwap_c}
\begin{tikzpicture}[>=latex,scale=0.8]
\tikzstyle{ideal}=[circle,draw=black,fill=black,inner sep=1.5pt]
\tikzstyle{quotient}=[circle,draw=black,thick,inner sep=1.5pt]

\node at (0,0) [ideal] {};
\node at (1,0) [ideal] {};
\node at (2,0) [ideal] {};
\node at (3,0) [ideal] {};
\node at (4,0) [ideal] {};

\node at (1,-1) [ideal] {};
\node at (2,-1) [ideal] {};
\node at (3,-1) [ideal] {};
\node (W_5) at (4,-1) [regular polygon,regular polygon sides=4,draw=black,inner sep=1.5pt] {};

\node at (2,-2) [ideal] {};
\node at (3,-2) [ideal] {};
\node (W_6) at (4,-2) [regular polygon,regular polygon sides=4,draw=black,fill=black,inner sep=1.5pt] {};

\draw [->] (W_6) -- (W_5);

\node at (3,-3) [quotient] {};
\node at (4,-3) [quotient] {};

\node at (4,-4) [quotient] {};

\draw[rounded corners=5pt,dashed] (-0.5,0.5) -- (4.5,0.5)  -- (4.5,-0.5) -- (3.5,-0.5) -- (3.5,-1.5) -- (4.5,-1.5) -- (4.5,-2.5) -- (1.5,-2.5) -- (1.5,-1.5) -- (0.5,-1.5) -- (0.5,-0.5)  -- (-0.5,-0.5) -- cycle;
\end{tikzpicture}
}
\\
\subfloat[][$\{I_4\}\setminus\{x_2x_1x_0^2\}\cup\{x_2x_0^3\}$.]{\label{fig:ExampleSwap_d}
\begin{tikzpicture}[>=latex,scale=0.8]
\tikzstyle{ideal}=[circle,draw=black,fill=black,inner sep=1.5pt]
\tikzstyle{quotient}=[circle,draw=black,thick,inner sep=1.5pt]

\node at (0,0) [ideal] {};
\node at (1,0) [ideal] {};
\node at (2,0) [ideal] {};
\node at (3,0) [ideal] {};
\node at (4,0) [ideal] {};

\node at (1,-1) [ideal] {};
\node at (2,-1) [ideal] {};
\node at (3,-1) [ideal] {};
\node  at (4,-1) [ideal] {};

\node at (2,-2) [ideal] {};
\node (W_5) at (3,-2) [regular polygon,regular polygon sides=4,draw=black,inner sep=1.5pt] {};
\node  at (4,-2) [quotient] {};

\node at (3,-3)(W_6) [regular polygon,regular polygon sides=4,draw=black,fill=black,inner sep=1.5pt] {};
\node at (4,-3) [quotient] {};

\node at (4,-4) [quotient] {};

\draw [->] (W_6) -- (W_5);

\draw[rounded corners=5pt,dashed] (-0.5,0.5) -- (4.5,0.5)  -- (4.5,-1.5) -- (2.5,-1.5) -- (2.5,-2.5) -- (1.5,-2.5) -- (1.5,-1.5) -- (0.5,-1.5) -- (0.5,-0.5)  -- (-0.5,-0.5) -- cycle;
\draw[rounded corners=5pt,dashed] (2.5,-2.5) -- (3.5,-2.5) -- (3.5,-3.5) -- (2.5,-3.5) -- cycle;
\end{tikzpicture}
}
\qquad\qquad\qquad
\subfloat[][$\{I_4\}\setminus\{x_2x_1x_0^2\}\cup\{x_2x_0^3\}$.]{\label{fig:ExampleSwap_e}
\begin{tikzpicture}[>=latex,scale=0.8]
\tikzstyle{ideal}=[circle,draw=black,fill=black,inner sep=1.5pt]
\tikzstyle{quotient}=[circle,draw=black,thick,inner sep=1.5pt]

\node at (0,0) [ideal] {};
\node at (1,0) [ideal] {};
\node at (2,0) [ideal] {};
\node at (3,0) [ideal] {};
\node at (4,0) [ideal] {};

\node at (1,-1) [ideal] {};
\node at (2,-1) [ideal] {};
\node at (3,-1) [ideal] {};
\node (W_5) at (4,-1) [regular polygon,regular polygon sides=4,draw=black,inner sep=1.5pt] {};

\node at (2,-2) [ideal] {};
\node  at (3,-2) [ideal] {};
\node  at (4,-2) [quotient] {};

\node at (3,-3)(W_6) [regular polygon,regular polygon sides=4,draw=black,fill=black,inner sep=1.5pt] {};
\node at (4,-3) [quotient] {};

\node at (4,-4) [quotient] {};

\draw[rounded corners=5pt] (-0.5,0.5) -- (4.5,0.5)  -- (4.5,-0.5) -- (3.5,-0.5) -- (3.5,-3.5) -- (2.5,-3.5) -- (2.5,-2.5) -- (1.5,-2.5) -- (1.5,-1.5) -- (0.5,-1.5) -- (0.5,-0.5)  -- (-0.5,-0.5) -- cycle;
\end{tikzpicture}
}
\end{center}
\caption[Attempts of exchange of pairs of monomials preserving the Borel property and the Hilbert polyomial.]{\label{fig:ExampleSwap}The graphical description of Example \ref{ex:ExampleSwap}. In Figure \ref{fig:ExampleSwap_a} there is the Borel set defined by the ideal initially considered. In the other figures, there are the sets that can be obtained swapping a pair of monomials composed by a minimal element of $\{I_4\}$ and a maximal element of the complement. With $\square$ and $\blacksquare$ we denote the monomials exchanged. The arrows correspond to the increasing elementary moves showing that the set of monomials obtained is not Borel.}
\end{figure}
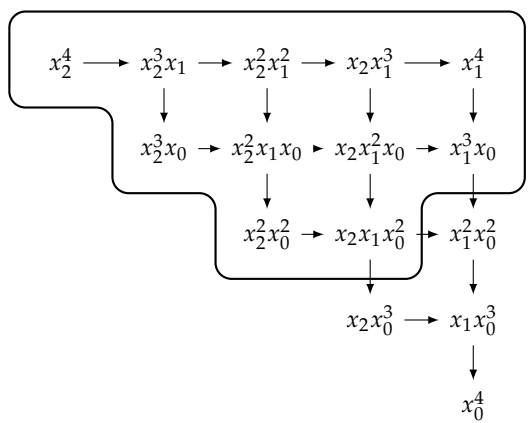
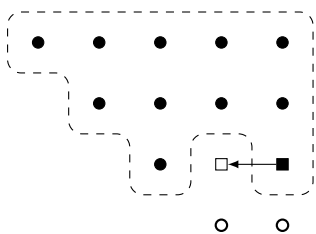
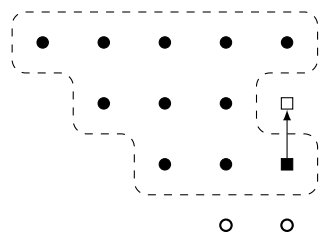
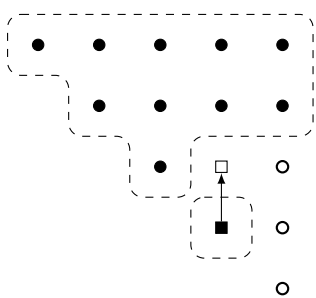
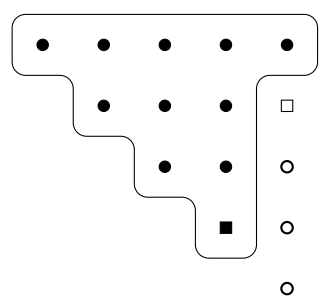

\begin{figure}[H]
\capstart
\begin{center}
\subfloat[][The Borel set $\{I_4\}$ defined by the ideal $I = (x_2^3,x_2^2x_1,x_2x_1^2)\subset {\K[x_0,x_1,x_2]}$ with Hilbert polynomial $p(t) = t+4$.]{\label{fig:wrongHP_a}
\begin{tikzpicture}[scale=0.75]
\tikzstyle{ideal}=[circle,draw=black,fill=black,inner sep=1.5pt]
\tikzstyle{quotient}=[circle,draw=black,thick,inner sep=1.5pt]
\node at (-1.5,0) [] {};
\node at (5.5,0) [] {};

\node at (0,0) [ideal] {};
\node at (1,0) [ideal] {};
\node at (2,0) [ideal] {};
\node at (3,0) [ideal] {};
\node at (4,0) [regular polygon,regular polygon sides=4,draw=black,inner sep=2pt] {};

\node at (1,-1) [ideal] {};
\node at (2,-1) [ideal] {};
\node at (3,-1) [regular polygon,regular polygon sides=4,draw=black,fill=black,inner sep=2pt] {};
\node at (4,-1) [quotient] {};

\node at (2,-2) [quotient] {};
\node at (3,-2) [quotient] {};
\node at (4,-2) [quotient] {};

\node at (3,-3) [quotient] {};
\node at (4,-3) [quotient] {};

\node at (4,-4) [quotient] {};

\draw [rounded corners=5pt] (-0.5,0.5) -- (3.5,0.5) -- (3.5,-1.5) -- (0.5,-1.5) -- (0.5,-0.5) -- (-0.5,-0.5) -- cycle;
\end{tikzpicture}
}
\qquad\qquad
\subfloat[][The Borel set $\{I_4\} \setminus\{x_2x_1^2 x_0\} \cup \{x_1^4\}$ corresponding to the ideal $(x_2^3,x_2^2x_1,x_2x_1^3,x_1^4) \subset {\K[x_0,x_1,x_2]}$ with Hilbert polynomial $\overline{p}(t) = 8$.]{\label{fig:wrongHP_b}
\begin{tikzpicture}[scale=0.75]
\tikzstyle{ideal}=[circle,draw=black,fill=black,inner sep=1.5pt]
\tikzstyle{quotient}=[circle,draw=black,thick,inner sep=1.5pt]

\node at (-1.5,0) [] {};
\node at (5.5,0) [] {};

\node at (0,0) [ideal] {};
\node at (1,0) [ideal] {};
\node at (2,0) [ideal] {};
\node at (3,0) [ideal] {};
\node at (4,0) [regular polygon,regular polygon sides=4,draw=black,fill=black,inner sep=2pt] {};

\node at (1,-1) [ideal] {};
\node at (2,-1) [ideal] {};
\node at (3,-1) [regular polygon,regular polygon sides=4,draw=black,inner sep=2pt] {};
\node at (4,-1) [quotient] {};

\node at (2,-2) [quotient] {};
\node at (3,-2) [quotient] {};
\node at (4,-2) [quotient] {};

\node at (3,-3) [quotient] {};
\node at (4,-3) [quotient] {};

\node at (4,-4) [quotient] {};

\draw [rounded corners=5pt] (-0.5,0.5) -- (4.5,0.5) -- (4.5,-0.5) -- (2.5,-0.5) -- (2.5,-1.5) -- (0.5,-1.5) -- (0.5,-0.5) -- (-0.5,-0.5) -- cycle;
\end{tikzpicture}
}
\end{center}
\caption[Example of an exchange of monomials that does not preserves the Hilbert polynomial.]{\label{fig:wrongHP} An example of an exchange of monomials that does not preserves the Hilbert polynomial.}
\end{figure}
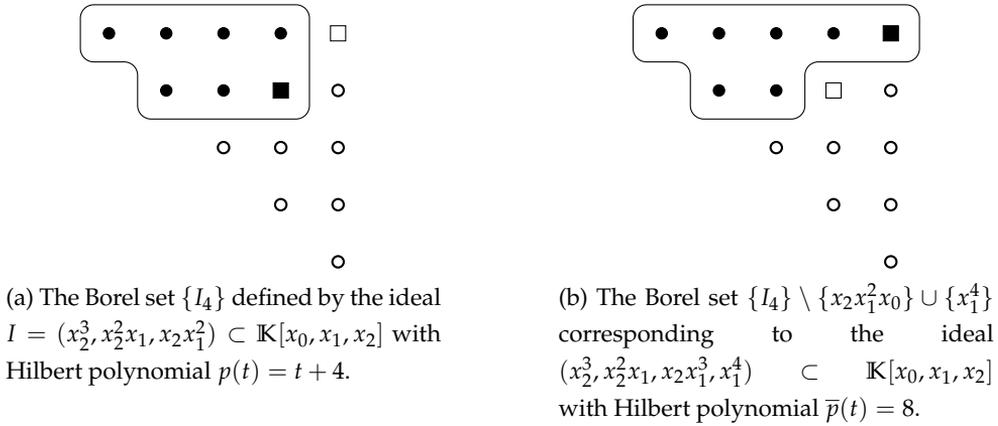

\begin{definition}\label{def:decreasingMoves}
Let $\mathscr{B} \subset \poset$ be a Borel set and let $x^\alpha$ be a minimal monomial of $\restrict{\mathscr{B}}{k}$ with $\min x^\alpha = k$. The set of monomials in $\mathscr{B}$ smaller than $x^\alpha$
\begin{equation}\label{eq:smallerMonomials}
\left\{ x^{\gamma} \in \mathscr{B}\ \vert\ x^\gamma \leq_B x^\alpha \right\}
\end{equation}
is uniquely associated to the composition of decreasing elementary moves
\begin{equation}\label{eq:smallerMoves}
\mathcal{F}_{\alpha} = \left\{ \mu_b \down{j_b}\circ \cdots\circ \mu_1 \down{j_1}\ \ \big\vert\ \ \mu_b \down{j_b}\circ \cdots\circ \mu_1 \down{j_1}(x^\alpha) \in \mathscr{B}\right\}.
\end{equation}
Note that being $x^\alpha$ minimal in the restriction $\restrict{\mathscr{B}}{k}$, the elementary moves appearing in the composed moves of $\mathcal{F}_{\alpha}$ have index smaller than or equal to $k$.
We call $\mathcal{F}_\alpha$ \emph{decreasing set} of $x^\alpha$ and we denote with $\mathcal{F}_{\alpha}(x^\alpha)$ the set of monomials in \eqref{eq:smallerMonomials}.\index{descreasing set of a monomial}

Given a maximal element $x^\beta$ of $\restrict{\mathscr{B}^{\mathcal{C}}}{k}$, we will say that the set $\mathcal{F}_\alpha$ is \emph{Borel-admissible}\index{descreasing set of a monomial!Borel-admissible} w.r.t. $x^\beta$ if every $F \in \mathcal{F}_\alpha$ is admissible also for $x^\beta$ (note that $\min x^\beta = k$) and if $\up{i}\big(F(x^\beta)\big) \in \mathscr{B}$ for each admissible elementary increasing move $\up{i},\ i \geqslant k$. We denote by $\mathcal{F}_{\alpha}(x^\beta)$ the set of monomials obtained from $x^\beta$ applying the moves in $\mathcal{F}_\alpha$.
\end{definition}

\begin{remark}\label{rk:borelAndDisjoint}
Note that, recalling Definition \ref{def:minimalMaximal}, the set $\mathcal{F}_{\alpha}(x^\beta)$ contains $k$-maximal elements (by construction) and $\mathcal{F}_{\alpha}(x^\alpha)$ $k$-minimal monomials, indeed for each $x^\gamma = F(x^\alpha),\ F \in \mathcal{F}_{\alpha}$ and any admissible $\down{j},\ j > k$, $\down{j}\big(F(x^\alpha)\big) = F \big(\down{j}(x^\alpha)\big) \in \mathscr{B}^{\mathcal{C}}$, because $\down{j}(x^\alpha) \in \mathscr{B}^{\mathcal{C}}$. Furthermore these two sets of monomials are always disjoint. Indeed we can write $x^\alpha = x^{\overline{\alpha}} x_k^a$ and $x^{\overline{\beta}}x_k^a$, so that $\min x^{\overline{\alpha}} \geqslant k$ and $\min x^{\overline{\beta}} \geqslant k$. The generic monomials in $\mathcal{F}_\alpha(x^\alpha)$ and $\mathcal{F}_\alpha(x^\beta)$ are of the type $x^{\overline{\alpha}}x^\gamma$ and $x^{\overline{\beta}}x^\delta$ where $\max x^\gamma \leqslant k$ and $\max x^\delta \leqslant k$, so since $x^{\overline{\alpha}}\neq x^{\overline{\beta}}$ we can not obtain the same monomial in the two sets.
\end{remark}

\begin{lemma}\label{lem:PreservingBorel}
 Let $\mathscr{B} \subset \pos{n}{m}$ be a Borel set and let $x^\alpha$ and $x^\beta$ be a minimal monomial of $\restrict{\mathscr{B}}{k}$ and a maximal monomial of $\restrict{\mathscr{B}^{\mathcal{C}}}{k}$ with $\min x^\alpha = \min x^\beta = k$. If there does not exist an elementary move $\down{j},\ j > k$ such that $\down{j}(x^\alpha) = x^\beta$ and if the decreasing set $\mathcal{F}_\alpha$ of $x^\alpha$ is Borel-admissible w.r.t. $x^\beta$, then $\overline{\mathscr{B}} = \mathscr{B} \setminus \mathcal{F}_\alpha(x^\alpha) \cup \mathcal{F}_\alpha(x^\beta)$ is a Borel set and $\left\vert \restrict{\mathscr{B}^{\mathcal{C}}}{i} \right\vert = \left\vert \restrict{\overline{\mathscr{B}}^{\mathcal{C}}}{i} \right\vert,\ \forall\ i$.
\begin{proof}
The set $\overline{\mathscr{B}} = \mathscr{B}\setminus \mathcal{F}_\alpha(x^\alpha) \cup \mathcal{F}_\alpha(x^\beta)$ turns out to be Borel, mainly for Remark \ref{rk:borelAndDisjoint} and because the monomials in the set $\mathcal{F}_\alpha(x^\alpha)$ and $\mathcal{F}_\alpha(x^\beta)$ are not comparable w.r.t. the Borel order $\leq_B$.

The condition $\big\vert \restrict{\mathscr{B}^{\mathcal{C}}}{i} \big\vert = \big\vert \restrict{\overline{\mathscr{B}}^{\mathcal{C}}}{i} \big\vert$ follows from the fact that for each $F \in \mathcal{F}_{\alpha}$, $\min F(x^\alpha) = \min F(x^\beta)$. \end{proof}
\end{lemma}

\begin{example}\label{ex:breakPreserveBorel}
Let us see with some example in the poset $\pos{3}{4}$ the problems we can meet in choosing sets of monomials to exchange preserving both the Borel property and the Hilbert polynomial.
\begin{description}
 \item[(Figure \ref{fig:borelAdmissible_a})] Let $\{I_3\}$ be the Borel set defined by the ideal $I = (x_3,x_2^3)$. $x_3x_1^2$ is a minimal monomial in $\restrict{\{I_3\}}{1}$ and $x_2^2x_1$ a maximal one in $\restrict{\{I_3\}^{\mathcal{C}}}{1}$. The decreasing set of $x_3x_1^2$ is $\mathcal{F}_{x_3x_1^2} = \{\mathrm{id},\down{1},2\down{1}\}$
and it turns out to be not Borel-admissible w.r.t. $x_2^2x_1$, because $2\down{1}$ is not admissible for $x_2^2 x_1$.
\item[(Figure \ref{fig:borelAdmissible_b})] Let now consider the Borel set $\{J_4\}$ defined by $(x_3^2,x_3x_2^2,x_2^3,x_3x_2x_1^2)$. $x_2^3x_1$ is a minimal monomial in $\restrict{\{J_4\}}{1}$ and $x_3x_1^3$ a maximal one in $\restrict{\{J_4\}^{\mathcal{C}}}{1}$. The decreasing set of $x_2^3 x_1$ is $\mathcal{F}_{x_2^3x_1} = \{\mathrm{id},\down{1}\}$ and each move in it is also admissible w.r.t. $x_3x_1^3$. But $\up{1} \big(\down{1}(x_3x_1^3)\big) = x_3x_2x_1x_0 \notin \{J_4\}$, so exchanging $\mathcal{F}_{x_2^3x_1}(x_2^3x_1)$ with $\mathcal{F}_{x_2^3x_1}(x_3x_1^3)$ would not respect the Borel property.
\item[] The same decreasing set arises in the Borel set defined by the ideal $(x_3^2,x_3x_2^2,x_2^3,$ $x_3x_2x_1)$, but in this case $\up{1} \big(\down{1}(x_3x_1^3)\big) = x_3x_2x_1x_0$ would belong to the ideal, so the decreasing set would be Borel-admissible w.r.t. $x_3x_1^3$.
\end{description}
\end{example}

After having understood how to move monomials out of and in to a Borel set preserving the Borel property and the Hilbert polynomial associated, we carry on trying to determine a flat deformation having among the fibers the Borel-fixed ideals associated to the Borel sets. We will use the following property.

\begin{proposition}[{\cite[Chapter 1 Section 3]{Artin}}]\label{prop:artinSyzygies}\index{flatness}
Let $A$ be a local $\K$-algebra and let us consider $M = \K[x]/\langle f_1,\ldots,f_s\rangle$ and $M_A = A[x]/\langle f_1',\ldots,f_s' \rangle$, where $f_i'$ is a lifting of $f_i$ in $A[x]$, i.e. tensoring the ideal $\langle f_1',\ldots,f_s'\rangle$ by $\K$, the residue field of $A$, leads to the ideal $\langle f_1,\ldots,f_s\rangle$. Then, $M_A$ is flat over $A$ if and only if any relation among $(f_1,\ldots,f_s)$ lifts to a relation among $(f_1',\ldots,f_s')$. 
\end{proposition}

\input{Figures/Ch3/exampleBorelAdmissible}

\newpage

\begin{theorem}\label{th:flatGeneralDim}
Let $\{I_r\}\subset\pos{n}{r}$ be a Borel set defined by a Borel-fixed ideal $I$. Moreover let us suppose that there exist two monomials: $x^\alpha$ minimal element of $\restrict{\{I_r\}}{k}$ and $x^\beta$ maximal element of $\restrict{\{I_r\}^{\mathcal{C}}}{k}$, such that $\min x^\alpha = \min x^\beta = k$, there does not exist a decreasing move $\down{j},\ j > k$ with $\down{j}(x^\alpha) =  x^\beta$ and such that the decreasing set $\mathcal{F}_{\alpha}$ of $x^\alpha$ is Borel-admissible w.r.t. $x^\beta$. Then, the family of subschemes of $\PP^n$ parametrized by the ideal 
\begin{equation}\label{eq:familyIdeals}
 \mathcal{I} = \left\langle \left\{I_r\right\}\setminus\mathcal{F}_\alpha(x^\alpha)\cup\{y_0 \, F(x^\alpha) + y_1 \, F(x^\beta) \ \vert\ F \in \mathcal{F}_\alpha\}\right\rangle
\end{equation}
is flat over $\PP^1$.
\end{theorem}
\begin{proof}
First of all let $\mathscr{B} = \{I_r\}$ and let us call $p(t)$ the Hilbert polynomial of $I$. By Lemma \ref{lem:PreservingBorel}, also the ideal $\mathcal{I}\vert_{[0:1]} = \big\langle \{I_r\}\setminus\mathcal{F}_\alpha(x^\alpha)\cup\mathcal{F}_\alpha(x^\beta)\big\rangle$ has Hilbert polynomial $p(t)$.

Let us suppose $y_0 \neq 0$ and set $z = \frac{y_1}{y_0}$ let us denote by $I_{z}$ the ideal $\mathcal{I}\vert_{[1:z]} \subset \K[z][x]$. By \cite[Chapter III Theorem 9.9]{Hartshorne} we know that to prove the flatness it suffices to check that for each point $\mathfrak{p} \in \Spec \K[z]$ the Hilbert polynomial of $(I_z)_{\mathfrak{p}}$ is $p(t)$. Let us start considering the presentation map of $I = I_{\geqslant r}$
\begin{equation}
\begin{split}
\psi_0 : \big(\K[x](-r)\big)^{\vert\mathscr{B}\vert} & {} \longrightarrow \K[x]\\
 \parbox{2cm}{\centering $\mathbf{e}_{\gamma}$} &{} \longmapsto \parbox{1cm}{\centering $x^\gamma$,}\quad \forall\ x^\gamma \in \mathscr{B}
\end{split}
\end{equation}
and the set of Eliahou-Kervaire syzygies generating $\ker \psi_0$
\begin{equation*}
\syz_{\textnormal{EK}}(I) = \left\{x_i \mathbf{e}_{\gamma} - \min x^\gamma \mathbf{e}_\delta\ \vert\ \forall\ x^\gamma \in \mathscr{B},\ \forall\ x_i >_B \min x^\gamma\text{ s.t. } x_i x^\gamma = \dec{x^\delta}{\min x^\gamma}{I}{}\right\}.
\end{equation*}
Now we will show that the kernel of the presentation map of $I_z$
\begin{equation}
\begin{split}
\psi : \big(\K[z][x](-r)\big)^{\vert\mathscr{B}\vert} & {} \longrightarrow\ \K[z][x]\\
 \parbox{2cm}{\centering $\mathbf{f}_{\gamma}$} &{} \longmapsto \begin{cases}
 x^\gamma &  \forall\ x^\gamma \in \mathscr{B}\setminus\mathcal{F}_{\alpha}(x^\alpha)\\
 F(x^\alpha) + zF(x^\beta), & \forall\ x^\gamma = F(x^\alpha), F \in \mathcal{F}_{\alpha}
 \end{cases}
\end{split}
\end{equation}
has a set of generators of syzygies lifted directly from $\syz_{\textnormal{EK}}(I)$. For any $x^\gamma \in \mathscr{B} \setminus \mathcal{F}_{\alpha}(x^\alpha)$, since $x^\gamma = \psi(\mathbf{f}_\gamma) = \psi_0(\mathbf{e}_\gamma)$, we consider the same syzygies of $I$
\[
x_i \mathbf{f}_{\gamma} - \min x^\gamma \mathbf{f}_\delta \qquad \text{in place of} \qquad x_i \mathbf{e}_{\gamma} - \min x^\gamma \mathbf{e}_\delta,
\]
indeed $x^{\delta} = \frac{x_i}{\min x^\gamma} x^\gamma = \up{i-1} \circ \cdots \circ \up{\min x^\gamma}(x^\gamma)$ belongs to $\mathscr{B}\setminus\mathcal{F}_{\alpha}(x^\alpha)$ as well. 

Let us now look at the generators $x^{\overline{\alpha}}+zx^{\overline{\beta}} = F(x^\alpha) + z F(x^\beta)$. For each $x_i >_B x_k = \min x^\alpha = \min x^\beta >_B \min x^{\overline{\alpha}} = \min x^{\overline{\beta}}$, both $x_i x^{\overline{\alpha}}$ and $x_i x^{\overline{\beta}}$ belong to $\mathscr{B}\setminus\mathcal{F}_{\alpha}(x^\alpha)$ because the identity $\frac{x_i}{\min x^{\overline{\alpha}}} = \up{i-1}\circ\cdots\circ\up{\min x^{\overline{\alpha}}}$ and from the fact that $x^\beta$ is $k$-maximal element and that in order for moving $x^{\overline{\alpha}}$ to any other monomial of $\mathcal{F}_{\alpha}$ we should only use moves $\up{j}$ with index $0\leqslant j < k$. Thus we consider the following syzygies
\[
x_i \mathbf{f}_{\overline{\alpha}} - \min x^{\overline{\alpha}} \mathbf{f}_{\gamma} - z\, \min x^{\overline{\beta}} \mathbf{f}_{\delta} \qquad(\text{in place of } x_i \mathbf{e}_{\overline{\alpha}} - \min x^{\overline{\alpha}} \mathbf{e}_{\gamma})
\]
where $ x_i x^{\overline{\alpha}} = \dec{x^\gamma}{\min x^{\overline{\alpha}}}{I_r}{},\ x_i x^{\overline{\beta}} = \dec{x^\delta}{\min x^{\overline{\beta}}}{I_r}{}$ and $x^\gamma,x^\delta \in \mathscr{B}\setminus\mathcal{F}_\alpha$. The last possibility to consider is when multiplying $x^{\overline{\alpha}} + zx^{\overline{\beta}}$ by a variable $\min x^{\overline{\alpha}} <_B x_i \leq_B x_k$. In this case $\frac{x_i}{\min x^{\overline{\alpha}}} = \up{i-1}\circ\cdots\circ\up{\min x^{\overline{\alpha}}}$ involves only elementary decreasing moves with index included in 0 and $k-1$, that is
\[
\up{i-1}\circ\cdots\circ\up{\min x^{\overline{\alpha}}}(x^{\overline{\alpha}}) = \up{i-1}\circ\cdots\circ\up{\min x^{\overline{\alpha}}}\big(F(x^\alpha)\big) = x^{\widetilde{\alpha}} \in \mathcal{F}_{\alpha}(x^\alpha)
\]
and
\[
\up{i-1}\circ\cdots\circ\up{\min x^{\overline{\beta}}}(x^{\overline{\beta}}) = \up{i-1}\circ\cdots\circ\up{\min x^{\overline{\beta}}}\big(F(x^\beta)\big) = x^{\widetilde{\beta}} \in \mathcal{F}_{\alpha}(x^\beta).
\]
In this case the sygygies are equal to the corresponding syzygies of $I$:
\[
x_i \mathbf{f}_{\overline{\alpha}} - \min x^{\overline{\alpha}} \mathbf{f}_{\widetilde{\alpha}} \qquad \text{in place of}\qquad x_i \mathbf{e}_{\overline{\alpha}} - \min x^{\overline{\alpha}} \mathbf{e}_{\widetilde{\alpha}}.
\]

For all the closed point $\mathfrak{p}$ of $\Spec \K[z]$, having lifted the syzygies, we have that $\dim \big((I_z)_{(\mathfrak{p})}\big)_r = \dim I_r$ and $\dim \big((I_z)_{(\mathfrak{p})}\big)_{r+1} = \dim I_{r+1}$, hence by Gotzmann's Persistence Theorem (Theorem \ref{th:PersistenceTheorem}) the Hilbert polynomial of the fiber of the point $\mathfrak{p}$ is equal to the Hilbert polynomial of $\Proj \K[x]/I$, i.e. $p(t)$. For the generic point, after having localized in $(I_z)_{(0)}$, the flatness is ensured by Proposition \ref{prop:artinSyzygies}.

The same reasoning works under the hypothesis $y_1 \neq 0$, considering as Borel set $\overline{\mathscr{B}} = \mathscr{B}\setminus\mathcal{F}_{\alpha}(x^\alpha)\cup\mathcal{F}_{\alpha}(x^\beta)$. $x^\beta$ will be a minimal monomial of $\restrict{\overline{\mathscr{B}}}{k}$, $x^\alpha$ a maximal monomial of $\restrict{\overline{\mathscr{B}}^\mathcal{C}}{k}$ and the decreasing set $\mathcal{F}_{\beta}$ will coincide with $\mathcal{F}_{\alpha}$.
\end{proof}

\begin{example}\label{ex:liftingSyzygies}
Let us consider the ideal $I = (x_3^2,x_3x_2,x_3x_1,x_2^3)_{\geqslant 4} \subset \K[x_0,x_1,x_2,x_3]$. Its Hilbert polynomial is $p(t) = 3t+1$ with Gotzmann number equal to $4$. $x_3x_1^3$ is a minimal monomial of $\restrict{\{I_4\}}{1}$ with decreasing set $\mathcal{F}_{x_3x_1^3} = \{\mathrm{id},\down{1},2\down{1}\}$ Borel-admissible w.r.t. $x_2^2 x_1^2$,  maximal element of $\restrict{\{I_4\}^{\mathcal{C}}}{1}$. Let us consider the monomial $x_3x_1^2 x_0 \in I$. The Eliahou-Kervaire syzygies of $I$ involving it are
\[
\begin{split}
&x_3 \cdot x_3x_1^2 x_0 = \dec{x_3^2x_1^2}{x_0}{I}{} \quad \Rightarrow \quad x_3 \mathbf{e}_{x_3x_1^2 x_0} - x_0 \mathbf{e}_{x_3^2x_1^2},\\
&x_2 \cdot x_3x_1^2 x_0 = \dec{x_3 x_2 x_1^2}{x_0}{I}{} \quad \Rightarrow \quad x_2 \mathbf{e}_{x_3x_1^2 x_0} - x_0 \mathbf{e}_{x_3 x_2 x_1^2},\\
&x_1 \cdot x_3x_1^2 x_0 = \dec{x_3x_1^3}{x_0}{I}{} \quad \Rightarrow \quad x_3 \mathbf{e}_{x_3x_1^2 x_0} - x_0 \mathbf{e}_{x_3x_1^3}.
\end{split}
\]
Set $J = \langle \{I_4\}\setminus\{x_3x_1^3,x_3x_1^2 x_0,x_3x_1x_0^2\}\cup\{x_3x_1^3+z\, x_2^2x_1^2,x_3x_1^2 x_0+z\, x_2^2x_1 x_0,x_3x_1x_0^2+z\, x_2^2x_0^2\}\rangle$, the lifted syzygies involving $x_3x_1^2 x_0+z\, x_2^2x_1 x_0$ are
\[
\begin{split}
x_3 \cdot (x_3x_1^2 x_0+z\, x_2^2x_1 x_0) &{} = \dec{x_3^2 x_1^2}{x_0}{I}{} + z\dec{x_3x_2^2x_1}{x_0}{I}{}\quad \Rightarrow\\
& x_3 \mathbf{f}_{x_3x_1^2 x_0} - x_0 \mathbf{f}_{x_3^2x_1^2} - z x_0 \mathbf{f}_{x_3x_2^2x_1},\\
x_2 \cdot (x_3x_1^2 x_0+z\, x_2^2x_1 x_0) &{} = \dec{x_2x_3 x_1^2}{x_0}{I}{} + z\dec{x_2^3x_1}{x_0}{I}{}\quad \Rightarrow\\
& x_2 \mathbf{f}_{x_3x_1^2 x_0} - x_0 \mathbf{f}_{x_2x_3x_1^2} - z x_0 \mathbf{f}_{x_2^3x_1},\\
x_1 \cdot (x_3x_1^2 x_0+z\, x_2^2x_1 x_0) &{} = x_0 \cdot (x_3x_1^3 + z x_2^2 x_1^2)\quad \Rightarrow\quad x_3 \mathbf{f}_{x_3x_1^2 x_0} - x_0 \mathbf{f}_{x_3x_1^3}.
\end{split}
\]
\end{example}

\begin{theorem}\label{th:curveGeneralDim}
Let $\Hilb{n}{p(t)}$ be the Hilbert scheme parametrizing subschemes of $\PP^n$ with Hilbert polynomial $p(t)$, whose Gotzmann number is $r$. Consider two Borel-fixed ideals $I$ and $J$ such that
\begin{enumerate}[(a)]
\item $I$ and $J$ define two $\K$-rational points of the Hilbert functor, i.e. $\Proj \K[x]/I,\Proj \K[x]/J$ in $\hilb{n}{p(t)}(\K)$;
\item there exist a minimal monomial $x^\alpha$ of $\restrict{\{I_r\}}{k}$ and a maximal monomial $x^\beta$ of $\restrict{\{I_r\}^{\mathcal{C}}}{k}$ such that $\min x^\alpha = \min x^\beta = k$, there does not exist any $\down{j}$ for which $\down{j}(x^\alpha) = x^\beta$, the decreasing set $\mathcal{F}_{\alpha}$ of $x^\alpha$ is Borel-admissible w.r.t. $x^\beta$ and
\[
\{J_r\} = \{I_r\} \setminus \mathcal{F}_{\alpha}(x^\alpha) \cup \mathcal{F}_{\alpha}(x^\beta). 
\]
\end{enumerate}
Then, there is a rational curve $\mathcal{C}: \PP^1 \rightarrow \Hilb{n}{p(t)}$ having two fibers corresponding to the $\K$-rational points defined by $I$ and $J$.

Moreover, if we consider the construction of the Hilbert scheme $\Hilb{n}{p(t)}$ as subscheme of the Grassmannian $\Grass{q(r)}{N}{\K}$, where $N = \binom{n+r}{n}$ and $q(r) = N-p(r)$, described in Chapter \ref{ch:HilbertScheme}, the degree of the curve $\mathcal{C}$ via the  Pl\"ucker embedding \eqref{eq:pluckerEmbedSubspace}  $\mathscr{P}: \Grass{q(r)}{N}{\K} \rightarrow \PP\big[\wedge^{q(r)} \K[x]_r\big]$ is $\left\vert\mathcal{F}_{\alpha}\right\vert$ and $\mathcal{C}$ is isomorphic to the rational normal curve in a convenient subspace of $\PP\big[\wedge^{q(r)} \K[x]_r\big]$.
\end{theorem}
\begin{proof}
The Borel sets $\{I_r\}$ and $\{J_r\}$ realize the hypothesis of Theorem \ref{th:flatGeneralDim}, hence the ideal
\[
\mathcal{I} = \left\langle \{I_r\} \setminus \mathcal{F}_{\alpha} \cup \{y_0 F(x^\alpha) + y_1 F(x^\beta) \ \vert\ F \in \mathcal{F}_{\alpha}\}\right\rangle \subset \K[y_0,y_1][x]
\]
parametrizes a flat family with Hilbert polynomial $p(t)$. By the universal property of the Hilbert scheme (Proposition \ref{prop:universalFamily}, see also \cite[Chapter III Remark 9.8.1]{Hartshorne}) the flat family $\Proj \K[y_0,y_1][x]/\mathcal{I} \rightarrow \PP^1$ determines uniquely a map $\mathcal{C}:\PP^1 \rightarrow \Hilb{n}{p(t)}$ and by construction the fibers over the points $[0:1]$ and $[1:0]$ correspond to the points defined by $I$ and $J$.

For the second part of the theorem, let us think to the submodule defined by $\mathcal{I}_r$ in $(\K[y_0,y_1])[x]_r$. As usual we can represent it by a matrix with $q(r)$ rows and $N$ columns with maximal rank; assuming to order the columns with the monomials of $\{I_r\}\setminus\mathcal{F}_{\alpha}(x^\alpha)$, the monomials of $\mathcal{F}_{\alpha}(x^\alpha)$, then the monomials $\mathcal{F}_{\alpha}(x^\beta)$ and finally the remaining ones, the subspace $\mathcal{I}_r$ is represented by the matrix
\small
\begin{align*}
& \qquad \quad\ \{I_r\}\setminus\mathcal{F}_{\alpha}(x^\alpha) \qquad \qquad \mathcal{F}_{\alpha}(x^\alpha) \qquad \qquad \mathcal{F}_{\alpha}(x^\beta) & \\
&\left(
\begin{array}{ccccc | ccc | ccc | cccccccc}
1 &  & \ldots & & 0 & & \vdots & & & \vdots & &   & & & & & & &\\
 & \ddots & & &  & & & & & & &   & & \vdots & & & & &\\
\vdots &  & 1 & & \vdots & \ldots & 0 & \ldots & \ldots & 0 & \ldots &   & \ldots  & 0 & \ldots & & & &\\
 &&& \ddots & & & & & & & &    & & \vdots & & & & &\\
 0& & \ldots & & 1 & & \vdots & & & \vdots & & & & & & & \vdots & &\\
 \cline{1-11}
  &  & \vdots & &  & y_0 & & & y_1 & & & & & & & \ldots & 0 & \ldots &\\
 \ldots &  & 0 & & \ldots & & \ddots & & & \ddots & & & & & & &\vdots &  &\\
 &  & \vdots & &  & & &  y_0 & & & y_1 & & & & & & & &\\
\end{array}
\right)&
\end{align*}
\normalsize
The Pl\"ucker coordinates via $\mathscr{P}$ are the $q(r)$-minors of this matrix. The non-vanishing coordinates correspond to the submatrices obtained by all columns associated to the monomials of $\{I_r\}\setminus\mathcal{F}_{\alpha}(x^\alpha)$ plus $\vert\mathcal{F}_{\alpha}\vert$ columns chosen among those associated to the monomials in $\mathcal{F}_{\alpha}(x^\alpha) \cup \mathcal{F}_{\alpha}(x^\beta)$, with the constraint of picking a monomial from each pair $(F(x^\alpha),F(x^\beta))$. Called $\HH$ the set of indices corresponding to the columns labeled with the monomials in $\{I_r\}\setminus\mathcal{F}_{\alpha}(x^\alpha)$, to any subset $\JJ$ of columns in $\mathcal{F}_{\alpha}(x^\alpha)$ the subset $\overline{\JJ}$ of columns of $\mathcal{F}_{\alpha}(x^\beta)$ giving a nonzero minor is uniquely associated and the Pl\"ucker coordinate $\Delta_{\HH\cup\JJ\cup\overline{\JJ}}$ is equal to $y_0^{\vert\JJ\vert} y_1^{\vert\overline{\JJ}\vert}$. Hence there are $2^{\vert\mathcal{F}_{\alpha}\vert}$ nonzero coordinates, and for any pair $\JJ,\JJ'$ of set of columns of $\mathcal{F}_\alpha(x^\alpha)$ such that $\vert\JJ\vert = \vert\JJ'\vert$, the Pl\"ucker coordinates $\Delta_{\HH\cup\JJ\cup\overline{\JJ}}$ and $\Delta_{\HH\cup\JJ'\cup\overline{\JJ}'}$ coincide, so that cutting with a suitable sets of hyperplanes we can obtain the curve $\PP^1 \rightarrow \PP^{\vert\mathcal{F}_{\alpha}\vert}$
\[
[y_0 : y_1] \longmapsto \left[y_0^{\vert\mathcal{F}_{\alpha}\vert}: y_0^{\vert\mathcal{F}_{\alpha}\vert-1}y_1 : \ldots : y_0 y_1^{\vert\mathcal{F}_{\alpha}\vert-1}: y_1^{\vert\mathcal{F}_{\alpha}\vert}\right]. \qedhere
\]
\end{proof}

\begin{example}
Let us consider again the ideal $I = (x_3^2,x_3x_2,x_3x_1,x_2^3)_{\geqslant 4}$ introduced in Example \ref{ex:liftingSyzygies} and the rational deformation defined by the ideal
\[
\begin{split}
\mathcal{I} = {}& \big\langle \{I_4\}\setminus\{x_3x_1^3,x_3x_1^2 x_0,x_3x_1x_0^2\}\\
&{}\cup\{y_0\,x_3x_1^3+y_1\, x_2^2x_1^2,y_0\, x_3x_1^2 x_0+y_1\, x_2^2x_1 x_0,y_0\, x_3x_1x_0^2+y_1\, x_2^2x_0^2\}\big\rangle.
\end{split}
\]
$I$ defines a point of the Hilbert scheme \gls{Hilb3tp1P3} that we embed in a projective space by means of the Pl\"ucker embedding of the Grassmannian $\mathscr{P}:\Grass{22}{35}{\K} \rightarrow \PP^{\binom{35}{22}-1}$. Considered the correspondence \eqref{eq:monomialCorrespondence} between indices from 1 up to $35$ and multiindices defining monomials in $\K[x_0,x_1,x_2,x_3]_4$, we call $\HH$ the set containing the indices associated to the monomials in $\{I_4\}\setminus\{x_3x_1^3,x_3x_1^2 x_0,x_3x_1x_0^2\}$ and moreover we have
\[
\begin{array}{lcl}
13 \leftrightarrow x_3x_1^3, & & 12 \leftrightarrow x_2^2 x_1^2,\\
23 \leftrightarrow x_3x_1^2 x_0, & & 22 \leftrightarrow x_2^2 x_1 x_0,\\
29 \leftrightarrow x_3x_1 x_0, & & 28 \leftrightarrow x_2^2 x_0^2.\\
\end{array}
\]
In order for obtaining a nonzero Pl\"ucker coordinate, we have to pick a monomial from each binomial, so there are $2^3$  non-vanishing Pl\"ucker coordinates:
\[
\begin{split}
& \Delta_{\HH \cup (13,23,29)} = y_0^3,\qquad \Delta_{\HH \cup (12,23,29)} = \Delta_{\HH \cup (13,22,29)} = \Delta_{\HH \cup (13,23,28)} = y_0^2 y_1,\\
& \Delta_{\HH \cup (12,22,29)} = \Delta_{\HH \cup (12,23,28)} =  \Delta_{\HH \cup (13,22,28)} = y_0 y_1^2,\qquad
 \Delta_{\HH\cup(12,22,28)} = y_1^3.
\end{split}
\]
Finally in the projective space $\PP^3$ obtained by $\PP^{\binom{35}{22}-1}$ cutting with the hyperplanes
\[
\begin{cases}
\Delta_{\HH \cup (13,22,29)} - \Delta_{\HH \cup (12,23,29)} = 0,\\
\Delta_{\HH \cup (13,23,28)} - \Delta_{\HH \cup (12,23,29)} = 0,\\
\Delta_{\HH \cup (12,23,28)} - \Delta_{\HH \cup (12,22,29)} = 0,\\
\Delta_{\HH \cup (13,22,28)} - \Delta_{\HH \cup (12,22,29)} = 0,\\
\Delta_\II = 0,\quad \II \text{ not appearing in the 8 coordinates listed above}
\end{cases}
\]
the curve is exactly the rational normal curves of degree 3.
\end{example}

\begin{definition}\index{Borel rational deformation}\index{Borel rational curve}\index{Borel degeneration}
Given two Borel-fixed ideals $I$ and $J$ that verify the hypothesis of Theorem \ref{th:flatGeneralDim} and Theorem \ref{th:curveGeneralDim}, we call \emph{Borel rational deformation} the deformation defined by the ideal \eqref{eq:familyIdeals} and \emph{Borel rational curve} the corresponding curve on the Hilbert scheme. Moreover we will call $J$ a \emph{Borel degeneration}\index{Borel degeneration} of $I$ (and viceversa).
\end{definition}

\begin{remark}\index{graph of monomials ideals}
A Borel rational deformation between two Borel-fixed ideals $I$ and $\overline{I}$ corresponds also to an edge connecting the vertices $I$ and $\overline{I}$ in the \emph{graph of monomial ideals} introduced by Altmann and Sturmfels in \cite{AltmannSturmfels}. The vertices of this graph are the monomial ideals in $\K[x]$ and two ideals $I,\overline{I}$ are connected by an edge if there exists an ideal $J$ such that the set of all its initial ideals (w.r.t. all term orderings) is $\{I,\overline{I}\}$.

Let $J = \mathcal{I}\vert_{[1:1]}$ be the ideal of the family described in \eqref{eq:familyIdeals} defining the deformation from $I$ to $\widetilde{I}$. By construction $x^\alpha$ and $x^\beta$ are not comparable w.r.t. the Borel partial order $\leq_B$, so for any term orderings $\sigma$, $x^\alpha >_{\sigma} x^\beta$ or $x^\alpha <_{\sigma} x^\beta$. 

Let $k = \min x^\alpha = \min x^\beta$ and let $s = \max \{\mu_k \ \vert\ \mu_1\down{1}\circ\cdots\circ\mu_k \down{k} = F \in \mathcal{F}_{\alpha}\}$. $x_k^s \mid \text{gcd}(x^\alpha,x^\beta)$, then $x^\alpha = x^{\overline{\alpha}}  x_k^s$, $x^\beta = x^{\overline{\beta}} x_k^s$ and $F(x^\alpha) = x^{\overline{\alpha}}\, F(x_k^s)$, $F(x^\beta) = x^{\overline{\beta}}\, F(x_k^s)$. Finally
\[
 x^\alpha = x^{\overline{\alpha}} x_k^s \gtrless_{\sigma} x^{\overline{\beta}} x_k^s = x^\beta \ \Rightarrow\ x^{\overline{\alpha}} \gtrless_{\sigma} x^{\overline{\beta}} \ \Rightarrow\ F(x^\alpha) = x^{\overline{\alpha}} F(x_k^s) \gtrless_{\sigma} x^{\overline{\beta}} F(x_k^s) = F(x^\beta),
\]
for each $F \in \mathcal{F}_\alpha$. So the ideal $J$ has $\{I,\overline{I}\}$ as set of possible initial ideals.
\end{remark}

An algorithm for computing all the possible Borel rational deformations (or degenerations) of a Borel-fixed ideal naturally arises from Theorem \ref{th:flatGeneralDim}. In Algorithm \ref{alg:allDeformations} and Algorithm \ref{alg:allDegenerations} there are their pseudocode descriptions: the key point is to look in any restriction of the Borel set $\{I_m\}$ for minimal and maximal elements with the required property.

\begin{algorithm}
\caption[Pseudocode description of the algorithm computing all the possible Borel rational deformations of a given Borel-fixed ideal.]{Pseudocode description of the algorithm computing all the possible Borel rational deformations of a given Borel-fixed ideal $I$. For more details see Appendix \ref{ch:HSCpackage}.}
\label{alg:allDeformations}
\begin{algorithmic}[1]
\STATE $\textsc{BorelRationalDeformations}(I,m)$
\REQUIRE $I\subset\K[x]$, Borel-fixed ideal.
\REQUIRE $m$, positive integer.
\ENSURE the set of Borel rational deformations involving $I$, constructed considering the Borel set $\{I_m\}$.
\STATE $\textsf{deformations} \leftarrow \emptyset$;
\FOR{$k=0,\ldots,n-1$}
\STATE $\textsf{minimalMonomials} \leftarrow \textsc{MinimalElements}\big(\restrict{\{I_m\}}{k}\big)$;
\STATE $\textsf{maximalMonomials} \leftarrow \textsc{MaximalElements}\big(\restrict{\{I_m\}^{\mathcal{C}}}{k}\big)$;
\FORALL{$x^\alpha \in \textsf{minimalMonomials}$}
\FORALL{$x^\beta \in \textsf{maximalMonomials}$}
\IF{$\nexists\ \down{j},\ j > k$ s.t. $\down{j}(x^\alpha)=x^\beta$}
\STATE $\mathcal{F}_{\alpha} \leftarrow \textsc{DecreasingSet}(\{I_m\},x^\alpha)$;
\IF{$\mathcal{F}_{\alpha}$ is Borel-admissible w.r.t. $x^\beta$}
\STATE $\mathcal{I} \leftarrow \left\langle \{I_m\}\setminus \mathcal{F}_{\alpha}(x^\alpha)\cup \left\{y_0F(x^\alpha)+y_1F(x^\beta)\ \vert\ F \in \mathcal{F}_{\alpha}\right\}\right\rangle$;
\STATE $\textsf{deformations} \leftarrow \textsf{deformations} \cup \left\{\mathcal{I}\right\}$;
\ENDIF
\ENDIF
\ENDFOR
\ENDFOR
\ENDFOR
\RETURN $\textsf{deformations}$;
\end{algorithmic}

\bigskip

\begin{algorithmic}
\STATE $\textsc{DecreasingSet}(\mathscr{B},x^\alpha)$
\REQUIRE $\mathscr{B}$, a Borel set.
\REQUIRE $x^\alpha$, a monomial of $\mathscr{B}$. 
\ENSURE the set of decreasing moves going from $x^\alpha$ to any other monomial $x^\gamma \in \mathscr{B}$, i.e. such that $x^\alpha >_B x^\gamma$.
\end{algorithmic}
\end{algorithm}

\begin{algorithm}
\caption[Pseudocode description of the algorithm computing all the possible Borel rational degeneration of a given Borel-fixed ideal $I$.]{Pseudocode description of the algorithm computing all the possible Borel rational degeneration of a given Borel-fixed ideal $I$. For more details see Appendix \ref{ch:HSCpackage}.}
\label{alg:allDegenerations}
\begin{algorithmic}[1]
\STATE $\textsc{BorelRationalDegenerations}(I,m)$
\REQUIRE $I\subset\K[x]$, Borel-fixed ideal.
\REQUIRE $m$, positive integer.
\ENSURE the set of Borel rational degenerations of $I$, constructed considering the Borel set $\{I_m\}$.
\STATE $\textsf{deformations} \leftarrow \textsc{BorelRationalDeformations}(I,m)$;
\STATE $\textsf{degenerations} \leftarrow \emptyset$;
\FORALL{$\mathcal{I} \in \textsf{deformations}$}
\STATE $\textsf{degenerations} \leftarrow \textsf{degenerations} \cup \left\{\mathcal{I}\vert_{[0:1]}\right\}$;
\ENDFOR
\RETURN $\textsf{deformations}$;
\end{algorithmic}
\end{algorithm}

\begin{example}\label{ex:multipleDef}
 Let us simulate an execution of Algorithm \ref{alg:allDeformations} on the ideal $I = (x_3^2,x_3 x_2^2, x_3 x_2 x_1,x_2^4,x_2^3 x_1,x_2^2 x_1^2) \subset \K[x_0,x_1,x_2,x_3]$ (see Figure \ref{fig:exampleDeformations_a} for its Green's diagram). The corresponding subscheme $\Proj \K[x_0,x_1,x_2,x_3]/I$ has Hilbert polynomial $p(t)=3t+5$, whose Gotzmann number is $8$, so we consider the Borel set $\{I_8\}$.
\begin{description}
\item[$k=0.$] The minimal monomials of $\{I_8\}$ are $\{x_3^2 x_0^6,x_3 x_2 x_1 x_0^5,x_2^2 x_1^2 x_0^4\}$ and the maximal elements of $\{I_8\}^{\mathcal{C}}$ are $\{x_3 x_2 x_0^6,x_2^3 x_0^5\}$. The pairs $(x_3x_2x_1x_0^5,x_3x_2x_0^6)$ and $(x_3^2x_0^6,x_3x_2x_0^6)$ are discarded because $\down{1}(x_3x_2x_1x_0^5)= x_3x_2x_0^6$ and $\down{3}(x_3^2x_0^6)=x_3x_2x_0^6$, so we have four possibilities. In the case $k=0$, the Borel-admissibility derives directly from the fact that the decreasing set contains only the identity move, so we do not need to check further conditions.
\begin{description}
\item[(Figure \ref{fig:exampleDeformations_b})] $x_3^2 x_0^6$ and $x_2^3 x_0^5$ define the Borel rational deformation
\[
 \left\langle \{I_8\} \setminus \{x_3^2 x_0^6\} \cup \{y_0 \, x_3^2 x_0^6 + y_1 x_2^3 x_0^5\} \right\rangle \subset \K[y_0,y_1][x_0,x_1,x_2,x_3]
\]
and the fiber over the point $[0:1]$ is the ideal $(x_3^3,x_3^2 x_2, x_3 x_2^2, x_2^3,x_3^2 x_1,$ $x_3 x_2 x_1,x_2^2 x_1^2)_{\geqslant 8}$.
\item[(Figure \ref{fig:exampleDeformations_c})] $x_3 x_2 x_1 x_0^5$ and $x_2^3 x_0^5$ define the deformation
\[
 \left\langle \{I_8\} \setminus \{x_3 x_2 x_1 x_0^5\} \cup \{y_0 \, x_3 x_2 x_1 x_0^5 + y_1 x_2^3 x_0^5\} \right\rangle \subset \K[y_0,y_1][x_0,x_1,x_2,x_3]
\]
with fiber over the point $[0:1]$ equal to $(x_3^2,x_3 x_2^2,x_2^3,x_3 x_2 x_1^2,x_2^2 x_1^2)_{\geqslant 8}$.
\item[(Figure \ref{fig:exampleDeformations_d})] $x_2^2 x_1^2 x_0^4$ and $x_3 x_2 x_0^6$ define the deformation
\[
\left\langle \{I_8\} \setminus \{x_2^2 x_1^2 x_0^4\} \cup \{y_0 \, x_2^2 x_1^2 x_0^4 + y_1 x_3 x_2 x_0^6\} \right\rangle \subset \K[y_0,y_1][x_0,x_1,x_2,x_3]
\]
and the fiber over the point $[0:1]$ is the ideal $(x_3^2,x_3 x_2,x_2^4,x_2^3 x_1, x_2^2 x_1^3)_{\geqslant 8}$.
\item[(Figure \ref{fig:exampleDeformations_e})] $x_2^2 x_1^2 x_0^4$ and $x_2^3 x_0^5$ define the deformation
\[
\left\langle \{I_8\} \setminus \{x_2^2 x_1^2 x_0^4\} \cup \{y_0 \, x_2^2 x_1^2 x_0^4 + y_1 x_2^3 x_0^5\} \right\rangle \subset \K[y_0,y_1][x_0,x_1,x_2,x_3]
\]
with the ideal $(x_3^2,x_3 x_2^2,x_2^3,x_3 x_2 x_1,x_2^2 x_1^3)_{\geqslant 8}$ as fiber over $[0:1]$.
\end{description}
\item[$k=1.$] $\restrict{\{I_8\}}{1}$ has only $x_2^2 x_1^6$ as minimal element and $x_3x_1^7$ is the only maximal element of $\restrict{\{I_8\}^{\mathcal{C}}}{1}$. The decreasing set of $x_2^2 x_1^6$ is
\[
\mathcal{F}_{x_2^2 x_1^6} = \{\mathrm{id},\down{1},2\down{1},3\down{1},4\down{1}\},
\]
Borel-admissible w.r.t. $x_3 x_1^7$. The Borel rational deformation defined by $x_2^2 x_1^6$ and $x_3x_1^7$ is
\[
\begin{split}
&\big\langle \{I_8\} \setminus \{x_2^2 x_1^6,x_2^2 x_1^5 x_0,x_2^2 x_1^4x_0^2,x_2^2 x_1^3 x_0^3,x_2^2 x_1^2 x_0^4\} \cup \{y_0x_2^2 x_1^6 + y_1 x_3 x_1^7,  \\
&\ y_0x_2^2 x_1^5 x_0 + y_1 x_3 x_1^6 x_0, y_0x_2^2 x_1^4 x_0^2 + y_1 x_3 x_1^5 x_0^3, y_0x_2^2 x_1^3 x_0^3 + y_1 x_3 x_1^4 x_0^3,\\
&\ y_0x_2^2 x_1^2 x_0^4 + y_1 x_3 x_1^3 x_0^4\} \big\rangle \subset \K[y_0,y_1][x_0,x_1,x_2,x_3]
\end{split}
\]
and the fiber over $[0:1]$ is the ideal $(x_3^2,x_3 x_2^2,x_3 x_2 x_1,x_2^4,x_2^3 x_1,x_3 x_1^3)_{\geqslant 8}$ (Figure \ref{fig:exampleDeformations_f}).
\end{description}
\end{example}

\vspace*{\stretch{1}}

\begin{figure}[H]
\begin{center}
\captionsetup[subfloat]{singlelinecheck=false,format=hang,justification=centering}
\subfloat[][$(x_3^2,x_3 x_2^2, x_3 x_2 x_1,x_2^4,x_2^3 x_1,x_2^2 x_1^2)$.]{\label{fig:exampleDeformations_a}
\begin{tikzpicture}[scale=0.5]
\tikzstyle{ideal}=[circle,draw=black,fill=black,inner sep=3.5pt]
\tikzstyle{quotient}=[circle,draw=black,thick,inner sep=3.5pt]
\tikzstyle{mixed}=[circle,draw=black,thick,inner sep=2pt]

\tikzstyle{empty}=[circle,draw=white,thick,inner sep=3.5pt]

\node at (-5,0) [] {};
\node at (5,0) [] {};

\node (r0) at (0,0) [quotient] {};

\node (r10) at (-1,-1.7) [quotient] {};
\node (r11) at (1,-1.7) [quotient] {};

\node (r20) at (-2,-3.4) [ideal] {};
\node (r21) at (0,-3.4) [mixed] {\tiny $1$};
\node (r22) at (2,-3.4) [mixed] {\tiny $2$};

\node (r32) at (1,-5.1) [ideal] {};
\node (r33) at (3,-5.1) [mixed] {\tiny $1$};

\node (r40) at (-4,-6.8) [empty] {};
\node (r44) at (4,-6.8) [ideal] {};
\end{tikzpicture}
}
\qquad\qquad
\subfloat[][$(x_3^3,x_3^2 x_2, x_3 x_2^2, x_2^3,x_3^2 x_1,x_3 x_2 x_1,x_2^2 x_1^2)$.]{\label{fig:exampleDeformations_b}
\begin{tikzpicture}[scale=0.5]
\tikzstyle{ideal}=[circle,draw=black,fill=black,inner sep=3.5pt]
\tikzstyle{quotient}=[circle,draw=black,thick,inner sep=3.5pt]
\tikzstyle{mixed}=[circle,draw=black,thick,inner sep=2pt]
\tikzstyle{empty}=[circle,draw=white,thick,inner sep=3.5pt]

\node at (-5,0) [] {};
\node at (5,0) [] {};

\node (r0) at (0,0) [quotient] {};

\node (r10) at (-1,-1.7) [quotient] {};
\node (r11) at (1,-1.7) [quotient] {};

\node (r20) at (-2,-3.4) [mixed] {\tiny $1$};
\node (r21) at (0,-3.4) [mixed] {\tiny $1$};
\node (r22) at (2,-3.4) [mixed] {\tiny $2$};

\node (r30) at (-3,-5.1) [ideal] {};
\node (r31) at (-1,-5.1) [ideal] {};
\node (r32) at (1,-5.1) [ideal] {};
\node (r33) at (3,-5.1) [ideal] {};

\node (r40) at (-4,-6.8) [empty] {};
\node (r44) at (4,-6.8) [empty] {};
\end{tikzpicture}
} \\

\bigskip

\subfloat[][$(x_3^2,x_3 x_2^2,x_2^3,x_3 x_2 x_1^2,x_2^2 x_1^2)$.]{\label{fig:exampleDeformations_c}
\begin{tikzpicture}[scale=0.5]
\tikzstyle{ideal}=[circle,draw=black,fill=black,inner sep=3.5pt]
\tikzstyle{quotient}=[circle,draw=black,thick,inner sep=3.5pt]
\tikzstyle{mixed}=[circle,draw=black,thick,inner sep=2pt]
\tikzstyle{empty}=[circle,draw=white,thick,inner sep=3.5pt]

\node at (-5,0) [] {};
\node at (5,0) [] {};

\node (r0) at (0,0) [quotient] {};

\node (r10) at (-1,-1.7) [quotient] {};
\node (r11) at (1,-1.7) [quotient] {};

\node (r20) at (-2,-3.4) [ideal] {};
\node (r21) at (0,-3.4) [mixed] {\tiny $3$};
\node (r22) at (2,-3.4) [mixed] {\tiny $2$};

\node (r32) at (1,-5.1) [ideal] {};
\node (r33) at (3,-5.1) [ideal] {};

\node (r40) at (-4,-6.8) [empty] {};
\node (r44) at (4,-6.8) [empty] {};
\end{tikzpicture}
}
\qquad\qquad
\subfloat[][$(x_3^2,x_3 x_2,x_2^4,x_2^3 x_1, x_2^2 x_1^3)$.]{\label{fig:exampleDeformations_d}
\begin{tikzpicture}[scale=0.5]
\tikzstyle{ideal}=[circle,draw=black,fill=black,inner sep=3.5pt]
\tikzstyle{quotient}=[circle,draw=black,thick,inner sep=3.5pt]
\tikzstyle{mixed}=[circle,draw=black,thick,inner sep=2pt]
\tikzstyle{empty}=[circle,draw=white,thick,inner sep=3.5pt]

\node at (-5,0) [] {};
\node at (5,0) [] {};

\node (r0) at (0,0) [quotient] {};

\node (r10) at (-1,-1.7) [quotient] {};
\node (r11) at (1,-1.7) [quotient] {};

\node (r20) at (-2,-3.4) [ideal] {};
\node (r21) at (0,-3.4) [ideal] {};
\node (r22) at (2,-3.4) [mixed] {\tiny $3$};

\node (r33) at (3,-5.1) [mixed] {\tiny $1$};

\node (r40) at (-4,-6.8) [empty] {};
\node (r44) at (4,-6.8) [ideal] {};
\end{tikzpicture}
} \\

\bigskip

\subfloat[][$(x_3^2,x_3 x_2^2,x_2^3,x_3 x_2 x_1,x_2^2 x_1^3)$.]{\label{fig:exampleDeformations_e}
\begin{tikzpicture}[scale=0.5]
\tikzstyle{ideal}=[circle,draw=black,fill=black,inner sep=3.5pt]
\tikzstyle{quotient}=[circle,draw=black,thick,inner sep=3.5pt]
\tikzstyle{mixed}=[circle,draw=black,thick,inner sep=2pt]
\tikzstyle{empty}=[circle,draw=white,thick,inner sep=3.5pt]

\node at (-5,0) [] {};
\node at (5,0) [] {};

\node (r0) at (0,0) [quotient] {};

\node (r10) at (-1,-1.7) [quotient] {};
\node (r11) at (1,-1.7) [quotient] {};

\node (r20) at (-2,-3.4) [ideal] {};
\node (r21) at (0,-3.4) [mixed] {\tiny $1$};
\node (r22) at (2,-3.4) [mixed] {\tiny $3$};

\node (r32) at (1,-5.1) [ideal] {};
\node (r33) at (3,-5.1) [ideal] {};

\node (r40) at (-4,-6.8) [empty] {};
\node (r44) at (4,-6.8) [empty] {};
\end{tikzpicture}
} 
\qquad\qquad
\subfloat[][$(x_3^2,x_3 x_2^2,x_3 x_2 x_1,x_2^4,x_2^3 x_1,x_3 x_1^3)$.]{\label{fig:exampleDeformations_f}
\begin{tikzpicture}[scale=0.5]
\tikzstyle{ideal}=[circle,draw=black,fill=black,inner sep=3.5pt]
\tikzstyle{quotient}=[circle,draw=black,thick,inner sep=3.5pt]
\tikzstyle{mixed}=[circle,draw=black,thick,inner sep=2pt]
\tikzstyle{empty}=[circle,draw=white,thick,inner sep=3.5pt]

\node at (-5,0) [] {};
\node at (5,0) [] {};

\node (r0) at (0,0) [quotient] {};

\node (r10) at (-1,-1.7) [mixed] {\tiny $3$};
\node (r11) at (1,-1.7) [quotient] {};

\node (r20) at (-2,-3.4) [ideal] {};
\node (r21) at (0,-3.4) [mixed] {\tiny $1$};
\node (r22) at (2,-3.4) [quotient] {};

\node (r32) at (1,-5.1) [ideal] {};
\node (r33) at (3,-5.1) [mixed] {\tiny $1$};

\node (r40) at (-4,-6.8) [empty] {};
\node (r44) at (4,-6.8) [ideal] {};
\end{tikzpicture}
}
\end{center}
\caption[Green's diagrams of the Borel-fixed ideals obtained by a Borel degeneration from ${(x_3^2,x_3 x_2^2, x_3 x_2 x_1,x_2^4,x_2^3 x_1,x_2^2 x_1^2) \subset \K[x_0,x_1,x_2,x_3]}$.]{\label{fig:exampleDeformations} Green's diagrams of the Borel-fixed ideals introduced in Example \ref{ex:multipleDef}.}
\end{figure}

\vspace*{\stretch{1}}

\newpage

\begin{remark}\index{Hilbert function}
We briefly discuss our method in relation with the Hilbert function and so with the results exposed by Mall in \cite{Mall}. Since our interest is oriented toward Hilbert schemes, where the crucial aspect is the Hilbert polynomial, we usually prefer to consider as ideal $I$ defining a subschemes with Hilbert polynomial $p(t)$ as generated in degree $r$, i.e. $I = I_{\geqslant r} = (I_r)$ such that the Hilbert function of $\K[x]/I$ is
\[
\textit{HF}_{\K[x]/I}(t) = \begin{cases}
\binom{n+t}{n},&\text{if}\ t < r,\\
p(t),&\text{if}\ t \geqslant r,
\end{cases}
\]
where $r$ is the Gotzmann number of the Hilbert polynomial $p(t)$. However, it is possible to adapt the technique to work also on the Hilbert function of the quotient modules defined by the saturation of such Borel ideals. We highlight this point starting from Example 3.6 and Example 3.9 of \cite{Mall}. Mall showed that the ideal $J=(x_3^2,x_3x_2,x_3x_1-x_2^2) \subset \K[x_0,x_1,x_2,x_3]$ has only two possible initial ideals (varying the term ordering): $I = (x_3^2,x_3x_2,x_3x_1,x_2^3)$ and $\overline{I}=(x_3^2,x_3x_2,x_2^2)$, so that these ideals are connected through two Gr\"obner deformations and moreover the quotient modules defined have the same Hilbert function $\textit{HF}_{\K[x]/I}(t) = \textit{HF}_{\K[x]/\overline{I}}(t) = \left(1,4,7,10,3t+1,\ldots\right)$.

In our perspective, since the Hilbert polynomial of $\K[x]/I$ is $p(t) = 3t+1$ with Gotzmann number $r=4$, we would consider the ideals $I_{\geqslant 4}$ and $\overline{I}_{\geqslant 4}$. Applying Algorithm \ref{alg:allDeformations} on the Borel set $\{\overline{I}_4\}$, we obtain the deformation
\[
\begin{split}
\overline{J} = &\big\langle\{\overline{I}_4\} \setminus \{x_2^2x_1^2,x_2^2x_1x_0,x_2^2x_0^2\} \\
&\cup \{y_0\, x_2^2x_1^2 + y_1\, x_3x_1^3,y_0\, x_2^2x_1x_0 + y_1\, x_3x_1^2x_2,y_0\, x_2^2x_0^2 + y_1\, x_3 x_1 x_0^2\}\big\rangle
\end{split}
\]
with fibers $I_{\geqslant 4}$ and $\overline{I}_{\geqslant 4}$.
Keeping in mind that to compute the saturation of a Borel-fixed ideal it is sufficient to put the smallest variable $x_0$ equal to 1, $I = (I_{\geqslant 4})^\sat $ and $\overline{I} = (\overline{I}_{\geqslant 4})^\sat$ must have the same Hilbert function because we are swapping pairs of monomials with the same power of the variable $x_0$. Furthermore we can deduce also the deformation between $I$ and $\overline{I}$: first of all we consider the saturation of the ideal $\left\langle\{\overline{I}_2\}_4 \setminus \{x_2^2x_1^2,x_2^2x_1x_0,x_2^2x_0^2\}\right\rangle = (x_3^2,x_3x_2,x_2^3)$ (the common part of $I$ and $\overline{I}$) and then we add to the generators the binomial $y_0\, x_2^2 + y_1\, x_3 x_1$
\[
\overline{J}^\sat = (x_3^2,x_3x_2,x_2^3,y_0\, x_2^2 + y_1\, x_3 x_1).
\]

Hence to have deformations preserving the Hilbert function of the saturation of the Borel-fixed ideals involved, the key point is to swap couples of monomials having the same power of the smallest variable (in our case $x_0$). This happens by construction whenever the deformation is defined by monomials $x^\alpha,x^\beta$ such that $\min x^\alpha = \min x^\beta > 0$, whereas it happens rarely if the deformation is ruled by a single couple of monomials as the following example shows. Let us consider the ideal $I_{\geqslant 4}$ and the ideal $\widetilde{I} = (x_3,x_2^4,x_2^3x_1)_{\geqslant 4}$, fibers of the deformation
\[
\left\langle\left\{I_4\right\} \setminus \left\{x_2^3 x_0\right\} \cup\left\{y_0 \, x_2^3 x_0 + y_1\, x_3x_0^3\right\}\right\rangle.
\]
The monomials involved in the deformation have not the same power of $x_0$, indeed
\[
\begin{split}
&\textit{HF}_{\K[x]/I}(t) = \left(1,4,7,10,3t+1,\ldots\right),\\
&\textit{HF}_{\K[x]/\widetilde{I}^\sat}(t) = \left(1,3,6,10,3t+1,\ldots\right).
\end{split}
\]
\end{remark}

\section{The connectedness of the Hilbert scheme}\label{sec:connected}
In this section, we will study consecutive deformations of Borel-fixed ideals and we want to introduced a way to control the \lq\lq direction\rq\rq\ toward which we move, for obtaining a technique similar to the one introduced by Peeva and Stillman \cite{PeevaStillman}. The deformations introduced in \cite{PeevaStillman} are affine and based on Gr\"obner basis tools, but anyhow the idea is to exchange monomials belonging to the ideal with monomials not belonging. The goal is to determine a sequence of deformations leading from any Borel-fixed ideal to the lexicographic ideal,\index{lexicographic ideal} so the choice of the monomials to exchange is governed by the $\mathtt{DegLex}$ term ordering.

Therefore we start by slightly modifying the strategy of Algorithm \ref{alg:allDeformations} and Algorithm \ref{alg:allDegenerations}, adding a term ordering $\sigma$, refinement of the Borel partial order $\leq_B$, that we will use to choose in a unique way one of the possible Borel rational deformations of an ideal. The point is to replace some monomials of the ideal with some others that are greater that them w.r.t. $\sigma$.

\begin{algorithm}[H]
\caption[Algorithm determining the Borel rational degeneration in the \lq\lq direction\rq\rq\ fixed by a term ordering.]{How to find a Borel rational degeneration with special \lq\lq direction\rq\rq\ of a fixed term ordering. For details see Appendix \ref{ch:HSCpackage}.}
\label{alg:deformation}
\begin{algorithmic}[1]
\STATE $\textsc{OrientedBorelRationalDegeneration}(I,\sigma)$
\REQUIRE $I\subset\K[x]$, Borel-fixed ideal.
\REQUIRE $\sigma$, term ordering, refinement of the Borel partial order $\leq_B$.
\ENSURE a saturated Borel-fixed ideal $J$, $\sigma$-Borel degeneration of $I$.\\ If $I$ has not a $\sigma$-Borel degeneration, the function returns the ideal itself.
\STATE $p(t) \leftarrow \textsc{HilbertPolynomial}(\K[x]/I)$;
\STATE $r \leftarrow \textsc{GotzmannNumber}\big(p(t)\big)$; 
\FOR{$k=0,\ldots,n-1$}
\STATE $x^\alpha \leftarrow \min_{\sigma} \restrict{\{I_r\}}{k}$;
\STATE $x^\beta \leftarrow \max_{\sigma} \restrict{\{I_r\}^{\mathcal{C}}}{k}$;
\IF{$x^\alpha <_{\sigma} x^\beta$}
\STATE $\mathcal{F}_{\alpha} \leftarrow \textsc{DecreasingSet}(\{I_r\},x^\alpha)$;
\IF{$\mathcal{F}_{\alpha}$ is Borel-admissible}
\RETURN $\left\langle\{I_r\}\setminus\mathcal{F}_{\alpha}(x^\alpha) \cup \mathcal{F}_{\alpha}(x^\beta)\right\rangle^{\sat}$;
\ENDIF
\ENDIF
\ENDFOR
\RETURN $I$;
\end{algorithmic}
\end{algorithm}

We remark that $x^\alpha$ and $x^\beta$ are surely a minimal element of $\restrict{\{I_r\}}{k}$ and a maximal element of $\restrict{\{I_r\}^{\mathcal{C}}}{k}$, by Remark \ref{rk:minMax}. Moreover the condition $x^\alpha <_\sigma x^\beta$ guarantees that $x^\beta$ can not be obtained from $x^\alpha$ by an elementary move $\down{j}$. In fact if such a move exists, $x^\alpha >_B x^\beta$ implies $x^\alpha >_\sigma x^\beta$ for any term ordering $\sigma$. 

The uniqueness of the deformation is imposed by making the algorithm returning the first deformation with the property $x^\alpha <_\sigma x^\beta$. If Algorithm \ref{alg:deformation} does not find any deformation, it returns the same ideal $I$ given as input; for instance this happens if we apply the algorithm to the lexicographic ideal with the $\mathtt{DegLex}$ term ordering. Indeed $x^\gamma >_{\mathtt{DegLex}} x^\delta, \forall\ x^\gamma \in \{I_m\},\ x^\delta \in \{I_m\}^{\mathcal{C}},\ \forall\ m$ by definition, thus the condition $x^\alpha <_\sigma x^\beta$ would always be false. 

\begin{definition}\label{def:preceqdef}
 Given a Borel-fixed ideal $I$ and a term ordering $\sigma$, we say that
\begin{itemize}
 \item\index{Borel degeneration!oriented} the ideal $J \neq I$ returned by $\textsc{OrientedBorelRationalDegeneration}(I,\sigma)$ (Algorithm \ref{alg:deformation}) is a \emph{$\sigma$-Borel rational degeneration} or simply a \emph{$\sigma$-Borel degeneration} of $I$ and we call \emph{$\sigma$-Borel rational deformation} of $I$ the Borel rational deformation having as fibers $I$ and $J$;
 \item\index{$\sigma$-endpoint} $I$ is a \emph{$\sigma$-endpoint} if $\textsc{OrientedBorelRationalDeformation}(I,\sigma)$ returns the same ideal $I$.
\end{itemize}
\end{definition}

After having determined a method similar to the one proposed in \cite{PeevaStillman}, we want to compare these two approaches, so we recall briefly the strategy and some notation of that paper. Given an ideal $I = I_{\geqslant r} \subset \K[x]$, Peeva and Stillman compute the monomials $x^\beta = \max_{\DegLex} \{I_r\}^{\mathcal{C}}$ (they call it \emph{first gap}) and $x^\alpha = \max_{\DegLex} \{x^\delta \in \{I_r\}\ \vert\ x^\delta <_{\DegLex} x^\beta\}$. If $I$ is the lexicographic ideal, the set $\{x^\delta \in \{I_r\}\ \vert\ x^\delta <_{\DegLex} x^\beta\}$ will be obviously empty. At this point they determine the set of monomials 
\[
T = \{x^{\overline{\alpha}} x^{\gamma}\in \{I_r\} \ \vert\ x^\alpha = x^{\overline{\alpha}}\cdot (\min x^\alpha)^{s},\ x^{\overline{\beta}} x^\gamma \notin \{I_r\} \text{ and } \max x^\gamma <_B \min x^{\overline{\beta}}\}   
\]
and they fix the monomial $x^{\overline{\alpha}}$ of minimal degree and the corresponding $x^{\overline{\beta}}$ of minimal degree. Let us denote them by $x^{\widetilde{\alpha}}$ and $x^{\widetilde{\beta}}$. They finally form the ideal
\[
 \widetilde{I} = \left\langle\{I_r\}\setminus\big\{x^{\widetilde{\alpha}} x^{\gamma} \in \{I_r\}\big\}\cup\big\{x^{\widetilde{\alpha}} x^{\gamma} - x^{\widetilde{\beta}} x^{\gamma} \ \vert\ x^{\widetilde{\alpha}} x^{\gamma} \in \{I_r\}\big\}\right\rangle,\quad \max x^\gamma \leq_B \min x^{\widetilde{\beta}}.
\]
The Borel ideal \lq\lq \DegLex-closer\rq\rq\ to the lexicographic ideal is $\GIN_{\DegLex}\big(\IN_{\DegLex}(\widetilde{I})\big)$.

The computation of a generic initial ideal\index{generic initial ideal} reveals an important difference between the two techniques: this choice of monomials involved in the substitution generally does not preserve the Borel condition, so that to restore it a computation of a $\GIN$ could be needed, as the following example shows.

\begin{example}\label{ex:PSbreakingBorel}
 Let us consider the ideal $I = (x_2^4,x_2^3x_1,x_2^2 x_1^2,x_2 x_1^3)_{\geqslant 7} \subset \K[x_0,x_1,x_2]$. The Hilbert polynomial is $p(t) = t+7$ with Gotzmann number 7. The first gap is $x_2^3x_0^4$ and the \DegLex-greatest monomial in $\{I_7\}$, \DegLex-smaller than the gap, is $x_2^2 x_1^5$.
\[
 T = \left\{x_2^2 x_1^2 x_0^3, x_2^2 x_1^3 x_0^2, x_2^2 x_1^4 x_0, x_2^2 x_1^5\right\}
\]
so that $x^{\widetilde{\alpha}} = x_2^2 x_1^2$ and $x^{\widetilde{\beta}} = x_2^3 x_0$. We construct the ideal 
\[
\widetilde{I} = \left\langle \{I_7\} \setminus\{x_2^2 x_1^2 x_0^3\}\cup\{x_2^2 x_1^2 x_0^3- x_2^3 x_0^4\} \right\rangle,
\]
and $\IN_{\DegLex}(\widetilde{I}) = (x_2^3,x_2x_1^3)_{\geqslant 7}$ is not Borel-fixed, because $x_2 x_1^3 x_0^3 \in \IN_{\DegLex}(\widetilde{I})$  and $\up{1}(x_2 x_1^3 x_0^3) = x_2^2 x_1^2 x_0^3 \notin \IN_{\DegLex}(\widetilde{I})$. Finally 
\[
\GIN_{\DegLex}\big(\IN_{\DegLex}(\widetilde{I})\big) = (x_2^3,x_2^2 x_1^2,x_2 x_1^4)_{\geqslant 7}
\]
is again Borel-fixed.

Using Algorithm \ref{alg:deformation} to compute the $\mathtt{DegLex}$-Borel degeneration of $I$, we have
\[
\min_{\DegLex} \{I_7\} =  x_2x_1^3x_0^3 <_{\DegLex} x_2^3 x_0^4 = \max_{\DegLex} \{I_7\}^{\mathcal{C}},
\]
and since the decreasing set of $x_2x_1^3x_0^3$ contains only the identity, the $\DegLex$-Borel degeneration of $I$ is the ideal $J = \left\langle \{I_7\}\setminus\{x_2x_1^3x_0^3\}\cup\{x_2^3 x_0^4\} \right\rangle = (x_2^3,x_2^2 x_1^2,x_2 x_1^4)_{\geqslant 7}$.
\end{example}

The choice of the first gap is very similar to the choice of the maximal monomial in $\{I_r\}^{\mathcal{C}}$, because in both cases we look for a \lq\lq greatest\rq\rq\ element. But the two techniques differ in the choice of the monomial inside the ideal: we look for minimal elements whereas Peeva and Stillman consider maximal elements lower than the gap. Hence we expect that in the cases in which the monomials chosen are the same, the deformation is almost equal, in the sense that the monomial ideal we reach is the same.

\begin{example}\label{ex:PSsameRational}
We consider the ideal $I = (x_3^2,x_3 x_2,x_2^3,x_2^2 x_1)_{\geqslant 5} \subset \K[x_0,x_1,x_2,x_3]$. The Hilbert polynomial is $p(t) = 3t+2$ with Gotzmann number 5. The first gap is $x_3 x_1^4$ and the greatest monomial of the ideal smaller than it is $x_2^5$. Then
\[
 T = \left\{ x_2^2 x_1^3, x_2^2 x_1^2 x_0, x_2^2 x_1 x_0^2, x_2^3 x_1^2, x_2^3 x_1 x_0, x_2^3 x_0^2, x_2^4 x_1, x_2^4 x_0, x_2^5\right\},
\]
$x^{\widetilde{\alpha}} = x_2^2$ and $x^{\widetilde{\beta}} = x_3 x_1$, so that 
\[
\begin{split}
 \widetilde{I} &{} = \Big\langle \{I_5\} \setminus\left\{x_2^2 x_1^3, x_2^2 x_1^2 x_0, x_2^2 x_1 x_0^2\right\}\\
 &\qquad{} \cup\left\{x_2^2 x_1^3 - x_3 x_1^4, x_2^2 x_1^2 x_0 - x_3 x_1^3 x_0, x_2^2 x_1 x_0^2 - x_3 x_1^2 x_0^2\right\} \Big\rangle
 \end{split}
\]
and $\IN_{\DegLex}(\widetilde{I}) = \GIN_{\DegLex}\big(\IN_{\DegLex}(\widetilde{I})\big) = (x_3^2, x_3 x_2, x_3 x_1^2,x_2^3)_{\geqslant 5}$.

Applying Algorithm \ref{alg:deformation} on $I$ with the $\mathtt{DegLex}$ term ordering, we have
\[
 \begin{split}
\min_{\DegLex} \{I_5\} = x_2^2 x_1 x_0^2 &{} >_{\DegLex} x_2^2 x_0^3 = \max_{\DegLex} \{I_5\}^{\mathcal{C}},\\
 \min_{\DegLex} \restrict{\{I_5\}}{1} = x_2^2 x_1^3 &{} <_{\DegLex} x_3 x_1^4 = \max_{\DegLex} \restrict{\{I_5\}^{\mathcal{C}}}{1}. 
 \end{split}
\]
Then the algorithm determines $\mathcal{F}_{x_2^2 x_1^3} = \{\textrm{id},\down{1},2\down{1}\}$ that is Borel-admissible w.r.t. $x_3 x_1^4$, and swapping $\mathcal{F}_{x_2^2 x_1^3}(x_2^2 x_1^3) = \{x_2^2 x_1^3,x_2^2 x_1^2 x_0,x_2^2 x_1 x_0^2\}$ with $\mathcal{F}_{x_2^2 x_1^3}(x_3 x_1^4) = \{x_3 x_1^4,$ $x_3 x_1^3 x_0,x_3 x_1^2 x_0^2\}$ we obtain again the ideal
\[
 J = (x_3^2, x_3 x_2, x_3 x_1^2,x_2^3)_{\geqslant 5} = \IN_{\DegLex}(\widetilde{I}).
\]
\end{example}

Our next goal is to prove that also with our method it is possible to reach from any point defined by a Borel-fixed ideal on the Hilbert scheme $\Hilb{n}{p(t)}$ some special point, by a sequence of consecutive Borel rational degenerations. Of course the lexicographic ideal defines a special point, but since we build the deformations working exclusively on the Borel sets defined by the piece of degree equal to $r$ of the ideals ($r$ is always the Gotzmann number of $p(t)$), we can generalize the property to the class of hilb-segment ideals\index{segment ideal!hilb-segment ideal} (Definition \ref{def:segments}).

\bigskip

Chosen a Hilbert scheme $\Hilb{n}{p(t)}$ and a term ordering $\sigma$, we want to have an overview on all $\sigma$-Borel rational degenerations among Borel-fixed ideals defining $\K$-rational points of $\Hilb{n}{p(t)}$.

\begin{definition}\index{degeneration graph|(}
 Let $\Hilb{n}{p(t)}$ be the Hilbert scheme parametrizing subschemes of the projective space $\PP^n$ with Hilbert polynomial $p(t)$ and let $\sigma$ be any term ordering, refinement of the Borel partial order $\leq_B$. We define the \emph{$\sigma$-degeneration graph} of $\Hilb{n}{p(t)}$ as the graph whose vertices correspond to all (saturated) Borel-fixed ideals defining points of $\Hilb{n}{p(t)}$ and whose edges represent the $\sigma$-Borel rational deformations, i.e. any edge goes from a Borel-fixed ideal to its $\sigma$-Borel degeneration.
The algorithm computing the $\sigma$-degeneration graph is described in Algorithm \ref{alg:deformationGraph}.
\end{definition}

\begin{algorithm}
\caption[Algorithm\hfill computing\hfill the\hfill $\sigma$-degeneration\hfill graph\hfill associated\hfill to\hfill a\newline Hilbert scheme.]{\label{alg:deformationGraph} How to compute the $\sigma$-degeneration graph associated to a Hilbert scheme. For details see Appendix \ref{ch:HSCpackage}.}
\begin{algorithmic}[1]
\STATE $\textsc{DegenerationGraph}(\Hilb{n}{p(t)},\sigma)$
\REQUIRE $\Hilb{n}{p(t)}$, Hilbert scheme.
\REQUIRE $\sigma$, term ordering, refinement of the Borel partial order $\leq_B$.
\ENSURE the $\sigma$-degeneration graph $(\textsf{vertices},\textsf{edges})$ of $\Hilb{n}{p(t)}$.
\STATE $\textsf{vertices} \leftarrow \textsc{BorelGeneratorDFS}\big(\K[x_0,\ldots,x_n],p(t)\big)$;
\STATE $\textsf{edges} \leftarrow \emptyset$;
\FORALL{$I \in \textsf{vertices}$}
\STATE $J \leftarrow \textsc{OrientedBorelRationalDeformation}(I,\sigma)$;
\IF{$J \neq I$}
\STATE $\textsf{edges} \leftarrow \textsf{edges} \cup \{(I,J)\}$;
\ENDIF
\ENDFOR
\RETURN $(\textsf{vertices},\textsf{edges})$;
\end{algorithmic}
\end{algorithm}

\begin{theorem}\label{prop:connectedTree}
Let $\Hilb{n}{p(t)}$ be a Hilbert scheme and let $\sigma$ be a term ordering. If there exists a Borel-fixed ideal $I$ defining a point of $\Hilb{n}{p(t)}$ which is a hilb-segment ideal w.r.t. $\sigma$, then the $\sigma$-degeneration graph of $\Hilb{n}{p(t)}$ is a rooted tree, with the ideal $I$ as root.
\end{theorem}

We recall briefly what we mean by rooted tree. A \emph{tree}\index{tree} is a connected graph $(V,E)$, such that $\vert E \vert = \vert V \vert -1$. A \emph{rooted tree}\index{tree!rooted} is a tree in which a fixed vertex (the root) determines a natural orientation of the edges, \lq\lq toward to\rq\rq\ and \lq\lq away from\rq\rq\ the root.

\begin{proof}
By definition $I$ is a $\sigma$-endpoint, because $I$ could not have $\sigma$-Borel degeneration, so $I$ is the natural root of the graph. To prove that the $\sigma$-degeneration graph is a rooted tree, it is sufficient to show that any other Borel ideal $J \neq I$ has a $\sigma$-Borel degeneration.

Let $r$ be the Gotzmann number of $p(t)$. For any ideal $J \neq I$, there exists a pair of monomials $(x^\alpha,x^\beta)$ such that $\{J_r\} \ni x^\alpha <_{\sigma} x^\beta \in \{J_r\}^{\mathcal{C}}$. So let $0\leqslant k < n$ be the integer such that
\[
\min_\sigma \restrict{\{J_r\}}{i} >_\sigma \max_\sigma \restrict{\{J_r\}^{\mathcal{C}}}{i}, \qquad \forall\ i=0,\ldots,k-1
\]
and
\[
\min_\sigma \restrict{\{J_r\}}{k} = x^\alpha <_\sigma x^\beta = \max_\sigma \restrict{\{J_r\}^{\mathcal{C}}}{k}.
\]
Let $\mathcal{F}_\alpha$ be the decreasing set of $x^\alpha$. If $\mathcal{F}_\alpha$ is Borel-admissible w.r.t. $x^\beta$, we finish because Algorithm \ref{alg:deformation} applied on $J$ and $\sigma$ does not return ideal $J$ itself. Then let us assume that $\mathcal{F}_\alpha$ is not Borel-admissible w.r.t. $x^\beta$: we want to show that there must exist a $\sigma$-Borel degeneration of $J$ determined by a pair of monomials $x^{\overline{\alpha}}$ and $x^{\overline{\beta}}$ such that $\min x^{\overline{\alpha}} = \min x^{\overline{\beta}} < k$, i.e. a $\sigma$-Borel degeneration of $J$ that Algorithm \ref{alg:deformation} should find before examining the restriction $\restrict{\{J_r\}}{k}$.

The first reason for which $\mathcal{F}_{\alpha}$ could be not Borel-admissible w.r.t. $x^\beta$ is the existence of some decreasing move in $\mathcal{F}_\alpha$ not admissible w.r.t. $x^\beta$. Let $\mathcal{G} \subset \mathcal{F}_\alpha$ be the set of the decreasing moves admissible on both monomials. Since $x^\alpha \in J$ could not belong to the hilb-segment ideal $I$, also every monomial in $\mathcal{F}_\alpha(x^\alpha)$ does not belong to $I$. To go back to $I$, the monomials in $\mathcal{G}(x^\alpha)$ could be replaced by the monomials in $\mathcal{G}(x^\beta)$ and the others in $\mathcal{F}_\alpha(x^\alpha)\setminus\mathcal{G}(x^\alpha)$ should to be replaced by monomials not obtained by decreasing moves from monomials in $\{J_r\}_k$, that is $\max_{\sigma} \restrict{\{J_r\}^{\mathcal{C}}}{i} >_\sigma \min_{\sigma} \restrict{\{J_r\}}{i}$, for some $i < s$.

The second reason for which $\mathcal{F}_{\alpha}$ could be not Borel-admissible is the existence of a move $F \in \mathcal{F}_{\alpha}$ such that the monomial $F(x^\beta)$ is not a $k$-maximal element, namely there exists $\up{j},\ j \geqslant k$, such that $\up{j} \big(F(x^\beta)\big) \notin \{J_r\}$. Let $i = \min F(x^\alpha) = \min F(x^\beta) < k$.
$\up{j} F(x^\beta)$ can not be obtained by a monomial in $\restrict{\{J_r\}^{\mathcal{C}}}{k}$ applying a composition of decreasing elementary moves $G$, because $G(x^\delta) = \up{j} F(x^\beta)$ implies $x^\delta >_B x^\beta$ in contradiction with the hypothesis $x^\beta = \max_{\sigma} \restrict{\{J_r\}}{k}$. Since $x^\alpha <_\sigma x^\beta \Rightarrow F(x^\alpha) <_\sigma F(x^\beta)$,
\[
\max_{\sigma} \restrict{\{J_r\}}{i} \geq_{\sigma} \up{j}\big(F(x^\beta)\big) >_\sigma F(x^\beta) >_{\sigma} F(x^\alpha) \geq_{\sigma} \min_{\sigma} \restrict{\{J_r\}}{k}.
\]

In both cases, there exist $i < k$ and a pair of monomials $x^{\overline{\alpha}},x^{\overline{\beta}}$
\[
  x^{\overline{\beta}} = \max_{\sigma} \restrict{\{J_r\}}{i} >_\sigma \min_{\sigma} \restrict{\{J_r\}}{i} = x^{\overline{\alpha}}.
\]
We compute again the decreasing set $\mathcal{F}_{\overline{\alpha}}$ and we check if it is Borel-admissible w.r.t. $x^{\overline{\beta}}$: if not we repeat the reasoning and we look for monomials involving more variables. Finally, we are sure to find a Borel-admissible set of decreasing moves because if we reach the smallest variable $x_0$, the decreasing set only contains the identity move. Note that it is not possible to have cycles thanks total order on the monomials, because each $\sigma$-Borel degeneration approaches to the hilb-segment ideal.
\end{proof}

\begin{corollary}\index{Hilbert scheme!connectedness of the}
 The Hilbert scheme $\Hilb{n}{p(t)}$ is connected.
 \end{corollary}
\begin{proof}
Let $I$ be any ideal defining a point on $\Hilb{n}{p(t)}$. As usual through an affine Gr\"obner degeneration the point defined by $I$ can be connected to the point defined by the Borel-fixed ideal $\GIN(I)$.

Of course on $\Hilb{n}{p(t)}$, there is the lexicographic point corresponding to the lexicographic ideal\index{lexicographic ideal} (Proposition \ref{prop:satLexIdeal}), hilb-segment ideal w.r.t. $\DegLex$. By Theorem \ref{prop:connectedTree}, the $\DegLex$-degeneration graph is a \emph{connected} rooted tree, so the point defined by any ideal $I$ can be connected to the lexicographic point by an initial affine degeneration and a sequence of $\DegLex$-Borel rational degenerations.
\end{proof}

We underline that if $\Hilb{n}{p(t)}$ does not contain a point defined by a hilb-segment ideal w.r.t. $\sigma$, the $\sigma$-degeneration graph could be not connected, as the following example shows.

\begin{example}\label{ex:defGraphs}
 Let us consider the Hilbert scheme \gls{Hilb6tm5P3}. In $\K[x_0,x_1,x_2,x_3]$ there are 11 saturated Borel-fixed ideals with Hilbert polynomial $6t-5$ (whose Gotzmann number is 10) and many of them are hilb-segment ideals. In the following list of the ideals, we specify the term ordering (computed with Algorithm \ref{alg:HilbRegSegment}) for which the corresponding ideal becomes (possibly) a hilb-segment ideal: 
\[
\begin{array}{l c l}
 I_{1} = (x_3,x_2^7,x_2^6x_1^4), & & \DegLex;\\
 I_{2} = (x_3,x_2^8,x_2^7x_1,x_2^6x_1^3), & & \omega_2 = (37,6,3,1);\\
 I_{3} = (x_3^2,x_3x_2,x_3x_1,x_2^7,x_2^6x_1^3), && \omega_3 = (21,4,2,1);\\
 I_{4} = (x_3^2,x_3x_2,x_3x_1,x_2^8,x_2^7x_1,x_2^6x_1^2); && \\
 I_{5} = (x_3^2,x_3x_2,x_3x_1^2,x_2^7,x_2^6x_1^2), && \omega_5 = (25,5,2,1);\\
 I_6 =(x_3^2,x_3x_2,x_3x_1^3,x_2^7,x_2^6x_1), && \omega_6 = (29,6,2,1);\\
 I_{7} = (x_3^2,x_3x_2^2,x_3x_2x_1,x_3x_1^2,x_2^7,x_2^6x_1); && \\
 I_{8} = (x_3^2,x_3x_2,x_3x_1^4,x_2^6), && \omega_8 = (33,7,2,1); \\
 I_{9} = (x_3^2,x_3x_2^2,x_3x_2x_1,x_3x_1^3,x_2^6); && \\
 I_{10} = (x_3^3,x_3^2x_2,x_3x_2^2,x_3^2x_1,x_3x_2x_1,x_3x_1^2,x_2^6); && \\
 I_{11} = (x_3^2,x_3x_2,x_2^5), && \omega_{11} = (33,11,2,1).
\end{array}
\]
The degeneration graphs of $\Hilb{n}{p(t)}$ w.r.t. all the term ordering listed above turn out to be rooted trees as shown in Figures \ref{fig:BorelGraph_1}--\ref{fig:BorelGraph_11}. 

\input{Figures/Ch3/HilbConnected}

The\hfill ideal\hfill generates\hfill by\hfill the\hfill greatest\hfill $\binom{3+10}{3} - p(10) = 231$\hfill monomials\hfill in\\ $\K[x_0,x_1,x_2,x_3]_{10}$ w.r.t. $\RevLex$ has constant Hilbert polynomial equal to $55$, so $\Hilb{3}{6t-5}$ does not contain a point defined by a hilb-segment ideal w.r.t. $\RevLex$. Applying Algorithm \ref{alg:deformationGraph} on $\Hilb{3}{6t-5}$ and $\RevLex$, we find that both $I_{10}$ and $I_{11}$ are $\RevLex$-endpoint, so that the $\RevLex$-degeneration graph is not connected (Figure \ref{fig:BorelGraph_RevLex}).
\end{example}

\subsection{The special case of constant Hilbert polynomials}

Let us finally consider the special case of the \DegLex-degeneration graph of Hilbert schemes of points. Given a Borel-fixed ideal $I$ defining a point of $\Hilb{n}{s}$, we highlight that among the minimal monomials of $\{I_s\}$ there will be surely a monomial of the type $x_1^{a} x_0^{s-a}$. In fact $x_1^s$ belongs to the ideal and applying the decreasing move $\down{1}$ repeatedly, we will find a monomial $x_1^{a} x_0^{s-a}$ such that $\down{1}(x_1^a x_0^{s-a}) \notin \{I_r\}$. It is easy to deduce that the power $a$ of the variable $x_1$ is equal to the regularity of the saturated ideal $I^{\sat}$, because $x_1^a$ is one of the generators of the saturation (Proposition \ref{prop:saturationBorelIdeal}) and it is the one of highest degree (Proposition \ref{prop:regularityDegreeGenerators}). 

\begin{proposition}\label{prop:regSteps}\index{regularity of an ideal}
Let $I \subset \K[x]$ be a Borel-fixed ideal defining a point of $\Hilb{n}{s}$ and let $\reg(I)$ be the regularity of its saturation $I^\sat$. $s-\reg(I)$ consecutive $\DegLex$-Borel rational degenerations lead from the point defined by $I$ to the lexicographic point.
\begin{proof}
For the lexicographic ideal $L$ associated to the constant Hilbert polynomial $p(t)=s$, the order set contains the monomials
\[
\{L_s\}^{\mathcal{C}} = \left\{ x_0^s,x_1 x_0^{s-1},\ldots,x_1^{s-1}x_0\right\}.
\]
For the Borel-fixed ideal $I$ defining a point of $\Hilb{n}{s}$, we can divide the order set as follows:
\[
 \{I_s\}^{\mathcal{C}} = \left\{x_0^s,\ldots,x_1^{\reg(I)-1} x_0^{s-\reg(I)+1}\right\} \cup \left\{x^{\gamma_1},\ldots,x^{\gamma_{s-\reg(I)}}\right\},
\]
where $\max x^{\gamma_i} > 1,\ \forall\ i = 1,\ldots,s-\reg(I)$.

The $\DegLex$-Borel degeneration of $I$, that we can obtain applying Algorithm \ref{alg:deformation}, is determined by the monomials
\[
\min_{\DegLex} \{I_s\} = x_1^{\reg(I)} x_0^{s-\reg(I)}\quad\text{and}\quad \max_{\DegLex} \{I_s\}^{\mathcal{C}} = \max_{\DegLex} \left\{x^{\gamma_1},\ldots,x^{\gamma_{s-\reg(I)}}\right\} = x^\beta,
\]
so that the $\DegLex$-Borel degeneration of $I$ is the ideal
\[
J = \left\langle \{I_s\}\setminus\left\{x_1^{\reg(I)} x_0^{s-\reg(I)}\right\} \cup \{x^\beta\}\right\rangle^{\sat}
\]
whose order set in degree $s$ is
\[
\{J_s\}^{\mathcal{C}} = \left\{ x_0^s,\ldots,x_1^{\reg(I)-1} x_0^{s-\reg(I)+1},x_1^{\reg(I)} x_0^{s-\reg(I)}\right\} \cup \left(\left\{x^{\gamma_1},\ldots,x^{\gamma_{s-\reg(I)}}\right\}\setminus\{x^\beta\}\right)
\]
and whose regularity is $\reg(J) = \reg(I)+1$. Repeating $s-\reg(I)$ times this process, we will obtain an ideal with regularity equal to $\reg(I) + \big(s-\reg(I)\big)= s$, i.e. the lexicographic ideal $L$.
\end{proof}
\end{proposition}

\begin{figure}[!ht]
\begin{center}
\captionsetup[subfloat]{justification=centering}
\subfloat[][$I = (x_2^3,x_2^2 x_1,x_2x_1^2,x_1^3)$, $\reg(I) = 3$.]{
\begin{tikzpicture}[scale=0.6]
\tikzstyle{ideal}=[circle,draw=black,fill=black,inner sep=1.5pt]
\tikzstyle{quotient}=[circle,draw=black,fill=white,thick,inner sep=1.5pt]


\filldraw [rounded corners = 5pt,fill=red!25,draw=white] (3.5,-3.5) -- (5.5,-3.5) -- (5.5,-5.5) -- (4.5,-5.5) -- (4.5,-4.5) -- (3.5,-4.5) -- cycle;

\filldraw [rounded corners = 5pt,fill=blue!25,draw=white] (5.5,-3.5) -- (6.5,-3.5) -- (6.5,-6.5) -- (5.5,-6.5) -- cycle;

\node at (-1.1,0) [] {};
\node at (7.1,0) [] {};

\node at (0,0) [ideal] {};
\node at (1,0) [ideal] {};
\node at (2,0) [ideal] {};
\node at (3,0) [ideal] {};
\node at (4,0) [ideal] {};
\node at (5,0) [ideal] {};
\node at (6,0) [ideal] {};

\node at (1,-1) [ideal] {};
\node at (2,-1) [ideal] {};
\node at (3,-1) [ideal] {};
\node at (4,-1) [ideal] {};
\node at (5,-1) [ideal] {};
\node at (6,-1) [ideal] {};

\node at (2,-2) [ideal] {};
\node at (3,-2) [ideal] {};
\node at (4,-2) [ideal] {};
\node at (5,-2) [ideal] {};
\node at (6,-2) [ideal] {};

\node at (3,-3) [ideal] {};
\node at (4,-3) [ideal] {};
\node at (5,-3) [ideal] {};
\node at (6,-3) [regular polygon,regular polygon sides=4,draw=black,fill=black,inner sep=1.5pt] {};

\node at (4,-4) [regular polygon,regular polygon sides=4,draw=black,fill=white,inner sep=1.5pt] {};
\node at (5,-4) [quotient] {};
\node at (6,-4) [quotient] {};

\node at (5,-5) [quotient] {};
\node at (6,-5) [quotient] {};

\node at (6,-6) [quotient] {};
\end{tikzpicture}
}
\qquad\qquad
\subfloat[][$\overline{I} = (x_2^2,x_2x_1^2,x_1^4)$, $\reg(\overline{I}) = 4$.]{
\begin{tikzpicture}[scale=0.6]
\tikzstyle{ideal}=[circle,draw=black,fill=black,inner sep=1.5pt]
\tikzstyle{quotient}=[circle,draw=black,fill=white,thick,inner sep=1.5pt]

\filldraw [rounded corners = 5pt,fill=red!25,draw=white] (4.5,-3.5) -- (5.5,-3.5) -- (5.5,-5.5) -- (4.5,-5.5) -- cycle;

\filldraw [rounded corners = 5pt,fill=blue!25,draw=white] (5.5,-2.5) -- (6.5,-2.5) -- (6.5,-6.5) -- (5.5,-6.5) -- cycle;

\node at (-1,0) [] {};
\node at (7,0) [] {};

\node at (0,0) [ideal] {};
\node at (1,0) [ideal] {};
\node at (2,0) [ideal] {};
\node at (3,0) [ideal] {};
\node at (4,0) [ideal] {};
\node at (5,0) [ideal] {};
\node at (6,0) [ideal] {};

\node at (1,-1) [ideal] {};
\node at (2,-1) [ideal] {};
\node at (3,-1) [ideal] {};
\node at (4,-1) [ideal] {};
\node at (5,-1) [ideal] {};
\node at (6,-1) [ideal] {};

\node at (2,-2) [ideal] {};
\node at (3,-2) [ideal] {};
\node at (4,-2) [ideal] {};
\node at (5,-2) [ideal] {};
\node at (6,-2) [regular polygon,regular polygon sides=4,draw=black,fill=black,inner sep=1.5pt] {};

\node at (3,-3) [ideal] {};
\node at (4,-3) [ideal] {};
\node at (5,-3) [ideal] {};
\node at (6,-3) [quotient] {};

\node at (4,-4) [ideal] {};
\node at (5,-4) [regular polygon,regular polygon sides=4,draw=black,fill=white,inner sep=1.5pt] {};
\node at (6,-4) [quotient] {};

\node at (5,-5) [quotient] {};
\node at (6,-5) [quotient] {};

\node at (6,-6) [quotient] {};
\end{tikzpicture}
}\\

\bigskip

\subfloat[][$\widetilde{I} = (x_2^2,x_2x_1,x_1^5)$, $\reg(\widetilde{I}) = 5$.]{
\begin{tikzpicture}[scale=0.6]
\tikzstyle{ideal}=[circle,draw=black,fill=black,inner sep=1.5pt]
\tikzstyle{quotient}=[circle,draw=black,fill=white,thick,inner sep=1.5pt]

\filldraw [rounded corners = 5pt,fill=red!25,draw=white] (4.5,-4.5) -- (5.5,-4.5) -- (5.5,-5.5) -- (4.5,-5.5) -- cycle;
\filldraw [rounded corners = 5pt,fill=blue!25,draw=white] (5.5,-1.5) -- (6.5,-1.5) -- (6.5,-6.5) -- (5.5,-6.5) -- cycle;

\node at (-1,0) [] {};
\node at (7,0) [] {};

\node at (0,0) [ideal] {};
\node at (1,0) [ideal] {};
\node at (2,0) [ideal] {};
\node at (3,0) [ideal] {};
\node at (4,0) [ideal] {};
\node at (5,0) [ideal] {};
\node at (6,0) [ideal] {};

\node at (1,-1) [ideal] {};
\node at (2,-1) [ideal] {};
\node at (3,-1) [ideal] {};
\node at (4,-1) [ideal] {};
\node at (5,-1) [ideal] {};
\node at (6,-1) [regular polygon,regular polygon sides=4,draw=black,fill=black,inner sep=1.5pt] {};

\node at (2,-2) [ideal] {};
\node at (3,-2) [ideal] {};
\node at (4,-2) [ideal] {};
\node at (5,-2) [ideal] {};
\node at (6,-2) [quotient] {};

\node at (3,-3) [ideal] {};
\node at (4,-3) [ideal] {};
\node at (5,-3) [ideal] {};
\node at (6,-3) [quotient] {};

\node at (4,-4) [ideal] {};
\node at (5,-4) [ideal] {};
\node at (6,-4) [quotient] {};

\node at (5,-5) [regular polygon,regular polygon sides=4,draw=black,fill=white,inner sep=1.5pt] {};
\node at (6,-5) [quotient] {};

\node at (6,-6) [quotient] {};
\end{tikzpicture}
}
\qquad\qquad
\subfloat[][$L = (x_2,x_1^6)$, $\reg(L) = 6$.]{
\begin{tikzpicture}[scale=0.6]
\tikzstyle{ideal}=[circle,draw=black,fill=black,inner sep=1.5pt]
\tikzstyle{quotient}=[circle,draw=black,fill=white,thick,inner sep=1.5pt]

\filldraw [rounded corners = 5pt,fill=blue!25,draw=white] (5.5,-0.5) -- (6.5,-0.5) -- (6.5,-6.5) -- (5.5,-6.5) -- cycle;

\node at (-1,0) [] {};
\node at (7,0) [] {};

\node at (0,0) [ideal] {};
\node at (1,0) [ideal] {};
\node at (2,0) [ideal] {};
\node at (3,0) [ideal] {};
\node at (4,0) [ideal] {};
\node at (5,0) [ideal] {};
\node at (6,0) [ideal] {};

\node at (1,-1) [ideal] {};
\node at (2,-1) [ideal] {};
\node at (3,-1) [ideal] {};
\node at (4,-1) [ideal] {};
\node at (5,-1) [ideal] {};
\node at (6,-1) [quotient] {};

\node at (2,-2) [ideal] {};
\node at (3,-2) [ideal] {};
\node at (4,-2) [ideal] {};
\node at (5,-2) [ideal] {};
\node at (6,-2) [quotient] {};

\node at (3,-3) [ideal] {};
\node at (4,-3) [ideal] {};
\node at (5,-3) [ideal] {};
\node at (6,-3) [quotient] {};

\node at (4,-4) [ideal] {};
\node at (5,-4) [ideal] {};
\node at (6,-4) [quotient] {};

\node at (5,-5) [ideal] {};
\node at (6,-5) [quotient] {};

\node at (6,-6) [quotient] {};
\end{tikzpicture}
}
\end{center}
\caption{\label{fig:defToLex} An example of the sequence of \texttt{DegLex}-Borel degeneration leading from the point defined by a generic Borel-fixed ideal to the lexicographic point.}
\end{figure}
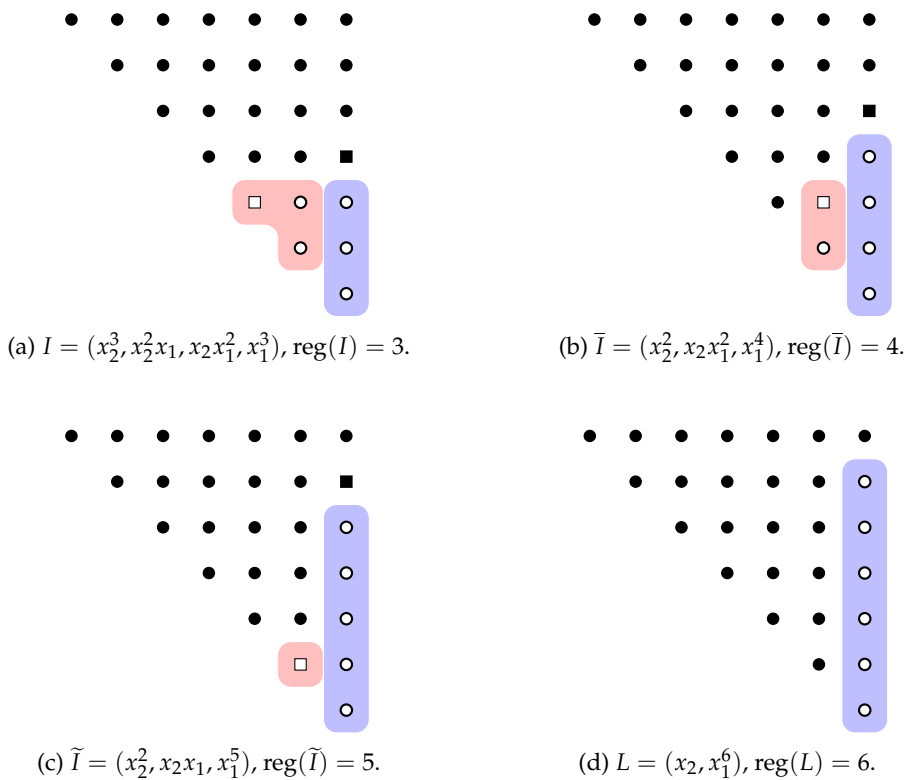

\begin{corollary}\label{cor:treeDepth}
 Let $\Hilb{n}{s}$ be a zero-dimensional Hilbert scheme. The \textnormal{\DegLex}-degeneration graph of $\Hilb{n}{s}$ has height equal to $s-a$, where $a$ is the smallest positive integer such that
\[
 s \leqslant \sum_{i=0}^{a-1} \binom{n-1 + i}{n-1}.
\]
\end{corollary}
We recall that the height of a rooted tree is the maximal distance between a vertex and the root, where the vertices connected to the root have distance 1, the vertices connected to vertices of distance 1 have distance 2 and so on.
\begin{proof}
By Proposition \ref{prop:regSteps}, to determine the height of the $\DegLex$-degeneration graph of $\Hilb{n}{s}$, we have to understand which is the lowest regularity of a Borel-fixed ideal with constant Hilbert polynomial $p(t) = s$. We saw that the regularity of such an ideal $I$ coincides with the minimal power of the variable $x_1$ in a monomial of the type $x_1^a x_0^{s-a}$ belonging  to $\{I_s\}$. So the question is how many monomials we can put in a order set not containing $x_1^a x_0^{s-a}$ and the answer is that we can put all the monomials with a power of the variable $x_0$ greater than $s-a$, i.e.
\[
\underbrace{\binom{n-1+a-1}{n-1}}_{\K[x_1,\ldots,x_n]_{a-1}\cdot x_0^{s-a+1}} + \ldots + \underbrace{\binom{n-1+1}{n-1}}_{\K[x_1,\ldots,x_n]_1 \cdot x_0^{s-1}} + \underbrace{\binom{n-1+0}{n-1}}_{x_0^s} = \sum_{i=0}^{a-1}\binom{n-1+i}{n-1}.
\]
Finally, to obtain a Borel-fixed ideal $I$ having $x_1^a$ as minimal generator, namely $\reg(I) = a$, we need $s \leqslant \sum_{i=0}^{a-1}\binom{n-1+i}{n-1}$.

Note that the choosing monomials to put in the order set w.r.t. a decreasing order on the power of the smallest variable $x_0$ agrees with choosing the smallest monomial w.r.t. $\RevLex$ among those still belonging to the Borel set, hence the lowest regularity $a$ of a Borel-fixed ideal with Hilbert polynomial $p(t) = s$ is always realized by the hilb-segment ideal w.r.t. $\RevLex$.
\end{proof}

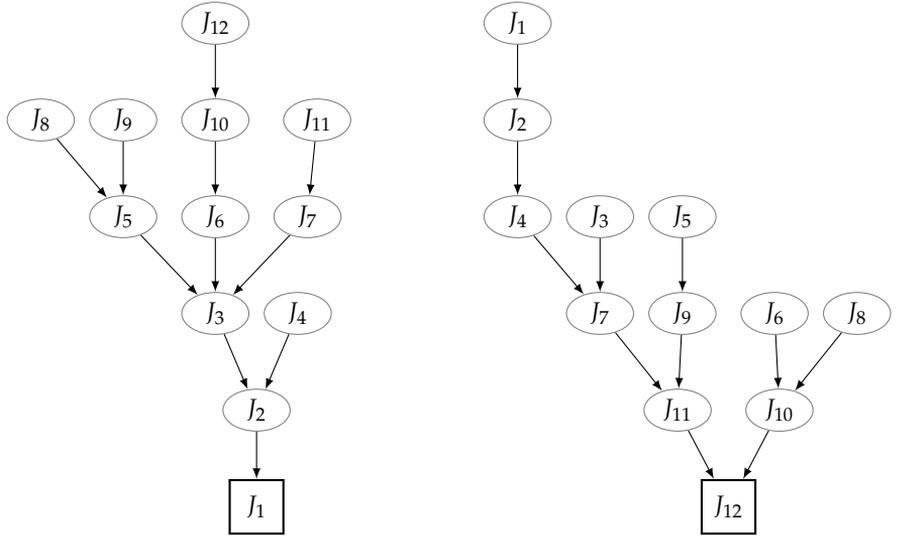
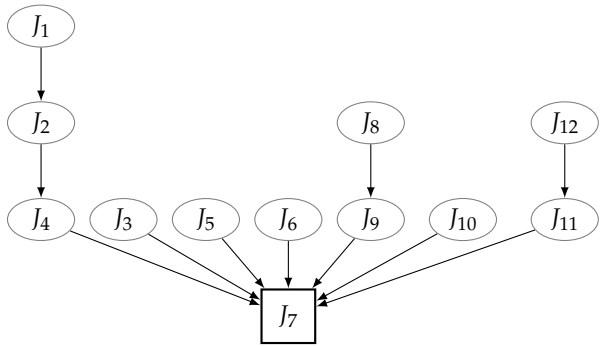
\begin{figure}[!ht]
\begin{center}
\subfloat[][The $\DegLex$-degeneration graph of $\Hilb{3}{8}$.]{\label{fig:heightGraphs_a}
\begin{tikzpicture}[>=latex,scale=0.9]
\tikzstyle{place}=[ellipse,draw=black!50,minimum width=25pt,inner sep=1.5pt]
\tikzstyle{place1}=[rectangle,draw=black,thick,minimum width=20pt,minimum height=20pt]

\node (11) at (122bp,167bp) [place] {\footnotesize $J_{11}$};
  \node (10) at (80bp,167bp) [place] {\footnotesize $J_{10}$};
  \node (12) at (80bp,207bp) [place] {\footnotesize $J_{12}$};
  \node (1) at (97bp,7bp) [place1] {\footnotesize $J_1$};
  \node (3) at (80bp,87bp) [place] {\footnotesize $J_3$};
  \node (2) at (97bp,47bp) [place] {\footnotesize $J_2$};
  \node (5) at (42bp,127bp) [place] {\footnotesize $J_5$};
  \node (4) at (114bp,87bp) [place] {\footnotesize $J_4$};
  \node (7) at (118bp,127bp) [place] {\footnotesize $J_7$};
  \node (6) at (80bp,127bp) [place] {\footnotesize $J_6$};
  \node (9) at (42bp,167bp) [place] {\footnotesize $J_9$};
  \node (8) at (8bp,167bp) [place] {\footnotesize $J_8$};
  \draw [->] (11) -- (7);
  \draw [->] (10) -- (6);
  \draw [->] (2) -- (1);
  \draw [->] (12) -- (10);
  \draw [->] (7) -- (3);
  \draw [->] (3) -- (2);
  \draw [->] (8) -- (5);
  \draw [->] (5) -- (3);
  \draw [->] (9) -- (5);
  \draw [->] (6) -- (3);
  \draw [->] (4) -- (2);
\end{tikzpicture}
}
\qquad\qquad
\subfloat[][The $\RevLex$-degeneration graph of $\Hilb{3}{8}$.]{\label{fig:heightGraphs_b}
\begin{tikzpicture}[>=latex,scale=0.9]
\tikzstyle{place}=[ellipse,draw=black!50,minimum width=25pt,inner sep=1.5pt]
\tikzstyle{place1}=[rectangle,draw=black,thick,minimum width=20pt,minimum height=20pt]

\node (11) at (74bp,47bp) [place] {\footnotesize $J_{11}$};
  \node (10) at (116bp,47bp) [place] {\footnotesize $J_{10}$};
  \node (12) at (95bp,7bp) [place1] {\footnotesize $J_{12}$};
  \node (1) at (8bp,207bp) [place] {\footnotesize $J_1$};
  \node (3) at (42bp,127bp) [place] {\footnotesize $J_3$};
  \node (2) at (8bp,167bp) [place] {\footnotesize $J_2$};
  \node (5) at (76bp,127bp) [place] {\footnotesize $J_5$};
  \node (4) at (8bp,127bp) [place] {\footnotesize $J_4$};
  \node (7) at (42bp,87bp) [place] {\footnotesize $J_7$};
  \node (6) at (114bp,87bp) [place] {\footnotesize $J_6$};
  \node (9) at (76bp,87bp) [place] {\footnotesize $J_9$};
  \node (8) at (148bp,87bp) [place] {\footnotesize $J_8$};
  \draw [->] (3) -- (7);
  \draw [->] (5) -- (9);
  \draw [->] (4) -- (7);
  \draw [->] (10) -- (12);
  \draw [->] (8) -- (10);
  \draw [->] (6) -- (10);
  \draw [->] (1) -- (2);
  \draw [->] (9) -- (11);
  \draw [->] (11) -- (12);
  \draw [->] (7) -- (11);
  \draw [->] (2) -- (4);
\end{tikzpicture}
}\\

\bigskip

\subfloat[][The $\omega_7$-degeneration graph of $\Hilb{3}{8}$.]{\label{fig:heightGraphs_c}
\begin{tikzpicture}[>=latex,scale=0.9]
\tikzstyle{place}=[ellipse,draw=black!50,minimum width=25pt,inner sep=1.5pt]
\tikzstyle{place1}=[rectangle,draw=black,thick,minimum width=20pt,minimum height=20pt]

\node (11) at (224bp,47bp) [place] {\footnotesize $J_{11}$};
  \node (10) at (182bp,47bp) [place] {\footnotesize $J_{10}$};
  \node (12) at (224bp,87bp) [place] {\footnotesize $J_{12}$};
  \node (1) at (8bp,127bp) [place] {\footnotesize $J_1$};
  \node (3) at (42bp,47bp) [place] {\footnotesize $J_3$};
  \node (2) at (8bp,87bp) [place] {\footnotesize $J_2$};
  \node (5) at (76bp,47bp) [place] {\footnotesize $J_5$};
  \node (4) at (8bp,47bp) [place] {\footnotesize $J_4$};
  \node (7) at (110bp,7bp) [place1] {\footnotesize $J_7$};
  \node (6) at (110bp,47bp) [place] {\footnotesize $J_6$};
  \node (9) at (144bp,47bp) [place] {\footnotesize $J_9$};
  \node (8) at (144bp,87bp) [place] {\footnotesize $J_8$};
  \draw [->] (11) -- (7);
  \draw [->] (3) -- (7);
  \draw [->] (10) -- (7);
  \draw [->] (5) -- (7);
  \draw [->] (4) -- (7);
  \draw [->] (6) -- (7);
  \draw [->] (12) -- (11);
  \draw [->] (9) -- (7);
  \draw [->] (1) -- (2);
  \draw [->] (8) -- (9);
  \draw [->] (2) -- (4);
\end{tikzpicture}
}
\end{center}
\caption{\label{fig:heightGraphs} The graphical representation of some degeneration graphs discussed in Example \ref{ex:depthDefGraph}. The graphs are drawn as trees to highlight their height.}
\end{figure}

\begin{example}\label{ex:depthDefGraph}
 Let us consider the Hilbert scheme \gls{Hilb8P3}. There are 12 Borel-fixed ideals, that we list again with the term ordering for which they (possibly) are hilb-segment:
\[
\begin{array}{l c l}
J_1 = (x_3,x_2,x_1^8), && \DegLex;\\
J_2 = (x_3,x_2^2,x_2x_1,x_1^7), && \omega_2 = (8,7,2,1);\\
J_3 = (x_3,x_2^2,x_2x_1^2,x_1^6), && \omega_3 = (7,5,2,1);\\
J_4 = (x_3^2,x_3x_2,x_2^2,x_3x_1,x_2x_1,x_1^6), && \omega_4 = (11,10,3,1);\\
J_5 = (x_3,x_2^2,x_2x_1^3,x_1^5), && \omega_5 = (11,6,3,1);\\
J_6 = (x_3,x_2^3,x_2^2x_1,x_2x_1^2,x_1^5); && \\
J_7 = (x_3^2,x_3x_2,x_2^2,x_3x_1,x_2x_1^2,x_1^5), && \omega_7 = (5,4,2,1);\\
J_8 = (x_3,x_2^3,x_2^2x_1,x_2x_1^3,x_1^4), && \omega_8 = (9,4,3,1);\\
J_9 = (x_3^2,x_3x_2,x_2^2,x_3x_1,x_2x_1^3,x_1^4); && \\
J_{10} = (x_3^2,x_3x_2,x_3x_1,x_2^3,x_2^2x_1,x_2x_1^2,x_1^4), && \omega_{10} = (6,4,3,1);\\
J_{11} = (x_3^2,x_3x_2,x_2^2,x_3x_1^2,x_2x_1^2,x_1^4), && \omega_{11} = (6,5,3,1);\\
J_{12} = (x_3^2,x_3x_2,x_2^3,x_2^2x_1,x_3x_1^2,x_2x_1^2,x_1^3), && \RevLex.
\end{array}
\]
The $\DegLex$-degeneration graph of $\Hilb{3}{8}$ (Figure \ref{fig:heightGraphs_a}) has height $5$, because 
\[
8 \leqslant \sum_{i=0}^{3-1} \binom{2+i}{2} = \binom{2}{2} + \binom{3}{2} + \binom{4}{2} = 10 \quad \text{and} \quad 8 > \sum_{i=0}^{2-1} \binom{2+i}{2} = 4,
\]
and the $\RevLex$-degeneration graph (Figure \ref{fig:heightGraphs_b}) has the same height, being the hilb-segment ideal $J_{12}$ w.r.t. $\RevLex$ the most distant vertix from the lexicographic ideal $J_1$.

Degeneration graphs with a lower height can be found considering hilb-segment ideals with an intermediate regularity (Figure \ref{fig:heightGraphs_c}).
\end{example}

\newpage

\vspace*{\stretch{1}}

\index{degeneration graph|)}

\vspace*{\stretch{2}}

\section[Borel-fixed ideals on a same component of the Hilbert scheme]{Borel-fixed ideals defining points lying on a same\\ component of the Hilbert scheme}

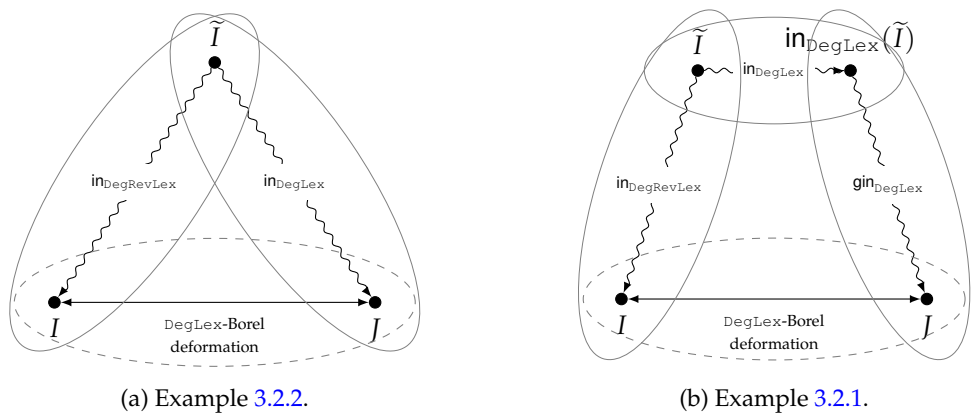
\begin{figure}[H] 
\begin{center}
\subfloat[][Example \ref{ex:PSsameRational}.]{\label{fig:ComponentGraph_a}
\begin{tikzpicture}[>=latex,scale=1.05]
\tikzstyle{place}=[circle,draw,fill,inner sep=0pt,minimum size=1.5mm]
\tikzstyle{place1}=[shape=ellipse,draw,minimum width=5.2cm,minimum height=1.8cm]
\node (I) at (0,0) [place,label=270:{\small $I$}] {};
\node (I2) at (2,3) [place,label=90:{\small $\widetilde{I}$}] {};
\node (J) at (4,0) [place,label=270:{\small $J$}] {};

\draw [->,thin,decorate,decoration={snake,amplitude=.3mm,segment length=2mm}] (I2) --node[fill=white]{\tiny $\IN_{\RevLex}$} (I);
\draw [->,thin,decorate,decoration={snake,amplitude=.3mm,segment length=2mm}] (I2) --node[fill=white]{\tiny $\IN_{\DegLex}$} (J);
\draw [<->,thin] (I) --node[below]{\tiny $\begin{array}{c} \text{\DegLex-Borel}\\ \text{deformation}\end{array}$} (J);

\node at (3,1.5) [place1,rotate=124,very thin,black!50] {};
\node at (1,1.5) [place1,rotate=56,very thin,black!50] {};
\draw [thin,dashed,black!50] (2,0) ellipse (2.5 and 0.8);
\end{tikzpicture}
}
\qquad\qquad
\subfloat[][Example \ref{ex:PSbreakingBorel}.]{\label{fig:ComponentGraph_b}
  \begin{tikzpicture}[>=latex]
\tikzstyle{place}=[circle,draw,fill,inner sep=0pt,minimum size=1.5mm]
\tikzstyle{place1}=[shape=ellipse,draw,minimum width=4.8cm,minimum height=1.6cm]

\node (I) at (0,0) [place,label=270:{\small $I$}] {};
\node (I2) at (1,3) [place,label=90:{\small $\widetilde{I}$}] {};
\node (I3) at (3,3) [place,label=90:{\small $\IN_{\DegLex}(\widetilde{I})$}] {};
\node (J) at (4,0) [place,label=270:{\small $J$}] {};

\draw [<->,thin] (I) --node[below]{\tiny $\begin{array}{c} \text{\DegLex-Borel}\\ \text{deformation}\end{array}$} (J);
\draw [->,thin,decorate,decoration={snake,amplitude=.3mm,segment length=2mm}] (I2) --node[fill=white]{\tiny $\IN_{\RevLex}$} (I);
\draw [->,thin,decorate,decoration={snake,amplitude=.3mm,segment length=2mm}] (I2) --node[fill=white]{\tiny $\IN_{\DegLex}$} (I3);
\draw [->,thin,decorate,decoration={snake,amplitude=.3mm,segment length=2mm}] (I3) --node[fill=white]{\tiny $\GIN_{\DegLex}$} (J);

\draw [thin,dashed,black!50] (2,0) ellipse (2.5 and 0.8);
\draw [very thin,black!50] (2,3) ellipse (1.7 and 0.7);
\node at (3.5,1.5) [place1,rotate=108,very thin,black!50] {};
\node at (0.5,1.5) [place1,rotate=72,very thin,black!50] {};

\end{tikzpicture}
}
\caption[An outline of what we can deduced about components of the Hilbert scheme from the deformations used to prove the connectedness.]{\label{fig:ComponentGraph} A Borel rational deformation ensures that the two Borel ideals defining it lie on a same component (drawn with the dashed line). With Peeva-Stillman method, in the best case (Fig. \ref{fig:ComponentGraph_a}), we can deduce that the components containing two Borel-fixed ideals have non-empty intersection but in the general case (Fig. \ref{fig:ComponentGraph_b}), we could get no information.}
\end{center}
\end{figure}

The last section of this chapter is devoted to the study of components of the Hilbert scheme, keeping in mind a question posed by Reeves in \cite{Reeves} \lq\lq Is the subset of Borel-fixed ideals on a component enough to determine the component?\rq\rq. In fact the points corresponding to two ideals connected by any Borel rational deformation lie on the same component of the Hilbert scheme. We underline that the technique introduced  by Peeva and Stillman in \cite{PeevaStillman} is slightly different, because in order for passing from a Borel-fixed ideal to another one they use at least two affine deformations (see Figure \ref{fig:ComponentGraph}).

We introduce a new graph related to Borel-fixed ideals and Hilbert scheme components, different from the incidence graph defined by Reeves in \cite{Reeves}, but also useful to understand the intersections among components and we will call it Borel incidence graph.

\begin{definition}\index{Borel incidence graph}
 Given the Hilbert scheme $\Hilb{n}{p(t)}$, we define the \emph{Borel incidence graph of $\Hilb{n}{p(t)}$} as the graph whose vertices correspond to the (saturated) Borel-fixed ideals of $\K[x]$ with Hilbert polynomial $p(t)$ and whose edges represents 
Borel rational deformations.
\end{definition}

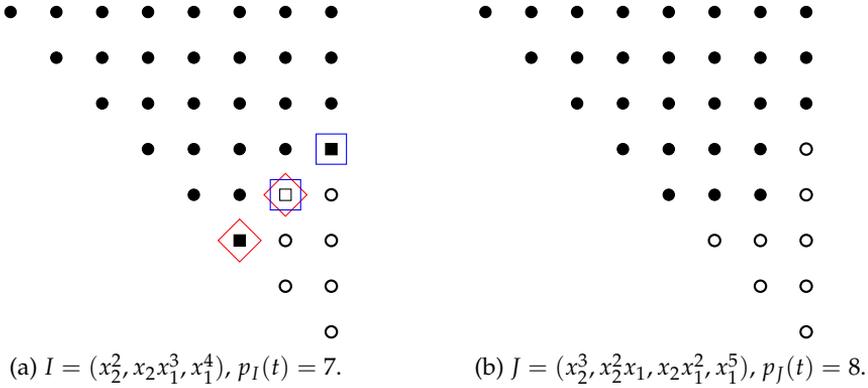
\begin{figure}[!ht]
\begin{center}
\subfloat[][$I = (x_2^2,x_2x_1^3,x_1^4)$, $p_I(t) = 7$.]{\label{fig:problem1_a}
\begin{tikzpicture}[scale=0.6]
\tikzstyle{ideal}=[circle,draw=black,fill=black,inner sep=1.5pt]
\tikzstyle{quotient}=[circle,draw=black,fill=white,thick,inner sep=1.5pt]

\node at (-1,1) [ideal] {};
\node at (0,1) [ideal] {};
\node at (1,1) [ideal] {};
\node at (2,1) [ideal] {};
\node at (3,1) [ideal] {};
\node at (4,1) [ideal] {};
\node at (5,1) [ideal] {};
\node at (6,1) [ideal] {};

\node at (0,0) [ideal] {};
\node at (1,0) [ideal] {};
\node at (2,0) [ideal] {};
\node at (3,0) [ideal] {};
\node at (4,0) [ideal] {};
\node at (5,0) [ideal] {};
\node at (6,0) [ideal] {};

\node at (1,-1) [ideal] {};
\node at (2,-1) [ideal] {};
\node at (3,-1) [ideal] {};
\node at (4,-1) [ideal] {};
\node at (5,-1) [ideal] {};
\node at (6,-1) [ideal] {};

\node at (2,-2) [ideal] {};
\node at (3,-2) [ideal] {};
\node at (4,-2) [ideal] {};
\node at (5,-2) [ideal] {};
\node at (6,-2) [regular polygon,regular polygon sides=4,draw=black,fill=black,inner sep=1.5pt] {};

\node at (3,-3) [ideal] {};
\node at (4,-3) [ideal] {};
\node at (5,-3) [regular polygon,regular polygon sides=4,draw=black,fill=white,inner sep=1.5pt] {};
\node at (6,-3) [quotient] {};

\node at (4,-4) [regular polygon,regular polygon sides=4,draw=black,fill=black,inner sep=1.5pt] {};
\node at (5,-4) [quotient] {};
\node at (6,-4) [quotient] {};

\node at (5,-5) [quotient] {};
\node at (6,-5) [quotient] {};

\node at (6,-6) [quotient] {};
\node at (4,-4) [regular polygon,regular polygon sides=4,draw=red,inner sep=4pt,rotate=45] {};
\node at (5,-3) [regular polygon,regular polygon sides=4,draw=red,inner sep=4pt,rotate=45] {};

\node at (5,-3) [regular polygon,regular polygon sides=4,draw=blue,inner sep=4pt] {};
\node at (6,-2) [regular polygon,regular polygon sides=4,draw=blue,inner sep=4pt] {};

\end{tikzpicture}
}
\qquad \qquad
\subfloat[][$J = (x_2^3,x_2^2x_1,x_2x_1^2,x_1^5)$, $p_J(t) = 8$.]{\label{fig:problem1_b}
\begin{tikzpicture}[scale=0.6]
\tikzstyle{ideal}=[circle,draw=black,fill=black,inner sep=1.5pt]
\tikzstyle{quotient}=[circle,draw=black,fill=white,thick,inner sep=1.5pt]

\node at (7,1) [] {};
\node at (-1,1) [ideal] {};
\node at (0,1) [ideal] {};
\node at (1,1) [ideal] {};
\node at (2,1) [ideal] {};
\node at (3,1) [ideal] {};
\node at (4,1) [ideal] {};
\node at (5,1) [ideal] {};
\node at (6,1) [ideal] {};

\node at (0,0) [ideal] {};
\node at (1,0) [ideal] {};
\node at (2,0) [ideal] {};
\node at (3,0) [ideal] {};
\node at (4,0) [ideal] {};
\node at (5,0) [ideal] {};
\node at (6,0) [ideal] {};

\node at (1,-1) [ideal] {};
\node at (2,-1) [ideal] {};
\node at (3,-1) [ideal] {};
\node at (4,-1) [ideal] {};
\node at (5,-1) [ideal] {};
\node at (6,-1) [ideal] {};

\node at (2,-2) [ideal] {};
\node at (3,-2) [ideal] {};
\node at (4,-2) [ideal] {};
\node at (5,-2) [ideal] {};
\node at (6,-2) [quotient] {};

\node at (3,-3) [ideal] {};
\node at (4,-3) [ideal] {};
\node at (5,-3) [ideal] {};
\node at (6,-3) [quotient] {};

\node at (4,-4) [quotient] {};
\node at (5,-4) [quotient] {};
\node at (6,-4) [quotient] {};

\node at (5,-5) [quotient] {};
\node at (6,-5) [quotient] {};

\node at (6,-6) [quotient] {};
\end{tikzpicture}
}
\caption[Example of Borel rational deformations that can not be performed simultaneously because involving not disjoint sets of monomials.]{\label{fig:problem1} Example of two Borel rational deformations that can not be performed simultaneously because involving not disjoint sets of monomials. Both deformations involve the monomial $x_2 x_1^2 x_0^4$, thus swapping the monomials $x_2^2 x_0^5,x_1^5x_0^2$ with $x_2 x_1^2 x_0^4$ preserves the Borel property but the Hilbert polynomial changes.}
\end{center}
\end{figure}

To construct the graph we can use Algorithm \ref{alg:BorelGeneratorDFS} to determine the vertices and then we can apply Algorithm \ref{alg:allDeformations} on every Borel-fixed ideal to compute the edges (discarding repetitions). But we would like to add more edges to the graph, so we now discuss if it is possible to perform two Borel rational deformations of the same ideal simultaneously. Let us consider a Borel-fixed ideal $I \subset \K[x]$ with Hilbert polynomial $p(t)$ having two Borel degenerations: $J_1 = \left\langle\{I_r\}\setminus\mathcal{F}_{\alpha_1}(x^{\alpha_1})\cup\mathcal{F}_{\alpha_1}(x^{\beta_1})\right\rangle$ and $J_2 = \big\langle\{I\}\setminus\mathcal{F}_{\alpha_2}(x^{\alpha_2})$ $\left.\cup\mathcal{F}_{\alpha_2}(x^{\beta_2})\right\rangle$. The point is to understand under which conditions swapping at the same time $\mathcal{F}_{\alpha_1}(x^{\alpha_1})$ with $\mathcal{F}_{\alpha_1}(x^{\beta_1})$ and $\mathcal{F}_{\alpha_2}(x^{\alpha_2})$ with $\mathcal{F}_{\alpha_2}(x^{\beta_2})$ preserves the Borel property.
The first problems that can arise are:

\begin{enumerate}
 \item\label{it:problem_1} if the sets of monomials are not disjoint, there are some problems in the definition of the deformation. For instance if $\mathcal{F}_{\alpha_1}(x^{\alpha_1}) \cap \mathcal{F}_{\alpha_2}(x^{\alpha_2}) = \emptyset$ and $\mathcal{F}_{\alpha_1}(x^{\beta_1}) \cap \mathcal{F}_{\alpha_2}(x^{\beta_2}) \neq \emptyset$, the number of monomials we remove from the Borel set would be different from the number of monomials we add (see Figure \ref{fig:problem1});
 \item\label{it:problem_2} if $x^{\beta_2}$ can be obtained by a decreasing move from $x^{\alpha_1}$ (or viceversa $x^{\beta_1}$ from $x^{\alpha_2}$), performing both exchanges we would not obtain a Borel set $\mathscr{B}$, because $\mathscr{B} \ni x^{\beta_2} <_B x^{\alpha_1} \notin \mathscr{B}$ (see Figure \ref{fig:problem2});
 \item\label{it:problem_3} assuming $\min x^{\alpha_1} > \min x^{\alpha_2}$, if there exists an admissible move $\down{j}$, $j \geqslant \min x^{\alpha_1}$, such that $\down{j}(x^{\alpha_2}) <_B x^{\beta_1}$, it could happen that $\down{j}(x^{\alpha_2}) \in \mathcal{F}_{\alpha_1}(x^{\beta_1})$ and again swapping both sets of monomials we would not obtain a Borel set $\mathscr{B}$ because $\mathscr{B} \ni \down{j}(x^{\alpha_2}) <_B x^{\alpha_2} \notin \mathscr{B}$ (see Figure \ref{fig:problem2}).
\end{enumerate}

\newpage

\begin{theorem}\label{th:compatible}
 Let $I \subset \K[x]$ be a Borel-fixed ideal with Hilbert polynomial $p(t)$ with Gotzmann number $r$. Let us assume that there are $s$ Borel rational deformations: 
\[
\mathcal{I}_k = \left\langle \{I_r\} \setminus \mathcal{F}_{\alpha_k}(x^{\alpha_k})\cup\{y_{k\,0} F(x^{\alpha_k}) + y_{k\,1} F(x^{\beta_k})\ \vert\ F \in \mathcal{F}_{\alpha_k}\}\right\rangle \subset \K[y_{k\, 0},y_{k\, 1}][x],
\]
$k=1,\ldots,s$, such that
\begin{enumerate}[(i)]
 \item\label{it:compatible_hp0} $\min x^{\alpha_1} \geqslant \ldots \geqslant \min^{\alpha_s}$;
 \item\label{it:compatible_hp1} $\forall\ i,j$, $x^{\alpha_i}$ can not be obtained from $x^{\beta_j}$ by applying a decreasing move;
 \item\label{it:compatible_hp2} $\forall\ i,j$ such that $\min x^{\alpha_i} > \min x^{\alpha_j}$, $\down{h}(x^{\alpha_j})$ cannot be obtained from $x^{\beta_i}$ through a decreasing move, for all admissible $\down{h},\ h \geqslant \min x^{\alpha_i}$;
 \item\label{it:compatible_hp3} $\forall\ i,j$, $\mathcal{F}_{\alpha_i}(x^{\alpha_i}) \cap \mathcal{F}_{\alpha_j}(x^{\alpha_j}) = \emptyset$ and $\mathcal{F}_{\alpha_i}(x^{\beta_i}) \cap \mathcal{F}_{\alpha_j}(x^{\beta_j}) = \emptyset$.
\end{enumerate}
Then the family of subschemes of $\PP^n$ parametrized by the ideal
\begin{equation}\label{eq:multipleDef}
 \mathcal{I} = \left\langle \left\{I\right\}\setminus \left(\bigcup_{k} \mathcal{F}_{\alpha_k}(x^{\alpha_k})\right)\cup\left( \bigcup_{k} \left\{y_{k\,0}\, F(x^{\alpha_k}) + y_{k\,1}\, F(x^{\beta_k})\ \big\vert\ F \in\mathcal{F}_{\alpha_k}\right\}\right)\right\rangle
\end{equation}
is flat over $(\PP^1)^{\times s}$ and has $2^s$ fibers corresponding to Borel-fixed ideals.
\end{theorem}

\newpage

\vspace*{\stretch{1}}

\input{Figures/Ch3/problemsComposition2}

\vspace*{\stretch{1}}

\newpage

\begin{proof}
The key point is that the hypotheses ensure that there are not linear syzygies between two monomials belonging to two different groups, so assuming $y_{k\, 0} \neq 0,\ \forall\ k$ and $z_k = \frac{y_{k\, 1}}{y_{k\, 0}}$, we can lift simultaneously the syzygies of Eliahou-Kervaire among the monomials of $\{I_r\}$ as done in the proof of Theorem \ref{th:flatGeneralDim} for each group of binomials $F(x^{\alpha_k}) + z_k F(x^{\beta_k}),\ \forall\ F \in \mathcal{F}_{\alpha_k},\ \forall\ k=1,\ldots,s$. For each localization we deduce the flatness\index{flatness} by Proposition \ref{prop:artinSyzygies} and by the symmetry between the sets of monomials $\mathcal{F}_{\alpha_k}(x^{\alpha_k})$ and $\mathcal{F}_{\alpha_k}(x^{\beta_k})$, the property holds for any other choice of non-vanishing variables $y_{k\, h}$. 

The Borel-fixed ideals appearing the family correspond to the points of $\left(\PP^1\right)^{\times s}$ such that for each $k$ $[y_{k\,0}:y_{k\,1}] = [0:1]$ or $[y_{k\,0}:y_{k\,1}] = [1:0]$, and so there are $2^s$ possibilities.
\end{proof}

\begin{corollary}\label{cor:curveCompatible}
 Let $I \subset \K[x]$ be a Borel-fixed ideal and let $\mathcal{I}_1,\ldots,\mathcal{I}_s$ be $s$ Borel rational deformations of $I$ as in Theorem \ref{th:compatible}. For any Borel-fixed ideal $J$ belonging to the family over $\left(\PP^1\right)^{\times s}$, there exists a rational deformation having both $I$ and $J$ as fibers. 
 Called $p(t)$ the Hilbert polynomial of $I$ and $r$ its Gotzmann number, the points defined by $I$ and $J$ on the Hilbert scheme $\Hilb{n}{p(t)}$ are connected by a rational curve $\mathcal{C}:\PP^1 \rightarrow \Hilb{n}{p(t)}$. Considering the construction of the Hilbert scheme $\Hilb{n}{p(t)}$ as subscheme of the Grassmannian $\Grass{q(r)}{N}{\K}$, where $N = \binom{n+r}{n}$ and $q(r) = N-p(r)$, the degree of the curve $\mathcal{C}$ via the  Pl\"ucker embedding \eqref{eq:pluckerEmbedSubspace}  $\mathscr{P}: \Grass{q(r)}{N}{\K} \rightarrow \PP\big[\wedge^{q(r)} \K[x]_r\big]$ is $q(r) - \vert\{I_r\}\cap\{J_r\}\vert$.
\end{corollary}
\begin{proof}
Let $\mathcal{I}$ be the family defined in \eqref{eq:multipleDef}. To obtain a rational deformation having $I$ and any other Borel-fixed ideal $J$ of $\mathcal{I}$, it suffices to specialize the variables $y_{k\,0},y_{k\,1}$ as follows:
\begin{itemize}
\item $y_{k\, 0}= 1,\ y_{k\, 1} = 0$, if $\mathcal{F}_{\alpha_k}(x^{\alpha_k}) \subset \{I_r\}$ and $\mathcal{F}_{\alpha_k}(x^{\alpha_k}) \subset \{J_r\}$;
\item $y_{k\, 0}= y_0,\ y_{k\, 1} = y_1$, if $\mathcal{F}_{\alpha_k}(x^{\alpha_k}) \subset \{I_r\}$ and $\mathcal{F}_{\alpha_k}(x^{\beta_k}) \subset \{J_r\}$.
\end{itemize}
As usual, the ideal $I$ corresponds to the fiber of the point $[1:0]$ and the ideal $J$ to the fiber of $[0:1]$.

To determine the degree of the curve through the Pl\"ucker embedding, we can use again the matrix representation of the submodule $\mathcal{I}_r \subset (\K[y_0,y_1])[x]_r$ used in the proof of Theorem \ref{th:curveGeneralDim}. In this case the columns of the first block correspond to the monomials belonging to both $\{I_r\}$ and $\{J_r\}$, the columns of the second block to the monomials in $\{I_r\} \setminus\left(\{I_r\}\cap\{J_r\}\right)$ and the columns of the third block to the monomials in $\{J_r\} \setminus\left(\{I_r\}\cap\{J_r\}\right)$. Thus the degree of the Pl\"ucker coordinates in the variables $y_0$ and $y_1$ is equal to $\vert\{I_r\} \setminus\left(\{I_r\}\cap\{J_r\}\right)\vert = \vert\{I_r\}\vert -\vert\{I_r\}\cap\{J_r\}\vert = q(r)-\vert\{I_r\}\cap\{J_r\}\vert$.
\end{proof}

\begin{definition}
 Let $I \subset \K[x]$ be a Borel-fixed ideal and let $\mathcal{I}_1,\ldots,\mathcal{I}_s$ be $s$ Borel rational deformations of $I$ as in Theorem \ref{th:compatible}. We say that $\{\mathcal{I}_1,\ldots,\mathcal{I}_s\}$ is a set of \emph{compatible} Borel rational deformations. Given a second Borel-fixed ideal $J$ belonging to the family of ideals $\mathcal{I}$ over $(\PP^1)^{\times s}$ but not to the set of Borel degenerations of $I$ computed with Algorithm \ref{alg:allDegenerations}, we call \emph{composed Borel rational deformation}\index{Borel rational deformation!composed} the deformation defined in Corollary \ref{cor:curveCompatible} having both $I$ and $J$ as fibers and \emph{(composed) Borel rational degeneration} of $I$ the ideal $J$.
\end{definition}

\begin{example}
The two Borel rational deformations
\[
 \begin{split}
  \mathcal{I}_1& {} = \left\langle\{I_8\}\setminus \{x_3^2 x_0^6\} \cup \{y_0 \, x_3^2 x_0^6 + y_1 x_2^3 x_0^5\}\right\rangle \\
  \mathcal{I}_2& {} = \left\langle\{I_8\}\setminus\mathcal{F}_{x_2^2 x_1^6}(x_2^2 x_1^6)\cup\left\{y_0 F(x_2^2 x_1^6) + y_1 F(x_3 x_1^7)\ \vert\ F \in \mathcal{F}_{x_2^2 x_1^6}\right\}\right\rangle
 \end{split}
\]
described in Example \ref{ex:multipleDef} are in the hypothesis of Theorem \ref{th:compatible}, so there is a flat family over $\PP^1 \times \PP^1$ defined by the ideal
\[
\begin{split}
\mathcal{I} = \Big\langle \{I_8\} &{}\setminus \left\{x_3^2 x_0^6\right\} \cup \left\{ y_{00}\, x_3^2 x_0^6 + y_{01}\, x_2^3 x_0^5\right\} \\
&{} \setminus \mathcal{F}_{x_2^2 x_1^6}(x_2^2 x_1^6) \cup\left\{y_{10} F(x_2^2 x_1^6) + y_{11} F(x_3 x_1^7)\ \vert\ F \in \mathcal{F}_{x_2^2 x_1^6}\right\}\Big\rangle.
\end{split}
\]
Let us assume $y_{00} \neq 0$, $y_{10} \neq 0$, $z_0 = \frac{y_{01}}{y_{00}}$, $z_1 = \frac{y_{11}}{y_{10}}$ and let us see how the Eliahou-Kervaire syzygies of $I_{\geqslant 8}$ lift to a set of syzygies for the family $\mathcal{I}\vert_{[1:z_0],[1:z_1]}$. The presentation map is
\begin{equation}
\begin{split}
& \psi : \big(\K[z_0,z_1][x](-8)\big)^{136}\ \longrightarrow\ \K[z_0,z_1][x]\\
 & \mathbf{f}_{\gamma}\ \longmapsto \begin{cases}
 x^\gamma &  \forall\ x^\gamma \in \{I_8\}\setminus \left\{x_3^2 x_0^6\right\}\setminus \mathcal{F}_{x_2^2 x_1^6}(x_2^2 x_1^6) \\
 x_3^2 x_0^6 + z_0\, x_2^3 x_0^5 & x^\gamma = x_3^2 x_0^6\\
 F(x_2^2 x_1^6) + z_1 F(x_3 x_1^7), & \forall\ x^\gamma = F(x_2^2 x_1^6),\ F \in \mathcal{F}_{x_2^2 x_1^6}.
 \end{cases}
\end{split}
\end{equation}
For the monomials in $\{I_8\}\setminus \left\{x_3^2 x_0^6\right\}\setminus \mathcal{F}_{x_2^2 x_1^6}(x_2^2 x_1^6)$ we lift the syzygies among the same monomials in $I_8 \subset \K[x]$. Let us look at the binomials.
\begin{itemize}
  \item $\mathbf{f}_{x_3^2 x_0^6} \mapsto x_3^2 x_0^6 + z_0\, x_2^3 x_0^5$
\begin{description}
  \item[$(\cdot x_1)$] $x_1 \mathbf{f}_{x_3^2 x_0^6} - x_0\mathbf{f}_{x_3^2 x_1 x_0^5} - z_0 x_0 \mathbf{f}_{x_2^3 x_1 x_0^4}$;
  \item[$(\cdot x_2)$] $x_2 \mathbf{f}_{x_3^2 x_0^6} - x_0\mathbf{f}_{x_3^2 x_2 x_0^5} - z_0 x_0 \mathbf{f}_{x_2^4 x_0^4}$;
  \item[$(\cdot x_3)$] $x_3 \mathbf{f}_{x_3^2 x_0^6} - x_0\mathbf{f}_{x_3^3 x_0^5} - z_0 x_0 \mathbf{f}_{x_3 x_2^3 x_0^4}$.
\end{description}
  \item $\mathbf{f}_{x_2^2 x_1^6} \mapsto x_2^2 x_1^6 + z_1 \, x_3 x_1^7$
\begin{description}
\item[$(\cdot x_2)$] $x_2 \mathbf{f}_{x_2^2 x_1^6} - x_1 \mathbf{f}_{x_2^3 x_1^5} - z_1 x_1 \mathbf{f}_{x_3x_2 x_1^6}$;
\item[$(\cdot x_3)$] $x_3 \mathbf{f}_{x_2^2 x_1^6} - x_1 \mathbf{f}_{x_3 x_2^2 x_1^5} - z_1 x_1 \mathbf{f}_{x_3^2 x_1^6}$.
\end{description}
  \item $\mathbf{f}_{x_2^2 x_1^5 x_0} \mapsto x_2^2 x_1^5 x_0 + z_1 \, x_3 x_1^6 x_0$
\begin{description}
\item[$(\cdot x_1)$] $x_1 \mathbf{f}_{x_2^2 x_1^5 x_0} - x_0 \mathbf{f}_{x_2^2 x_1^6}$;
\item[$(\cdot x_2)$] $x_2 \mathbf{f}_{x_2^2 x_1^5 x_0} - x_0 \mathbf{f}_{x_2^3 x_1^5} - z_1 x_0 \mathbf{f}_{x_3x_2 x_1^6}$;
\item[$(\cdot x_3)$] $x_3 \mathbf{f}_{x_2^2 x_1^5 x_0} - x_0 \mathbf{f}_{x_3 x_2^2 x_1^5} - z_1 x_0 \mathbf{f}_{x_3^2 x_1^6}$.
\end{description}
  \item $\mathbf{f}_{x_2^2 x_1^4 x_0^2} \mapsto x_2^2 x_1^4 x_0^2 + z_1 \, x_3 x_1^5 x_0^2$
\begin{description}
\item[$(\cdot x_1)$] $x_1 \mathbf{f}_{x_2^2 x_1^4 x_0^2} - x_0 \mathbf{f}_{x_2^2 x_1^5 x_0}$;
\item[$(\cdot x_2)$] $x_2 \mathbf{f}_{x_2^2 x_1^4 x_0^2} - x_0 \mathbf{f}_{x_2^3 x_1^4 x_0} - z_1 x_0 \mathbf{f}_{x_3x_2 x_1^5 x_0}$;
\item[$(\cdot x_3)$] $x_3 \mathbf{f}_{x_2^2 x_1^4 x_0^2} - x_0 \mathbf{f}_{x_3 x_2^2 x_1^4 x_0} - z_1 x_0 \mathbf{f}_{x_3^2 x_1^5 x_0}$.
\end{description}
  \item $\mathbf{f}_{x_2^2 x_1^3 x_0^3} \mapsto x_2^2 x_1^3 x_0^3 + z_1 \, x_3 x_1^4 x_0^3$
\begin{description}
\item[$(\cdot x_1)$] $x_1 \mathbf{f}_{x_2^2 x_1^3 x_0^3} - x_0 \mathbf{f}_{x_2^2 x_1^4 x_0^2}$;
\item[$(\cdot x_2)$] $x_2 \mathbf{f}_{x_2^2 x_1^3 x_0^3} - x_0 \mathbf{f}_{x_2^3 x_1^3 x_0^2} - z_1 x_0 \mathbf{f}_{x_3x_2 x_1^4 x_0^2}$;
\item[$(\cdot x_3)$] $x_3 \mathbf{f}_{x_2^2 x_1^3 x_0^3} - x_0 \mathbf{f}_{x_3 x_2^2 x_1^3 x_0^2} - z_1 x_0 \mathbf{f}_{x_3^2 x_1^4 x_0^2}$.
\end{description}
  \item $\mathbf{f}_{x_2^2 x_1^2 x_0^4} \mapsto x_2^2 x_1^2 x_0^4 + z_1 \, x_3 x_1^3 x_0^4$
  \begin{description}
\item[$(\cdot x_1)$] $x_1 \mathbf{f}_{x_2^2 x_1^2 x_0^4} - x_0 \mathbf{f}_{x_2^2 x_1^3 x_0^3}$;
\item[$(\cdot x_2)$] $x_2 \mathbf{f}_{x_2^2 x_1^2 x_0^4} - x_0 \mathbf{f}_{x_2^3 x_1^2 x_0^3} - z_1 x_0 \mathbf{f}_{x_3x_2 x_1^3 x_0^3}$;
\item[$(\cdot x_3)$] $x_3 \mathbf{f}_{x_2^2 x_1^2 x_0^4} - x_0 \mathbf{f}_{x_3 x_2^2 x_1^2 x_0^3} - z_1 x_0 \mathbf{f}_{x_3^2 x_1^3 x_0^3}$.
\end{description}
\end{itemize}

Thus the Borel-fixed ideals belonging to this family are $4$:
\[
\begin{split}
&\quad I = (x_3^2,x_3 x_2^2, x_3 x_2 x_1,x_2^4,x_2^3 x_1,x_2^2 x_1^2)\\
& \left.\begin{array}{l}
\, (x_3^3,x_3^2 x_2, x_3 x_2^2, x_2^3,x_3^2 x_1,x_3 x_2 x_1,x_2^2 x_1^2)\\
\, (x_3^2,x_3 x_2^2,x_3 x_2 x_1,x_2^4,x_2^3 x_1,x_3 x_1^3) \\
\end{array}\right\} \text{ simple Borel degenerations}\\
&\quad (x_3^3,x_3^2 x_2,x_3 x_2^2,x_2^3,x_3^2 x_1, x_3 x_2 x_1,x_3 x_1^3)\ \text{ composed Borel degeneration}
\end{split}
\]
and for any pair among them, the points defined by the Borel-fixed ideals on the Hilbert scheme $\Hilb{3}{3t+5}$ are connected by a rational curve.
\end{example}

\begin{algorithm}[H]
 \caption[Pseudocode description of the method computing the Borel incidence graph.]{ The pseudocode description of the method computing the Borel incidence graph. For details see Appendix \ref{ch:HSCpackage}.}
 \label{alg:BorelIncidenceGraph}
 \begin{algorithmic}[1]
\STATE $\textsc{BorelIncidenceGraph}(\Hilb{n}{p(t)})$
\REQUIRE $\Hilb{n}{p(t)}$, a Hilbert scheme.
\ENSURE the Borel incidence graph $(\textsf{vertices},\textsf{edges})$ of $\Hilb{n}{p(t)}$.
\STATE $r \leftarrow \textsc{GotzmannNumber}\big(p(t)\big)$;
\STATE $\textsf{vertices} \leftarrow \textsc{BorelGeneratorDFS}\big(\K[x_0,\ldots,x_n],p(t)\big)$;
\STATE $\textsf{edges} \leftarrow \emptyset$;
\FORALL{$I \in \textsf{vertices}$}
\STATE $\textsf{simpleDeformations} \leftarrow \textsc{BorelRationalDeformations}(I,r)$;
\FORALL{$\{\mathcal{I}_1,\ldots,\mathcal{I}_s\} \subset \textsf{simpleDeformations}$}
\IF{$\{\mathcal{I}_1,\ldots,\mathcal{I}_s\}$ compatible set of Borel deformations}
\STATE $J \leftarrow \left\langle \{I_r\} \setminus \left(\bigcup_k \mathcal{F}_{\alpha_k}(x^{\alpha_k})\right) \cup \left(\bigcup_k \mathcal{F}_{\alpha_k}(x^{\beta_k})\right) \right\rangle^{\sat}$;
\IF{$(I,J) \notin \textsf{edges}$ \AND $(J,I) \notin \textsf{edges}$}
\STATE $\textsf{edges} \leftarrow \textsf{edges} \cup \{(I,J)\}$;
\ENDIF
\ENDIF
\ENDFOR
\ENDFOR
\RETURN $(\textsf{vertices},\textsf{edges})$;
 \end{algorithmic}
\end{algorithm}

Theorem \ref{th:compatible} and Corollary \ref{cor:curveCompatible} introduce other Borel rational deformations between pairs of Borel-fixed ideals. Hence in the algorithm computing the Borel incidence graph (Algorithm \ref{alg:BorelIncidenceGraph}) we add these composed Borel rational deformations to those obtained applying Algorithm \ref{alg:allDeformations} on every ideal.

Furthermore Theorem \ref{th:compatible} gives us a new criterion to detect points lying on a common component of the Hilbert scheme, indeed the points defined by Borel-fixed ideals belonging to any family over $\left(\PP^{1}\right)^{\times s}$ are on the same component. With the following examples we show that this criterion is a sufficient condition that does not overlap Reeves criterion \cite{Reeves} based on the hyperplane section.  

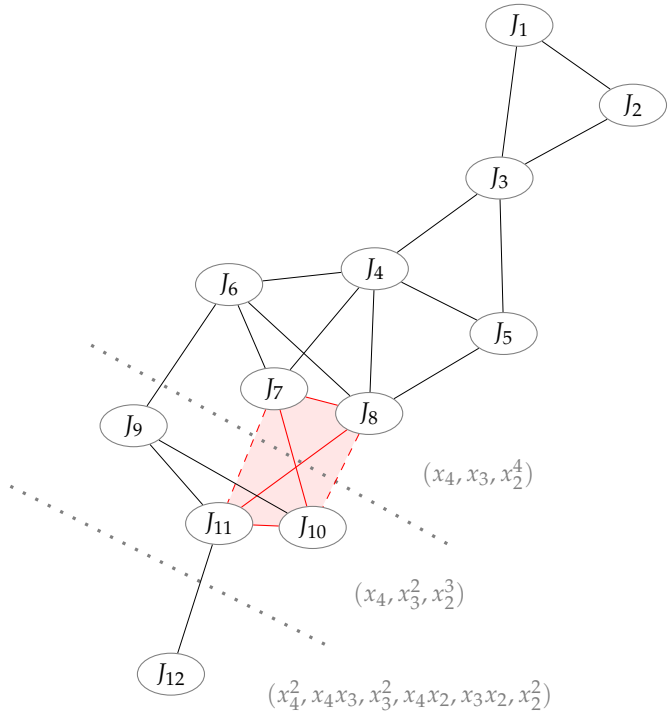
\begin{figure}[!ht]
\begin{center}
 \begin{tikzpicture}[>=latex,scale=0.75]
\tikzstyle{place}=[ellipse,draw=black!50,fill=white,minimum width=25pt,inner sep=1.5pt]
\filldraw [fill=red!10,draw=white] (79.058bp,152.12bp) -- (126.09bp,139.86bp) -- (98.117bp,83.314bp) -- (51.839bp,85.239bp) -- cycle;


\draw [very thick,black!50,loosely dotted] (-11.5bp,173.78bp) -- (168.117bp,73.314bp);

\draw [very thick,black!50,loosely dotted] (-51.5bp,103.78bp) -- (108.117bp,23.314bp);

\node at (180bp,110bp) [black!50] {\footnotesize $(x_4,x_3,x_2^4)$};
\node at (146.09bp,50bp) [black!50] {\footnotesize $(x_4,x_3^2,x_2^3)$};
\node at (146.09bp,0bp) [black!50] {\footnotesize $(x_4^2,x_4x_3,x_3^2,x_4x_2,x_3x_2,x_2^2)$};

\node (11) at (51.839bp,85.239bp) [place] {\footnotesize $J_{11}$};
  \node (10) at (98.117bp,83.314bp) [place] {\footnotesize $J_{10}$};
  \node (12) at (27.971bp,10.5bp) [place] {\footnotesize $J_{12}$};
  \node (1) at (200.4bp,332.5bp) [place] {\footnotesize $J_{1}$};
  \node (3) at (190.56bp,256.71bp) [place] {\footnotesize $J_{3}$};
  \node (2) at (256.64bp,292.87bp) [place] {\footnotesize $J_{2}$};
  \node (5) at (192.61bp,179.55bp) [place] {\footnotesize $J_{5}$};
  \node (4) at (128.86bp,211.69bp) [place] {\footnotesize $J_{4}$};
  \node (7) at (79.058bp,152.12bp) [place] {\footnotesize $J_{7}$};
  \node (6) at (56.737bp,203.84bp) [place] {\footnotesize $J_{6}$};
  \node (9) at (9.5bp,133.78bp) [place] {\footnotesize $J_{9}$};
  \node (8) at (126.09bp,139.86bp) [place] {\footnotesize $J_{8}$};
  
  \draw [very thin] (5) -- (8);
  \draw [very thin] (6) -- (7);
  \draw [dashed,red] (8) -- (10);
  \draw [very thin] (3) -- (4);
  \draw [very thin] (9) -- (11);
  \draw [dashed,red] (7) -- (11);
  \draw [very thin] (2) -- (3);
  \draw [red] (8) -- (11);
  \draw [very thin] (9) -- (10);
  \draw [red] (7) -- (10);
  \draw [very thin] (6) -- (8);
  \draw [very thin] (4) -- (6);
  \draw [very thin] (4) -- (5);
  \draw [very thin] (11) -- (12);
  \draw [very thin] (6) -- (9);
  \draw [very thin] (4) -- (8);
  \draw [red] (7) -- (8);
  \draw [very thin] (4) -- (7);
  \draw [very thin] (3) -- (5);
  \draw [very thin] (1) -- (2);
  \draw [red] (10) -- (11);
  \draw [very thin] (1) -- (3);
 \end{tikzpicture}
\caption[The Borel incidence graph of $\Hilb{4}{4t+1}$.]{\label{fig:BorelIncidence} The Borel incidence graph of $\Hilb{4}{4t+1}$. The light red quadrangle corresponds to the family over $\PP^1 \times \PP^1$ and the dashed lines correspond to composed Borel rational deformations. The dotted lines divide the ideals w.r.t. the hyperplane section.} 
\end{center}
\end{figure}

\begin{example}
Let us consider the Hilbert scheme containing the rational normal curve of degree 4, that is \gls{Hilb4tp1P4}. On it there are $12$ points defined by Borel-fixed ideals:
\small
\[
\begin{split}
&J_1 = (x_4,x_3,x_2^5,x_2^4x_1^3),\\ 
&J_2 = (x_4,x_3,x_2^6,x_2^5x_1,x_2^4x_1^2), \\
&J_3 = (x_4,x_3^2,x_3x_2,x_3x_1,x_2^5,x_2^4x_1^2),\\
&J_4 = (x_4,x_3^2,x_3x_2,x_3x_1^2,x_2^5,x_2^4x_1), \\
&J_5 = (x_4^2,x_4x_3,x_3^2,x_4x_2,x_3x_2,x_4x_1,x_3x_1,x_2^5,x_2^4x_1), \\
&J_6 = (x_4,x_3^2,x_3x_2,x_2^4,x_3x_1^3),\\
& J_7 = (x_4,x_3^2,x_3x_2^2,x_3x_2x_1,x_3x_1^2,x_2^4),\\
& J_8 = (x_4^2,x_4x_3,x_3^2,x_4x_2,x_3x_2,x_4x_1,x_3x_1^2,x_2^4), \\
& J_9 = (x_4,x_3^2,x_3x_2,x_2^4,x_2^3x_1),\\
& J_{10} = (x_4,x_3^2,x_3x_2^2,x_2^3,x_3x_2x_1),\\
& J_{11} = (x_4^2,x_4x_3,x_3^2,x_4x_2,x_3x_2,x_4x_1,x_2^3),\\
& J_{12} = (x_4^2,x_4x_3,x_3^2,x_4x_2,x_3x_2,x_2^2).\\
\end{split}
\]
\normalsize

The result by Reeves \cite{Reeves} says that there exists a component of the Hilbert scheme containing the points defined by Borel-fixed ideals with the same hyperplane section, i.e. in the case of $\Hilb{4}{4t+1}$ there are components (even not different) containing:
\[
\begin{split}
&\bullet\ J_1,J_2,J_3,J_4,J_5,J_6,J_7,J_8\text{ having hyperplane section }(x_4,x_3,x_2^4);\\
&\bullet\ J_9,J_{10},J_{11}\text{ having hyperplane section } (x_4,x_3^2,x_3x_2,x_2^3);\\
&\bullet\ J_{12}.
\end{split}
\]

Computing with Algorithm \ref{alg:BorelIncidenceGraph} the Borel incidence graph of $\Hilb{4}{4t+1}$ (Figure \ref{fig:BorelIncidence}), we find one composed deformation over $\PP^1 \times \PP^1$, containing the ideals $J_7,J_8,$ $J_{10},J_{11}$, so that the corresponding points lie on a same component.
\end{example}

\vspace*{\stretch{1}}

\vspace*{\stretch{1}}

\newpage

\begin{example}
On the Hilbert scheme \gls{Hilb6tm3P3} containing the complete intersections $(2,3)$ in $\PP^3$, there are 31 points defined by Borel-fixed ideals, with three different admissible hyperplane sections (of 6 points in $\PP^2$):\\
$(x_3,x_2^6)$

{\footnotesize
\[
\begin{array}{l l}
 I_{1} = (x_3,x_2^7,x_2^6x_1^6), & I_{14} = (x_3^2,x_3x_2,x_3x_1^4,x_2^7,x_2^6x_1^2), \\
 I_{2} = (x_3,x_2^8,x_2^7x_1,x_2^6x_1^5), & I_{15} = (x_3^2,x_3x_2^2,x_3x_2x_1,x_3x_1^3,x_2^7,x_2^6x_1^2), \\
  I_{3} = (x_3^2,x_3x_2,x_3x_1,x_2^7,x_2^6x_1^5), & I_{16} = (x_3^3,x_3^2x_2,x_3x_2^2,x_3^2x_1,x_3x_2x_1,x_3x_1^2,x_2^7,x_2^6x_1^2), \\
  I_{4} = (x_3,x_2^8,x_2^7x_1^2,x_2^6x_1^4), & I_{17} = (x_3^2,x_3x_2,x_3x_1^5,x_2^7,x_2^6x_1), \\
 I_{5} = (x_3^2,x_3x_2,x_3x_1,x_2^8,x_2^7x_1,x_2^6x_1^4) , & I_{18} = (x_3^2,x_3x_2^2,x_3x_2x_1,x_3x_1^4,x_2^7,x_2^6x_1), \\
 I_{6} = (x_3^2,x_3x_2,x_3x_1^2,x_2^7,x_2^6x_1^4), & I_{19} = (x_3^2,x_3x_2^2,x_3x_2x_1^2,x_3x_1^3,x_2^7,x_2^6x_1), \\
 I_{7} = (x_3,x_2^9,x_2^8x_1,x_2^7x_1^2,x_2^6x_1^3), & I_{20} = (x_3^3,x_3^2x_2,x_3x_2^2,x_3^2x_1,x_3x_2x_1,x_3x_1^3,x_2^7,x_2^6x_1), \\
 I_{8} = (x_3^2,x_3x_2,x_3x_1,x_2^8,x_2^7x_1^2,x_2^6x_1^3), & I_{21} = (x_3^2,x_3x_2,x_2^6,x_3x_1^6), \\
 I_{9} = (x_3^2,x_3x_2,x_3x_1^2,x_2^8,x_2^7x_1,x_2^6x_1^3), & I_{22} = (x_3^2,x_3x_2^2,x_3x_2x_1,x_2^6,x_3x_1^5), \\
 I_{10} = (x_3^2,x_3x_2,x_3x_1^3,x_2^7,x_2^6x_1^3), & I_{23} = (x_3^2,x_3x_2^2,x_3x_2x_1^2,x_3x_1^4,x_2^6), \\
 I_{11} = (x_3^2,x_3x_2^2,x_3x_2x_1,x_3x_1^2,x_2^7,x_2^6x_1^3), & I_{24} = (x_3^3,x_3^2x_2,x_3x_2^2,x_3^2x_1,x_3x_2x_1,x_3x_1^4,x_2^6), \\
 I_{12} = (x_3^2,x_3x_2,x_3x_1^3,x_2^8,x_2^7x_1,x_2^6x_1^2), & I_{25} = (x_3^2,x_3x_2^3,x_3x_2^2x_1,x_3x_2x_1^2,x_3x_1^3,x_2^6), \\
 I_{13} = (x_3^2,x_3x_2^2,x_3x_2x_1,x_3x_1^2,x_2^8,x_2^7x_1,x_2^6x_1^2), & I_{26} = (x_3^3,x_3^2x_2,x_3x_2^2,x_3^2x_1,x_3x_2x_1^2,x_3x_1^3,x_2^6); \\
\end{array}
\]}
$(x_3^2,x_3x_2,x_2^5)$
\footnotesize
\[
\begin{array}{ll}
 I_{27} = (x_3^2,x_3x_2,x_2^6,x_2^5x_1^2), & I_{29} = (x_3^2,x_3x_2^2,x_3x_2x_1^2,x_2^5), \\
I_{28} = (x_3^2,x_3x_2^2,x_3x_2x_1,x_2^6,x_2^5x_1), & I_{30} = (x_3^3,x_3^2x_2,x_3x_2^2,x_3^2x_1,x_3x_2x_1,x_2^5); \\
\end{array}
\]
\normalsize
$(x_3^2,x_3 x_2^2,x_2^4)$
\footnotesize
\[
I_{31} = (x_3^2,x_3x_2^2,x_2^4)
\]
\normalsize

In this case, Algorithm \ref{alg:BorelIncidenceGraph} finds 5 families over $\PP^1 \times \PP^1$.
\begin{itemize}
\item Starting from $I_9$ and its Borel deformations 
\[
\begin{split}
&\mathcal{I}_1\quad\text{defined by}\quad x_3x_1^2x_0^9,\ x_2^6x_1^2x_0^4,\ \mathcal{F}_{x_3x_1^2x_0^9}=\{\mathrm{id}\},\\
&\mathcal{I}_2\quad\text{defined by}\quad x_3x_1^2x_0^9,\ x_2^7x_0^5,\ \mathcal{F}_{x_3x_1^2x_0^9}=\{\mathrm{id}\},\\
&\mathcal{I}_3\quad\text{defined by}\quad x_3x_2x_0^{10},\ x_2^7x_0^5,\ \mathcal{F}_{x_3x_2x_0^{10}}=\{\mathrm{id}\},\\
&\mathcal{I}_4\quad\text{defined by}\quad x_3x_2x_0^{10},\ x_2^6x_1^2x_0^4,\ \mathcal{F}_{x_3x_2x_0^{10}}=\{\mathrm{id}\},\\
\end{split}
\]
we can determine 2 families over $\PP^1 \times \PP^1$, because both $\{\mathcal{I}_1,\mathcal{I}_3\}$ and $\{\mathcal{I}_2,\mathcal{I}_4\}$ are sets of compatible deformations:
\[
\begin{split}
\mathcal{I}_{1,3} = \Big\langle \left\{(I_{9})_{12}\right\}&{} \setminus \{x_3x_1^2x_0^9\} \cup \{ y_{00}x_3x_1^2x_0^9 + y_{01}x_2^6x_1^2x_0^4\}\\
&{} \setminus \{x_3x_2x_0^{10}\} \cup\{y_{10}x_3x_2x_0^{10} + y_{11}x_2^7x_0^5\}\Big\rangle,\\
\mathcal{I}_{2,4} = \Big\langle \left\{(I_{9})_{12}\right\}&{} \setminus \{x_3x_1^2x_0^9\} \cup \{ y_{00}x_3x_1^2x_0^9 + y_{01}x_2^7x_0^5\}\\
&{} \setminus \{x_3x_2x_0^{10}\} \cup\{y_{10}x_3x_2x_0^{10} + y_{11}x_2^6x_1^2x_0^4\}\Big\rangle.\\
\end{split}
\]
The Borel-fixed ideals belonging to $\mathcal{I}_{1,3}$ are $I_9$, $I_{11}$, $I_{12}$ and $I_{15}$ and the ideals belonging to $\mathcal{I}_{2,4}$ are $I_9$, $I_{10}$, $I_{13}$ and $I_{15}$. Note that $I_{15}$ turns out to be in both cases a composed Borel degeneration of $I_9$ but the rational deformations having $I_9$ and $I_{15}$ as fibers that we can construct applying Corollary \ref{cor:curveCompatible} on $\mathcal{I}_{1,3}$ and $\mathcal{I}_{2,4}$ are different (see Figure \ref{fig:doubleDef}).
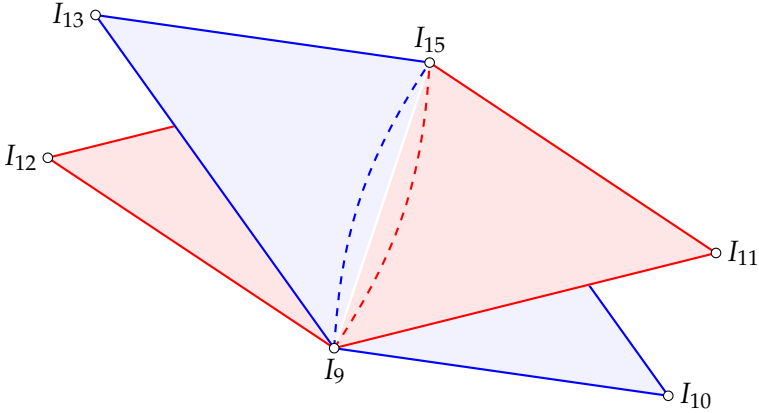
\begin{figure}
\begin{center}
\begin{tikzpicture}[scale=1.25]
\filldraw [draw = white,thick,fill=blue!5] (0,0) -- (3.5,-0.5) -- (1,3) -- cycle;
\draw [blue,thick] (0,0) -- (3.5,-0.5) -- (1,3);
\filldraw [draw = white,thick,fill=red!10] (0,0) -- (4,1) -- (1,3) -- cycle;
\draw [red,thick] (0,0) -- (4,1) -- (1,3);
\filldraw [draw = white,thick,fill=red!10] (0,0) -- (1,3) -- (-3,2) -- cycle;
\draw [red,thick] (1,3) -- (-3,2) -- (0,0);
\filldraw [draw = white,thick,fill=blue!5] (0,0) -- (1,3) -- (-2.5,3.5) -- cycle;
\draw [blue,thick] (1,3) -- (-2.5,3.5) -- (0,0);

\draw [red,dashed,thick] (0,0) to [bend right=15] (1,3);
\draw [blue,dashed,thick] (0,0) to [bend left=15] (1,3);

\node at (0,0) [below] {$I_{9}$};
\node at (1,3) [above] {$I_{15}$};

\node at (4,1) [right] {$I_{11}$};
\node at (-3,2) [left] {$I_{12}$};

\node at (3.5,-0.5) [right] {$I_{10}$};
\node at (-2.5,3.5) [left] {$I_{13}$};

\draw [black,fill=white] (0,0) circle (0.05);
\draw [black,fill=white] (1,3) circle (0.05);
\draw [black,fill=white] (4,1) circle (0.05);
\draw [black,fill=white] (-3,2) circle (0.05);
\draw [black,fill=white] (3.5,-0.5) circle (0.05);
\draw [black,fill=white] (-2.5,3.5) circle (0.05);
\end{tikzpicture}
\caption{\label{fig:doubleDef} The two family over $\PP^1 \times \PP^1$ that give rise to two composed Borel rational deformations containing the same pairs of Borel-fixed ideals, even if they do not coincide.}
\end{center}
\end{figure}
\item Considered the ideal $I_{18} = (x_3^2,x_3x_2^2,x_3x_2x_1,x_3x_1^4,x_2^7,x_2^6x_1)$ and its Borel rational deformations
\[
\begin{split}
&\overline{\mathcal{I}}_{1} \quad \text{defined by} \quad x_3x_2x_1x_0^9,\ x_2^6x_0^6,\ \mathcal{F}_{x_3x_2x_1x_0^9} = \{\mathrm{id}\},\\
&\overline{\mathcal{I}}_{2} \quad \text{defined by} \quad x_3x_2x_1x_0^9,\ x_3x_1^3x_0^8,\ \mathcal{F}_{x_3x_2x_1x_0^9} = \{\mathrm{id}\},\\
\end{split}
\]
\[
\begin{split}
&\overline{\mathcal{I}}_{3} \quad \text{defined by} \quad x_3^2x_0^{10},\ x_3x_1^3x_0^8,\ \mathcal{F}_{x_3^2x_0^{10}} = \{\mathrm{id}\},\\
&\overline{\mathcal{I}}_{4} \quad \text{defined by} \quad x_3^2x_0^{10},\ x_2^6x_0^6,\ \mathcal{F}_{x_3^2x_0^{10}} = \{\mathrm{id}\},\\
\end{split}
\]
we obtain again 2 families over $\PP^1 \times \PP^1$ by the compatibility of $\{\overline{\mathcal{I}}_{1},\overline{\mathcal{I}}_{3}\}$ and $\{\overline{\mathcal{I}}_{2},\overline{\mathcal{I}}_{4}\}$:
\[
\begin{split}
\overline{\mathcal{I}}_{1,3} = \Big\langle \left\{(I_{18})_{12}\right\}&{} \setminus \{x_3x_2x_1x_0^9\} \cup \{ y_{00}x_3x_2x_1x_0^9 + y_{01}x_2^6x_0^6\}\\
&{} \setminus \{x_3^2x_0^{10}\} \cup\{y_{10}x_3^2x_0^{10} + y_{11} x_3x_1^3x_0^8\}\Big\rangle,\\
\overline{\mathcal{I}}_{2,4} = \Big\langle \left\{(I_{18})_{12}\right\}&{} \setminus \{x_3x_2x_1x_0^9\} \cup \{ y_{00}x_3x_2x_1x_0^9 + y_{01} x_3x_1^3x_0^8\}\\
&{} \setminus \{x_3^2x_0^{10}\} \cup\{y_{10}x_3^2x_0^{10} + y_{11}x_2^6x_0^6\}\Big\rangle.\\
\end{split}
\]
$\overline{\mathcal{I}}_{1,3}$ contains $I_{18}$, $I_{20}$, $I_{23}$, $I_{26}$, $\overline{\mathcal{I}}_{2,4}$ contains $I_{18}$, $I_{19}$, $I_{24}$, $I_{26}$, and as before $I_{26}$ is a composed Borel degeneration of $I_{18}$ that can be obtained by means of 2 different composed Borel rational deformations.

\item Consider the ideal $I_{23}=(x_3^2,x_3x_2^2,x_3x_2x_1^2,x_3x_1^4,x_2^6)$ and the Borel rational deformations 
\[
\begin{split}
&\widetilde{\mathcal{I}}_1\quad\text{defined by}\quad x_3^2x_0^{10},\ x_3x_2x_1x_0^9, \ \mathcal{F}_{x_3^2x_0^{10}} = \{\mathrm{id}\},\\
&\widetilde{\mathcal{I}}_2\quad\text{defined by}\quad x_3x_1^{11},\ x_2^5x_1^7, \ \mathcal{F}_{x_3x_1^{11}} = \{\mathrm{id},\down{1},2\down{1},\ldots,7\down{1}\}.
\end{split}
\]
They are compatible and give rise to the family over $\PP^1 \times \PP^1$
\[
\begin{split}
\mathcal{I}_1 = \Big\langle \left\{(I_{23})_{12}\right\}&{} \setminus \{x_3^2x_0^{10}\} \cup \{ y_{00}x_3^2x_0^{10} + y_{01}x_3x_2x_1x_0^9\}\\
&{} \setminus \mathcal{F}_{x_3x_1^{11}}(x_3x_1^{11})\cup\left\{y_{10}F(x_3x_1^{11}) + y_{11}F(x_2^5x_1^7)\ \vert\ F \in \mathcal{F}_{x_3x_1^{11}}\right\}\Big\rangle
\end{split}
\]
containing 4 Borel-fixed ideal $I_{23},I_{24},I_{29}$ and $I_{30}$. Since one of the simple Borel deformations is defined by monomials in $\restrict{\{(I_{23})_{12}\}}{1}$, the hyperplane section of the 4 ideals could not be the same, so we deduce that the corresponding four points, that we can not discuss with Reeves criterion, lie on a common component of $\Hilb{3}{6t-3}$. 
\end{itemize}

\end{example}

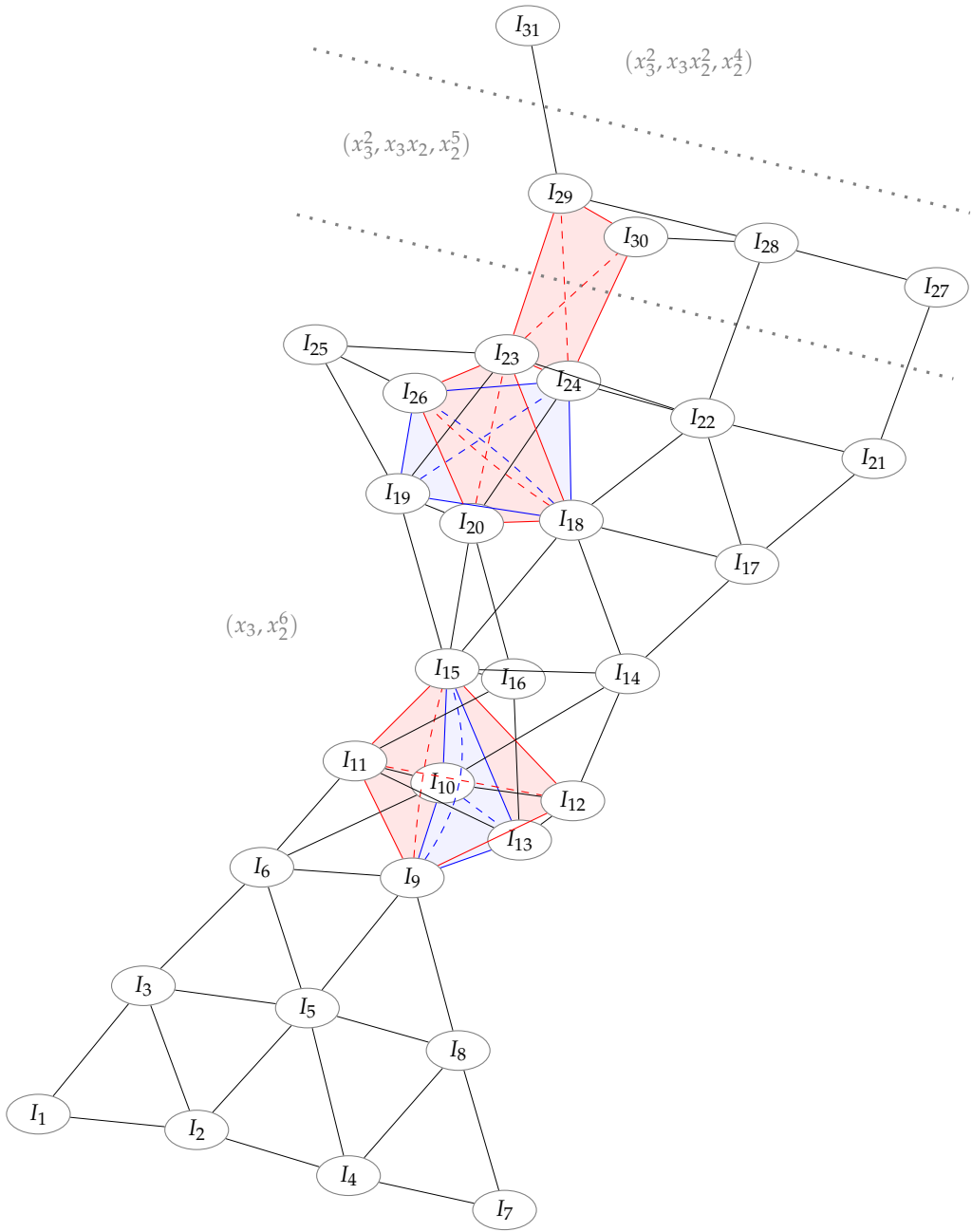
\begin{figure}[!ht]
\begin{center}
\begin{tikzpicture}[>=latex,scale=0.65]
\tikzstyle{place}=[ellipse,draw=black!50,fill=white,minimum width=25pt,inner sep=1.5pt]

\filldraw [draw=white,fill=red!10] (234.44bp,248.88bp) -- (200bp,319.39bp) -- (255.47bp,374.96bp) -- (331.05bp,295.34bp)  --cycle;
\filldraw [draw=white,fill=blue!5] (234.44bp,248.88bp) -- (252.76bp,306.11bp) -- (255.47bp,374.96bp) -- (299.05bp,271.69bp)  --cycle;

\filldraw [draw=white,fill=blue!5] (330.16bp,464.62bp) -- (225.78bp,480.64bp) -- (236.02bp,541.33bp) -- (328.62bp,548.69bp)  -- cycle;
\filldraw [draw=white,fill=red!10] (330.16bp,464.62bp) -- (270.16bp,462.73bp) -- (236.02bp,541.33bp)  -- (291.36bp,564.59bp)-- cycle;

\filldraw [draw=white,fill=red!10] (291.36bp,564.59bp) -- (328.62bp,548.69bp) -- (369.02bp,635.53bp) -- (323.52bp,661.79bp) -- cycle;

\draw [loosely dotted,black!50,-,very thick] (175bp,749bp) -- (570bp,650bp);

\draw [loosely dotted,black!50,-,very thick] (165bp,649bp) -- (560bp,550bp);

\node at (400.96bp,740.99bp) [black!50] {\footnotesize $(x_3^2,x_3x_2^2,x_2^4)$};
\node at (230.96bp,690.99bp) [black!50] {\footnotesize $(x_3^2,x_3x_2,x_2^5)$};
\node at (143.96bp,400.99bp) [black!50] {\footnotesize $(x_3,x_2^6)$};

  \node (24) at (328.62bp,548.69bp) [place] {\footnotesize $I_{24}$};
  \node (25) at (176.16bp,570.58bp) [place] {\footnotesize $I_{25}$};
  \node (26) at (236.02bp,541.33bp) [place] {\footnotesize $I_{26}$};
  \node (27) at (549.85bp,605.15bp) [place] {\footnotesize $I_{27}$};
  \node (20) at (270.16bp,462.73bp) [place] {\footnotesize $I_{20}$};
  \node (21) at (512.19bp,501.83bp) [place] {\footnotesize $I_{21}$};
  \node (22) at (409.14bp,526.23bp) [place] {\footnotesize $I_{22}$};
  \node (23) at (291.36bp,564.59bp) [place] {\footnotesize $I_{23}$};
  \node (28) at (447.54bp,631.92bp) [place] {\footnotesize $I_{28}$};
  \node (29) at (323.52bp,661.79bp) [place] {\footnotesize $I_{29}$};
  \node (1) at (9.5bp,106.34bp) [place] {\footnotesize $I_{1}$};
  \node (3) at (72.7bp,183.79bp) [place] {\footnotesize $I_{3}$};
  \node (2) at (104.76bp,96.843bp) [place] {\footnotesize $I_{2}$};
  \node (5) at (171.35bp,170.35bp) [place] {\footnotesize $I_{5}$};
  \node (4) at (196.1bp,69.246bp) [place] {\footnotesize $I_{4}$};
  \node (7) at (289.99bp,48.5bp) [place] {\footnotesize $I_{7}$};
  \node (6) at (143.95bp,255.27bp) [place] {\footnotesize $I_{6}$};
  \node (9) at (234.44bp,248.88bp) [place] {\footnotesize $I_{9}$};
  \node (8) at (262.06bp,144.72bp) [place] {\footnotesize $I_{8}$};
  \node (11) at (200bp,319.39bp) [place] {\footnotesize $I_{11}$};
  \node (10) at (252.76bp,306.11bp) [place] {\footnotesize $I_{10}$};
  \node (13) at (299.05bp,271.69bp) [place] {\footnotesize $I_{13}$};
  \node (12) at (331.05bp,295.34bp) [place] {\footnotesize $I_{12}$};
  \node (15) at (255.47bp,374.96bp) [place] {\footnotesize $I_{15}$};
  \node (14) at (363.95bp,372bp) [place] {\footnotesize $I_{14}$};
  \node (17) at (435.77bp,438.13bp) [place] {\footnotesize $I_{17}$};
  \node (16) at (295.3bp,368.95bp) [place] {\footnotesize $I_{16}$};
  \node (19) at (225.78bp,480.64bp) [place] {\footnotesize $I_{19}$};
  \node (18) at (330.16bp,464.62bp) [place] {\footnotesize $I_{18}$};
  \node (31) at (303.96bp,762.99bp) [place] {\footnotesize $I_{31}$};
  \node (30) at (369.02bp,635.53bp) [place] {\footnotesize $I_{30}$};
  \draw [very thin] (12) -- (13);
  \draw [red] (23) -- (29);
  \draw [very thin] (5) -- (8);
  \draw [blue] (19) -- (26);
  \draw [very thin] (14) -- (18);
  \draw [very thin] (29) -- (31);
  \draw [blue] (24) -- (26);
  \draw [blue] (18) -- (19);
  \draw [very thin] (14) -- (15);
  \draw [blue] (10) -- (15);
  \draw [dashed,red] (9) to [bend left=5] (15);
  \draw [dashed,blue] (9) to [bend right=25] (15);
  \draw [very thin] (23) -- (25);
  \draw [very thin] (19) -- (23);
  \draw [very thin] (6) -- (10);
  \draw [red] (20) -- (26);
  \draw [very thin] (1) -- (2);
  \draw [red] (9) -- (11);
  \draw [very thin] (18) -- (22);
  \draw [very thin] (10) -- (11);
  \draw [very thin] (12) -- (14);
  \draw [very thin] (13) -- (16);
  \draw [very thin] (2) -- (3);
  \draw [red] (23) -- (24);
  \draw [blue] (13) -- (15);
  \draw [red,dashed] (18) to [bend left=5] (26);
  \draw [blue,dashed] (18) to [bend right = 5] (26);
  \draw [very thin] (2) -- (5);
  \draw [very thin] (8) -- (9);
  \draw [very thin] (15) -- (20);
  \draw [dashed,red] (24) -- (29);
  \draw [red] (29) -- (30);
  \draw [very thin] (21) -- (22);
  \draw [red] (24) -- (30);
  \draw [very thin] (6) -- (11);
  \draw [red] (23) -- (26);
  \draw [dashed,red] (23) -- (30);
  \draw [blue] (9) -- (10);
  \draw [very thin] (19) -- (25);
  \draw [very thin] (28) -- (30);
  \draw [very thin] (21) -- (27);
  \draw [red] (12) -- (15);
  \draw [very thin] (11) -- (13);
  \draw [very thin] (2) -- (4);
  \draw [dashed,blue] (10) -- (13);
  \draw [very thin] (15) -- (16);
  \draw [very thin] (10) -- (12);
  \draw [dashed,red] (20) -- (23);
  \draw [very thin] (15) -- (19);
  \draw [very thin] (14) -- (17);
  \draw [very thin] (17) -- (18);
  \draw [dashed,blue] (19) -- (24);
  \draw [very thin] (4) -- (5);
  \draw [very thin] (11) -- (16);
  \draw [very thin] (22) -- (24);
  \draw [very thin] (17) -- (22);
  \draw [blue] (9) -- (13);
  \draw [red] (18) -- (20);
  \draw [very thin] (22) -- (23);
  \draw [dashed,red] (11) -- (12);
  \draw [very thin] (6) -- (9);
  \draw [very thin] (17) -- (21);
  \draw [very thin] (20) -- (24);
  \draw [very thin] (3) -- (6);
  \draw [very thin] (4) -- (8);
  \draw [very thin] (5) -- (6);
  \draw [very thin] (22) -- (28);
  \draw [very thin] (15) -- (18);
  \draw [very thin] (5) -- (9);
  \draw [very thin] (4) -- (7);
  \draw [very thin] (25) -- (26);
  \draw [very thin] (28) -- (29);
  \draw [very thin] (7) -- (8);
  \draw [very thin] (3) -- (5);
  \draw [very thin] (16) -- (20);
  \draw [very thin] (10) -- (14);
  \draw [red] (9) -- (12);
  \draw [very thin] (27) -- (28);
  \draw [red] (11) -- (15);
  \draw [blue] (18) -- (24);
  \draw [very thin] (19) -- (20);
  \draw [very thin] (1) -- (3);
  \draw [red] (18) -- (23);
\end{tikzpicture}
\caption{\label{fig:incidenceBig} The Borel incidence graph of $\Hilb{3}{6t-3}$.}
\end{center}
\end{figure}

\chapter{Borel open covering of Hilbert schemes}\label{ch:openSubsets}

In Chapter \ref{ch:HilbertScheme}, we showed that describing explicitly a Hilbert scheme is a very hard task, even in the easiest (from a geometric point of view) cases, because the costruction of $\Hilb{n}{p(t)}$ as subscheme of a Grassmannian requires a huge number of variables. The first idea to reduce the complexity is to study the Hilbert scheme locally, i.e. to consider the affine open covering of the Grassmannian and to look at the intersection between the Hilbert scheme and each open subset. Indeed to study the Grassmannian $\Grass{q}{N}{\K}$ globally, we have to use the Pl\"ucker embedding \eqref{eq:pluckerEmbedSubspace} that requires $\binom{N}{q}$ (projective) variables, whereas if we want to study $\Grass{q}{N}{\K}$ locally, we can consider its open covering having open subsets described by $q(N-q)$ (affine) variables. 

Since the Grassmannians considered in the study of Hilbert schemes actually parametrize ideals in a polynomial ring, the second idea is to exploit the algorithmic tools developed by the computational algebra, particularly the theory of Gr\"obner basis. The application of Gr\"obner bases to the study of Hilbert schemes was already introduced by Carr\`a Ferro in \cite{CarraFerro} and but our interest originates mainly in the ideas exposed in the paper \cite{NotariSpreafico} by Notari and Spreafico.

\medskip

The results I will expose in this chapter belong to several joint papers \cite{LellaRoggero,CioffiRoggero,BCLR,BLR} with M. Roggero, F. Cioffi and C. Bertone.

\section{Gr\"obner strata}\label{sec:GrobnerStrata}

\begin{definition}\label{def:tail}
Let us consider any term ordering $\sigma$ and a monomial ideal $J \subset \K[x]$ (not even Borel-fixed). We define the \emph{homogeneous tail}\index{homogeneous tail} (in the following tail for short) of $x^{\alpha} \in J$ as the set of monomials:
\begin{equation}
\gls{tail} = \left\{ x^\beta \notin J \text{ s.t. } \vert\beta\vert=\vert\alpha\vert \ \text{and}\ x^\beta <_{\sigma} x^{\alpha}\right\}
\end{equation}
\end{definition}

Every ideal $I \subset \K[x]$, having $J=(x^{\alpha_1},\ldots,x^{\alpha_s})$ as initial ideal w.r.t. $\sigma$, has a reduced Gr\"obner basis of $\{f_1,\dots,f_s\}$ where:
\begin{equation}\label{eq:polibasepartenza1}
  f_i = x^{\alpha_i} + \sum_{x^{\beta} \in \tail{x^{\alpha_i}}{\sigma}{J}} {c}_{\alpha_i \beta} x^{\beta}, \qquad {c}_{\alpha_i \beta}\in \K.
\end{equation}
Thus it is very natural to parameterize the family of all the ideals $I$ such that $\IN_\sigma(I) = J$ by the coefficients ${c}_{\alpha_i \beta}$; in this way the fanily corresponds to a subset of $\K^{\vert\tail{J}{\sigma}{}\vert}$, where $\tail{J}{\sigma}{}=  \tail{x^{\alpha_1}}{\sigma}{J} \times \cdots \times \tail{x^{\alpha_s}}{\sigma}{J}$.

\begin{definition} \label{def:codaridotta}  
Given a monomial ideal $J = (x^{\alpha_1},\ldots,x^{\alpha_s}) \subset \K[x]$ and a term ordering $\sigma$, let us fix a subset $T_i$ for each tail $\tail{x^{\alpha_i}}{\sigma}{J}$. Set $T = \{T_1,\ldots,T_s\}$, we will denote by $\St{\sigma}(J,T)$ the family of all ideals $I$ in $\K[x]$ such that $\IN_{\sigma}(I)=J$ and whose reduced Gr\"obner basis ${f}_1,\ldots,{f}_s$ is of the type:
\begin{equation}\label{eq:polibasepartenza2}
  {f}_i = x^{\alpha_i} + \sum_{x^{\beta} \in T_{i}} {c}_{\alpha_i \beta} x^{\beta},
\end{equation}
and we will call it \emph{$T$-Gr\"obner stratum} of $J$.
Moreover we will use \gls{GBstratum}, and we will call it \emph{Gr\"obner stratum}\index{Gr\"obner stratum} of $J$, whenever we consider the complete tail of every generator of $J$, i.e. $T_i = \tail{x^{\alpha_i}}{\sigma}{J},\ \forall\ i$.
\end{definition} 

\begin{remark}\label{rk:termordersignificance}
It will be clearer later that the term ordering affects the construction of a Gr\"obner stratum only because it determines which monomials can belong to the tails; indeed two different term orderings giving the same tails will lead to the same Gr\"obner stratum.
\end{remark}

Every ideal $I$ in the family $\St{\sigma}(J,T)$  is uniquely determined by a point in the affine space $\AA^N$ ($N=\sum_i \vert T_i \vert$) where we fix coordinates ${C}_{\alpha_i\beta}$ corresponding to the coefficients $c_{\alpha_i \beta}$ that appear in \eqref{eq:polibasepartenza2}. The subset of $\AA^N$ corresponding to $\St{\sigma}(J,T)$ turns out to be a closed algebraic set. More precisely, we will see how it can be endowed in a very natural way with a structure of affine subscheme, possibly reducible or non reduced, that is we will see that it can be obtained as the subscheme of $\AA^N$ defined by an ideal $\mathfrak{h}(J,T)$ in $\K[C]$, where $C$ is the set of variables $C_{\alpha_i \beta}$.

\begin{definition} \label{def:procedura} 
We will denote by $\mathfrak{h}(J,T)$ and  $\mathfrak{L}(J,T)$ the ideals in $\K[C]$ generated by the following procedure.
\begin{description}
\item[Step 1] Consider the set of polynomials $\mathcal{B} = \{F_1,\ldots,F_s\}$ such that
\begin{equation}\label{eq:polibasegenerica}
  F_i = x^{\alpha_i} + \sum_{x^{\beta} \in T_i}  {C}_{\alpha_i \beta} x^\beta \in \K[C][x].
\end{equation}
\item[Step 2] Consider a set $\text{\glshyperlink[$\syz(J)$]{syz}} = \left\{ x^{\gamma} \mathbf{e}_i - x^\delta \mathbf{e}_j\ \vert\ x^\gamma x^{\alpha_i} - x^{\delta}x^{\alpha_j}\right\}$ that generates the syzygies of $J$.
\item[Step 3a] For every $x^{\gamma} \mathbf{e}_i - x^\delta \mathbf{e}_j \in \syz(J)$, compute a complete reduction w.r.t. $\mathcal{B}$ of the $S$-polynomial $\glslink{Spoly}{S(F_i,F_j)} = x^{\gamma} F_i - x^\delta F_j$: $S(F_i,F_j) \stackrel{\mathcal{B}}{\longrightarrow} R_{ij}$.
\item[Step 3b] For every $x^{\gamma} \mathbf{e}_i - x^\delta \mathbf{e}_j \in \syz(J)$, compute a complete reduction w.r.t. $J$ of the $S$-polynomial $S(F_i,F_j) = x^{\gamma} F_i - x^\delta F_j$: $S(F_i,F_j) \stackrel{J}{\longrightarrow} M_{ij}$.
\item[Step 4a] Call $\mathfrak{h}(J,T)$ the ideal of $\K[C]$ generated by the coefficients (polynomials in $\K[C]$) of the reduced polynomials $R_{ij}$ computed at \textbf{Step 3a}.
\item[Step 4b] Call $\mathfrak{L}(J,T)$ the ideal of $\K[C]$ generated by the coefficients of the reduced polynomials $M_{ij}$ computed at \textbf{Step 3b}. Note that by construction the coefficients in $M_{ij}$ are linear, so actually $\mathfrak{L}(J,T)$ turns out to be a vector subspace of the vector space $\langle C \rangle$ spanned by the variables $C$.
\end{description}
\end{definition}

It is almost evident, that the definition of $\mathfrak{h}(J,T)$ is nothing else than Buchberger's characterization of Gr\"obner basis if we think to the $C_{\alpha_i\beta}$'s as constant in $\K$ instead of variables. In fact the variables $C_{\alpha_i \beta}$ do not appear in the leading terms w.r.t. $\sigma$ of $F_i$ and so their specialization in $\K$ commutes with reduction with respect to $\mathcal{B}$. Thus $(\ldots,{c}_{\alpha_ \beta},\ldots)$ is a closed point in the support of \glshyperlink[$Z \big(\mathfrak{h}(J,T)\big)$]{ZariskiClosed} in $\AA^N$ if and only if it corresponds to polynomials $f_1,\dots,f_s$ in $\K[x]$ that form a Gr\"obner basis. 
Then the support of \glshyperlink[$Z \big(\mathfrak{h}(J,T)\big)$]{ZariskiClosed} is uniquely defined; however a priori the ideal $\mathfrak{h}(J,T)$ could depend on the choices we perform computing it, that is (1) on the choice of the set $\syz(J)$ of generators of the syzygies and (2) on the choices did during the reduction of any $S$-polynomial $S(F_i,F_j)$ (which in general is not uniquely determined).  

Thanks again to Buchberger's criterion, we can prove that indeed $\mathfrak{h}(J,T)$ only depends on $J$, $T$ and of course on $\sigma$, because it can be defined in an equivalent intrinsic way. 

\begin{proposition}\label{th:unica} Let  $J \subset \K[x]$ be a monomial ideal and let $\sigma$ be any term ordering. Consider the set $\mathcal{B}=\{F_1,\dots, F_s\},\ F_i \in \K[C][x]$ as in \eqref{eq:polibasegenerica} and an ideal $\mathfrak{a}$ in $\K[C]$ with  Gr\"obner basis $\mathcal{A}$. The following conditions are equivalent: 
\begin{enumerate}[(i)]
\item\label{it:unica_1} $\mathcal{B} \cup \mathcal{A} $ is a Gr\"obner basis in $\K[C,x]$;
\item\label{it:unica_1'} $\mathfrak{a}$ contains the coefficients (polynomials in $\K[C]$) of all the polynomials in the ideal  $(F_1,\dots, F_s)$ that are reduced modulo $J$;
\item\label{it:unica_2} $\mathfrak{a}$  contains  all the coefficients of every complete reduction of $S(F_i,F_j)$ with respect to $\mathcal{B}$ for every $i,j$;  
\item\label{it:unica_3} $\mathfrak{a}$ contains  all the coefficients of some (even partial) reduction  with respect to $\mathcal{B}$ of $S(F_i,F_j)$ for every $i,j$;
\item\label{it:unica_4} $\mathfrak{a}$ contains  all the coefficients of some (even partial) reduction  with respect to $\mathcal{B}$ of $S(F_i,F_j)$, for every $(i,j)$ corresponding to a set $\syz(J)$ of generators of the syzygies of $J$.
\end{enumerate}
\end{proposition}
\begin{proof} $\emph{(\ref{it:unica_1})} \Rightarrow \emph{(\ref{it:unica_1'})}$. Let $G$ be a polynomial in $(F_1,\dots, F_s)$ which is reduced modulo $J$. By  hypothesis, $G$ must be reducible to 0 through $\mathcal{B} \cup \mathcal{A}$, so that the further steps of reduction have to be performed just using $\mathcal{A}$. Any step of reduction through $\mathcal{A}$ does not change the monomials in $\K[x]$ but only modifies their coefficients (in $\K[C]$), then  $G \stackrel{\mathcal{A}}{\longrightarrow} 0$, that is every coefficient in $\K[C]$ of $G$ can be reduced to $0$ using $\mathcal{A}$: this shows that all the coefficients in $\K[x]$ of $G$ belong to $\mathfrak{a}$.

$\emph{(\ref{it:unica_1'})} \Rightarrow \emph{(\ref{it:unica_2})}$, $\emph{(\ref{it:unica_2})} \Rightarrow \emph{(\ref{it:unica_3})}$ and $\emph{(\ref{it:unica_3})} \Rightarrow \emph{(\ref{it:unica_4})}$ are obvious.

$\emph{(\ref{it:unica_4})} \Rightarrow \emph{(\ref{it:unica_1})}$. We can check that $\mathcal{B} \cup \mathcal{A} $ is a Gr\"obner basis using the refined Buchberger criterion (see for instance \cite[Theorem 9, p. 104]{CLOiva}). If $\mathcal{A} =\{a_1,\dots,a_r\}$, a set of generators for the syzygies of the ideal $(x^{\alpha_1}, \dots, x^{\alpha_s}, \IN_{\sigma}(a_1), \dots, \IN_{\sigma}(a_r))$ can be obtained as the union of a set of generators of the syzygies of $(x^{\alpha_1}, \dots, x^{\alpha_s})$, a set of generators of the syzygies of $\IN_{\sigma}(\mathfrak{a}) = \big(\IN_{\sigma}(a_1), \dots, \IN_{\sigma}(a_r)\big)$ and the trivial syzygies $\big(x^{\alpha_i}, \IN_{\sigma}(a_j)\big)$. Then: 
\begin{itemize}
 \item $S(a_i,a_j) \xrightarrow{\mathcal{B} \cup \mathcal{A}} 0$, since $\mathcal{A}$ is a Gr\"obner basis and $\mathcal{A} \subset \mathcal{B}\cup\mathcal{A}$;
 \item $S(F_i,a_j) \xrightarrow{\mathcal{B} \cup \mathcal{A}} 0$, since the leading terms of $F_i$ and $a_j$ are coprime and $F_i,a_j \in \mathcal{B}\cup\mathcal{A}$;
 \item $S(F_i,F_j) \xrightarrow{\mathcal{B}\cup\mathcal{A}} 0$ in at least one way, by hypothesis. \qedhere
\end{itemize} 
\end{proof}

There are many ideals $\mathfrak{a}$ fulfilling the equivalent conditions of Proposition \ref{th:unica}: for instance we can consider the irrelevant maximal ideal in $\K[C]$ or any ideal obtained accordingly with condition \emph{(\ref{it:unica_3})}. Moreover, if $\mathfrak{a}$ satisfies those conditions and $\mathfrak{a}' \supset \mathfrak{a}$, then also $\mathfrak{a}'$  does, and if several ideals $\mathfrak{a}_l$ satisfy the conditions, then also their intersection $\bigcap_l \mathfrak{a}_l$ does. As a consequence of these remarks we obtain the proof of the uniqueness of the ideal $\mathfrak{h}(J,T)$ given by Definition \ref{def:procedura}.

\begin{theorem}\label{th:cor-unica}  Let $J$ be a monomial ideal and $T$ be the list of subsets of the tails of $J$ as above. Then:   
\begin{enumerate}[(i)]
\item\label{it:cor-unica_i} $\mathfrak{h}(J,T)$ is uniquely defined; indeed $\mathfrak{h}(J,T)=\bigcap_l \mathfrak{a}_l$, $\mathfrak{a}_l$  satisfying the equivalent conditions of Proposition \ref{th:unica};
\item\label{it:cor-unica_ii} $\mathfrak{L}(J,T)$ is uniquely defined.
\end{enumerate}
\end{theorem}
\begin{proof} 
\emph{(\ref{it:cor-unica_i})} $\mathfrak{h}(J,T)$ is one of the ideals $\mathfrak{a}_l$, because it satisfies condition \emph{(\ref{it:unica_4})} of Proposition \ref{th:unica}; on the other hand, if $\mathfrak{a}_l$ satisfies condition \emph{(\ref{it:unica_2})} of Proposition \ref{th:unica}, then clearly $\mathfrak{a} \supseteq \mathfrak{h}(J,T)$.

\emph{(\ref{it:cor-unica_ii})} It suffices to observe that the generators for $\mathfrak{L}(J,T)$ are the degree 1 homogeneous components of the generators of $\mathfrak{h}(J,T)$ given in its construction (Definition \ref{def:procedura}).
\end{proof}

By abuse of notation we will denote by the same symbol $\St{\sigma}(J,T)$ the family of ideals and the subscheme in $\AA^N$ given by the ideal $\mathfrak{h}(J,T)$.  Note that $\mathfrak{h}(J,T)$ is not always a prime ideal and so $\St{\sigma}(J,T)$  is not necessarily irreducible nor reduced, as the following trivial example shows.

\begin{example} Let $J=(x_2^2,x_2x_1)\subset \K[x_0,x_1,x_2]$ and $\sigma$ be any term ordering. Let us choose $T=\left\{T_{x_2^2}=\emptyset,T_{x_2x_1}=\{x_1x_0\}\right\}$ and  construct the ideal of the $T$-Gr\"obner stratum $\St{\sigma}(J,T)$ according to Definition \ref{def:procedura}:  
\[
\begin{split}
& F_{x_2^2}=x_2^2, \quad F_{x_2x_1}=x_2x_1+C x_1x_0,\\
& S(F_{x_2^2},F_{x_2x_1})=x_1F_{x_2^2}-x_2 F_{x_2x_1}={}\\
& \phantom{S(F_{x_2^2},F_{x_2x_1})} = -C x_2 x_1 x_0 \xrightarrow{F_{x_2^2},F_{x_2x_1}} -C x_2 x_1 x_0 + C x_0 F_{x_2 x_1 } = C^2 x_1x_0^2.
\end{split}
\] 
Then $\mathfrak{h}(J,T)=(C^2)$ that is $\St{\sigma}(J,T)$ is a double point on the affine line $\AA^1$.
\end{example}

\subsection{Gr\"obner strata are homogeneous varieties}

Let us denote by \gls{monomialGroup} the group of the monomials in the field of fraction $\K(x)$ of the ring $\K[x]$. Any term ordering $\sigma$ makes $\overline{\mathbb{T}}_x$ a totally ordered group in the obvious way.

\begin{definition}\label{def:lambda} Let us denote by $\ell$ the grading induced on $\K[C]$ by any term ordering $\sigma$ on $\K[x]$ through the map
\begin{equation}
\begin{split}
\ell:\K[C]&{} \longrightarrow \overline{\mathbb{T}}_x,\\
C_{\alpha_i \beta} &{} \longmapsto \frac{x^{\alpha_i}}{x^\beta}. 
\end{split}
\end{equation}
\end{definition}

 As we will use also the usual grading over $\ZZ$ where all the variables have degree 1, we will always write explicitly   the symbol $\ell$ when the above defined grading is concerned (so, $\ell$-degree $\gamma$ with $\gamma\in \ZZ^{n+1}$, $\ell$-homogeneous of degree $\gamma$  etc.) and we will leave the simple terms when the usual grading is involved (so, degree $r$ with $r\in \ZZ$, homogeneous of degree $r$ etc.).
 
\begin{proposition}[{\cite[Lemma 2.8]{RoggeroTerracini}}]\label{th:lambda1}  
\begin{enumerate}[(i)]
\item\label{it:lambda1_i} The grading $\ell$ is positive. 
\item\label{it:lambda1_ii} $\mathfrak{h}(J,T)$ is a $\ell$-homogeneous ideal.
\end{enumerate}
\end{proposition}
\begin{proof} \emph{(\ref{it:lambda1_i})}  Let us observe that all the variables have $\ell$-degree higher than that of the constant 1. Indeed
\[
\ell(C_{\alpha_i \beta}) >_\ell \ell(1) \quad \Longleftrightarrow\quad \frac{x^{\alpha_i}}{x^\beta} >_{\sigma} 1 \quad \Longleftrightarrow\quad x^{\alpha_i} >_\sigma x^\beta.
\]
As well known, this condition is equivalent to the positivity of the grading (see \cite[Chapter 4]{KreuzerRobbiano2}).

\emph{(\ref{it:lambda1_ii})} Let us consider the grading on $\K[C,x]$ induced by the map $\ell: \K[C,x] \rightarrow \overline{\mathbb{T}}_x$ sending $\ell(x_j) = x_j$ and $\ell(C_{\alpha_i \beta}) = \frac{x^{\alpha_i}}{x^\beta}$ and note that it coincides with the grading introduced in Definition \ref{def:lambda} on the restriction to $\K[C]$. Every monomial that appears in $F_i$ is of the type $C_{\alpha_i \beta}x^{\beta}$ and so its $\ell$-degree is $\ell(C_{\alpha_i \beta}x^{\beta})=\ell(C_{\alpha_i \beta})\cdot\ell(x^{\beta})= \frac{x^{\alpha_i}}{x^\beta} x^\beta = x^{\alpha_i}$.
Thus all the polynomials $F_i$ are $\ell$-homogeneous and then also the $S$-polynomials $S(F_i,F_j)$ and their reductions are $\ell$-homogeneous. Finally, the coefficients of any monomial $x^\gamma$ (which are polynomials in $\K[C]$) in such reductions are $\ell$-homogeneous. 
\end{proof}

\begin{proposition}[{\cite[Theorem 3.2]{FerrareseRoggero}}]\label{th:richiami} 
The linear space $Z\big(\mathfrak{L}(J,T)\big)$ can be naturally identified with the Zariski tangent space to $\St{\sigma}(J,T)$ at the origin.

If $C'\subset C$ is any subset of $\dim_{\K} Z\big(\mathfrak{L}(J,T)\big)$ variables such that  $\mathfrak{L}(J,T)\oplus\langle C' \rangle=\langle C\rangle$, then $\mathfrak{h}(J,T)\cap \K[C']$ defines a $\ell$-homogeneous subvariety in $\AA^{\vert C'\vert}$ isomorphic to $\St{\sigma}(\mathfrak{h},T)$.
\end{proposition}
\begin{proof}
By definition there exist $e = \vert C \setminus C'\vert$ $\ell$-homogeneous linear form $l_1,\ldots,l_e$ in $\mathfrak{L}(J,T)$ such that $\{l_1,\ldots,l_e\} \cup C'$ is a basis for the $\K$-vector space of linear forms in $\K[C]$. Then $\mathfrak{h}(J,T)$ has a set of $\ell$-homogeneous generators of the type
\begin{equation}
l_1 + q_1, \ldots, l_e + q_e, q_{e+1}, \ldots, q_{e+s}
\end{equation}
where $q_1,\ldots,q_{e+s} \in \K[C']$ so that the inclusion
\begin{equation}
\K[C']/(q_{e+1},\ldots,q_{e+s}) \hookrightarrow \K[C]/\mathfrak{h}(J,T)
\end{equation}
is indeed an isomorphism (see also \cite[Proposition 2.4]{RoggeroTerracini}). The hypothesis $\vert C'\vert = \dim_{\K} \mathfrak{L}(J,T)$ ensures that $l_1,\ldots,l_e$ generate $\mathfrak{L}(J,T)$ so that $q_{e+1},\ldots,q_{e+s}$ belong to $(C')^2\K[C']$ and the tangent space at the origin of $Z((q_{e+1},\ldots,q_{e+s}))$  is a linear space of dimension $\vert C '\vert$, i.e. $\AA^{\vert C '\vert}$ itself.
\end{proof}

We may summarize the previous result saying that \emph{$\St{\sigma}(J,T)$ can be embedded in its Zariski tangent space at the origin}. This explains the following terminology.

\begin{definition}\label{def:eliminabili} We call \emph{embedding dimension}\index{Gr\"obner stratum!embedding dimension of a}\index{embedding dimension|see{Gr\"obner stratum}} of $\St{\sigma}(J,T)$ the dimension of the affine space $\AA^{\vert C'\vert}$ defined in Proposition \ref{th:richiami} and we will denote it by $\ed \St{\sigma}(J,T)$, i.e. $\gls{ed} = \vert C'\vert$. The complement  $C'' = C\setminus C'$ is a \emph{maximal set of eliminable variables} for $\mathfrak{h}(J,T)$.
\end{definition}

\begin{corollary}\label{th:liscio} In the above notation, the following statements are equivalent:
\begin{enumerate}[(i)]
 \item $\St{\sigma}(J,T)\simeq \AA^{\ed \St{\sigma}(J,T)}$;
 \item $\St{\sigma}(J,T)$ is smooth;
 \item the origin is a smooth point for $\St{\sigma}(J,T)$;
 \item $\ed \St{\sigma}(J,T) \leqslant \dim_\K \St{\sigma}(J,T)$.
\end{enumerate}
\end{corollary}

Note that in general a maximal set of eliminable variables (and so its complementary) is not uniquely determined. However, if $C_{\alpha_i \beta}\in \mathfrak{L}(J,T)$, then $C_{\alpha_i \beta}$  belongs to any set of eliminable variables; on the other hand, if $C_{\alpha_i \beta}$ does not appear in any element of $\mathfrak{L}(J,T)$, then $C_{\alpha_i \beta}$ does not belong to any set of eliminable variables.

There is an easy criterion that allows us to decide if a variable is eliminable or not.
\begin{criterion}\label{criterio} Let us consider two polynomials $F_i$ and $F_j$ among those defined in \eqref{eq:polibasegenerica} such that $\IN_{\sigma}(F_i)= x^{\alpha_i}$, $\IN_{\sigma}(F_j)=x^{\alpha_j}$ and let $C_{\alpha_i\beta}$ be a variable appearing in the tail of $F_i$. Using the reduction with respect to $J$ of a $\ell$-homogeneous polynomial $x^{\gamma}F_i-x^{\delta}F_j$ we can see that:
\begin{enumerate}[(a)]
\item\label{it:criterio_i} if $x^{\gamma+\beta}\notin J$ and $x^{\gamma+\beta -\delta}$ is not a monomial that appears in $F_j$, then $C_{\alpha_i\beta}\in \mathfrak{L}(J,T)$;
\item\label{it:criterio_ii} if $x^{\gamma+\beta}\notin J$ and $x^{\beta'}=x^{\gamma +\beta -\delta}$ is a monomial that appears in $F_j$, then $C_{\alpha_i\beta}-C_{\alpha_j\beta'}\in \mathfrak{L}(J,T)$.
\end{enumerate}
Moreover if $C_{\alpha_i\beta}-C_{\alpha_j\beta'} \in \mathfrak{L}(J,T)$, then every maximal set of eliminable variables must contain at least either one of them.
\end{criterion}

In most cases the number $N = |C|$ is very big and $\mathfrak{h}(J,T)$ needs a lot of generators so that finding it explicitly is a very heavy computation. On the contrary $\mathfrak{L}(J,T)$ is very fast to compute and so we can easily obtain a set of eliminable variables $C''$; a forgoing knowledge of $C'$ allows a simpler computation of the ideal  $\mathfrak{h}(J,T)\cap \K[C\setminus C'']$ that gives $\St{\sigma}(J,T)$ embedded in the affine space of minimal dimension $\AA^{\ed \St{\sigma}(J,T)}$.

Furthermore, in many interesting cases we can greatly bring down the number of involved variables thanks to another kind of argument.

\begin{theorem}\label{th:saturato} Let $J \subset \K[x]$ be a Borel-fixed saturated monomial ideal with basis $G$, $m$ any  integer and $\mathfrak{h}(J_{\geqslant m})$  the ideal  of $\St{\sigma}(J_{\geqslant m})$  as in Definition \ref{def:procedura}. 
\begin{enumerate}[(i)]
\item\label{it:saturato_i}  There is a set of eliminable variables for $\mathfrak{h}(J_{\geqslant m})$ 
that contains all variables except at most those appearing in polynomials $F_i$ whose leading term is either  $x^\alpha \in G_{\geq m}$ or  $x^\alpha x_0^{m-\vert \alpha \vert }$, where $x^\alpha \in G_{<m}$.
\item\label{it:saturato_ii}  $\St{\sigma}(J_{\geqslant m-1})$ is a closed subscheme of $ \St{\sigma}(J_{\geqslant m})$. More precisely   $\St{\sigma}(J_{\geqslant m-1})$ is isomorphic to $\St{\sigma}(J_{\geqslant m},T)$,  where $T$ contains the complete tail of a monomial in the basis of $J_{\geqslant m}$ if it is not divided by $x_0$, and a tail containing only  monomials divided by $x_0$ otherwise.
\item\label{it:saturato_iii}  If  $x_{1}$ does not appear in any monomial of degree  $m+1$ in the monomial basis of $J$, then $\St{\sigma}(J_{\geqslant m-1})\simeq \St{\sigma}(J_{\geqslant m})$. 
\item\label{it:saturato_v}  If  $x_{1}$ appears in $N$ monomials of degree $m+1$ in the monomial basis $G$ of $J$, then $\ed \St{\sigma}(J_{\geqslant m})\geqslant \ed \St{\sigma}(J_{\geqslant m-1}) + NM$, where $M$ is the number of monomials of the basis of $J$ of degree smaller than $m$. 
\item\label{it:saturato_vv} $\St{\sigma}(J_{\geqslant m-1})\not \simeq \St{\sigma}(J_{\geqslant m})$ if and only if $x_{1}$ appears in monomials of degree $m+1$ in the monomial basis of $J$ and $J_{\geqslant m-1}\neq  J_{\geqslant m}$.
\item\label{it:saturato_iv} If $s$ is the maximal degree of a monomial divided by $x_{1}$ in the monomial basis of $J$, then $\St{\sigma}(J_{\geqslant s-1})\simeq \St{\sigma}(J_{\geqslant m})$ for every $m\geqslant s$.
\end{enumerate}
\end{theorem}
\begin{proof} \emph{(\ref{it:saturato_i})} Let us consider any monomial $x^\eta$ in the monomial basis of $J_{\geqslant m}$ which does not belong to $G_{\geqslant m}$ and such it that could be written as $x^\eta=x^\alpha x^\epsilon$ where $x^\alpha \in G_{<m}$  and $x^\epsilon$ is a monomial of degree $m-\vert\alpha\vert$,  $x^\epsilon \neq x_0^{m-\vert\alpha\vert}$. Then among the polynomials $F_i$ there are:
\[
 \begin{split}
&  F = x^\alpha x_0^{m-\vert\alpha\vert} + \sum C_{\beta} x^\beta,\\
&  F' = x^{\alpha+\varepsilon} + \sum C'_{\delta} x^\delta. 
 \end{split}
\]
We have to prove that all the variables $C'$ that appear in $F'$ can be eliminated. The $S$-polynomial of $F$ and $F'$ is:
\[
 S(F,F') = x_0^p F' - x^{\varepsilon'} F = 
\sum C'_{\delta} x^\delta x_0^p - \sum C_{\beta} x^{\beta+\varepsilon'}.
\]
No monomial $x^\delta x_0^p$ in the first summand belongs to $J_{\geqslant m}$ because   $x^\delta \notin J$  and $J$ is saturated and Borel-fixed. Thus, the linear part of the coefficient of $x^\delta x_0^p$ in  the reduction of this $S$-polynomial will be either $C'_{\delta}$ or $C'_{\delta}-C_{\beta}$.  Then $C'$ is a set of eliminable variables for $J_{\geqslant m}$.

\medskip

\emph{(\ref{it:saturato_ii})} The first part of this statement is a special case of general facts proved in \cite[Section 3]{HaimanSturmfels}. We directly prove the  second part (which  implies the first one). Here  we  denote by $x^\alpha$ and  $x^\gamma$ the monomials in the basis of $J_{\geqslant m-1}$ of degree $m-1$ and $\geq m$ respectively, and we set:
\[
\begin{split} F_{\alpha} & {}= x^\alpha  + \sum C_{\alpha\beta} x^\beta, \qquad  \vert\alpha\vert = m-1  \\
  F_{\gamma} & {}= x^\gamma  + \sum C_{\gamma\eta} x^\eta, \qquad \vert\gamma\vert \geqslant m  
  \end{split}
\]
  where  $ x^\beta $ varies among all   monomials of degree $m-1$ in the tail of $ x^\alpha$ and $x^\eta $ among those of the same degree of $ x^\gamma$  in its tail.   Applying the  procedure described in Definition \ref{def:procedura} on such set of polynomials we  define $\St{\sigma}(J_{\geqslant m-1})$  by means of an ideal $\mathfrak{h}(J_{\geqslant m-1}) \subset \K[C]$.
  
The basis of $J_{\geqslant m}$ is made by  monomials of the following three types:
\begin{itemize}
\item monomials $x^\gamma$ of degree $\geqslant m$, that  also belong to the basis of $J_{\geqslant m-1}$;
\item monomials $x^\alpha x_0$ such that  $x^\alpha $ is any monomial of degree $m-1$ in the basis of $J_{\geqslant m-1}$;
\item  monomials $x^\alpha x_i$ of degree $m$ such that $x^\alpha $ is as above and $ \min(x^\alpha)> x_i \neq x_0$.
\end{itemize} 
We set:
 \begin{equation} \label{listaF}
 \begin{split}
 \overline{F}_{\alpha 0} & {}= x^\alpha x_0 + \sum C_{\alpha\delta} x^\delta x_0, \\
 \overline{F}_{\alpha i} & {}= x^\alpha x_i + \sum C'_{\alpha i \tau} x^\tau  \quad \vert \tau\vert= m    \quad x^\tau < x^\alpha x_i, \\
 \overline{F}_{\gamma} & {} = x^\gamma  + \sum C_{\gamma\eta} x^\eta.
 \end{split}
 \end{equation}
  Note that we use  the same  names   for some of the coefficients that appears in polynomials $F$ and $\overline{F}$, so that $\overline{F}_{\alpha 0}=x_0 F_\alpha$ and $\overline{F}_{\gamma}=F_{\gamma}$.
  
 Applying the  procedure described in Definition \ref{def:procedura} on the set of polynomials $\overline{F}$  we obtain an ideal  $\mathfrak{h}' \subset \K[C,C']$ defining $\St{\sigma}(J_{\geqslant m},T)$. 
  
Thanks  to \emph{(\ref{it:saturato_i})} we know that $C'$ is a set of eliminable variables for $\mathfrak{h}'$ and so $\St{\sigma}(J_{\geqslant m},T)$ is also defined by $\overline{\mathfrak{h}}=\mathfrak{h}'\cap \K[C]$. The statement follows once we show that $\overline{\mathfrak{h}} =\mathfrak{h}(J_{\geqslant m-1})$.
  
In order to eliminate the variables $C'$ we consider every monomial  $x^\alpha x_i =\IN (\overline{F}_{\alpha i} )$ and reduce it   using the polynomials $F$. In this way we obtain  a polynomial $H_{\alpha i} \in (F)\K[X,C]$ such that  $x^\alpha x_i+H_{ \alpha i} $ is completely reduced w.r.t. $J$.  Then also  $x^\alpha x_i x_0+H_{ \alpha i}x_0 +\sum C'_{\alpha i \tau} x^\tau x_0   $ (i.e. $\overline{F}_{\alpha i}x_0+H_{ \alpha i}x_0  $) is  reduced modulo $J$ and moreover it belongs to $(\overline{F})\K[X,C,C']$  because $x_0 F\subseteq (\overline{F})\K[X,C,C']$.  The coefficients of the monomials in the variables $x$ belong to $\mathfrak{h}'$, because the ideal $\mathfrak{h}'$  is generated by the coefficient of monomials in $x$ in the polynomials in $(F)\K[X,C,C']$ that are reduced modulo $J_{m-1}$ or modulo $J$, which is the same (Proposition \ref{th:unica}\emph{(\ref{it:unica_1'})} and Theorem \ref{th:cor-unica}). The coefficients of the monomials in $x$ of $F_{\alpha i}x_0+H_{ \alpha i}x_0  $ are also the coefficients of the monomials in $x$ of $F_{\alpha i}+H_{ \alpha i} $, and are precisely the set of polynomials of the type $C'_{ \alpha i \tau}-\phi_{ \alpha i \tau}(C)$ that allow us to   eliminate   the variables $C'$. So  the elimination of   $C'$ is obtained simply  putting $C'_{ \alpha i \tau}=\phi_{ \alpha i \tau}(C)$. In this way $F_{\alpha i}$ becomes $-H_{ \alpha i}$ that belongs to $(G)\K[X,C]$.

  The ideal $\overline{\mathfrak{h}}$, obtained from $\mathfrak{h}'$ eliminating $C'$, can  also be obtained first eliminating $C'$ and after taking the coefficients of the monomials in $x$, because the procedure of eliminating   $C'$ and that of taking coefficients.  So  $\overline{\mathfrak{h}}$ is generated by the coefficients of monomials in $x$ of polynomials in  $(x_0G_\alpha,-H_{ \alpha i}, G_\gamma)\K[X,C]$ that are reduced modulo $J$.  Hence $\overline{\mathfrak{h}}\subseteq \mathfrak{h}(J_{\geqslant m-1})$ because $(x_0G_\alpha,-H_{ \alpha i}, G_\gamma)\K[X,C]\subset (G)\K[X,C]$.
 
  On the other hand, $x_0(G)\K[X,C]=(x_0G_\alpha, x_0G_\gamma)\K[X,C]\subset (x_0G_\alpha, G_\gamma)\K[X,C]$. Moreover two polynomials $Q$ and $x_0Q$ have the same coefficients of the monomials in $x$  and either one is   reduced modulo $J$ if and only the other is. Hence  we obtain the opposite inclusion  $\mathfrak{h}(J_{\geqslant m-1}) \subseteq \overline{\mathfrak{h}}$.
  
  \medskip
  
\emph{(\ref{it:saturato_iii})} We use \emph{(\ref{it:saturato_ii})} and prove that in the present hypothesis, $\St{\sigma}(J_{\geqslant m})\simeq \St{\sigma}(J_{\geqslant m},T)$, where $T$ is defined as in  \emph{(\ref{it:saturato_ii})}.    Following Definition \ref{def:procedura}, we obtain the ideal  $\mathfrak{h}(J_{\geqslant m})$ of $\St{\sigma}(J_{\geqslant m})$  using:
\[
 \begin{split}
 F''_{\alpha 0} & {}= x^\alpha x_0 + \sum C_{\alpha\delta} x^\delta x_0 +\sum C''_{\alpha\mu} x^\mu \quad, \quad x_0 \nmid x^\mu\\
 F_{\alpha i} & {}= x^\alpha x_i + \sum C'_{\alpha i \tau} x^\tau   \\
 F_{\gamma} & {} = x^\gamma  + \sum C_{\gamma\eta} x^\eta= G_\gamma.
 \end{split}
 \]
 Note that $ F_{\alpha i}$ and $ F_{\gamma}$ are as in  \emph{(\ref{it:saturato_ii})}, but all the degree $m$ monomials of  the  tail of $x_0x^\alpha$ appear in $F''_{\alpha 0}$, and not only those  divided by $x_0$.
  
For every monomial $x^\alpha$ of degree $m-1$ in the basis of $J_{m-1}$, let us consider  the $S$-polynomial:
\[
S(F''_{\alpha 0} ,F_{\alpha 1})= \sum C_{\alpha\delta} x^\delta x_{1}x_0 +\sum C''_{\alpha\mu} x^\mu x_{1} - \sum C'_{\alpha i \tau} x^\tau x_0.
\]
 By hypothesis no monomial appearing in it belongs to $J_{m}$. In fact  $x^\mu x_1 \in J$ if and only if it is a minimal generator of $J$, which is excluded by hypothesis because its degree is $m+1$, or it is of the type $x^\alpha x_a$ with $x^\alpha$ minimal generator of $J_m$  and  $x_a=\min (x^\mu x_1)=x_{1}$, while $x^\mu \notin J_m$. Then $S(F''_{\alpha 0} ,F_{\alpha 1})$ is already reduced with respect to $J_m$ and so the coefficients of the monomials in $x$  belong to $\mathfrak{h}(J_{\geqslant m})$. Especially, as both $x^\delta x_1x_0$ and   $x^\tau x_0$ are multiple of $x_0$, while $x^\mu x_1$ is not, the coefficient of $x^\mu x_1$ is simply  $C''_{\alpha\sigma} $ so that each $C''_{\alpha\sigma} $ belongs to $\mathfrak{h}(J_{\geqslant m})$.  Hence we can eliminate all the variables $C''$, just putting them equal to 0. In this way    $F''_{\alpha 0}$ becomes $\overline{F}_{\alpha 0}$ as in \eqref{listaF} and $\St{\sigma}(J_{\geqslant m})\simeq \St{\sigma} (J_{\geqslant m},T)$, where $T$ is as in \emph{(\ref{it:saturato_ii})}, and we conclude  because $\St{\sigma} (J_{\geqslant m},T) \simeq \St{\sigma} (J_{\geqslant m-1}) $. 
 
 \medskip
 
\emph{(\ref{it:saturato_v})} By \emph{(\ref{it:saturato_ii})}, we know that $\ed \St{\sigma}(J_{\geqslant m}) \geqslant \ed \St{\sigma}(J_{\geqslant m},T) = \St{\sigma}(J_{\geqslant m-1})$, where the tails defined in $T$ contain only monomials divided by $x_0$. Let us now consider a monomial $x^\alpha$ among the generators of $J$ of degree smaller than $m$ and a generator $x^\gamma$ of degree $m+1$ divided by $x_1$. Computing the stratum $\St{\sigma}(J_{\geqslant m})$, in the tail of $x^\alpha x_0^{m-\vert\alpha\vert}$ there is the monomial $x^\beta = x^\gamma/x_1$ not belonging to $T$. Let us call $D$ the coefficient of $x^\beta$, that is
\[
F = x^\alpha x_0^{m-\vert\alpha\vert} + \ldots	 + D x^\beta + \ldots.
\]
Thinking about the Eliahou-Kervaire syzygies of the ideal $J$, it is easy to see that in any $S$-polynomial, $F$ is surely multiplied by a monomial $x^\delta,\ \min x^\delta > \min  x^\alpha x_0^{m-\vert\alpha\vert} = 0$. Therefore in every $S$-polynomial the monomial $x^\beta x^\delta = (x^\beta x_i) x^{\delta'}$ belongs to $J$ because of the Borel-fixed hypothesis, so that it can be reduced. Finally there is no equation involving the variable $D$, so it is free and it cannot be eliminated. Repeating the reasoning for the $M$ minimal generators of degree smaller than $m+1$ and for the $N$ generators divided by $x_1$ of degree $m+1$, we obtain the thesis.
  
  \medskip
  
\emph{(\ref{it:saturato_vv})} and \emph{(\ref{it:saturato_iv})} straightforward applying \emph{(\ref{it:saturato_v})} and \emph{(\ref{it:saturato_iii})}.
\end{proof}

With the following examples, we want to underline again the not so crucial role played by term ordering in this construction (Example \ref{ex:strato1}) and we want to show (Example \ref{ex:strato2} and Example \ref{rk:stratoLex}) that the estimate of growth of the embedding dimension of the stratum introduced in Theorem \ref{th:saturato}\emph{(\ref{it:saturato_v})} is a lower bound.

\begin{example}\label{ex:strato1}
Let us consider the ideals $I = (x_3,x_2^2,x_2 x_1)$ and $I_{\geqslant 2} = (x_3^2,x_3 x_2,x_3 x_1,$ $x_3 x_0,x_2^2,x_2 x_1)$ in the ring $\K[x_0,x_1,x_2,x_3]$ and the Gr\"obner strata of the ideal $I_{\geqslant 2}$ according to two different term orderings: $\St{\DegLex}(I_{\geqslant 2})$ and $\St{\RevLex}(I_{\geqslant 2})$. In the first case there are $24$ monomials in the complete tails, i.e. 24 new variables $C$, whereas in the second case they are $23$, so we may guess that the family of the ideals with initial ideal $I_{\geqslant 2}$ w.r.t. $\DegLex$ could be different from the family of the ideals with initial ideal $I_{\geqslant 2}$ w.r.t. $\RevLex$.

However applying Theorem \ref{th:saturato}, we can see that $\St{\DegLex}(I_{\geqslant 2}) \simeq \St{\DegLex}(I)$ and $\St{\RevLex}(I_{\geqslant 2}) \simeq \St{\RevLex}(I)$. Now the tails of the 3 monomials that generate $I$ are the same w.r.t. both term orders and then (see Remark \ref{rk:termordersignificance})
\[
 \St{\DegLex}(I_{\geqslant 2}) \simeq \St{\DegLex}(I) = \St{\RevLex}(I) \simeq \St{\RevLex}(I_{\geqslant 2}).
\]
\end{example}

\begin{example}\label{ex:strato2}
Let us consider the polynomial ring $\K[x_0,x_1,x_2,x_3]$, the ideal $J =(x_3^2,x_3 x_2,x_3 x_1^4,x_2^7,x_2^6 x_1^2)$ and the term ordering $\omega$ given by the matrix \eqref{eq:termOrdering} with second row equal to $(23,5,2,1)$. By the previous theorem we know that
\[
\begin{split}
& \St{\omega}(J) \simeq \St{\omega}(J_{\geqslant 3}),\\
& \St{\omega}(J_{\geqslant 4}) \simeq \St{\omega}(J_{\geqslant 5}) \simeq \St{\omega}(J_{\geqslant 6}),\\
& \St{\omega}(J_{\geqslant 7}) \simeq \St{\omega}(J_{\geqslant m}),\ \forall\ m \geqslant 8
\end{split}
\]
and
\[
\ed \St{\omega}(J_{\geqslant 4}) \geqslant \ed \St{\omega}(J) + 2 \qquad \ed \St{\omega}(J_{\geqslant 7}) \geqslant \ed \St{\omega}(J_{\geqslant 4}) + 3
\]
By an explicit computation, $\ed \St{\omega}(J) = 46$,  $\ed \St{\omega}(J_{\geqslant 4}) = 50$ and  $\ed \St{\omega}(J_{\geqslant 7}) = 56$.
\end{example}

\begin{example}\label{rk:stratoLex}\index{lexicographic ideal}
For any term ordering $\sigma$, there are at most two possible classes of isomorphism for the strata $\St{\sigma}(L_{\geqslant m})$, where $L \subset \K[x]$ is a saturated lexicographic ideal: $\St{\sigma}(L)$ and $\St{\sigma}(L_{\geqslant r-1})$, where $r$ is the maximal degree of a minimal generator, in fact the variable $x_1$ appears (if it does) only in the generator of degree $r$. Called $b$ the number of generators of degree $r$, applying Theorem \ref{th:saturato}\emph{(\ref{it:saturato_v})}, we have
\[
\ed \St{\sigma}(L_{\geqslant r-1}) \geqslant \ed \St{\sigma}(L) + n-b.
\]

If the monomial of maximal degree in the basis does not contain the variable $x_1$, we have $\St{\sigma}(L_{\geqslant m}) \simeq \St{\sigma}(L),\ \forall\ m$. 
\end{example}

We conclude this section with a result similar to the one stated in Theorem \ref{th:saturato} that concerns only the case of Gr\"obner strata w.r.t. $\RevLex$.
\begin{proposition}
Let $J$ be a Borel-fixed saturated ideal and let us consider the $\RevLex$ term ordering. Then
\begin{equation*}
 \St{\RevLex}(J) \simeq \St{\RevLex}(J_{\geqslant m}),\quad \forall\ m.
\end{equation*} 
\end{proposition}
\begin{proof}
The arguments to achieve the proof are very similar to the arguments used in the proof of Theorem \ref{th:saturato}. First of all let us consider the monomials 
\[
 F_\alpha = x^\alpha + \sum C_{\alpha\beta} x^\beta
\]
corresponding to the monomial basis $G$ of $J$ and the ideal $\mathfrak{h}(J) \subset \K[C]$ of the stratum $\St{\RevLex}(J)$. 

In order to compute $\St{\RevLex}(J_{\geqslant m})$, we have to consider again polynomials $F_\alpha$ as before if $x^\alpha \in G_{\geqslant m}$ and new polynomials $\overline{F}_{\alpha\varepsilon}$ such that $\IN(\overline{F}_{\alpha\varepsilon}) = x^{\alpha+\varepsilon}$, $\forall\ x^\alpha \in G_{< m}$ and $\forall\ x^\varepsilon$ of degree $m-\vert\alpha\vert$, especially $x^\alpha x_0^{m-\vert\alpha\vert}$.
Then by the definition itself of $\RevLex$, the tail of $x^\alpha x_0^{m-\vert\alpha\vert}$ contains exactly the monomials in the tail of $x^\alpha$ multiplied by $x_0^{m-\vert\alpha\vert}$. So we can write
\[
\overline{F}_{\alpha\varepsilon} = \begin{cases}
                         x^{\alpha+\varepsilon} + \sum E_{\alpha\delta}^\varepsilon x^\delta,& \forall\ x^\varepsilon \neq x_0^{m-\vert\alpha\vert},\\
                         x^\alpha x_0^{m-\vert\alpha\vert} + \sum C_{\alpha\beta} x^\beta x_0^{m-\vert\alpha\vert} = x_0^{m-\vert\alpha\vert} F_\alpha,&
                         \text{if } x^\varepsilon = x_0^{m-\vert\alpha\vert}\\
                        \end{cases}
\]
hence $\mathfrak{h}(J_{\geqslant m}) \subset \K[C,E]$ (note that in the present case variables $D$ do not appear by construction).

By Theorem \ref{th:saturato}\emph{(\ref{it:saturato_i})}, we know that all the variables $E$ can be eliminated. By the same reasoning used in the proof of Theorem \ref{th:saturato}\emph{(\ref{it:saturato_ii})}, the ideal $\overline{\mathfrak{h}} = \mathfrak{h}(J_{\geqslant m}) \cap \K[C]$ contains the coefficients of the monomials in $x$ of a set of S-polynomials corresponding to a set of the $S$-polynomials of the monomial basis of $J$: so $\St{\RevLex}(J) \simeq \St{\RevLex}(J_{\geqslant m})$.
\end{proof}

\section{Open subsets of the Hilbert scheme I}\label{sec:openSubsetsI}\index{Hilbert scheme}

We now discuss the relation between Gr\"obner strata and Hilbert schemes, so in this section we will use again the main related notation. Given a Hilbert polynomial $p(t)$, $r$ will denote its Gotzmann number, $N(t) = \binom{n+t}{n}$ and $q(t) = N(t) - p(t)$.

\begin{lemma}\label{th:sett} Let $J \subset \K[x]$ be any monomial ideal defining a subscheme $\Proj \K[x]/J \subset \PP^n$ with Hilbert polynomial $p(t)$ and let $\sigma$ be any term ordering. Then (at least set-theoretically)  $\St{\sigma}(J) \subseteq \Hilb{n}{p(t)}$. 
\end{lemma}
\begin{proof} Let $I$ be any ideal in $\St{\sigma}(J)$. By hypothesis $\IN_{\sigma}(I)=J$ and then $I$ and $J$ share the same Hilbert function. Therefore $I$ defines a subscheme with Hilbert polynomial $p(t)$, i.e. a point of $\Hilb{n}{p(t)}$.
\end{proof}

Now we will see that with some restriction the set-theoretic inclusions are in fact algebraic maps and that for some special ideals they are indeed open injections. The crucial point is that the stratum structure (and so its injection in the Hilbert scheme) depends on the ideal $J$ and not on the the corresponding subscheme $Z=\Proj \K[x]/J$. 
For instance the stratum of the saturated lexicographic ideal $L$ associated to the Hilbert polynomial $p(t)$ is not in general isomorphic to an open subset of $\Hilb{n}{p(t)}$ (see \cite{RoggeroTerracini} and Example \ref{rk:stratoLex}), whereas, as we will see,  the stratum of its truncation $L_{\geqslant r}$ is an open subset of the Reeves-Stillman component of $\Hilb{n}{p(t)}$.

\bigskip

Let $J = J_{\geqslant r}$ be a monomial ideal with Hilbert polynomial $p(t)$ (and Gotzmann number $r$) and let $\sigma$ be any term ordering. As seen in the previous section,  every ideal $I$ such that $\IN_{\sigma}(I)=J$ has a (unique) reduced Gr\"obner basis $\{f_1, \dots, f_{q(r)}\}$ where $f_i$ is as in Definition \ref{def:codaridotta}. Not every ideal generated by $q(r)$ polynomials of such a type has $J$ as initial ideal. In order to obtain equations for $\St{\sigma}(J)$ we consider again the coefficients $c_{\alpha_i \beta}$ appearing in the $f_i$ as new variables; more precisely let  $C=\{ C_{\alpha_i \beta},\ i=1, \dots, q(r), \ x^{\beta}\notin J_r \text{ and } x^{\alpha_i} >_{\sigma} x^{\beta} \}$ be new variables and consider $q(r)$ polynomials in $\K[C][x]$ of the following type:
\begin{equation}\label{eq:costF}
  F_i = x^{\alpha_i} + \sum_{x^{\beta}\in \tail{x^{\alpha_i}}{\sigma}{J}}  {C}_{\alpha_i \beta} x^\beta.  
\end{equation}
We obtain the ideal $\mathfrak{h}(J)$ of $\St{\sigma}(J)$ collecting the coefficients (polynomials in $\K[C]$) of the monomials in $x$ of some complete reduction with respect to $F_1, \dots, F_{q(r)}$ of all the $S$-polynomials \glshyperlink[$S(F_i,F_j)$]{Spoly}, corresponding to a set of generators for \glshyperlink[$\syz(J)$]{syz} (see Theorem \ref{th:cor-unica} and Proposition \ref{th:unica}\emph{(\ref{it:unica_4})}).

\begin{proposition}\label{th:matriceA} 
In the above notation, let $J=J_{\geqslant r}$ be a monomial ideal with Hilbert polynomial $p(t)$ and let $\mathcal{A}$ be the $q(r)(n+1)\times N(r+1)$ matrix whose columns correspond to the monomials in $\K[x]_{r+1}$ and whose rows correspond to the polynomials $x_jF_i$, for all $j = 0,\ldots,n$ and $i = 1,\dots,q(r)$ (and the entries are polynomials in $\K[C]$ corresponding to the coefficients of the monomials in $x$). Then the ideal $\mathfrak{h}(J)$ of the Gr\"obner stratum $\St{\sigma}(J)$ is generated by the minors of $\mathcal{A}$ of dimension $q(r+1)+1$.
\end{proposition}
\begin{proof} By abuse of notation we write in the same way a polynomial and the corresponding row in the matrix $\mathcal{A}$. It is quite evident by elementary arguments of linear algebra, that the ideal $\mathfrak{a}\subseteq \K[C]$, generated by all minors of dimension $q(r+1)+1$ does not change if we perform some row reduction on $\mathcal{A}$. Let $\mathcal{G}$ be a set of $q(r+1)$ rows whose leading terms are a basis of $J_{r+1}$. If $x_hF_i \notin \mathcal{G}$, then it has the same leading term than one in $\mathcal{G}$, say $x_kF_j$; we can substitute $x_hF_i$ with $x_hF_i-x_kF_j$. In this way the rows not in $\mathcal{G}$ become precisely all the $S$-polynomials $S(F_i,F_j)$ that have $r+1$ w.r.t. the variables $x$.

At the end of this sequence of row reductions, we can write the matrix as follows:
\begin{equation}\label{eq:matriceA1lineare}
\left(\begin{array}{c|c}
D & E \\
\hline
S & L
\end{array}\right)
\end{equation}
where $D$ is a $q(r+1) \times q(r+1)$ upper triangular matrix with 1's along the main diagonal, whose rows correspond to $\mathcal{G}$ and whose columns correspond to monomials in $J_{r+1}$.

Using rows in $\mathcal{G}$, we now perform a sequence of rows reductions on the following ones, in order to annihilate all the coefficients of monomials in $J_{r+1}$, that is the entries of the submatrix $S$: if $a(C)$ is the first non-zero entry in a row not in $\mathcal{G}$ and its column corresponds to the monomial $x^{\gamma} \in J_{r+1}$, we add to this row $-a(C)x_kF_j$, where $x_kF_j \in \mathcal{G}$ and $\IN_{\sigma}(x_kF_j)=x^{\gamma}$. This is nothing else than a step of reduction with respect to $\{F_1, \dots, F_{q(r)}\}$. At the end of this second turn of rows reductions, we can write the matrix as follows:
\begin{equation}\label{eq:matriceA1strato}
\left( \begin{array}{c|c}
D & E \\
\hline
0 & R 
\end{array}\right) 
\end{equation} 
where the rows in $(D\ |\ E)$ are unchanged whereas the rows in $(0\ |\ R)$ are the coefficients of the monomials in $x$ of the complete reductions of $S$-polynomials of degree $r+1$ w.r.t. variables $x$. Then $\mathfrak{a}$ is generated by the entries of $R$ and so $\mathfrak{a} \subset \mathfrak{h}(J)$.

We can see that this inclusion is in fact an equality taking in mind Theorem \ref{th:EKsyzygies} and Proposition \ref{th:unica}\emph{(\ref{it:unica_4})}: the first one says that $\syz(J)$ is generated in degree $r+1$ and the second one that in this case $\mathfrak{h}(J)$ is generated by the coefficients of the monomials in $x$ of complete reductions of the $S$-polynomials $S(F_i,F_j)$ of degree $r+1$ w.r.t. the variables $x$.
\end{proof}

The following corollary just express in an explicit way two properties contained in the proof of Proposition \ref{th:matriceA}.

\begin{corollary}\label{cor:matriceA} In the above notation:
\begin{itemize}
\item the ideal $\mathfrak{h}(J)$ is generated by the entries of the submatrix $R$ in \eqref{eq:matriceA1strato};
\item the vector space $\mathfrak{L}(J)$ is generated by the entries of the submatrix $L$ in \eqref{eq:matriceA1lineare}.
\end{itemize}
\end{corollary}

\bigskip

Let us now look briefly at the construction of the Hilbert scheme from a local perspective, that is finding equations of the open subset $\mathcal{U}_{\Delta_{\II}} \cap \Hilb{n}{p(t)}$ for any open affine subset\index{Hilbert scheme!open subset of the} $\mathcal{U}_{\Delta_\II}$ of the Grassmannian $\Grass{q(r)}{N(r)}{\K}$\index{Grassmannian} given by the non-vanishing of a Pl\"ucker coordinate $\Delta_\II$. 
Since the base vector space is that one spanned by the monomials of degree $r$ in $\K[x]$, we can associate to any Pl\"ucker coordinate a monomial ideal $J$ generated by $q(r)$ monomials of degree $r$. Therefore we will denote by $\mathcal{U}_{J}$ and $\mathcal{H}_{J}$ respectively the open subsets of $\Grass{q(r)}{N(r)}{\K}$ and of $\Hilb{n}{p(t)}$ corresponding to the the Pl\"{u}cker coordinate associated to the ideal $J$. 

In a natural way $\mathcal{U}_{J}$ is isomorphic to the affine space  $\AA^{q(r)\big(N(r)-q(r)\big)}$. Indeed, if $J=(x^{\alpha_1}, \dots, x^{\alpha_{q(r)}})$, every point in $\mathcal{U}_{J}$ is uniquely identified by the reduced, ordered set of generators $\langle g_1, \dots,g_{q(r)} \rangle $  of the type  $g_i = x^{\alpha_i}+ \sum c_{\alpha_1 \beta} x^\beta$, where  $c_{\alpha_i \beta}\in \K$ and $x^{\beta}$ is any monomial in $ \K[x]_r\setminus J$. Then we consider on $\AA^{q(r)\big(N(r)-q(r)\big)}$ the coordinates $C_{\alpha_i \beta}$. Note that each $C_{\alpha_i \beta}$  naturally corresponds to the  Pl\"{u}cker coordinate associated to the ideal $J'=( x^{\alpha_1}, \dots,x^{\alpha_{i-1}},x^{\beta} ,x^{\alpha_{i+1}},\dots,  x^{\alpha_{q(r)}})$ (but of course not all the Pl\"ucker coordinates are of this type).

Now we can mimic the construction of Gr\"obner strata and obtain the  defining  ideal  of $\mathcal{H}_J$ as a subscheme of  $\AA^{q(r)\big(N(r)-q(r)\big)}$. Let us consider the set of variables $\overline{C}=\{ C_{\alpha_i \beta}, i=1, \dots, q(r), \ x^{\beta}\in \K[x]_r \setminus J\}$ and $q(r)$ polynomials $\overline{F}_1, \dots, \overline{F}_{q(r)}$ in $\K[\overline{C}][x]$ of the type:
\begin{equation}\label{eq:costG} 
\overline{F}_i = x^{\alpha_i}+\sum_{\begin{subarray}{c}\vert\beta\vert = r\\ x^\beta \notin J\end{subarray}} C_{\alpha_i \beta} x^\beta   
\end{equation}
and let $\overline{\mathcal{A}}$ be the $(n+1)q(r)\times N(r+1)$ matrix representing the polynomials $x_j\overline{F}_i$ (as the matrix $\mathcal{A}$ represents the polynomials $x_j F_i$ in the proof of Proposition \ref{th:matriceA}). Then consider the ideal $\overline{\mathfrak{a}}(J) \subset \K[\overline{C}]$  generated by the minors of $\overline{\mathcal{A}}$ of dimension $q(r+1)+1$. 

\begin{proposition}\label{th:rkB}  
$\overline{\mathfrak{a}}(J)$ is the ideal of $\mathcal{H}_J$ as a closed subscheme of $\AA^{q(r)p(r)}$.
\end{proposition}
\begin{proof}
Every ideal $I \in \mathcal{U}_J$ can be obtained from $(\overline{F}_1, \dots, \overline{F}_{q(r)})$ specializing (in a unique way) the variables $C_{\alpha_i \beta}$ to $c_{\alpha_1 \beta}\in \K$. Obviously not all the specializations give ideals $I \in \mathcal{H}_J$, that is  with Hilbert polynomial $p(t)$, because we have to ask both $\dim_\K I_r =q(r)$ and $\dim_\K I_{r+1} = q(r+1)$: thanks to Gotzmann's Persistence Theorem we know that these two necessary conditions are also sufficient. 

In the open subset $\mathcal{U}_{J}$ the first condition always holds and the rank of every specialization of $\overline{\mathcal{A}}$ is $\geqslant q(r+1)$ by Macaulay's Estimate on the Growth of Ideals.\index{Macaulay's Estimate on the Growth of Ideals} Therefore $\mathcal{H}_J$ is given by the condition $\rank \overline{\mathcal{A}} \leqslant q(r+1)$.
\end{proof}

Let us know consider a special ordering of Pl\"ucker coordinates. We write the $q(r)$ monomials generating the ideal associated to the Pl\"ucker coordinate in decreasing order with respect to $\sigma$; then given $J_1=(x^{\alpha_1}>_\sigma \dots >_\sigma x^{\alpha_{q(r)}})$ and $J_2=(x^{\gamma_1}>_{\sigma}  \dots >_{\sigma} x^{\gamma_{q(r)}})$,  then   $J_1 > J_2$ if $x^{\alpha_i} = x^{\gamma_i}$ for every $i = 1,\ldots,s-1$ and $x^{\alpha_s}>_{\sigma} x^{\gamma_s}$ (i.e. a lexicographic order on the \lq\lq alphabet\rq\rq of the monomials of degree $r$ ordered w.r.t. $\sigma$).  

It is now easy to compare, for the same monomial ideal $J = J_{\geqslant r}$ with Hilbert polynomial $p(t)$, the Gr\"obner stratum $\St{\sigma}(J)$ and the open subset $\mathcal{H}_{J}$.  We underline that for our purpose it will be sufficient to consider the open subsets $\mathcal{H}_{J}$ corresponding to monomial ideals  $J$ defining points of $\Hilb{n}{p(t)}$,  because (scheme-theoretically) they  cover  $\Hilb{n}{p(t)}$. Indeed, if $I$ has Hilbert polynomial $p(t)$, also $\IN_{\sigma}(I)$ does and so $I \in \mathcal{H}_{\IN_{\sigma}(I)}$.

\begin{theorem}\label{th:aperto} Let $p(t)$ be any admissible Hilbert polynomial in $\PP^n$ with Gotzmann number $r$ and let $\sigma$ be any term ordering. 
\begin{enumerate}[(i)]
\item\label{it:aperto_i}  If $J=J_{\geqslant r}$ is a monomial ideal with Hilbert polynomial $p(t)$, then $\St{\sigma}(J)$ is naturally isomorphic to the locally closed subscheme of $\Hilb{n}{p(t)}$ given by the conditions that the Pl\"{u}cker coordinate corresponding to $J$ does not vanish and the preceding ones vanish.
\item\label{it:aperto_ii} For every isolated, irreducible component $H$ of $\Hilb{n}{p(t)}$, there is a monomial ideal $J= J_{\geqslant r}$ with Hilbert polynomial $p(t)$ such that an irreducible component of \glslink{supp}{$\Supp \St{\sigma}(J)$} is an open subset of $\Supp H$. Then $\Supp H$ has an open subset which is a homogeneous affine variety with respect to a non-standard positive grading.
\item\label{it:aperto_iii} Every smooth irreducible component $H$ of $\Hilb{n}{p(t)}$ is rational. 
The same holds for every smooth, irreducible component of $\Supp \Hilb{n}{p(t)}$.
\end{enumerate}
\end{theorem}
\begin{proof} 
\emph{(\ref{it:aperto_i})} We obtain the two affine varieties $\St{\sigma}(J)$ and $\mathcal{H}_{J}$ in a quite similar way. The only difference comes from the definition of the set of polynomials $F_1, \dots, F_{q(r)}$ given in \eqref{eq:costF}, leading to equations for $\St{\sigma}(J)$, and the set of polynomials $\overline{F}_1, \dots, \overline{F}_{q(r)}$ given in \eqref{eq:costG}, leading to equations for $\mathcal{H}_{J}$: in $\overline{F}_i$ the sum is over all the degree $r$ monomials $x^{\beta}\notin J$ whereas in $F_i$ we also assume that $x^{\beta} <_\sigma \IN_{\sigma}(F_i) = x^{\alpha_i}$. Therefore we can think of $\St{\sigma}(J)$ as the affine subscheme defined by the ideal $\overline{\mathfrak{h}}(J)$ in the ring $\K[\overline{C}][X]$, where $\overline{C}=\{C_{\alpha_i \beta} \ \vert\ i=0, \dots, q(r),\ x^{\beta} \in \K[x]_r \setminus J \}$ generated by $\mathfrak{h}(J)$ and by $\big(C_{\alpha_i \beta}\ \vert\  x^{\beta} >_\sigma \IN_{\sigma}(F_i) \big)$, namely  $\overline{\mathfrak{h}}(J)=\mathfrak{h}(J)\K[\overline{C}]+( \overline{C}\setminus C)$. Now we can  conclude because  the Pl\"{u}cker coordinates before that associated to $J$ vanish if and only if all the $C_{\alpha_i \beta}$ such that $x^{\beta} >_\sigma \IN_{\sigma}(F_i)$ vanish. 

\emph{(\ref{it:aperto_ii})}\hfill As\hfill $J$\hfill varies\hfill among\hfill the\hfill finite\hfill set\hfill of\hfill the\hfill monomial\hfill ideals\hfill in\hfill $\Hilb{n}{p(t)}$,\hfill the\\ Gr\"obner strata $\St{\sigma}(J)$ give a set theoretical covering of $\Hilb{n}{p(t)}$ by locally closed subschemes. Then there is a suitable ideal $J$ such that an irreducible component of $\Supp \St{\sigma}(J)$ is an open subset of $H$.
The support and the irreducible components of the support of $\St{\sigma}(J)$ are homogeneous (see \cite[Section IV.3.3]{BourbakiCA} and \cite[Corollary 2.7]{FerrareseRoggero}), having $\St{\sigma}(J)$ a structure of homogeneous affine scheme with respect to a non-standard positive grading $\ell$.

\emph{(\ref{it:aperto_iii})} If $H$ is a smooth, irreducible component of either $\Hilb{n}{p(t)}$ or $\Supp \Hilb{n}{p(t)}$, then it is also reduced. Thanks to the previous item we know that an open subset of $H$ is an affine homogeneous variety with respect to a positive grading. Moreover this open subset is also smooth and so  it is isomorphic to an affine space, by Corollary \ref{th:liscio}. 
\end{proof}

\begin{remark} Let $J$ be a monomial ideal in $\Hilb{n}{p(t)}$ and let $\mathfrak{a}(J) \subset \K[\overline{C}]$ the ideal of $\mathcal{H}_{J}$. It is possible to define a grading $\ell'$ on $\K[\overline{C}]$ such that $\mathfrak{a}(J)$ becomes homogeneous, by the analogous definition: $\ell'(C_{\alpha_i \beta})=\frac{x^{\alpha_i}}{x^\beta}$ if $C_{\alpha_i \beta}$ appears in $\overline{F}_i$ \eqref{eq:costG}. However  this grading $\ell'$ is not necessarily   positive and so it gives less interesting consequences.

If an irreducible  component $H$ of $\Hilb{n}{p(t)}$ is also reduced, Theorem \ref{th:aperto} insures that there is an open subset of $H$ which has the structure of homogeneous variety with respect to a positive grading induced  from that of a suitable Gr\"obner stratum $\St{\sigma}(J)$. On the other hand, in the case of  a non-reduced component we only know that the support of a suitable open subset  is   homogeneous with respect to a positive grading, but this does not imply that the open subset itself is homogeneous. 
\end{remark}

The proof of Theorem \ref{th:aperto} suggests that whenever there do not exist monomials $x^\beta >_{\sigma} \IN_{\sigma}(F_i)$ not belonging to the ideal $J$, the constructions of $\St{\sigma}(J)$ and $\mathcal{H}_J$ are substantially equal, hence the hilb-segment ideals\index{segment ideal!hilb-segment ideal} introduced in Definition \ref{def:segments} assume a great importance.

\begin{corollary}\label{cor:vuoto} Let $p(t)$ be a Hilbert polynomial with Gotzmann number $r$ and let $J = J_{\geqslant r}$ be a hilb-segment ideal (and so Borel-fixed) such that $J_r$ defines a point in the Grassmannian $\Grass{q(r)}{N(r)}{\K}$. If $\K[x]/J$ has Hilbert polynomial different from $p(t)$, then the open subset $\mathcal{H}_J$ of $\Hilb{n}{p(t)}$ is empty.
\end{corollary}
\begin{proof} 
Called  $\sigma$ the term ordering realizing $J$ as hilb-segment ideal, any point $I \in \mathcal{H}_{J}$ should belong to the Gr\"obner stratum $\St{\sigma}(J)$, that is it should share the same Hilbert polynomial of $J$, which is not $p(t)$.
\end{proof}

The first of the following examples highlights both that Theorem \ref{th:aperto} does not hold for a monomial ideal $J$ defining a point of $\Grass{q(r)}{N(r)}{\K}$ but not of $\Hilb{n}{p(t)}$ and that Corollary \ref{cor:vuoto} does not hold  for a monomial ideal $J$ in $\Grass{q(r)}{N(r)}{\K}$ which is not a hilb-segment ideal. Moreover Example \ref{ex:got} presents a concrete case of empty $\mathcal{H}_{J}$ as discussed in the previous corollary.

\begin{example}\label{ex:wrongStratumHilb2P2}
Let us consider the constant Hilbert polynomial $p(t)=2$ in $\PP^2$. \gls{Hilb2P2} is irreducible of dimension 4 (see \cite{IarrobinoRed}). The monomial ideal $J =(x_2^2,x_2x_1,x_1^2,$ $x_0^2)$ is generated by $4$ monomials of degree $2$, but does not belong to $\Hilb{2}{2}$ because its radical is the irrelevant maximal ideal. However, $\mathcal{H}_J$ is non-empty because it contains for instance all the reduced subschemes given by  couples  of points $[1:a:b],\, [1:a':b']\in \PP^2$ such that $ab'\neq a'b$. By the way, $\St{\sigma}(J)$ cannot have any common point with $\Hilb{2}{2}$.
\end{example}

\begin{example}\label{ex:got}  Let us consider the Hilbert polynomial $p(t) = 4t$ and the projective space $\PP^3$. The Gotzmann number is $6$, so that the Hilbert scheme \gls{Hilb4tP3} is embedded in the Grassmannian $\Grass{60}{84}{\K}$. The ideal $J$ generated by the greatest 60 monomials in $\K[x_0,x_1,x_2,x_3]_6$ w.r.t. $\RevLex$ defines by construction a segment $\{J_r\} \subset \pos{3}{6}$ and has constant Hilbert polynomial $\overline{p}(t) = 24$. Thus $J$ defines a point of $\Grass{60}{84}{\K}$ not belonging to $\Hilb{3}{4t}$, therefore $\mathcal{H}_{J}$ is empty.
\end{example}

\begin{corollary}\label{th:segment01} Let $p(t)$ be an admissible Hilbert polynomial in $\PP^n$ and let $H$ be an isolated and irreducible component of $\Hilb{n}{p(t)}$.  
If $H$ contains a point defined by a hilb-segment ideal $J$ w.r.t. some term ordering $\sigma$, then $\St{\sigma}(J)$ is an open subset of $H$, so that $H$ has an open subset  which is an homogeneous affine variety with respect to a non-standard positive grading.
\end{corollary}
\begin{proof} If $J$ is a hilb-segment ideal, then there are no Pl\"{u}cker coordinates preceding that one associated to $J$. Thus $\St{\sigma}(J)\simeq \mathcal{H}_J$ (see Theorem \ref{th:aperto}) and so $\mathcal{H}_J$ is an affine  homogeneous scheme with respect to a positive grading.   
\end{proof}

\begin{corollary}\label{th:segment}
Let $J \subset \K[x]$ be a hilb-segment ideal w.r.t. some term ordering $\sigma$ defining a point of the Hilbert scheme $\Hilb{n}{p(t)}$ and let $H$ be an irreducible component of $\Hilb{n}{p(t)}$ containing the point defined by $J$. If either of the following condition holds: 
\begin{enumerate}
\item  $\St{\sigma}(J)$ is an affine space,
\item  $J$ is a smooth point of $\St{\sigma}(J)$,
\item  $J$ is a smooth point of $\Hilb{n}{p(t)}$,
\end{enumerate}
then $H$ is rational.
\end{corollary}
\begin{proof}
Straightforward consequence of the previous result and of Corollary \ref{th:liscio}.
\end{proof}

\subsection{The Reeves and Stillman component of $\Hilb{n}{p(t)}$} \label{sec:lexsegment} 

A nice application of the results just proved concerns the component of the Hilbert scheme $\Hilb{n}{p(t)}$ containing the point defined by the lexicographic ideal $L \subset \K[x]$ associated to the Hilbert polynomial $p(t)$.\index{lexicographic ideal} This component is unique and it is usually denoted by $H_{RS}$ and called Reeves ans Stillman component of $\Hilb{n}{p(t)}$ because in \cite{ReevesStillman} they prove that the point of $\Hilb{n}{p(t)}$ corresponding to $\Proj \K[x]/L$ (the lexicographic point) is smooth.\index{lexicographic point}

Putting together the smoothness of the lexicographic point and Corollary \ref{th:segment}, we obtain the following property.

\begin{corollary}\label{th:lexlex}
The Reeves and Stillman component $H_{RS}$ of $\Hilb{n}{p(t)}$ is rational.
\end{corollary}

Reeves and Stillman get the proof by a computation of the Zariski tangent space dimension; however we are able to prove the same result applying our technique as a theoretical tool. Mimicking the notation used in \cite{ReevesStillman}, we denote by $L(b_0,\ldots,b_{n-1})$ the truncation of the lexicographic ideal \eqref{eq:satLexIdeal} in degree $r = \sum_i b_i$, i.e. the hilb-segment ideal w.r.t. $\DegLex$
\[
L(b_0,\ldots,b_{n-1}) = \left\langle x^\alpha\ \big\vert\ x^\alpha \geq_{\DegLex} x_n^{b_{n-1}} \cdots x_1^{b_0} \right\rangle.
\]

\begin{theorem}\label{th:lexstratumopenrational}
The Gr\"obner stratum $\St{\DegLex}\big(L(b_0,\ldots,b_{n-1})\big)$ of the lexicographic ideal associated to the Hilbert polynomial $p(t)$ is isomorphic to an affine space. Therefore the component $H_{RS}$ of $\Hilb{n}{p(t)}$ is rational.
\end{theorem}
\begin{proof}
Thanks to Corollary \ref{th:segment} we can obtain the complete statement proving that the Gr\"obner stratum $\St{\DegLex}(L(b_0,\dots,b_{n-1}))$ is an affine space, that is showing that a same number is both a lower-bound for its dimension and an upper-bound for its embedding dimension; the first part corresponds to Theorem 4.1 (here in terms of initial ideals) and the second one corresponds to Theorem 3.3 of \cite{ReevesStillman}.

We proceed by induction on the number $n$ of variables and on the Gotzmann number $r$.
In order to obtain an upper-bound for the embedding dimension we look for a maximal set of eliminable variables $C''\subset C$, using Criterion \ref{criterio}. If $\{x^{\alpha_1},\ldots,x^{\alpha_n}\}$ is the monomial basis of the saturation $L$ of $L(b_0,\dots,b_{n-1})$, then we can assume to order the polynomials $F_1,\dots,F_{q(r)}\in \K[C][x]$ (that we use to construct $\St{\DegLex}(L(b_0,\dots,b_{n-1}))$) so that $\IN_{\DegLex}(F_1)=x^{\alpha_1}x_0^{r-\vert \alpha_1 \vert} >_{\DegLex} \cdots >_{\DegLex} \IN_{\DegLex}(F_n)=x^{\alpha_n}x_0^{r-\vert \alpha_n \vert}$. So by Theorem \ref{th:saturato} we can initially consider a set of eliminable variables $C''$, containing all the variables $C_{\alpha_j \beta}$ appearing in $F_j$ for every $j>n$.

\medskip
\textbf{Step 1} The zero-dimensional case: $ \St{\DegLex}\big(L(b_0,\dots,0)\big)\simeq \AA^{n b_0}$.

\noindent \textbf{Claim 1i.}  $\dim \St{\DegLex}\big(L(b_0,\dots,0)\big)\geqslant nb_0$.

\noindent The zero-dimensional scheme $Z$ of $b_0$ general points in $\PP^n$ has Gotzmann number $b_0$ and Hilbert polynomial $p(t)=b_0$. Moreover $\IN_{\DegLex}\big(I(Z)_{\geqslant b_0}\big)\supseteq L(b_0,\ldots,b_{n-1})$, because for every monomial $x^{\gamma} \geq_{\DegLex} x_1^{b_0}$ we can find some homogeneous polynomial of the type $x^{\gamma}-\sum_{j=1}^{b_0} c_j x_1^{b_0-j} x_0^j$ vanishing in the $b_0$ points of $Z$: we can find the $c_j$'s solving a $b_0\times b_0$ linear system with a Vandermonde associated matrix. As both $\IN_{\DegLex}\big(I(Z)_{\geqslant b_0}\big)$ and $L(b_0,0,\ldots,0)$ are generated in degree $b_0$, they coincide; so $I(Z)_{\geqslant r} \in \St{\DegLex}\big(L(b_0,0,\ldots,0)\big)$ and we conclude since we can choose $Z$ in a family of dimension $nb_0$.

\noindent \textbf{Claim  1ii.} $\ed \St{\DegLex}\big(L(b_0,\dots,0)\big)\leqslant nb_0$. 

\noindent The saturation of $L(b_0,0,\ldots,0)$ is the ideal $(x_n,x_{n-1},\ldots,x_{2},x_1^{a_0})$, which is generated by $n$ monomials; moreover there are only $b_0$ monomials of degree $b_0$ not contained in $L(b_0,0,\ldots,0)$: Theorem \ref{th:saturato} leads to the conclusion.

\medskip

\textbf{Step 2}  If $\St{\DegLex}\big(L(0,b_1,\dots,b_{n-1})) \simeq \AA^K$ then $ \St{\DegLex}(L(b_0,b_1, \dots, b_{n-1}))\simeq \AA^{K+nb_0}$.

\noindent \textbf{Claim 2i.} $\dim \St{\DegLex}\big(L(b_0,b_1, \dots, b_{n-1})\big) \geqslant \dim  \St{\DegLex}\big(L(0,b_1,\dots,b_{n-1})\big)+nb_0=K+nb_0$.

\noindent
Let $Y$ be any closed subscheme in $\PP^n$ such that $I(Y)_{\geqslant r}\in \St{\DegLex}\big(L(0,b_1,\ldots,b_{n-1})\big)$ and consider the set $Z$ of $b_0$ points in $\PP^n$. If we choose the $b_0$ points in $Z$ general enough, then $I(Z\cup Y)=I(Z)\cdot I(Y)$. Then we conclude thanks to the previous step, as $\IN_{\DegLex}\big(I(Z)\big)=L(b_0,0,\dots, 0)$ and $L(b_0,\ldots,b_{n-1})=L(b_1,\ldots,b_{n-1})\cdot L(b_0,\dots,0)$.

\noindent \textbf{Claim 2ii.} $\ed \St{\DegLex}\big(L(b_0,b_1, \dots, b_{n-1})\big) \leqslant \ed \St{\DegLex}\big(L(0,b_1, \dots, b_{n-1})\big) +  nb_0=K+nb_0$.

First of all, let us consider all the polynomials $F_i$ such that $x_0^{r-b_0} \mid \IN_{\DegLex}(F_i)$ and the set of variables $C_{\alpha_i\beta}$ appearing in them such that $x^{\beta}=x^{\beta_1}x_0^{r-b_1} $ for some monomial $x^{\beta_1}\notin L(0,b_1,\ldots,b_{n-1})$: a multiple of  $x^\beta$ belongs to $L(b_0,\ldots,b_{n-1})$ if and only if the corresponding multiple of $x^{\beta_1}$  belongs to $L(0,b_1,\ldots,b_{n-1})$. Then $F_i=x_0^{r-b_0}F^{(1)}_{i}+\dots $, where the $F^{(1)}_{i}$'s are the polynomials  that appear in the definition of $\St{\DegLex}\big(L(0,b_1,\ldots,b_{n-1})\big)$. Using the $S$-polynomials involving pairs of such polynomials we see that $\mathfrak{L}\big(L(0,b_1,\ldots,b_{n-1})\big)\subseteq \mathfrak{L}\big(L(b_0,\ldots,b_{n-1})\big)$; thus all the variables $C_{\alpha_i\beta}$ of this type are eliminable, except at most $K=\ed \St{\DegLex}\big(L(0,b_1,\ldots,$ $b_{n-1})\big)$ of them.

Moreover, for every $i \leqslant n$ there are $b_0$ variables $C_{\alpha_i\beta}$ such that $x^{\beta}\notin L(b_0,\ldots,b_{n-1})$, $x^{\beta}\in L(0,\ldots,b_{n-1})$: there are  $x_n^{b_{n-1}}\cdots x_{2}^{b_1}x_1^{b_0-j}x_0^j$, $j=1, \dots,b_0$. If we specialize to 0 all the variables of the two above considered types, the embedding dimension drops at most by $\ed \St{\DegLex}\big(L(0,b_1,\ldots,b_{n-1})\big)+nb_0=K+nb_0$. 

Now it will be sufficient to verify that all the remaining variables $C_{\alpha_i\beta}$ are eliminable, using Criterion \ref{criterio}. Assume that $ x^{\beta}<_{\DegLex} x_n^{b_{n-1}}\cdots x_{2}^{b_1}$ and $x_0^{r-b_0}\nmid x^{\beta}$.
\begin{itemize}
\item If $i>n$, all the variables are eliminable using those appearing in $F_1, \dots, F_n$, thanks to Corollary \ref{th:saturato}.
\item If $i<n$, using $S(F_i,F_j)$, where $\IN_{\DegLex}(F_j)=x^{\alpha_i}x_1^{r-\vert \alpha_i \vert}$, we see that $C_{\alpha_i\beta} \in \mathfrak{L}\big(L(b_0,\ldots,b_{n-1})\big)$.  
\item If $i=n$,  using $S(F_n,F_{n-1})=x_{2}x_0^{b_0-1}F_n-x_1^{b_0}F_{n-1}$, we see that $C_{\alpha_n\beta} \in \mathfrak{L}\big(L(b_0,$ $\ldots,b_{n-1})\big)$ (note that by the previous item  $C_{\alpha_{n-1}\beta'} \in \mathfrak{L}(L(b_0,\ldots,b_{n-1}))$).	
\end{itemize}

\medskip 

\textbf{Step 3:} If $L(0,a_1,\dots,a_{n-1}) \simeq \AA^{K_1}$ then $L(0,a_1,\dots, a_d)\simeq \AA^{K_2}$ where $d$ is the maximal index $< n$ such that $a_d\neq 0$ (the degree of the Hilbert polynomial) and $ K_2=K_1+(n-d)(d+1)+\binom{b_{n-1}+n}{n}-1$ (or $K_2=K_1+\binom{b_{n-1}+n}{n}-1$ if $d$ does not exist).

\noindent Here we compare the ideal $L(0,a_1,\dots,a_{n-1})$ in $\K[x]$ and the ideal $L(0,a_1,\dots, a_d)$ in $\K[x_0, \dots,x_d]$.
Observe that both $L(0,a_1,\dots,a_{n-1})^{\sat}$ and $L(0,a_1,\dots, a_d)^{\sat}$ fulfill the hypothesis of Theorem \ref{th:saturato}\emph{(\ref{it:saturato_iv})} (see also Example \ref{rk:stratoLex}); then it holds
\begin{eqnarray*}
&\St{\DegLex}\big(L(b_0,\ldots,b_{n-1})\big)\simeq \St{\DegLex}\left(L(0,a_1,\dots,a_{n-1})^{\sat}\right),&\\
&\St{\DegLex}\big(L(0,a_1,\dots, a_d)\big)\simeq \St{\DegLex}\left(L(0,a_1,\dots, a_d)^{\sat}\right).& 
\end{eqnarray*}
The statement for the saturated ideals $L(0,a_1,\dots, a_{n-1})^{\sat}$ and $L(0,a_1,\dots, a_d)^{\sat}$ is proved using the same technique as above in \cite[Corollary 5.5]{RoggeroTerracini}.
\end{proof}

\newpage

I think that the main advantageous aspects of this approach for a local study of the Hilbert schemes are:
\begin{itemize}
\item[($+$1)] the possibility of exploiting Gr\"obner bases tools. Through both theoretical improvements of algorithms and the constant development of computers, this theory provides very efficient methods to study many non-trivial cases;
\item[($+$2)] Theorem \ref{th:saturato} allows to reduce the number of parameters required for the description of open subsets of $\Hilb{n}{p(t)}$, that in general is very large, being the Hilbert scheme embedded in a suitable Grassmannian.
\end{itemize}

On the other hand, there are also some critical aspects:
\begin{itemize}
\item[($-$1)] the ideals for which the Gr\"obner stratum is an open subset of $\Hilb{n}{p(t)}$ are hilb-segment, i.e. Borel-fixed, but as seen in Section \ref{sec:segments} not every Borel-fixed ideal is a hilb-segment ideal, hence there could be components of $\Hilb{n}{p(t)}$ not containing hilb-segment ideals, i.e. components that can not be studied by means of Gr\"obner strata;
\item[($-$2)] by Lemma \ref{th:sett} we know that considering \lq\lq all\rq\rq\ the monomial ideals with a fixed Hilbert polynomial $p(t)$ and the associated Gr\"obner strata we can cover set-theoretically the Hilbert scheme $\Hilb{n}{p(t)}$. Two questions that arise immediately are 1. what about the scheme structure? 2. is there a smaller set of monomial ideals sufficient to cover all the Hilbert scheme?
\item[($-$3)] Theorem \ref{th:saturato} says also that many of the parameters considered in the costruction of the Gr\"obner stratum are eliminable, so from an algorithmic point of view it would be preferable to avoid to introduce them.
\end{itemize}

Thus the first goal has been to generalize the construction of Gr\"obner strata to a generic Borel-fixed ideal. The main problem is that avoiding the use of a term ordering we give up some of the basic tools we used, primarily the Buchberger's algorithm and the associated noetherian reduction of polynomials. In next section we introduce a new noetherian reduction procedure which does not use any term ordering.

\section{Cioffi and Roggero's results}\label{sec:buch}

In this section we recall the main results exposed in the paper by F. Cioffi and M. Roggero \lq\lq Flat families by Borel-fixed ideals and a generalization of Gr\"obner bases\rq\rq\ \cite{CioffiRoggero} adapting the notation to that used so far. 

\bigskip

From now on, for any monomial ideal $J \subset \K[x]$ we will denote by $G_J$ the set of minimal generators of $J$ and by \glslink{sousEscalier}{$\mathcal{N}(J)$} its \emph{sous-escalier},\index{sous-escalier} that it the set of monomials not belonging to $J$.

\begin{lemma} Let $J$ be a Borel-fixed ideal in $\K[x]$. Then:
\begin{enumerate}[(i)]
\item $x^\alpha \in J\setminus G_J\ \Rightarrow\ \dfrac{x^\alpha}{\min x^\alpha} \in J$;
\item $x^\beta \notin J$ and  $x_i x^\beta \in J \ \Rightarrow$ either $x_ix^\beta \in G_J$ or $x_i > \min x^\beta$.
\end{enumerate}
\end{lemma}
\begin{proof}
Both properties follow from the combinatorial characterization of Borel-fixed ideals.
\end{proof}

For any monomial $x^\alpha \in \K[x]$, we will denote by $x^{\underline{\alpha}}$ the monomial obtained from $x^\alpha$ with the substitution $x_0 = 1$, i.e. $\alpha=(\alpha_0,\ldots,\alpha_n) \rightarrow \overline{\alpha} = (0,\alpha_1,\ldots,\alpha_n)$. Analogously for any monomial ideal $J \subset \K[x]$, \glslink{saturateUnderline}{$\underline{J}$} will be the ideal generated by $\{x^{\underline{\alpha}}\ \vert\ x^\alpha \in G_J\}$. Note that if $J$ is a Borel-fixed ideal $\underline{J}$ coincides with $J^\sat$.

\begin{definition}
For any non-zero homogeneous polynomial $f \in \K[x]$, the \emph{support} of $f$ is the set \gls{suppPoly} of monomials that appear in $f$ with a non-zero coefficient. 
\end{definition}

\begin{definition}[\cite{ReevesSturmfels}]
A \emph{marked polynomial}\index{marked polynomial} is a polynomial $f\in \K[x]$ together with a specified monomial of its support $\Supp f$ that will be called \emph{head term} of $f$ and denoted by \gls{HT}. 
\end{definition} 

\begin{remark}
Although we mainly use the word \lq\lq monomial\rq\rq, we say \lq\lq head term\rq\rq\ for coherency with the notation introduced in \cite{ReevesSturmfels}. Anyway, there will be no possible ambiguity on the meaning of \lq\lq head term of $f$\rq\rq, because we will always use marked polynomials $f$ such that the coefficient of $\Ht(f)$ in $f$ is 1.
\end{remark}

\begin{definition}[\cite{CioffiRoggero}]
Given a monomial ideal $J$ and an ideal $I$, a \emph{$J$-reduced form modulo $I$} of a polynomial $h$ is a polynomial $h_0$ such that $h-h_0\in I$ and $\Supp h_0 \subseteq \cN(J)$. A polynomial is \emph{$J$-reduced} if its support is contained in $\cN(J)$. If {there is a unique} $J$-reduced form modulo $I$ of $h$, we call it \emph{$J$-normal form modulo $I$} and denote it by $\Nf_J(h)$.
\end{definition}

Note that every polynomial $h$ has a unique $J$-reduced form modulo an ideal $I$ if and only if $\cN(J)$ is a $\K$-basis for the quotient $\K[x]/I$ or, equivalently, $\K[x]=I\oplus \langle \cN (J) \rangle $ as a $\K$-vector space. If moreover $I$ is homogeneous, the $J$-reduced form modulo $I$ of a homogeneous polynomial is supposed to be homogeneous too. These facts motivate the following definitions.

\begin{definition}\label{def:defJbase}  
A finite set $G$ of homogeneous marked polynomials $f_\alpha=x^\alpha-\sum c_{\alpha\beta} x^\beta$, with  $\Ht(f_\alpha)=x^\alpha$, is called $J$-\emph{marked set}\index{marked set} if the head terms $\Ht(f_\alpha)$ are pairwise different, form the monomial basis $G_J$ of a monomial ideal $J$ and every $x^\beta$ belongs to $\cN(J)$, i.e. $\vert\Supp f \cap J \vert =1$. We call the polynomial $\Ht(f_\alpha)-f_\alpha$ \emph{tail} of $f_\alpha$ and we denote it by \glslink{tailMP}{$\tail{f_\alpha}{}{}$}, so that $\Supp \tail{f_\alpha}{}{} \subseteq \cN(J)$.  A $J$-marked set $G$ is a $J$-\emph{marked basis}\index{marked basis} if $\cN(J)$ is a basis of $\K[x]/(G)$ as a $\K$-vector space. 
\end{definition}

\begin{definition}
The family of all homogeneous ideals $I$ such that $\cN(J)$ is a basis of the quotient $\K[x]/I$ as a $\K$-vector space will be denoted by \gls{markedFamily} and called $J$-\emph{marked family}.\index{marked family} If $J$ is a Borel-fixed ideal, then $\Mf(J)$ can be endowed with a natural structure of scheme (see \cite[Section 4]{CioffiRoggero}) that we call $J$-\emph{marked scheme}.\index{marked scheme}
\end{definition} 

\begin{remark}\label{Jmarkedbasis}
\begin{enumerate}[(i)]
\item \label{rk:Jmarkedbasis_i} The ideal $(G)$ generated by a $J$-marked basis $G$ has the same Hilbert function of $J$,\index{Hilbert function} hence $\dim_\K J_t = \dim_\K (G)_t$, by the definition of $J$-marked basis itself. Moreover, note that a $J$-marked basis is unique for the ideal that it generates, by the uniqueness of the $J$-normal forms modulo $I$ of the monomials in $G_J$.
\item \label{rk:Jmarkedbasis_ii} $\Mf(J)$ contains every homogeneous ideal having $J$ as initial ideal w.r.t. some term order, but it can also contain other ideals (see \cite[Example 3.18]{CioffiRoggero}).
\item \label{rk:Jmarkedbasis_iii} When $J$ is a Borel-fixed ideal, every homogeneous polynomial has a $J$-reduced form modulo any ideal generated by a $J$-marked set $G$ (\cite[Theorem 2.2]{CioffiRoggero}). 
\end{enumerate}
\end{remark}
 
\begin{proposition}\label{cor1} Let $J$ be a Borel-fixed ideal, $I$ be a homogeneous ideal generated by a $J$-marked set $G$. The following facts are equivalent:  
\begin{enumerate}[(i)]
\item\label{it:cor1_i}  $I\in \Mf(J)$
\item\label{it:cor1_ii} $G$ is a $J$-marked basis;
\item\label{it:cor1_iii} $\dim_\K I_{t}=\dim_\K J_{t}$, for every integer $t$;
\item \label{it:cor1_iv} if $h \in I$ and $h$ is $J$-reduced, then $h=0$.
\end{enumerate}
\end{proposition}

\begin{proof} For the equivalence among the first three statements, see \cite[Corollaries 2.3, 2.4, 2.5]{CioffiRoggero}. For the equivalence among \emph{(\ref{it:cor1_i})} and \emph{(\ref{it:cor1_iv})}, observe that if $I\in \Mf(J)$, then every polynomial has a unique $J$-reduced form; so, the $J$-reduced form of a polynomial of $I$ must to be null. Vice versa, it is enough to show that every polynomial $f$ has a unique $J$-reduced form. Let $\overline{f}$ and $\widetilde{f}$ be two $J$-reduced form of $f$. Then, $\overline{f} -\widetilde{f}$ is a $J$-reduced polynomial of $I$ because $f-\overline{f}$ and $f-\widetilde{f}$ belong to $I$ by definition. We have done, because $\overline{f} -\widetilde{f}$ is null by the hypothesis.
\end{proof}

Thinking for a moment about Hilbert schemes, we want to observe that two different ideals $I_1$ and $\mathfrak{b}$ of the same $J$-marked scheme $\Mf(J)$ give rise to different projective schemes of the same Hilbert scheme $\Hilb{n}{p(t)}$. Indeed, by the uniqueness of the reduced form, there is a monomial $x^\alpha\in G_J$ such that the corresponding polynomials $f_\alpha^{1}$ and $f_\alpha^{2}$ of the $J$-marked bases of $I_1$ and $I_2$, respectively, are different and moreover such that $f_\alpha^{1}\not\in I_2,\ f_\alpha^{2}\not\in {I_1}$. 
If $I_1$ and $I_2$ defined the same projective scheme, we would have $(I_1)_r=(I_2)_r$ for some $r\gg 0$. Hence $x_0^{r-m}f_\alpha^{1}\in I_2$ with normal form modulo $I_2$ given by $x_0^{r-m} f_\alpha^{1} - x_0^{r-m} f_\alpha^{2}$, that is impossible because of Proposition \ref{cor1}\emph{(\ref{it:cor1_iv})}.

\bigskip

Now we come back to deal with Borel-fixed ideals,\index{Borel-fixed ideal} exposing special properties of $J$-marked families in this case.
So from now on $J$ will be always considered Borel-fixed and $G$ will be a $J$-marked set.

\begin{definition}\label{def:$V^J_s$}
If $m_J$ is the initial degree of $J$, we set $V^J_{m_J} = G_{m_J}$; so, for every term $x^\alpha\in G_J$ of degree $m_J$, there is a unique polynomial $g_\alpha \in V^J_{m_J}$ such that $\Ht(g_\alpha)=x^\alpha$. For every $m > m_J$ and for every $x^\alpha\in J_m\setminus G_J$, we set $g_\alpha = x_i g_\epsilon$, where $x_i = \min x^\alpha$ and $g_\epsilon$ is the unique polynomial of $V^J_{s-1}$ with head term $x^\epsilon=\frac{x^\alpha}{\min x^\alpha}$, and we let $V^J_s = G_m \cup \{g_\alpha \text{ s.t. } \ x^\alpha \in J_m\setminus G_J\}$.
\end{definition}

In the following we let $V^J=\cup_s V^J_s$. Moreover, $\langle V^J\rangle$ denotes the vector space generated by the polynomials in $V^J$ and $\xrightarrow{V^J_s}$ is the reduction relation defined in the usual sense of Gr\"obner basis theory. 

\begin{remark}
An ideal $I$ belongs to $\Mf(J)$ if and only if $I=\langle V^J\rangle $ as a vector space. Indeed, for every integer $m$, the number of elements in $V^J_{m}$ is equal to the number of monomials in $J_{m}$.
\end{remark}

\begin{example}
Let us consider the Borel-fixed ideal $J=(x_2^2,x_2x_1,x_1^3) \subset \K[x_0,x_1,x_2]$ and the $J$-marked set
\[
\begin{split}
 f_{x_2^2} &{}= x_2^2 + 2x_1^2 - \frac{1}{3} x_1x_0,\\
 f_{x_2x_1} &{}= x_2x_1 - x_0^2,\\
 f_{x_1^3} &{}= x_1^3 + \frac{1}{5}x_1^2 x_0 - x_1x_0^2 - 3x_0^3.\\
\end{split}
\]
The initial degree of $J$ is 2, so we have:
\begin{description}
\item[$(m=2)$] $J_2 = \left\langle x_2^2,x_2x_1 \right\rangle$
\[
V_2^J = \left\{g_{x_2^2}=f_{x_2^2},g_{x_2x_1}=f_{x_2x_1}\right\};
\]
\item[$(m=3)$] $J_3 = \left\langle x_2^3,x_2^2 x_1,x_2 x_1^2, x_1^3, x_2^2 x_0, x_2x_1x_0\right\rangle$
\[
\begin{split}
V_3^J = \big\{& g_{x_2^3} = x_2 g_{x_2^2}, g_{x_2^2 x_1} = x_1 g_{x_2^2}, g_{x_2 x_1^2} = x_1 g_{x_2x_1},\\
               & g_{x_2^2 x_0} = x_0 g_{x_2^2}, g_{x_2 x_1 x_0} = x_0 g_{x_2x_1}\big\} \cup \left\{ g_{x_1^3}=f_{x_1^3}\right\};
\end{split}
\]
\item[$(m=4)$] $J_4 = \left\langle x_2^4,x_2^3 x_1,x_2^2 x_1^2, x_2x_1^3,x_1^4, x_2^3x_0,x_2^2 x_1x_0,x_2 x_1^2x_0, x_1^3x_0, x_2^2 x_0^2, x_2x_1x_0^2\right\rangle$
\[
\begin{split}
V_4^J = \big\{& g_{x_2^4} = x_2 g_{x_2^3} = x_2^2 f_{x_2^2}, g_{x_2^3 x_1} = x_1 g_{x_2^3} = x_2x_1 f_{x_2^2},g_{x_2^2 x_1^2} = x_1 g_{x_2^2 x_1} = x_1^2 f_{x_2^2},\\
& g_{x_2x_1^3} = x_1 g_{x_2 x_1^2} = x_1^2 f_{x_2x_1}, g_{x_1^4} = x_1 g_{x_1^3} = x_1 f_{x_1^3},\\
& g_{x_2^3x_0} = x_0 g_{x_2^3} = x_2x_0 f_{x_2^2}, g_{x_2^2 x_1x_0} = x_0 g_{x_2^2 x_1} = x_1x_0 f_{x_2^2},\\
& g_{x_2 x_1^2x_0} = x_0 g_{x_2 x_1^2} = x_1x_0 f_{x_2x_1}, g_{x_1^3x_0} = x_0 g_{x_1^3} = x_0 f_{x_1^3},\\
& g_{x_2^2 x_0^2} = x_0 g_{x_2^2 x_0} = x_0^2 f_{x_2^2}, g_{x_2x_1x_0^2} = x_0 g_{x_2x_1x_0} = x_0^2 f_{x_2x_1}\big\}.
\end{split}
\]
\end{description}
\end{example}

Keeping in mind the canonical decomposition and the decomposition map of a Borel-fixed ideal introduced in Definition \ref{def:canonicalDec}, we prove the following lemma.

\begin{lemma}\label{descLex}
Let $J$ be a Borel-fixed ideal. If $x^\epsilon$ belongs to $\cN(J)$ and $x^{\epsilon}\cdot x^\delta=x^{\epsilon+\delta}$ belongs to $J$ for some $x^\delta$, then $x^{\epsilon+\delta}= \dec{x^\alpha}{x^\eta}{J}{}$ with $x^{\eta}<_{\Lex}x^\delta$.  Furthermore:
\begin{enumerate}[(i)]
\item if $\vert \delta\vert=\vert\eta\vert$, then $x^{\eta}<_{B}x^\delta$; and
\item $x^{\underline{\eta}}<_{\Lex}x^{\underline{\delta}}$.
\end{enumerate}
\end{lemma}

\begin{proof}
We can assume that $x^{\delta}$ and $x^{\eta}$ are coprime; indeed, if this is not the case, we can divide the involved equalities of monomials by  $\gcd(x^{\delta},x^{\eta})$. If $x^\eta=1$, all the statements are obvious. If $x^\eta\neq 1$, then $\min x^\delta \mid x^\alpha$ because $x^{\delta}$ and $ x^{\eta}$ are coprime, hence $ \min x^\delta \geq \min x^\alpha \geq \max x^\eta $ and so $\min x^\delta > \max x^\eta $ because they cannot coincide. This inequality implies both $x^{\eta}<_{\Lex}x^\delta$ and $x^{\underline{\eta}}<_{\Lex}x^{\underline{\delta}}$. Moreover, if $\vert \delta\vert=\vert\eta\vert$, this is also sufficient to conclude that $x^{\eta}<_{B}x^\delta$.
\end{proof}

\begin{remark} \label{rm:V}
Observe that if $g_\beta=x^\delta f_\alpha$ belongs to $V^J_s$, then $x^\beta=\dec{x^\alpha}{x^\delta}{J}{}$.
\end{remark}

We have already recalled that, when $J$ is a Borel-fixed ideal, every homogeneous polynomial has a $J$-reduced form modulo an ideal generated by a $J$-marked set $G$ (Remark \ref{Jmarkedbasis}\ref{rk:Jmarkedbasis_iii}). Further, a $J$-reduced form of a homogeneous polynomial can be constructed by a suitable reduction relation, as it is recalled by next Proposition. 

\begin{proposition}[{\cite[Proposition 3.6]{CioffiRoggero}}] \label{construction of $J$-normal form}
With the above notation, every monomial $x^\beta\in J_m$ can be reduced to a $J$-reduced form modulo $(G)$ in a finite number of reduction steps, using only polynomials of $V^J_s$. Hence, the reduction relation $\xrightarrow{V^J_s}$ is Noetherian. 
\end{proposition}

The Noetherianity of the reduction relation $\xrightarrow{V^J_s}$ provides an algorithm that reduces every homogeneous polynomial of degree $m$ to a $J$-reduced form modulo $(G)$ in a finite number of steps. We note that on one hand it is convenient to substitute the polynomials in $V^J_s$ by their $J$-reduced normal forms for an efficient implementation of a reduction algorithm, but on the other hand it is convenient to use in the proofs the polynomials of $V^J_s$ as constructed in Definition \ref{def:$V^J_s$}. 

Using the Noetherianity of the reduction relation $\xrightarrow{V^J_s}$, we can recognize when a $J$-marked set is a $J$-marked basis by a Buchberger-like criterion \cite[Theorem 3.12]{CioffiRoggero}. To this aim we need to pose an order on the set $W^J_m =\{x^\delta f_\alpha \text{ s.t. } f_\alpha \in G \text{ and } \vert \delta+\alpha\vert=m\}$ that becomes a set of marked polynomials by letting $\Ht(x^\delta f_\alpha)=x^{\delta+\alpha}$. Note that $I_m$ is generated by $W^J_m$ as a vector space. We set $W^J=\cup_m W^J_m$.

\begin{remark}
We point out that Definition 3.9 of \cite{CioffiRoggero} does not work well for our purpose. Hence, we introduce the following Definition \ref{order2} and provide a new proof of \cite[Lemma 3.10]{CioffiRoggero} by the following Lemma \ref{reduction2}.
\end{remark}

\begin{definition}\label{order2}
Let $\leq$ be any order on $G$ and $x^\delta f_\alpha$, $x^{\delta'}f_{\alpha'}$ be two elements of $W^J_m$. We set 
\begin{equation}
x^\delta f_\alpha \geq_m x^{\delta'}f_{\alpha'} \Leftrightarrow 
x^\delta >_{\Lex} x^{\delta'} \text{ or } x^\delta = x^{\delta'} \text{ and } f_\alpha \geq f_{\alpha'}.
\end{equation}
\end{definition}

\begin{lemma}\label{reduction2}
\begin{enumerate}[(i)]
\item \label{reduction2-i} For every two elements $x^\delta f_\alpha$, $x^{\delta'}f_{\alpha'}$ of $W^J_m$ we get 
\[
x^\delta f_\alpha \geq_m x^{\delta'}f_{\alpha'} \Rightarrow \forall x^\eta :\  x^{\delta+\eta} {f_\alpha} \geq_{m'} x^{\delta'+\eta}f_{\alpha'},
\]
 where $m'=\vert \delta +\eta + \alpha\vert$.
\item \label{reduction2-ii} Every polynomial $g_\beta\in V^J_s$ is the minimum w.r.t. $\leq_m$ of the subset $W^J_\beta$ of $W^J_m$ containing all polynomials of $W^J_m$ with $x^\beta$ as head term.
\item \label{reduction2-iii}If $x^\delta f_\alpha$ belongs to $W^J_m\setminus G_m \text{ and } x^\beta$ belongs to $\Supp x^\delta f_\alpha \setminus \{x^\delta x^\alpha\} \text{ with } g_\beta\in V^J_s$, then $x^\delta f_{\alpha} \succ_m g_\beta$.
\end{enumerate}
\end{lemma}
\begin{proof} 
\emph{(\ref{reduction2-i})} This follows by the analogous property of the term order $>_{\Lex}$.

\emph{(\ref{reduction2-ii})} Let $g_\beta=x^{\delta'}f_{\alpha'}$ be the polynomial of $V^J$ such that $x^\beta= \dec{x^{\alpha'}}{x^{\delta'}}{J}{}$ and $x^\delta f_\alpha$ be another polynomial of $W^{J}_{\beta}$. We can assume that $x^{\delta}$ and $x^{\delta'}$ are coprime; otherwise, we can divide the involved inequalities of monomials by $\gcd(x^{\delta},x^{\delta'})$. By Remark \ref{rm:V}, we have that $\max x^{\delta'} \leq \min x^{\alpha'}$ and $\max x^{\delta} > \min x^{\alpha}$. Then, we get $\max x^\delta > \max x^{\delta'} $ because $x^{\alpha'}\nmid x^\alpha$ and $x^{\alpha}\nmid x^{\alpha'}$. Thus, $x^\delta >_{\Lex} x^{\delta'}$.

\emph{(\ref{reduction2-iii})} If $x^\beta$ belongs to $G_J$ we are done. Otherwise, let $x^\beta=\dec{x^{\alpha'}}{x^{\delta'}}{J}{}$ and note that every term of $\Supp x^\delta f_\alpha$ is a multiple of $x^\delta$, in particular $x^{\beta}=x^{\delta+\gamma}$ for some $x^\gamma \in \cN(J)$. By Lemma \ref{descLex}, we get $x^{\delta'}<_{\Lex}x^\delta$. 
\end{proof}

\begin{definition}
The \emph{S-polynomial} of two elements $f_\alpha$, $f_{\alpha'}$ of a $J$-marked set $G$ is the polynomial $S(f_\alpha, f_{\alpha'}) =x^\gamma f_\alpha-x^{\gamma'}f_{\alpha'}$, where $x^{\gamma+\alpha}=x^{\gamma'+\alpha'}= \gls{lcm}(x^\alpha,x^{\alpha'})$.
\end{definition}

\begin{theorem}[Buchberger-like criterion]\label{BuchCrit1}
Let $J$ be a Borel-fixed ideal and $I$ the homogeneous ideal generated by a $J$-marked set $G$. With the above notation, the following statements are equivalent:
\begin{enumerate}[(i)]
\item \label{buchcrit1-i} $I\in \Mf(J)$;
\item \label{buchcrit1-ii} $\forall f_\alpha, f_{\alpha'} \in G,\  S(f_\alpha, f_{\alpha'}) \stackrel{V^J_s}\longrightarrow 0$;
\item \label{buchcrit1-iii}$ \forall f_\alpha, f_{\alpha'} \in G,\  S(f_\alpha, f_{\alpha'})=x^\gamma f_\alpha-x^{\gamma'}f_{\alpha'} =\sum a_j x^{\eta_j}f_{\alpha_j}$, with $x^{\eta_j}<_{\Lex}\max_{\Lex}\{x^\gamma,x^{\gamma'}\}$ and $f_{\alpha_j}\in G$.
\end{enumerate}
\end{theorem}
\begin{proof}
$\emph{(\ref{buchcrit1-i})} \Rightarrow \emph{(\ref{buchcrit1-ii})}$ Recall that $I \in \Mf(J)$ if and only if $G$ is a $J$-marked basis, so that every polynomial has a unique $J$-normal form modulo $I$. Since $S(f_\alpha, f_{\alpha'})$ belongs to $I$ by construction, its $J$-normal form modulo $I$ is null.

$\emph{(\ref{buchcrit1-ii})} \Rightarrow \emph{(\ref{buchcrit1-iii})}$ Straightforward by the definition of the reduction relation $\xrightarrow{V^J_s}$ and by Lemma \ref{reduction2}\emph{(\ref{reduction2-iii})}.

$\emph{(\ref{buchcrit1-iii})} \Rightarrow \emph{(\ref{buchcrit1-i})}$ We want to prove that $I=\langle V^J \rangle$ or, equivalently, that $\langle V^J\rangle=\langle W^J\rangle$.  It is sufficient to prove that $x^\eta\cdot V^J\subseteq \langle V^J\rangle$, for every monomial  $x^\eta$. We proceed by induction on the monomials $x^\eta$, ordered according to $\Lex$. The thesis is obviously true for $x^\eta=1$. We then assume that the thesis holds for any monomial $x^{\eta'}$ such that $x^{\eta'}<_{\Lex}x^\eta$. 

If $\vert\eta\vert> 1$, we can consider any product $x^{\eta}=x^{\eta_1}\cdot x^{\eta_2}$, $x^{\eta_1}$ and $x^{\eta_2}$ non-constant. Since $x^{\eta_i}<_{\Lex}x^\eta, i=1,2$, we immediately obtain by induction
\[
x^\eta\cdot V^J= x^{\eta_1}\cdot (x^{\eta_2}\cdot V^J)\subseteq x^{\eta_1}\langle V^J\rangle\subseteq\langle V^J\rangle.
\]

If $\vert \eta\vert=1$, then we need to prove that $x_i\cdot V^J\subseteq\langle V^J\rangle$. Since $x_0V^J\subseteq  V^J$, it is then sufficient to prove the thesis for $x^\eta=x_i$, assuming that the thesis holds for every $x^{\eta'}<_{\Lex}x_i$. We consider $g_\beta=x^\delta f_\alpha \in V^J$, where $\max x^\delta\leq \min x^\alpha$. If $x_ig_\beta$ does not belong to $V^J$, then $\max (x_i\cdot x^\delta)>\min x^\alpha $, so $x_i>\min x^\alpha $. In particular, $x_i>\min x^\alpha \geq \max x^\delta $, so $x_i >_{\Lex} x^\delta$: by induction, it is now sufficient to prove the thesis for $x_i f_\alpha$.

We consider an $S$-polynomial $S(f_\alpha,f_{\alpha'})=x_i f_\alpha-x^\gamma f_{\alpha'}$ such that $x^\gamma<_\Lex x_i$. Such $S$-polynomial always exists: for instance, we can consider $x_ix^\alpha= \dec{x^{\alpha'}}{x^{\eta'}}{J}{}$. 
By the hypothesis $x_i f_\alpha-x^{\eta'}f_{\alpha'} =\sum a_j x^{{\eta'}_j}f_{\alpha_j}$ where $x^{\eta'}f_{\alpha'},  x^{{\eta'}_j}f_{\alpha_j}\in \langle V^J\rangle$ by induction since $x^{\eta'}, x^{{\eta'}_j}$ are lower than $x_i$ w.r.t. $\Lex$ and then $x_i f_\alpha$ belongs to $\langle V^J\rangle$.
\end{proof}

\begin{remark}
In \cite[Section 3]{CioffiRoggero} some results about syzygies of the ideal $I$ generated by a $J$-marked basis are proposed, by using the order on $W^J_m$ defined in Definition 3.9 of \cite{CioffiRoggero} that does not work well. Anyway, the order defined in Definition \ref{order2} works well also in that context of syzygies.
\end{remark}

\begin{definition}
We call \emph{Eliahou-Kervaire couple} of the $J$-marked set $G$ any couple of polynomials $f_\alpha, f_\beta$, $\Ht(f_\alpha)=x^\alpha$, $\Ht(f_\beta)=x^\beta$, such that
\[
x_jx^\alpha=\dec{x^\beta}{x^\eta}{J}{} \text{ for some } x_j>\min x^\alpha.
\]
We call \emph{Eliahou-Kervaire $S$-polynomial} (EK-polynomial, for short) of $G$ an $S$-polynomial among an Eliahou-Kervaire couple of polynomials $f_\alpha$ and $f_\beta$. We denote such $S$-polynomial by $S^{\text{EK}}(f_\alpha,f_\beta)$. Observe that, thanks to the definition, an EK-polynomial is of kind
\[
S^{\text{EK}}(f_\alpha,f_\beta)=x_jf_\alpha-x^\eta f_\beta, \text{ for some } x_j>\min x^\alpha, \text{ with } x_jx^\alpha=\dec{x^\beta}{x^\eta}{J}{}.
\]
\end{definition}

\begin{remark}\label{Buch-proof3}
We underline that in the proof of Theorem \ref{BuchCrit1} the crucial point is the existence of an $S$-polynomial of kind $x_if_\alpha-x^\eta f_\alpha$ with $x^\eta<_{\Lex}x_i$, and we use an EK-polynomial. An analogous argument will be used in the proof of Theorem \ref{nostrobuch} and of Theorem \ref{paolo}.
\end{remark}

As pointed out in Remark \ref{Buch-proof3}, in the proof of Theorem \ref{BuchCrit1} we just need to assume that \emph{(\ref{buchcrit1-iii})} holds for EK-polynomials. Hence, we have the following result.

\begin{corollary}\label{sizEK}
With the same notation of Theorem \ref{BuchCrit1},
\[
I\in \Mf(J)\Leftrightarrow \text{ for every EK-polynomial, } S^{\text{EK}}(f_\alpha,f_\beta)\xrightarrow{V^J_s}0.
\]
\end{corollary}

\section[Superminimal generators and a new Noetherian reduction]{Superminimal generators and a new Noetherian\\ reduction}

Form this section onwards, we consider a Borel-fixed ideal $J$ obtained as truncation in degree $m$ of a saturated ideal, that is $J = \underline{J}_{\geqslant m}$ for some integer $m$. In this case we have $(G_J)_{>m}\subset G_{\underline{J}}$.

\begin{definition}\label{defsupermin}
The \emph{set of superminimal generators}\index{superminimal generators} of $J$ is
\[
sG_J=\{x^{\alpha} \in G_J   \ \vert\  x^{\underline{\alpha}}  \in G_{\underline{J}}\}.
\] 
Hence, for every $x^\alpha \in sG_J$, we have an integer $t_\alpha = \alpha_0$ such that $x^\alpha=x_0^{t_\alpha} x^{\underline{\alpha}}$; more precisely, 
\[
t_\alpha=\begin{cases}
0, &\text{ if } \vert\underline{\alpha}\vert \geq m\\
m-\vert\underline{\alpha}\vert, &\text{ otherwise}
\end{cases}.
\]
Given a $J$-marked set $G$, the set $sG$ of \emph{superminimal generators} of $G$ is
\[
sG = \left\{f_{\alpha}= x^{\alpha}-\tail{f_{\alpha}}{}{}\in G\ \vert x^{\alpha} \in  sG_J\right\}.
\]
\end{definition}

\begin{definition}\label{sminred} Given a $J$-marked set $G$ and two polynomials $h$ and $h_1$, we say that $h$ is in \emph{$sG$-relation} with $h_1$ if there is a monomial $x^\gamma \in \Supp h \cap  J$, $\Coeff(x^\gamma)=c$, such that $x^\gamma$ is divisible by a superminimal generator $x^\alpha$ of $J$ with $x^\gamma=x^{\alpha}\cdot x^\epsilon = \dec{x^{\underline{\alpha}}}{x^\eta}{\underline{J}}{}$ and $h_1=h-c\cdot x^\epsilon f_\alpha$, that is $h_1$ is obtained by replacing in $h$ the monomial $x^\gamma$ by $x^\epsilon \cdot \tail{f_{\alpha}}{}{}$.
We call \emph{superminimal reduction}\index{superminimal reduction} the transitive closure of the above relation and denote it by $\xrightarrow{sG}$. Moreover, we say that:
\begin{itemize} 
\item \emph{$h$ can be reduced to $h_1$} by $\xrightarrow{sG}$ if $h \xrightarrow{sG}h_1$;
\item \emph{$h$ is reduced w.r.t. $sG$} if no monomial in $\Supp h$ is divisible by a monomial of $sG_J$;
\item \emph{$h$ is strongly reduced} if no monomial in $\Supp h$ is divisible by a monomial of $G_{\underline{J}}$, that is $h$ is $\underline{J}$-reduced. In other words, $h$ is strongly reduced if for every $t$, $x_0^t\cdot h$ is reduced w.r.t. $sG$.
\end{itemize}
\end{definition}

\begin{remark}\label{code}
Given a polynomial $f_\alpha$ of a $J$-marked set $G$ and any positive integer $t$, then $\Supp (x_0^t\cdot \tail{f_\alpha}{}{})\subseteq \cN(J)$. Furthermore, if $G$ is a $J$-marked basis, then we also have $x_0^t\cdot \tail{f_\alpha}{}{}=\Nf(x_0^t\cdot x^\alpha)$ since in this case the tail of $f_\alpha$ is indeed  the $J$-normal form modulo $(G)$ of $x^\alpha$ (see Remark \ref{Jmarkedbasis}\ref{rk:Jmarkedbasis_i}). More generally, if $h$ is a homogeneous polynomial of degree $\deg(h)\geqslant m$, then $\Nf(x_0^t\cdot h)= x_0^t\cdot \Nf(h)$.
\end{remark}

It is interesting to notice that the subset $sG$ of Definition \ref{defsupermin} is a subset of $V$, but not every step of reduction by $\xrightarrow{sG}$ is also a step of reduction by $\xrightarrow{V^J_s}$, as shown in the following example.

\begin{example} 
Consider the Borel-fixed ideal
\[
J = \underline{J}_{\geqslant 2} = (x_2,x_1^2)_{\geqslant 2} = (x_2^2,x_2x_1,x_1^2,x_2x_0)
\]
and let $G$ be a $J$-marked set. Consider the monomial $x_2^2 x_0$. The only way to reduce $x_2^2 x_0$ via $\xrightarrow{V_3^J}$ leads to $x_0\cdot \tail{g_{x_2^2}}{}{}$, where $g_{x_2^2} = f_{x_2^2} \in V^J_2$. Moreover, $x_0\cdot \tail{g_{x_2^2}}{}{}$ is not further reducible, because all the monomials of its support belong to $\cN(J)$. On the other hand, according to Definition \ref{sminred}, a first step of reduction of the monomial $x_2^2 x_0$ via  $\xrightarrow{sG}$ is $ x_2^2 x_0\xrightarrow{sG} x_2 \cdot \tail{f_{x_2x_0}}{}{}$, where $f$ is the polynomial in $sG$ with $\Ht(f)=x_2 x_0$. Since $x_2$ is a monomial of $G_{\underline{J}}$, every monomial appearing in $\Supp \big(x_2\cdot \tail{f_{x_2x_0}}{}{}\big)$ belongs to $J=\underline{J}_{\geqslant 2}$, and so we will need further steps of reduction via $\xrightarrow{sG}$ to compute a polynomial reduced w.r.t. $sG$.
\end{example}

\begin{theorem}\label{ridsm} With the above notation:
\begin{enumerate}[(i)]
\item\label{ridsm-i} $\xrightarrow{sG}$ is Noetherian;
\item\label{ridsm-ii} for every homogeneous polynomial $h$ there exist $\overline{t}= \overline{t}(h)$ and a polynomial $\overline{h}$ strongly reduced such that $x_0^{\overline{t}}\cdot h\xrightarrow{sG} \overline{h}$;
\item\label{ridsm-iii} If moreover $G$ is a $J$-marked basis and $\deg h \geqslant m$, then  $\overline{h}=\Nf (x_0^{\overline{t}}\cdot h)= x_0^{\overline{t}} \cdot \Nf(h)$ where $\Nf(h)$ is the unique $J$-normal form modulo $(G)$ of $h$.
\end{enumerate}
\end{theorem}

\begin{proof}
\emph{(\ref{ridsm-i})}  If $\xrightarrow{sG}$ was not Noetherian, by Lemma \ref{descLex} and Lemma \ref{reduction2}, we would be able to find infinite descending chains of monomials w.r.t. $<_{\Lex}$.

\emph{(\ref{ridsm-ii})}  It is sufficient to prove the thesis for monomials $x^\gamma$ in $\underline{J}$. Let $x^\gamma= \dec{x^{\underline{\alpha}}}{x^\eta}{\underline{J}}{}$. If $x^\eta=1$, then $x^{\alpha}= x_0^{t_\alpha}\cdot x^{\underline{\alpha}}$ is in $sG_J$, $f_{\alpha}$ belongs to $sG$ and $x_0^{t_\alpha} \cdot x^{\underline\alpha} \xrightarrow{sG}\tail{f_{\alpha}}{}{}$, where $\Supp \tail{f_{\alpha}}{}{} \subseteq \cN(J)$. In this case $\overline{h}=\tail{f_{\alpha}}{}{}$ and $\overline{t}=t_\alpha$. 

If $x^\eta\neq 1$, we can assume that the thesis holds for any monomial $x^{\gamma'}= \dec{x^{\underline{\beta}}}{x^{\eta'}}{\underline{J}}{}$, such that $x^{\eta'} <_{\Lex} x^{\eta}$. We perform a first reduction $x_0^{t_{\alpha}} \cdot x^\gamma \xrightarrow{sG}  x^\eta\cdot \tail{f_\alpha}{}{}$. If $x^\eta\cdot \tail{f_\alpha}{}{}$ is strongly reduced, we are done. Otherwise, we have $x^\eta \neq x_0^{\vert \eta \vert}$. For every monomial  $x^{\gamma'}\in \Supp \big(x^\eta\cdot \tail{f_\alpha}{}{}\big)\cap \underline{J}$ we have $x^{\gamma'}=\dec{x^{\underline{\beta}'}}{x^{\eta'}}{\underline{J}}{}$, with  $x^{\eta'}<_{\Lex}x^\eta$ by Lemma \ref{descLex}. So, we have also $x_0^t  \cdot x^{\eta'} <_{\Lex} x^\eta$, for every $t$. By the inductive hypothesis we can find a suitable power $t$ of $x_0$ such that every monomial in $x_0^t\cdot x^\eta \cdot \tail{f_\alpha}{}{}$ can be reduced by $\xrightarrow{sG}$ to a strongly reduced polynomial. Thus  $x_0^t\cdot x^\eta \cdot \tail{f_\alpha}{}{} \xrightarrow{sG}\overline{h}$ with $\overline{h}$ strongly reduced. In this case $\overline{t}(x^\gamma)=t_\alpha+ t=t_\alpha + \overline {t}\big(x^\eta \cdot \tail{f_\alpha}{}{}\big)$.

\emph{(\ref{ridsm-ii})} Since $\deg h \geqslant m$, we observe that $x_0^{\overline{t}}\cdot h-\overline{h}\in (G)$ and  $\Supp \overline{h}\subseteq \cN(\underline{J})\subseteq \cN(J)$, hence  $\overline{h}$ is a $J$-reduced form modulo $(G)$ of $x_0^{\overline{t}}\cdot h$. Therefore, if $G$ is a $J$-marked basis, $\overline{h}$ is the unique $J$-normal form of $x_0^{\overline{t}}\cdot h$. Moreover, $\overline{h}=\Nf(x_0^{\overline{t}}\cdot h)= x_0^{\overline{t}}\cdot \Nf(h)$ because $\deg h\geqslant m$ (see Remark \ref{code}).
\end{proof}

If we consider the more general setting used in Section \ref{sec:buch}, we are not able to generalize the properties of the reduction $\xrightarrow{V^J_s}$ to $\xrightarrow{sG}$. Indeed, in our proofs we will often need that polynomials $f_\alpha\in G$ have the following property:
\begin{equation}\label{propJ}
\forall \  x^\beta\in \Supp \tail{f_\alpha}{}{}, \ x^\beta \in \cN(J_0).
\end{equation}
Assume that $J$ is Borel-fixed, with initial degree $m_J$, that $\underline{J}$ is its saturation and $m$ is such that $J_{\geqslant m}=\underline{J}_{\geqslant m}$.
Condition \eqref{propJ} is necessary for the properties we will prove (see Example \ref{jnontagl}) and it obviously holds for $J=\underline{J}_{\geqslant {m_J}}$, that is when $m_J=m$. If we work in a more  general setting, in which $m>m_J$, then we  need to assume that \eqref{propJ} holds for the $J$-marked set we consider; however, in this way, we are able to characterize only a subset of the $J$-marked scheme.

Theorem \ref{ridsm} has some interesting consequences.

\begin{corollary}\label{hsconsm}
Let $I$ be an ideal generated by a $J$-marked set $G$. Then:  
\[
I\in \Mf(J) \Longleftrightarrow \forall\ h\in I, \exists\ t \text{ s.t. } x_0^t\cdot h\xrightarrow{sG}0.
\]
\end{corollary}
\begin{proof} Let $h\in I$. If  $I$ belongs to $\Mf(J)$, then $G$ is a $J$-marked basis and the equivalent condition of Proposition \ref{cor1}\emph{(\ref{it:cor1_iv})} holds. Especially $\Nf(h)=0$. Moreover, by Theorem \ref{ridsm}\emph{(\ref{ridsm-iii})}, we have $x_0^t \cdot h\xrightarrow{sG} x_0^t\cdot \Nf(h)$ for a suitable $t$, and we conclude.

Vice versa, we use again Proposition \ref{cor1}. For every  $h\in I$ such that $h$ is $J$-reduced modulo $I$, then $h$ is also strongly reduced w.r.t. $sG$, that is  $x_0^t\cdot h$ is not further reducible through $\xrightarrow{sG}$ for every $t$ and by the hypothesis $x_0^t\cdot h=0$. Thus $h=0$ and we conclude.
\end{proof}

\begin{corollary} \label{unicsG}
Given a set of marked polynomials $\Gamma=\{f_\beta \ \vert \ \Ht(f_\beta)=x^\beta \in sG_J,$ $\Supp (x^\beta-f_\beta)\subset \cN(J)\}$, there is at most one ideal $I$ of $\Mf(J)$ such that $\Gamma$ is the set of superminimals of $I$.
\end{corollary}
\begin{proof}
Suppose that $I$ is an ideal in $\Mf(J)$ such that $\Gamma$ is its set of superminimals and let $G$ be its $J$-marked basis. If $f_\alpha$ is a polynomial of $G$, then $f_\alpha=x^\alpha-\Nf(x^\alpha)$, where $\Nf(x^\alpha)$ is uniquely determined by $\Gamma$ in the following way: for every $x^\alpha\in G_J\setminus sG_J$, we consider an integer $t$ such that $x_0^t\cdot x^\alpha \xrightarrow{\Gamma}\overline{h}$ with $\overline{h}$ strongly reduced, that exists by Theorem \ref{ridsm}, and set $\Nf(x^\alpha)=\frac{\overline{h}}{ x_0^t}$. Note that if $x_0^t$ does not divide $\overline{h}$, then such an ideal $I$ does not exist.
\end{proof}

Now we will prove that the Buchberger-like criterion introduced for $\xrightarrow{V^J_s}$ can be rephrased in terms of the $\xrightarrow{sG}$ reduction, showing an analogous of Theorem \ref{BuchCrit1} for $\xrightarrow{sG}$. 

\begin{lemma}\label{x0V}
Let $h$ be a homogeneous polynomial of degree $s\geqslant m$.
\[
h\in \langle V^J_s\rangle \Longleftrightarrow x_0\cdot h\in \langle V^J_{s+1} \rangle.
\]
\end{lemma}
\begin{proof}
If $h \in \langle V^J_s\rangle$, then $x_0\cdot h\in \langle V_{s+1}\rangle $ by definition of $V$.

Vice versa, assume that $x_0\cdot h \in \langle V_{s+1}\rangle$. This is equivalent to $x_0\cdot h\xrightarrow{V^J_{s+1}}0$. Every monomial in $\Supp (x_0\cdot h)$ can be written as $x_0\cdot x^\epsilon$; observe that $x_0\cdot x^\epsilon \notin sG_J$, because $\deg x_0\cdot x^\epsilon > m$. Then, if $x_0\cdot x^\epsilon$ belongs to $J$, we can decompose it as $x_0\cdot x^\epsilon = \dec{x^\alpha}{x^\eta}{J}{}$, $x^\alpha \in G_J$ and $x^\eta\neq 1$. Since $\min x^\alpha \geq \max x^\eta$, we have that $x^\eta$ is divisible by $x_0$. So $x^\eta=x_0\cdot x^{\eta'}$. 

Summing up, in order to reduce the monomial $x_0\cdot x^\epsilon$ of $x_0\cdot h$ using $V^J$, we use the polynomial $x_0\cdot x^{\eta'}\cdot f_\alpha \in V^J$,  $\Ht(f_\alpha)=x^\alpha$. If the coefficient of $x_0\cdot x^\epsilon$ in $x_0\cdot h$ is $a$, we obtain
\[
x_0\cdot h\xrightarrow{V^J_{s+1}} x_0\cdot(h-a\cdot x^{\eta'}f_\alpha).
\] 
At every step of reduction, we obtain a polynomial which is divisible by $x_0$. In particular,
\[
x_0\cdot h \in \langle V^J_{s+1}\rangle \Rightarrow x_0\cdot h=x_0\cdot \sum a_i x^{\eta_i}f_{\alpha_i},\text{ where } x_0\cdot x^{\eta_i}f_{\alpha_i}\in V^J_{s+1}.
\]
Then we have that $h= \sum a_i x^{\eta_i}f_{\alpha_i}$ and $x^{\eta_i}f_{\alpha_i}\in V^J_{s}$, that is $h\in \langle V^J_{s}\rangle$.
\end{proof}

Consider $f_\alpha,f_{\alpha'}\in G$, the $S$-polynomial $S(f_\alpha, f_{\alpha'})=x^\gamma f_\alpha-x^{\gamma'}f_{\alpha'}$ and assume that $x^{\gamma'}<_{\Lex}x^{\gamma}$. By Lemma \ref{reduction2}\emph{(\ref{reduction2-iii})}, if $S(f_\alpha,f_\alpha')\xrightarrow{V^J_s} h$, then $S(f_\alpha,f_\alpha')-h=\sum a_jx^{\delta_j}f_{\beta_j}$
with $x^{\delta_j}f_{\beta_j} \in V^J_s$, $x^{\delta_j}<_{\Lex}x^\gamma$. Now we show that a similar result holds for the superminimal reduction $\xrightarrow{sG}$.

\begin{lemma}\label{lexrid}
Consider $f_\alpha,f_{\alpha'}\in G$, the $S$-polynomial $S(f_\alpha, f_{\alpha'})=x^\gamma f_\alpha-x^{\gamma'}f_{\alpha'}$ and assume that $x^{\gamma'}<_{\Lex}x^{\gamma}$. If $x_0^t\cdot S(f_\alpha,f_\alpha')\xrightarrow{sG} h$, then $x_0^t\cdot S(f_\alpha,f_\alpha')-h=\sum a_jx^{\eta_j}f_{\beta_j}$
with $f_{\beta_j} \in sG$, $x^{\eta_j}<_{\Lex}x^\gamma$ and $x^{\underline{\eta_j}}<_{\Lex}x^{\underline{\gamma}}$. 
\end{lemma}

\begin{proof}
Consider a monomial $x_0^t\cdot x^\gamma\cdot x^\epsilon$ in $\Supp \big(x_0^t\cdot x^\gamma\cdot \tail{f_\alpha}{}{}\big)\cap J$. Such a monomial decomposes as $x_0^t\cdot x^\gamma\cdot x^\epsilon= \dec{x^{\underline{\beta}}}{x^{\eta}}{\underline{J}}{}$, $x^{\eta}<_{\Lex}x_0^t x^\gamma$ and $x^{\underline\eta}<_{\Lex}x^{\underline\gamma}$ by Lemma \ref{descLex},  because $x^\epsilon\in \cN(J)$. The same holds for any further reduction and the same argument applies to monomials appearing in $\Supp \big(x_0^t\cdot x^{\gamma'}\cdot \tail{f_{\alpha'}}{}{}\big)$.
\end{proof}

If we just consider $J$ Borel-fixed, with initial degree $m_J$, $\overline{J}$ its saturation and $m$ such that $J_{\geqslant m}=\underline{J}_{\geqslant m}$, then Lemma \ref{lexrid} is no longer true when $m>m_J$: as we already pointed out we need to have $m_J=m$.

\begin{example}\label{jnontagl}
In $\K[x_0,x_1,x_2,x_3]$, consider $\underline{J}=(x_3^2,x_2x_3,x_1x_3,x_2^2)$ and 
\[
\begin{split}
J&{} =(x_3^2)\cdot (x_0,x_1,x_2,x_3)_{\geqslant 2}+(x_3x_2)\cdot (x_0,x_1,x_2,x_3)_{\geqslant 2}+\\
&{} +(x_3x_1)\cdot (x_0,x_1,x_2,x_3)_{\geqslant 2}+ (x_2^2)\cdot (x_0,x_1,x_2,x_3)_{\geqslant 4}.
\end{split}
\]
The ideal $\underline{J}$ is the saturation of $J$, but $J\neq \underline{J}_{\geqslant m}$ for any integer $m$.
Consider a $J$-marked set $G$ and $f_\alpha,f_\beta \in G$ such that $\Ht(f_\alpha)=x_3x_2 x_0^2$ and $\Ht(f_\beta)=x_3 x_1 x_0^2$ and consider $x_2^4\in \Supp \tail{f_{\beta}}{}{}$. Then $S(f_\alpha,f_\beta)=x_1f_\alpha-x_2f_\beta$. If we apply Definition \ref{sminred} and Theorem \ref{ridsm}, we reduce $x_2^4\in \Supp S(f_\alpha,f_\beta)$ by $\xrightarrow{sG}$, pre-multiplying by $x_0^4$. We get that $x_2^4x_0^4$ belongs to $\Supp x_0^4S(f_\alpha,f_\beta)$ and $x_0^4x_2^4=\dec{x_2^2}{x_0^4x_2^2}{\underline{J}}{}$. But $x_2^2x_0^4>_{\Lex} x_2$.
\end{example}

\begin{theorem}\label{nostrobuch}
With the fixed notation, let $I$ be the homogeneous ideal generated by a $J$-marked set $G$. The following statements are equivalent:
\begin{enumerate}[(i)]
\item\label{nostrobuch_i} $I\in \Mf(J)$;
\item\label{nostrobuch_ii} $\forall\ f_\alpha, f_{\alpha'} \in G,\ \exists\ t \text{ s.t. } x_0^t\cdot S(f_\alpha, f_{\alpha'}) \xrightarrow{sG} 0$;
\item\label{nostrobuch_iii} $\forall\ f_\alpha, f_{\alpha'} \in G,\ \exists\ t \text{ s.t. } x_0^t\cdot S(f_\alpha, f_{\alpha'})=x_0^t(x^\gamma f_\alpha-x^{\gamma'}f_{\alpha'}) =\sum a_j x^{\eta_j}f_{\alpha_j}$, with $x^{\eta_j}<_{\Lex}\max_{\Lex}\{x^\gamma,x^{\gamma'}\}$ and $f_{\alpha_j}\in sG$.
\end{enumerate}
\end{theorem}
\begin{proof}
$\emph{(\ref{nostrobuch_i})} \Rightarrow \emph{(\ref{nostrobuch_ii})}$ If $I\in \Mf(J)$, we can apply Corollary \ref{hsconsm} because any $S$-polynomial among elements in $G$ belongs to $I$.

$\emph{(\ref{nostrobuch_i})} \Rightarrow \emph{(\ref{nostrobuch_ii})}$ Straightforward by Lemma \ref{lexrid}.

$\emph{(\ref{nostrobuch_i})} \Rightarrow \emph{(\ref{nostrobuch_ii})}$ Assuming \emph{(\ref{nostrobuch_iii})}, we prove that $\langle V^J\rangle=\langle W^J\rangle$ by an argument analogous to that applied in the proof of Theorem \ref{BuchCrit1}. It is sufficient to prove that $x^\eta\cdot V^J\subseteq \langle V^J\rangle$, for every monomial $x^\eta$. We proceed by induction on the monomials $x^\eta$, ordered according to $\Lex$. The thesis is obviously true for $x^\eta=1$. We then assume that the thesis holds for any monomial $x^{\eta'}$ such that $x^{\eta'}<_{\Lex}x^\eta$.
 
If $\vert\eta\vert> 1$, we can consider any product $x^{\eta}=x^{\eta_1}\cdot x^{\eta_2}$, $x^{\eta_1}$ and $x^{\eta_2}$ non-constant. Since $x^{\eta_i}<_{\Lex}x^\eta, i=1,2$, we immediately obtain by induction 
\[
x^\eta\cdot V^J = x^{\eta_1}\cdot \left(x^{\eta_2}\cdot V^J \right)\subseteq x^{\eta_1}\langle V^J\rangle\subseteq\langle V^J\rangle.
\]

If $\vert \eta\vert=1$, then we need to prove that $x_i\cdot V^J\subseteq\langle V^J\rangle$. Since $x_0 V^J\subseteq  V^J$, it is then sufficient to prove the thesis for $x^\eta=x_i$, $i\geqslant 1$, assuming that the thesis holds for every $x^{\eta'}<_{\Lex} x_i$. We consider $g_\beta=x^\delta f_\alpha \in V^J$, where $\max x^\delta \leq \min x^\alpha$. If $x_ig_\beta$ does not belong to $V^J$, then $\max (x_i\cdot x^\delta) > \min x^\alpha$, so $x_i>\min x^\alpha$ because $\max x^\delta \leq \min x^\alpha$ by construction. In particular, $x_i>\min x^\alpha \geq \max x^\delta$, so $x_i>_\Lex x^\delta$ and it is sufficient to prove the thesis for $x_i f_\alpha$.

We consider an $S$-polynomial $S(f_\alpha,f_{\alpha'})=x_i f_\alpha-x^\gamma f_{\alpha'}$ such that $x^\gamma<_\Lex x_i$. Such $S$-polynomial always exists: for instance, we can consider $x_ix^\alpha=\dec{x^{\alpha'}}{x^{\eta'}}{J}{}$.

By hypothesis there is $t$ such that $x_0^t S(f_\alpha,f_{\alpha'})=x_0^t(x_i f_\alpha-x^{\eta'}f_{\alpha'}) =\sum a_j x^{{\eta'}_j}f_{\alpha_j}$ where $x_0^tx^{\eta'}, x^{{\eta'}_j}$ are lower than $x_i$ w.r.t. $\Lex$. Then $x^{\eta'}f_{\alpha'}$,  $x^{{\eta'}_j}f_{\alpha_j}$  belong to $\langle V^J\rangle$ by induction and by Lemma \ref{x0V}, and we conclude that $x_i f_\alpha \in \langle V^J \rangle$, by Lemma \ref{x0V} again.
\end{proof}

As pointed out in Remark \ref{Buch-proof3} concerning the proof of Theorem \ref{BuchCrit1}, also in the proof of Theorem \ref{nostrobuch} it would be sufficient to assume statement \emph{(\ref{nostrobuch_iii})} only for EK-polynomials. We then have the following result.

\begin{corollary}\label{cor:nostrobuch}
With the same notations of Theorem \ref{nostrobuch},
\[
I\in \Mf(J) \Leftrightarrow\  \forall\ \text{EK-polynomial}\ \exists \ t \text{ s.t. }x_0^t S^{\text{EK}}(f_\alpha,f_{\alpha'})\xrightarrow{sG}0.
\]
\end{corollary}

We have just showed that it is sufficient to work with superminimal reduction for testing if an ideal $I=(G)$ is in $\Mf(J)$. Then, one may think that it is enough to reduce $S$-polynomials among elements in $sG$ (either with $\xrightarrow{V^J_s}$ or with $\xrightarrow{sG}$). The following example clearly shows that this is not true.

\begin{example}
We consider the Borel-fixed ideal 
\[
J= \underline{J}_{\geqslant 2} = (x_3,x_2^2)_{\geqslant 2} = (x_3^2,x_3x_2,x_2^2,x_3x_1,x_3x_0)\subseteq \K[x_0,x_1,x_2,x_3].
\]
 In this case, $sG_J$ contains only two monomials, $x_3x_0$ and $x_2^2$.

If $G$ is any $J$-marked set, then $sG=\{f_{x_3x_0},f_{x_2^2}\}$. The unique $S$-polynomial among superminimal elements is
\[
S(f_{x_3x_0},f_{x_2^2})=x_2^2f_{x_3x_0}-x_3x_0f_{x_2^2}=x_3x_0\cdot\Nf(x_2^2)-x_2^2\cdot \Nf(x_3x_0).
\]
Any monomial appearing in $\Supp \Nf(x_3x_0)$ is in $\cN(J)_2=\K[x_0,x_1,x_2]_2\setminus\{x_2^2\}$. Then any monomial appearing in $\Supp \big(x_2^2\cdot\Nf(x_3x_0)\big)$ is further reduced by $f_{x_2^2}$, obtaining by $\xrightarrow{V^J_s}$ or $\xrightarrow{sG}$
\[
S(f_{x_3x_0},f_{x_2^2})=x_3x_0\cdot\Nf(x_2^2)-\Nf(x_2^2)\cdot\Nf(x_3x_0)=\Nf(x_2^2)\cdot f_{x_3x_0}\rightarrow 0.
\]
Nevertheless, even if the only $S$-polynomial among superminimal generators reduces to 0, we need to impose also other conditions in order to get a $J$-marked basis $G$.

If we consider $sG=\{f_{x_3x_0},f_{x_2^2}\}$ with $f_{x_3x_0}=x_3x_0+x_1^2$ and $f_{x_2^2}=x_2^2$, then for any choice of $f_{x_3^2},f_{x_3x_2},f_{x_3x_1}$, the $S$-polynomial among $f_{x_3x_1}$ and $f_{x_3x_0}$ does not reduce to 0:
\[
S(f_{x_3x_1},f_{x_3x_0})=x_0f_{x_3x_1}-x_1f_{x_3x_0}=\sum_{x^{\alpha_i}\in \cN(J)_2}a_i x^{\alpha_i}x_0-x_1^3.
\]
The monomials $x^{\alpha_i} x_0$ are in $\cN(J)_3$ and are not further reducible (neither by $\xrightarrow{V^J_{s}}$ nor by $\xrightarrow{sG}$). Furthermore, $x_1^3$ does not appear among monomials $m_ix_0$, so it is not canceled. So, for any choice of coefficients in the tail of $f_{x_3x_1}$, we have an $S$-polynomial which is not reducible to 0, and so any $J$-marked set containing $f_{x_3x_0}=x_3x_0+x_1^2$ is not a $J$-marked basis.
\end{example}

The previous example shows that it is not enough to impose condition \emph{(\ref{nostrobuch_ii})} of Theorem \ref{nostrobuch} to $S$-polynomials among elements of $sG$. We need to consider some other $S$-polynomial in order to get a $J$-marked basis.

\begin{theorem}\label{paolo}
With the fixed notations, consider the following sets of $S$-polynomials:
\begin{eqnarray}
\label{eq:nonLinSyz} & L_1 = \left\{ x^{\gamma} f_{\alpha}-x^{\gamma'} f_{\alpha'}\ \big\vert\
x^\alpha,x^{\alpha'} \in sG_J\right\},&\\
\label{eq:LinSyz} & L_2= \left\{x_i f_{\alpha'}-x_0 f_{\alpha}\ \big\vert \ x_i = \min\limits_{j>0} \{x_j \text{ s.t. } x_j \mid x^\alpha\},\ \vert\alpha\vert = \vert\alpha'\vert = m\right\}.&
\end{eqnarray}
Then: 
\[
 I\in\Mf(J)\Longleftrightarrow\ \forall\  S(f_\alpha,f_{\alpha'})\in L_1\cup L_2, \ \exists \ t \text{ s.t. }x_0^t\cdot S(f_\alpha,f_{\alpha'})\xrightarrow{sG}0.
\]
\end{theorem}
\begin{proof}
($\Rightarrow$) If $I$ belongs to $\Mf(J)$, then it is enough to apply Theorem \ref{BuchCrit1}.

($\Leftarrow$) Vice versa, we want to prove that $\langle V^J\rangle=\langle W^J\rangle$, that is $x_i\cdot V^J\subseteq \langle V^J\rangle$ for every $i=0,\dots, n$. We proceed by induction on the variables. By construction we have $x_0\cdot V^J \subseteq \langle V^J\rangle$. We now assume that $(x_0,\dots,x_{i-1})V^J\subseteq \langle V^J\rangle$ and we prove that $x_i\cdot V^J\subseteq \langle V^J\rangle$. Consider $x^\delta f_\beta \in V^J$. The thesis is that $x_i\cdot x^\delta f_\beta$ is contained in $\langle V^J\rangle$. If $x_i x^\delta f_\beta$ does not belong to $V^J$, then $\max (x_i\cdot x^\delta) > \min x^\beta$, so $x_i > \min x^\beta$ because $\max x^\delta \leq \min x^\beta$ by construction. In particular, $x_i > \min x^\beta \geq \max x^\delta$, so that it is sufficient to prove the thesis for $x_i f_\beta$, because by induction then we have $x^\delta x_i f_\beta\in \langle V^J\rangle$. 

Consider $x^\beta= \dec{x^{\underline{\alpha}}}{x^\eta}{\underline{J}}{}$. We have a first case when $x^\eta=1$. Then $x^\beta=x^{\underline{\alpha}}$ and $f_\beta$ belongs to $sG$. We consider $x^{\underline{\alpha}} x_i= \dec{x^{\underline{\alpha}'}}{x^{\eta'}}{\underline{J}}{}$. Observe that since $x_i > \min x^{\underline{\alpha}}$ then $x_i$ does not divide $x^{\eta'}$ and $\max x^{\eta'} <x_i$. Consider $x^{{\alpha}'}=x^{\underline{\alpha}'}\cdot x_0^{t_{\alpha'}}$, so that we can take the polynomial $f_{{\alpha}'}\in sG$. We construct the $S$-polynomial between $f_\beta$ and $f_{{\alpha}'}$
\[
S(f_\beta,f_{{\alpha}'})=x_0^{t_{{\alpha}'}}x_i f_\beta-x^{\eta'}f_{{\alpha}'}
\]
that belongs to $L_1$. Thus, by the hypothesis and by Lemma \ref{lexrid}, there is $k$ such that
\[
x_0^k S(f_\beta,f_{{\alpha}'})=x_0^k (x_0^{t_{{\alpha}'}}x_i f_\beta-x^{\eta'}f_{{\alpha}'})=\sum a_j x^{\eta_j}f_{\alpha_j},
\]
with $x^{\eta_j}<_{\Lex}x_i$ and $f_{\alpha_j}\in sG$. Hence we obtain that both  $x^{\eta_j}f_{\alpha_j}$ and  $x^{\eta'}f_{{\alpha}'}$ belong to $\langle V^J\rangle$ by induction on the variables, and so $x_if_\beta$ belongs to $\langle V^J\rangle$ (by Lemma \ref{x0V}).

We have a second case when $x^\eta=x_0^t$, $t>0$. Then, $\vert\beta\vert=m$ and $f_\beta$ belongs to $sG$. Let $x_ix^\beta = \dec{x^{\underline{\alpha}'}}{x^{\eta'}}{\underline{J}}{}$. If $x_i> \min x^{\underline{\alpha}'}$, then $x^{\eta'}$ is not divisible by $x_i$ and we repeat the argument above. Otherwise, $x_i \leq \min x^{\underline{\alpha}'})$ and $x_i$ does not divide $x^{\eta'}$, so that $x_i=\min x^{\underline{\alpha}'}$ and $x^{\eta'}<_{\Lex} x_i$. Then, we take $x^{\beta'}=\frac{x^\beta}{x_0}\cdot x_i$ that belongs to $G_J$ because it has degree $m$. So, $x_i f_\beta-x_0f_{\beta'}$ belongs to $L_2$ and we repeat the same reasoning above.
 
We now assume the thesis holds for every $f_{\beta'}$ such that $x^{\beta'}= \dec{x^{\underline{\alpha}'}}{x^{\eta'}}{\underline{J}}{}$ with $x^{\eta'}<_{\Lex}x^\eta$. By the base of the induction, we can suppose that $x^\eta \geq_{\Lex}x_1$; so, $f_\beta$ does not belong to $sG$ and hence it has degree $m$. Let $x_j = \min_{l>0} \{x_l \text{ s.t. } x_l\mid x^\beta\}$. Observe that if $x_0$ does not divide $x^\beta$, then $x_j= \min x^\beta$; in this case, we have  $x_i > x_j$ because  $x_i > \min x^\beta$. Anyway, first we suppose that $x_i\leq x_j$; so, $x_j> \min x^\beta$ and $x_0$ divides $x^\beta$. We consider $x^{\beta'}=\frac{x^\beta}{x_0}\cdot x_i$ and the following $S$-polynomial
\[
S(f_{\beta},f_{\beta'})=x_if_{\beta}-x_0f_{\beta'},
\]
that belongs to $L_2$ and we repeat the argument of the previous case.
  
We now assume that $x_i>x_j$ and consider $x^{\beta'}=\frac{x^\beta}{x_j}\cdot x_0=\dec{x^{\underline{\alpha}'}}{x^{\eta'}}{\underline{J}}{}$. Observe that $x^{\eta'}<_{\Lex}x^\eta$ because $x^{\eta'}=\frac{x^\eta}{x_j}\cdot x_0$. We consider the $S$-polynomial:
  \[
  S(f_{\beta'},f_\beta)=x_j f_{\beta'}-x_0 f_\beta
  \]
that belongs to $L_2$; so, by the hypothesis and by Lemma \ref{lexrid}, there is an integer $t$ such that
\begin{equation}\label{spoly}
x_0^t S(f_{\beta'},f_\beta)=x_0^t(x_j f_{\beta'}-x_0 f_\beta)=\sum a_l x^{\eta_l}f_{\alpha_l}
\end{equation}
with $x^{\eta_l}<_{\Lex}x_j$, $f_{\alpha_l}\in sG$. We now multiply \eqref{spoly} by $x_i$. We observe that $x_i f_{\alpha_l}$ belongs to $\langle V^J\rangle$, because $f_{\alpha_l}\in sG$ and by the first two cases. Also $x_i f_{\beta'}$ belongs to $\langle V^J\rangle$ because $x^{\eta'}<_{\Lex}x^\eta$. Moreover, $x_jx_i f_{\beta'}$ belongs to $\langle V^J\rangle$ by induction on the variables. So, $x_i f_\beta$ belongs to $\langle V^J\rangle$ thanks to Lemma \ref{x0V}.
\end{proof}

\begin{remark}
Theorem \ref{paolo} is an improvement, from the computational point of view, of Theorem \ref{nostrobuch}. Indeed, it gives a criterion to establish if $I$ belongs to $\Mf(J)$ in\hfill which\hfill we\hfill  compute\hfill the\hfill reduction\hfill by\hfill $\xrightarrow{s{G}}$\hfill of\hfill a\hfill subset\hfill of\hfill the\hfill whole\hfill set\hfill of\\ $S$-polynomials among elements in $G$.
Actually,  in the proof of Theorem \ref{paolo} we do not need all the $S$-polynomials of the set $L_1$, but only those of type $x_0^{t_{{\alpha}'}}x_if_\beta-x^{\eta'}f_{{\alpha}'}$. This fact can lead us to a further improvement of the efficiency of our algorithms.
\end{remark}

\section{Explicit construction of marked families}\label{ctildate}

In this section, we define an affine scheme whose points correspond to the all the ideals belonging to the $J$-marked family $\Mf(J)$, being $J = \underline{J}_{\geqslant m}$ a truncated Borel-fixed ideal as before, as done in Section \ref{sec:GrobnerStrata} for families sharing the same initial ideal.

\begin{definition}[{Cf. with Definition \ref{def:procedura}}]\label{def:procedureMf}\index{marked family!costruction of a} Let $J = \underline{J}_{\geqslant m} \subset \K[x]$ be a Borel-fixed ideal.
\begin{description}
\item[Step 1] Consider the set of polynomials $\mathcal{G} = \{F_{\alpha}\}_{x^\alpha \in G_J}$
\begin{equation}\label{eq:JbaseC}
F_{\alpha} = x^\alpha + \sum_{x^\beta \in \cN(J)_{\vert\alpha\vert}} C_{\alpha\beta} x^\beta \in \K[C][x]
\end{equation}
where $C$ denotes the whole set of variables $C_{\alpha \beta},\ \forall\ x^\alpha \in G_J,\ \forall\ x^\beta \in \cN(J)_{\vert\alpha\vert}$. Moreover we look at $F_{\alpha}$'s as marked polynomial with $\Ht(F^\alpha) = x^\alpha$.
\item[Step 2] Consider the analogous of $V^J_s$ and $W^J_s$ denoting them by $\mathcal{V}^J_s$ and $\mathcal{W}^J_s$. 
\item[Step 3a] For any pairs $F_{\alpha},F_{\alpha}'$, compute the $J$-reduced forms by $\xrightarrow{\mathcal{V}_s}$ of the $S$-polynomial $S(F_{\alpha},F_{\alpha}')$: $S(F_{\alpha},F_{\alpha}') \xrightarrow{\mathcal{V}_s} H_{\alpha\alpha'}$.
\item[Step 3b] For any EK-pairs $F_{\alpha},F_{\alpha}'$, compute the $J$-reduced forms by $\xrightarrow{\mathcal{V}_s}$ of the EK-polynomial $S^{\textnormal{EK}}(F_{\alpha},F_{\alpha}')$: $S^{\textnormal{EK}}(F_{\alpha},F_{\alpha}') \xrightarrow{\mathcal{V}_s} H^{\textnormal{EK}}_{\alpha\alpha'}$.
\item[Step 4a] Call $\mathfrak{A}_J$ the ideal in $\K[C]$ generated by the coefficients (polynomials of $\K[C]$) of the monomials in the variables $x$ appearing in $H_{\alpha\alpha'}$.
\item[Step 4b] Call $\mathfrak{A}'_J$ the ideal in $\K[C]$ generated by the coefficients (polynomials of $\K[C]$) of the monomials in the variables $x$ appearing in $H^{\textnormal{EK}}_{\alpha\alpha'}$.
\end{description}
\end{definition}

Cioffi and Roggero in \cite[Section 4]{CioffiRoggero} prove that the ideal $\mathfrak{A}_J$ does not depend on the reduction ${\xrightarrow{\ \mathcal{V}^J_{s}}\ }$ and defines the subscheme structure of $\Mf(J)$ in the affine space $\AA^{\vert C\vert}$. By definition $\mathfrak{A}'_J\subseteq\mathfrak A_J$. Anyway, we will prove that $\mathfrak{A}'_J$ and $\mathfrak{A}_J$ are the same ideal, although $\mathfrak{A}_J$ is defined by a set of generators bigger than the set of generators of $\mathfrak{A}'_J$. More precisely, we prove that the ideal $\mathfrak{A}'_J$ contains the coefficients of every $J$-reduced polynomial in $(\mathcal{G}) \subset \K[C][x]$.

\begin{lemma}\label{formule} 
\begin{enumerate}[(i)]
\item\label{formule_i} For every monomial $x^\beta = \dec{x^\alpha}{x^\delta}{J}{} \in J$, we have a formula of type 
\[
x^\beta=\sum a_i x^{\gamma_i} F_{\alpha_i}+H_\beta,
\]
with $a_i\in K[C]$, $x^{\gamma_i} F_{\alpha_i}\in \mathcal{V}^J$,  $x^{\gamma_i}<_{\Lex} x^\delta$ and $\Supp H_\beta\subset \cN(J)$.
\item\label{formule_ii} For every polynomial $x_iF_\alpha\in \mathcal{W}^J\setminus\mathcal{V}^J$, we have a formula of type 
\[
x_i F_\alpha=\sum b_j x^{\eta_j} F_{\alpha_j}+H_{i,\alpha},
\] 
with $b_j\in K[C]$, $x^{\eta_j} F_{\alpha_j}\in \mathcal{V}^J$, $x^{\eta_j}<_{\Lex} x_i$, $\Supp H_{i,\alpha} \subset \cN(J)$ and the coefficients appearing in $H_{i,\alpha}$ belong to $\mathfrak{A}'_J$.
\end{enumerate}
\end{lemma}
\begin{proof}
Statement \emph{(\ref{formule_i})} follows from the existence of $J$-reduced forms obtained by ${\xrightarrow{\mathcal{V}^J_{s}}}$. Statement \emph{(\ref{formule_ii})} follows also from the definition of $\mathfrak A'_J$.
\end{proof}

\begin{proposition} \label{prop:formula}
For every polynomial $x^\delta F_\alpha\in \mathcal{W}^J\setminus \mathcal{V}^J$, we have
\begin{equation} \label{eq:formula}
x^\delta F_\alpha=\sum b_j x^{\eta_j} F_{\alpha_j}+H_{\delta\alpha},
\end{equation}
with $b_j\in K[C]$, $x^{\eta_j} F_{\alpha_j}\in \mathcal{V}^J$, $x^{\eta_j}<_{\Lex} x^\delta$, $\Supp H_{\delta\alpha}\subset \cN(J)$ and the coefficients appearing in $H_{\delta\alpha}$ belong to $\mathfrak{A}'_J$.
\end{proposition}
\begin{proof}
For $\vert\delta\vert=1$ it is enough to use Lemma \ref{formule}\emph{(\ref{formule_ii})}. Assume that $\vert\delta\vert>1$ and that the thesis holds for every $x^{\delta'} <_{\Lex} x^{\delta}$. Let $x_i=\min x^\delta$ and $x^{\delta'}=\frac{x^\delta}{x_i}$, so that $x^{\delta'} F_\alpha$ belongs to $\mathcal{W}^J\setminus \mathcal{V}^J$. 

By the inductive hypothesis, we have $x^{\delta'}F_\alpha=\sum b'_jx^{\eta'_j} F_{\alpha_j}+H_{\delta'\alpha}$, with $x^{\eta'_j} <_{\Lex} x^{\delta'}$. So, multiplying by $x_i$, we obtain $x^{\delta}F_\alpha=\sum b'_jx_i x^{\eta'_j} F_{\alpha_j}+x_i H_{\delta'\alpha}$ and the thesis holds for every polynomial $x_i x^{\eta'_j} F_{\alpha_j}$ that belongs to $\mathcal{W}^J\setminus \mathcal{V}^J$ because $x_ix^{\eta'_j} <_{\Lex} x_ix^{\delta'}=x^\delta$. Then, we substitute such polynomials by formulas of type \eqref{eq:formula} and obtain 
\[
x^{\delta}F_\alpha=\sum b_s x^{\eta_s} F_{\alpha_s}+H'+x_i H_{\delta'\alpha}
\]
where the first sum satisfies the conditions of \eqref{eq:formula} and $H'$ is $J$-reduced with $\Supp H'$ contained in $\cN(J)$ and the coefficients of $H'$ are in $\mathfrak{A}'_J$.

Note that $x_i H_{\delta'\alpha}$ and $H_{\delta'\alpha}$ have the same coefficients belonging to $\mathfrak{A}'_J$, but we do not know if $\Supp (x_i H_{\delta'\alpha}) \subset \cN(J)$. If $x^{\beta'}\in \Supp H_{\delta'\alpha}$ has coefficient $b$ in $H_{\delta'\alpha}$ and $x^\beta=x_i x^{\beta'}$ belongs to $J$, then we can use Lemma \ref{formule}\emph{(\ref{formule_i})} obtaining $b x^\beta=\sum b a_k x^{\gamma_i} F_{\alpha_k}+ b H_\beta$. Moreover, if $x^\beta=\dec{x^{\alpha'}}{x^\epsilon}{J}{}$, then $x^{\gamma_i}<_{\Lex} x^\epsilon <_{\Lex} x_i <_{\Lex} x^\delta$ and all coefficients of $H_\beta$ belong to $\mathfrak{A}'_J$ because they are divisible by $b$. Substituting all such monomials $x^\beta$, we obtain the thesis and $H_{\delta\alpha}$ is $J$-reduced with coefficients in $\mathfrak{A}'_J$, because it is the sum of $J$-reduced polynomials with coefficients in $\mathfrak{A}'_J$.
\end{proof}

\begin{corollary} \label{cor:formula}
Every polynomial of $(\mathcal{G})$ can be written in a unique way as $\sum b_j x^{\eta_j} F_{\alpha_j}+H$, with $b_j\in K[C]$, $x^{\eta_j} F_{\alpha_j}\in \mathcal{V}^J$ and $H$ $J$-reduced. Moreover, we obtain also that the coefficients of $H$ belongs to $\mathfrak{A}'_J$.
\end{corollary}
\begin{proof}
By definition, every polynomial of $(\mathcal{G})$ is a linear combination of polynomials of $\mathcal{V}^J\cap (\mathcal{W}\setminus \mathcal{V})$ with coefficients in $K[C]$ and, by Proposition \ref{prop:formula}, every such polynomial can be written has described in the statement. Hence, we have only to prove the uniqueness of this writing.

Let $\sum b_jx^{\eta_j} F_{\alpha_j}+H=0$ be the difference between two writings of the same polynomial of $(\mathcal{G})$, with $b_j\neq 0$, $x^{\eta_j} F_{\alpha_j}\in \mathcal{V}$ pairwise different and $H$ $J$-reduced. Let $x^{\eta_1} x^{\alpha_1}$ the maximum of the monomials w.r.t. the order for which $x^{\eta_i} x^{\alpha_i}$ is lower than $x^{\eta_j} x^{\alpha_j}$ if $x^{\eta_i} <_{\Lex} x^{\eta_j}$ or $x^{\eta_i} = x^{\eta_j}$ and $x^{\alpha_i} < x^{\alpha_j}$, where $<$ is any order fixed on $G_J$. By definition of $\mathcal{V}^J$, the unique polynomial of $\mathcal{V}^J$ with head term $x^{\eta_1} x^{\alpha_1}$ is $x^{\eta_1} F_{\alpha_1}$. Moreover, the monomial $x^{\eta_1} x^{\alpha_1}$ does not appear with a non-null coefficient in any polynomial of the sum because every other monomial belongs to $\cN(J)$ or is lower than it, by construction and by Lemma \ref{descLex}. Further, $x^{\eta_1} x^{\alpha_1}$ does not belong to $\Supp H$ because $\Supp H\subset \cN(J)$ and $x^{\eta_1} x^{\alpha_1}\in J$. Thus, we obtain a contradiction to the fact that $b_j\neq 0$.
\end{proof}

\begin{corollary} \label{last}
The ideal $\mathfrak{A}'_J$ contains the coefficients of every $J$-reduced polynomial of $(\mathcal{G})$. In particular, $\mathfrak{A}'_J=\mathfrak{A}_J$.
\end{corollary}
\begin{proof} 
Let $F$ be a $J$-reduced polynomial of $(\mathcal{G})$ and let $F=\sum b_jx^{\eta_j} F_{\alpha_j}+H$ as in Corollary \ref{cor:formula}. Since $F$ itself is $J$-reduced, also $F=0+F$ is a formula as described in Corollary \ref{cor:formula} and we obtain that $F=H$, by the uniqueness of this formula. Hence, we have that the coefficients of $F$ and $H$ are the same and are in $\mathfrak A'_J$. The last assertion is due to the definition of $\mathfrak{A}_J$.
\end{proof}

\begin{remark}\label{rem5}
Actually, for every ideal $\widehat{\mathfrak{A}}_J\subseteq \mathfrak{A}_J\subseteq \K[C]$ such that condition \emph{(\ref{formule_ii})} of Lemma \ref{formule} holds, also Corollary \ref{last} holds. We are then allowed to choose different sets of $S$-polynomials of $\mathcal{G}$ in order to obtain generators of the ideal $\mathfrak{A}_J$.
\end{remark}

We know recall the construction of the $J$-marked schemes using matrices to underline the close relation between them and the open affine subsets of the Grasmmannians. By Gotzmann's Persistence Theorem, a specialization $C \rightarrow c\in \K^{\vert C\vert}$ transforms the $J$-marked set $\mathcal{G}$ in a $J$-marked basis $G$ if and only if $\dim_\K (G)_s = \dim_\K J_s$, for every degree $s$. Thus, for each $s$, consider the matrix $\mathcal{A}_s$ whose columns correspond to the terms of degree $s$ in $\K[x]$ and whose rows contain the coefficients of the terms in every polynomial of degree $s$ of type $x^\delta F_\alpha$. Hence, every entry of the matrix $\mathcal{A}_s$ is $1$, $0$ or one of the variables $C$. Let $\mathfrak{a}$ be the ideal of $\K[C]$ generated by the minors of order $\dim_\K J_s + 1$ of $\mathcal{A}_s$, for every $s$.  

\begin{proposition}[{\cite[Lemma 4.2]{CioffiRoggero}}]\label{prop:MftrhoughMatrices}
The ideal $\mathfrak{a}$ is equal to the ideal $\mathfrak{A}$.
\end{proposition}
\begin{proof}
Let $a_s = \dim_\K J_s$. We consider in $\mathcal{A}_s$ the $a_s \times a_s$ submatrix  whose columns corresponds to the terms in $J_s$ and whose rows are given  by the polynomials $x^\beta F_\alpha$ in $\mathcal{V}_s$. Up to a permutation of rows and columns, this submatrix is upper-triangular with 1 on the main diagonal. We may also assume that it corresponds to the first $a_m$ rows and columns in $\mathcal{A}_s$. Then the ideal $\mathfrak{a}$ is generated by the determinants of $(a_s+1)\times (a_s+1)$ submatrices containing that above considered.  Moreover the Gaussian row-reduction of $\mathcal{A}_s$ with respect to the first $a_m$ rows is nothing else than the $\mathcal{V}_s$-reduction  of the $S$-polynomials of the special type considered defining $\mathfrak{A}$.
\end{proof}

As the superminimal reduction uses less polynomials than $\xrightarrow{V^J_s}$, we now exploit it to embed $\Mf(J)$ in an affine subspace of $\AA^{\vert C\vert}$ of lower dimension. 

\begin{definition} If $\mathcal{G}$ is the set of marked polynomials given in \eqref{eq:JbaseC}, we will call \emph{set of superminimal generators},\index{superminimal generators} and denote it by $s\mathcal{G}$,  the subset of $\mathcal{G}$
\begin{equation}\label{eq:JsuperminC}
s\mathcal{G} = \left\{ F_{\alpha} \in \mathcal{G}\ \vert\ \Ht(F_{\alpha}) = x^\alpha \in sG_J\right\}.
\end{equation}
 We will denote by $\widetilde{C} \subset C$ the set of variables appearing in the tails of the polynomials in $s\mathcal{G}$.
\end{definition}

Note that the $J$-marked basis $G$ of every $I\in \Mf(J)$ is obtained by specializing in a suitable way the variables $C$ in $\mathcal{G}$ and that the set of superminimal generators $sG$ of $I$ is obtained in the same way by $s\mathcal{G}$ through the same specialization of the variables $\widetilde{C}$. 

\begin{definition}\label{elim}
Let $x^\alpha \in G_J$ and $t$ be an integer such that $x_0^t \cdot x^\alpha \xrightarrow{s\mathcal{G}} H_\alpha$, with $H_\alpha$ strongly reduced (the integer $t$ exists by Theorem \ref{ridsm}). We can write $H_\alpha= H'_\alpha  + x_0^t \cdot H''_\alpha$, where no monomial appearing in $H'_\alpha$ is divisible by $x_0^t$. We will denote by:
\begin{itemize} 
\item $\mathfrak{B}=\left \{ C_{\alpha \gamma}- \phi_{\alpha \gamma} \text{ s.t. } x^\alpha \in G_J \setminus sG_J,\ x^\gamma \in \cN(J)_{\vert \alpha \vert}\right\}$ the set of the coefficients of $\tail{F_\alpha}{}{}-H''_\alpha$ for every $x^\alpha \in G_J$;
\item $\mathfrak{D}_1\subset \K[\widetilde{C}]$ the set of the coefficients of $H'_\alpha$ for every $x^\alpha \in G_J\setminus sG_J$;
\item $\mathfrak{D}_2$\hfill the\hfill set\hfill of\hfill the\hfill coefficients\hfill of\hfill the\hfill strongly\hfill reduced\hfill polynomials\hfill in\\ $(s\mathcal{G})\K[\widetilde{C}][x]$.
\end{itemize} 
\end{definition}

\begin{theorem}\label{smallaffine} The $J$-marked scheme $\Mf(J)$ is defined by the ideal $\widetilde{\mathfrak{A}}_J=\mathfrak{A}_J \cap \K[\widetilde{C}]$ as  subscheme of the affine space $\AA^{\vert \widetilde{C}\vert}$. 
Moreover $\mathfrak{A}_J=(\mathfrak{B} \cup \mathfrak{D}_1 \cup \mathfrak{D}_2)\K[C]$ and $\widetilde{\mathfrak{A}}_J=(\mathfrak{D}_1 \cup \mathfrak{D}_2)K[\widetilde{C}]$.
\end{theorem}
\begin{proof} 
For the first part it suffices to prove that $\mathfrak{A}_J$ contains $\mathfrak{B}$ and so it contains an element of the type $C_{\alpha\gamma}- \phi_{\alpha \gamma}$, for every $C_{\alpha\gamma}\in C\setminus\widetilde{C}$, where $\phi_{\alpha\gamma}\in \K[\widetilde{C}]$, that allows the elimination of the variables $C_{\alpha\gamma}\in C\setminus \widetilde{C}$.

It is clear by the construction in Definition \ref{elim} that $H_\alpha$ belongs to $\K[\widetilde{C}][x]$ and that both $x_0^t \cdot \tail{F_\alpha}{}{}$ and $H_\alpha$ are strongly reduced. Thus their difference $x_0^t \cdot \tail{F_\alpha}{}{}-H_\alpha$ is strongly reduced and moreover it belongs to $(\mathcal{G})$, because $x_0^t \cdot \tail{F_\alpha}{}{}-H_\alpha= -x_0^t \cdot F_\alpha+(x_0^t \cdot x^\alpha-H_\alpha)$. Hence, by Corollary \ref{last}, its coefficients belong to $\mathfrak{A}_J$ and in particular the coefficient of $x_0^t \cdot x^\gamma$ is of the type $C_{\alpha \gamma}- \phi_{\alpha \gamma}$, with $\phi_{\alpha \gamma}\in \K[\widetilde{C}]$. Then $\mathfrak{A}_J \supseteq \mathfrak{B}$ and $\mathfrak{A}_J$ is generated by $\mathfrak{B} \cup \widetilde{\mathfrak{A}}_J$.

To prove the second part, it is sufficient to show that $\mathfrak{A}_J \cap \K[\widetilde{C}]=(\mathfrak{D}_1 \cup \mathfrak{D}_2)\K[\widetilde{C}]$.

($\supseteq$) Taking the coefficients in $x_0^t \cdot \tail{F_\alpha}{}{}-H_\alpha$ of monomials that are not divisible by $x_0^t$, we see that $\mathfrak{A}_J$ contains the coefficients of $H'_\alpha$. Then $\mathfrak{A}_J \cap \K[\widetilde{C}]  \supseteq \mathfrak{D}_1$, because $H'_\alpha \in \K[\widetilde{C}][x]$. Moreover we recall that $\mathfrak{A}_J$ is made by all the coefficients in the polynomials of $(\mathcal{G})$ that are strongly reduced. Indeed, $\mathfrak{A}_J$ is made by all the coefficients of the polynomials of $(\mathcal{G})$ that are $J$-reduced. But the degree of the monomials in the variables $x$ of every polynomial in $(\mathcal{G})$ is $\geqslant m$ and then \lq\lq$J$-reduced\rq\rq \ is equivalent to \lq\lq$\underline{J}$-reduced\rq\rq, that it is strongly reduced. Then $\mathfrak{A}_J \cap \K[\widetilde{C}]  \supseteq \mathfrak{D}_2$, because $(s\mathcal{G})\K[\widetilde{C}][x] \subset (\mathcal{G})$.

($\subseteq$)  For every polynomial $F\in \K[C,x]$, let us denote by $F^\phi$  the polynomial in $\K[\widetilde{C},x]$ obtained substituting every $C_{\alpha \gamma}\in C \setminus \widetilde{C}$ by $\phi_{\alpha \gamma}$. 
Observe that for ever $x^\alpha \in G_J$ we have $x_0^t\cdot F_\alpha ^\phi= x_0^t(x^\alpha-H''_\alpha )+H'_\alpha$ and moreover $x_0^t(x^\alpha-H'')+H'_\alpha \in (s\mathcal{G})\K[\widetilde{C},x]$. It remains to prove that every element $w\in \mathfrak{A}_J \cap K[\widetilde{C}]$ can be obtained modulo $\mathfrak{D}_1$ as a coefficient in some strongly reduced polynomial of the ideal $(s\mathcal{G})\subset \K[\widetilde{C}]$. We know that $w$ is a coefficient in a strongly reduced polynomial $D\in (\mathcal{G})$.
  
If $D=\sum D_\alpha F_\alpha\in (\mathcal{G})$, then for a suitable $t$, $$x_0^t\cdot D^\phi =\sum D^\phi_\alpha\cdot \left( x_0^t\cdot (x^\alpha -H''_\alpha)+H'_\alpha\right) \in (s\mathcal{G})\K[\widetilde{C}][x]$$ and $w$ is still one of the coefficients of $D^\phi$ because it does not contain any variable in $C\setminus \widetilde{C}$ and so it remains unchanged. Moreover if $D$ is strongly reduced, also $D^\phi$ is strongly reduced and so $w\in (\mathfrak{D}_1 \cup \mathfrak{D}_2)\K[\widetilde{C}]$.
\end{proof}

\begin{proposition}\label{formuleconsuperminimals}Let $\mathfrak{U}\subseteq \widetilde{\mathfrak{A}}_J$ be any ideal in $\K[\widetilde{C}]$ such that:
\begin{enumerate}[(i)]
\item\label{formuleconsuperminimals1} for every monomial $x^\beta = \dec{x^{\underline{\alpha}}}{x^\delta}{\underline{J}}{} \in J$, there exists $t$ such that we have a formula of type 
\[
x_0^t\cdot x^\beta=\sum b_ix^{\eta_i} F_{\alpha_i}+H_\beta,
\] 
with $a_i\in \K[\widetilde{C}]$, $F_{\alpha_i}\in s\mathcal{G}$,  $x^{\eta_i}<_{\Lex}x^\delta$, $x^{\underline{\eta}_j+\underline{\alpha}_j}= \dec{x^{\underline{\alpha}_j}}{x^{\underline{\eta}_j}}{\underline{J}}{}$ and $H_\beta=x_0^t\cdot H_\beta'+H_\beta''$, with $\Supp H_\beta \subset \cN(J)$, $x_0^t$ does not divide $H_\beta''$ and the coefficients of $H_\beta''$ belong to $\mathfrak{U}$;
\item\label{formuleconsuperminimals2} for every polynomial $F_\alpha\in s\mathcal{G}$ and for every $x_i>\min x^{\underline \alpha}$ there exists $t$ such that we have a formula of type
\[
x_0^t x_i F_\alpha=\sum b_j x^{\eta_j}F_{\alpha_j}+H_{i,\alpha}
\]
where $b_j\in \K[\widetilde{C}]$, $F_{\alpha_j}\in s\mathcal{G}$, $x^{\eta_j}<_{\Lex}x_i$, $x^{\underline{\eta}_j+\underline{\alpha}_j}= \dec{x^{\underline{\alpha}_j}}{x^{\underline{\eta}_j}}{\underline{J}}{}$, $\Supp H_{i,\alpha}\subseteq \cN(J)$ and the coefficients of $H_{i,\alpha}$ belongs to $\mathfrak{U}$.
\end{enumerate}
Then $\mathfrak{U}=(\mathfrak{D}_1\cup\mathfrak{D}_2)$.
\end{proposition}
\begin{proof}
Thanks to \emph{(\ref{formuleconsuperminimals1})}, we immediately have that $\mathfrak{D}_1\subseteq \mathfrak{U}$.

For the inclusion $\mathfrak{D}_2\subseteq \mathfrak{U}$, observe that if \emph{(\ref{formuleconsuperminimals1})} and \emph{(\ref{formuleconsuperminimals2})} hold for $\mathfrak{U}$, then we can use the same arguments of Proposition \ref{prop:formula} and obtain that for every $F_\alpha \in s\mathcal{G}$, for every $x^\delta$, there exists $t$ such that
\begin{equation}\label{eqsuper}
x_0^t x^\delta F_\alpha=\sum b_j x^{\eta_j}F_{\alpha_j}+H
\end{equation}
with $b_j\in \K[\widetilde{C}]$, $F_{\alpha_j}\in s\mathcal{G}$, $x^{\eta_j}<_{\Lex}x^\delta$, $x^{\underline{\eta}_j+\underline{\alpha}_j}= \dec{x^{\underline{\alpha}_j}}{x^{\underline{\eta}_j}}{\underline{J}}{}$, $\Supp H_{\delta,\alpha} \subseteq \cN(J)$ and the coefficients of $H_{\delta,\alpha}$ belong to $\mathfrak{U}$.

We can also prove the uniqueness of such a rewriting: thanks to the uniqueness of the $\underline{J}$-canonical decomposition (Lemma \ref{lem:monomialDecomposition}), the polynomials $x^{\eta_j}F_{\alpha_j}$ that can appear in \eqref{eqsuper} have pairwise different head terms. So an analogous of Corollary \ref{cor:formula} holds for this setting. Thanks to this uniqueness, as in Corollary \ref{last}, we get the non trivial inclusion of the thesis.
\end{proof}

Proposition \ref{formuleconsuperminimals} is very important from the computational point of view: indeed, 
different choices of sets of $S$-polynomials to reduce give different sets of generators for $\widetilde{\mathfrak{A}}_J$. For instance we can get a set of generators for $\widetilde{\mathfrak{A}}_J$ starting from EK-polynomials among polynomials in $\mathcal{G}$ or starting from $\mathcal{L}_1$ and $\mathcal{L}_2$ in $\K[C][x]$, corresponding to $L_1$ and $L_2$ as defined in Theorem \ref{paolo}. However, a good choice of the set of $S$-polynomials can strongly influence the efficiency of an algorithm computing equations for $\Mf(J)$. 

\medskip

As seen in Section \ref{sec:GrobnerStrata}, a Gr\"obner stratum $\St{\sigma}(J)$ can be isomorphically projected in its Zariski tangent space at the origin $T_0\big(\St{\sigma}(J)\big)$ and moreover if the origin is a smooth point, then the stratum is isomorphic to this tangent space. In general, if we do not consider a term ordering we cannot  project isomorphically $\Mf(J)$ into $T_0\big(\Mf(J)\big)$, but in any case, the dimension of this tangent space plays an interesting role in the following theorem.

\begin{remark}\label{tangente}
Let $\mathfrak{L}(J)$ be the ideal generated in $\K[C]$ by the linear components of the generators of $\mathfrak{A}_J$. Then, the Zariski tangent space $T_0\big(\Mf(J)\big)$ of the $J$-marked scheme at the origin can be naturally identified to the linear space of $\AA^{\vert C\vert}$ defined as the set of zeros of $\mathfrak{L}(J)$.
\end{remark}

We now prove the analogous of Theorem \ref{th:saturato} for Gr\"obner strata in the case of marked families. Using at the same time several truncations of a saturated Borel-fixed ideal $\underline{J}$, we introduce the following notation: 
\begin{itemize}
\item $s\mathcal{G}^{(m)}$ will denote the superminimal generators associated to $\underline{J}_{\geqslant m}$ and $\widetilde{C}^{(m)}$ the corresponding variables;
\item $\widetilde{\mathfrak{A}}^{(m)}$ will denote the ideal defining the affine subscheme $\Mf(\underline{J}_{\geqslant m})$ in the ring $\K[\widetilde{C}^{(m)}]$ (as in Theorem  \ref{smallaffine}). 
\end{itemize}

\begin{theorem}[Cf. with Theorem \ref{th:saturato}]\label{schisom}
Let $\underline{J}$ be a saturated Borel-fixed ideal and let $m$ be any integer. 
With the previous notations, the followings hold:
\begin{enumerate}[(i)]
\item\label{schisom_i} $\Mf(\underline{J}_{\geqslant m-1})$ is a closed subscheme of $\Mf(\underline{J}_{\geqslant m})$ cut out by a suitable linear space.
\item\label{schisom_iiii} Let $N$ be the number of monomials $x^\alpha \in G_{\underline{J}}$ of degree $m+1$ divisible by $x_1$ and $M =\vert G_{\underline{J}}\cap \K[x]_{\leqslant m-1}\vert$; then, 
\[
\dim_\K T_0\big(\Mf(\underline{J}_{\geqslant m})\big)\geqslant \dim_\K T_0\big(\Mf(\underline{J}_{\geqslant m-1})\big)+ NM.
\]
\item \label{schisom_iiiii} $\Mf(\underline{J}_{\geqslant m-1})\simeq\Mf(\underline{J}_{\geqslant m})$ if and only if either $\underline{J}_{\geqslant m-1}= \underline{J}_{\geqslant m}$ or no monomial of degree $m+1$ in $G_{\underline{J}}$ is divisible by $x_1$.
\end{enumerate}
In particular: 
\begin{equation}
\Mf(\underline{J}_{\rho-1}) \simeq \Mf(\underline{J}_{\reg(\underline{J})-1}) \simeq \Mf(\underline{J}_{\reg(\underline{J})})\simeq \Mf(\underline{J}_{r}),
\end{equation}
where $\rho$ is the maximal degree of monomials divisible by $x_1$ in $G_{\underline{J}}$ and $r$ is the Gotzmann number of the Hilbert polynomial $p(t)$ of $\K[x]/\underline{J}$.

Moreover, $\Mf(\underline{J}_{\reg(\underline{J})})$  can be embedded in an affine space of dimension $\vert G_{\underline{J}}\vert \cdot p\big( \reg(\underline{J})\big)$, and the same holds for $\Mf(\underline{J}_{\geqslant m})$, for every $m\geqslant \reg(\underline{J})$.
\end{theorem}
\begin{proof} 
\emph{(\ref{schisom_i})} Thanks to Theorem \ref{smallaffine}, a marked scheme is defined by an ideal generated by polynomials of $\K[\widetilde{C}]$ that are constructed using only the superminimals.
So, now it is enough to prove that the set of superminimals $s\mathcal{G}^{(m-1)}$ corresponds to $s\mathcal{G}^{(m)}$ modulo a subset of the variables $\widetilde{C}^{(m)}$, in the following sense. 

Consider $x^\alpha \in sG_{\underline{J}_{\geqslant m-1}}$. If $\vert \alpha \vert\geqslant m$, then $x^\alpha$ belongs to $sG_{\underline{J}_{\geqslant m}}$ and we can identify $F^{(m)}_\alpha \in s\mathcal{G}^{(m)}$ and $F^{(m-1)}_\alpha \in s\mathcal{G}^{(m-1)}$ (and in particular the variables in their tails: $\widetilde{C}^{(m)}_{\alpha \gamma}=\widetilde{C}^{(m-1)}_{\alpha \gamma}$).

If $\vert\alpha\vert=m-1$, then we can consider the corresponding superminimal element $F^{(m)}_\beta \in s\mathcal{G}^{(m)}$, with $x^\beta=x_0\cdot x^\alpha$. Then we identify the variable $\widetilde{C}_{\beta\delta'}^{(m)}$, which is the coefficient of a monomial in $\Supp F_\beta^{(m)}$ of kind $x^{\delta'}=x_0\cdot x^\delta$, with the variable $\widetilde{C}_{\alpha \delta}^{(m-1)}$ which is the coefficient of the monomial $x^\delta$ in $\Supp F_\alpha^{(m-1)}$. 

We repeat this identifications for all $x^\alpha \in sG_{\underline{J}_{\geqslant m-1}}$ and we denote by $\overline{C}^{(m)}$ the subset of $\widetilde{C}^{(m)}$ containing the variables non-identified with variables of $\widetilde{C}^{(m-1)}$, that is the variables appearing as coefficients of monomials not divisible by $x_0$ in the tails of polynomials in $s\mathcal{G}^{(m)}\setminus s\mathcal{G}^{(m-1)}$. Now, every polynomial in $s\mathcal{G}^{(m)}\mod (\overline{C}^{(m)})$ either belongs to $s\mathcal{G}^{(m-1)}$ or is a polynomials of $s\mathcal{G}^{(m-1)}$ multiplied by $x_0$. Thanks to Theorem \ref{smallaffine}, we have that 
\begin{equation*}
\widetilde{\mathfrak{A}}^{(m)}+\left(\overline{C}^{(m)}\right)\simeq \widetilde{\mathfrak{A}}^{(m-1)}.
\end{equation*}

\emph{(\ref{schisom_iiii})} We now consider $x^\gamma \in G_{\underline{J}}$, $\vert \gamma\vert=m+1$, $x^\gamma$ divisible by $x_1$. We define $x^\beta=x^\gamma/x_1$; observe that $x^\beta \notin \underline{J}$. Furthermore, $x^\beta$ is not divisible by $x_0$, otherwise $x^\gamma$ would be too.

Then, for every $x^{\underline{\alpha}} \in G_{\underline{J}}$ with $\vert\underline{\alpha}\vert\leqslant m-1$, there is $F_{\alpha}=x^{\underline{\alpha}} x_0^{m-\vert\underline{\alpha}\vert}-\tail{F_{\alpha}}{}{} \in s\mathcal{G}^{(m)}$ such that $x^\beta \in \Supp \tail{F_{\alpha}}{}{}$. We focus on the coefficient $\widetilde{C}_{\alpha\beta}^{(m)}$ of $x^\beta$. Since $x^\beta$ is not divisible by $x_0$, $\widetilde{C}_{\alpha\beta}^{(m)}$ cannot be identified with a coefficient appearing in $F_{\alpha}^{(m-1)}=x^{\underline{\alpha}} x_0^{m-\vert\underline{\alpha}\vert-1}-\tail{F_{\alpha}^{(m-1)}}{}{}\in s\mathcal{G}^{(m-1)}$. So $\widetilde{C}_{\alpha\beta}^{(m)}$ belongs to the subset of variables $\overline{C}^{(m)}$ defined in the proof of \emph{(\ref{schisom_i})}.

We now use the construction of $T_0\big(\Mf(\underline{J}_{\geqslant m})\big)$ of Remark \ref{tangente}. If we think about syzygies of the ideal $\underline{J}_{\geqslant m}$, we can see that in any $S$-polynomial, $F_{\alpha}^{(m)}$ is multiplied by a monomial $x^\delta$ divisible by $x_i$, $i>0$. In particular, $x^\delta\cdot x^\beta\in \underline{J}_{\geqslant m}$; indeed, if $x_i=x_1$ we are done by construction, otherwise we apply the Borel-fixed property because $\frac{x^\gamma}{x_i}\cdot x_1 \cdot x^\beta$ belongs to $\underline{J}$. This means that the coefficient $\widetilde{C}_{\alpha\beta}^{(m)}$ does not appear in any equation defining $T_0\big(\Mf(\underline{J}_{\geqslant m})\big)$.

Applying this argument to the $N$ monomials in $G_{\underline{J}}$ of degree $m+1$ which are divisible by $x_1$ and to the $N$ monomials in $G_{\underline{J}}$ of degree $\leqslant m-1$, we obtain the result.

\emph{(\ref{schisom_iiiii})} If $\underline{J}_{\geqslant m}=\underline{J}_{\geqslant m-1}$, obviously $\Mf(\underline{J}_{\geqslant m})=\Mf(\underline{J}_{\geqslant m-1})$.  
We now assume that $\underline{J}_{\geqslant m}\neq \underline{J}_{\geqslant m-1}$ and no monomial of degree $m+1$ in the monomial basis of $\underline{J}$  is divisible by $x_1$; we prove that every polynomial in $s\mathcal{G}^{(m)}$ either belong to $s\mathcal{G}^{(m-1)}$ or it is the product of the \lq\lq corresponding\rq\rq\  polynomial in $s\mathcal{G}^{(m-1)}$ by $x_0$.

If $x^\alpha \in sG_{\underline{J}_{\geqslant m-1}}$ and $\vert\alpha\vert \geqslant m$, then $F^{(m)}_\alpha \in s\mathcal{G}^{(m)}$ and $F^{(m-1)}_\alpha \in s\mathcal{G}^{(m-1)}$ have the same shape and we can identify them letting $\widetilde{C}^{(m)}_{\alpha \gamma}=\widetilde{C}^{(m-1)}_{\alpha \gamma}$, as done in the proof of \emph{(\ref{schisom_i})}. If $\vert\alpha\vert =m-1$, then  $x^\beta=x_0\cdot x^\alpha \in sG_{\underline{J}_{\geqslant m}}$ and all the monomials in the support of $x_0\cdot F^{(m-1)}_\alpha$  appear in the support of $F^{(m)}_\beta$ (and we identify their coefficients as above). In the support of $F^{(m)}_\beta$ there are also some more monomials that are not divisible by $x_0$. We will prove now that the coefficients of these last monomials indeed belong to $\widetilde{\mathfrak{A}}^{(m)}$.

Consider the monomial $x_0\cdot x_1\cdot x^\alpha$. If we perform its reduction using $s\mathcal{G}^{(m)}$, the first step of reduction will lead to
\[
x_0\cdot x_1\cdot x^\alpha \xrightarrow{s\mathcal{G}^{(m)}} x_1 \tail{F^{(m)}_\beta}{}{}.
\]
Let $x^\gamma$ be a monomial of $\Supp \tail{F^{(m)}_\beta}{}{}$. If $x_1\cdot x^\gamma \in \underline{J}_{\geqslant m}$, then $x_1\cdot x^\gamma = \dec{x^{\underline{\alpha}'}}{x^\eta}{\underline{J}}{}$, with $x^{\underline{\alpha}'}\in G_{\underline{J}}$ and $x^\eta<_{\Lex} x_1$. 
If $x^\eta=1$, then $\vert\underline{\alpha}'\vert=m+1$ and $x^{\underline{\alpha}'}$ is divisible by $x_1$, against the hypothesis.
Then $x^\eta=x_0^{\overline{t}}$, with $t>0$, and so the monomial $x_1\cdot x^\gamma\in \underline{J}_{\geqslant m}$ is actually divisible by $x_0$.
If $x_1\cdot x^\gamma \in \cN(\underline{J}_{\geqslant m})$, then this monomial is not further reducible, so that its coefficient belongs to $\widetilde{\mathfrak{A}}^{(m)}$.

Vice versa, by absurd suppose now that $\underline{J}_{\geqslant m-1}\neq \underline{J}_{\geqslant m}$ and that there exists $x^\alpha \in G_{\overline{J}}$ divisible by $x_1$, $\vert\alpha\vert=m+1$.
Using \emph{(\ref{schisom_iiii})}, we have that $T_0\big(\Mf(\underline{J}_{\geqslant m-1})\big)\not\simeq T_0\big(\Mf(\underline{J}_{\geqslant m})\big)$\hfill because\hfill $\dim_\K T_0\big(\Mf(\underline{J}_{\geqslant m-1})\big) < \dim_\K T_0\big(\Mf(\underline{J}_{\geqslant m})\big)$,\hfill and\hfill so\\ $\Mf(\underline{J}_{\geqslant m-1}) \not\simeq \Mf(\underline{J}_{\geqslant m})$.

For the last part of the statement, note that if $\rho$ is the maximal degree of a monomial divisible by $x_1$ in the monomial basis of $\underline{J}$, for every $m \geqslant \rho$, applying iteratively \emph{(\ref{schisom_iiiii})} we obtain
\[
\Mf(\underline{J}_{\geqslant \rho-1})\simeq \Mf(\underline{J}_{\geqslant m}).
\]
If especially $m \geqslant \reg(\underline{J})$, the Hilbert function and the Hilbert polynomial $p(t)$ of $\K[x]/\underline{J}$ surely coincide, hence $\Mf(\underline{J}_{\geqslant m})$ can be embedded in an affine space of dimension 
\[
\vert \widetilde{C}^{(m)}\vert= \sum_{x^\alpha \in sG_{\underline{J}_{\geqslant m}}} p(\vert\alpha\vert)
\]
 and in this case every monomial in $sG_{\underline{J}_{\geqslant m}}$ has degree $m$.
\end{proof}

\begin{remark}
By Theorem \ref{schisom}, we can embed $\Mf(\underline{J}_{\geqslant m})$ in an affine space of dimension $\vert\widetilde{C}^{(m)}\vert$, for every $m\geqslant \rho-1$. Anyway, we cannot always compute the dimension of this affine space using the Hilbert polynomial $p(t)$, since $m$ may be strictly lower than the regularity of the Hilbert function. 
For $m< \reg(\underline{I})$, we have that
\[
\vert\widetilde{C}^{(m)}\vert=\sum_{x^{\alpha}\in sG_{\underline{J}_{\geqslant m}}}\vert \cN(\underline{J}_{\geqslant m})_{\vert\alpha\vert}\vert.
\]
\end{remark}

\subsection{The pseudocode description of the algorithm}

We will now expose the pseudocode of the algorithm for computing the equations defining the affine scheme that describes a $J$-marked family $\Mf(J)$ for a truncated Borel-fixed ideal $J = \underline{J}_{\geqslant m}$, mainly based on Theorem \ref{paolo}.

\begin{algorithm}[H]
\capstart
\begin{algorithmic}
\STATE $\textsc{MinimalGenerators}(J)$
\REQUIRE $J$, a truncation of a saturated Borel-fixed ideal, i.e. $J = \underline{J}_{\geqslant m}$ for some $m$.
\ENSURE the set $G_J$ of minimal generators of $J$.
\end{algorithmic}

\bigskip

\begin{algorithmic}
\STATE $\textsc{SuperminimalGenerators}(J)$
\REQUIRE $J$, a truncation of a saturated Borel-fixed ideal, i.e. $J = \underline{J}_{\geqslant m}$ for some $m$.
\ENSURE the set $sG_J$ of superminimal generators of $J$.
\end{algorithmic}

\bigskip

\begin{algorithmic}
\STATE $\textsc{Reduce}(H,\textsf{nfs})$
\REQUIRE $H$, a polynomial in $\K[C][x]$;
\REQUIRE $\textsf{nfs}$, a set of marked polynomials $\Ht(F_\alpha) - \tail{F_{\alpha}}{}{}$ of $\K[C][x]$ such that $\Ht(F_{\alpha}) \in \K[x]$ and $\tail{F_{\alpha}}{}{} = \Nf\big(\Ht(F_\alpha)\big)$.
\ENSURE the polynomial $\overline{H}$ computed by replacing each monomial $x^\delta \in \Supp H$ with $x^\eta \tail{F_{\alpha}}{}{}$ with $F_{\alpha} \in \textsf{nfs}$.
\end{algorithmic}
\end{algorithm}

\begin{algorithm}[H]
\begin{algorithmic}
\STATE $\textsc{SuperminimalReduction}(H,s\mathcal{G})$
\REQUIRE $H$, a polynomial in $\K[C][x]$;
\REQUIRE $s\mathcal{G}$, the set of superminimal generators for some marked family.
\ENSURE $\overline{H}$ such that there exists $t$ for which $x_0^t H \xrightarrow{s\mathcal{G}} \overline{H}$.
\end{algorithmic}

\bigskip

\begin{algorithmic}
\STATE $\textsc{SuperminimalSyzygies}(J)$
\REQUIRE $J$, a truncation of a saturated Borel-fixed ideal, i.e. $J = \underline{J}_{\geqslant m}$ for some $m$.
\ENSURE the set of syzygies between pairs of superminimal generators of $J$ corresponding to the set $L_1$ of Theorem \ref{paolo}.
\end{algorithmic}

\bigskip

\begin{algorithmic}
\STATE $\textsc{Coeff}(H,x^\beta)$
\REQUIRE $H$, a polynomial in $\K[C][x]$.
\REQUIRE $x^\beta$, a monomial in $\K[x]$.
\ENSURE the coefficient of the monomial $x^\beta$ in $H$ (obviously 0 if $x^\beta \notin \Supp H$).
\end{algorithmic}
\caption{Auxiliary methods for the algorithm computing the affine scheme that describes a marked family.}
\label{alg:auxMf}
\end{algorithm}

\begin{algorithm}[H]
\capstart
\label{alg:markedFamily}
\begin{algorithmic}[1]
\STATE $\textsc{MarkedFamily}(J)$
\REQUIRE $J$, a truncation of a saturated Borel-fixed ideal, i.e. $J = \underline{J}_{\geqslant m}$ for some $m$.
\ENSURE the ideal of the affine scheme describing $\Mf(J)$.
\STATE $sG_J \leftarrow \textsc{SuperminimalGenerators}(J)$;
\STATE $s\mathcal{G}\leftarrow \emptyset$;
\FORALL{$x^\alpha \in sG_J$}
\STATE $F_{\alpha} \leftarrow x^\alpha$;
\FORALL{$x^\beta \in \cN(J)_{\vert\alpha\vert}$}
\STATE $F_{\alpha} \leftarrow F_{\alpha} + C_{\alpha\beta} x^\beta$;
\ENDFOR
\STATE $s\mathcal{G} \leftarrow s\mathcal{G} \cup \{F_{\alpha}\}$;
\ENDFOR
\STATE $G \leftarrow \textsc{MinimalGenerators}(J) \setminus sG_J$;
\STATE $\textsf{knownNF} \leftarrow s\mathcal{G}$;
\end{algorithmic}
\end{algorithm}

\begin{algorithm}[H]
\caption[Algorithm for computing the affine scheme that describe a marked family.]{The algorithm for computing the affine scheme that describe a marked family.}
\label{alg:markedFamily2}
\begin{algorithmic}[1]\setcounter{ALC@line}{12}
\STATE $\textsf{equations} \leftarrow \emptyset$;
\WHILE{$G \neq \emptyset$}
\STATE $x^\alpha \leftarrow \min_{\RevLex} G$;
\STATE $x_i \leftarrow \min_{j>0} \{ x_j \mid x^\alpha\}$;
\STATE $x^\gamma \leftarrow \frac{x_0}{x_i}x^\alpha$; \hfill \COMMENT{This is a sygygy of the set $L_2$ of Theorem \ref{paolo}}
\STATE $H \leftarrow \textsc{Reduce}\big(x_i \Nf(x^\gamma),\textsf{knownNF}\big)$;
\STATE $(Q,R) \leftarrow {}\text{polynomials such that } H = Q \cdot x_0 + R$;
\FORALL{$x^\delta \in \Supp R$}
\STATE $\textsf{equations} \leftarrow \textsf{equations} \cup \{\textsc{Coeff}(R,x^\delta)\}$;\hfill \COMMENT{We are imposing $R = 0$}
\ENDFOR
\STATE $\textsf{knownNF} \leftarrow \textsf{knownNF} \cup \left\{x^\alpha - Q \right\}$; \hfill \COMMENT{Because $H = \Nf(x_0 x^\gamma) = x_0\Nf(x^\gamma)$}
\STATE $G \leftarrow G \setminus \{x^{\alpha}\}$;
\ENDWHILE
\STATE $\textsf{syzygies} \leftarrow \textsc{SuperminimalSyzygies}(J)$;
\FORALL{$x^{\gamma}\mathbf{e}_{\alpha} - x^{\gamma'}\mathbf{e}_{\alpha'} \in \textsf{syzygies}$}
\STATE $S(F_{\alpha},F_{\alpha'}) \leftarrow x^{\gamma}F_{\alpha} - x^{\gamma'}F_{\alpha'}$;
\STATE $H \leftarrow \textsc{SuperminimalReduction}\big(S(F_{\alpha},F_{\alpha'}),s\mathcal{G}\big)$;
\FORALL{$x^\delta \in \Supp H$}
\STATE $\textsf{equations} \leftarrow \textsf{equations} \cup \{\textsc{Coeff}(H,x^\delta)\}$;
\ENDFOR
\ENDFOR
\RETURN $\big\langle\textsc{equations}\big\rangle$;
\end{algorithmic}
\end{algorithm}
It is very useful also to have a method that computes the dimension of the tangent space at the origin of a marked family (i.e. the number of monomials in the tails of the superminimal generators) avoiding to generate the complete equations of the ideal of the family. For instance if we know a lower bound of the dimension of the tangent space (determined by geometric arguments or for any other reason) and the marked family realizes this bound, we can conclude directly that the marked family is an affine space and that there are no relations among the variables $C$.

\begin{algorithm}[H]
\caption[Algorithm computing the dimension of the tangent space at the origin of a marked family.]{The algorithm computing the dimension of the tangent space at the origin of a marked family.}
\label{alg:dimTang}
\begin{algorithmic}[1]
\STATE $\textsc{TangentSpaceDimension}(J)$
\REQUIRE $J$, a truncation of a saturated Borel-fixed ideal, i.e. $J = \underline{J}_{\geqslant m}$ for some $m$.
\ENSURE the dimension of the tangent space at the origin of $\Mf(J)$.
\STATE $sG_J \leftarrow \textsc{SuperminimalGenerators}(J)$;
\STATE $N \leftarrow 0$;
\FORALL{$x^\alpha \in sG_J$}
\STATE $N \leftarrow N + \left\vert \cN(J)_{\vert\alpha\vert}\right\vert$;
\ENDFOR
\RETURN $N$;
\end{algorithmic}
\end{algorithm}

\begin{algorithm}[H]
\caption[Algorithm computing the Gr\"obner stratum of a gen-segment ideal.]{The algorithm computing the Gr\"obner stratum of a gen-segment ideal.}
\label{alg:grobnerStrata}
\begin{algorithmic}[1]
\STATE $\textsc{Gr\"obnerStratum}(J,\sigma)$
\REQUIRE $J$, a truncation of a saturated Borel-fixed ideal, i.e. $J = \underline{J}_{\geqslant m}$ for some $m$.
\REQUIRE $\sigma$, a term ordering such that $J$ is a gen-segment ideal w.r.t. $\sigma$.
\ENSURE the ideal of the affine scheme describing $\St{\sigma}(J)$.
\STATE $\textsf{ideal} \leftarrow \textsc{MarkedFamily}(J)$;
\STATE $\textsf{equations} \leftarrow {}$ generators of $\textsf{ideal}$;
\FORALL{$f \in \textsf{equations}$}
\IF{$f \notin (C)^2$}
\STATE $C_{\alpha\beta} = \max_{\ell} \Supp f$;\hfill\COMMENT{$\ell$ is the positive grading induced on $C$ by $\sigma$ (Definition \ref{def:lambda})}
\STATE $\phi_{\alpha\beta} = \frac{f}{\textsc{Coeff}(C_{\alpha\beta},f)} + C_{\alpha\beta}$;
\STATE $\textsf{ideal} \leftarrow \textsc{Substitute}(C_{\alpha\beta} \leftarrow \phi_{\alpha\beta},\textsf{ideal})$;
\STATE $\textsf{ideal} \leftarrow \textsc{Substitute}(C_{\alpha\beta} \leftarrow \phi_{\alpha\beta},\textsf{equations})$;
\ENDIF
\ENDFOR
\RETURN \textsf{ideal};
\end{algorithmic}
\end{algorithm}

We can start from the same algorithm also for computing the Gr\"obner stratum of a gen-segment ideal, indeed by Theorem \ref{th:saturato}\emph{(\ref{it:saturato_i})} we know that also for Gr\"obner strata the relevant variables $C$ are those in the tails of superminimal generators and in the case of a segment ideal $J$ the tail of a monomial $x^\alpha \in J$ contains all the monomials of $\cN(J)_{\vert\alpha\vert}$. After having computed the equations of the marked family, we can exploit the term ordering for further eliminated other variables $C_{\alpha\beta}$ appearing in degree 1 in any of the equation of the ideal of the scheme defining the marked family (Algorithm \ref{alg:grobnerStrata}).

Always starting from Algorithm \ref{alg:markedFamily2}, we can determine the procedure computing the tangent space dimension at the origin of a Gr\"obner stratum, i.e. the its embedding dimension (Algorithm \ref{alg:ed}).

\vspace*{\stretch{1}}

\begin{algorithm}[H]
\capstart
\begin{algorithmic}[1]
\STATE $\textsc{EmbeddingDimension}(J,\sigma)$
\REQUIRE $J$, a truncation of a saturated Borel-fixed ideal, i.e. $J = \underline{J}_{\geqslant m}$ for some $m$.
\REQUIRE $\sigma$, a term ordering such that $J$ is a gen-segment ideal w.r.t. $\sigma$.
\ENSURE the embedding dimension of $\St{\sigma}(J)$.
\STATE $sG_J \leftarrow \textsc{SuperminimalGenerators}(J)$;
\STATE $s\mathcal{G}\leftarrow \emptyset$;
\STATE $N \leftarrow 0$;
\FORALL{$x^\alpha \in sG_J$}
\STATE $N \leftarrow N + \left\vert\cN(J)_{\vert\alpha\vert}\right\vert$;
\STATE $F_{\alpha} \leftarrow x^\alpha$;
\FORALL{$x^\beta \in \cN(J)_{\vert\alpha\vert}$}
\STATE $F_{\alpha} \leftarrow F_{\alpha} + C_{\alpha\beta} x^\beta$;
\ENDFOR
\STATE $s\mathcal{G} \leftarrow s\mathcal{G} \cup \{F_{\alpha}\}$;
\ENDFOR
\STATE $G \leftarrow \textsc{MinimalGenerators}(J) \setminus sG_J$;
\STATE $\textsf{linearEquations} \leftarrow \emptyset$;
\end{algorithmic}
\end{algorithm}

\begin{algorithm}[H]
\caption[Algorithm\hfill computing\hfill the\hfill embedding\hfill dimension\hfill of\hfill a\hfill Gr\"obner\newline stratum.]{The algorithm computing the embedding dimension of a Gr\"obner stratum.}
\label{alg:ed}
\begin{algorithmic}[1]\setcounter{ALC@line}{14}
\WHILE{$G \neq \emptyset$}
\STATE $x^\alpha \leftarrow \min_{\RevLex} G$;
\STATE $x_i \leftarrow \min_{j>0} \{ x_j \mid x^\alpha\}$;
\STATE $x^\gamma \leftarrow \frac{x_0}{x_i}x^\alpha$; 
\STATE $H \leftarrow \textsc{Reduce}\big(x_i \Nf(x^\gamma),G_J\big)$;\hfill \COMMENT{We delete all the monomials in $J$}
\STATE $(Q,R) \leftarrow {}\text{polynomials such that } H = Q \cdot x_0 + R$;
\FORALL{$x^\delta \in \Supp R$}
\STATE $\textsf{equations} \leftarrow \textsf{equations} \cup \{\textsc{Coeff}(R,x^\delta)\}$;
\ENDFOR
\STATE $G \leftarrow G \setminus \{x^{\alpha}\}$;
\ENDWHILE
\STATE $\textsf{syzygies} \leftarrow \textsc{SuperminimalSyzygies}(J)$;
\FORALL{$x^{\gamma}\mathbf{e}_{\alpha} - x^{\gamma'}\mathbf{e}_{\alpha'} \in \textsf{syzygies}$}
\STATE $S(F_{\alpha},F_{\alpha'}) \leftarrow x^{\gamma}F_{\alpha} - x^{\gamma'}F_{\alpha'}$;
\STATE $H \leftarrow \textsc{Reduce}\big(S(F_{\alpha},F_{\alpha'}),G_J\big)$;
\FORALL{$x^\delta \in \Supp H$}
\STATE $\textsf{equations} \leftarrow \textsf{equations} \cup \{\textsc{Coeff}(H,x^\delta)\}$;
\ENDFOR
\ENDFOR
\RETURN $N-\dim_\K \left\langle \textsf{equation} \right\rangle$;
\end{algorithmic}
\end{algorithm}

\begin{example}\label{ex:Hilb4tP3}
Let us consider the saturated Borel-fixed ideals in $\K[x_0,x_1,x_2,x_3]$ with Hilbert polynomial $p(t) = 4t$:
\[
\begin{split}
& J_1 = (x_3,x_2^5,x_2^4x_1^2),\\
& J_2 = (x_3^2,x_3x_2,x_3x_1,x_2^5,x_2^4x_1),\\
& J_3 = (x_3^2,x_3x_2,x_3x_1^2,x_2^4),\\
& J_4 = (x_3^2,x_3x_2,x_2^3).\\
\end{split}
\]
They are all hilb-segment ideals and by Proposition $\ref{prop:segmentsImplications}$ any truncation will be a gen-segment ideal, so that any marked family coincides with the Gr\"obner stratum w.r.t. the term ordering making any ideal a gen-segment ideal. We now discuss how the computational complexity of the costruction of such families of ideals decreases applying the results introduced up to this point.
\begin{enumerate}
\item[$J_1.$] This is the lexicographic ideal associated to $p(t) = 4t$. As seen in Example \ref{rk:stratoLex}, there are only two possible marked family structures:
\[
\Mf(J_1) \simeq \Mf\left((J_1)_{\geqslant m}\right),\ m = 2,3,4, \qquad \Mf\left((J_1)_{\geqslant 5}\right) \simeq \Mf\left((J_1)_{\geqslant m}\right),\ m > 5.
\]
Applying Algorithm \ref{alg:dimTang} we have that
\[
\dim_{\K} T_0\left(\Mf(J_1)\right) = 47 \qquad \text{and} \qquad \dim_{\K} T_0\left(\Mf( (J_1)_{\geqslant 5})\right) = 64.
\]
Furthermore, if we use the homogeneous positive grading induced by the term\hfill ordering\hfill $\DegLex$,\hfill we\hfill obtain\hfill that\hfill both\hfill $\Mf(J_1) \simeq \St{\DegLex}(J_1)$\hfill and\\ $\Mf( (J_1)_{\geqslant 5}) \simeq \St{\DegLex}( (J_1)_{\geqslant 5})$ are affine spaces and
\[
\ed \St{\DegLex}(J_1) = 21 \qquad \text{and} \qquad \ed \St{\DegLex}( (J_1)_{\geqslant 5}) = 23.
\] 
\item[$J_2.$] Again there are two possible marked scheme structure up to isomorphism
\[
\begin{split}
&\Mf(J_2) \simeq \Mf\left((J_2)_{\geqslant 3}\right),\qquad \dim_{\K} T_0 \left( \Mf(J_2)\right) = 61,\\ 
&\Mf\left((J_2)_{\geqslant4}\right) \simeq \Mf\left((J_2)_{\geqslant m}\right),\ \forall\ m > 4,\qquad \dim_{\K} T_0\left(\Mf( (J_2)_{\geqslant 4})\right) = 88.
\end{split}
\]
By computing the Gr\"obner strata w.r.t. $\omega_2 = (9,3,2,1)$, we find that the variables we need to describe these families of ideals are
\[
\ed \St{\omega_2}(J_2) = 24 \qquad \text{and} \qquad \ed \St{\omega_2}( (J_2)_{\geqslant 4}) = 27.
\]
\item[$J_3.$] In this case $\rho = 3$, thus the marked families of the truncations of $J_3$ are all isomorphic to $\Mf(J_3) \simeq \St{\omega_3}(J_3)$, where $\omega_3 = (7,3,2,1)$:
\[
\Mf(J_3) \simeq \Mf\left((J_3)_{\geqslant m}\right),\ \forall\ m,\quad \dim_{\K} T_0 \left( \Mf(J_3)\right) = 44,\quad \ed \St{\omega_3}(J_3) = 24.
\]
\item[$J_4.$] In this case the ideal is ACM, so $x_1$ do not appear in any generator, so we can consider again the saturated ideal itself.
\[
\Mf(J_4) \simeq \Mf\left((J_4)_{\geqslant m}\right),\ \forall\ m,\quad \dim_{\K} T_0 \left( \Mf(J_4)\right) = 28,\quad \ed \St{\omega_4}(J_4) = 16,
\]
where $\omega_4 = (6,4,2,1)$.
\end{enumerate}
\end{example}

\section{Open subsets of the Hilbert scheme II}\label{sec:openSubsetsII}\index{Hilbert scheme!open subset of the}

In this final section of the chapter we will use marked families to cover the Hilbert scheme. The results we will expose belong to the submitted paper \lq\lq Borel open covering of Hilbert schemes\rq\rq\ \cite{BLR} written in collaboration with C. Bertone and M. Roggero.

\bigskip

We consider again the embedding of the Hilbert scheme\index{Hilbert scheme} in a suitable projective space through the Pl\"ucker embedding $\Hilb{n}{p(t)} \subset \Grass{q(r)}{N(r)}{\K} \hookrightarrow \PP^{\binom{N(r)}{p(r)}-1}$\index{Grassmannian} and we will use the same notation introduced in Section \ref{sec:openSubsetsI}. As we seen, each open affine subset $\mathcal{U}_J$ of the Grassmannian and the corresponding open subset $\mathcal{H}_J = \mathcal{U}_J \cap \Hilb{n}{p(t)}$ of the Hilbert scheme can be uniquely identified with a monomial ideal of $\K[x]$ generated by $q(r)$ monomials of degree $r$ that fixes the Pl\"ucker coordinate $\Delta_J \neq 0$. We will denote this set of ideals with $\mathcal{M}^n$, with $\mathcal{B}^n$ its subset composed by Borel-fixed ideals and with $\mathcal{B}^n_{p(t)}$ the subset of Borel-fixed ideals with Hilbert polynomial $p(t)$, i.e.
\begin{equation}
\glshyperlink{BorelNpoly} \subset \gls{BorelN} \subset \gls{MonIdealN}.
\end{equation}

In the following results, we state some close relation between open subsets of Grassmannians\hfill and\hfill properties\hfill of\hfill the\hfill ideals\hfill that\hfill correspond\hfill to\hfill the\hfill points\hfill of\\ $\Grass{q(r)}{N(r)}{\K}$. 
 
\begin{lemma}\label{lem:pluckerlocali} 
Let $J$ and $I$ be ideals in $\Grass{q(r)}{N(r)}{\K}$,  with $J\in \mathcal{M}^n$ and $G_J$ its monomial basis. Then the following statements are equivalent:
\begin{enumerate}[(i)]
\item\label{it:pluckerlocali_i} $\Delta_J(I) \neq 0$;
\item\label{it:pluckerlocali_ii} $I_r$ can be represented by a matrix $\mathfrak{M}(I_r)$ of the form $\left(\begin{array}{l|r} \mathrm{Id}  & R\end{array}\right)$, where the left block is the $q(r)\times q(r)$ identity matrix and corresponds to the monomials in $G_J$ and the entries of the right block $R$ are constants $-c_{\alpha\beta}$, where $x^\alpha \in G_J$ and $x^\beta \in \cN(J)_r$;
\item\label{it:pluckerlocali_iii} $I$ is generated by a $J$-marked set:
\begin{equation}
G=\left \{f_\alpha= x^\alpha-\sum c_{\alpha\beta} x^\beta\ \vert\  \Ht(f_\alpha)=x^\alpha\in G_J \right\}.
\end{equation}
\end{enumerate}
If the previous conditions hold and moreover $J'$ is another monomial ideal in $\mathcal{M}^n$, then: 
\begin{enumerate}[(i)] \setcounter{enumi}{3}
\item\label{it:pluckerlocali_1} $\Delta_{J'}(I)/\Delta_J(I)$  can  be  expressed  as  a  polynomial  in  the   $c_{\alpha \beta}$'s of  degree  $\vert G_J \setminus (G_{J} \cap G_{J'})\vert$;
\item\label{it:pluckerlocali_2} especially, if $G_{J'}=G_J\setminus \{x^{\alpha}\} \cup \{x^{\beta}\}$, then (up to the sign) $\Delta_{J'}(I)/\Delta_J(I)= c_{\alpha \beta}$;
\item\label{it:pluckerlocali_3} we can fix an isomorphism $\AA^{p(r)q(r)}\simeq\mathcal{U}_J $ such that the constants $c_{\alpha \beta}$ are the coordinates of $I$ in $\AA^{p(r)q(r)}$. 
\end{enumerate}
\end{lemma}
\begin{proof}
 $\emph{(\ref{it:pluckerlocali_i})} \Rightarrow \emph{(\ref{it:pluckerlocali_ii})}$ It suffices to multiply any matrix $\mathfrak{M}(I_r)$ by the inverse of its submatrix corresponding to the columns fixed by $J$, since its determinant is $\Delta_J(I)\neq 0$.

$\emph{(\ref{it:pluckerlocali_ii})} \Rightarrow \emph{(\ref{it:pluckerlocali_i})}$ is obvious. 

$\emph{(\ref{it:pluckerlocali_ii})} \Rightarrow \emph{(\ref{it:pluckerlocali_iii})}$ The generators of $I$ given by the rows of $\mathfrak{M}(I_r)$ are indeed a $J$-marked set and, vice versa, the matrix containing the coefficients of the polynomials $f_{\alpha}$ has precisely the shape required in \emph{(\ref{it:pluckerlocali_ii})}.
 
 Finally \emph{(\ref{it:pluckerlocali_1})}, \emph{(\ref{it:pluckerlocali_2})} and \emph{(\ref{it:pluckerlocali_3})} are easy consequences of \emph{(\ref{it:pluckerlocali_ii})}.
\end{proof}

As before we will denote by $\mathcal{G}$ the $J$-marked set:
\begin{equation}\label{JbaseC} 
\mathcal{G}= \left\{F_\alpha= x^\alpha-\sum C_{\alpha\beta} x^\beta\ \big\vert\ \Ht(F_\alpha)=x^\alpha\in G_J, x^\beta \in \cN(J)\right\}\end{equation}
and by $(\mathcal{G})$ the ideal generated in the ring  $\K[C][x]$,  where $C$ is as usual the compact notation for the set of new variables $C_{\alpha \beta}$, $x^\alpha \in G_J$, $x^\beta \in \cN(J)_r$.

\begin{corollary}\label{cor:newVariables}  
In the hypothesis of Lemma \ref{lem:pluckerlocali}, $\mathcal{U}_J$ is isomorphic to the affine space $\AA^{p(r)q(r)} = \Spec \K[C]$.
 The (closed) points in $\mathcal{U}_J$ correspond to all ideals that we obtain from $(\mathcal{G})$ specializing the variables $C_{\alpha \beta}$ to $c_{\alpha \beta}\in \K$.
\end{corollary}

\begin{remark}\label{rk:MarkedSetNotBasis} 
If $J \in \mathcal{M}^n$, then the open subset $\mathcal{U}_J$ is a parameter space for the set of $J$-marked sets. However it is not in general isomorphic to the $J$-marked scheme $\Mf(J)$ because, for instance, the Hilbert polynomial is not necessarily constant on $\mathcal{U}_J$ (see \cite[Example 1.10]{CioffiRoggero}).
\end{remark}

\begin{remark}\label{gradodelta}  
Let $J,J'$ be any couple of monomial ideals in $\mathcal{M}^n$ and let $\delta=\vert G_{J}\setminus G_{J'}\vert $. The localization of $\Delta_{J'}$ in $\K[C]$, the coordinate ring of  $\AA^{p(r)q(r)}\simeq \mathcal{U}_J$, gives a polynomial of degree $\delta$ as shown in Lemma \ref{lem:pluckerlocali}\emph{(\ref{it:pluckerlocali_1})}. In the \lq\lq worst\rq\rq\ case, if we consider $J'$ such that $G_{J}\cap G_{J'}=\emptyset$, then $\delta=q(r)$.
\end{remark} 
  
As $J$ varies in $\mathcal{M}^n$, the open sets $\mathcal{U}_J$ cover the Grassmannian, and so the subsets $\mathcal{H}_J=\mathcal{U}_J\cap \Hilb{n}{p(t)},\ J \in \mathcal{M}^n$ give an open covering of $\Hilb{n}{p(t)}$ by affine subschemes.  It is quite obvious that every open subset $\mathcal{U}_J$ is non-empty, because it contains   the point corresponding to $J$ itself.  If $J$ has Hilbert polynomial $p(t)$ also $\mathcal{H}_J$ is non-empty because it contains $J$ itself, but when $\Proj \K[x]/J \notin \hilb{n}{p(t)}(\K)$, it is not easy to understand general properties of $\mathcal{H}_J$ or even to decide if it is empty or not. For this reason we prefer to consider a slightly different open covering for $\Grass{q(r)}{N(r)}{\K}$ and $\Hilb{n}{p(t)}$, obtained considering only the set of Borel-fixed ideals $\mathcal{B}^n\subset\mathcal{M}^n$, that we will prove to be more convenient for our purposes.
  
\begin{definition}
Given the Grassmannian $\Grass{q(r)}{N(r)}{\K}$ and the set $\mathcal{B}^n$ of all the Borel-fixed ideals $J \subset \K[x]$ such that $\dim_\K J_r = q(r)$, we call \emph{Borel region}\index{Borel region|see{Grassmannian!Borel region of the}}\index{Grassmannian!Borel region of the} of $\Grass{q(r)}{N(r)}{\K}$ the union
\begin{equation}\label{eq:BorelMarkedRegion}
\mathcal{U}=\bigcup_{J\in  \mathcal{B}^n}\mathcal{U}_{J}.
\end{equation}
\end{definition}

Another key role will be played by the linear group \glshyperlink{GL} and its induced action on the Grassmannian $\Grass{q(r)}{N(r)}{\K}$. The action of an element of $\GL(n+1)$ on $\PP^n$ corresponds to a different choice of the basis for $\K[x]_1$ and therefore to a different choice of the basis for $\K[x]_r$. So $g \in \GL(n+1)$ induces a linear change of Pl\"ucker coordinates in the projective space $\PP^{\binom{N(r)}{p(r)}-1}$ in which $\Grass{q(r)}{N(r)}{\K}$ is embedded. Note that not all the linear changes of Pl\"ucker coordinates can be obtained by the action of some element of $\GL(n+1)$ on $\PP^n$.

\begin{lemma}\label{lem:ricgrass} The action of $\GL(n+1)$ on the Borel region $\mathcal{U}$ gives an open covering of $\Grass{q(r)}{N(r)}{\K}$, that is:
\begin{equation}\label{eq:ActionBorelMarkedRegion}
\bigcup_{g\in \GL(n+1)} g\centerdot\mathcal{U} =\bigcup_{\begin{subarray}{c} J\in \mathcal{B}^n\\ g\in \GL(n+1)\end{subarray}} g\centerdot\mathcal{U}_{J} =\Grass{q(r)}{N(r)}{\K}.
\end{equation}
\end{lemma}
\begin{proof} 
Let $I \in \Grass{q(r)}{N(r)}{\K}$ be any ideal and let $\sigma$ be any term order on the monomials of $\K[x]$. Due to Galligo's Theorem \cite{Galligo}, in generic coordinates the initial ideal $J'$ of $I$ is Borel-fixed, and then $J =(J'_r)$ is Borel-fixed too. Moreover, by construction  $J$ is generated by $q(r)$ monomials of degree $r$ and so $J\in \mathcal{B}^n$. Hence for a general $g\in \GL(n+1)$ we have $\Delta_J(g\centerdot I)\neq 0$ that is $g\centerdot I \in \mathcal{U}_{J}$, so that $I\in g^{-1}\centerdot\mathcal{U}_{J}$.   
\end{proof}

\begin{remark}
In the proof of Lemma \ref{lem:ricgrass} we deal with the generic initial ideal $J'$,  which may have minimal generators of degree $>r$. To avoid this problem, we consider $J'_r$, which is Borel-fixed  and is generated by $q(r)$ monomials in degree $r$.
\end{remark}

The new covering of the Grassmannian $\Grass{q(r)}{N(r)}{\K}$ shown in Lemma \ref{lem:ricgrass}  will turn out to be more suitable to study local properties of Hilbert schemes. 
\begin{definition}\label{def:BorelOpenCoveringGrass}
Given the Grassmannian $\Grass{q(r)}{N(r)}{\K}$, we define the \emph{Borel covering}\index{Grassmannian!Borel covering of the}\index{Borel covering of the Grassmannian|see{Grassmannian, Borel covering of the}} of $\Grass{q(r)}{N(r)}{\K}$ as the family of all open subsets of the type $g\centerdot \mathcal{U}_J$ where $J\in \mathcal{B}^n$ and $g\in \GL(n+1)$.
\end{definition}

What happens to the restriction to the Hilbert scheme of the Borel covering of the Grassmannian? First we investigate whether  $\mathcal{H}_{J}=\mathcal U_J \cap \Hilb{n}{p(t)}$ is empty or not, for $J\in \mathcal{B}^n$. Of course if $J$ belongs to $\Hilb{n}{p(t)}$, then $ \mathcal{H}_{J}$ cannot be empty because it contains at least $J$. Moreover, since it is an open subset of $\Hilb{n}{p(t)}$, if it contains a point, it also contains an open subset of at least one irreducible component of $\Hilb{n}{p(t)}$. 

\begin{proposition}\label{fuori} 
 If $J\in \mathcal{B}^n$, then:
\begin{equation}\label{hilbvuoto}
 \mathcal{H}_J \neq  \emptyset\quad\Longleftrightarrow\quad \Proj \K[x]/J \in \hilb{n}{p(t)}(\K).
\end{equation}
As a consequence, if we define the \emph{Borel region} $\mathcal{H}$ of $\Hilb{n}{p(t)}$ as $\mathcal {H}=\mathcal{U}\cap \Hilb{n}{p(t)}= \cup_{J\in B^n_{p(t)}} \mathcal{H}_J$, we get:
\begin{equation*}
\Hilb{n}{p(t)}=\bigcup_{\begin{subarray}{c}g\in \GL(n+1) \\ J \in \mathcal{B}^n_{p(t)}\end{subarray}} g \centerdot\mathcal{H}_J=\bigcup_{g\in \GL(n+1) } g \centerdot\mathcal{H}.
\end{equation*} 
\end{proposition} 
\begin{proof} 
We prove only the non-trivial part $(\Rightarrow)$ of the first statement. Assume that $\Proj \K[x]/J \notin \hilb{n}{p(t)}(\K)$. By Gotzmann's Persistence Theorem,\index{Gotzmann's Persistence Theorem} this is equivalent to $\dim_\K J_{r+1} > q(r+1)$. If $I$ is any ideal in $\mathcal{U}_J$, then it has a set of generators as those given in Lemma \ref{lem:pluckerlocali}\emph{(\ref{it:pluckerlocali_iii})}, so that  $\dim_\K I_{r+1} \geqslant \dim_\K J_{r+1} > q(r+1)$ (see \cite[Corollary 2.3]{CioffiRoggero}). Hence $\Proj \K[x]/I \notin \hilb{n}{p(t)}(\K)$.
The other  statement is a direct consequence of the first one and of Lemma \ref{lem:ricgrass}.
\end{proof}
 
 We point out that in Proposition \ref{fuori}, the hypothesis $J\in \mathcal{B}^n$ is necessary, as shown  in the following example.
 \begin{example}[{Cf. Example \ref{ex:wrongStratumHilb2P2}}]\label{open0dim}
  Let us consider again the Hilbert scheme \gls{Hilb2P2}: in this case $r=2$ and $q(2)=4$. The monomial ideal $J=(x_{2}^2,x_2x_{1},x_1^2,x_{0}^2)$ is generated by $4$ monomials of degree 2, it is not Borel-fixed and  obviously does not belong to $\Hilb{2}{2}$ because it is a primary ideal over the irrelevant maximal ideal $(x_{2},x_1,x_{0})$. Nevertheless, $\mathcal{H}_J$ contains the ideal $(x_{2}^2-x_{2}x_{0},x_2x_1,x_1^2 - x_{1}x_0,x_{0}^2-x_{1}x_0-x_{2}x_{0})$ corresponding to the set of points $\{[1:0:1], [1:1:0]\}$ and, more generally, all the ideals corresponding to  pair of  distinct points outside the line $x_{0}=0$ and not on the same line through $[1:0:0]$.
 \end{example}
  
 \begin{corollary}
Set-theoretically we have that:
\begin{equation}\label{eq:hyperplanes}
\Hilb{n}{p(t)} \subseteq \bigcap_{\begin{subarray}{c}g\in \GL(n+1) \\ J \in \mathcal{B}^n\setminus \mathcal B^n_{p(t)}\end{subarray}}g \centerdot \Pi_J
\end{equation}
where $\Pi_J$ is the hyperplane in $\PP^{\binom{N(r)}{p(r)}-1}$ given by $\Delta_J=0$.
\end{corollary}

\begin{example}\label{nonBorel}
  Let us consider the Hilbert polynomial $p(t)=3t$ in $\PP^3$. The closed points of \gls{Hilb3tP3} corresponds to curves in $\PP^3$ of degree 3 and arithmetic genus 1, hence it contains all the smooth plane elliptic curves and also some singular or reducible or non-reduced curve.  The  Gotzmann number of $p(t)=3t$ is $r=3$ and so $q(3)=11$ and $\binom{20}{11}-1=167959$.
  
The only Borel-fixed ideal defining points on $\Hilb{3}{3t}$ is the lexicographic ideal: 
\[
L= \underline{L}_{\geqslant 3} = (x_3,x_2^3)_{\geqslant 3}.
\]
The Borel region of $\Hilb{3}{3t}$ is then  equal to the open subset $\mathcal{U}_{L}\cap \Hilb{3}{3t}$. The Grassmannian $\Grass{11}{20}{\K}$ in which  $\Hilb{3}{3t}$  is embedded has dimension $q(3)\cdot p(3)=99$.  
Using Algorithm \ref{alg:BorelGeneratorDFS}, we compute the complete list of Borel-fixed ideals in $\Grass{11}{20}{\K}$ that do not belong to $\Hilb{3}{3t}$:
\[
\begin{array}{r cc l}
2t+3 &&& J_1 = (x_3^2,x_3x_2,x_3x_1,x_2^3,x_2^2x_1)_{\geqslant 3},\\
&&& J_2 = (x_3^2,x_3x_2,x_2^2,x_3x_1^2)_{\geqslant 3},\\
t+6 &&& J_3 = (x_3^2,x_3x_2,x_2^3,x_2^2x_1,x_3x_1^2,x_2x_1^2)_{\geqslant 3},\\
9 &&& J_4 = (x_3^2,x_3x_2^2,x_2^3,x_3x_2x_1,x_2^2x_1,x_3x_1^2,x_2x_1^2,x_1^3)_{\geqslant 3}
\end{array}
\]
Then $\Hilb{3}{3t}$ as a subscheme of  $\Grass{11}{20}{\K} \hookrightarrow \PP^{167959}$  is contained set-theoretically in the intersection of the hyperplanes $\Pi_{J_i}$ given by $\Delta_{J_i}=0$, $i=1,2,3,4$ (and in all the hyperplanes obtained from these by the action of $\GL(4)$).
\end{example}
  
\begin{definition} 
 The \emph{Borel covering} of $\Hilb{n}{p(t)}$\index{Hilbert scheme!Borel covering of the}\index{Borel covering of the Hilbert scheme|see{Hilbert scheme, Borel covering of the}} will be the family of all the open subsets of $g\centerdot\mathcal{H}_J$ where $J\in \mathcal{B}^n_{p(t)}$ and $g\in  \GL(n+1)$.
 
 For any Borel ideal $J$ in $\mathcal{B}^n_{p(t)}$, the open subset $\mathcal{H}_J=\mathcal{U}_J\cap \Hilb{n}{p(t)}$ will be called the \emph{$J$-marked region}\index{marked region} of $\Hilb{n}{p(t)}$.
 \end{definition} 
 
 The name \lq\lq $J$-marked region\rq\rq\ comes from Lemma \ref{lem:pluckerlocali} and its connection with the $J$-marked scheme will be clearer with Theorem \ref{isomHM}.

\begin{remark}\label{unosolo}
\begin{enumerate}
\item If the Hilbert polynomial $p(t)$ is the constant $r$,  then every Borel ideal $J\in \Grass{q(r)}{N(r)}{\K}$ belongs to $\Hilb{n}{r}$ i.e. $\mathcal{B}^n = \mathcal{B}^n_{p(t)}$.  Then in the zero-dimensional case the family of hyperplanes $\Pi_J$ considered in Proposition \ref{fuori} is indeed empty.
\item If $\deg p(t) = d\geqslant 1$, $\mathcal{B}^n\setminus \mathcal{B}^n_{p(t)}$  in general is not empty and its elements define subschemes of $\PP^n$ of dimension equal to or lower than the one of the subschemes parametrized by $\Hilb{n}{p(t)}$.
Indeed, if $I\in \mathcal{B}^n$ has Hilbert polynomial $\widetilde{p}(t)\neq p(t)$, then being $\widetilde{q}(r) = q(r)$
for Gotzmann's Persistence Theorem, $\dim_\K I_t > q(t)$ for $t\geqslant r$. Hence for $t\gg0$, $\widetilde{q}(t)>q(t)$ and $\widetilde{p}(t)< p(t)$. So $\deg \widetilde p(t)\leqslant d$.
\item \label{unosolo_iii}If $\mathcal{B}^n_{p(t)}$ contains  only one ideal, then $\Hilb{n}{p(t)}$ is a  smooth rational projective variety. Indeed, we know that $\mathcal{B}^n_{p(t)}$ contains at least the lexicographic ideal $L$, i.e. the ideal generated in degree $r$ by the $q(r)$ maximal monomials w.r.t. the term order $\DegLex$.  In Section \ref{sec:lexsegment} we proved that $\mathcal{H}=\mathcal{H}_{L}$ is isomorphic to an affine space. By Proposition \ref{fuori}, as $g$ varies in $\GL(n+1)$, the open subsets  $g \centerdot \mathcal{H}$ cover $\Hilb{n}{p(t)}$. Thus $\Hilb{n}{p(t)}$ is smooth and rational as claimed. 
\end{enumerate}
\end{remark}

The open subset $\mathcal{H}_J = \mathcal{U}_J\cap \Hilb{n}{p(t)}$ of $\Hilb{n}{p(t)}$ is then a closed subscheme in the affine space in $\AA^{p(r)q(r)}\simeq \mathcal{U}_J$. Moving from Lemma \ref{lem:pluckerlocali} and Corollary \ref{cor:newVariables} we can determine the scheme structure of $\mathcal{H}_J$ in $\AA^{p(r)q(r)}$, starting from the set of $J$-marked polynomials $\mathcal{G}$ as in Definition \ref{JbaseC}.

\begin{definition} 
We will denote by $\mathfrak{A}_J$ the ideal in $\K[C]$ defining $\mathcal{H}_J$ as an affine subscheme of $\AA^{p(r)q(r)}$ through the isomorphism of Lemma \ref{lem:pluckerlocali}\emph{(\ref{it:pluckerlocali_3})}.
\end{definition}

\begin{remark}\label{eleminA} 
We  obtain every ideal $I\in \mathcal{U}_J$   specializing (in a unique way) the variables $C_{\alpha \beta}$ in $(\mathcal{G})$ to $c_{\alpha\beta}\in \K$, but not every specialization gives rise to an ideal $I$ in $\mathcal{H}_J$, that is to an ideal with Hilbert polynomial $p(t)$. This last condition holds for an ideal $I$ if and only if every  polynomial has an unique $J$-normal form modulo $I$, that is if and only if every $J$-reduced polynomial  in $I$ vanishes. Hence, the ideal $\mathfrak{A}_J$  is made by  the coefficients w.r.t. the variables $x$ of all the polynomials  $(\mathcal{G})\subset \K[C][x]$ that are $J$-reduced.
\end{remark}

Due to Macaulay's Estimate on the Growth of Ideals\index{Macaulay's Estimate on the Growth of Ideals} we know that if $I$ is generated by a $J$-marked set, then $\dim_\K I_t \geqslant q(t)$ for every $t\geqslant r$. Moreover, by Gotzmann's Persistence Theorem,\index{Gotzmann's Persistence Theorem}  $\Proj \K[x]/I\in \hilb{n}{p(t)}(\K)$ if the equality holds for $r+1$, that is if $\dim_\K I_{r+1}=q(r+1)$.

Then let $(\mathcal{G})\subseteq \K[C][x]$, $\mathfrak{M}\big((\mathcal{G})_r\big)$ and $\mathfrak{M}\big((\mathcal{G})_{r+1}\big)$ be respectively the matrices whose columns correspond to the monomials in $\K[x]_r$ and $\K[x]_{r+1}$ and whose rows contain the coefficients of monomials in the polynomials $F_\alpha$ and $x_iF_\alpha$ respectively. Thus, a set of generators for the ideal $\mathfrak{A}_J$  is given by the minors of order $q(r+1)+1$ of the matrix $\mathfrak{M}\big((\mathcal{G})_{r+1}\big)$.  

An easy computation shows that this is in general a very large set of polynomials! Indeed $\mathfrak{M}\big((\mathcal{G})_{r+1}\big)$ is a $(n+1) q(r) \times N(r)$ matrix and the number of its minors of order $q(r+1)+1$ is $\binom{(n+1)q(r)}{q(r+1)+1} \cdot \binom{N(r)}{q(r+1)+1}$ and their degree is up to $q(r+1)+1$. Looking at the special form of $\mathfrak{M}\big((\mathcal{G})_{r+1}\big)$, we will show in Theorem \ref{gradod+2} that the number of minors of $\mathfrak{M}\big((\mathcal{G})_{r+1}\big)$  that are sufficient  to impose the condition on the the rank  can be drastically reduced  and the degree of the involved determinants is bounded by $d+2$. 

\begin{example}\label{numerominori}
Let us consider for instance the Hilbert scheme \gls{Hilb2P2}. For every monomial ideal $J\in \mathcal{M}^n$, we have $r=2$, $q(2)=4$, $q(3)=8$, $N(3)=10$. Then $\mathfrak{M}\big((\mathcal{G}_J)_{3}\big)$ is a $12\times 10$ matrix and the number of its minors of order 9 (with degree up to 9) is $\binom{3\cdot 4}{9}\cdot \binom{10}{9}=2200$.
\end{example}

In order to obtain a better set of generators for $\mathfrak{A}_J$, we now prove that the open subset $\mathcal{H}_J$ for a Borel-fixed ideal $J$ defining a point of $\Hilb{n}{p(t)}$ is nothing else but the $J$-marked scheme $\Mf(J)$. 

\begin{theorem}\label{isomHM}
There is a scheme theoretic isomorphism: 
\begin{equation}
\mathcal{H}_J \simeq \Mf(J).
\end{equation}
 \end{theorem}  
 \begin{proof}
 The thesis directly follows from the two constructions of $\Mf(J)$ and $\mathcal{H}_J$. Both constructions start from a $J$-marked set $\mathcal{G}\subseteq \K[C][x]$ (as in Definition \ref{JbaseC}). As shown in Proposition \ref{prop:MftrhoughMatrices}, we can obtain a set of generators for the ideal defining $\Mf(J)$ imposing conditions on the rank of some matrices. In the present  hypothesis,  we can  consider  only one matrix, the one corresponding to the degree $r+1$,  and impose that its rank is  $\leqslant \dim_\K J_{r+1}$. This matrix turns out to be indeed  $\mathfrak{M}\big((\mathcal{G})_{r+1}\big)$ and $\dim_\K J_{r+1}=q(r+1)$. Then in both cases, a set of generators  is  given by    the minors of order $q(r+1)+1$ of the matrix  $\mathfrak{M}\big((\mathcal{G})_{r+1}\big)$.
 \end{proof}  

Thanks to this last result, $\mathfrak{A}_J$ is the ideal in $\K[C]$ defining $\mathcal{H}_J$ or equivalently $\Mf(J)$ as an affine subscheme in $\AA^{p(r)q(r)}$. The isomorphism between a $J$-marked region of $\Hilb{n}{p(t)}$ and the corresponding $J$-marked scheme allows us to embed $\mathcal{H}_J$ in affine linear spaces of \lq\lq low\rq\rq\ dimension by Theorem \ref{schisom}. We can choose linear spaces of different dimension, depending on whether we want to keep control on the degree of the equations defining the scheme structure or not.

\subsection{Equations defining $\mathcal{H}_J$ in local Pl\"ucker coordinates}\label{sec:localEqd+2}

In our reasoning the matrix $\mathfrak{M}\big((\mathcal{G})_{r+1}\big)$ has a major role,  therefore we now look closer at its shape. Remind that the ideal $J = J_{\geqslant r} = (J_r)$ belongs to $\mathcal{B}^n_{p(t)}$ and that we are assuming that the degree of the Hilbert polynomial $p(t)$ is equal to $d$.

\begin{lemma}\label{nice}  Up to permutations on rows and columns, $\mathfrak{M}\big((\mathcal{G})_{r+1}\big)$  assumes the following simple form:
   \begin{equation}
 \mathfrak{M}\big((\mathcal{G})_{r+1}\big)=\left( \begin{array}{cccc|ccc}
  \mathtt{Id}(n,\dots,d+1)   &\bullet &\bullet &\bullet &\bullet   &\dots&\bullet\\ 
  0  &\mathtt{Id}(d)   &\bullet&\bullet&\bullet   &\dots&\bullet\\
 \vdots   &  \vdots & \ddots &\bullet  &\vdots  & \vdots&\vdots\\
   0   &0 &0 & \mathtt{Id}(0)  &\bullet   &\dots&\bullet\\
   \hline 
    \star  &\bullet   &\bullet&\bullet&\bullet   &\dots&\bullet\\
 \vdots   &  \vdots & \vdots &\vdots &\vdots  & \vdots&\vdots\\
   \star  &\bullet   &\bullet&\bullet&\bullet   &\dots&\bullet
   \end{array}\right)   
  \end{equation}
    where 
\begin{itemize}
     \item the columns on the left of the vertical line correspond to monomials in $J_{r+1}$;
     \item the columns on the right of the vertical line correspond to monomials in $\cN(J)_{r+1}$;
      \item $\mathrm{Id}(n,\dots,d+1) $ is an identity matrix of order $\binom{n-d+r}{r+1}$, corresponding to the monomials in $\K[x_{d+1},\dots,x_n]_{r+1}$;
     \item  $\mathrm{Id}(d) $, \dots ,$\mathtt{Id}(0) $ are identity matrices of suitable dimensions $\leqslant q(r)$, corresponding to monomials in $J_{r+1}$ with minimal variable $x_d,\dots, x_0$ respectively;
     \item   \lq\lq $\star$\rq\rq\ stands for  entries that are all 0, except at most one entry equal to 1 in each row;
     \item    \lq\lq $\bullet$\rq\rq\  stands for entries that are either 0 or   coefficients $-C_{\alpha\beta}$.    
  \end{itemize}
\end{lemma}
\begin{proof} 
We consider the $\K[C]$-module of polynomials in $(\mathcal{G})$ of degree $r+1$ with respect to the variables $x$ and its set of generators $\{x_i F_\alpha\ \vert\ F_\alpha \in\mathcal{G}, i=0,\dots,n\}$. We write inside the matrix $\mathfrak{M}\big((\mathcal{G})_{r+1}\big)$ the coefficients of the monomials in $\K[x]_{r+1}$ appearing in these polynomials $x_iF_\alpha$.

 First of all we order the columns writing first the monomials in $J_{r+1}$, listed in decreasing order w.r.t. $\RevLex$, and then the monomials in $\cN(J)_{r+1}$. In this way the first monomial is $x_n^{r+1}$, the only one with minimal variable $x_n$, after this there are   the monomials whose minimal variable is $x_{n-1}$, and so on. 
  
  The rows are ordered in a similar way. 
   Every monomial in $ x^\gamma J_{r+1}$ can be written as a product  $\dec{x^\alpha}{x_i}{J}{}$  such that $x_i = \min x_i x^\alpha \leq \min x^\alpha$ and $x^\alpha \in J_r$ (Lemma \ref{lem:monomialDecomposition}).  The  first rows (those above the horizontal line in the picture) correspond  to polynomials  $x_iF_\alpha$ such that $x_i = \min x_ix^\alpha$ ordered w.r.t. $\RevLex$ on the initial  monomials $x_ix^\alpha$. The first row   corresponds to $x_n F_{x_n^r}$, after there are the rows corresponding to polynomials of the type $x_{n-1}F_{\alpha}$ with $x^\alpha \in \K[x_{n-1},x_n]$ and so on.  Below the horizontal line  we list the rows corresponding to the remaining polynomials $x_iF_\alpha$ such that $\min x_i x^\alpha < x_i$.
   
 The top left submatrix, let us call it $\mathcal{D}$, is an upper triangular matrix of order $q(r+1)$. In fact, as $J\in \mathcal{B}^n_{p(t)}$, then $J_{r+1}$ contains $\K[x_{d+1},\dots,x_{n}]_{r+1}$ (see Proposition \ref{prop:degreeHPpowerVars}) and so each monomial in $J_{r+1}$ corresponds to one and only one column and  row in $\mathcal{D}$.
 
  Moreover in  the row of  $\mathcal{D}$ corresponding to a polynomial $x_iF_\alpha$ with initial monomial $x_ix^\alpha$, the entry on the main diagonal  is the coefficient of  $x_ix^\alpha$ in $x_iF_\alpha$, i.e. $1$. If $x_ix^\beta$ is any  monomial appearing in $x_i \tail{F_\alpha}{}{}$,  then either $x_ix^\beta\notin J$, hence its coefficient $-C_{\alpha\beta}$ is written on the right of the vertical line, or  $x_ix^\beta\in J$, that is $x_ix^\beta=x_jx^{\alpha'}$ for some $x^{\alpha'}\in J$ and  $x_j =\min x_i x^\beta < x_i$,  hence its coefficient $-C_{\alpha\beta}$  is written in one of the columns corresponding to monomials with minimal variable $x_j$ lower than $x_i$. Thus in $\mathcal{D}$ there are identity blocks $\mathrm{Id}(i)$ corresponding to monomials in $J_{r+1}$ with minimal variable $x_i$. 
    
    Furthermore,  the minimal variable in every monomial $x^\beta \in \cN(J)_r$ is lower than or equal to $x_{d}$: hence the first block of $\mathcal{D}$   is a big  identity matrix $\mathtt{Id}(n,\dots,d+1)$  of order $\binom{n-d+r}{r+1}$, corresponding to monomials in $J_{r+1}$ with minimal variable $x_{d+1}, \dots, x_n$.
    The same arguments holds for  the \lq\lq $\star$\rq\rq\ under the horizontal line.
 \end{proof}

We will now determine the dimension of a linear affine space in which $\mathcal{H}_J$ can be embedded and furthermore to study in which cases we can control the degree of the defining equations, bounding it using only $d$.  As $\mathfrak{A}_J$  is the localization in the open subset $\mathcal{U}_J$ of the ideal defining the Hilbert scheme $\Hilb{n}{p(t)}$ in \glshyperlink{KDelta}, we can derive a bound on a set of generators of $\mathfrak{A}_J$  from the known bounds for analogous global results. In Chapter \ref{ch:HilbertScheme} we showed that Iarrobino and Kleiman proved that $\Hilb{n}{p(t)}$ is generated in degree $q(r+1)+1$. Later on, Haiman and Sturmfels proved the Bayer's conjecture saying that $\Hilb{n}{p(t)}$ is generated in the far lower degree $n+1$. Unluckily, global results do not give a satisfying bound in the local case, because the global Pl\"ucker coordinate $\Delta_{J'}$, when localized in $\mathcal{U}_J$, corresponds to a polynomial in $\K[C]$ whose  degree can vary between 1  and $q(r)$ (Remark \ref{gradodelta}).

\begin{example}
We consider again \gls{Hilb3tP3} as in Example \ref{nonBorel}.
\begin{itemize}
	\item Localizing $\Delta_{L}$ at $\mathcal{U}_{J_1}$, we obtain a polynomial of degree 1, since the monomial basis of $L$ is $G_{J_1}\setminus \{x_2 x_1^2\}\cup\{x_3^2x_0\}$ (Lemma \ref{lem:pluckerlocali}\emph{(\ref{it:pluckerlocali_2})});
	\item for $\mathcal{U}_{J_i}$, $i=2,3,4$, we count the monomials in $G_{J_i}\setminus G_{L}$. We then obtain that localizing at $\mathcal U_{J_i}$, the Pl\"ucker coordinate $\Delta_{L}$ becomes a polynomial of degree $i$, $i=2,3,4$ in the $C_{\alpha\gamma}$ (Lemma \ref{lem:pluckerlocali}\emph{(\ref{it:pluckerlocali_1})}).
\end{itemize}
\end{example}

We now prove that the equations defining $\mathcal{H}_J$ in $\AA^{p(r)q(r)}$, that is in the local case, are of degree $\leqslant d+2$.

\begin{theorem}\label{gradod+2} 
Let $p(t)$ be an admissible Hilbert polynomial in $\PP^n$,  of degree $d$ and Gotzmann number $r$.
If $J\in \mathcal{B}^n_{p(t)}$, then the ideal   $\mathfrak{A}_J$ defining $\mathcal{H}_J$  as a subscheme of $\AA^{p(r)q(r)}$, that is in \lq\lq local\rq\rq\ Pl\"ucker coordinates, is generated  in degree smaller than or equal to $d+2$ by $p(r+1)\cdot\left((n+1)q(r)-q(r+1)\right)$ polynomials.
\end{theorem}
 We give two (equivalent) proofs that $\mathfrak{A}_J$ is generated in degree $\leqslant d+2$: the first one uses minors of the matrix $\mathfrak{M}\big((\mathcal{G})_{r+1}\big)$, and in this way we also count the number of generators; the second one uses the Buchberger-like criterion on the reduction of $S$-polynomials.
 
\begin{proof}[First proof]
The ideal $\mathfrak{A}_J$ of $\mathcal{H}_J$ is generated by  the minors of order $q(r+1)+1$ of $\mathfrak{M}\big((\mathcal{G})_{r+1}\big)$, that we think to be written like in Lemma \ref{nice}.  As the determinant of the top left submatrix  of order $q(r+1)$ (called $\mathcal{D}$ in the proof of Lemma \ref{nice}) is 1, we can just consider the minors of order $q(r+1)+1$ containing $\mathcal{D}$:
\begin{equation}
  \det \left( \begin{array}{cccc|c}
  \mathtt{Id}(n,\dots,d+1)   &\bullet &\bullet &\bullet &\bullet   \\
  0  &   \mathtt{Id}(d)   &\bullet&\bullet&\bullet   \\
 \vdots   &  \vdots & \ddots &\bullet  &\vdots  \\
   0   &0 &0 &   \mathtt{Id}(0)  &\bullet   \\
   \hline 
    \star  &\bullet   &\bullet&\bullet&\bullet    
   \end{array}\right).  
\end{equation}
   We perform Gaussian reduction on the last rows. In $\star$ there is at most a non-zero element, which is 1; if necessary, we perform a first row reduction, to make it a 0. At the end of this first step of reduction, the degree of $\bullet$ in the last row remains at most 1 in $\K[C]$.
    
With the second row reduction,  we obtain that the above determinant is equal to the following:
\begin{equation*}
  \det \left( \begin{array}{cccc|c}
  \mathtt{Id}(n,\dots,d+1)   &\bullet &\bullet &\bullet &\bullet   \\
  0  &   \mathtt{Id}(d)   &\bullet&\bullet&\bullet   \\
 \vdots   &  \vdots & \ddots &\bullet  &\vdots  \\
   0   &0 &0 &   \mathtt{Id}(0)  &\bullet   \\
   \hline 
    0 & 0   &\circ_2&\circ_2&\circ_2    
   \end{array}\right).  
 \end{equation*} 
where $\circ_2 $ stands for polynomials in $\K[C]$  of degree at most $2$. 

Going on with Gaussian reduction, the determinant is equal to the  element appearing in the last line and last column, which is a   polynomial in $\K[C]$ of degree $\leqslant d+2$. 

For the number of polynomials that generate $\mathfrak{A}_J$, we simply count the number of minors of $\mathfrak{M}\big((\mathcal{G})_{r+1}\big)$ of order $q(r+1)+1$ containing the matrix $\mathcal{D}$.
\end{proof}
  
 \begin{proof}[Second proof] 
  As shown in Theorem \ref{nostrobuch} and Corollary \ref{cor:nostrobuch}, we can obtain a set of generators for  $\mathfrak{A}_J$ also using a special procedure of reduction of $S$-polynomials of elements in $\mathcal{G}$ with a  Buchberger-like criterion analogous to the one for Gr\"obner bases. The only difference when a term order is not defined is that reductions must be chosen in a careful way in order to have a Noetherian reduction. In the present  hypothesis,  we consider only EK-polynomials  of degree $r+1$ with respect to the variables $x$, that is of the type $x_iF_\alpha -x_jF_{\alpha'}$ with $x_i > \min x^\alpha$ and $x_i x^\alpha= \dec{x^{\alpha'}}{x_j}{J}{}$, that  correspond to a basis of the syzygies of $J$ in degree $r+1$.
 
If $x_i x^{\beta}$  is a monomial of $J$ that appear in $x_iF_\alpha -x_jF_{\alpha'}$, then  $x^\beta \in \cN(J)_r$ and  $x_i x^{\beta}=\dec{x^\gamma}{x_h}{J}{}$, i.e. $x_h = \min x_ix^{\beta} <x_i$ and $x^\beta \in J_r$.  Then we can perform a step of reduction $\xrightarrow{\mathcal{V}^J_{r+1}}$ of $x_h x^\gamma$ rewriting it by $x_h \tail{F_\gamma}{}{}$. If some monomial of $x_h x^\beta-x_h F_\gamma$ belongs to $J$, then again we can reduce it using some polynomial $x_{h'}F_{\gamma'}$ with $x_{h'}<x_h$.
  
At every step of reduction  a monomial is replaced by a sum of other monomials multiplied by one of the variables $C$. Then at every step of reduction the  degree of the coefficients   directly involved  increases  by 1.  If $x^{\eta_0},x^{\eta_1}, \dots, x^{\eta_s}$ is a sequence of monomials in $J_{m}$ such that $x^{\eta_{i+1}}$ appears in the tail of the reduction of $x^{\eta_{i}}$, then $\min x^{\eta_{i+1}} < \min x^{\eta_{i}}$.  As the minimal variable of any monomial in $\cN(J)_r$ is lower than or equal to $x_{d}$,    the   length of any such  chain   is at most $d+1$. Thus, the final degree of the coefficients is at most $1+(d+1)=d+2$.
 \end{proof}

\begin{example}
We consider again \gls{Hilb2P2}, already investigated in Example \ref{numerominori}. If we consider all the minors of $\mathfrak{M}\big((\mathcal{G})_{r+1}\big)$ of order $q(r+1)+1$, we obtain a set of generators for $\mathfrak{A}_J$ of cardinality 2200.
Using Theorem \ref{gradod+2}, we see that actually in order to define $\mathfrak{A}_J$ we just need 8 minors of degree 2.
\end{example}

By  Theorem \ref{isomHM}, we can exploit the techniques presented  for $J$-marked schemes to embed $\mathcal{H}_J$ in linear affine spaces of lower dimension. We consider again the notation introduced for Theorem \ref{schisom}: for an ideal $J \in \mathcal{B}^n_{p(t)}$, let $\underline{J}$ be its saturation and $\rho$ the maximal degree of a monomial divisible by $x_1$ in $G_{\underline{J}}$; if there are no such monomials in $G_{\underline{J}}$, we set $\rho=0$. Moreover if $x^{\alpha}\in G_{\underline{J}}$, we write $x^{\overline{\alpha}}$ for the monomial $x^{\alpha}x_0^{m-\vert\alpha\vert}\in G_{\underline{J}_{\geqslant m}}$, if $\vert\alpha\vert<m$; otherwise  $x^{\overline{\alpha}}=x^\alpha$. Finally, we will denote by $\varphi_{J,r}$ the embedding $\mathcal{H}_J\hookrightarrow \AA^{p(r)q(r)}$ given by Theorem \ref{isomHM} and Theorem \ref{gradod+2}.

\begin{theorem}\label{isomH}
In the established setting, the followings hold:
\begin{enumerate}[(i)]
\item if $m \geqslant r$, then  $\Mf(\underline{J}_{\geqslant m})\simeq \mathcal{H}_J$;
 \item if $m < r$, then  $\Mf(\underline{J}_{\geqslant m})$ is a closed subscheme of $\mathcal{H}_J$, (eventually equal). If we consider the embedding $\phi_{J,r}(\mathcal{H}_J) \subset \AA^{p(r)q(r)}$, then $\Mf(\underline{J}_{\geqslant m})$ is cut out by a suitable linear space;
  \item $\mathcal{H}_J\simeq\Mf(\underline{J}_{\geqslant m})$ if  and only if either $\underline{J}_{\geqslant m}= J$ or $m\geqslant \rho-1$.
  \end{enumerate}
 In particular, if $\rho >0$, then  $\rho -1$ is the smallest integer $m$ such that: 
\[ 
\mathcal{H}_{J}\simeq\Mf(\underline{J}_{\geqslant m}).
\]
   Especially, the isomorphism  $\mathcal{H}_{J}\simeq \Mf(\underline{J}_{\geqslant \reg(\underline{J})})$ induces an embedding $\phi_{J,\reg(\underline{J})}$ of $\mathcal{H}_{J}$ in an  affine space of dimension   $\vert G_{\underline{J}}\vert \cdot p\big(\reg(\underline{J})\big)$ and the isomorphism  $\mathcal{H}_{J}\simeq \Mf(\underline{J}_{\geqslant \rho-1})$ induces an embedding  $\phi_{J,\rho -1}$ of  $\mathcal{H}_J$ in an  affine space of dimension 
   \[
\sum_{x^{\overline{\alpha}}\in G_{\underline{J}_{\geqslant \rho-1}}} \left\vert \cN(\underline{J}_{\geqslant \rho-1})_{\vert\overline{\alpha}\vert}\right\vert.
\]
\end{theorem}
\begin{proof}
Thanks to the isomorphism $\mathcal{H}_J\simeq \Mf(\underline{J}_{\geqslant r})$ of Theorem \ref{isomHM}, the statements are straightforward consequences of Theorem \ref{schisom}.
\end{proof}

The embedding  $\phi_{J,\rho-1}$ (or more generally $\phi_{J,m}$ with $\rho-1 \leqslant m < \reg(\underline{J})$) of  $\mathcal{H}_J$ in affine spaces defined in Theorem \ref{isomH} are computationally  advantageous, because in order to compute equations for $\mathcal{H}_J$ we   deal with a small number of variables, namely smaller than $p(r)q(r)$; however, using these embedding   we do not have any control on the degree of the equations defining $\mathcal{H}_J$. 

If we can do computations for an embedding in a bigger affine space, considering $\mathcal{H}_J$ in $\AA^{p(\reg(\underline{J}))q(\reg(\underline{J}))}$, then the equations defining $\mathcal{H}_J$ as a subscheme of $\AA^{p(\reg(\underline{J}))q(\reg(\underline{J}))}$ are bounded, as we show in the following theorem.  Furthermore  we can compare computationally two open subsets of this kind.

\begin{theorem}\label{smallemb} 
Consider $J\in \mathcal{B}^n_{p(t)}$.
\begin{enumerate}[(i)]
	\item\label{smallemb-i} $\mathcal{H}_{J}$  can be embedded as a closed subscheme in  $\AA^{p(m)q(m)}$ where  $m$ is any integer $\geqslant \reg(\underline{J})$, by an ideal  generated in degree $\leqslant d+2$.
	\item \label{smallemb-ii} If $\phi_{J_i,r}: \mathcal{H}_{J_i} \rightarrow \AA^{p(r)q(r)} $ are the embedding  for the open subsets corresponding to two Borel-fixed ideals  $J_1$ and $J_2$  belonging to $\mathcal{B}^n_{p(t)}$, then  in general $\phi_{J_1,r}(\mathcal{H}_{J_1}\cap \mathcal{H}_{J_2}) \neq \phi_{J_2,r}(\mathcal{H}_{J_1}\cap \mathcal{H}_{J_2}) $. More precisely:
	\[
	\phi_{J_1,r}(\mathcal{H}_{J_1}\cap \mathcal{H}_{J_2}) =\phi_{J_1,r}(\mathcal{H}_{J_1})\setminus F_1 , \quad \phi_{J_2,r}(\mathcal{H}_{J_1}\cap \mathcal{H}_{J_2}) =\phi_{J_2,r}(\mathcal{H}_{J_2})\setminus F_2,
	\]
	 where $F_1$ e $F_2$ are hypersurfaces of the same degree $\left\vert G_{J_1} \setminus ( G_{J_1} \cap G_{J_2})\right\vert$ in $\AA^{p(r)q(r)}$.
		\item\label{smallemb-iii}   If we consider  $\overline{m}\geqslant \max \{ \reg(\underline{J}_1), \reg(\underline{J}_2)\}$
	 then statement (\ref{smallemb-ii}) holds considering the embedding $\phi_{J_i,\overline{m}}: \mathcal{H}_{J_i} \hookrightarrow \AA^{p(\overline{m})q(\overline{m})} $.
	\end{enumerate}
\end{theorem}
\begin{proof} 
\emph{(\ref{smallemb-i})} Using Theorem \ref{isomHM}, $\mathcal{H}_{J}=\Mf({J})$. Furthermore, thanks to Theorem \ref{isomH} , we have that $\mathcal{H}_{J}\simeq\Mf(\underline{J}_{\geqslant \reg(\underline{J})})$. Applying Theorem \ref{nostrobuch}, it is sufficient to consider reductions of $S$-polynomials in degree $\reg(\underline{J})+1$ and we conclude as in Theorem \ref{gradod+2}.

\emph{(\ref{smallemb-ii})} $F_1$ is defined by the equation of $\frac{\Delta_{J_2}}{\Delta_{J_1}}$ in $\AA^{p(r)q(r)}$  and so its degree corresponds to $\left\vert G_{J_1} \setminus ( G_{J_1} \cap G_{J_2})\right\vert$, by Lemma \ref{lem:pluckerlocali}\emph{(\ref{it:pluckerlocali_2})}.

\emph{(\ref{smallemb-iii})} is a straightforward consequence of \emph{(\ref{smallemb-i})} and \emph{(\ref{smallemb-ii})}.
\end{proof}

\begin{example}
We consider the Hilbert scheme of $s$ points in $\PP^n$, $\Hilb{n}{s}$. For any monomial ideal $J$ , we have that  the open subset of the Grassmannian $\mathcal{U}_J$ is isomorphic to $\AA^{s q(s)}$, where $q(s)=\binom{n+s}{n}-s$.
The saturated lexicographic ideal $\underline{L}=(x_n,\dots,x_2,x_1^{s})$ has regularity $s$ the open subset $\mathcal{H}_{\underline{L}_{\geqslant s}} \subset \Hilb{n}{s}$, the Reeves-Stillman component, contains all the subschemes of $\PP^n$ made up of $s$ distinct points, so it has dimension $\geqslant n s$.  Using Theorem  \ref{isomH}, $\mathcal{H}_{\underline{L}_{\geqslant s}}$ is embedded in an affine space of dimension $\vert G_{\underline{L}}\vert \cdot p(s)= n s$. Then, $\mathcal{H}_{\underline{L}_{\geqslant s}} \simeq \AA^{ns}$.
\end{example}

\begin{example}
We can now easily study some features of \gls{Hilb3tP3}, that we have already investigated in Example \ref{nonBorel}. The Borel region of $\Hilb{3}{3t}$  is made up of one open subset only, corresponding to the lexicographic ideal $L= \underline{L}_{\geqslant 3} =(x_3,x_2^3)_{\geqslant 3}$, as already pointed out in Example \ref{nonBorel}, using Proposition \ref{fuori} and  Theorem \ref{isomHM}. Since no monomial in the basis of $(x_3,x_2^3)$ is divisible by $x_1$, using Theorem  \ref{isomH}, we have that $\Mf(\underline{L})\simeq\mathcal{H}_{L}$ and we can embed $\mathcal{H}_{L}$ in $\AA^{12}$.

Furthermore, since the monomial basis of $\underline{L}$ is made up of two coprime monomials, we have that
every ideal $I$ in $\Mf(\underline{L})$ corresponds to the complete intersection of a plane and a cubic; we then have that $\mathcal{H}_{L}$ has dimension $\geqslant 12$, and so $\mathcal{H}_{L}\simeq \AA^{12}$. Every point of $\Hilb{3}{3t}$ is, up to a change of coordinates,  a point of $\mathcal{H}_{L}$, hence every scheme in $\PP^3$ with Hilbert polynomial $3t$ is a $(1,3)$-complete intersection.  
\end{example}

\begin{example}
Let us consider the Hilbert scheme \gls{Hilb4tP3}. Continuing the computation of marked schemes started in Example \ref{ex:Hilb4tP3}, we obtain that the Gr\"obner stratum $\St{\omega_4}(J_4)$ of $J_{4}=(x_3^2,x_3 x_2,x_2^3)$ is isomorphic to the affine space $\AA^{16}$, namely $\Hilb{3}{4t}$ has a rational component of dimension 16 corresponding to $(2,2)$-complete intersection (called Vainsencher-Avritzer component \cite{VainsencherAvritzer}). A second component is the Reeves-Stillman one, containing the lexicographic point. Being such a point smooth \cite{ReevesStillman}, the Gr\"obner stratum $\St{\DegLex}\big((J_1)_{\geqslant 5}\big)$ has to be an affine scheme, thus isomorphic to $\AA^{23}$ (cf. with \cite{GotzmannHilb4tP3}). 

We know that the point defined by $J_3 = (x_3^2,x_3x_2,x_3x_1^2,x_2^4)$ lies on the intersection of these two components: indeed it belongs to the Reeves-Stillman component because it has the same hyperplane section of the lexicographic ideal (Reeves criterion) and it belongs to the Vainsencher-Avritzer component because there is a Borel rational deformation having $(J_3)_{\geqslant 6}$ and $(J_{4})_{\geqslant 6}$ as fibers. By explicit computation (see Example	\ref{ex:Hilb4tP3_M2}) we find that the ideal of the Gr\"obner stratum $\St{\omega_3}(J_3)$ can be decomposed in an ideal generated by a single variable (of 24), that correspond to an affine scheme of dimension 23, and an ideal defining a affine subscheme isomorphic to $\AA^{16}$. Hence the two components of $\Hilb{3}{4t}$ intersect transversely.
\end{example}

\chapter{Low degree equations defining the Hilbert scheme}\label{ch:LowDegreeEquations}

In this chapter, we will introduce a new set of equations defining the Hilbert scheme $\Hilb{n}{p(t)}$ as subscheme of the suitable Grassmannian $\Grass{q(r)}{N(r)}{\K}$ as already done in Chapter \ref{ch:HilbertScheme}. This topic is placed here, after a long discussion on Borel-fixed ideals, because ideas behind the main result we will prove come from the construction of the Borel open covering of the Hilbert scheme discussed in Chapter \ref{ch:openSubsets}, in particularly Section \ref{sec:localEqd+2}. We remark that it is not trivial to extend the bound on the degree of equations defining locally the Hilbert scheme to the global equations, because as seen in Lemma \ref{lem:pluckerlocali} global coordinates can correspond to polynomials of high degree in the local coordinates. The results of this chapter belong to the preprint \lq\lq Low degree equations defining the Hilbert scheme\rq\rq\ \cite{BLMR}, joint paper with J. Brachat, B. Mourrain and M. Roggero.

\section{BLMR equations}\index{Hilbert scheme!equations of the|(}
Given an admissible Hilbert polynomial $p(t)$ on $\PP^n$ with Gotzmann number $r$ and degree $d$ and the associated volume polynomial $q(t) = \binom{n+t}{n}-p(t)$, we set
\begin{equation}
q'(t)=q(t) - \binom{n-d-1+t}{n-d-1}=\binom{n+t}{n}- \binom{n-d-1+t}{n-d-1} -p(t),
\end{equation}
so that $q(t)-q'(t)=\dim_\K \K[x_{d+1},\ldots, x_{n}]_t$.

\begin{proposition}\label{forsechiara} 
Let  $\mathcal{U}'$ be the set of all the elements $I \in \Grass{q(r)}{\K[x]_r}{\K}$ such that $I_r$ has a set of generators of the type:
\begin{equation}\label{eq:Jbase}
\begin{split} 
G_I^r&{}=\big\{x^{\alpha}+f_{\alpha} \ \vert\ x^\alpha \in \K[x_{d+1},\ldots,x_{n}]_r \text{ and } f_\alpha \in (x_d,\ldots,x_0) \big\}\\
&{}\quad \cup \big\{g_j\ \vert\  g_{j}\in (x_d, \ldots, x_0)\big\}  
\end{split}
\end{equation}
where the first set of generators contains by construction $q(r) - q'(r)$ elements and the second set has $q'(r)$ polynomials, so that $\dim_\K I_r = q(r)$.

Then $\mathcal{U}'$ is a non-empty open subset in $\Grass{q(r)}{\K[x]_r}{\K}$ and $I_{r+1}$ has a set of generators $G_I^{r+1}$ that can be represented by a matrix of the type:  
 \begin{equation}\label{matrice}
 \mathcal{A}_{r+1}=\left( \begin{array}{ccc|ccc}
 & \mathtt{Id} & &\bullet &\bullet &\bullet \\ 
\hline
&  0 &   & & \mathcal{D}_1 &  \\
   \hline 
 & 0  & &    & \mathcal{D}_2 &   
   \end{array}\right)   
  \end{equation}
where:
\begin{itemize}
\item the columns belonging to the left part of the matrix correspond to the monomials in $\K[x_{d+1},\ldots,x_{n}]_{r+1}$ and the columns on the right to the monomials in $(x_0, \ldots, x_d)_{r+1}$;
\item the top-left submatrix  $\mathrm{Id}$ is the identity matrix of order $q(r+1)-q'(r+1)$;  
\item the rows of $\mathcal{D}_1$ contain the coefficients of all the generators multiplied by a variable $x_h,\ h=0,\ldots,d$; 
\item the rows of $\mathcal{D}_2$ contain the coefficients of the generators $g_{j}$ multiplied by a variable $x_h,\ h=d+1,\ldots,n$ and the coefficients of the polynomials $x_{i'}f_{\alpha'}- x_i f_\alpha$ such that $x_{i'}x^{\alpha'}=x_i x^{\alpha}$ and $i,i' \geqslant d+1$.
\end{itemize}

Moreover the subset $\mathcal{U} \subset \mathcal{U}'$ of all the ideals $I_r$ such that $ \rank \mathcal{D}_1\geqslant q'(r+1)$ is open and $\mathcal{U}^{\GL} = \{ g\centerdot\mathcal{U} \ \vert\  g\in \GL(n+1)\}$  is an open covering of $\Grass{q(r)}{\K[x]_r}{\K}$.
\end{proposition}
\begin{proof}  Let us consider the canonical projection
\[
\pi: \K[x_{0}, \ldots, x_n]_r\ \longrightarrow\ \left(\K[x_{0}, \ldots, x_n]/(x_0, \ldots x_d)\right)_r\simeq \K[x_{d+1},\ldots, x_n]_r.
\]
The subset $\mathcal{U}'$ of $\Grass{q(r)}{\K[x]_r}{\K}$ is open because $\mathcal{U}' = \pi^{-1}(\K[x_{d+1}, \ldots, x_n]_r)$ ${} \cap \Grass{q(r)}{\K[x]_r}{\K}$. Moreover $\mathcal{U}'$ is non-empty because any Borel-fixed ideal $J$ defining a point of $\Grass{q(r)}{\K[x]_r}{k}$ (i.e. $\dim_\K J_r = q(r)$) belongs to $\mathcal{U}'$. Indeed, $\dim_\K J_t\geqslant q(t),\ \forall\ t\geqslant r$ (by Macaulay's Estimate on the Growth of Ideals) implies that the Hilbert polynomial of $\K[x_0,\ldots,x_n]/J$ has degree smaller than or equal to $\deg p(t)=d$  and so $\K[x_{d+1},\ldots,x_{n}]_{\geqslant r} \subseteq J$  (see Proposition \ref{prop:degreeHPpowerVars}). Therefore the ideal $J_{r}$ has a basis $G_I^r$ (as $\K$-vector space) as the one described in \eqref{eq:Jbase} so that $J_r \in \mathcal{U}'$.

For any $I_r \in \mathcal{U}'$, a set of generators of $I_{r+1}$ is $\K[x]_1 \cdot G_I^r = \bigcup_i \left\{x_i\, G_I^r\right\}$. The set of generators $G_I^{r+1}$ we are looking for can be easily obtained from $\K[x]_1\cdot G_I^r $ just modifying few elements: for every monomial $x^\gamma$ in $ \K[x_{d+1},\ldots,x_n]_{r+1}$ we choose only one product $x_i(x^{\alpha}+f_\alpha)$ such that  $x^\gamma=x_i\,x^{\alpha}$, to be left in ${G}_I^{r+1}$ (and corresponding to a row in the first block of $\mathcal{A}_{r+1}$), whereas we replace any other polynomial $x_{i'}(x^{\alpha'}+f_{\alpha'})$, such that $x^\gamma=x_{i'}x^{\alpha'}$, by $x_{i'}f_{\alpha'}-x_i f_{\alpha}$ (which belongs to $(x_0, \ldots, x_d)$ and corresponds to a row of $\mathcal{D}_2$).

Obviously the condition $\rank \mathcal{D}_1 \geqslant q'(r+1)$ is an open condition and we call $\mathcal{U} \subset \mathcal{U}'$ the corresponding open subset. Again this open subset is not empty because it contains for instance all the subspaces $J_r$ defined by a Borel-fixed ideal $J\in \Grass{q(r)}{\K[x]_r}{\K}$.

\medskip

To prove the last statement, we consider any term ordering $\sigma$, refinement of $\leq_B$ (i.e. $x_n >_\sigma \cdots >_{\sigma} x_0$). Then for a generic $g \in \GL(n+1)$, $J =(\IN( g\centerdot I)_r)$ is Borel-fixed.  Note
 that $J$ belongs to $\Grass{q(r)}{\K[x]_r}{\K}$, but if $ I\notin \hilb{n}{p(t)}(\K)$, $J$ can  differ from $\IN(g\centerdot I)$. As $J$ is Borel-fixed
\[
\begin{split}
\dim_\K (x_0J_{r}+\ldots+x_dJ_{r})&{}=\dim_\K  J_{r+1} \cap (x_0,\ldots,x_d) = \\
                                 &{}=\dim_\K J_{r+1} - \dim_\K \K[x_{d+1}, \ldots, x_n]_{r+1} \geqslant q'(r+1),
\end{split}
\]
hence
\[
\begin{split}
\dim_\K \left(x_0(g\centerdot I)_{r}+\ldots +x_d(g\centerdot I)_{r}\right) & {}\geqslant \dim_\K \left(x_0\IN(g\centerdot I)_r + \ldots +x_d\IN(g\centerdot I)_r\right) \geqslant {}\\
&{}\geqslant\dim_\K (x_0J_r+\cdots +x_dJ_r )\geqslant q'(r+1).
\end{split}
\] 
Finally we can conclude that $g\centerdot I\in \mathcal{U}$ because a set of generators of the vector space $x_0(g\centerdot I)_{r}+\cdots +x_d(g\centerdot I)_{r}$ corresponds to the rows of $\mathcal{D}_1$.
\end{proof}

Making reference to the matrix $ \mathcal{A}_{r+1}$ in \eqref{matrice}, let us call
\[
\mathcal{D} = \left( \begin{array}{c}
   \mathcal{D}_1  \\
     \mathcal{D}_2    
   \end{array}\right). 
\]

\begin{corollary}\label{necessaria} 
Let $I \subset \K[x]$ be any ideal belonging to $\mathcal{U}$. Then
\begin{enumerate}[(i)]
\item\label{condizione1} $I\in \mathcal{U} \cap \hilb{n}{p(t)}(\K) \quad \Longleftrightarrow \quad \rank \mathcal{D} = \rank \mathcal{D}_1 = q'(r+1)$;
\item\label{condizione2} $I\notin \mathcal{U} \cap \hilb{n}{p(t)}(\K) \quad \Longleftrightarrow\quad \begin{array}{l}\text{either }\rank \mathcal{D}_1 > q'(r+1)\\ \text{or }\rank \mathcal{D}_1 =q'(r+1)<\rank \mathcal{D}\end{array}$.
\end{enumerate}
\end{corollary}
\begin{proof}  \emph{(\ref{condizione1})} follows from the fact that $I\in \hilb{n}{p(t)}(\K)$ if and only if $\dim_\K I_{r+1}=\rank \mathcal{A}_{r+1} =q(r+1)$ and from the special form of $\mathcal{A}_{r+1}$. 

\emph{(\ref{condizione2})} is another way to write \emph{(\ref{condizione1})}. 
\end{proof}

We can now describe the set of equations $\mathfrak{H}$ that we will prove to describe the Hilbert scheme $\Hilb{n}{p(t)}$. As usual $r$ is the Gotzmann number of $p(t)$ and let $I$ be an element of $\Grass{q(r)}{\K[x]_r}{\K}$. We will describe the equations of $\mathfrak{H}$ in terms of the Pl\"ucker coordinates of $I$ via the  Pl\"ucker embedding \eqref{eq:pluckerEmbedSubspace}  $\mathscr{P}: \Grass{q(r)}{N}{\K} \rightarrow \PP\big[\wedge^{q(r)} \K[x]_r\big]$.

A first subset $\mathfrak{H}'$ of the equations is construct as follows.
Let us choose in all the possible ways a set of $d+1$ elements of the type $x_i \Lambda^{(s_i)}_{\JJ_i}$ such that $i = 0,\ldots,d$ and $\sum_i s_i = q'(r+1)$. Moreover let us consider any variable $x_j,\ j=d+1,\ldots,n$ and any multiindex $\KK = \{\kk_1,\ldots,\kk_{q(r)-1}\}$, such that all the monomials belonging to $\K[x_{d+1}, \ldots, x_n]_r$ are contained in the set $\{x^{\alpha(\kk_1)},\ldots,x^{\alpha(\kk_{q(r)-1})}\}$ (in order for $\Lambda^{(1)}_{\KK}$ to contain only monomials in $(x_0,\ldots,x_d)$).
The first part $\mathfrak{H}'$ of the polynomials in $\mathfrak{H}$ is represented by all the coefficients of all exterior products of the type
\begin{equation}\label{secondo}
\left(\bigwedge_{i=0}^d x_i\, \Lambda^{(s_i)}_{\JJ_i} \right) \wedge x_j\, \Lambda^{(1)}_{\KK}.
\end{equation}

The second part of the equations $\mathfrak{H}''$ is constructed as follows. 
We choose in all the possible ways a multiindex $\HH = \{\hh_1,\ldots,\hh_{q(r)}\}$, such that again the corresponding set of monomials $\{x^{\alpha(\hh_1)},\ldots,x^{\alpha(\hh_{q(r)})}\}$ of degree $r$ contains the monomials in $\K[x_{d+1}, \ldots, x_n]_r$. For any $x^\gamma \in \K[x_{d+1}, \ldots, x_n]_{r-1}$ and any couple of variables $x_j,x_{j'} \in \{x_{d+1},\ldots,x_n\}$, let $\overline{\hh}$ and $\overline{\hh}'$ be the indices such that $x_j x^{\gamma} = x^{\alpha(\overline{\hh})}$ and $x_{j'} x^{\gamma} = x^{\alpha(\overline{\hh}')}$.  The second part $\mathfrak{H}''$ of the polynomials in $\mathfrak{H}$ can be obtained collecting the coefficients of all the exterior products of the type
\begin{equation}\label{terzo} 
\left(\bigwedge_{i=0}^d x_i\, \Lambda^{(s_i)}_{\JJ_i} \right) \wedge \left(x_{j'}\, \Lambda^{(1)}_{\HH_1} +x_{j}\, \Lambda^{(1)}_{\HH_2} \right) 
\end{equation}
where $\HH_1=\HH\setminus \{\overline{\hh}\}$ and $\HH_2=\HH\setminus \{\overline{\hh}'\}$. Note that we need to sum the two vectors in order for the monomials $x_{j'}x^{\alpha(\overline{h})}$ and $x_j x^{\alpha(\overline{h}')}$ to delete, because $\Delta_{\HH_1 \vert \overline{h}} = - \Delta_{\HH_2 \vert \overline{h}'}$.

\begin{theorem}\label{th:degreeEqField} Let $p(t)$ be an admissible Hilbert polynomial in $\PP^n$ of degree $d$ and Gotzmann number $r$. The set $\hilb{n}{p(t)}(\K)$ of $\K$-valued point of $\hilb{n}{p(t)}$ can be described by a closed subscheme of the Grassmannian $\Grass{q(r)}{\K[x]_r}{\K}$ defined by the equations $\mathfrak{H}$ given in \eqref{secondo}, \eqref{terzo} of degree smaller than or equal to $d+2$. 
\end{theorem}

\begin{proof}  
We divide the proof in two steps.

\noindent\textbf{Step 1.} Firstly, we consider an ideal $I \in \mathcal{U}$ and show that every polynomial in $\mathfrak{H}$ vanishes on the Pl\"ucker coordinates $[\ldots:\Delta_{\II}(I):\ldots]$ of $I$ if and only if $I \in \mathcal{U}\cap \hilb{n}{p(t)}(\K)$. 

Note that the conditions given by the vanishing of the polynomials in $\mathfrak{H} = \mathfrak{H}' \cup \mathfrak{H}''$  mean that $q'(r+1)$ rows in the matrix $\mathcal{D}_1$ and one rows in the matrix $\mathcal{D}_2$ are linearly dependent (directly by the construction of the matrix $\mathcal{A}_{r+1}$ in Proposition \ref{forsechiara}). These conditions ensure also that $q'(r+1)+1$ rows of the matrix $\mathcal{D}_1$ are dependent, because of the well-known property of vector spaces saying that $s+1$ vectors $v_1,\ldots,v_{s+1}$, such that every subset of $s$ elements is linearly dependent with any other vector $u \neq 0$, are dependent. The only delicate issue, that we will discuss later in Remark \ref{rk:gNotEmpty}, is checking that $\mathcal{D}_2$ is not a zero matrix.
 
By Corollary \ref{necessaria}, $I \in \mathcal{U}$ belongs to $\hilb{n}{p(t)}(\K)$ if and only if the polynomials of $\mathfrak{H}$ vanish on $[\ldots:\Delta_\II(I):\ldots]$. Note that the coefficients of the exterior products in \eqref{secondo} and \eqref{terzo} are polynomials of degree $\leqslant d+2$ (more precisely the degree is the number of non-zero $s_i$).

\medskip

\noindent\textbf{Step 2.}
Let $I$ be an element of $\Grass{q(r)}{\K[x]_r}{\K}$ and $g = (g_{i,j})$ be an element of $\GL(n+1)$. 
The Pl\"ucker coordinates $[\cdots:\Delta_{\II}(g\centerdot I):\cdots]$ of $g\centerdot I \in \Grass{q(r)}{\K[x]_r}{\K}$ are bi-homogeneous polynomials of degree 1 in the Pl\"ucker coordinates $[\cdots:\Delta_\II(I):\cdots]$ and of degree $q(r)\cdot r$ in the coefficients $g_{i,j}$ of the matrix $g$.
So given a homogeneous polynomial $P$ of degree $s \leqslant d+2$ in $\mathfrak{H}$, $P([\cdots:\Delta_{\II}(g\centerdot I):\cdots])$ is a bi-homogeneous polynomial of degree $s$ in $[\cdots:\Delta_{\II}(I):\cdots]$ and of degree $q(r)\cdot r \cdot s$ in $g_{i,j}$. 
At this point we collect, and denote by $C_P$, the homogeneous polynomials of degree $s \leqslant d+2$ in the Pl\"ucker coordinates $[\cdots:\Delta_{\II}(I):\cdots]$, that spring up as coefficients of $P([\cdots:\Delta_{\II}(g\centerdot I):\cdots])$, viewed as a homogeneous polynomial of degree $q(r)\cdot r \cdot s$ in the variables $g_{i,j}$.

From Proposition \ref{forsechiara} and Corollary \ref{necessaria}, $I$ belongs to $\hilb{n}{p(t)}(\K)$ if and only if for a generic changes of variables $g \in \GL(n+1)$
\[
g\centerdot I \in \mathcal{U} \cap \hilb{n}{p(t)}(\K)
\]
i.e.
\[
\text{all the homogeneous polynomials } P \in \mathfrak{H} \text{ vanish at } [\cdots:\Delta_{\II}(g \centerdot I):\cdots],
\] 
or equivalently 
\[
\text{all the coefficients}\ C_{P} \text{ for } P \in \mathfrak{H} \text{ vanish at } [\cdots:\Delta_{\II}(I):\cdots].
\] 
We finally proved that $I \in \Grass{q(r)}{\K[x]_r}{\K}$ belongs to $\hilb{n}{p(t)}(\K)$ if and only if $[\cdots:\Delta_{\II}(I):\cdots]$ satisfies all the equations of the set 
\begin{equation}
\bigcup_{P \in \mathfrak{H}} C_{P}
\end{equation}
which consists of homogeneous polynomial of degree smaller than or equal to $d+2$.
\end{proof}

\begin{remark}\label{rk:gNotEmpty}
Note that a necessary condition in order that $\mathcal{D}_2$ is empty is that $I_r$ has no generators belonging to $G^r_I$ of the type $g_j$ and now we prove that it is not possible. For the sake of simplicity, we can think about the monomial ideal obtained in the case $f_\alpha = g_j = 0,\ \forall\ \alpha,\ \forall\ j$. Such an ideal should have a Hilbert polynomial $\widetilde{p}(t)$ such that $\widetilde{p}(r) = \binom{n+r}{n}-\binom{n-d+r}{n-d}$. Let us show that $\widetilde{p}(r)$ can not be equal to $P(r)$ with $P$ a Hilbert polynomial with Gotzmann number equal to $r$ and of degree $d$.

The first point is to compute the maximal value in degree $r$ of a Hilbert polynomial of degree $d$ and Gotzmann number $r$. By the Gotzmann decomposition of Hilbert polynomials \eqref{eq:GotzmannDecomposition}, we know that among the Hilbert polynomials of degree $d$ and Gotzmann number $r$ there is\index{Hilbert polynomial!Gotzmann's representation of a}
\[
p(t) = \binom{t+d}{d} + \binom{t+d-1}{d} + \cdots + \binom{t+d-(r-1)}{d},
\]
and\hfill that\hfill any\hfill other\hfill Hilbert\hfill polynomial\hfill has\hfill at\hfill least\hfill one\hfill binomial\hfill coefficient\\ $\binom{t+(d-i)-j}{d-i}$ replacing $\binom{t+d-j}{d}$ $(i,j > 0)$. Because of $\binom{r+d-j}{d} > \binom{r+(d-i)-j}{d-i}$, the maximal value reached is
\[
P = \max \left\{P(r)\ \vert\ P(t)\text{ of degree } d \text{ and Gotzmann number } r\right\} = \sum_{i=1}^r \binom{d+i}{d}.
\]

Finally, starting from the decomposition
\[
\K[x_0,\ldots,x_n]_r = \bigcup_{i=0,\ldots,r} \K[x_0,\ldots,x_d]_i \cdot \K[x_{d+1},\ldots,x_n]_{r-i}, 
\]
we have that
\[
\begin{split}
\widetilde{p}(r) & {} = \binom{n+r}{n} - \binom{n-d+r+1}{n-d+1} =\\ &{} = \sum_{i=1}^r \binom{d+i}{d}\binom{n-d+1+r-i}{n-d+1} > \sum_{i=1}^r \binom{d+i}{d} = P.
\end{split}
\]
\end{remark}

\section{Extension of the coefficient ring}

We will now prove that the subscheme of $\Grass{q(r)}{\K[x]_r}{\K}$ described in Theorem \ref{th:degreeEqField} representing the set $\hilb{n}{p(t)}(\K)$\index{Hilbert functor} is indeed the Hilbert scheme we are looking for, so let us denote it by $\Hilb{n}{p(t)}$. We have seen in Chapter \ref{ch:HilbertScheme} that any element $Z \in \hilb{n}{p(t)}(A)$ defines an element of $\grass{q(r)}{N(r)}(A)$, that is, being $\grass{q(r)}{N(r)}$\index{Grassmann functor} represented by $\Grass{q(r)}{N(r)}{\K}$, a morphism $f_{Z} : \Spec A \rightarrow \Grass{q(r)}{N(r)}{\K}$. To prove that $\Hilb{n}{p(t)}$ represents $\hilb{n}{p(t)}$, we have to show that to any element $Z \in \hilb{n}{p(t)}(A)$ we can associate a morphism $\overline{f}_Z: \Spec A \rightarrow \Hilb{n}{p(t)}$. To accomplish this task we prove that the morphism $f_Z$ factors through the inclusion $\Hilb{n}{p(t)} \hookrightarrow \Grass{q(r)}{N(r)}{\K}$ as subschemes of $\PP \big[ \wedge^{q(r)} \K[x]_r\big]$:
\begin{center}
\begin{tikzpicture}[scale=0.9]
\node (A) at (0,0) [] {$\Spec A$};
\node (B) at (5,0) [] {$\Grass{q(r)}{N(r)}{\K}$};
\node (C) at (2.5,-2) [] {$\Hilb{n}{p(r)}$};
\draw [->] (A) -- node[rectangle,fill=white,inner sep=1pt]{\footnotesize $f_Z$} (B);
\draw [->,dashed] (A) -- node[rectangle,fill=white,inner sep=1pt]{\footnotesize $\overline{f}_Z$} (C);
\draw [right hook->] (C) -- (B);
\node at (2.5,-0.8) [] {$\circlearrowleft$};
\end{tikzpicture}
\end{center}

\begin{theorem}\label{th:eqLocalRing}
The Hilbert scheme $\Hilb{n}{p(t)}$ can be defined as subscheme of the Grassmannian $\Grass{q(r)}{N(r)}{\K}$ by equations of degree equal to or less than $d+2$.
\end{theorem}
\begin{proof}
For any element of $Z \in \hilb{n}{p(t)}(A)$ we consider the submodule $ I(Z)_r \subset A[x_0,\ldots,x_n]_r$ that belongs to $\grass{q(r)}{N(r)}(A)$. To prove that the map $f_Z : \Spec A \rightarrow \Grass{q(r)}{N(r)}{\K}$ factors through $\Hilb{n}{p(t)}$, we now prove that an element of $I \in \grass{q(r)}{N(r)}(A)$ belongs to $\hilb{n}{p(t)}(A)$ if and only if the equations \eqref{secondo} and \eqref{terzo} extended to the projective space $\PP^{\binom{N(r)}{q(r)}-1}_A$ by means of the tensor product for the $\K$-algebra $A$ $\K[x] \otimes_{\K} A = A[x]$ are satisfied by the Pl\"ucker coordinates of $I$ in the Grasmmannian $\Grass{q(r)}{N(r)}{A} = \Grass{q(r)}{N(r)}{\K} \times_{\Spec \K} \Spec A$.

In Chapter 1 we discussed about the fact that the flatness\index{flatness} is a local property, therefore we can consider $A$ a local $\K$-algebra with maximal ideal $\mathfrak{m}$ and residue field $k(\mathfrak{m})$.

\noindent\textbf{Step 1.} Given $I_{r} \in \grass{q(r)}{N(r)}(A)$, firstly let us prove that if equations given in Theorem \ref{th:degreeEqField} are satisfied, then $I_{r}$ belongs to $\hilb{n}{p(t)}(A)$ (i.e. $A[x]_{1}\cdot I_{r}$ belongs to $\grass{q(r+1)}{N(r+1)}(A)$ according to Gotzmann's Persistence Theorem.

Let us consider $I_{r} \in \grass{q(r)}{N(r)}(A)$ satisfying the extension of equations \eqref{secondo} and \eqref{terzo}. Tensoring by the residue field $k(\mathfrak{m})$ and using Nakayama's Lemma, we can determine a free submodule $J_{r} \subset A[x]_{r}$ generated by $q(r)$ monomials having the Borel-fixed property, such that with a generic change of coordinates, the monomials $N(J_r) = \{ x^\beta \in A[x]_r\ \vert\ x^\beta \notin J_r\}$ form a basis of $A[x]_{r}/I_{r}$ as a free $A$-module (see \cite[Chapter 15]{Eisenbud}). Now we consider the exact sequence
\[
0\ \longrightarrow\ I_{r}\ \longrightarrow\ A[x]_{r}\ \longrightarrow\ A[x]_{r}/I_{r}\ \longrightarrow\ 0,
\]
and we tensor it by the residue field $\resField{m}$, obtaining
\[
I_{r} \otimes_{\K} \resField{m}\ \longrightarrow\ \resField{m}[x]_{r} = \K[x]_r \otimes_{\K} \resField{m}\ \longrightarrow\ A[x]_{r}/I_{r}\otimes_{\K} \resField{m}\ \longrightarrow\ 0.
\]
Called $I_{r}^{\resField{m}}$ the image of $I_{r} \otimes_{\K} \resField{m}$ in $\resField{m}[x]$, by the assumptions and by Theorem \ref{th:degreeEqField} we deduce that $I_{r}^{\resField{m}}$ belongs to $\hilb{n}{p(t)}(\resField{m})$. Consequently, $J_{r}^{\resField{m}}$ (resp. $J_r$) also belongs to $\hilb{n}{p(t)}(\resField{m})$ (resp. $\hilb{n}{p (t)} (A)$) and thus $(J_{r}^{\resField{m}})$ (resp. $(J_{r})$) defines a Borel-fixed ideal with Hilbert polynomial $p(t)$ in $\resField{m}[x]$ (resp. in $A[x]$).

For a generic change of coordinates, as $N(J_r)$ is a basis of $A[x]_r/I_{r}$, $I_{r}$ is also a free $A$-module of rank $q(r)$ with a basis of the form:
\begin{equation}\label{eq:genIr}
\left\{x^{\alpha} - \sum_{x^{\beta}\in N(J_r)}c_{\alpha\beta}x^{\beta}\quad\Big\vert\quad x^{\alpha} \in J_{r}\right\}.
\end{equation}
Therefore we can choose a system of generators for $A[x]_{1}\cdot I_{r}$ equal to that one described with the matrix $\mathcal{A}_{r+1}$ \eqref{matrice} in Proposition \ref{forsechiara}. Up to a change of coordinates, the $A$-module $\langle \mathcal{D}_{1} \rangle$ generated by the lines of $\mathcal{D}_{1}$ (and by definition equal to $x_{0}I_{r}+\cdots+x_{d}I_{r}$) contains a family $\mathcal{F}$ of $q'(r+1)$ polynomials of the form:
\begin{equation}\label{eq:genD1}
\mathcal{F}=\left\{x^{\gamma} - \sum_{x^{\eta}\in N(J_{r+1})}c_{\gamma\eta}x^{\eta}\quad\Big\vert\quad x^{\gamma} \in x_{0}J_{r}+\cdots+x_{d}J_{r}\right\}.
\end{equation}
Let us prove it by induction on $0 \leqslant i \leqslant d$ for $x^{\alpha} \in x_{i}J_{r}$.
Let $x^{\gamma}=x_{0} x^{\alpha}$ with $x^{\alpha}$ in $J_{r}$. Among the generators \eqref{eq:genIr} of $I_r$ there is
\[
x^{\alpha} - \sum_{x^{\beta}\in N(J_r)} c_{\alpha\beta}x^{\beta},
\]
thus
\[
x_0 x^{\alpha} - \sum_{x^{\beta}\in N(J_r)}c_{\alpha\beta}x_0 x^{\beta}
\]
belongs to $x_{0}I_{r} \subset x_{0}I_{r}+\ldots+x_{d}I_{r}$. Moreover because of $(J_{r})$ is Borel-fixed and $(J_{r+1}:A[x]_1) = J_r$ ($J_r \in \hilb{n}{p(t)}(A)$), it is easy to check that $x_{0}N(J_r) \subset N(J_{r+1})$. Hence the assertion is proved for $i=0$ and let us suppose that it holds for all $0\leqslant j < i$. Considered $x^{\gamma}=x_{i}x^{\alpha}$ with $x^{\alpha} \in J_{r}$, again
\[
x_i x^{\alpha} - \sum_{x^{\beta}\in N(J_r)}c_{\alpha\beta}x_i x^{\beta} \in x_{i}I_{r} \subset x_{0}I_{r}+\ldots+x_{d}I_{r}.
\]
If $x_{i}x^{\beta}\ (x^{\beta} \in N(J_r))$ does not belong to $N(J_{r+1})$, then there exists $x^\epsilon \in J_r$ such that 
\[
x_{i}x^{\beta} = x_{j} x^{\epsilon}
\]
and $j < i$ because of the Borel-fixed property. Then, by induction, we can replace $x_{i}x^{\beta} = x_{j}x^{\epsilon}$ with an element of the $A$-module generated by $N(J_{r+1})$ modulo $x_{0} I_{r}+\cdots+x_{d}I_{r}$, finally proving that the family $\mathcal{F}$ described in \eqref{eq:genD1} belongs to  $\langle \mathcal{D}_{1} \rangle$ (i.e. to $x_{0}I_{r}+\ldots+x_{d}I_{r}$).

As equations of Theorem \ref{th:degreeEqField} are satisfied, equations \eqref{secondo} and \eqref{terzo} are also satisfied for a generic change of coordinates, so that we can assume without loss of generality that there exist $J_{r}$ Borel-fixed and $\mathcal{F}$ as in \eqref{eq:genD1}, such that $I_{r}$ satisfies \eqref{secondo} and \eqref{terzo}.

Now we want to show that equations \eqref{secondo} and \eqref{terzo} imply that $\mathcal{F}$ generates the $A$-module $\langle \mathcal{D}_{2} \rangle$ spanned by the lines of $\mathcal{D}_{2}$. As a matter of fact, equations \eqref{secondo} and \eqref{terzo} imply that the exterior product between $q'(r+1)$ polynomials in $\langle \mathcal{D}_{1} \rangle$ and one polynomial in $\langle \mathcal{D}_{2} \rangle$ always vanishes. In particular, the exterior product between the $q'(r+1)$ polynomials that belong to $\mathcal{F}$ and any polynomial $g$ in $\langle \mathcal{D}_{2} \rangle$ is equal to zero. We deduce easily that $g$ belongs to $\langle \mathcal{F} \rangle$ and that $\mathcal{F}$ generates $\langle \mathcal{D}_{2} \rangle $. 

Moreover $\mathcal{F}$ generates $\langle \mathcal{D}_{1} \rangle$. With the same reasoning used in the proof of Theorem \ref{th:degreeEqField} and in Remark \ref{rk:gNotEmpty}, it is easy to prove that any exterior product between $q'(r+1)+1$ polynomials in $\langle \mathcal{D}_{1} \rangle$ is equal to zero. In particular, the exterior product between the $q'(r+1)$ polynomials that belong to $\mathcal{F}$ and any polynomial $g$ in $\langle \mathcal{D}_{1} \rangle$ is equal to zero. So again $g$ belongs to $\langle \mathcal{F} \rangle$ and $\mathcal{F}$ generates $\langle \mathcal{D}_{1} \rangle$ (keeping in mind that the free $A$-module $A[x]_r$ has a basis that contains $\mathcal{F}$).

Finally, we conclude that $I_{r+1}$ is a free $A$-module with basis $\mathcal{F}$ plus the polynomials represented by the lines in the first rows of $\mathcal{A}_{r+1}$ and rewriting this family of polynomials using linear combinations of elements in $\mathcal{F}$ we can obtain a basis of the form
\begin{equation}\label{eq:genIr1}
\left\{x^{\gamma} - \sum_{x^{\eta}\in N(J_{r+1})}c_{\gamma\eta}x^{\eta} \quad\Big\vert\quad x^{\gamma} \in J_{r+1} \right\}.
\end{equation}
$A[x]_{r+1}/I_{r+1}$ turns out to be an $A$-module with basis $N(J_{r+1})$, so $I_{r} \in \hilb{n}{p (t)}(A)$.

\medskip

\noindent\textbf{Step 2.}
Let us suppose that $I_{r} \in \grass{q(r)}{N(r)}(A)$ belongs to $\hilb{n}{p(t)}(A)$ and let us prove that it satisfies the extension of equations given in Theorem \ref{th:degreeEqField}. This is equivalent to prove that equations \eqref{secondo} and \eqref{terzo} are satisfied for a generic changes of coordinates. From the Generic Initial Theorem \cite[Theorem 15.20]{Eisenbud} and Nakayama's Lemma, there exists a Borel-fixed monomial ideal $J$ with Hilbert polynomial $p(t)$, such that for a generic changes of coordinates, $N(J_{r})$ and $N(J_{r+1})$ are a basis of respectively $A[x]_{r}/I_{r}$ and $A[x]_{r+1}/I_{r+1}$ as free $A$-modules. As mentioned in \textbf{Step 1}, we can represent $I_{r+1}$ with the matrix $\mathcal{A}_{r+1}$ introduced in Proposition \ref{forsechiara} and find a family $\mathcal{F}$ in $\langle \mathcal{D}_{1} \rangle$ of the form \eqref{eq:genD1}.

As $N(J_{r+1})$ is a basis of $A[x]_{r+1}/I_{r+1}$ as a free $A$-module, every polynomial given by a line in $\mathcal{D}_{1}$ or $\mathcal{D}_{2}$ belongs to  $\langle \mathcal{F} \rangle $. Therefore equations \eqref{secondo} and \eqref{terzo} are satisfied for a generic change of coordinates and equations $\mathfrak{H}$ of Theorem \ref{th:degreeEqField} are satisfied.
\end{proof}
\index{Hilbert scheme!equations of the|)}

\begin{example}\label{ex:ex:MainExampleBLMR}
Let us apply Theorem \ref{th:degreeEqField} to \gls{Hilb2P2}, already considered in Example \ref{ex:MainExampleGotzmann}, Example \ref{ex:MainExampleIarrobinoKleiman} and Example \ref{ex:MainExampleBayerHaimanSturmfels}. Since the Hilbert polynomial is constant, to compute the first set of equations $\mathfrak{H}'$ \eqref{secondo}, we have to consider the wedge product between $x_0\Lambda^{(4)}_{\emptyset}$ and the 2 elements $x_2\Lambda^{(1)}_{123},x_1\Lambda^{(1)}_{123}$ not containing monomials in $\K[x_2,x_1]_2$. We obtain 12 polynomials:
\small
\[
\begin{array}{l l l}
\bullet\ -\Delta_{26}\Delta_{46}+\Delta_{45}\Delta_{46}+\Delta_{16}\Delta_{56}, & \qquad & \bullet\ -\Delta_{45}^2+\Delta_{25}\Delta_{46}-\Delta_{15}\Delta_{56},\\
\bullet\ -\Delta_{24}\Delta_{46}+\Delta_{14}\Delta_{56}, & & \bullet\ \Delta_{34}\Delta_{45}+\Delta_{23}\Delta_{46}-\Delta_{13}\Delta_{56},\\
\bullet\ -\Delta_{24}\Delta_{45}+\Delta_{12}\Delta_{56}, & & \bullet\ \Delta_{14}\Delta_{45}-\Delta_{12}\Delta_{46},\\\bullet\ -\Delta_{36}\Delta_{46}+\Delta_{26}\Delta_{56}+\Delta_{45}\Delta_{56}, & & \bullet\ \Delta_{35}\Delta_{46}-\Delta_{25}\Delta_{56},\\
\bullet\ -\Delta_{45}^2-\Delta_{34}\Delta_{46}+\Delta_{24}\Delta_{56}, & & \bullet\ \Delta_{35}\Delta_{45}-\Delta_{23}\Delta_{56},\\
\bullet\ -\Delta_{25}\Delta_{45}+\Delta_{23}\Delta_{46}, & & \bullet\ \Delta_{15}\Delta_{45}-\Delta_{13}\Delta_{46}+\Delta_{12}\Delta_{56}.\\
\end{array}
\]
\normalsize
To compute the second set of equations $\mathfrak{H}''$ \eqref{terzo}, we have to consider the coefficients of the wedge product between $x_0\Lambda^{(4)}_{\emptyset}$ and one of the  following elements
\[
\begin{array}{l c l}
x_1\Lambda^{(1)}_{234}+x_2\Lambda^{(1)}_{134}, && x_1\Lambda^{(1)}_{134}+x_2\Lambda^{(1)}_{124},\\
x_1\Lambda^{(1)}_{235}+x_2\Lambda^{(1)}_{135}, && x_1\Lambda^{(1)}_{135}+x_2\Lambda^{(1)}_{125},\\
x_1\Lambda^{(1)}_{236}+x_2\Lambda^{(1)}_{136}, && x_1\Lambda^{(1)}_{136}+x_2\Lambda^{(1)}_{126}.
\end{array}
\]
We obtain other 36 generators:
\small
\[
\begin{array}{lll}
\bullet\ \Delta_{16}\Delta_{25}-\Delta_{15}\Delta_{26}-\Delta_{24}\Delta_{26}+\Delta_{14}\Delta_{36}, & & \bullet\ \Delta_{24}\Delta_{25}-\Delta_{14}\Delta_{35},\\
\bullet\ -\Delta_{15}\Delta_{24}-\Delta_{24}^2+\Delta_{14}\Delta_{25}+\Delta_{14}\Delta_{34}, && \bullet\ \Delta_{15}\Delta_{23}+\Delta_{23}\Delta_{24}-\Delta_{13}\Delta_{25},\\
\bullet\ -\Delta_{14}\Delta_{23}+\Delta_{12}\Delta_{25}, && \bullet\ \Delta_{13}\Delta_{14}-\Delta_{12}\Delta_{15}-\Delta_{12}\Delta_{24},\\
\bullet\ -\Delta_{25}\Delta_{26}-\Delta_{26}\Delta_{34}+\Delta_{16}\Delta_{35}+\Delta_{24}\Delta_{36},
&& \bullet\ \Delta_{25}^2+\Delta_{25}\Delta_{34}-\Delta_{15}\Delta_{35}-\Delta_{24}\Delta_{35},\\
\bullet\ -\Delta_{24}\Delta_{25}+\Delta_{14}\Delta_{35},
&& \bullet\ \Delta_{23}\Delta_{25}+\Delta_{23}\Delta_{34}-\Delta_{13}\Delta_{35},\\
\bullet\ -\Delta_{23}\Delta_{24}+\Delta_{12}\Delta_{35},
&& \bullet\ \Delta_{13}\Delta_{24}-\Delta_{12}\Delta_{25}-\Delta_{12}\Delta_{34},\\
\bullet\ -\Delta_{24}\Delta_{46}+\Delta_{14}\Delta_{56},
&& \bullet\ \Delta_{16}\Delta_{25}-\Delta_{15}\Delta_{26}+\Delta_{24}\Delta_{45},\\
\end{array}
\]
\[
\begin{array}{lll}
\bullet\ -\Delta_{16}\Delta_{24}+\Delta_{14}\Delta_{26}-\Delta_{14}\Delta_{45},
&& \bullet\ \Delta_{16}\Delta_{23}-\Delta_{13}\Delta_{26}-\Delta_{24}\Delta_{34}+\Delta_{14}\Delta_{35},\\
\bullet\ \Delta_{24}^2-\Delta_{14}\Delta_{25}+\Delta_{12}\Delta_{26},
&& \bullet\ \Delta_{14}\Delta_{15}-\Delta_{12}\Delta_{16}-\Delta_{14}\Delta_{24},\\
\bullet\ -\Delta_{26}^2+\Delta_{16}\Delta_{36}-\Delta_{34}\Delta_{46}+\Delta_{24}\Delta_{56},
&& \bullet\ \Delta_{25}\Delta_{26}-\Delta_{15}\Delta_{36}+\Delta_{34}\Delta_{45},\\
\bullet\ -\Delta_{24}\Delta_{26}+\Delta_{14}\Delta_{36}-\Delta_{24}\Delta_{45},
&& \bullet\ \Delta_{23}\Delta_{26}-\Delta_{34}^2+\Delta_{24}\Delta_{35}-\Delta_{13}\Delta_{36},\\
\bullet\ -\Delta_{24}\Delta_{25}+\Delta_{24}\Delta_{34}+\Delta_{12}\Delta_{36},
&& \bullet\ \Delta_{15}\Delta_{24}-\Delta_{12}\Delta_{26}-\Delta_{14}\Delta_{34},\\
\bullet\ \Delta_{26}^2-\Delta_{16}\Delta_{36}-\Delta_{25}\Delta_{46}+\Delta_{15}\Delta_{56},
&& \bullet\ -\Delta_{25}\Delta_{26}+\Delta_{16}\Delta_{35}+\Delta_{25}\Delta_{45},\\
\bullet\ \Delta_{24}\Delta_{26}-\Delta_{16}\Delta_{34}-\Delta_{15}\Delta_{45},
&& \bullet\ -\Delta_{23}\Delta_{26}-\Delta_{25}\Delta_{34}+\Delta_{15}\Delta_{35},\\
\bullet\ \Delta_{16}\Delta_{23}-\Delta_{15}\Delta_{25}+\Delta_{24}\Delta_{25},
&& \bullet\ \Delta_{15}^2-\Delta_{13}\Delta_{16}-\Delta_{14}\Delta_{25}+\Delta_{12}\Delta_{26},\\
\bullet\ -\Delta_{35}\Delta_{46}+\Delta_{25}\Delta_{56},
&& \bullet\ \Delta_{26}\Delta_{35}-\Delta_{25}\Delta_{36}+\Delta_{35}\Delta_{45},\\
\bullet\ -\Delta_{26}\Delta_{34}+\Delta_{24}\Delta_{36}-\Delta_{25}\Delta_{45},
&& \bullet\ \Delta_{25}\Delta_{35}-\Delta_{34}\Delta_{35}-\Delta_{23}\Delta_{36},\\
\bullet\ -\Delta_{25}^2+\Delta_{23}\Delta_{26}+\Delta_{24}\Delta_{35},
&& \bullet\ \Delta_{15}\Delta_{25}-\Delta_{13}\Delta_{26}-\Delta_{14}\Delta_{35}+\Delta_{12}\Delta_{36}.\\
\end{array}
\] 

\normalsize
Now we need to introduce the action of $\GL_{\K}(3)$ on $\K[x_0,x_1,x_2]$ and to understand how the induced action on $\Grass{4}{6}{\K}$ works. Given an element $g = (g_{ij}) \in \GL_{\K}(3)$ and its action
\[
\left(\begin{array}{c c c}
g_{00} & g_{01} & g_{02}\\ g_{10} & g_{11} & g_{12}\\ g_{20} & g_{21} & g_{22}
\end{array}\right) \centerdot \left(\begin{array}{c}
x_0 \\ x_1\\ x_2
\end{array}\right),
\]
the induced action on the $\K$-vector space $\K[x_0,x_1,x_2]_2$ is represented by the matrix
\footnotesize
\[
\left(\begin{array}{c c c c c c}
g_{00}^2 & 2g_{00}g_{01} & 2g_{00}g_{02} & g_{01}^2  & 2g_{01}g_{02} & g_{02}^2\\
g_{00}g_{10} & g_{01}g_{10} + g_{00}g_{11} & g_{02}g_{10} + g_{00}g_{12} & g_{01}g_{11}  & g_{02}g_{11} + g_{01}g_{12} & g_{02}g_{12}\\
g_{00}g_{20} & g_{01}g_{20} + g_{00}g_{21} & g_{02}g_{20} + g_{00}g_{22} & g_{01}g_{21}  & g_{02}g_{21} + g_{01}g_{22} & g_{02}g_{22}\\
g_{10}^2 & 2g_{10}g_{11} & 2g_{10}g_{12} & g_{11}^2 & 2g_{11}g_{12} & g_{12}^2\\
g_{10}g_{20} & g_{11}g_{20} + g_{10}g_{21} & g_{12}g_{20} + g_{10}g_{22}  & g_{11}g_{21} & g_{12}g_{21} + g_{11}g_{22} & g_{12}g_{22}\\
g_{20}^2 & 2g_{20}g_{21} & 2g_{20}g_{22} & g_{21}^2 & 2g_{21}g_{22} & g_{22}^2\\
\end{array}\right) \centerdot \left(\begin{array}{c} x_0^2 \\ x_0 x_1 \\ x_0 x_2 \\ x_1^2 \\ x_1 x_2 \\ x_2^2 \end{array}\right). 
\]
\normalsize
To write explicitly the action of $g$ on the Pl\"ucker coordinates, we can consider the element $\Lambda^{(4)}_{123456}$, substitute each element $x^\beta$ of the basis of $\K[x_0,x_1,x_2]_2$ with $g \centerdot x^\beta$ obtaining the exterior product $g\centerdot \Lambda^{(4)}_{123456}$ and then the action is determined by looking at the coefficients of the same element of the basis of $\wedge^4 \K[x_0,x_1,x_2]_2$ in $\Lambda^{(4)}_{123456}$ and $g\centerdot \Lambda^{(4)}_{123456}$. For instance the Pl\"ucker coordinate $\Delta_{25}$, coefficient of $x_2^2 \wedge x_1^2 \wedge x_2x_0 \wedge x_0^2$ in $\Lambda^{(4)}_{123456}$, will be send by the action of $g$ to the coefficient of $x_2^2 \wedge x_1^2 \wedge x_2x_0 \wedge x_0^2$ in $g\centerdot \Lambda^{(4)}_{123456}$.

Finally, for every polynomial $P$ contained in $\mathfrak{H}$, we have to compute the action of $g$, that is substituting each $\Delta_{ab}$ with $g\centerdot \Delta_{ab}$, and then we have to collect all the coefficients (polynomials in the Pl\"ucker coordinates of degree 2) of $g \centerdot P$, viewed as polynomial in the variables $g_{ij}$. For instance collecting the coefficients of the polynomial $g\centerdot(\Delta_{35}\Delta_{46}-\Delta_{25}\Delta_{56}),\ \Delta_{35}\Delta_{46}-\Delta_{25}\Delta_{56} \in \mathfrak{H}'$, we obtain 3495 polynomials that give some of the equations defining the Hilbert scheme $\Hilb{2}{2}$.

\end{example}

\chapter[On the connectedness of Hilbert schemes of lcm curves in $\PP^3$]{On the connectedness of\\ Hilbert schemes of locally Cohen-Macaulay curves in $\PP^3$}\label{ch:lcmHilb}

In this chapter I expose the first results obtained in collaboration with E. Schlesinger about the question of the connectedness of the Hilbert scheme of locally Cohen-Macaulay (lcm for short) curves in $\PP^3$. Basically this is an expanded version of the paper \lq\lq{The Hilbert schemes of locally Cohen-Macaulay curves in $\PP^3$ may after all be connected\rq\rq\ \cite{LellaSchlesinger}.

This chapter wants to be a further evidence of the potential of the ideas we are proposing to study many problems of geometric nature. Indeed looking at locally Cohen-Macaulay curves in $\PP^3$ with a computational and combinatorial eye, we prove quite easily the connectedness of the Hilbert scheme of (locally Cohen-Macaulay) curves of degree $4$ and genus $-3$, considered one of the best candidates to be a non-connected Hilbert scheme. 

Throughout this chapter, unless otherwise specified we will say Hilbert scheme meaning Hilbert scheme of locally Cohen-Macaulay curves.

\section{Introduction to the problem}

By the term \emph{curve} we will mean a one dimensional subscheme without isolated or embedded zero dimensional components; in the literature such an object is usually called a \emph{locally Cohen-Macaulay} curve.\index{curve!locally Cohen-Macaulay} One of the major tool in the study of Hilbert schemes of codimension two subschemes of a fixed projective scheme is the liaison theory. We refer to  \cite{HartshorneToulouse} and \cite{Migliore} for a modern treatment of this topic. Let us recall the main results of biliaison theory in the context of curves in $\PP^3$. We say that a curve $D$ is obtained from a curve $C$ by \emph{elementary biliaison}\index{biliaison} of \emph{height} $h$ if $D$ is linearly equivalent to $C+hH$ on a surface
$S$ that contains both $C$ and $D$ and has $H$ as its plane section \cite{HartshorneGD} (see also \cite[Chapter III]{MDP}).
For example, an effective divisor $D$ of bidegree $(a,b)$ on a smooth quadric surface $Q \simeq \PP^1 \times \PP^1 \subset \PP^3$ is obtained by an elementary biliaison of nonnegative height from a curve $C$ which is the disjoint union of $d=|a-b|$ lines on the quadric.
Biliaison is the equivalence relation generated by elementary biliaisons. Rao \cite{Rao} proved that there
is an invariant that distinguishes biliaison equivalence classes: the finite length graded module
$M_C = \bigoplus_{n \in \ZZ} H^1 (\PP^3, \mathcal{I}_C(n))$, which is now commonly referred to as
the Rao (or Hartshorne-Rao) module of $C$; two curves are in the same biliaison
class if and only if their Rao modules are isomorphic up to a twist. The structure of a biliaison
class is well understood. A curve $C$ in a biliaison class is said to be \emph{minimal} if for every other
curve $D$ of the class $M_D \simeq M_C (-h)$ with $h \geqslant 0$; this means the Rao module of $D$ is obtained shifting
the Rao module of $C$ to the right. The main result of the theory is:
\begin{enumerate}
\item[(i)]\label{property1}
Every\hfill biliaison\hfill class\hfill contains\hfill minimal\hfill curves;\hfill the\hfill family\hfill of\hfill the\hfill minimal\\ curves of the class is irreducible, and
any two minimal curves of the class can be joined by a finite number of elementary biliaisons of height zero.
\item[(ii)]\label{property2}
If $D$ is a non minimal curve, then $D$ is obtained from a minimal curve $C$ of the class by a finite
sequence of elementary biliaisons of positive height. This is known as the \emph{Lazarsfeld-Rao property}.
\end{enumerate}
Note that from \ref{property2} it follows that a curve is minimal in its biliaison class if and only if it has minimal degree among curves of the class. This result has a long history. Lazarsfeld and Rao \cite{LazarsfeldRao} proved the Lazarsfeld Rao property under a cohomological condition on $C$ that guarantees $C$ is minimal; they only considered \emph{trivial} elementary biliaisons ($C_i=C_{i-1}+nH$, no linear
equivalence allowed), but at the end of the process they needed a deformation with constant cohomology.
The existence of minimal curves and the Lazarsfeld-Rao property were proven independently in \cite{MDP}
and \cite{BallicoBolondiMigliore}. Strano \cite{Strano} showed the deformation at the end of the Lazarsfeld-Rao process is not needed if
one allows linear equivalence in the definition of elementary biliaison.  The version of the theorem we have given is due to Hartshorne \cite[Theorem 3.4]{HartshorneToulouse}, where the precise conditions on the ambient projective scheme are determined for the Lazarsfeld-Rao property to hold for biliaison classes of codimension $2$ subschemes.

We would like to stress the fact that for this theory to work it is necessary to consider locally Cohen-Macaulay curves:
even if one starts with a smooth irreducible curve, the minimal curve of the class may fail to be reduced or irreducible;
it may even not be generically a  local complete intersection (as it was assumed by Rao in \cite{Rao}). So
let \gls{lcmHilb} denote the Hilbert scheme parametrizing (locally Cohen-Macaulay) curves of degree $d$ and arithmetic genus $g$ in $\PP^3$; it is an open
subscheme of the full Hilbert scheme $\Hilb{3}{dt+1-g}$. In  \cite{MDPbounds} and \cite{MDPextremal} Martin-Deschamps and Perrin have shown that, whenever nonempty,
$\lcmHilb{d}{g}$ has a component $E_{d,g}$ whose closed points correspond to curves that have maximal cohomology, in the sense that
\begin{equation}\label{extr-in}
\dim H^{i} (\PP^3,\mathcal{I}_C(n)) \leqslant \dim H^{i} (\PP^3,\mathcal{I}_E(n))
\end{equation}
for every $C \in \lcmHilb{d}{g}$, $E \in E_{d,g}$, $n \in \ZZ$ and $i=0,1,2$. Note that by \eqref{extr-in} there is no obstruction coming from the semicontinuity of cohomology to specializing any curve in $\lcmHilb{d}{g}$ to an extremal curve. This remark raised the question whether every curve can be specialized to an extremal curve.\index{curve!extremal}\index{extremal curve|see{curve, extremal}} This is known to be false for curves that are not generically of embedding dimension two \cite{NolletSchlesingerDegree4}, but it is open for, say, smooth curves.
A weaker version of this question, proposed by Hartshorne in \cite{HartshorneConnectedness} and \cite{HartshorneConnectedness2},  is whether $\lcmHilb{d}{g}$ is connected, that is, if every curve in $\lcmHilb{d}{g}$ belongs to the connected component containing the extremal curves.
It is an interesting problem because the Hilbert scheme of smooth curves is not connected \cite{Hartshorne}, while the full Hilbert scheme is connected (as seen in Chapter \ref{ch:deformations}), but through schemes with zero dimensional components.

Let us review what is known about the problem of connectedness of $\lcmHilb{d}{g}$. No example of a nonconnected $\lcmHilb{d}{g}$ has been found so far. The Hilbert scheme $\lcmHilb{d}{g}$ is connected when $g \geqslant \binom{d-3}{2}-1$  (see \cite{AitAmrane,HMDP,Sabadini}) and when $d \leqslant 4$ (see \cite{NolletDegree3,NolletSchlesingerDegree4}; this is non trivial, because when $g$ is very negative the Hilbert scheme has a large number of irreducible components).
Hartshorne \cite{HartshorneConnectedness} has shown that smooth irreducible nonspecial curves and arithmetically Cohen-Macaulay (ACM) curves are in the connected component containing the extremal curves,
and so are Koszul curves by \cite{PerrinKoszul}. Building on the results of \cite{HartshorneConnectedness}, E. Schlesinger has shown
\cite{SchlesingerFootnote} that, if a minimal curve $C$ can be connected to an extremal curve by flat families
lying on surfaces of degree $s$, where $s$ is the least degree of a surface containing $C$, then \emph{every} curve in the biliaison class of $C$ is in the connected component of  extremal curves in its Hilbert scheme.

By \cite{HartshorneGD} and \cite{HartshorneSchlesinger2H}, any curve contained in a singular quadric surface, including a double plane, is in the connected component of  extremal curves \cite{HartshorneSchlesinger2H}. On the other hand, the case of curves on a smooth quadric surface
has been open so far. By biliaison it is enough to deal with divisors of bidegree $(d,0)$. The case $d \leqslant 2$ is then trivial,
and Nollet \cite{NolletDegree3} has shown that is possible  to specialize a divisor of bidegree $(3,0)$ to an extremal curve, but
it has been an open question whether one could specialize four (or more) disjoint lines on a smooth quadric
to an extremal curve of the same genus; consequently, the Hilbert scheme of curves of degree $\lcmHilb{10}{12}$, which has an irreducible component whose general member is a divisor of bidegree $(7,3)$ on a smooth quadric surface, was proposed as a  candidate for an example of a nonconnected $\lcmHilb{d}{g}$: see \cite[Ex. 4.2]{HartshorneConnectedness}, \cite[Section 4]{HartshorneConnectedness2}, and the open problems list of the 2010 AIM workshop on \emph{Components of  Hilbert Schemes}  available at \href{http://aimpl.org/hilbertschemes}{\texttt{aimpl.org/hilbertschemes}}.

\section{Extremal curves} \label{two}
In this section we establish notation and terminology and review some known results that we will need later.
We work over an algebraically closed field $\K$ of arbitrary characteristic. We denote with the symbol
$\mathcal{I}_X$ the ideal sheaf of a subscheme $X \subset \PP^3$. Given a coherent
sheaf $\mathcal{F}$ on $\PP^3$, we define $h^i(\mathcal{F})= \dim  H^{i}(\PP^3,\mathcal{F})$ and
$H^i_{\ast} (\mathcal{F})= \bigoplus_{n \in \ZZ} H^{i}(\PP^3,\mathcal{F}(n))$. We write
$\K[x]=\K[x_0,x_1,x_2,x_3]$ for the homogeneous coordinate ring $H^{0}_{\ast}(\OO_{\PP^3})$ of $\PP^3$.

\begin{definition}
A \emph{curve} in $\PP^3$ (or more precisely a locally Cohen-Macaulay curve)\index{curve!locally Cohen-Macaulay}
is a one dimensional subscheme $C \subset \PP^3$ without zero dimensional associated points; this means
that all irreducible components of $C$ have dimension $1$, and that $C$ has no embedded points.
\end{definition}
We denote by $\lcmHilb{d}{g}$ the Hilbert scheme parametrizing curves of degree $d$ and arithmetic
genus $g$ in $\PP^3$ \cite{MDP}. This is an open subscheme of the Hilbert scheme $\Hilb{3}{dt+1-g}$ parametrizing
all one dimensional subschemes of $\PP^3$ with Hilbert polynomial $dt+1-g$.

\begin{definition}
 We say that a curve $E \subset \PP^3$ of degree $d$ and genus $g$ is \emph{extremal} if
\begin{equation}\label{extr_def}\index{curve!extremal}
h^{i} (\mathcal{I}_C(n)) \leqslant h^{i} (\PP^3,\mathcal{I}_E(n))
\end{equation}
for $i=0,1,2$, for every curve $C$ of degree $d$ and genus $g$, and for every $ n \in \ZZ$,
\end{definition}
Thus a curve is extremal if it has the largest possible cohomology.
One knows that $\lcmHilb{d}{g}$ is nonempty if  and only if either $g=\binom{d-1}{2}$, in which case it is
irreducible and consists of plane curves, or $g \leqslant \binom{d-2}{2}$. Whenever nonempty,
$\lcmHilb{d}{g}$ contains extremal curves \cite{MDPbounds}; in fact, the extremal curves form an irreducible
$E_{d,g}$ of $\lcmHilb{d}{g}$ \cite{MDPextremal}.

\begin{remark}
Our definition of extremal curves is equivalent to the one given by Hartshorne \cite{HartshorneConnectedness,HartshorneConnectedness2}.
Martin-Deschamps and Perrin did not include ACM curves with maximal cohomology among extremal curves;
the difference comes up only when $g=\binom{d-1}{2}$ and $g=\binom{d-2}{2}$.
In \cite{MDPbounds} the functions $h^{i} (\mathcal{I}_E (n))$ are computed explicitly for an extremal curve $E$, and the
bounds \eqref{extr_def} are proven, under  the assumption that the field $\K$ has characteristic zero;
this assumption on the characteristic is not necessary \cite{NolletSubextremal}.
\end{remark}

Martin-Deschamps and Perrin also
compute the Rao module $M_E = H^{1}_{\ast} (\mathcal{I}_E)$ of an extremal curve:

\begin{theorem}[{\cite{MDPbounds,MDPextremal}}]
Let $(d,g)$ be two integers satisfying $d \geqslant 2$ and
$g \leqslant \binom{d-2}{2}-1$.
A curve $E$ of degree $d$ and genus $g$ is extremal if and only if
\[
M_E \simeq R/(l_1,l_2,F,G) (b)
\]
where $(l_1,l_2,F,G)$ is a regular sequence, $\deg l_1 =\deg l_2 =1$,
$\deg F= \binom{d-2}{2}-g$, $\deg G= \binom{d-1}{2}-g$, and
$b=\deg F-1$.
\end{theorem}

The following proposition describes extremal curves supported on a line. It is a special case of
\cite[Proposition 0.6]{MDPextremal};  we state it in the form needed later in the chapter.

\begin{proposition} \label{extr_mls}
Let $(d,g)$ be a pair of integers satisfying $g \leqslant \binom{d-2}{2}-1$. Let $F$ and $G$
be two forms of degrees $\deg F= \binom{d-2}{2}-g$ and $\deg G= \binom{d-1}{2}-g$ in $\K[x_1,x_0]$ with no common zeros.
The surface $S$ of equation $x_3 G-x_2^{d-1}F=0$ is irreducible and generically smooth along the line
$L$ of equations $x_3=x_2=0$. It therefore contains a unique curve $E$ of degree $d$ supported on $L$.
The curve $E$ is extremal of degree $d$ and genus $g$, and its Rao module is
\[
M_E \simeq \K[x_0,x_1,x_2,x_3]/(x_3,x_2,F,G) (b) \simeq \K[x_0,x_1]/(F,G) (b)
\]
where $b=\deg(F)+1$. The homogeneous ideal of $E$ is generated by $x_3^2,x_3x_2,x_2^d$ and $x_3G-x_2^{d-1}F$.
\end{proposition}

\begin{proof}
The surface $S$ is irreducible because $F$ and $G$ have no common zeros, and it is smooth at points
of $L$ where $G$ is different from zero. Therefore  the ideal of $L$ in the local ring
$\OO_{S,\xi}$ of the generic
point $\xi$ of $L$ is generated by one function $t$, and the ideal of a curve of degree $d$ supported on $L$
must be $t^d \OO_{S,\xi}$. Since a locally Cohen-Macaulay curve supported on $L$ is determined by its ideal
at the generic point of $L$, we see that there is a unique curve $D_m \subset S$ supported on $L$ of degree $m$
for every $m \geqslant 1$. For $m=d-1$, the curve $P=D_{d-1}$ is the planar multiple structure of equations
$x_3=x_2^{d-1}=0$.  We note that
$\mathcal{I}_P \otimes \OO_L \simeq \OO_L (-1) \oplus \OO_L (1-d)$  where the two generators are the images
of $x_3$ and $x_2^{d-1}$. The two forms $F$ and $G$ define a surjective map $\OO_L (-1) \oplus \OO_L (1-d) \rightarrow \OO_L (b)$;
composing this with the natural map  $\mathcal{I}_P \rightarrow \mathcal{I}_P \otimes \OO_L$  we obtain
a surjection $\phi: \mathcal{I}_P \rightarrow \OO_L (b)$. We let $E$ be the subscheme of $\PP^3$ whose ideal sheaf is
the kernel of $\phi$. By construction we have an exact sequence
 \begin{equation}\label{extr_seq}
    0\ \longrightarrow\ \mathcal{I}_E \ \longrightarrow\ \mathcal{I}_P \ \longrightarrow\ \OO_L (b)\ \longrightarrow\ 0.
 \end{equation}
This sequence shows that $E$ is a (locally Cohen-Macaulay) curve of degree $d$ and genus $g$, and that its homogeneous
ideal is generated by $x_3^2$,$x_3x_2$,$x_2^d$ and $x_3G-x_2^{d-1}F$. Therefore $E=D_d$ is the unique curve of degree $d$ contained in $S$ and supported on $L$. Finally, the long exact cohomology sequence associated to \eqref{extr_seq} shows that the Rao module of $E$ is
\[
M_E= \K[x_0,x_1]/(F,G) (b) = \K[x_0,x_1,x_2,x_3]/(x_3,x_2,F,G) (b).
\]
Hence $E$ is an extremal curve.
\end{proof}

Extremal curves seem to have a special role on the study of the connectedness of $\lcmHilb{d}{g}$ as the lexicographic ideal for the connectedness of the full Hilbert scheme $\Hilb{n}{p(t)}$. The lexicographic ideal is so important for two reasons:
\begin{enumerate}
\item every Hilbert scheme $\Hilb{n}{p(t)}$ contains a point defined by the lexicographic ideal associated to $p(t)$;
\item being the lexicographic ideal a segment ideal w.r.t. $\DegLex$, this term ordering fixes the \lq\lq direction\rq\rq\ to follow to define flat deformations and specializations of ideals in order to approach the lexicographic point (polarizations for Hartshorne \cite{HartshorneThesis}, Gr\"obner degenerations for Peeva and Stillman \cite{PeevaStillman} and Borel degenerations in Chapter \ref{ch:deformations}).
\end{enumerate}

Each Hilbert scheme of locally Cohen-Macaulay curves $\lcmHilb{d}{g}$ contains extremal curves \cite{MDP}, therefore they are a good candidate to play the analogous special role as the lexicographic ideal. The point we want to discuss is how to detect the \lq\lq direction\rq\rq\ that allows to approaching them, basically applying Gr\"obner degenerations. This direction can not be given by a term ordering, mainly for two reasons:
\begin{itemize}
\item there is not a unique extremal curve, but we have to consider a wide class of curves, precisely a component of $\lcmHilb{d}{g}$;
\item the ideal defining an extremal curve is not a monomial ideal.
\end{itemize}
Thus the idea is to consider a weight order $\omega$ such that
\begin{enumerate}
\item it does not distinguish $x_1$ and $x_0$;
\item the generator $x_3 G - x_2^{d-1} F$ is $\omega$-homogeneous.
\end{enumerate}
The first consequence is that the order defined by $\omega$ on the monomials of $\K[x]$ is not a total order, but the Gr\"obner machinery still works. Indeed given a polynomial $P(x_0,x_1,x_2,x_3)$, the \emph{initial form} $\IN_{\omega}(P)$ of $P$ with respect to $\omega$ is the sum of all the terms $c_{\alpha} x^\alpha$ in $P$ for which the scalar product
\[
\omega \cdot \alpha= \omega_3 \alpha_3 + \omega_2 \alpha_2 + \omega_1 \alpha_1 + \omega_0 \alpha_0
\]
is maximal. The initial ideal $\IN_\omega (I)$ of an ideal $I$ is the ideal generated by the initial forms  $\IN_{\omega} (P)$ as $P$ varies in $I$. Also in this case, there is a flat family over the affine line $\AA^1 = \Spec \K[t]$ whose fibers over $t \neq 0$ are isomorphic to $\Proj \K[x]/I$, while the special fiber over zero is the subscheme of $\PP^3$ defined by the $\IN_\omega (I)$: see for example \cite{BayerMumford} and \cite[Theorem 15.17]{Eisenbud}. Roughly, this family is defined letting the one dimensional torus act on $\PP^3$ by $t \centerdot [x_3:x_2:x_1:x_0]=[t^{\omega_3}x:t^{\omega_3}x_2:t^{\omega_1}x_1:t^{\omega_0}x_0]$ and taking the limit as $t$ goes to zero, so that set theoretically we are projecting $C$ onto the line $L$ of equations $x_3=x_2=0$, but what is interesting is the scheme theoretic structure of the limit.

\medskip

Let us suppose 
\begin{equation}\label{eq:semiBorel}
x_3 >_{\omega} x_2 >_\omega x_1 =_{\omega} x_0
\end{equation}
and $\deg_\omega x_1 = \deg_{\omega} x_0 = 1$, i.e. $\omega = (\omega_3,\omega_2,1,1)$. The order induced by the hypothesis \eqref{eq:semiBorel} on the monomials of a fixed degree can be represented as a planar graph in similar way to that one used to represent the Borel partial order $\leq_B$. 

\begin{definition}
We call $\Omega(m)$ the graph defined as follows:
\begin{itemize}
\item the\hfill vertices\hfill correspond\hfill to\hfill the\hfill classes\hfill of\hfill polynomials\hfill $x_3^{\alpha_3} x_2^{\alpha_2} F$,\hfill with\hfill $F$\hfill in\\ $\K[x_0,x_1]$ such that $\alpha_3 + \alpha_2 + \deg F = m$;
\item the edges correspond to the two classes of transformations $\lcmdown{3} = \frac{x_2}{x_3}$ and $\lcmdown{2}=\frac{l}{x_2}$ with $l \in \K[x_0,x_1]_1$.
\end{itemize}
\end{definition}

\begin{proposition}
The\hfill graph\hfill $\Omega(m)$\hfill is\hfill isomorphic\hfill to\hfill the\hfill graph\hfill representing\hfill the\hfill poset\\ $\pos{2}{m}$.
\end{proposition}
\begin{proof}
Assumed that the monomials of $\pos{2}{m}$ belong to the ring $\K[x_1,x_2,x_3]$, it suffices to consider the map $\psi: \Omega(m) \rightarrow \pos{2}{m}$
\[
\begin{array}{cccl}
x_3^{\alpha_3} x_2^{\alpha_2} F & \longmapsto & x_3^{\alpha_3} x_2^{\alpha_2} x_1^{\deg F}, & F\in \K[x_1,x_0]\\
\lcmdown{i} & \longmapsto & \down{i}, & i=2,3. 
\end{array}\qedhere
\]
\end{proof}

\begin{definition}
Let $I \subset \K[x]$ be an ideal defining a locally Cohen-Macaulay curve $C$ in $\lcmHilb{d}{g},\ d \geqslant 3$ and let $r$ be the maximal degree of a generator of $I$. We will say that $I$ is a \emph{lcm-segment} ideal\index{segment ideal!lcm-segment ideal}\index{lcm-segment ideal|see{segment ideal}} if there exists a weight order $\omega = (\omega_3,\omega_2,1,1)$ such that
\begin{enumerate}[(i)]
\item every polynomial in the basis of $\langle I_r\rangle$ as $\K$-vector space is $\omega$-homogeneous;
\item the image through $\psi$ of the set of all monomials involved in a basis of $\langle I_r\rangle$ is a segment of $\pos{2}{r}$ w.r.t. $\overline{\omega} = (\omega_3,\omega_2,1)$.
\end{enumerate}
\end{definition}

\begin{proposition}\label{prop:weightExtremal}
The ideal $I_E \subset \K[x_0,x_1,x_2,x_3]$ defining any extremal curve $E \in \lcmHilb{d}{g}$ supported on a line is a lcm-segment ideal w.r.t. $\omega = (d,2,1,1)$.
\end{proposition}
\begin{proof}
By Proposition \ref{extr_mls}, we know that the ideal of an extremal curve in of degree $d$ and genus $g$ supported on a line is
\[
I_E = (x_3^2,x_3x_2,x_2^d,x_3 G - x_2^{d-1} F) \quad \text{with} \quad \deg (x_3 G - x_2^{d-1} F) = \binom{d-1}{2} - g +1.
\]
Set $r = \binom{d-1}{2} - g +1$, the vector space $\left\langle \left(I_E\right)_r\right\rangle$ as only one non-monomial generator: $x_3 G - x_2^{d-1} F$. So we need a weight order $\omega$ such that
\[
\deg_\omega x_3 G = \omega_3 + \binom{d-1}{2} - g = (d-1)\omega_2 + \binom{d-2}{2} - g = \deg_\omega x_2^{d-1} F
\]
i.e. $\omega_3 = (d-1)\omega_2 - \binom{d-2}{1} = (d-1)\omega_2 - d + 2$.

Then the image through $\psi$ of monomials involved in the basis of $\left\langle \left(I_E\right)_r\right\rangle$ corresponds to the Borel set\index{Borel set} $\mathscr{B} = \{(x_3^2,x_3 x_2,x_2^d)_r\} \cup \{x_3x_1^{r-1},x_2^{d-1} x_1^{r-d+1}\} \subset \pos{2}{r}$ (see Figure \ref{fig:extremalCurve}), so that 
\[
\mathscr{N} = \pos{2}{r} \setminus \mathscr{B} = \left\{x_2^{d-2} x_1^{r-d+2},\ldots,x_1^r\right\}.
\]
The two monomials of $\mathscr{B}$ images of the generator $x_3G - x_2^{d-1}F$ are minimal elements and $\mathscr{N}$ has as unique maximal monomial $x_2^{d-2} x_1^{r-d+2}$, hence in order for $\mathscr{B}$ to be a segment w.r.t. $\omega = (\omega_3,\omega_2,\omega_1)$ we need
\[
\left\{\begin{array}{l}
\omega_3 > \omega_2\\
\omega_2 > \omega_1\\
\omega_3 + (r-1) \omega_1 > (d-2)\omega_2 + (r-d+2) \omega_1
\end{array}\right.
\] 
Finally adding the hypothesis $\omega_1 = 1$ and the requirement of homogeneity for $x_3G - x_2^{d-1} F$, we have to solve the system
\[
\left\{\begin{array}{l}
\omega_3 > \omega_2\\
\omega_2 > 1\\
\omega_3 + (r-1) > (d-2)\omega_2 + (r-d+2)\\
\omega_3 = (d-1)\omega_2 - d + 2
\end{array}\right.
\]
Replacing $\omega_3$ in the third inequality using the equality we obtain $\omega_2 > 1$, so that a solution of the system is
\[
\left\{\begin{array}{l}
\omega_3 > \omega_2\\
\omega_2 > 1\\
\omega_3 = (d-1)\omega_2 - d + 2
\end{array}\right. \quad\Rightarrow\quad
\left\{\begin{array}{l}
\omega_2 = 2\\
\omega_3 = 2(d-1) - d + 2 = d > 2
\end{array}\right.
\qedhere
\]
\end{proof}
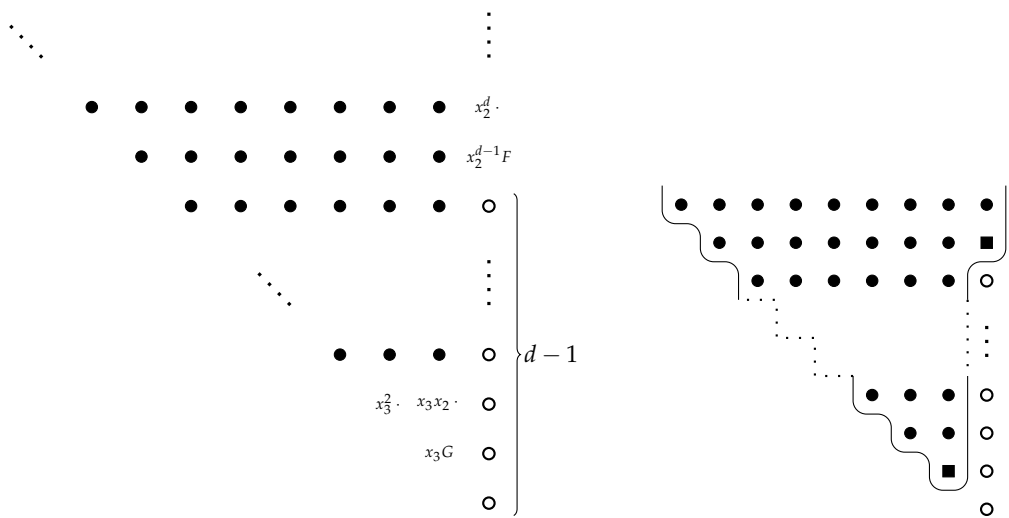
\begin{figure}[!ht]
\begin{center}
\subfloat[][The ideal of an extremal curve drawn in the graph $\Omega(r)$.]{\label{fig:extremalCurve_a}
\begin{tikzpicture}[scale=0.65,decoration=brace]
\tikzstyle{ideal}=[circle,draw=black,fill=black,inner sep=1.5pt]
\tikzstyle{quotient}=[circle,draw=black,thick,inner sep=1.5pt]

\node at (-4,4) [ideal] {};
\node at (-3,4) [ideal] {};
\node at (-2,4) [ideal] {};
\node (n)  at (-1,4) [ideal] {};
\node (n0) at (0,4) [ideal] {};
\node (n1) at (1,4) [ideal] {};
\node (n2) at (2,4) [ideal] {};
\node (n3) at (3,4) [ideal] {};
\node (n4) at (4,4) [] {\tiny $x_2^d \cdot{}$};

\node at (-3,3) [ideal] {};
\node at (-2,3) [ideal] {};
\node (n)  at (-1,3) [ideal] {};
\node (n0) at (0,3) [ideal] {};
\node (n1) at (1,3) [ideal] {};
\node (n2) at (2,3) [ideal] {};
\node (n3) at (3,3) [ideal] {};
\node (n4) at (4,3) [] {\tiny $x_2^{d-1}F$};

\node at (-2,2) [ideal] {};
\node (n)  at (-1,2) [ideal] {};
\node (n0) at (0,2) [ideal] {};
\node (n1) at (1,2) [ideal] {};
\node (n2) at (2,2) [ideal] {};
\node (n3) at (3,2) [ideal] {};
\node (n4) at (4,2) [quotient] {};

\node (11) at (1,-1) [ideal] {};
\node (12) at (2,-1) [ideal] {};
\node (13) at (3,-1) [ideal] {};
\node (14) at (4,-1) [quotient] {};

\node (22) at (2,-2) [] {\tiny $x_3^2 \cdot{}$};
\node (23) at (3,-2) [] {\tiny $x_3 x_2 \cdot{}$};
\node (24) at (4,-2) [quotient] {};

\node (33) at (3,-3) [] {\tiny $x_3G$};
\node (34) at (4,-3) [quotient] {};

\node (44) at (4,-4) [quotient] {};

\draw [loosely dotted,very thick] (4,0) -- (4,1); 
\draw [loosely dotted,very thick] (0,0) -- (-0.7,0.7); 

\draw [loosely dotted,very thick] (4,5) -- (4,6); 
\draw [loosely dotted,very thick] (-5,5) -- (-5.7,5.7);


\draw [decorate] (4.5,2.25) --node[right]{\footnotesize $d-1$} (4.5,-4.25);
\end{tikzpicture}
}
\qquad
\subfloat[][The Borel set in $\pos{2}{r}$ associated to an extremal curve.]{\label{fig:extremalCurve_b}
\begin{tikzpicture}[scale=0.5,decoration=brace]
\tikzstyle{ideal}=[circle,draw=black,fill=black,inner sep=1.5pt]
\tikzstyle{quotient}=[circle,draw=black,thick,inner sep=1.5pt]

\node at (-4,4) [ideal] {};
\node at (-3,4) [ideal] {};
\node at (-2,4) [ideal] {};
\node (n)  at (-1,4) [ideal] {};
\node (n0) at (0,4) [ideal] {};
\node (n1) at (1,4) [ideal] {};
\node (n2) at (2,4) [ideal] {};
\node (n3) at (3,4) [ideal] {};
\node (n4) at (4,4) [ideal] {};

\node at (-3,3) [ideal] {};
\node at (-2,3) [ideal] {};
\node (n)  at (-1,3) [ideal] {};
\node (n0) at (0,3) [ideal] {};
\node (n1) at (1,3) [ideal] {};
\node (n2) at (2,3) [ideal] {};
\node (n3) at (3,3) [ideal] {};
\node (n4) at (4,3) [regular polygon,regular polygon sides=4,draw=black,fill=black,inner sep=1.5pt] {};

\node at (-2,2) [ideal] {};
\node (n)  at (-1,2) [ideal] {};
\node (n0) at (0,2) [ideal] {};
\node (n1) at (1,2) [ideal] {};
\node (n2) at (2,2) [ideal] {};
\node (n3) at (3,2) [ideal] {};
\node (n4) at (4,2) [quotient] {};

\node (11) at (1,-1) [ideal] {};
\node (12) at (2,-1) [ideal] {};
\node (13) at (3,-1) [ideal] {};
\node (14) at (4,-1) [quotient] {};

\node (22) at (2,-2) [ideal] {};
\node (23) at (3,-2) [ideal] {};
\node (24) at (4,-2) [quotient] {};

\node (33) at (3,-3) [regular polygon,regular polygon sides=4,draw=black,fill=black,inner sep=1.5pt] {};
\node (34) at (4,-3) [quotient] {};

\node (44) at (4,-4) [quotient] {};

\draw [loosely dotted,very thick] (4,0) -- (4,1); 

\draw [loosely dotted,thick] (3.5,1.5) -- (3.5,-0.5);
\draw [loosely dotted,thick] (0.5,-0.5) -- (-0.5,-0.5) -- (-0.5,0.5) -- (-1.5,0.5) -- (-1.5,1.5) -- (-2.5,1.5);
\draw [rounded corners=5pt] (3.5,-0.5) -- (3.5,-3.5) -- (2.5,-3.5) -- (2.5,-2.5) -- (1.5,-2.5) -- (1.5,-1.5) -- (0.5,-1.5) -- (0.5,-0.5);
\draw [rounded corners=5pt] (3.5,1.5) -- (3.5,2.5) -- (4.5,2.5) -- (4.5,4.5);
\draw [rounded corners=5pt] (-2.5,1.5) -- (-2.5,2.5) -- (-3.5,2.5) -- (-3.5,3.5) -- (-4.5,3.5) -- (-4.5,4.5);
\end{tikzpicture}
}
\caption{\label{fig:extremalCurve} The ideal defining an extremal curve represented in $\Omega(r)$ and the associated Borel set in $\pos{2}{r}$.}
\end{center}
\end{figure}

\begin{remark}
The weight order making the ideal of an extremal curve $E \in \lcmHilb{d}{g}$ a lcm-segment ideal does not depend on the genus of the curve, but only on the degree.  
\end{remark}

\section[The Hilbert scheme of curves of degree 4 and genus -3 is connected]{\gls{lcmHilb4m3} is connected} \label{three}

In this section we construct a flat family of curves whose general member is a disjoint union
of lines on a smooth quadric surface $Q$ and whose special member is an extremal multiple line.
This specialization is obtained considering the initial ideal with respect to the weight vector
$\omega=(d,2,1,1)$. In case $d=3$, such a specialization was constructed by Nollet \cite{NolletDegree3}
for a triple structure $3L$ on the line $L$ on $Q$ (without using the language of weight vectors and initial ideals).
For $d \geqslant 4$, if we begin with the $d$-uple structure $dL$ on the quadric, we obtain a limit
with embedded points. We now show that, if we take a sufficiently general divisor of bidegree $(d,0)$,
for any $d \geqslant 3$, then we obtain an extremal curve as a limit.

\begin{theorem} \label{main}
Let $Q$ be the quadric surface of equation
\begin{equation}
q(x_0,x_1,x_2,x_3)=x_3(x_3+x_0)-x_2x_1=0.
\end{equation}
For every $a \in \K$ let $L_a \subset Q$ denote the line of equations $x_3-a x_1=x_2-a(ax_1+x_0)=0$.
Given $d \geqslant 3$ and $a_1, \ldots, a_d \in \K$, consider the divisor
\begin{equation}
C=L_{a_1}+ \cdots + L_{a_d}
\end{equation}
on $Q$. If the sums $a_{i} + a_{j}$  for $1 \leqslant i<j \leqslant d$ are all distinct, then there is a flat family of pairs $\mathcal{C} \subset \mathcal{Q} \rightarrow \AA^1$, whose fiber over $1$ is $(C,Q)$, whose fiber over $t \neq 0$ consists of $d$ disjoint lines on a smooth quadric surface, and whose fiber over $0$ is an extremal curve in the double plane of equation $x_3^2=0$.
\end{theorem}
\begin{proof}
Let $C_0$ denote the subscheme defined by (the saturation of) $\IN_\omega(I_C)$,  where $I_C$ denotes the homogeneous ideal of $C$. By Gr\"obner bases theory, there is a flat specialization from $C$ to $C_0$; since $\IN_{\omega}(q)=x_3^2$, the smooth quadric $Q$ specializes to the double plane $x_3^2=0$ as $C$ specializes to $C_0$.
We let $l_i=x_3-a_i x_1$ and $m_i=x_2-a_i(a_ix_1+x_0)$ denote the given equations for the line $L_{a_i}$. Since
$I_C$ contains the product of the ideals of the lines $L_{a_i}$,
\[
x_2^d= \IN_{\omega} (m_1 \cdot \ldots \cdot m_d) \in I_{C_0}.
\]
Therefore $C_0$ is contained in the complete intersection $x_3^2=x_2^d=0$ and so it is supported on the line $L$.
We will show that $I_C$ contains a polynomial $A(x_0,x_1,x_2,x_3)$ such that $\IN_\omega(A)=x_3G-x_2^{d-1}F$
where $F$ and $G$ are homogeneous forms in $\K[x_0,x_1]$, $\deg G= \binom{d}{2}$, $\deg F= \binom{d-1}{2}+1$,
and $F$ and $G$ have no common zero in $\PP^1=\Proj \K[x_0,x_1]$. It follows that $C_0$ is contained
in the surface $S$ of equation $x_3 G-x_2^{d-1}F$. By flatness, the Hilbert polynomial of $C_0$ coincides
with that of $C$, so $C_0$ is a one dimensional subscheme of $\PP^3$ of degree $d$ and genus $1-d$.
Let $E$ the largest Cohen-Macaulay curve contained in $C_0$: it is the curve of degree $d$ obtained
from $C_0$ throwing away its embedded points. By Proposition \ref{extr_mls} $E$ is the unique
curve of degree $d$ contained in $S$ and supported on the line $L$; it is an extremal curve of degree
$d$ and genus $1-d$. Since $E \subset C_0$ and the two schemes have the same Hilbert polynomial, we conclude
$E=C_0$. Thus the limit is an extremal curve, and the statement is proven.

To conclude the proof, we need to find $A \in I_C$ with $\IN_\omega(A)=x_3 G-x_2^{d-1}F$. In principle,
this is a Gr\"obner basis calculation, but luckily we can bypass such a calculation because we were
able to find, with the help of some computation performed with \textit{Macaulay2} \cite{M2}, a
determinantal formula for $A$ (note that, while $I_C$ is generated by forms of degree $\leqslant d$,
the degree of $A$ is much larger than $d$). Let $A=x_3 G-x_1 B$ denote the determinant
\begin{equation*}
\begin{vmatrix}
l_1         & m_1     & m_1^2  & \ldots & m_1^{d-1} \\
l_2         & m_2     & m_2^2  & \ldots & m_2^{d-1} \\
\vdots    & \vdots  & \vdots & \ddots & \vdots    \\
l_d         & m_d     & m_d^2  & \ldots & m_d^{d-1} \\
\end{vmatrix}.
\end{equation*}
Since $l_i = x_3 - a_i \, x_1$, by linearity
\begin{equation}\label{ad}
A=x_3
\begin{vmatrix}
1         & m_1     & m_1^2  & \ldots & m_1^{d-1} \\
1         & m_2     & m_2^2  & \ldots & m_2^{d-1} \\
\vdots    & \vdots  & \vdots & \ddots & \vdots    \\
1         & m_d     & m_d^2  & \ldots & m_d^{d-1} \\
\end{vmatrix}
-x_1
\begin{vmatrix}
a_1         & m_1     & m_1^2  & \ldots & m_1^{d-1} \\
a_2         & m_2     & m_2^2  & \ldots & m_2^{d-1} \\
\vdots    & \vdots  & \vdots & \ddots & \vdots    \\
a_d         & m_d     & m_d^2  & \ldots & m_d^{d-1} \\
\end{vmatrix}.
\end{equation}
The linear forms $l_i$ and $m_i$ are the equations of the line $L_{a_i}$, hence the polynomial $A$
belongs to the ideal of $C$. In the expansion $A=xG-zB$
the polynomial $G$ is a Vandermonde determinant, i.e.
\begin{equation} \label{eqg}
G= \prod_{1 \leqslant i < j \leqslant d} (m_j-m_i) =  \prod_{1 \leqslant i < j \leqslant d} (a_i-a_j) \big((a_i+a_j)x_1+x_0\big)
\end{equation}
We note that $G$ is nonzero: the hypothesis that the sums $a_i+a_j$ be all distinct implies
that $a_i-a_j \neq 0$ for every $i < j$. Furthermore, the zeros of $G$ in $\PP^1=\Proj \K[x_0,x_1]$ are the points $[1:-a_i-a_j]$.

The polynomial $B$ is
\begin{equation} \label{eqh}
\begin{split}
B&{}= \begin{vmatrix}
a_1         & m_1     & m_1^2  & \ldots & m_1^{d-1} \\
a_2         & m_2     & m_2^2  & \ldots & m_2^{d-1} \\
\vdots    & \vdots  & \vdots & \ddots & \vdots    \\
a_d         & m_d     & m_d^2  & \ldots & m_d^{d-1} \\
\end{vmatrix}=\\
&{} = \sum_{j=1}^d (-1)^{j-1} a_j\, m_1  \cdots \widehat{m_j} \cdots m_d \prod_{\begin{subarray}{c}1 \leqslant h < k \leqslant d\\ h\neq j, k \neq j \end{subarray}} (m_k-m_h).
\end{split}
\end{equation}
Since $m_k - m_h=  (a_h-a_k) \big((a_h+a_k)x_1+x_0\big)$ and the initial term of $m_j$ is $x_2$, the initial term of $B$ is the polynomial
\begin{equation} \label{eqlt}
\sum_{j=1}^d (-1)^{j-1} a_j\, x_2^{d-1} \prod_{\begin{subarray}{c}1 \leqslant h < k \leqslant d\\ h\neq j, k \neq j \end{subarray}}(a_h-a_k) \big((a_h+a_k)x_1+x_0\big)=x_2^{d-1}P
\end{equation}
provided $P =P(x_0,x_1) \neq 0$. We will prove not only that $P$ is not zero, but also that it has no zero in common with $G$, that is, it does not vanish at the points $[1:-a_i-a_j]$. By symmetry, it is enough to show that
$P(1,- a_1 - a_2) \neq 0$.
For this, we write $P$ as a determinant: if we set $p_i=a_i^2 x_1+a_i x_0$, then
\begin{equation} \label{eqf2}
P=
\sum_{j=1}^d (-1)^{j-1} a_j \prod_{\begin{subarray}{c}1 \leqslant h < k \leqslant d\\ h\neq j, k \neq j \end{subarray}}(p_k-p_h) =
\begin{vmatrix}
a_1     &1       & p_1        & \ldots   & p_1^{d-2} \\
a_2     &1       & p_2        & \ldots   & p_2^{d-2} \\
\vdots  &\vdots  & \vdots     & \ddots   &\vdots     \\
a_d     &1       & p_d        & \ldots   & p_d^{d-2} \\
\end{vmatrix}
\end{equation}
Now note that $p_1(1,-a_1-a_2)=p_2 (1,-a_1-a_2)=-a_1 a_2$ and $p_i(1,-a_1-a_2)=a_i^2-(a_1+a_2)a_i$ for $j \geqslant 3$. For simplicity we write $p_i$ in place of $p_i(1,-a_1-a_2)$.
Then
\begin{equation} \label{eqf3}
\begin{split}
P(1,&-a_1-a_2)=
\begin{vmatrix}
a_1     &1       & (-a_1a_2)        & \ldots   & (-a_1a_2)^{d-2} \\
a_2     &1       & (-a_1a_2)       & \ldots   & (-a_1a_2)^{d-2} \\
\vdots  &\vdots  & \vdots     & \ddots   &\vdots     \\
a_d     &1       & p_d        & \ldots   & p_d^{d-2} \\
\end{vmatrix}
=\\
&{}=
(a_1-a_2)
\begin{vmatrix}
1       & (-a_1a_2)        & \ldots   & (-a_1a_2)^{d-2} \\
1       & p_3        & \ldots   & p_{3}^{d-2} \\
\vdots  & \vdots     & \ddots   &\vdots     \\
1       & p_d        & \ldots   & p_d^{d-2} \\
\end{vmatrix}=
\\
&{} =
(a_1-a_2) \left(\prod_{j=3}^d (p_j+a_1a_2)\right)\left(  \prod_{3 \leqslant h < k \leqslant d} (p_k-p_h) \right)= \\
&{}=
(a_1-a_2) \left(\prod_{j=3}^d (a_j-a_1)(a_j-a_2) \right)\left(  \prod_{3 \leqslant h < k \leqslant d} (a_k\!-\!a_h)(a_h+a_k\!-\!a_1\!-\!a_2) \right).
\end{split}
\end{equation}
This shows $P(1,-a_1-a_2) \neq 0$ because of the assumption that sums $a_i+a_j$ be all distinct.

To finish, we let $F=x_1P$. Then $\IN_{\omega}(A)=x_3G-x_2^{d-1}F$, and $F$ and $G$ are homogeneous forms in $\K[x_0,x_1]$, $\deg G= \binom{d}{2}$, $\deg F= \binom{d-1}{2}+1$,
and $F$ and $G$ have no common zero in $\PP^1$. 
\end{proof}

We would like to make some remarks on the polynomial $P$. It is divisible by $x_1^{d-2}$. Indeed,
$P$ is a form in $x_1$ and $x_0$ of degree $ \binom{d-1}{2}$
with coefficients in $\K[a_1,\ldots,a_d]$. It is divisible by the Vandermonde determinant $V(a_1,\ldots,a_d)$
because it is antisymmetric in the $a_i$'s. Furthermore, the coefficient of $x_1^{\alpha_1} x_0^{\alpha_0}$ in $P$ is an antisymmetric polynomial
of degree $2\alpha_1+\alpha_0+1$ in the $a_i$'s: in order for it to be nonzero, it is necessary that
$2 \alpha_1 + \alpha_1 +1\geqslant \deg V= \binom{d}{2}$. Since $ \alpha_1+\alpha_0 = \binom{d-1}{2}$, we deduce
$ \alpha_1 \geqslant d-2$. This means that $P$ is divisible by $x_1^{d-2}$, and the coefficient of $x_1^{d-2}x_0^{\binom{d-1}{2}}$ is  $ V(a_1, \ldots, a_d)$ times a constant $-c_d$ that depends on $d$ but not on the $a_i$'s.  
To compute $c_d$, we eliminate the highest power of $x_0$ in each column by subtracting the correct linear combination of the previous columns: we obtain $c_2 = 1$ and
\begin{equation}\label{eq:cd}
c_d = \sum_{k=1}^{d-2}(-1)^{k+1} \binom{d-1-k}{k}c_{d-k},\qquad d \geqslant 3.
\end{equation}
The first few values are $c_3=1$, $c_4=2$, $c_5=5$, $c_6=14$, $c_7=42$, and at the end of the chapter we will prove that $c_d$ coincides with  the $(d-2)$-th Catalan number.\index{Catalan numbers}

\begin{remark}
We can give a geometric interpretation of the condition that the sums
$a_i+a_j$ be all distinct. In the family constructed in the proof of the theorem,
the union of the two lines $L_{a_i}$ and $L_{a_j}$ specializes to the planar double line
$x_3=x_2^2=0$ plus the embedded point  $x_3=x_2= (a_i+a_j)x_1+x_0=0$.
Thus the condition means that these embedded points are all distinct.
\end{remark}

\begin{example}
If the condition that the $a_i+a_j$ be all distinct is not satisfied, we expect the limit to acquire embedded
points. The reason is that in this case the proof of Theorem \ref{main} shows
that the polynomials $F$ and $G$ have a common zero, and so $x_3 G-x_2^{d-1}F$ is no longer irreducible.
For a specific example, we take $d=4$ and $a_1=0$, $a_2=1$, $a_3=2$ and $a_4=3$ (in characteristic $\neq 2,3$),
so that $a_1+a_4=a_2+a_3=3$.  In this case,
\[
\begin{split}
\frac{x_3 G-x_2^3 F}{12}  = {} &  x_3(x_1+x_0)(x_1+2x_0)(x_1+3x_0)^2(x_1+4x_0)(x_1+5x_0) \\
&{} - 2 x_2^3 x_1^3 (3x_1+x_0)= (3x_1+x_0)(x_3 G_1-x_2^3F_1).
\end{split}
\]
The initial ideal of $I_C$ is
\[
\begin{split}
\IN_{(4,2,1,1)}(I_C) = \big( & x_3^{2}, 6 x_3 x_2 x_1^{2}+2 x_3 x_2 x_1 x_0, x_2^{4}, x_3 x_2^{2} x_0,6 x_3 x_2^{2} x_1,\\
 &  6 x_3 x_2 x_1 x_0^{2}+2 x_3 x_2 x_0^{3}, x_3G_1-x_2^3F_1, 6 x_3 x_2^{3}\big)
 \end{split}
\]
The saturation of this ideal is
\[
\left( x_3^{2},3 x_3 x_2 x_1+x_3 x_2 x_0,x_3 x_2^{2},x_2^{4},x_3G_1-x_2^3F_1 \right),
\]
therefore the limit $C_0$ consists of the unique $4$ structure supported on $L$
contained in the surface $x_3G_1-x_2^3F_1$, which by Proposition \ref{extr_mls} is an extremal curve
of genus $-2$, plus an embedded point (of equation $3x_1+x_0=0$ on $L$).
\end{example}

\begin{corollary}
Let $C$ be an effective divisor of bidegree $(d,0)$ on a smooth quadric surface. Then every curve in the biliaison class
of $C$ is in the connected component of the extremal curves in its Hilbert scheme. In particular,
this holds for every curve on a smooth quadric surface.
\end{corollary}
\begin{proof}
The case $0 \leqslant d \leqslant 2$ has been proven by Hartshorne \cite{HartshorneConnectedness}.  When $d\geqslant 3$, the statement follows from Theorem \ref{main} and \cite[Theorem 2.3]{SchlesingerFootnote} because $C$ is a minimal curve.
\end{proof}

\subsection*{Catalan numbers}
\index{Catalan numbers|(}
\begin{definition}
Let us define for any $n > k \geqslant 2$
\begin{equation}
\displaystyle
\dCoef{n}{k} = \begin{cases} 
1, &\text{if } k=2\\
\sum\limits_{i=k+1}^d \tCoef{i}{k-1}, & \text{otherwise}.
\end{cases}
\end{equation}
\end{definition}

\begin{table}[!ht]
\small
\begin{center}
\begin{tabular}{ l | cccccccccccccc }
 & $n=3$& $n=4$& $n=5$& $n=6$& $n=7$& $n=8$& $n=9$& $n=10$& $n=11$ \\
 \hline
$k=2$ & 1 & 1 & 1 & 1 & 1 & 1 & 1 & 1 & 1  \\
$k=3$ &  & 1 & 2 & 3 & 4 & 5 & 6 & 7 & 8\\
$k=4$ &  &  & 2 & 5 & 9 & 14 & 20 & 27 & 35\\
$k=5$ &  &  &  & 5 & 14 & 28 & 48 & 75 & 110\\
$k=6$ &  &  &  &  & 14 & 42 & 90 & 165 & 275\\
$k=7$ &  &  &  &  &  & 42 & 132 & 297 & 572\\
$k=8$ &  &  &  &  &  & & 132 & 429 & 1001\\
$k=9$ &  &  &  &  &  & & & 429 & 1430\\
$k=10$ &  &  &  &  &  & & & & 1430\\
\end{tabular}
\end{center}
\normalsize
\caption{First values of $\tCoef{n}{k}$.}
\end{table}

\begin{proposition}
The number $\tCoef{d}{d-1}$ coincides with the Catalan number $C(d-3)$.
\end{proposition}
\begin{proof}
The definition of $\tCoef{n}{k}$ is the same as the construction of the Catalan numbers through the generalized Pascal triangle (see \cite[Example 3.5.5]{StanleyEnumComb}).
\end{proof}

We recall a binomial identity that will be very useful in the following. For any $n, m$
\begin{equation}
\sum	_{i=0}^m (-1)^i \binom{n}{i} = (-1)^m\binom{n-1}{m} 
\end{equation}
and equivalently for $n > 0$
\begin{equation}\label{eq:binomialIdentity}
\sum_{i=1}^m (-1)^{i+1} \binom{n}{i} = 1+ (-1)^{m+1} \binom{n-1}{m}.
\end{equation}

\begin{lemma}
For the coefficients $c_d$ described in \eqref{eq:cd},
\begin{equation}\label{eq:newPartialExpression}
c_d = \sum_{h=0}^{a-2} \dCoef{d-h}{a-h}\, c_{2+h}  + \sum_{k=1}^{d-a-1} (-1)^{k+1}\binom{d-k-a}{k} c_{d-k},\quad \forall\ d > a \geqslant 2.
\end{equation}
\end{lemma}
\begin{proof}
We proceed by induction on $a$. Let us consider any $c_d$ with $d > 2$. The sum
\[
\sum_{i=3}^d c_i = \sum_{i=3}^d \left(\sum_{k=1}^{i-2} (-1)^k \binom{i-k-1}{k} c_{i-k}\right).
\]
can be rewritten as
\[
\sum_{i=3}^d c_i = \sum_{k=1}^{d-2} \left( \sum_{i=1}^{k} (-1)^{i+1} \binom{d-k-1}{i} \right) c_{d-k}
\]
as suggested by the following diagram
\[
\begin{array}{l c c r}
\phantom{{}={}} c_d &&& \binom{d-2}{1} c_{d-1} - \binom{d-3}{2}c_{d-2} + \binom{d-4}{3} c_{d-3} - \binom{d-5}{4}c_{d-4} + \ldots + \binom{1}{d-2} c_2\\
{}+c_{d-1} &&&  +\binom{d-3}{1}c_{d-2} - \binom{d-4}{2} c_{d-3} + \binom{d-5}{3}c_{d-4} + \ldots + \binom{1}{d-3}c_2 \\
{}+c_{d-2} &&&  + \binom{d-4}{1} c_{d-3} - \binom{d-5}{2}c_{d-4} + \ldots + \binom{1}{d-4}c_2\\
\qquad\vdots&&& \vdots\quad\ \,\\
{}+c_3 &&& \binom{1}{1} c_2
\end{array}
\]
Since $d - k - 1$ is greater that 0 for each $k$, we can apply the binomial identity \eqref{eq:binomialIdentity}, obtaining
\[
\begin{split}
\sum_{i=3}^d c_i&{} = \sum_{k=1}^{d-2} \left( 1 + (-1)^{k+1} \binom{d-k-2}{k} \right) c_{d-k} = \\
&{} = \sum_{i=2}^{d-1} c_i + \sum_{k=1}^{d-2} (-1)^{k+1} \binom{d-k-2}{k} c_{d-k}.
\end{split}
\]
The binomial coefficient corresponding to $k=d-2$ surely vanishes and $1 = \tCoef{d}{2}$ so we proved
\[
c_{d} = \dCoef{d}{2}\,c_2 + \sum_{k=1}^{d-3} (-1)^{k+1} \binom{d-k-2}{k} c_{d-k}.
\]
Let us now suppose that the statement is true for $a-1$, that is for every $d > a-1$
\[
c_d = \sum_{h=0}^{a-3} \dCoef{d-h}{a-1-h} \, c_{2+h} + \sum_{k=1}^{d-(a-1)-1} (-1)^{k+1}\binom{d-k-(a-1)}{k} c_{d-k}.
\]
With the same reasoning applied before we compute
\[
\begin{split}
\sum_{i={a+1}}^d c_i &{} = \sum_{i=a+1}^d \left(\sum_{h=0}^{a-3} \dCoef{d-h}{a-1-h} \, c_{2+h}\right) + \\
&{}\qquad + \sum_{i=a+1}^d \left( \sum_{k=1}^{d-(a-1)-1} (-1)^{k+1}\binom{d-k-(a-1)}{k} c_{d-k} \right).
\end{split}
\]
For the first part of the sum, we can change the order of the summations obtaining
\[
\begin{split}
\sum_{i=a+1}^d \left(\sum_{h=0}^{a-3} \dCoef{d-h}{a-1-h} \, c_{2+h}\right) = \sum_{h=0}^{a-3} \left(\sum_{i=a+1}^d \dCoef{i-h}{a-1-h}\right) c_{2+h} = {}&\\
=\sum_{h=0}^{a-3} \left(\sum_{j=(a-h)+1}^{d-h} \dCoef{j}{(a-h)-1}\right) c_{2+h} = \sum_{h=0}^{a-3} \dCoef{d-h}{a-h}\, c_{2+h}.&
\end{split}
\]
For the second part of the sum, we change the order and we use again the binomial identity \eqref{eq:binomialIdentity}, since $d-k-(a-1) > 0,\ \forall\ k$:
\[
\begin{split}
&\sum_{i=a+1}^d \left( \sum_{k=1}^{d-(a-1)-1} (-1)^{k+1}\binom{d-k-(a-1)}{k} c_{d-k} \right) ={}\\
&\quad{} = \sum_{k=1}^{d-(a-1)-1} \left( \sum_{i=1}^{k} (-1)^{i+1} \binom{d-k-(a-1)}{i} \right) c_{d-k} = \\
&\quad{} = \sum_{k=1}^{d-(a-1)-1} \left( 1 + (-1)^{k+1} \binom{d-k-a}{k} \right) c_{d-k} =\\
&\quad{} = \sum_{i=a}^{d-1} c_i + \sum_{k=1}^{d-(a-1)-1} (-1)^{k+1} \binom{d-k-a}{k} c_{d-k}.
\end{split}
\]
For $k=d-(a-1)-1$, the binomial coefficient vanishes and $1=\tCoef{d-(a-2)}{a-(a-2)}$ so
\[
\begin{split}
c_d&{} = \sum_{h=0}^{a-3} \dCoef{d-h}{a-h}\, c_{2+h} + c_{a} + \sum_{k=1}^{d-a-1} (-1)^{k+1} \binom{d-k-a}{k} c_{d-k} =\\
&{} = \sum_{h=0}^{a-2} \dCoef{d-h}{a-h}\, c_{2+h} + \sum_{k=1}^{d-a-1} (-1)^{k+1} \binom{d-k-a}{k} c_{d-k}. \qedhere
\end{split}
\] 
\end{proof}

\begin{proposition}
\begin{equation}\label{eq:explicitCd}
c_d = \dCoef{d+1}{d} = C(d-2).
\end{equation}
\end{proposition}
\begin{proof}
We prove the equality by induction on $d$. Obviously $c_2 = 1 = \tCoef{3}{2} = C(0)$ and let us assume the statement true for any $c_k,\ 2\leqslant k < d$.
Considering the equality \eqref{eq:newPartialExpression} with $a=d-1$ the second sum is empty, so that
\[
c_d = \sum_{h=0}^{d-3} \dCoef{d-h}{d-h-1} c_{2+h}.
\]
By the inductive hypothesis
\[
\sum_{h=0}^{d-3} \dCoef{d-h}{d-h-1} c_{2+h} = \sum_{h=0}^{d-3} C(d-h-3) C(2+h-2) = \sum_{h=0}^D C(D-h) C(h).
\]
This is the recursive definition of Catalan numbers (see \cite[Section 7.5]{ConcreteMathematics}) so 
\[
\sum_{h=0}^D C(D-h) C(h) = C(D+1) \quad\Rightarrow\quad c_d = C(d-2) = \dCoef{d+1}{d}. \qedhere
\]
\end{proof}

\begin{corollary}
\begin{equation}
c_d = \frac{1}{d-1}\binom{2d-4}{d-2}.
\end{equation}
\end{corollary}

\begin{corollary}[A new recurrence relation for Catalan numbers]
\begin{equation}
C(0) = 1, \qquad C(n) = \sum_{k=1}^{n}(-1)^{k+1} \binom{n-k+1}{k}C(n-k).
\end{equation}
\end{corollary}
\index{Catalan numbers|)}

\appendix

\chapter{The \textit{Macaulay2} package \texttt{HilbertSchemesEquations}}\label{ch:Hilb2P2}

This\hfill chapter\hfill is\hfill supposed\hfill to\hfill be\hfill a\hfill handbook\hfill for\hfill the\hfill \textit{Macaulay2}\hfill \cite{M2}\hfill package\\ \texttt{HilbertSchemeEquations.m2}. We will introduce and explain the main functions of the package through the complete computation of the examples introduced in Chapter \ref{ch:HilbertScheme}, that is computing various sets of equations for the Hilbert scheme $\Hilb{2}{2}$.

\section{Basic features}
Firstly we need a method for computing the Gotzmann number of a Hilbert polynomial. \texttt{HilbertSchemesEquations.m2} provides two methods:
\begin{itemize}
\item \texttt{gotzmannDecomposition}, taking as input a single variable polynomial $p(t)$ and returning the Gotzmann decomposition of the polynomial as list of pairs $\{\ldots,(a_i,b_i),\ldots\}$ such that $p(t) = \sum_i \binom{t+a_i}{b_i}$:

\begin{code}
\begin{verbatim}
gotzmannDecomposition = method(TypicalValue => List)
--  INPUT: p, polynomial (one variable).
-- OUTPUT: a list of pairs {...(a_i,b_i)...} containing the 
--         Gotzmann decomposition of p.
--  ERROR: if numgens(ring(p)) > 1.
--         if coefficientRing(ring(p)) =!= ZZ and
--            coefficientRing(ring(p)) =!= QQ.
--         if p is not an admissible Hilbert polynomial.
\end{verbatim}
\end{code}
\item \texttt{gotzmannNumber}, taking as input a single variable polynomial $p(t)$ and returning the Gotzmann number:

\begin{code}
\begin{verbatim}
gotzmannNumber = method(TypicalValue => ZZ)
--  INPUT: p, polynomial (one variable).
-- OUTPUT: the Gotzmann number of p.
--  ERROR: if numgens(ring(p)) > 1.
--         if coefficientRing(ring(p)) =!= ZZ and 
--            coefficientRing(ring(p)) =!= QQ.
--         if p is not an admissible Hilbert polynomial.
\end{verbatim}
\end{code}
\end{itemize}

\begin{example}
We test these two methods on the polynomials $p(t) = 4t+1$ and $q(t) = \frac{1}{3}t^2 - \frac{2}{3}t + 1$. 

\begin{code}
\begin{verbatim}
Macaulay2, version 1.4
i1 : loadPackage "HilbertSchemesEquations";
i2 : R = QQ[t];
i3 : p = 4*t+1;
i4 : q = (t^2-2*t+3)/3;
i5 : gotzmannNumber p
o5 = 7
i6 : gotzmannDecomposition p
o6 = {(1, 1), (0, 1), (-1, 1), (-2, 1), (-4, 0), (-5, 0), (-6, 0)}
o6 : List
i7 : gotzmannDecomposition q
stdio:7:1:(3): error: argument 1: not admissible Hilbert polynomial
\end{verbatim}
\end{code}
Hence $q(t)$ is not admissible (as expected) and $p(t)$ can be decomposed as
\[
\binom{t+1}{1}+\binom{t}{1}+\binom{t-1}{1}+\binom{t-2}{1} + \binom{t-4}{0} + \binom{t-5}{0} + \binom{t-6}{0}.
\]
\end{example}

To compute the Pl\"ucker relations defining the Grassmannian, we will use the function \texttt{Grassmannian(ZZ,ZZ)} provided by \textit{Macaulay2}. We remark that to compute the ideal defining $\Grass{q}{N}{\K} \subset \PP_{\K}^{\binom{N}{q}-1}$ we need to call \texttt{Grassmannian($q-1$,$N-1$)}.

The package \texttt{HilbertSchemesEquations} gives methods to compute the generic generators of a subspace (and of its exterior powers) parametrized by $\Grass{q}{N}{\K}$ according to Definition \ref{def:genericGens}.
\begin{itemize}
\item \texttt{genericSubspaceGen}. This method requires 3 arguments: 
\begin{enumerate}
\item the exterior algebra generated by the basis $\{\vv_1,\ldots,\vv_N\}$ of the base vector space $V \simeq \K^N$, with coefficient in the ring of the Grassmannian $\Grass{q}{N}{\K}$;
\item the dimension $q$ of the subspaces parametrized by $\Grass{q}{N}{\K}$;
\item a multiindex $\JJ$; 
\end{enumerate}
and it returns the element $\Lambda^{(q-\vert \JJ\vert)}_{\JJ}$ associated to $\Grass{q}{N}{\K}$.
\begin{code}
\begin{verbatim}
genericSubspaceGen = method(TypicalValue => RingElement)
--  INPUT: A, an exterior algebra with coefficientRing(A) 
--            corresponding to a ring defining a grassmannian. 
--            gens(A) has to be a basis of the base vector 
--            space of the grassmannian.
--         q, the dimension of the subspaces parametrized by 
--            the grassmannian.
--         J, a set of indices.
-- OUTPUT: the generator Lambda^{(q-#J)}_J.
--  ERROR: if q < 1 or q > numgens(A).
--         if numgens(coefficientRing(A)) != binomial(numgens(A),q).
--         if #J >= q,
--         if min(J) < 0 or max(J) > numgens(A)-1.
\end{verbatim}
\end{code}
\item \texttt{genericSetOfSubspaceGens}. Also this method requires 3 arguments: the first 2 are the same required by \texttt{genericSubspaceGen} and the third is the order $s$ of the exterior power. It returns the set of generators $\Gamma^{(s)}$ associated to $\Grass{q}{N}{\K}$.

\begin{code}
\begin{verbatim}
genericSetOfSubspaceGens = method(TypicalValue => List)
--  INPUT: A, an exterior algebra with coefficientRing(A) 
--            corresponding to a ring defining a grassmannian.
--            gens(A) has to be a basis of the base vector 
--            space of the grassmannian.
--         q, the dimension of the subspaces parametrized by
--            the grassmannian.
--         s, the order of the exterior power.
-- OUTPUT: list containing the set of generators Gamma^{(s)}.
--  ERROR: if q < 1 or q > numgens(A).
--         if numgens(coefficientRing(A)) != binomial(numgens(A),q).
--         if s < 1 or s > q.
\end{verbatim}
\end{code}
\end{itemize}

\begin{example}\label{ex:MainExampleGrassCODE}
We introduce these two new methods applying them to the case of $\Grass{4}{6}{\K}$ discussed in Example \ref{ex:MainExampleGrass} and Example \ref{ex:MainExampleGens}.

\begin{code}
\begin{verbatim}
Macaulay2, version 1.4
i1 : loadPackage "HilbertSchemesEquations";
i2 : grass = Grassmannian (3,5,CoefficientRing=>QQ);
o2 : Ideal of QQ[p       , p       , p       , p       , p       , p       , 
                  0,1,2,3   0,1,2,4   0,1,3,4   0,2,3,4   1,2,3,4   0,1,2,5   
     --------------------------------------------------------------------------
     p       , p       , p       , p       , p       , p       , p       , 
      0,1,3,5   0,2,3,5   1,2,3,5   0,1,4,5   0,2,4,5   1,2,4,5   0,3,4,5 
     --------------------------------------------------------------------------
     p       , p       ]
      1,3,4,5   2,3,4,5
i3 : pluckerRelations = first entries gens grass;
i4 : #pluckerRelations 
o4 = 15
i5 : for i from 0 to #pluckerRelations-1 do (print(pluckerRelations#i););
     p       p        - p       p        + p       p
      1,2,4,5 0,3,4,5    0,2,4,5 1,3,4,5    0,1,4,5 2,3,4,5
     p       p        - p       p        + p       p
      1,2,3,5 0,3,4,5    0,2,3,5 1,3,4,5    0,1,3,5 2,3,4,5
     p       p        - p       p        + p       p
      1,2,3,4 0,3,4,5    0,2,3,4 1,3,4,5    0,1,3,4 2,3,4,5
     p       p        - p       p        + p       p
      1,2,3,5 0,2,4,5    0,2,3,5 1,2,4,5    0,1,2,5 2,3,4,5
     p       p        - p       p        + p       p
      1,2,3,4 0,2,4,5    0,2,3,4 1,2,4,5    0,1,2,4 2,3,4,5
     p       p        - p       p        + p       p
      1,2,3,5 0,1,4,5    0,1,3,5 1,2,4,5    0,1,2,5 1,3,4,5
     p       p        - p       p        + p       p
      0,2,3,5 0,1,4,5    0,1,3,5 0,2,4,5    0,1,2,5 0,3,4,5
     p       p        - p       p        + p       p
      1,2,3,4 0,1,4,5    0,1,3,4 1,2,4,5    0,1,2,4 1,3,4,5
     p       p        - p       p        + p       p
      0,2,3,4 0,1,4,5    0,1,3,4 0,2,4,5    0,1,2,4 0,3,4,5
     p       p        - p       p        + p       p
      1,2,3,4 0,2,3,5    0,2,3,4 1,2,3,5    0,1,2,3 2,3,4,5
     p       p        - p       p        + p       p
      1,2,3,4 0,1,3,5    0,1,3,4 1,2,3,5    0,1,2,3 1,3,4,5
     p       p        - p       p        + p       p
      0,2,3,4 0,1,3,5    0,1,3,4 0,2,3,5    0,1,2,3 0,3,4,5
     p       p        - p       p        + p       p
      1,2,3,4 0,1,2,5    0,1,2,4 1,2,3,5    0,1,2,3 1,2,4,5
     p       p        - p       p        + p       p
      0,2,3,4 0,1,2,5    0,1,2,4 0,2,3,5    0,1,2,3 0,2,4,5
     p       p        - p       p        + p       p
      0,1,3,4 0,1,2,5    0,1,2,4 0,1,3,5    0,1,2,3 0,1,4,5
\end{verbatim}
\end{code}
Then we introduce the exterior algebra generated by the basis $\{\vv_0,\ldots,\vv_5\}$ of the base vector space $V$ and we compute the sets of generators of any exterior power of a subspace parametrized by $\Grass{4}{6}{\K}$.

\begin{code}
\begin{verbatim}
i6 : G = ring(grass)/grass;
i7 : A = G[v_0..v_5,SkewCommutative=>true];
i8 : genericSubspaceGen (A,4,{0,4,5})

o8 = p       v  + p       v  + p       v
      0,1,4,5 1    0,2,4,5 2    0,3,4,5 3
o8 : A
i9 : Gamma1 = genericSetOfSubspaceGens (A,4,1);
i10 : #Gamma1
o10 = 20
i11 : genericSubspaceGen (A,4,{1,3})
o11 = - p       v v  + p       v v  - p       v v  + p       v v  - p       v v 
         0,1,2,3 0 2    0,1,3,4 0 4    1,2,3,4 2 4    0,1,3,5 0 5    1,2,3,5 2 5
      --------------------------------------------------------------------------
      + p       v v
         1,3,4,5 4 5
o11 : A
i12 : Gamma2 = genericSetOfSubspaceGens (A,4,2);
i13 : #Gamma2
o13 = 15
i14 : genericSubspaceGen (A,4,{2})
o14 = p       v v v  + p       v v v  - p       v v v  - p       v v v  +
       0,1,2,3 0 1 3    0,1,2,4 0 1 4    0,2,3,4 0 3 4    1,2,3,4 1 3 4  
      --------------------------------------------------------------------------
      p       v v v  - p       v v v  - p       v v v  - p       v v v  -
       0,1,2,5 0 1 5    0,2,3,5 0 3 5    1,2,3,5 1 3 5    0,2,4,5 0 4 5  
      --------------------------------------------------------------------------
      p       v v v  + p       v v v
       1,2,4,5 1 4 5    2,3,4,5 3 4 5
o14 : A
i15 : Gamma3 = genericSetOfSubspaceGens (A,4,3);
i16 : #Gamma3
o16 = 6
i17 : genericSubspaceGen (A,4,{})
o17 = p       v v v v  + p       v v v v  + p       v v v v  + p       v v v v 
       0,1,2,3 0 1 2 3    0,1,2,4 0 1 2 4    0,1,3,4 0 1 3 4    0,2,3,4 0 2 3 4
      --------------------------------------------------------------------------
      + p       v v v v  + p       v v v v  + p       v v v v  +
         1,2,3,4 1 2 3 4    0,1,2,5 0 1 2 5    0,1,3,5 0 1 3 5  
      --------------------------------------------------------------------------
      p       v v v v  + p       v v v v  + p       v v v v  + p       v v v v 
       0,2,3,5 0 2 3 5    1,2,3,5 1 2 3 5    0,1,4,5 0 1 4 5    0,2,4,5 0 2 4 5
      --------------------------------------------------------------------------
      + p       v v v v  + p       v v v v  + p       v v v v  +
         1,2,4,5 1 2 4 5    0,3,4,5 0 3 4 5    1,3,4,5 1 3 4 5  
      --------------------------------------------------------------------------
      p       v v v v
       2,3,4,5 2 3 4 5
o17 : A
\end{verbatim}
\end{code}
\end{example}

\section{Hilbert scheme equations}\label{sec:HilbEquationsCode}

\subsection*{Gotzmann equations}
The\hfill package\hfill \texttt{HilbertSchemesEquations}\hfill provides\hfill the\hfill function\\ \texttt{GotzmannHilbEquations}:

\begin{code}
\begin{verbatim}
GotzmannHilbEquations = method(TypicalValue => List, 
          Options => {PluckerRelations => true,SingleGrassmannian => true})
--  INPUT: p, admissible Hilbert polynomial.
--         n, dimension of the projective space.
-- OUTPUT: a list containing:
--           - #0 the ring in which the Hilbert scheme is embedded;
--           - #1 the ideal defining the Hilbert scheme.
-- OPTION: PluckerRelations (Boolean), default value true.
--            If true, the ideal of Plucker relations is computed.
--            If false, the Plucker coordinates are considered
--            without relations among them.
--         SingleGrassmannian (Boolean), default value true.
--            If true, the ideal is embedded in a single 
--            grassmannian. If false, the ideal is embedded
--            in the product of two grassmannians.
--  ERROR: if numgens(ring(p)) > 1.
--         if coefficientRing(ring(p)) != ZZ and 
--            coefficientRing(ring(p)) != QQ
--         if p is not admissible.
--         if n < 1.
--         if first degree (p) >= n.
\end{verbatim}
\end{code}
This method requires two inputs: a Hilbert polynomial $p(t)$ and a dimension $n$ of a projective space $\PP^n_{\K}$. It returns a sequence with two elements: the first is the ring in which the ideal of the Hilbert scheme is computed and the second is just the ideal.

There are two options:
\begin{itemize}
\item \texttt{PluckerRelations}, since computing the ideal of the Pl\"ucker relations is in general a hard task, it is possible to tell the function to ignore Pl\"ucker relations;
\item \texttt{SingleGrassmannian}, with this option it is possible to choose the embedding of the Hilbert scheme: if in the single Grassmannian $\Grass{q(r)}{N(r)}{\K}$ or in the product $\Grass{q(r)}{N(r)}{\K} \times \Grass{p(r+1)}{N(r+1)}{\K}$. To compute the image of the Hilbert scheme $\Hilb{n}{p(t)} \subset \Grass{q(r)}{N(r)}{\K} \times \Grass{p(r+1)}{N(r+1)}{\K}$ by the projection $\pi$ on the first factor, we use the projective elimination theory (see \cite[Chapter 8 Section 5]{CLOiva}). Given the map $\pi: \Proj \K[\Delta] \times \Proj \K[\nabla] \rightarrow \Proj \K[\Delta]$ and the ideal $\mathcal{I}_{\mathcal{H}}$ defining $\Hilb{n}{p(t)} \subset \Proj \K[\Delta] \times \Proj \K[\nabla]$, the ideal $\widehat{\mathcal{I}}_{\mathcal{H}}$ defining $\pi(\Hilb{n}{p(t)})$ can be computed applying the standard elimination algorithm on the affine open subset of the product of projective spaces. Denoted by $U_{\II}$ the open subset of $\Proj \K[\Delta]$ where $\Delta_{\II} \neq 0$ and by $V_{\JJ}$ the open subset of $\Proj \K[\nabla]$ where $\nabla_{\JJ} \neq 0$, we dehomogenize the ideal $\mathcal{I}_{\mathcal{H}}$, then we eliminate the variables $\nabla$ and finally we homogenize the ideal obtained with the variable $\Delta_\II$. Repeating this procedure for any pair $\II,\JJ$ we recover the ideal $\widetilde{\mathcal{I}}_{\mathcal{H}}$.
\end{itemize}

To recap, let us consider the Hilbert polynomial $p(t)$ with Gotzmann number $r$ and the projective space $\PP^n_{\K}$. Moreover let $\Grass{q(r)}{N(r)}{\K} \subset \PP^{\binom{N(r)}{q(r)}-1}_{\K} = \Proj \K[\Delta_{\II}]$ be defined by the ideal $\mathcal{Q}_1$ and $\Grass{p(r)}{N(r+1)}{\K} \subset \PP^{\binom{N(r+1)}{p(r+1)}-1}_{\K} = \Proj \K[\nabla_{\JJ}]$ defined by the ideal $\mathcal{Q}_2$. Finally let $\mathcal{B}_{\mathcal{H}}$ be the ideal generated by the bilinear equations introduced in the proof of Theorem \ref{th:GotzmannEquations}. 

Calling the method \texttt{GotzmannHilbEquations} on the pairs $(p(t),n)$ and varying the option, there are the four possibilities illustrated in the following table:
\begin{center}
\begin{tikzpicture}
\draw [-] (-1.5,0)--(10,0);
\draw [-] (0,0.75)--(10,0.75);
\draw [-] (0,1.5) -- (10,1.5);
\draw [-] (10,-4) -- (10,1.5);
\draw [-] (0,-4) -- (0,1.5);
\draw [-] (5,-4) -- (5,0.75);
\node at (5,1.125) [] {\texttt{SingleGrassmannian}};
\node at (7.5,0.375) [] {\textbf{\texttt{true}}};
\node at (2.5,0.375) [] {\textbf{\texttt{false}}};
\node at (2.5,-1) [] {$\begin{array}{c}\mathcal{I}_{\mathcal{H}} = (\mathcal{B}_{\mathcal{H}},\mathcal{Q}_1,\mathcal{Q}_2) \\ \Hilb{n}{p(t)} \simeq \Proj \frac{\K[\Delta,\nabla]}{\mathcal{I}_{\mathcal{H}}}\end{array}$};
\node at (7.5,-1) [] {$\begin{array}{c} \widehat{\mathcal{I}}_{\mathcal{H}} \subset \K[\Delta] \\ \Hilb{n}{p(t)} \simeq \Proj \frac{\K[\Delta]}{\widehat{\mathcal{I}}_{\mathcal{H}}}\end{array}$};
\node at (2.5,-3) [] {$\mathcal{B}_{\mathcal{H}} \subset \K[\Delta,\nabla]$};
\node at (7.5,-3) [] {$\widehat{\mathcal{B}}_{\mathcal{H}} \subset \K[\Delta]$};
\draw [-] (-1.5,-4)--(10,-4);
\node at (-0.375,-1) [rotate=90] {\textbf{\texttt{true}}};
\node at (-0.375,-3) [rotate=90] {\textbf{\texttt{false}}};
\node at (-1.125,-2) [rotate=90] {\texttt{PluckerRelations}}; 
\draw [-] (-0.75,-2)--(10,-2);
\draw [-] (-0.75,0)--(-0.75,-4);
\draw [-] (-1.5,0)--(-1.5,-4);
\end{tikzpicture}
\end{center}

\begin{example}\label{ex:MainExampleGotzmannCODE} 
We compute entirely Gotzmann equations of the Hilbert scheme \gls{Hilb2P2}, completing Example \ref{ex:MainExampleGotzmann}.

\begin{code}
\begin{verbatim}
Macaulay2, version 1.4
i1 : loadPackage "HilbertSchemesEquations";
i2 : R = QQ[t];
i3 : time Hilb = GotzmannHilbEquations (2_R,2);
     -- 600 bilinear equations
     -- used 10943.3 seconds
i4 : gensHilb = first entries gens (Hilb#1);
i5 : #gensHilb 
o5 = 30376
i6 : gbHilb = first entries gens gb (Hilb#1);
i7 : #gbHilb 
o7 = 50
i8 : hilbertPolynomial (Hilb#1,Projective=>false)
     7 4   15 3   45 2   15
o8 = -i  + --i  + --i  + --i + 1
     8      4      8      4
o8 : QQ[i]
\end{verbatim}
\end{code}

The computation is very long because to compute the projection of the ideal given by the bilinear equations on a single Grassmannian, the method has to eliminate the variables $\nabla$ in each open subset $U_{\II} \times V_{\JJ}$ of the open covering of $\PP^{14} \times \PP^{44}$. Each elimination correspond to a computation of a Gr\"obner basis and there are $15 \cdot 45 = 675$ possible open subsets. 

If we want to embed $\Hilb{2}{2}$ in the product $\Grass{4}{6}{\K}\times\Grass{2}{10}{\K}$, we have to switch to \texttt{false} the option \texttt{SingleGrassmannian} and the computation turns out to be very quick.

\begin{code}
\begin{verbatim}
i9 : time Hilb = GotzmannHilbEquations (2_P,2,SingleGrassmannian=>false);
     -- 600 bilinear equations
     -- used 0.542746 seconds
i10 : gensHilb = first entries gens (Hilb#1);
i11 : #gensHilb 
o11 = 735
i12 : gbHilb = first entries gens gb (Hilb#1);
i13 : #gbHilb 
o13 = 1992
i14 : hilbertPolynomial (Hilb#1,Projective=>false)
           1      16       1     15       739    14      239   13  
o14 = -----------i   + ---------i   + ----------i   + --------i   +
      14631321600      182891520      3657830400      52254720     
      --------------------------------------------------------------------------
         74153   12     41759  11    24703913  10    3192863 9    3190527833 8  
      ----------i   + --------i   + ----------i   + --------i  + -----------i  +
      1045094400      52254720      3657830400      73156608     14631321600    
      --------------------------------------------------------------------------
       8866763 7   674339683 6   30564427 5   1668104321 4   7824563 3  
      --------i  + ---------i  + --------i  + ----------i  + -------i  +
      10450944     261273600      3265920      101606400      508032    
      --------------------------------------------------------------------------
      6898529 2   5477
      -------i  + ----i + 1
       705600     1260
o14 : QQ[i]
\end{verbatim}
\end{code}
\end{example}

\subsection*{Iarrobino-Kleiman equations}
The\hfill package\hfill \texttt{HilbertSchemesEquations}\hfill provides\hfill also\hfill the\hfill code\hfill to\hfill compute\\ Iarrobino-Kleiman global equations for the Hilbert scheme, even if this method is totally unusable because the huge number of product to be computed as shown in Example \ref{ex:MainExampleIarrobinoKleiman}.

\begin{code}
\begin{verbatim}
IKHilbEquations = method(TypicalValue => List, 
                         Options => {PluckerRelations => true})
--  INPUT: p, admissible Hilbert polynomial.
--         n, dimension of the projective space.
-- OUTPUT: a list containing:
--           - #0 the ring in which the Hilbert scheme is embedded;
--           - #1 the ideal defining the Hilbert scheme.
-- OPTION: PluckerRelations (Boolean), default value true.
--            If true, the ideal of Plucker relations is computed.
--            If false, the Plucker coordinates are considered
--            without relations among them.
--  ERROR: if numgens(ring(p)) > 1.
--         if coefficientRing(ring(p)) != ZZ and 
--            coefficientRing(ring(p)) != QQ
--         if p is not admissible.
--         if n < 1.
--         if first degree (p) >= n.
\end{verbatim}
\end{code}

\subsection*{Bayer-Haiman-Sturmfels equations}
The method implementing the strategy introduced in the proof of Theorem \ref{th:BayerHaimanSturmfelsEquations} for computing the equations of $\Hilb{n}{p(t)}$ is called \texttt{BHSHilbEquations}.

\begin{code}
\begin{verbatim}
BHSHilbEquations = method(TypicalValue => List, 
                          Options => {PluckerRelations => true})
--  INPUT: p, admissible Hilbert polynomial.
--         n, dimension of the projective space.
-- OUTPUT: a list containing:
--           - #0 the ring in which the Hilbert scheme is embedded;
--           - #1 the ideal defining the Hilbert scheme.
-- OPTION: PluckerRelations (Boolean), default value true.
--            If true, the ideal of Plucker relations is computed.
--            If false, the Plucker coordinates are considered
--            without relations among them.
--  ERROR: if numgens(ring(p)) > 1.
--         if coefficientRing(ring(p)) != ZZ and 
--            coefficientRing(ring(p)) != QQ
--         if p is not admissible.
--         if n < 1.
--         if first degree (p) >= n.
\end{verbatim}
\end{code}
This function requires as input the same objects as the previous methods and it also has the option \texttt{PluckerRelation} to choose if considering the Pl\"ucker relations or not. If \texttt{PluckerRelations => true} then the relations are used during the computation of the exterior products, in order to reduce any coefficient in the Pl\"ucker coordinates.

\begin{example} 
We test the method \texttt{BHSHilbEquations} on the Hilbert scheme \gls{Hilb2P2}, completing Example \ref{ex:MainExampleBayerHaimanSturmfels}.

\begin{code}
\begin{verbatim}
Macaulay2, version 1.4
i1 : loadPackage "HilbertSchemesEquations";
i2 : R = QQ[t];
i3 : time Hilb = BHSHilbEquations (2_R,2);
     -- 577 exterior products
     -- used 1.67346 seconds
i4 : gensHilb = first entries gens (Hilb#1);
i5 : #gensHilb 
o5 = 3976
i6 : gbHilb = first entries gens gb (Hilb#1);
i7 : #gbHilb 
o7 = 272
i8 : hilbertPolynomial (Hilb#1,Projective=>false)
     7 4   15 3   45 2   15
o8 = -i  + --i  + --i  + --i + 1
     8      4      8      4
o8 : QQ[i]
i9 : time Hilb = BHSHilbEquations (2_R,2,PluckerRelations=>false);
     -- 617 exterior products
     -- used 1.01678 seconds
i10 : gensHilb = first entries gens (Hilb#1);
i11 : #gensHilb 
o11 = 1642
\end{verbatim}
\end{code}
We underline that by this computation the Hilbert scheme $\Hilb{2}{2}$ embedded in $\PP^{14}$ has dimension 4 and degree 21, i.e. $\frac{7}{8} = \frac{21}{4!}$, confirming what stated by Haiman and Sturmfels in \cite[page 756]{HaimanSturmfels}.
\end{example}

\chapter[The \texttt{HSC} java library]{The \texttt{H}(ilbert) \texttt{S}(cheme) \texttt{C}(omputation) java library}\label{ch:HSCpackage}

\section{The description of the library}

The \texttt{HSC} library contains several packages that we now briefly describe.

\begin{description}
\item[\texttt{HSC.math}] It contains the implementations of the rational numbers and some basic mathematical operations and functions.

\item[\texttt{HSC.hilbpoly}] It contains the implementations of Hilbert polynomials and some basic operations on them. For instance given a Hilbert polynomial $p(t)$, it is possible to compute $\Delta p(t)$, its Gotzmann number and Gotzmann decomposition and the saturated lexicographic ideal associated to it.

\item[\texttt{HSC.monomials}] It contains the implementations of monomials, monomial ideals and term orderings. All the basic operation on these objects are made available (elementary moves and evaluation on maximal and minimal variable, regularity of a Borel-fixed ideal, etc.).

\item[\texttt{HSC.borelfixed}] This package contains the implementation of the most important algorithms introduce in the thesis to work on Borel-fixed ideals. Let us look closer at its classes:
\begin{description}
\item[\textnormal{\texttt{PosetGraph}}] contains the implementation of the poset of monomials of a fixed degree with the Borel partial order $\leq_B$ represented as planar graph;
\item[\textnormal{\texttt{BorelGenerator}}] contains the algorithm for computing all the saturated Borel-fixed with chosen number of variables and Hilbert polynomial;
\item[\textnormal{\texttt{BorelInequalitiesSystem}}] contains the algorithms determining whenever a Borel-fixed ideal is a gen/reg/hilb-segment ideal (and computing the relative term ordering);
\item[\textnormal{\texttt{ConnectingPath}}] contains all the algorithms about Borel rational deformations of Borel-fixed ideals.
\end{description}
\item[\texttt{HSC.inequality}] It contains the implementations of the linear inequalities with integer coefficients needed for the simplex algorithm used in the computation of segment ideals.
\item[\texttt{HSC.utilities}] It contains methods to manage the input/output.
\end{description}

\section{Borel-fixed ideals and segment ideals}
Algorithm \ref{alg:BorelGeneratorDFS}, computing all the saturated Borel-fixed ideals with fixed number of variables and Hilbert polynomial, can be tested by the applet \texttt{Borel} \texttt{Generator} (Figure \ref{fig:BorelGenerator})\hfill available\hfill at\hfill \href{http://www.personalweb.unito.it/paolo.lella/HSC/borelGenerator.html}{\texttt{www.personalweb.unito.it/paolo.lella/HSC/}\\ \texttt{borelGenerator.html}}. There are two field to fill:
\begin{description}
\item[\textnormal{\textsf{Projective space}}] requires the dimension of the projective space;
\item[\textnormal{\textsf{Hilbert polynomial}}] requires the Hilbert polynomial as list of coefficients enclosed in square brackets. Any rational coefficient has to be enclosed in round brackets. As example
\[
\begin{split}
 p(t) = 4t \qquad &{}\leadsto \qquad [0,4], \\
 p(t) = \frac{3}{2}t^2 + \frac{5}{2}t - 1 \qquad &{}\leadsto \qquad [-1,(5/2),(3/2)].
\end{split}
\]
\end{description}

In the output window (Figure \ref{fig:BorelGeneratorOutput}), there will be the list of all Borel-fixed saturated ideals, with first element always the lexicographic ideal.

\begin{figure}[!ht]
\begin{center}
\includegraphics[width=12cm]{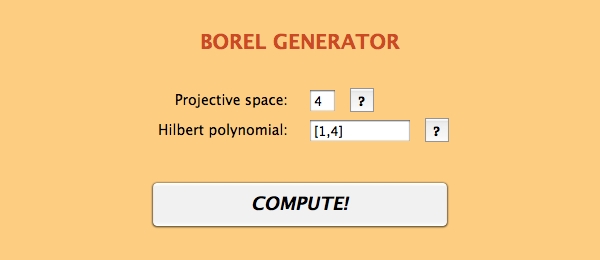}
\caption{\label{fig:BorelGenerator} The applet \texttt{Borel} \texttt{Generator}.}
\end{center}
\end{figure}

\begin{figure}[!ht]
\begin{center}
\includegraphics[width=10cm]{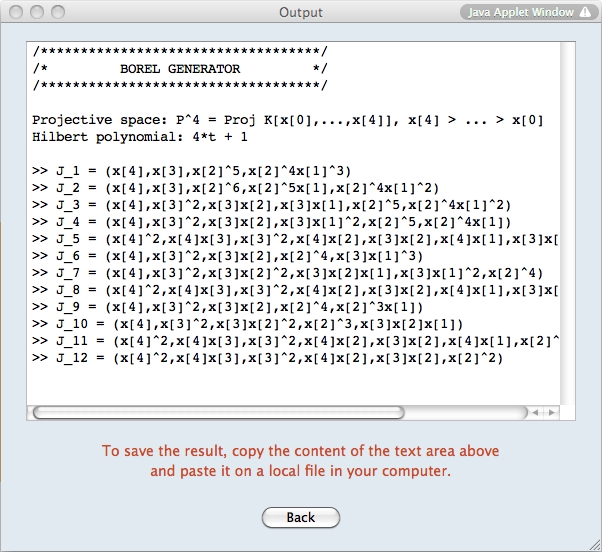}
\caption{\label{fig:BorelGeneratorOutput} The output window of the applet \texttt{Borel} \texttt{Generator}.}
\end{center}
\end{figure}

The applet \texttt{Segment} \texttt{Ideals} (Figure \ref{fig:Segment}) makes available the algorithms to determine whenever a Borel-fixed ideals is hilb/reg/gen-segment ideal (Algorithm \ref{alg:HilbRegSegment}\hfill and\hfill Algorithm\hfill \ref{alg:genSegment}).\hfill  It\hfill can\hfill be\hfill found\hfill at\hfill \href{http://www.personalweb.unito.it/paolo.lella/HSC/segment.html}{\texttt{www.personalweb.unito.it/}\\ \texttt{paolo.lella/HSC/segment.html}} requires three arguments as input:
\begin{description}
\item[\textnormal{\textsf{Projective space}}] needs a positive integer declaring the dimension of the projective space (i.e. $n$ means that the polynomial ring used will be $\K[x_0,\ldots,x_n]$);
\item[\textnormal{\textsf{Borel-fixed ideal}}] needs a string describing a Borel-fixed ideal, with the same syntax used in the output window of the applet \texttt{Borel} \texttt{Generator};
\item[\textnormal{\textsf{Truncation degree}}] needs a positive integer, which will be used to compute if the ideal truncated in such degree is a gen-segment ideal.
\end{description}

\begin{figure}[H]
\begin{center}
\includegraphics[width=12cm]{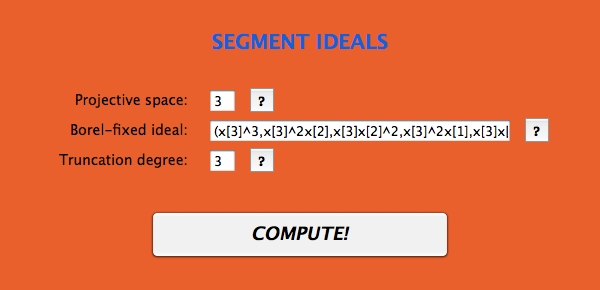}
\caption{\label{fig:Segment} The applet \texttt{Segment} \texttt{Ideals}.}
\end{center}
\end{figure}

In Figure \ref{fig:SegmentOutput}, there is the output window of this applet, with the results of the computation on the ideal discussed in Example \ref{ex:notRegbutGen}. If the ideal is some segment ideal, the applet specifies the term ordering with the first solution (that one with smallest coefficients) of the system of constraints given by the symplex algorithm.

\begin{figure}[H]
\begin{center}
\includegraphics[width=8cm]{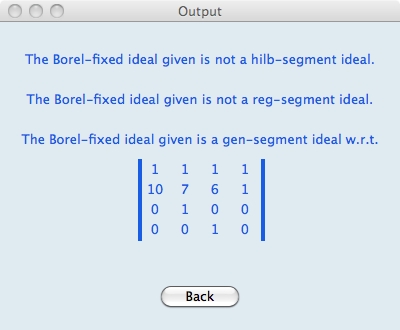}
\caption{\label{fig:SegmentOutput} The output window of the applet \texttt{Segment} \texttt{Ideals}.}
\end{center}
\end{figure}

\section{Borel rational deformations}

The\hfill applet\hfill computing\hfill all\hfill the\hfill Borel\hfill rational\hfill deformations\hfill of\hfill a\hfill Borel-fixed\hfill ideal\\ (Algorithm \ref{alg:allDeformations}) is called \texttt{Borel} \texttt{Rational} \texttt{Deformations} (Figure \ref{fig:Deformations}) and it can be\hfill found\hfill at\hfill \href{http://www.personalweb.unito.it/paolo.lella/HSC/deformations.html}{\texttt{www.personalweb.unito.it/paolo.lella/HSC/}\\ \texttt{deformations.html}}. 

It requires two arguments:
\begin{description}
\item[\textnormal{\textsf{Projective space}}] is the dimension of the projective space (as before).
\item[\textnormal{\textsf{Borel-fixed ideal}}] needs the string of the Borel-fixed ideals that we want to deform, again with the systax used in the applet \texttt{Borel} \texttt{Generator}.
\end{description}
The algorithm computes both simple and composed Borel rational deformations, considering the Borel set defined by the homogeneous piece of the ideal of degree equal to the Gotzmann number of its Hilbert polynomial, so that all the possible Borel rational curves on the Hilbert scheme passing through the given ideal are determined.

\begin{figure}[H]
\begin{center}
\includegraphics[width=12cm]{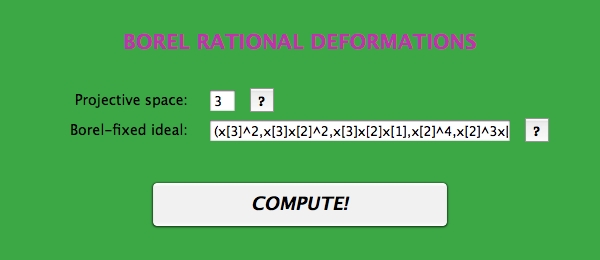}
\caption{\label{fig:Deformations} The applet \texttt{Borel} \texttt{Rational} \texttt{Deformations}.}
\end{center}
\end{figure}

\begin{figure}[H]
\begin{center}
\includegraphics[width=10cm]{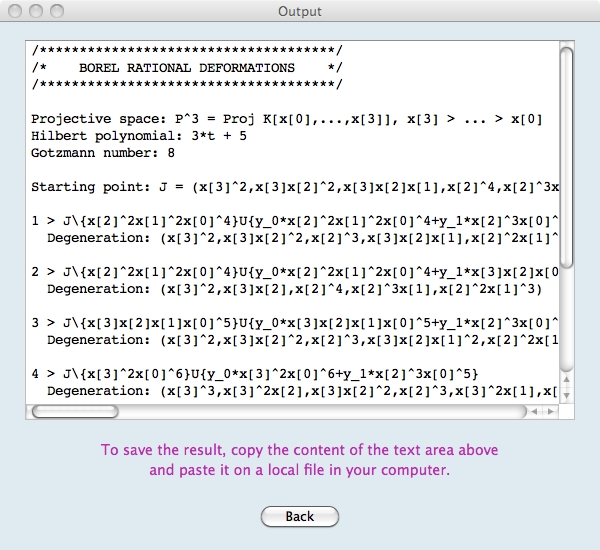}
\caption{\label{fig:DeformationsOutput} The output window of the applet \texttt{Borel} \texttt{Rational} \texttt{Deformations}.}
\end{center}
\end{figure}

\newpage 

To have a global glance of Borel rational curves on a Hilbert scheme, we can use the applet \texttt{Borel} \texttt{Incidence} \texttt{Graph} (Figure \ref{fig:IncidenceGraph}) that computes simple and composed Borel rational deformations of every Borel-fixed ideal with number of variables and Hilbert polynomial fixed (Algorithm \ref{alg:BorelIncidenceGraph}). It requires the same two arguments of the applet \texttt{Borel} \texttt{Generator} (\textsf{Projective space} and \textsf{Hilbert polynomial}) and\hfill it\hfill is\hfill available\hfill at\hfill \href{http://www.personalweb.unito.it/paolo.lella/HSC/borelIncidenceGraph.html}{\texttt{www.personalweb.unito.it/paolo.lella/HSC/}\\ \texttt{borelIncidenceGraph.html}}

\begin{figure}[H]
\begin{center}
\includegraphics[width=12cm]{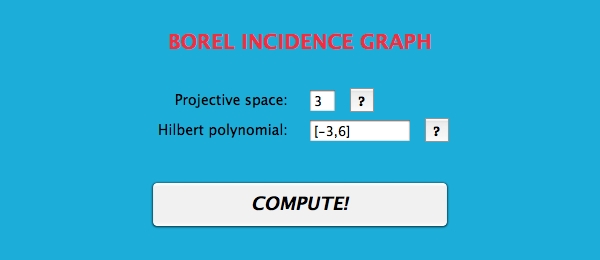}
\caption{\label{fig:IncidenceGraph} The applet \texttt{Borel} \texttt{Incidence} \texttt{Graph}.}
\end{center}
\end{figure}

The output window (Figure \ref{fig:IncidenceGraphOutput}) is divided in two part. The upper part contains the details about the deformations: for each Borel-fixed ideals there is the list of all the simple and composed Borel rational deformations involving it. In the lower part of the window there is the description of the graph (vertices and edges) described with the language of the free software \href{http://www.graphviz.org/}{\texttt{Graphviz}}. Copying the code and pasting it in a file with \texttt{.dot} extension and compiling with the \texttt{neato} processor it is possible to get the picture of the Borel incidence graph (see for instance Figure \ref{fig:BorelIncidence} and Figure \ref{fig:incidenceBig}).

\begin{figure}[H]
\begin{center}
\includegraphics[width=10cm]{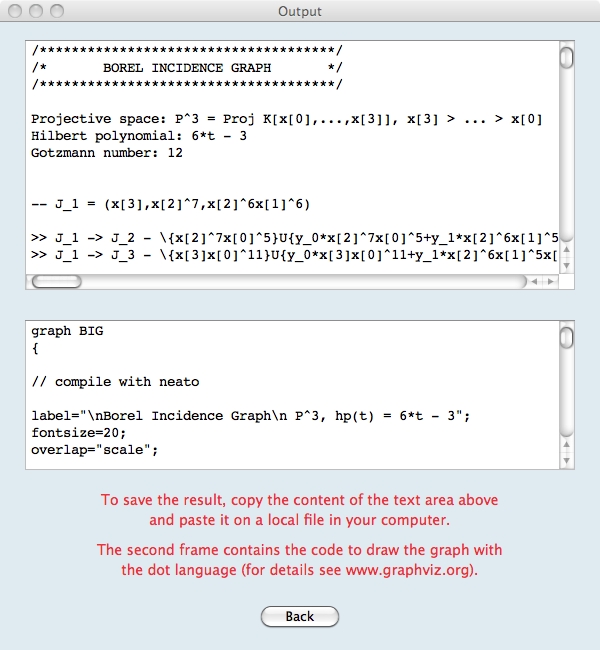}
\caption{\label{fig:IncidenceGraphOutput} The output window of the applet \texttt{Borel} \texttt{Incidence} \texttt{Graph}.}
\end{center}
\end{figure}

\section{$\sigma$-Borel degenerations}

The algorithms dealing with oriented Borel rational deformations are made available in the applets \texttt{Oriented} \texttt{Borel} \texttt{Rational} \texttt{Degeneration} (Figure \ref{fig:Degeneration}) and \texttt{Degeneration} \texttt{Graph} (Figure \ref{fig:DegenerationGraph}). 

\texttt{Oriented} \texttt{Borel} \texttt{Rational} \texttt{Degeneration} (Algorithm \ref{alg:deformation}) is available at \href{http://www.personalweb.unito.it/paolo.lella/HSC/TOdeformation.html}{\texttt{www.personalweb.unito.it/paolo.lella/HSC/TOdeformation.html}}\\ and requires three arguments: \textsf{Projective space} and \textsf{Borel-fixed ideal} are as in the applet \texttt{Borel} \texttt{Rational} \texttt{Deformations}. Moreover a term ordering is required: the default term order is $\DegLex$. To change it, we need to click on the button \lq\lq\textsf{Change}\rq\rq: in the dialog window that opens (Figure \ref{fig:termOrderWindow}) there are two options, \textsf{Graded} and \textsf{Reverse}, and the fields to fill with a sequence of rational coefficients determining the term ordering. The $\DegLex$ term ordering corresponds to \textsf{Graded: yes}, \textsf{Reverse: no}, \textsf{Weights: 1,0,\ldots,0}, whereas $\RevLex$ corresponds to \textsf{Graded: yes}, \textsf{Reverse: yes}, \textsf{Weights: 0,\ldots,0,-1}. To fix the term ordering described in \eqref{eq:termOrdering}, it suffices to choose \textsf{Graded: yes}, \textsf{Reverse: no} and to insert the vector $(\omega_n,\ldots,\omega_0)$.

\begin{figure}[H]
\begin{center}
\includegraphics[width=12cm]{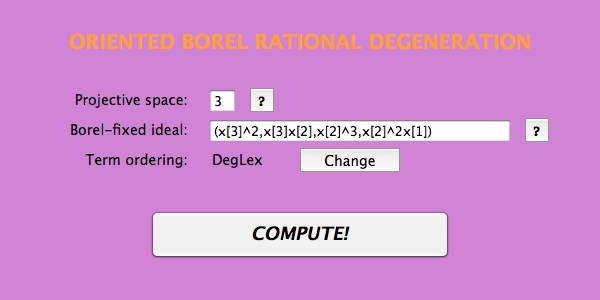}
\caption{\label{fig:Degeneration} The applet \texttt{Oriented }\texttt{Borel} \texttt{Rational} \texttt{Degeneration}.}
\end{center}
\end{figure}

\begin{figure}[H]
\begin{center}
\includegraphics[width=6cm]{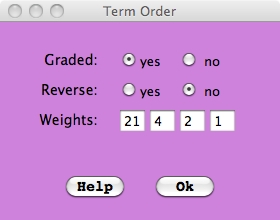}
\caption{\label{fig:termOrderWindow} The dialog window for changing the term ordering.}
\end{center}
\end{figure}

In the output window (Figure \ref{fig:DegenerationOutput}), there is the description of the degeneration computed in degree equal to the Gotzmann number of the Hilbert polynomial of the ideal, with the monomials exchanged specified and the Borel-fixed ideal obtained.

\begin{figure}[H]
\begin{center}
\includegraphics[width=10cm]{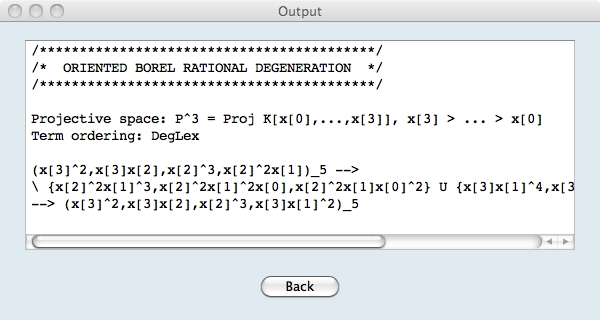}
\caption[The output window of the applet \texttt{Oriented} \texttt{Borel} \texttt{Rational}\newline {\texttt{Degeneration}}.]{\label{fig:DegenerationOutput} The output window of the applet \texttt{Oriented} \texttt{Borel} \texttt{Rational} \texttt{Degeneration}.}
\end{center}
\end{figure}

\begin{figure}[H]
\begin{center}
\includegraphics[width=12cm]{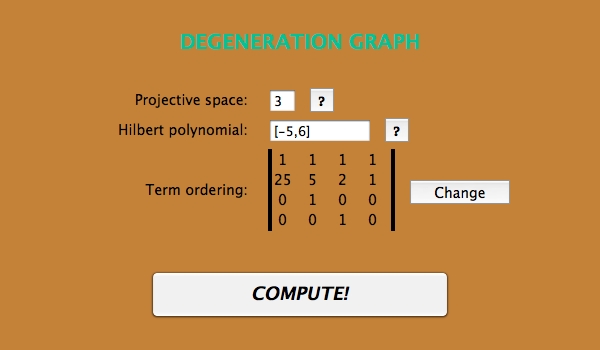}
\caption{\label{fig:DegenerationGraph} The applet \texttt{Degeneration} \texttt{Graph}.}
\end{center}
\end{figure}

\noindent The\hfill applet\hfill \texttt{Degeneration}\hfill \texttt{Graph}\hfill (available\hfill at\hfill \href{http://www.personalweb.unito.it/paolo.lella/HSC/deformationGraph.html}{\texttt{www.personalweb.unito.it/}\\ \texttt{paolo.lella/HSC/deformationGraph.html}})\hfill realizes\hfill Algorithm\hfill \ref{alg:deformationGraph}\hfill and\\ needs three arguments: \textsf{Projective space}, \textsf{Hilbert polynomial} and \textsf{Term ordering}. As for the applet computing the Borel incidence graph, its output window (Figure \ref{fig:DegenerationGraphOutput}) splits in two parts: the upper one contains the explicit description of the Borel-fixed ideals defining points on the chosen Hilbert polynomial with the relative Borel degeneration prescribed by the fixed term ordering, while the lower one contains the description of the direct graph again with the \href{http://www.graphviz.org/}{\texttt{Graphviz}} code. Copying the code and pasting it in a file with \texttt{.dot} extension and compiling with the \texttt{dot} processor it is possible to get the picture of the forest representing the degeneration graph (see for instance Figure \ref{fig:BorelGraph} and Figure \ref{fig:heightGraphs}).

\begin{figure}[H]
\begin{center}
\includegraphics[width=9.9cm]{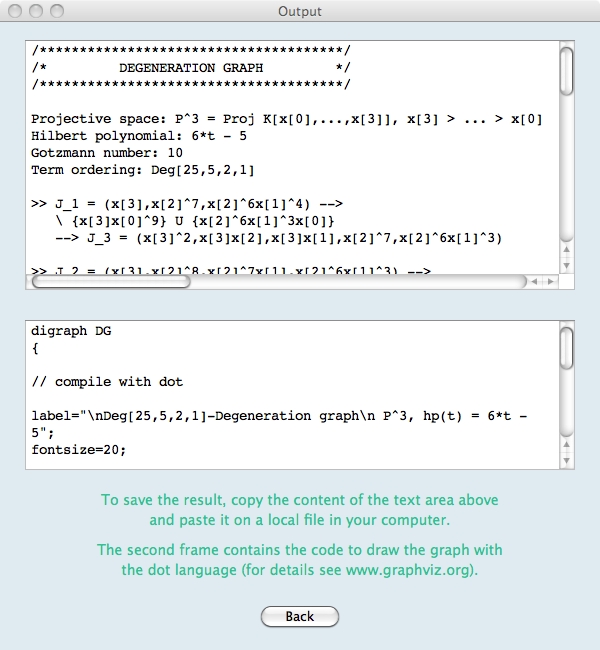}
\caption{\label{fig:DegenerationGraphOutput} The output window of the applet \texttt{Degeneration} \texttt{Graph}.}
\end{center}
\end{figure}

\chapter{The \textit{Macaulay2} package \texttt{MarkedSchemes}}\label{ch:markedSchemes}

This\hfill chapter\hfill is\hfill supposed\hfill to\hfill be\hfill a\hfill handbook\hfill for\hfill the\hfill \textit{Macaulay2}\hfill \cite{M2}\hfill package\\ \texttt{MarkedSchemes.m2}. We will introduce and explain the main functions of the package that allow to compute affine schemes associated to marked families and Gr\"obner strata.

\section{Basic features}
Firstly there are some basic methods to manipulate monomials (in the context of Borel-fixed ideals).
\begin{itemize}
\item \texttt{moveUp} and \texttt{moveDown} implement the elementary moves and require two arguments: a monomial and the index of the elementary move.

\begin{code}
\begin{verbatim}
moveUP = method(TypicalValue => RingElement)
--  INPUT: m monomial
--         i variable index
-- OUTPUT: the monomial m*(ring(m)_(i-1)/ring(m)_i);
--  ERROR: if m is not a monomial 
--         if m is a constant
--         if i < 1 or i > numgens(ring(m))-1
--         if the i-th variable does not divide m

moveDOWN = method(TypicalValue => RingElement)
--  INPUT: m monomial
--         i variable index
-- OUTPUT: the monomial m*(ring(m)_(i+1)/ring(m)_i);
--  ERROR: if m is not a monomial 
--         if m is a constant
--         if i < 0 or i >= numgens(ring(m))-1
--         if the i-th variable does not divide m
\end{verbatim}
\end{code}
\item \texttt{minimum} and \texttt{maximum} require as argument a monomial and return the index of the minimum/maximum variable dividing the monomial.

\begin{code}
\begin{verbatim}
minimum = method(TypicalValue => ZZ)
--  INPUT: m monomial
-- OUTPUT: the index of the smallest variable dividing m
--  ERROR: if m is not a monomial
--         if m is a constant

maximum = method(TypicalValue => ZZ)
--  INPUT: m monomial
-- OUTPUT: the index of the greatest variable dividing m
--  ERROR: if m is not a monomial
--         if m is a constant
\end{verbatim}
\end{code}
\item \texttt{canonicalDecomposition} is a method for computing the canonical decomposition of a monomial w.r.t. a Borel-fixed ideal containing it. It returns a sequence with two entries, the first one is the unique generator of the ideal giving the decomposition.

\begin{code}
\begin{verbatim}
canonicalDecomposition = method(TypicalValue => Sequence)
--  INPUT: m monomial 
--         J Borel-fixed ideal
-- OUTPUT: a sequence of two elements that represent the 
--           canonical decomposition of m over J
--  ERROR: if m is not a monomial 
--         if m is constant
--         if J is not a monomial ideal
--         if J is not a Borel-fixed ideal
--         if m does not belong to J
\end{verbatim}
\end{code}
\end{itemize}

\begin{example}
With this example, we want to point out that \textit{Macaulay2} does not permit to work with decreasing indexed variables, so we will have to consider $x_0 > \ldots > x_n$. 

\begin{code}
\begin{verbatim}
Macaulay2, version 1.4
i1 : loadPackage "MarkedSchemes";
i2 : R := QQ[x_0..x_3];
i3 : m = x_0*x_1
o3 = x x
      0 1
o3 : QQ[x , x , x , x ]
         0   1   2   3
i4 : mUP = moveUP(m,1)
      2
o4 = x
      0
o4 : QQ[x , x , x , x ]
         0   1   2   3
i5 : mDOWN = moveDOWN(m,1)
o5 = x x
      0 2
o5 : QQ[x , x , x , x ]
         0   1   2   3
i6 : minimum mDoWN
o6 = 2
i7 : maximum m
o7 = 0
i8 : J = ideal(mUP,m,mDoWN,x_1^5,x_1^4*x_2)
              2               5   4
o8 = ideal (x , x x , x x , x , x x )
              0   0 1   0 2   1   1 2
o8 : Ideal of QQ[x , x , x , x ]
                   0   1   2   3
i9 : f = x_0*x_1^5*x_2
         5
o9 = x x x
       0 1 2
o9 : QQ[x , x , x , x ]
          0   1   2   3
i10 : canonicalDecomposition (f,J)
              4
o10 = (x x , x x )
        0 1   1 2
o10 : Sequence
\end{verbatim}
\end{code}
\end{example}

\section{Marked families and Gr\"obner strata}

The basic method implementing Algorithm \ref{alg:markedFamily2} is called \texttt{markedScheme}.

\begin{code}
\begin{verbatim}
markedScheme = method(TypicalValue => Sequence, 
                      Options => {MonomialOrder => GRevLex,
                                  Segment => false,
                                  DescribeFamily => false,
                                  EliminateParameters => false});
--   INPUT: J Borel-fixed saturated ideal
--          s degree in which J has to be truncated
-- OPTIONS: Segment, if true there exists a term ordering for which the  
--                   ideal is a gen-segment ideal. 
--          MonomialOrder, if the option Segment is true, the option  
--                         MonomialOrder contains the term ordering for 
--                         which the ideal is a gen-segment ideal.
--          DescribeFamily, if true the function returns the complete 
--                          list of polynomial generators with coefficients
--                          in the parameters.
--          EliminateParameters, if true and Segment is true some 
--                               parameters willbe eliminated.
--  OUTPUT: sequence containting
--             - #0 the ring of parameters
--             - #1 the ideal of the marked scheme
--             - #2 the ring of the parameters and the starting variables
--             - #3 the polynomial generators of the family 
--                  (if opts.DescribeFamily all the polynomials, 
--                  otherwise only the superminimal generators)
--   ERROR: if J is not a monomial ideal
--          if J is not a Borel-fixed ideal
--          if s < 0
-- WARNING: if J is not saturated, the ideal will be saturated.
\end{verbatim}
\end{code}

With the default choices for the options of the method, the scheme of the $J$-marked family is computed avoiding any term ordering. If we are interested in the open\hfill subset\hfill of\hfill the\hfill Hilbert\hfill scheme\hfill defined\hfill by\hfill $J$,\hfill there\hfill is\hfill the\hfill method\\ \texttt{openSubsetHilb} that computes the maximal degree $\rho$ of a monomial generator of $J$ divided by the last but one variable and then calls $\texttt{markedScheme}$ with arguments $J$ and $\rho-1$.

\begin{example}\label{ex:appC1}
Let us consider the ideal $(x_0^2,x_0x_1^2,x_1^5) \subset \K[x_0,x_1,x_2],\ x_0 > x_1 > x_2$ and its truncation in degree 3.

\begin{code}
\begin{verbatim}
Macaulay2, version 1.4
i1 : loadPackage "MarkedSchemes";
i2 : R := QQ[x_0..x_2];
i3 : J = ideal(x_0^2,x_0*x_1^2,x_1^5)
             2     2   5
o3 = ideal (x , x x , x )
             0   0 1   1
o3 : Ideal of QQ[x , x , x ]
                  0   1   2
i4 : Mf = markedScheme(J,2);
     -- Computation of normal forms in progress...
     -- Completed in .073413
     -- Computation of S-polynomials in progress...
     -- Completed in .106108
i5 : Mf#0
o5 = QQ[p , p , p , p , p , p , p , p , p , p , p  , p  , p  , p  , p  , p  ,
         0   1   2   3   4   5   6   7   8   9   10   11   12   13   14   15 
     ---------------------------------------------------------------------------
         p  , p  ]
          16  17
o5 : PolynomialRing
i6 : Mf#1
                       2           2                                        2   
o6 = ideal (- p p   + p   - p , p p   - p p  + p p  - p p   + 2p p   - p , p p  
               4 10    10    3   9 10    4 8    3 9    2 10     8 10    1   9 10
     --------------------------------------------------------------------------- 
                                                          2                     
     + p p  + p p   + p , p p p   - p p  + p p  - p p  + p  + p p  + p p   - p ,
        8 9    7 10    6   8 9 10    4 6    3 7    2 8    8    1 9    6 10    0 
     ---------------------------------------------------------------------------
                                                                             
     p p p   + p p  + p , p p p   - p p  - p p  + p p  + p p  + p p  + p p  ,
      7 9 10    7 8    5   6 9 10    4 5    2 6    1 7    6 8    0 9    5 10 
     ---------------------------------------------------------------------------
                                      2                            3     
     p p p   - p p  + p p  + p p , - p p   - p p  - p p   - p , - p p   +
      5 9 10    2 5    0 7    5 8     9 10    8 9    7 10    6     9 10  
     ---------------------------------------------------------------------------
      2            2                          2                                
     p p  p   - p p  - 2p p p   + p p  p   - p  p   + p p p   + p p  p   - p p 
      9 10 17    8 9     7 9 10    4 10 15    10 15    8 9 17    7 10 17    7 8
     ---------------------------------------------------------------------------
                                     4    3          2                        
     - p p  + p p   + p p   - p , - p  + p p   - 3p p  + p p p   - 2p p  p   +
        6 9    3 15    6 17    5     9    9 17     7 9    4 9 15     9 10 15  
     ---------------------------------------------------------------------------
      2                  2                                                   
     p p   + 2p p p   - p  + p p   - p  p   + p p   + p p   - p p   + p p   +
      9 16     7 9 17    7    4 13    10 13    9 14    2 15    8 15    7 16  
     ---------------------------------------------------------------------------
               3      2                   2                                  
     p  , - p p  + p p p   - 2p p p  - p p  - p p  p   + p p p   - p p  p   +
      12     8 9    8 9 17     7 8 9    6 9    9 10 14    4 8 15    8 10 15  
     ---------------------------------------------------------------------------
                                                                               
     p p p   + p p p   + p p p   - p p  - p p  - p  p   + p p   + p p   + p p  
      8 9 16    7 8 17    6 9 17    6 7    5 9    10 12    3 13    1 15    6 16
     ---------------------------------------------------------------------------
                   3      2        2                                            
     + p p  , - p p  + p p p   - 2p p  - p p  p   + p p p   - p p  p   + p p p  
        5 17     7 9    7 9 17     7 9    9 10 13    4 7 15    7 10 15    7 9 16
     ---------------------------------------------------------------------------
        2                                        3      2                   2  
     + p p   + p p   - p p   + p p   + p  , - p p  + p p p   - 2p p p  - p p  -
        7 17    2 13    8 13    7 14    11     6 9    6 9 17     6 7 9    5 9  
     ---------------------------------------------------------------------------
                                                                                
     p p  p   + p p p   - p p  p   + p p p   + p p p   + p p p   - p p  - p  p  
      9 10 12    4 6 15    6 10 15    6 9 16    6 7 17    5 9 17    5 7    10 11
     ---------------------------------------------------------------------------
                                                   3      2
     - p p   + p p   + p p   + p p   + p p  , - p p  + p p p   - 2p p p  -
        8 12    1 13    6 14    0 15    5 16     5 9    5 9 17     5 7 9  
     ---------------------------------------------------------------------------
     p p  p   + p p p   - p p  p   + p p p   + p p p   - p p   + p p   + p p  )
      9 10 11    4 5 15    5 10 15    5 9 16    5 7 17    8 11    0 13    5 14

o6 : Ideal of QQ[p , p , p , p , p , p , p , p , p , p , p  , p  , p  , p  ,
                  0   1   2   3   4   5   6   7   8   9   10   11   12   13 
     ---------------------------------------------------------------------------
                  p  , p  , p  , p  ]
                   14   15   16   17
i7 : Mf#2

o7 = QQ[p , p , p , p , p , p , p , p , p , p , p  , p  , p  , p  , p  , p  ,
         0   1   2   3   4   5   6   7   8   9   10   11   12   13   14   15 
         p  , p  ][x , x , x ]
          16   17   0   1   2
o7 : PolynomialRing
i8 : Mf#3
       2               2                        2     2       3             
o8 = {x  - p x x  - p x  - p x x  - p x x  - p x , x x  - p  x  - p x x x  -
       0    4 0 1    3 1    2 0 2    1 1 2    0 2   0 1    10 1    9 0 1 2  
     ---------------------------------------------------------------------------
        2          2        2      3   5       4         3 2           3  
     p x x  - p x x  - p x x  - p x , x  - p  x x  - p  x x  - p  x x x  -
      8 1 2    7 0 2    6 1 2    5 2   1    17 1 2    16 1 2    15 0 1 2  
     ---------------------------------------------------------------------------
         2 3         4         4       5
     p  x x  - p  x x  - p  x x  - p  x }
      14 1 2    13 0 2    12 1 2    11 2
o8 : List
\end{verbatim}
\end{code}
Since the highest degree of a generator divisible by $x_1$ is 5, the optimal degree to compute the open subset of \gls{Hilb7P2} defined by $J$ is 4. Note that there are more parameters because there are more monomials in the tails of the superminimal generators.

\begin{code}
\begin{verbatim}
i9 : HJ = openSubsetHilb (J);
     -- Computation of normal forms in progress...
     -- Completed in .363012
     -- Computation of S-polynomials in progress...
     -- Completed in .152566
i10 : HJ#0
o10 = QQ[p , p , p , p , p , p , p , p , p , p , p  , p  , p  , p  , p  , p  ,
          0   1   2   3   4   5   6   7   8   9   10   11   12   13   14   15 
      --------------------------------------------------------------------------
          p  , p  , p  , p  , p  ]
           16   17   18   19   20
o10 : PolynomialRing
i11 : HJ#3
            4      3      2 2          2      2 2        3        3      4   
o11 = {- p x  - p x x  + x x  - p x x x  - p x x  - p x x  - p x x  - p x , -
          6 1    5 1 2    0 2    4 0 1 2    3 1 2    2 0 2    1 1 2    0 2   
      --------------------------------------------------------------------------
          4      2         3             2       2 2        3        3      4 
      p  x  + x x x  - p  x x  - p  x x x  - p  x x  - p x x  - p x x  - p x ,
       13 1    0 1 2    12 1 2    11 0 1 2    10 1 2    9 0 2    8 1 2    7 2 
      --------------------------------------------------------------------------
       5       4         3 2           3       2 3         4         4       5
      x  - p  x x  - p  x x  - p  x x x  - p  x x  - p  x x  - p  x x  - p  x }
       1    20 1 2    19 1 2    18 0 1 2    17 1 2    16 0 2    15 1 2    14 2
o11 : List
\end{verbatim}
\end{code}
\end{example}

There are also the methods computing the dimension of the tangent space at the origin.

\begin{code}
\begin{verbatim}
EDmarkedScheme = method(TypicalValue => ZZ, 
                        Options => {MonomialOrder => GRevLex,
                        Segment => false});
--   INPUT: J Borel-fixed saturated ideal
--          s degree in which J has to be truncated
-- OPTIONS: Segment, if true there exists a term ordering for which the  
--                   ideal is a gen-segment ideal. 
--          MonomialOrder, if the option Segment is true, the option  
--                         MonomialOrder contains the term ordering for 
--                         which the ideal is a gen-segment ideal.
--  OUTPUT: the number of parameters necessary to describe the marked family
--   ERROR: if J is not a monomial ideal
--          if J is not a Borel-fixed ideal
--          if s < 0
-- WARNING: if J is not saturated, the ideal will be saturated.

EDopenSubsetHilb = method(TypicalValue => ZZ, 
                          Options => {MonomialOrder => GRevLex,
                          Segment => false});
--   INPUT: J Borel-fixed saturated ideal
-- OPTIONS: Segment, if true there exists a term ordering for which the  
--                   ideal is a gen-segment ideal. 
--          MonomialOrder, if the option Segment is true, the option  
--                         MonomialOrder contains the term ordering for 
--                         which the ideal is a gen-segment ideal.
--  OUTPUT: the number of parameters necessary to describe the marked scheme
--   ERROR: if J is not a monomial ideal
--          if J is not a Borel-fixed ideal
--          if s < 0
-- WARNING: if J is not saturated, the ideal will be saturated.
\end{verbatim}
\end{code}

Considering Example \ref{ex:appC1}, we have that

\begin{code}
\begin{verbatim}
i12 : EDmarkedScheme(J,3)
o12 = 19
i13 : EDopenSubsetHilb(J)
o13 = 21
\end{verbatim}
\end{code}

\bigskip

To compute the Gr\"obner stratum of a gen-segment ideal, i.e. to add to the marked family the structure of homogeneous variety w.r.t. a positive grading, we have to use the options putting \texttt{Segment = true}, giving the term ordering for which the ideal is a gen-segment ideal and requiring the elimination of the parameters (\texttt{EliminateParameters = true}).

\begin{example}\label{ex:Hilb4tP3_M2}
Let\hfill us\hfill consider\hfill the\hfill Borel-fixed\hfill ideal\hfill $J = (x_0^2,x_0 x_1,x_0x_2^2,x_1^4)$\hfill in\\ $\K[x_0,x_1,x_2,x_3]\ (x_0 > x_1 > x_2 > x_3)$, that defines a point on the Hilbert scheme \gls{Hilb4tP3} (cf. Example \ref{ex:Hilb4tP3}).

\begin{code}
\begin{verbatim}
Macaulay2, version 1.4
i1 : loadPackage "MarkedSchemes";
i2 : R := QQ[x_0..x_3];
i3 : J = ideal(x_0^2,x_0*x_1,x_0*x_2^2,x_1^4)
             2           2   4
o3 = ideal (x , x x , x x , x )
             0   0 1   0 2   1
o3 : Ideal of QQ[x , x , x , x ]
                  0   1   2   3
i4 : EDopenSubsetHilb (J)
o4 = 44
i5 : EDopenSubsetHilb (J,Segment=>true,MonomialOrder=>{Weights=>{7,3,2,1}})
o5 = 24
i6 : St = openSubsetHilb (J,Segment=>true,
                            MonomialOrder=>{Weights=>{7,3,2,1}},
                            EliminateParameters=>true);
     -- Computation of normal forms in progress...
     -- Completed in .116961
     -- Computation of S-polynomials in progress...
     -- Completed in .395245
     -- Elimination of parameters in progress...
     -- Completed in 5.42916
i7 : A := St#0
o7 = QQ[p  , p  , p  , p  , p  , p  , p  , p  , p  , p  , p  , p  , p  , p  ,          
         12   13   16   17   19   20   23   24   25   26   28   30   31   32          
     --------------------------------------------------------------------------
        p  , p  , p  , p  , p  , p  , p  , p  , p  , p  ]
         34   35   36   37   38   39   40   41   42   43
o7 : PolynomialRing
i8 : numgens A
o8 = 24
i9 : idealSt = St#1;
i10 : gensIdealSt := first entries gens idealSt;
i11 : #gensIdealSt
o11 = 40
i12 : (factor(gensIdealSt#0))#0
o12 = p
       39
i13 : (factor(gensIdealSt#1))#0
o13 = p
       39
i14 : (factor(gensIdealSt#2))#0
o14 = p
       39
i15 : dec1 := ideal(p_39);
i16 : RScomponent := Spec (A/dec1);
i17 : dim RScomponent
o17 = 23
i18 : dec2 := saturate (idealSt,p_39);
i19 : VAcomponent := Spec (A/dec2);
i20 : dim VAcomponent
o20 = 16
\end{verbatim}
\end{code}
\end{example}

\backmatter

\renewcommand{\glossarypreamble}{The list of all Hilbert schemes discussed in the examples with a brief description of the geometric objects parametrized.\par}
 \renewcommand{\entryname}{Hilbert scheme}
\renewcommand{\descriptionname}{Geometric objects parametrized}
\renewcommand{\pagelistname}{pages}

\printglossary[type=hilbSchemes,title={Catalog of Hilbert schemes},toctitle={Catalog of Hilbert schemes},style=long3colheaderBIS]

\renewcommand{\glossarypreamble}{}
\renewcommand{\entryname}{Symbol}
\renewcommand{\descriptionname}{Typical usage or definition}
\renewcommand{\pagelistname}{page}

\printglossary[title={Symbols and Notation},toctitle={Symbols and Notation}]

\clearemptydoublepage

\bibliographystyle{amsplain}

\providecommand{\bysame}{\leavevmode\hbox to3em{\hrulefill}\thinspace}
\providecommand{\MR}{\relax\ifhmode\unskip\space\fi MR }
\providecommand{\MRhref}[2]{%
  \href{http://www.ams.org/mathscinet-getitem?mr=#1}{#2}
}
\providecommand{\href}[2]{#2}
\begin{thebibliography}{100}

\bibitem{AitAmrane}
Samir A{\"{\i}}t~Amrane, \emph{Sur le sch{\'e}ma de {H}ilbert des courbes de
  degr{\'e} {$d$} et genre {$(d-3)(d-4)/2$} de {$\bold P^3_k$}}, C. R. Acad.
  Sci. Paris S{\'e}r. I Math. \textbf{326} (1998), no.~7, 851--856.

\bibitem{AltmannSturmfels}
Klaus Altmann and Bernd Sturmfels, \emph{The graph of monomial ideals}, J. Pure
  Appl. Algebra \textbf{201} (2005), no.~1-3, 250--263.

\bibitem{Artin}
Michael Artin, \emph{{L}ectures on deformations of singularities}, {T}ata
  {I}nstitute on {F}undamental {R}esearch, Bombay, 1976.

\bibitem{AtiyahMacdonald}
Michael~F. Atiyah and Ian~G. MacDonald, \emph{Introduction to commutative
  algebra}, Addison-Wesley Publishing Co., Reading, Mass.-London-Don Mills,
  Ont., 1969.

\bibitem{BallicoBolondiMigliore}
Edoardo Ballico, Giorgio Bolondi, and Juan~Carlos Migliore, \emph{The
  {L}azarsfeld-{R}ao problem for liaison classes of two-codimensional
  subschemes of {${\bf P}^n$}}, Amer. J. Math. \textbf{113} (1991), no.~1,
  117--128.

\bibitem{BarnabeiBriniRota}
Marilena Barnabei, Andrea Brini, and Gian-Carlo Rota, \emph{On the exterior
  calculus of invariant theory}, J. Algebra \textbf{96} (1985), no.~1,
  120--160.

\bibitem{BayerThesis}
David Bayer, \emph{{T}he division algorithm and the {H}ilbert scheme}, Ph.D.
  thesis, Harvard University, 1982, p.~163.

\bibitem{BayerMorrison}
David Bayer and Ian Morrison, \emph{Standard bases and geometric invariant
  theory. {I}. {I}nitial ideals and state polytopes}, J. Symbolic Comput.
  \textbf{6} (1988), no.~2-3, 209--217. Computational aspects of commutative
  algebra.

\bibitem{BayerMumford}
David Bayer and David Mumford, \emph{What can be computed in algebraic
  geometry?}, ArXiv e-prints (1992), Available at
  \href{http://arxiv.org/abs/alg-geom/9304003}{\texttt{arxiv.org/abs/alg-geom/9304003}}.

\bibitem{BayerStillmanREG}
David Bayer and Michael Stillman, \emph{A criterion for detecting
  {$m$}-regularity}, Invent. Math. \textbf{87} (1987), no.~1, 1--11.

\bibitem{BayerStillmanGIN}
\bysame, \emph{A theorem on refining division orders by the reverse
  lexicographic order}, Duke Math. J. \textbf{55} (1987), no.~2, 321--328.

\bibitem{BCLR}
Cristina Bertone, Francesca Cioffi, Paolo Lella, and Margherita Roggero,
  \emph{Upgraded methods for the effective computation of marked schemes on a
  strongly stable ideal}, ArXiv e-prints (2011). Available at
  \href{http://arxiv.org/abs/1110.0698}{\texttt{arxiv.org/abs/1110.0698}}.

\bibitem{BLR}
Cristina Bertone, Paolo Lella, and Margherita Roggero, \emph{Borel open
  covering of {H}ilbert schemes}, ArXiv e-prints (2011). Available at
  \href{http://arxiv.org/abs/0909.2184}{\texttt{arxiv.org/abs/0909.2184}}.

\bibitem{Bertram}
Aaron Bertram, \emph{Construction of the {H}ilbert scheme}, Available at
  \href{http://www.math.utah.edu/~bertram/courses/hilbert/}{\texttt{www.math.utah.edu/\textasciitilde
  bertram/courses/hilbert/}}, 1999. Notes for a course at University of Utah.

\bibitem{BigattiRobbiano}
A.~M. Bigatti and L.~Robbiano, \emph{Borel sets and sectional matrices}, Ann.
  Comb. \textbf{1} (1997), no.~3, 197--213.

\bibitem{BourbakiCA}
Nicolas Bourbaki, \emph{Commutative algebra. {C}hapters 1--7}, Elements of
  Mathematics (Berlin), Springer-Verlag, Berlin, 1989. Translated from the
  French. Reprint of the 1972 edition.

\bibitem{BLMR}
Jerome Brachat, Paolo Lella, Bernard Mourrain, and Margherita Roggero,
  \emph{Low degree equations defining the {H}ilbert scheme}, ArXiv e-prints
  (2011). Available at
  \href{http://arxiv.org/abs/1104.2007}{\texttt{arxiv.org/abs/1104.2007}}.

\bibitem{Burden}
Richard~L. Burden, J.~Douglas Faires, and Albert~C. Reynolds, \emph{Numerical
  analysis}, Prindle, Weber \& Schmidt, Boston, Mass., 1978.

\bibitem{CarraFerro}
Giuseppa Carr{{\`a}}~Ferro, \emph{Gr{\"o}bner bases and {H}ilbert schemes.
  {I}}, J. Symbolic Comput. \textbf{6} (1988), no.~2-3, 219--230. Computational
  aspects of commutative algebra.

\bibitem{CilibertoSernesi}
Ciro Ciliberto and Edoardo Sernesi, \emph{Families of varieties and the
  {H}ilbert scheme}, Lectures on {R}iemann surfaces ({T}rieste, 1987), World
  Sci. Publ., Teaneck, NJ, 1989, pp.~428--499.

\bibitem{CLMR}
Francesca Cioffi, Paolo Lella, Maria~Grazia Marinari, and Margherita Roggero,
  \emph{Segments and {H}ilbert schemes of points}, Discrete Math. \textbf{311}
  (2011), no.~20, 2238--2252.

\bibitem{CioffiRoggero}
Francesca Cioffi and Margherita Roggero, \emph{Flat families by strongly stable
  ideals and a generalization of {G}r{\"o}bner bases}, J. Symbolic Comput.
  \textbf{46} (2011), no.~9, 1070--1084.

\bibitem{ConcaSidman}
Aldo Conca and Jessica Sidman, \emph{Generic initial ideals of points and
  curves}, J. Symbolic Comput. \textbf{40} (2005), no.~3, 1023--1038.

\bibitem{CLOiva}
David Cox, John Little, and Donal O'Shea, \emph{Ideals, varieties, and
  algorithms}, third ed., Undergraduate Texts in Mathematics, Springer, New
  York, 2007. An introduction to computational algebraic geometry and
  commutative algebra.

\bibitem{CLOuag}
David~A. Cox, John Little, and Donal O'Shea, \emph{Using algebraic geometry},
  second ed., Graduate Texts in Mathematics, vol. 185, Springer, New York,
  2005.

\bibitem{Deery}
Todd Deery, \emph{Rev-lex segment ideals and minimal {B}etti numbers}, The
  {C}urves {S}eminar at {Q}ueen's, {V}ol.\ {X} ({K}ingston, {ON}, 1995),
  Queen's Papers in Pure and Appl. Math., vol. 102, Queen's Univ., Kingston,
  ON, 1996, pp.~193--219.

\bibitem{Eisenbud}
David Eisenbud, \emph{Commutative algebra}, Graduate Texts in Mathematics, vol.
  150, Springer-Verlag, New York, 1995. With a view toward algebraic geometry.

\bibitem{EisenbudGoto}
David Eisenbud and Shiro Goto, \emph{Linear free resolutions and minimal
  multiplicity}, J. Algebra \textbf{88} (1984), no.~1, 89--133.

\bibitem{EisenbudHarris}
David Eisenbud and Joe Harris, \emph{The geometry of schemes}, Graduate Texts
  in Mathematics, vol. 197, Springer-Verlag, New York, 2000.

\bibitem{EliahouKervaire}
Shalom Eliahou and Michel Kervaire, \emph{Minimal resolutions of some monomial
  ideals}, J. Algebra \textbf{129} (1990), no.~1, 1--25.

\bibitem{Ellia}
Philippe Ellia, \emph{D'autres composantes non r{\'e}duites de {${\rm
  Hilb}\,{\bf P}^3$}}, Math. Ann. \textbf{277} (1987), no.~3, 433--446.

\bibitem{FerrareseRoggero}
Giorgio Ferrarese and Margherita Roggero, \emph{Homogeneous varieties for
  {H}ilbert schemes}, Int. J. Algebra \textbf{3} (2009), no.~9-12, 547--557.

\bibitem{Galligo}
Andr{{\'e}} Galligo, \emph{\`{A} propos du th{\'e}or{\`e}me de-pr{\'e}paration
  de {W}eierstrass}, Fonctions de plusieurs variables complexes ({S}{\'e}m.
  {F}ran\c cois {N}orguet, octobre 1970--d{\'e}cembre 1973; {\`a} la
  m{\'e}moire d'{A}ndr{\'e} {M}artineau), Lecture Notes in Math., vol. 409,
  Springer, Berlin, 1974. Th{{\`e}}se de 3{{\`e}}me cycle soutenue le 16 mai
  1973 {{\`a}} l'Institut de Math{{\'e}}matique et Sciences Physiques de
  l'Universit{{\'e}} de Nice, pp.~543--579.

\bibitem{Gotzmann}
Gerd Gotzmann, \emph{Eine {B}edingung f{\"u}r die {F}lachheit und das
  {H}ilbertpolynom eines graduierten {R}inges}, Math. Z. \textbf{158} (1978),
  no.~1, 61--70.

\bibitem{GotzmannHilb4tP3}
Gerd Gotzmann, \emph{The irreducible components of $\text{{H}ilb}^{4n}(p^3)$},
  ArXiv e-prints (2008). Available at
  \href{http://arxiv.org/abs/0811.3160}{\texttt{arxiv.org/abs/0811.3160}}.

\bibitem{ConcreteMathematics}
Ronald~L. Graham, Donald~E. Knuth, and Oren Patashnik, \emph{Concrete
  mathematics}, second ed., Addison-Wesley Publishing Company, Reading, MA,
  1994. A foundation for computer science.

\bibitem{M2}
Daniel~R. Grayson and Michael~E. Stillman, \emph{Macaulay2, a software system
  for research in algebraic geometry}. Available at
  \href{http://www.math.uiuc.edu/Macaulay2/}{\texttt{www.math.uiuc.edu/Macaulay2/}}.

\bibitem{GreenGIN}
Mark~L. Green, \emph{Generic initial ideals}, Six lectures on commutative
  algebra, Mod. Birkh{\"a}user Class., Birkh{\"a}user Verlag, Basel, 2010,
  pp.~119--186.

\bibitem{GriffithsHarris}
Phillip Griffiths and Joseph Harris, \emph{Principles of algebraic geometry},
  Wiley-Interscience [John Wiley \& Sons], New York, 1978, Pure and Applied
  Mathematics.

\bibitem{GrothendieckFGA}
Alexander Grothendieck, \emph{Techniques de construction et th{\'e}or{\`e}mes
  d'existence en g{\'e}om{\'e}trie alg{\'e}brique. {IV}. {L}es sch{\'e}mas de
  {H}ilbert}, S{\'e}minaire {B}ourbaki, {V}ol.\ 6, Soc. Math. France, Paris,
  1995, pp.~Exp.\ No.\ 221, 249--276.

\bibitem{HaimanSturmfels}
Mark Haiman and Bernd Sturmfels, \emph{Multigraded {H}ilbert schemes}, J.
  Algebraic Geom. \textbf{13} (2004), no.~4, 725--769.

\bibitem{HartshorneThesis}
Robin Hartshorne, \emph{Connectedness of the {H}ilbert scheme}, Inst. Hautes
  {\'E}tudes Sci. Publ. Math. (1966), no.~29, 5--48.

\bibitem{Hartshorne}
\bysame, \emph{Algebraic geometry}, Springer-Verlag, New York, 1977, Graduate
  Texts in Mathematics, No. 52.

\bibitem{HartshorneGD}
\bysame, \emph{Generalized divisors on {G}orenstein schemes}, $K$-Theory
  \textbf{8} (1994), no.~3, 287--339.

\bibitem{HartshorneConnectedness}
\bysame, \emph{On the connectedness of the {H}ilbert scheme of curves in
  {${\Bbb P}^3$}}, Comm. Algebra \textbf{28} (2000), no.~12, 6059--6077.
  Special issue in honor of Robin Hartshorne.

\bibitem{HartshorneToulouse}
\bysame, \emph{On {R}ao's theorems and the {L}azarsfeld-{R}ao property}, Ann.
  Fac. Sci. Toulouse Math. (6) \textbf{12} (2003), no.~3, 375--393.

\bibitem{HartshorneConnectedness2}
\bysame, \emph{Questions of connectedness of the {H}ilbert scheme of curves in
  {$\Bbb P^3$}}, Algebra, arithmetic and geometry with applications ({W}est
  {L}afayette, {IN}, 2000), Springer, Berlin, 2004, pp.~487--495.

\bibitem{HMDP}
Robin Hartshorne, Mireille Martin-Deschamps, and Daniel Perrin, \emph{Triades
  et familles de courbes gauches}, Math. Ann. \textbf{315} (1999), no.~3,
  397--468.

\bibitem{HartshorneSchlesinger2H}
Robin Hartshorne and Enrico Schlesinger, \emph{Curves in the double plane},
  Comm. Algebra \textbf{28} (2000), no.~12, 5655--5676. Special issue in honor
  of Robin Hartshorne.

\bibitem{IarrobinoRed}
Anthony~A. Iarrobino, \emph{Reducibility of the families of {$0$}-dimensional
  schemes on a variety}, Invent. Math. \textbf{15} (1972), 72--77.

\bibitem{IarrobinoPunctual}
\bysame, \emph{Punctual {H}ilbert schemes}, Mem. Amer. Math. Soc. \textbf{10}
  (1977), no.~188, viii+112.

\bibitem{IarrobinoKanev}
Anthony~A. Iarrobino and Vassil Kanev, \emph{Power sums, {G}orenstein algebras,
  and determinantal loci}, Lecture Notes in Math., vol. 1721, Springer-Verlag,
  Berlin, 1999. Appendix C by Iarrobino and Steven L. Kleiman.

\bibitem{Kleppe}
Jan~O. Kleppe, \emph{Nonreduced components of the {H}ilbert scheme of smooth
  space curves}, Space curves ({R}occa di {P}apa, 1985), Lecture Notes in
  Math., vol. 1266, Springer, Berlin, 1987, pp.~181--207.

\bibitem{KreuzerRobbiano1}
Martin Kreuzer and Lorenzo Robbiano, \emph{Computational commutative algebra.
  1}, Springer-Verlag, Berlin, 2000.

\bibitem{KreuzerRobbiano2}
\bysame, \emph{Computational commutative algebra. 2}, Springer-Verlag, Berlin,
  2005.

\bibitem{LazarsfeldRao}
Robert Lazarsfeld and Prabhakar~A. Rao, \emph{Linkage of general curves of
  large degree}, Algebraic geometry---open problems ({R}avello, 1982), Lecture
  Notes in Math., vol. 997, Springer, Berlin, 1983, pp.~267--289.

\bibitem{LellaDeformations}
Paolo Lella, \emph{A network of rational curves on the {H}ilbert scheme}, ArXiv
  e-prints (2010). Available at
  \href{http://arxiv.org/abs/1006.5020}{\texttt{arxiv.org/abs/1006.5020}}.

\bibitem{LellaRoggero}
Paolo Lella and Margherita Roggero, \emph{{R}ational components of {H}ilbert
  schemes}, {R}end. {S}emin. {M}at. {U}niv. {P}adova \textbf{126} (2011),
  11--45.

\bibitem{LellaSchlesinger}
Paolo Lella and Enrico Schlesinger, \emph{The {H}ilbert schemes of locally
  {C}ohen-{M}acaulay curves in $\mathbb{P}^3$ may after all be connected},
  ArXiv e-prints (2011). Available at
  \href{http://arxiv.org/abs/1110.2611}{\texttt{arxiv.org/abs/1110.2611}}.

\bibitem{Macaulay}
F.~S. Macaulay, \emph{{Some properties of enumeration in the theory of modular
  systems.}}, Proc. London Math. Soc. \textbf{26} (1927), 531--555.

\bibitem{Mall}
Daniel Mall, \emph{Connectedness of {H}ilbert function strata and other
  connectedness results}, J. Pure Appl. Algebra \textbf{150} (2000), no.~2,
  175--205.

\bibitem{MarinariBUMI}
Maria~Grazia Marinari, \emph{On {B}orel ideals}, Boll. Unione Mat. Ital. Sez. B
  Artic. Ric. Mat. (8) \textbf{4} (2001), no.~1, 207--237.

\bibitem{MarinariRamella1999}
Maria~Grazia Marinari and Luciana Ramella, \emph{Some properties of {B}orel
  ideals}, J. Pure Appl. Algebra \textbf{139} (1999), no.~1-3, 183--200.
  Effective methods in algebraic geometry (Saint-Malo, 1998).

\bibitem{MarinariRamella2005}
\bysame, \emph{A characterization of stable and {B}orel ideals}, Appl. Algebra
  Engrg. Comm. Comput. \textbf{16} (2005), no.~1, 45--68.

\bibitem{MarinariRamella2006}
\bysame, \emph{Borel ideals in three variables}, Beitr{\"a}ge Algebra Geom.
  \textbf{47} (2006), no.~1, 195--209.

\bibitem{MDP}
Mireille Martin-Deschamps and Daniel Perrin, \emph{Sur la classification des
  courbes gauches}, Ast{\'e}risque \textbf{184-185} (1990), 208.

\bibitem{MDPbounds}
\bysame, \emph{Sur les bornes du module de {R}ao}, C. R. Acad. Sci. Paris
  S{\'e}r. I Math. \textbf{317} (1993), no.~12, 1159--1162.

\bibitem{MDPextremal}
\bysame, \emph{Le sch{\'e}ma de {H}ilbert des courbes gauches localement
  {C}ohen-{M}acaulay n'est (presque) jamais r{\'e}duit}, Ann. Sci. {\'E}cole
  Norm. Sup. (4) \textbf{29} (1996), no.~6, 757--785.

\bibitem{MatsumuraCA}
Hideyuki Matsumura, \emph{Commutative algebra}, second ed., Mathematics Lecture
  Note Series, vol.~56, Benjamin/Cummings Publishing Co., Inc., Reading, Mass.,
  1980.

\bibitem{MatsumuraCRT}
\bysame, \emph{Commutative ring theory}, second ed., Cambridge Studies in
  Advanced Mathematics, vol.~8, Cambridge University Press, Cambridge, 1989.
  Translated from the Japanese by M. Reid.

\bibitem{Migliore}
Juan~C. Migliore, \emph{Introduction to liaison theory and deficiency modules},
  Progress in Mathematics, vol. 165, Birkh{\"a}user Boston Inc., Boston, MA,
  1998.

\bibitem{MillerSturmfels}
Ezra Miller and Bernd Sturmfels, \emph{Combinatorial commutative algebra},
  Graduate Texts in Mathematics, vol. 227, Springer-Verlag, New York, 2005.

\bibitem{MoraI}
Teo Mora, \emph{Solving polynomial equation systems. {I}}, Encyclopedia of
  Mathematics and its Applications, vol.~88, Cambridge University Press,
  Cambridge, 2003. The Kronecker-Duval philosophy.

\bibitem{MoraII}
\bysame, \emph{Solving polynomial equation systems. {II}}, Encyclopedia of
  Mathematics and its Applications, vol.~99, Cambridge University Press,
  Cambridge, 2005. Macaulay's paradigm and Gr{{\"o}}bner technology.

\bibitem{MoraRobbiano}
Teo Mora and Lorenzo Robbiano, \emph{The {G}r{\"o}bner fan of an ideal}, J.
  Symbolic Comput. \textbf{6} (1988), no.~2-3, 183--208. Computational aspects
  of commutative algebra.

\bibitem{Mourrain}
Bernard Mourrain, \emph{A new criterion for normal form algorithms}, Applied
  algebra, algebraic algorithms and error-correcting codes ({H}onolulu, {HI},
  1999), Lecture Notes in Comput. Sci., vol. 1719, Springer, Berlin, 1999,
  pp.~430--443.

\bibitem{MumfordPathologies}
David Mumford, \emph{Further pathologies in algebraic geometry}, Amer. J. Math.
  \textbf{84} (1962), 642--648.

\bibitem{Mumford}
\bysame, \emph{Lectures on curves on an algebraic surface}, With a section by
  G. M. Bergman. Annals of Mathematics Studies, No. 59, Princeton University
  Press, Princeton, N.J., 1966.

\bibitem{NolletDegree3}
Scott Nollet, \emph{The {H}ilbert schemes of degree three curves}, Ann. Sci.
  {\'E}cole Norm. Sup. (4) \textbf{30} (1997), no.~3, 367--384.

\bibitem{NolletSubextremal}
\bysame, \emph{Subextremal curves}, Manuscripta Math. \textbf{94} (1997),
  no.~3, 303--317.

\bibitem{NolletSchlesingerDegree4}
Scott Nollet and Enrico Schlesinger, \emph{Hilbert schemes of degree four
  curves}, Compositio Math. \textbf{139} (2003), no.~2, 169--196.

\bibitem{NotariSpreafico}
Roberto Notari and Maria~Luisa Spreafico, \emph{A stratification of {H}ilbert
  schemes by initial ideals and applications}, Manuscripta Math. \textbf{101}
  (2000), no.~4, 429--448.

\bibitem{PeevaStillman}
Irena Peeva and Mike Stillman, \emph{Connectedness of {H}ilbert schemes}, J.
  Algebraic Geom. \textbf{14} (2005), no.~2, 193--211.

\bibitem{PerrinKoszul}
Daniel Perrin, \emph{Un pas vers la connexit{\'e} du sch{\'e}ma de {H}ilbert:
  les courbes de {K}oszul sont dans la composante des extr{\'e}males}, Collect.
  Math. \textbf{52} (2001), no.~3, 295--319.

\bibitem{Rao}
Prabhakar~A. Rao, \emph{Liaison among curves in {${\bf P}^{3}$}}, Invent. Math.
  \textbf{50} (1978/79), no.~3, 205--217.

\bibitem{Reeves}
Alyson~A. Reeves, \emph{The radius of the {H}ilbert scheme}, J. Algebraic Geom.
  \textbf{4} (1995), no.~4, 639--657.

\bibitem{ReevesStillman}
Alyson~A. Reeves and Mike Stillman, \emph{Smoothness of the lexicographic
  point}, J. Algebraic Geom. \textbf{6} (1997), no.~2, 235--246.

\bibitem{ReevesSturmfels}
Alyson~A. Reeves and Bernd Sturmfels, \emph{A note on polynomial reduction}, J.
  Symbolic Comput. \textbf{16} (1993), no.~3, 273--277.

\bibitem{RobbianoBorderBasis}
Lorenzo Robbiano, \emph{On border basis and {G}r{\"o}bner basis schemes},
  Collect. Math. \textbf{60} (2009), no.~1, 11--25.

\bibitem{RoggeroTerracini}
Margherita Roggero and Lea Terracini, \emph{Ideals with an assigned initial
  ideals}, Int. Math. Forum \textbf{5} (2010), no.~53-56, 2731--2750.

\bibitem{Sabadini}
Irene Sabadini, \emph{On the {H}ilbert scheme of curves of degree {$d$} and
  genus {$\binom{d-3}{2}-1$}}, Matematiche (Catania) \textbf{55} (2000), no.~2,
  517--531 (2002). Dedicated to Silvio Greco on the occasion of his 60th
  birthday (Catania, 2001).

\bibitem{SchlesingerFootnote}
Enrico Schlesinger, \emph{Footnote to a paper by {R}. {H}artshorne: \lq\lq{O}n
  the connectedness of the {H}ilbert scheme of curves in {${\Bbb P}^3$}\rq\rq\
  [{C}omm.\ {A}lgebra {\bf 28} (2000), no.\ 12, 6059--6077; {MR}1808618
  (2002d:14003)]}, Comm. Algebra \textbf{28} (2000), no.~12, 6079--6083.
  Special issue in honor of Robin Hartshorne.

\bibitem{Sernesi}
Edoardo Sernesi, \emph{Deformations of algebraic schemes}, Grundlehren der
  Mathematischen Wissenschaften [Fundamental Principles of Mathematical
  Sciences], vol. 334, Springer-Verlag, Berlin, 2006.

\bibitem{Shafarevich1}
Igor~R. Shafarevich, \emph{Basic algebraic geometry. 1}, second ed.,
  Springer-Verlag, Berlin, 1994. Varieties in projective space. Translated from
  the 1988 Russian edition and with notes by Miles Reid.

\bibitem{Shafarevich2}
\bysame, \emph{Basic algebraic geometry. 2}, second ed., Springer-Verlag,
  Berlin, 1994. Schemes and complex manifolds. Translated from the 1988 Russian
  edition by Miles Reid.

\bibitem{ShermanPersistence}
Morgan Sherman, \emph{A local version of Gotzmann's persistence}, ArXiv
  e-prints (2007). Available at
  \href{http://arxiv.org/abs/0710.0186}{\texttt{arxiv.org/abs/0710.0186}}.

\bibitem{ShermanBorelFixed}
Morgan Sherman, \emph{On an extension of {G}alligo's theorem concerning the
  {B}orel-fixed points on the {H}ilbert scheme}, J. Algebra \textbf{318}
  (2007), no.~1, 47--67.

\bibitem{Stanley}
Richard~P. Stanley, \emph{Combinatorics and commutative algebra}, second ed.,
  Progress in Mathematics, vol.~41, Birkh{\"a}user Boston Inc., Boston, MA,
  1996.

\bibitem{StanleyEnumComb}
\bysame, \emph{Enumerative combinatorics. {V}ol. 1}, Cambridge Studies in
  Advanced Mathematics, vol.~49, Cambridge University Press, Cambridge, 1997.
  With a foreword by Gian-Carlo Rota. Corrected reprint of the 1986 original.

\bibitem{Strano}
Rosario Strano, \emph{Biliaison classes of curves in {${\bf P}^3$}}, Proc.
  Amer. Math. Soc. \textbf{132} (2004), no.~3, 649--658.

\bibitem{VainsencherAvritzer}
Israel Vainsencher and Dan Avritzer, \emph{Compactifying the space of elliptic
  quartic curves}, Complex projective geometry ({T}rieste, 1989/{B}ergen,
  1989), London Math. Soc. Lecture Note Ser., vol. 179, Cambridge Univ. Press,
  Cambridge, 1992, pp.~47--58.

\bibitem{VallaBellaterra}
Giuseppe Valla, \emph{Problems and results on {H}ilbert functions of graded
  algebras}, Six lectures on commutative algebra, Mod. Birkh{\"a}user Class.,
  Birkh{\"a}user Verlag, Basel, 2010, pp.~293--344.

\bibitem{Vanderbei}
Robert~J. Vanderbei, \emph{Linear programming}, second ed., International
  Series in Operations Research \& Management Science, 37, Kluwer Academic
  Publishers, Boston, MA, 2001. Foundations and extensions.

\end{thebibliography}
\nocite{*}
\references

\cleardoublepage
\printindex

\end{document}